\documentclass[twoside,a4paper,11pt]{book}

	\usepackage[utf8]{inputenc}
	\usepackage[T1]{fontenc}
	\usepackage{indentfirst}
	\usepackage{textcomp}
	\usepackage{yfonts}
	\usepackage[sc,osf]{mathpazo}
	\usepackage[english,french]{babel}
	
	\usepackage{amsmath,amsfonts,amssymb,amsopn,amscd,amsthm}
	\usepackage{esint,mathrsfs,units}
	\usepackage[bbgreekl]{mathbbol}
	\usepackage[enableskew]{youngtab}
	  
	\usepackage{graphicx}
	\usepackage[dvipsnames]{xcolor}
	\usepackage{pstricks,pst-plot,pst-3dplot,pst-vue3d,pst-grad}
	\usepackage{xy}  \xyoption{all}
	\usepackage{pb-diagram,pb-xy}
	\usepackage{rotating}

	\usepackage{enumerate}
	\usepackage{fancyhdr}
	\usepackage[colorlinks=true,citecolor=DarkOrchid,linkcolor=NavyBlue,naturalnames=true]{hyperref}
	
		\newcommand{\N}{\mathbb{N}}
		\newcommand{\Z}{\mathbb{Z}}
		\newcommand{\Q}{\mathbb{Q}}
		\newcommand{\R}{\mathbb{R}}
		\newcommand{\C}{\mathbb{C}}
		\newcommand{\CP}{\mathbb{C}\mathbb{P}}
		\renewcommand{\emptyset}{\varnothing}
		\newcommand{\E}{\mathrm{e}}
		\newcommand{\I}{\mathrm{i}}
		\newcommand{\eps}{\varepsilon}
		
		\newcommand{\scal}[2]{\left\langle #1\vphantom{#2}\,\right|
						      \left.#2 \vphantom{#1}\right\rangle}
		\renewcommand{\leq}{\leqslant}
		\renewcommand{\geq}{\geqslant}
		\newcommand{\lle}{\left[\!\left[}
		\newcommand{\rre}{\right]\!\right]}
		\newcommand{\vacuum}{\underline{\emptyset}}
		
		\newcommand{\alg}{\mathscr{A}}
		\newcommand{\blg}{\mathscr{B}}
		\newcommand{\clg}{\widetilde{\mathscr{A}}}
		\newcommand{\dlg}{\mathscr{D}}
		\newcommand{\For}{\mathbb{F}}
		\newcommand{\sym}{\mathfrak{S}}
		\newcommand{\wsym}{\mathfrak{W}}
		\newcommand{\Part}{\mathfrak{P}}
		\newcommand{\comp}{\mathfrak{C}}
		\newcommand{\IH}{\mathscr{H}}
		\newcommand{\Herm}{\mathrm{H}}
		\newcommand{\Unit}{\mathrm{U}}
		\newcommand{\GL}{\mathrm{GL}}
		\newcommand{\liegl}{\mathfrak{gl}}
		\newcommand{\GB}{\mathrm{B}}
		\newcommand{\Sp}{\mathrm{Sp}}
		\newcommand{\BSp}{\mathrm{BSp}}
		\newcommand{\Ind}{\mathrm{R}}
		\renewcommand{\hom}{\mathrm{Hom}}
		\newcommand{\hendo}{\mathrm{End}}
		\newcommand{\tr}{\mathrm{tr}\,}
		\newcommand{\proj}{\mathrm{pr}}
		\newcommand{\id}{\mathrm{id}}
		\newcommand{\Gal}{\mathfrak{G}}
		\newcommand{\mB}{\mathrm{B}}
		\newcommand{\dec}{\mathrm{Dec}}
		\newcommand{\IS}{\mathrm{IS}}
		\newcommand{\codim}{\mathrm{codim}\,}

		\newcommand{\proba}{\mathbb{P}}
		\newcommand{\esper}{\mathbb{E}}
		\newcommand{\Gel}{\mathbb{G}}
		\newcommand{\card}{\mathrm{card}\,}
		\newcommand{\leb}{\mathrm{L}}
		
		\newcommand{\ym}{\mathscr{Y}}
		\newcommand{\cym}{\mathscr{C}\!\mathscr{Y}}
		\newcommand{\obs}{\mathscr{O}}
		\newcommand{\imaj}{\mathrm{imaj}}
		\newcommand{\arra}{\mathcal{A}}
		\newcommand{\tilp}{\widetilde{p}}
		\newcommand{\orb}{\mathrm{orb}}
		
		\newcommand{\figcapt}[3]{\begin{figure}[ht] \begin{center} 
							{\footnotesize{#1}} 
							\caption[#3]{#2} \end{center} \end{figure}}
		\newcommand{\comment}[1]{}
		\newtheorem{theorem}{Théorème}[chapter]
		\newtheorem{proposition}[theorem]{Proposition}
		\newtheorem{corollary}[theorem]{Corollaire}
		\newtheorem{conjecture}[theorem]{Conjecture}
		\newtheorem{lemma}[theorem]{Lemme}
		\newtheorem{algorithm}[theorem]{Algorithme}
		\newtheorem{definition}[theorem]{Définition}
		\theoremstyle{remark}
		\newtheorem*{example}{Exemple}
		\newtheorem*{examples}{Exemples}
		\newtheorem*{remark}{Remarque}
\renewcommand{\cleardoublepage}{\clearpage}	

\begin{document}
\selectlanguage{french}
\title{Partitions aléatoires et théorie asymptotique des groupes symétriques, des algèbres d'Hecke et des groupes de Chevalley finis}
\author{Pierre-Loïc Méliot}
\date{18 novembre 2010}

\frontmatter

\pagestyle{empty}
\begin{center}
~\vspace{-1cm}~\\
\includegraphics{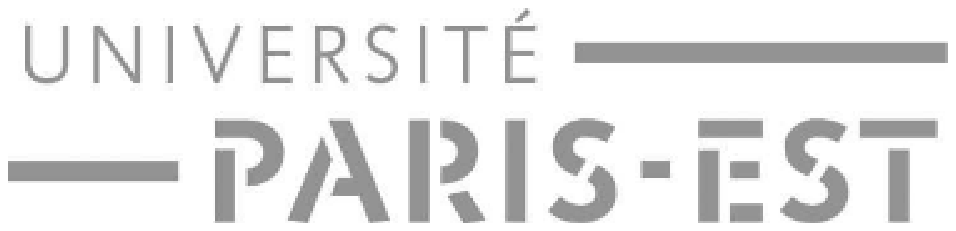}\\
\vspace{5mm}
{\LARGE\textsc{Thèse de doctorat}}\\
\vspace{2mm}
pour l'obtention du grade de\\
\vspace{2mm}
{\Large Docteur de l'Université Paris-Est}\\
\vspace{2mm}
{\large\textbf{Spécialité Mathématiques}}\\
\vspace{2mm}
\emph{au titre de l'\'Ecole Doctorale de Mathématiques et des Sciences \\et Techniques de l'Information et de la Communication.}\\
\vspace{1cm}
Présentée et soutenue publiquement par \\
\vspace{2mm}
{\Large Pierre-Loïc Méliot}\\
\vspace{2mm}
{le vendredi 17 décembre 2010.}\\
\vspace{1cm}
{\huge
\rule{16cm}{2pt}
\vspace{-2mm}\\
Partitions aléatoires et théorie asymptotique des groupes symétriques, des algèbres d'Hecke  et des groupes de Chevalley finis\\
\vspace{5mm}
\rule{16cm}{2pt}}
\\
\vspace{1cm}
{\large \textbf{Sous la direction de }}\\
\vspace{2mm}
{\large Philippe Biane,}\\
\vspace{1cm}
{\large \textbf{et devant le jury composé par }}\\
\vspace{2mm}
\begin{tabular}{ll}
\emph{Rapporteurs :}&Alexei Borodin,\\
&Cédric Lecouvey,\\
&\\
\emph{Examinateurs :}&Philippe Bougerol,\\
&Ashkan Nikeghbali,\\
&Jean-Yves Thibon,\\
&\\
\emph{Directeur :}&Philippe Biane.
\end{tabular}
\end{center}
\clearpage

\begin{quotation}
\phantom{~}$\,$\vspace{-1.5cm}

{\footnotesize\foreignlanguage{english}{
\noindent \textsc{Kit meanwhile had begun} to frequent the Applied Mechanics Institute. Since Prandtl's recent discovery of the boundary layer, things over there had been hopping, with intense inquiry into matters of lift and drag, powered flight poised like a new-feathered bird at the edge of history. Kit had not thought much about aerodynamics since his brainless sojourn in the Vibe embrace, when in the course of golfing parties out on Long Island he had become acquainted with the brambled guttie, a gutta-percha ball systematically roughened away from the perfectly spherical by molding little knobs all over the surface area. What he could not help noticing then, even though he was not all that crazy for the game, so inordinately populated by the likes of Scarsdale Vibe, was a particular mystery of flight---the undeniable lift of heart in seeing a struck ball---a tee shot especially---suddenly go into a steep ascent, an exhilarated denial of gravity you didn't have to be a golfer to appreciate. There being enough otherworldliness out on the links already. Finding himself more and more drawn to the microcosm on the other side of the B\"urgerstra\ss{}e, Kit soon understood that the brambling of the golf-ball surface had been a way to keep the boundary layer from detaching and falling apart into turbulence which would tend to drag the ball down, denying it its destiny in the sky. When he mentioned this in conversations at the saloons along the Brauweg frequented by engineering and physics students, some immediately suggested implications for the Earth, a brambled spheroid on the grand scale, in its passage through the \AE{}ther, being lifted not in the third dimension but on a euphoric world-line through Minkowski's ``four-dimensional physics.''\\
\phantom{~}$\,\,\,$``What happened to vectorism?'' Yashmeen teased.\\
\phantom{~}$\,\,\,$``There are vectors,'' Kit replied, ``and vectors. Over in Dr. Prandtl's shop, they're all straightforward lift and drift, velocity and so forth. You can draw pictures, of good old three-dimensional space if you like, or on the Complex plane, if Zhukovsky's Transformation is your glass of tea. Flights of arrows, teardrops. In Geheimrat Klein's shop, we were more used to expressing vectors without pictures, purely as an array of coefficients, no relation to anything physical, not even space itself, and writing them in any number of dimensions---according to Spectral Theory, up to infinity.''\\
\phantom{~}$\,\,\,$``And beyond,'' added G\"unther, nodding earnestly.
\\
~\\
 \noindent \textsc{In Hilbert's class} one day, she raised her hand. He twinkled at her to go ahead. ``\emph{Herr Geheimrat---}''\\
\phantom{~}$\,\,\,$``'Herr Professor' is good enough.''\\
\phantom{~}$\,\,\,$``The nontrivial zeroes of the $\zeta$-function...''\\
\phantom{~}$\,\,\,$``Ah.''\\
\phantom{~}$\,\,\,$She was trembling. She had not had much sleep. Hilbert had seen that sort of thing before, and rather a good deal of it since the turn of the century---since his own much-noted talk at the Sorbonne, he supposed, in which he had listed the outstanding problems in mathematics which would be addressed in the coming century, among them that of the zeroes of the $\zeta$-function.\\
\phantom{~}$\,\,\,$``Might they be correlated with eigenvalues of some Hermitian operator yet to be determined?''\\
\phantom{~}$\,\,\,$The twinkle, as some reported later, modulated to a steady pulsation. ``An intriguing suggestion, Fr\"aulein Halfcourt.'' Usually he addressed her as ``my child.'' ``Let us consider why this should be so.'' He peered, as if she were an apparition he was trying to see more clearly. ``Apart from eigenvalues, by their nature, being zeroes of \emph{some} equation,'' he prompted gently.\\
\phantom{~}$\,\,\,$``There is also this... spine of reality.'' Afterward she would remember she actually said ``\emph{R\"uckgrat von Wirklichkeit.}'' ``Though the members of a Hermitian may be complex, the eigenvalues are real. The entries on the main diagonal are real. The $\zeta$-function zeroes which lie along Real part = $\nicefrac{1}{2}$, are symmetrical about the real axis, and so...'' She hesitated. She had \emph{seen} it, for the moment, so clearly.\\
\phantom{~}$\,\,\,$``Let us apply some thought,'' said Hilbert. ``We will talk about this further.'' But she was to leave G\"ottingen shortly after this, and they would never have the chance to confer. As years passed, she would grow dim for Hilbert, her words those of an inner sprite too playful to frame a formal proposition, or to qualify as a fully habilitated Muse. And the idea itself would evolve into the celebrated Hilbert-P\'olya Conjecture.
}}
\end{quotation}
{\small $\qquad\qquad\qquad\qquad\qquad\qquad\qquad\qquad\qquad\qquad\qquad$ T. Pynchon, \emph{Against the Day}, p. 678-679, 2006.\phantom{\cite{Pyn06}}}\clearpage

\chapter*{Remerciements}
Mes premiers remerciements sont bien sûr adressés à mon directeur de thèse, Philippe Biane, qui m'a proposé il y a trois ans de travailler sur ce sujet. Il m'a laissé une liberté totale dans mes travaux, tout en prodiguant de précieux conseils. Le meilleur qu'il ait pu me donner fut sans doute de regarder des exemples de petite taille, et de faire des calculs explicites dans ces cas ;  à peu près tous les résultats de cette thèse ont été découverts à partir du moment où j'ai suivi ce simple précepte. Au-delà de l'aspect mathématique, j'ai eu un grand plaisir à travailler sous la direction d'une personne si agréable et joyeuse.\\

À l'autre bout de la frise chronologique, je remercie Alexei Borodin et Cédric Lecouvey qui ont accepté le rôle de rapporteur et qui ont pris la peine de lire ce mémoire et de m'aider à le corriger et à l'améliorer. Je remercie également Philippe Bougerol, Ashkan Nikeghbali et Jean-Yves Thibon, qui ont accepté de se joindre aux rapporteurs pour constituer mon jury de thèse, et de se déplacer pour assister à ma soutenance.
\bigskip

Au cours de ma thèse, l'équipe de combinatoire algébrique de Marne-La-Vallée m'a con\-stamment soutenu, et les séminaires du vendredi matin m'ont beaucoup appris. Les discussions qui s'y déroulent n'ont rien à envier à l'avalanche scientifico-baroque décrite dans le passage de \cite{Pyn06} cité ci-contre ; et de nombreuses idées me sont venues aux cours de ces échanges. Je remercie donc tout particulièrement Florent Hivert, Alain Lascoux, Jean-Christophe Novelli,  Jean-Yves Thibon, Nicolas Thiéry, et tous les autres participants de cette formidable assemblée. Je remercie également tous les thésards de l'équipe, qui ont grandement contribué à faire de ma thèse un moment agréable : Adrien Boussicault, Jean-Paul Bultel, Hayat Cheballah, Valentin Féray, Samuele Giraudo, Viviane Pons et Marc Sage. Un remerciement particulier doit être adressé à Valentin, avec qui j'ai coécrit mon premier article (\cite{FM10}). C'est lui qui m'a suggéré l'utilisation des observables de diagrammes pour l'étude asymptotique des diagrammes de Young, et qui m'a montré l'importance de l'algèbre des permutations partielles ; et ma thèse aurait sans doute été bien moins fructueuse sans ces suggestions.\\

Une partie de ma thèse qui n'est quasiment pas évoquée dans ce mémoire est l'écriture de programmes permettant de tirer au hasard des partitions selon les mesures de Plancherel et leurs diverses déformations, et de faire des calculs dans diverses algèbres combinatoires. Dans ce cadre, j'ai profité pleinement de l'émergence et du développement du système de calcul formel \texttt{sage} ; j'en remercie donc les développeurs, et plus particulièrement Florent Hivert et Nicolas Thiéry, qui m'ont aidé à corriger certains de mes programmes.\bigskip

Durant ces trois années, j'ai eu plusieurs fois la possibilité de voyager pour assister à des conférences et rencontrer d'autres chercheurs : c'est une autre partie de la thèse n'est pas relatée ici, et ce fut pourtant l'une des plus agréables, puisque j'ai eu la chance de découvrir (dans cet ordre) Montréal, Cargese, Z\"urich, San Francisco, Bielefeld, et en chacun de ces lieux des personnalités étonnantes et des mathématiques rafraichissantes (ou réciproquement). Je remercie ainsi K\"ursat Aker, Gérard Ben Arous, Marek Bozejko, Maciej Do{\l}ega, Ashkan Nikeghbali, Eric Nordenstam, Piotr \'Sniady et Anatoli Vershik pour les discussions que j'ai eues avec eux lors de ces voyages. Je remercie également infiniment Sylvie Cach et Line Fonfrede pour leur aide efficace dans l'organisation de ces voyages. J'ai également pu participer à quelques conférences ou séminaires tenus à Paris même ; concernant ceux-ci, je remercie Florent Benaych-Georges, qui m'a laissé exposer les résultats d'Okounkov (\cite{Oko00}) à l'Institut Henri Poincaré, et Paul Bourgade, avec qui j'ai eu d'intéressantes discussions.\bigskip

Finalement, mes derniers remerciements sont adressés à mes parents et à mes amis, pour leur soutien indéfectible. La plupart des idées mathématiques ne surviennent pas au cours de calculs apocalyptiques ; au contraire, elles émergent le plus souvent lorsque le mathématicien est au repos, et fait toute autre chose. Ainsi, les résultats présentés dans ce mémoire ont probablement  vu le jour au cours de la main de poker suivante
\begin{center}\psset{unit=1mm}
\pspicture(0,0)(100,30)
\rput(10,13){P}\rput(10,17){9$\diamondsuit\,$2$\spadesuit$}
\rput(90,17){K$\heartsuit\,$K$\clubsuit$}
\rput(90,13){PL}
\rput(50,1){Cl}\rput(50,5){A$\diamondsuit\,$A$\clubsuit$}
\rput(50,25){$\underline{9\heartsuit\,9\spadesuit\,7\diamondsuit}\,\,\,\,\underline{3\clubsuit}\,\,\,\,\underline{9\clubsuit}$}
\endpspicture
\end{center}
ou bien au détour d'un chemin corse,  à moins que ce ne fut en écoutant pour la onzième fois l'histoire du cosaque. Les héros de ces aventures se reconnaîtront sans doute, et ils méritent ces derniers remerciements. 
~\clearpage

\chapter{Avant-propos}
Si $G$ est un groupe fini et si $V$ est une représentation linéaire complexe et de dimension finie de $G$, la théorie de Schur permet de décomposer de manière unique $V$ en somme directe de modules irréductibles. Ainsi,
$$V \simeq_{\C G} \bigoplus_{\lambda \in \widehat{G}}\,n_{\lambda}V^{\lambda}\,,$$
où $\widehat{G}$ est l'ensemble (fini) des classes d'isomorphismes de $\C G$-modules simples, et où les $n_{\lambda}$ sont des entiers positifs ou nuls. Ce résultat bien connu peut être réinterprété de manière probabiliste : toute représentation (réductible et non réduite à $0$) $V$ de $G$ fournit une mesure de probabilité sur $\widehat{G}$ --- la \textbf{mesure de Plancherel} de $V$, voir la définition \ref{defplancherel} --- définie par 
$$\proba_{V}[\lambda \in \widehat{G}]=\frac{n_{\lambda}\,\dim V^{\lambda}}{\dim V}\,.$$
Considérons alors une famille croissante de groupes finis $(G_{n})_{n \in \N}$, et une famille <<~naturelle~>> de représentations $(V_{n})_{n \in \N}$ de ces groupes --- par exemple, on peut prendre $V_{n}=\C G_{n}$. La \textbf{théorie asymptotique de ces représentations} est l'ensemble des résultats qui permettent de répondre aux questions suivantes :\vspace{2mm}
\begin{enumerate}
\item A-t-on des résultats de convergence pour les variables aléatoires $\lambda \in \widehat{G}_{n}$ tirées suivant les mesures $\proba_{V_{n}}$ ? Par exemple, si $g$ est un élément d'un groupe $G_{n}$ et si $\lambda$ est tiré au hasard suivant la mesure $\proba_{V_{N \geq n}}$, que dire de la variable aléatoire 
$$\chi^{\lambda}(g)\,,$$
où $\chi^{\lambda}$ est le caractère irréductible associé à la classe d'isomorphismes de $\C G_{N}$-modules simples $\lambda$ ?\vspace{2mm}
\item Plus généralement, sans fixer de famille de représentations $(V_{n})_{n \in \N}$, comment exprimer simplement un caractère irréductible $\chi^{\lambda}(g)$ avec $\lambda \in \widehat{G}_{N}$ et $N$ grand, et $g$ fixé dans un groupe $G_{n \leq N}$ ? Quels sont les liens entre les représentations des <<~grands~>> groupes $G_{N}$, $N$ tendant vers l'infini, et les représentations de la limite inductive $G_{\infty}=\lim_{n \to \infty} G_{n}$ ?\vspace{2mm}
\item Les centres $Z(\C G_{n})$ des algèbres des groupes $G_{n}$ ont-ils des propriétés communes ? En particulier, les classes de conjugaison vérifient-elles des relations <<~génériques~>>, c'est-à-dire indépendantes de la taille $n$ du groupe ?\vspace{2mm}
\end{enumerate}
À partir des années 60, et suite au développement de l'analyse harmonique abstraite sur les groupes infinis et aux articles précurseurs \cite{FH59} et \cite{Tho64}, une théorie asymptotique des représentations a été démontrée pour la famille des \textbf{groupes symétriques} $(\sym_{n})_{n\in \N}$. Les éléments de $\widehat{\sym_{n}}$ sont les \textbf{partitions} de taille $n$ (voir le chapitre \ref{permutation}), c'est-à-dire les suites décroissantes d'entiers $\lambda=(\lambda_{1}\geq \cdots \geq \lambda_{r})$ telles que 
$$\lambda_{1}+\lambda_{2}+\cdots+\lambda_{r}=n\,\,;$$
par conséquent, toute représentation de $\sym_{n}$ fournit un modèle de \textbf{partition aléatoire}, \emph{i.e.}, une mesure de probabilité sur l'ensemble des partitions de taille $n$. Les travaux de Kerov et Vershik (\cite{KV77,KV81,Ker93a}) ont montré que sous les mesures de Plancherel des représentations régulières des groupes $\sym_{n}$, les partitions d'entiers étaient concentrées gaussiennement autour d'une <<~forme limite~>> (théorèmes \ref{firstasymptoticplancherel} et \ref{secondasymptoticplancherel}).\bigskip

Les motivations de la théorie asymptotique des représentations sont multiples, et au cours des deux dernières décennies, des liens ont été établis entre cette théorie et la théorie des matrices aléatoires (\emph{cf.} \cite{BOO00,Oko00}), le problème d'Ulam des plus longs sous-mots croissants dans un mot aléatoire (\cite{BDJ99,BDJ00,Joh01}), les probabilités libres (voir \cite{Bia98}), les pavages et les surfaces aléatoires (voir \cite{OR03}), et même la théorie de Gromov-Witten (\cite{Oko03b}) et les équations aux dérivées partielles hydrodynamiques (\cite{Ker93b,Ker99}). Retenons trois motivations essentielles :\vspace{2mm}
\begin{enumerate}
\item La mesure de Plancherel d'un $G$-module est évidemment importante pour l'\textbf{analyse harmonique abstraite} sur le groupe $G$. Ainsi, en connaître les propriétés apporte un point de vue nouveau (probabiliste) sur la théorie des représentations de $G$. Par exemple, dans le cas du groupe symétrique, ce nouveau paradigme a poussé Ivanov, Kerov et Olshanski à construire l'algèbre des permutations partielles (\cite{IK99,IO02}), que l'on peut aussi interpréter comme algèbre d'observables des partitions d'entiers (\emph{cf.} les chapitres \ref{tool} et \ref{badbeat}), et qui apporte une explication simple à de nombreux résultats sur la structure des (centres des) algèbres des groupes symétriques. Dans ce mémoire, outre le groupe symétrique, nous serons amené à étudier les groupes de Lie classiques sur les corps finis ; leur théorie des représentations est nettement plus complexe (voir le chapitre \ref{general}), et on peut espérer que l'approche probabiliste simplifie sa compréhension, ou du moins apporte de nouvelles idées dans ce domaine (voir par exemple les résultats du chapitre \ref{badbeat}, et le paragraphe \ref{polyobs}). \vspace{2mm}
\item Les mesures de Plancherel des représentations des groupes classiques ont souvent des interprétations combinatoires. En particulier, en utilisant l'algorithme RSK (\emph{cf.} \S\ref{ulamrsk}), Logan, Shepp, Kerov et Vershik (\cite{LS77,KV77}) ont ramené le problème des \textbf{plus longs sous-mots croissants} dans une permutation aléatoire à l'étude asymptotique des mesures de Plancherel des représentations régulières des groupes symétriques, et ainsi résolu le problème d'Ulam. On peut donc espérer résoudre des problèmes purement combinatoires à l'aide de la théorie asymptotique des représentations ; par exemple, les résultats que nous avons obtenus sur les $q$-mesures de Plancherel peuvent être énoncés en termes de statistiques des permutations aléatoires distribuées suivant un potentiel proportionnel à leur indice majeur (corollaires \ref{firstknuth} et \ref{secondknuth}). \vspace{2mm}
\item Les méthodes employées et les résultats obtenus en théorie asymptotique des représentations sont à rapprocher des méthodes et des résultats de la théorie des \textbf{matrices aléatoires}. Ces similarités seront évoquées en détail dans les chapitres \ref{determinantal} et \ref{matrix} ; il est très difficile d'établir un lien direct entre ces deux sujets (la meilleure tentative est probablement l'article d'Okounkov \cite{Oko00}), mais on retrouve les mêmes lois limites, les mêmes processus ponctuels déterminantaux, les mêmes processus gaussiens, \emph{etc.} de part et d'autre. Ainsi, les modèles de partitions aléatoires issus de la théorie (asymptotique) des représentations doivent être envisagés comme des analogues discrets des modèles de matrices aléatoires, et une avancée dans l'un des domaines peut donner de nouvelles idées pour résoudre les problèmes de l'autre.
\end{enumerate}

Au cours de cette thèse, nous avons établi des résultats analogues aux résultats asymptotiques classiques pour les mesures de Plancherel des représentations régulières des groupes symétriques, mais pour d'autres familles de groupes finis, et pour d'autres familles de représentations des groupes symétriques.  Compte tenu de la classification des groupes finis simples, les familles intéressantes qui viennent immédiatement après la famille $(\sym_{n})_{n \in \N}$ sont les familles de groupes de Chevalley classiques
$$(\GL(n,\For_{q}))_{n \in \N},\,\,\,(\mathrm{U}(n,\For_{q^{2}}))_{n \in \N},\,\,\,(\Sp(2n,\For_{q}))_{n \in \N},\,\,\,(\mathrm{O}(n,\For_{q}))_{n \in \N},\,\,\,\text{etc.}$$
Ainsi, nous nous sommes essentiellement penchés sur le cas des groupes finis de matrices $\GL(n,\For_{q})$, et nous avons étudié des familles de représentations de ces groupes qu'on peut voir comme $q$-déformations des représentations régulières $\sym_{n}\curvearrowright \C\sym_{n}$ (chapitre \ref{general} à \ref{qplancherelmeasure}). Plus précisément, nous avons établi des résultats de concentration gaussienne pour les \textbf{q-mesures de Plancherel} (théorèmes \ref{firstasymptoticqplancherel}, \ref{secondasymptoticqplancherel} et \ref{thirdasymptoticqplancherel}), qui sont les mesures de probabilité associées aux représentations
$$\GL(n,\For_{q}) \curvearrowright \C[\mathcal{F}(n,\For_{q})]\,,$$
où $\mathcal{F}(n,\For_{q})=\GL(n,\For_{q})/\GB(n,\For_{q})$ désigne la variété des drapeaux complets de taille $n$ sur le corps $\For_{q}$ (voir \cite{FM10,Mel10a}). En réalité, on retrouve ces phénomènes de concentration gaussienne dans un contexte beaucoup plus général, à savoir celui des modules sur un groupe de Chevalley fini obtenus par induction parabolique à partir d'un caractère cuspidal d'un sous-groupe de Lévi rationnel (chapitre \ref{arikikoike}). On a en particulier établi un résultat de concentration gaussienne pour l'analogue en type B du cas précédemment exposé (théorèmes \ref{asymptoticbqplancherel} et \ref{bqmix}), c'est-à-dire pour les représentations
$$\Sp(2n,\For_{q})\curvearrowright \C[\mathcal{F}^{\mB}(n,\For_{q})]\,,$$
où $\mathcal{F}^{\mB}(n,\For_{q})=\Sp(2n,\For_{q})/\BSp(2n,\For_{q})$ désigne la variété des drapeaux totalement anisotropes complets de taille $n$ sur le corps $\For_{q}$ (dans un espace symplectique de dimension $2n$). \bigskip

Incidemment, cette théorie asymptotique des représentations des groupes de Chevalley finis peut être réinterprétée en termes de représentations des groupes symétriques, ou plus précisément en termes de représentations \textbf{des algèbres d'Hecke} des groupes symétriques (et des autres groupes de Coxeter classiques). C'est ce point de vue qui nous a permis d'adapter les arguments de Kerov et Vershik, et d'établir nos résultats de concentration gaussienne. Ainsi, les modèles algébriques présentés ci-dessus constituent de nouveau des modèles de partitions aléatoires. Les partitions d'entiers sont des objets sont de nature géométrique planaire, et elles peuvent être étudiées à l'aide d'une \textbf{algèbre d'observables} modelée sur l'algèbre des fonctions symétriques (\cite{IO02}; chapitre \ref{tool}). Cette algèbre $\obs$ est l'analogue pour l'ensemble $\ym$ de toutes les partitions d'entiers de l'algèbre des fonctions rationnelles pour une variété algébrique ; ainsi, ses éléments sont les fonctions <<~polynomiales~>> des partitions. Elle permet de traiter avec le même formalisme :\vspace{2mm}
\begin{itemize}
\item[-] les représentations régulières des groupes symétriques (chapitre \ref{plancherel}),\vspace{2mm}
\item[-] les représentations des groupes $\GL(n,\For_{q})$ sur les modules $\C[\mathcal{F}(n,\For_{q})]$ (chapitre \ref{qplancherelmeasure}),\vspace{2mm}
\item[-] les représentations de Schur-Weyl des groupes symétriques sur des produits tensoriels (chapitre \ref{schurweylmeasure}),\vspace{2mm}
\item[-] et les modèles de Gelfand des groupes symétriques (chapitre \ref{gelfandmeasure}).\vspace{2mm}
\end{itemize}
Concernant les représentations de Schur-Weyl, nous sommes parvenus à préciser les résultats de P. Biane (\cite{Bia01a}), et nous avons montré que le théorème central limite de Kerov était également valable dans ce contexte (\emph{cf.} \cite{Mel10b} et le chapitre \ref{schurweylmeasure}), à une translation près le long de l'axe des abscisses. Ce résultat de concentration gaussienne des formes des diagrammes de Young est aussi valable pour les mesures de Gelfand des groupes symétriques (\emph{cf.} \cite{Mel10} et le chapitre \ref{gelfandmeasure}), et possède ainsi un caractère universel.  \bigskip

Nous avons mentionné ci-dessus que l'outil principal de notre étude asymptotique était l'algèbre $\obs$ des observables, dont les espérances jouent un rôle similaire pour les partitions aléatoires au rôle joué par les moments pour des variables aléatoires réelles. Or, l'algèbre $\obs$ se révèle être isomorphe à une sous-algèbre commutative de l'\textbf{algèbre d'Ivanov-Kerov des permutations partielles}, que l'on peut voir comme une limite projective des algèbres des groupes symétriques. Cette algèbre permet de démontrer des \textbf{identités génériques} dans les algèbres des groupes symétriques, c'est-à-dire des identités indépendantes de la taille $n$ du groupe symétrique. Nous avons tenté de généraliser cette construction au cas d'autres algèbres, notamment les algèbres d'Hecke $\IH_{q}(\sym_{n})$ et les algèbres des groupes linéaires finis $\GL(n,\For_{q})$ ; en particulier, nous avons démontré l'analogue pour les centres des algèbres d'Hecke d'un théorème de Farahat et Higman (\emph{cf.} \cite{Mel10c}). Toutes ces constructions rentrent dans le cadre extrêmement général et quelque peu bourbakiste des \textbf{fibrés de semi-groupes par des semi-treillis}. Dans ce même contexte, on peut formaliser des problèmes combinatoires sur les permutations, par exemple le problème des nombres de Hurwitz.\medskip
\begin{center}
$\ast\quad\ast\quad\ast$
\end{center}
\bigskip

La première partie de ce mémoire rappelle les résultats connus en théorie asymptotique des représentations pour les groupes symétriques. Nous avons pris la liberté de rédiger un exposé complet de cette théorie\footnote{Un exposé semblable est donné par A. Hora dans \cite{Hora07}.} ; ainsi, le lecteur pourra comprendre l'origine des outils employés pour l'étude asymptotique des algèbres d'Hecke et des groupes linéaires finis, et il aura une idée de l'étendue des résultats qu'on peut espérer obtenir dans ce nouveau cadre. Les chapitres \ref{permutation} à \ref{matrix} ne contiennent donc aucun résultat nouveau, mais leur lecture rendra nettement plus naturels les raisonnements des chapitres ultérieurs.\bigskip

La seconde partie du mémoire est consacrée à la théorie asymptotique des représentations des groupes de Chevalley finis sur leurs variétés de drapeaux. La combinatoire des représentations de ces groupes est hautement non triviale (notamment pour un lecteur probabiliste) ; nous rappelons donc cette théorie dans le chapitre \ref{general}, et nous y traitons en détail le cas des groupes $\GL(n,\For_{q})$. Les chapitres \ref{iwahori} à \ref{arikikoike} sont consacrés au coeur du problème, c'est-à-dire l'étude asymptotique des $q$-mesures de Plancherel ; en particulier, on y démontre les résultats asymptotiques \ref{firstasymptoticqplancherel}, \ref{secondasymptoticqplancherel}, \ref{thirdasymptoticqplancherel}, \ref{asymptoticbqplancherel} et \ref{bqmix}, qui sont les analogues des résultats de Kerov et Vershik pour le groupe symétrique. \bigskip

Dans la troisième partie, nous adaptons les outils utilisés pour l'étude des mesures de Plancherel et des $q$-mesures de Plancherel au cas des mesures de Schur-Weyl et des mesures de Gelfand. Ainsi, le chapitre \ref{schurweylmeasure} est consacré aux mesures de Schur-Weyl, qui en fonction des valeurs de leurs paramètres exhibent trois régimes asymptotiques distincts ; nous démontrons en particulier les théorèmes \ref{schurweylplus} et \ref{schurweylminus}. Dans le chapitre \ref{gelfandmeasure}, nous étudions une autre famille de représentations des groupes symétriques, les modèles de Gelfand ; leurs propriétés combinatoires sont liées à l'énumération des racines carrées dans les groupes symétriques. 
\bigskip

Enfin, la dernière partie est consacrée aux généralisations de la construction d'Ivanov-Kerov (\cite{IK99}) et du théorème de Farahat et Higman (\cite{FH59}). Dans le chapitre \ref{badbeat}, on démontre un analogue de ce théorème pour les éléments de Geck-Rouquier dans les centres des algèbres d'Iwahori-Hecke de type A (\cite{Mel10c}). Dans le chapitre \ref{bundle}, on montre que les constructions à la Ivanov-Kerov rentrent dans le cadre des fibrés de semi-groupes, et on applique ces notions au problème des nombres de Hurwitz. Il est vraisemblable qu'un théorème de type Farahat-Higman soit valable dans le contexte des centres des algèbres des groupes $\GL(n,\For_{q})$ ; on explique cette conjecture à la fin du chapitre \ref{bundle}. Le schéma qui suit résume le plan que nous venons de dresser, et indique les prérequis logiques à la lecture de chaque chapitre.

\begin{center}\psset{unit=1mm,arrowsize=3pt 2}\pspicture(-50,-175)(100,0)
\psframe(-40,-14)(40,0)
\rput(0,-5){\textbf{Chapitre 1.} Permutations, partitions,}
\rput(0,-9){représentations du groupe symétrique.}
\psline{->}(0,-14)(0,-25)
\psframe(-41,-35)(41,-25)
\rput(0,-30){\textbf{Chapitre 2.} Observables de diagrammes.}
\psline{->}(0,-35)(0,-60)(-20,-60)(-20,-65)
\psline{->}(0,-60)(24,-60)(24,-140)(45,-140)
\psline(15,-70)(24,-70)
\psframe(-55,-75)(15,-65)
\rput(-20,-70){\textbf{Chapitre 3.} Mesure de Plancherel.}
\psframe(-58,-99)(18,-85)
\psline{->}(-20,-75)(-20,-85)
\rput(-20,-90){\textbf{Chapitre 4.} Groupe symétrique infini,}
\rput(-20,-94){processus ponctuels déterminantaux.}
\psframe(-45,-124)(5,-110)
\psline{->}(-20,-99)(-20,-110)
\rput(-20,-115){\textbf{Chapitre 5.} Matrices et}
\rput(-20,-119){permutations aléatoires.}
\psframe(37,-55)(103,-45)
\rput(70,-50){\textbf{Chapitre 6}. Groupe linéaire fini.}
\psline{->}(68,-55)(68,-65)
\psframe(37,-75)(99,-65)
\psline(37,-50)(2,-50)
\psline{->}(-2,-50)(-64,-50)(-64,-170)(-52,-170)
\rput(68,-70){\textbf{Chapitre 7}. Algèbres d'Hecke.}
\psline(99,-70)(113,-70)
\psline{->}(68,-75)(68,-85)
\psframe(29,-95)(107,-85)
\psline{->}(24,-90)(29,-90)
\rput(68,-90){\textbf{Chapitre 8}. $q$-mesure de Plancherel.}
\psframe(42,-119)(94,-105)
\psline{->}(68,-95)(68,-105)
\rput(68,-110){\textbf{Chapitre 9}. Mesures}
\rput(68,-114){d'induction parabolique.}
\psframe(45,-149)(91,-135)
\rput(-19,-140){\textbf{Chapitre 10}. Mesures}
\rput(-19,-144){de Schur-Weyl.}
\psline{->}(24,-140)(3,-140)
\psframe(-41,-149)(3,-135)
\rput(68,-140){\textbf{Chapitre 11}. Mesures}
\rput(68,-144){de Gelfand.}
\psline{->}(40,-5)(113,-5)(113,-170)(92,-170)
\psline{->}(44,-170)(31,-170)
\psframe(44,-179)(92,-165)
\rput(68,-170){\textbf{Chapitre 12}. Algèbre}
\rput(68,-174){d'Hecke-Ivanov-Kerov.}
\psframe(-52,-179)(31,-165)
\rput(-10,-170){\textbf{Chapitre 13}. Fibrés de semi-groupes}
\rput(-10,-174){et limites projectives d'algèbres de groupes.}
\endpspicture
\psset{arrowsize=1.5pt 2}
\end{center}

\chapter*{Foreword}
\setcounter{footnote}{0}
\foreignlanguage{english}{If $G$ is a finite group and $V$ is a finite-dimensional complex linear representation of $G$, the Schur theory allows to decompose $V$ in a unique way as a direct sum of irreducible modules. Hence, 
$$V \simeq_{\C G} \bigoplus_{\lambda \in \widehat{G}}\,n_{\lambda}V^{\lambda}\,,$$
where $\widehat{G}$ is the set of isomorphism classes of simple $\C G$-modules, and the $n_{\lambda}$'s are non-negative integers. This well-known result can be restated with a probabilistic point of view: any (reducible and non-zero) representation $V$ of $G$ yields a probability measure on $\widehat{G}$ --- the so-called \textbf{Plancherel measure} of $V$, see Definition \ref{defplancherel} --- that is defined by
$$\proba_{V}[\lambda \in \widehat{G}]=\frac{n_{\lambda}\,\dim V^{\lambda}}{\dim V}\,.$$
Then, let us consider an increasing family of groups $(G_{n})_{n \in \N}$, and a ``natural'' family of representations of these groups $(V_{n})_{n \in \N}$ --- for instance, one can take $V_{n}=\C G_{n}$. The \textbf{asymptotic theory of these representations} is the set of results related to the following questions:\vspace{2mm}
\begin{enumerate}
\item Does one observe convergence phenomena for the random variables $\lambda \in \widehat{G}_{n}$ picked according to the measures $\proba_{V_{n}}$? For instance, if $g$ is in some group $G_{n}$ and $\lambda$ is chosen randomly according to the distribution $\proba_{V_{N \geq n}}$, how does the random variable
$$\chi^{\lambda}(g)$$
behave, where $\chi^{\lambda}$ is the irreducible character associated to the isomorphism class of simple $\C G_{N}$-modules $\lambda$?
\item More generally and without fixing a family of representations $(V_{n})_{n \in \N}$, can one exprime in a simple way the irreducible character value $\chi^{\lambda}(g)$, with $g$ fixed in some group $G_{n}$, and $\lambda \in \widehat{G}_{N}$ with $N \geq n$? Can one relate the representations of the ``big'' groups $G_{N}$ where $N$ goes to infinity, and the representations of the inductive limit $G_{\infty}=\varinjlim_{n \to \infty}G_{n}$?
\item Do the centers $Z(\C G_{n})$ of the group algebras of the $G_{n}$'s share some properties? In particular, do the conjugacy classes satisfy some ``generic'' identities, that is to say, identities that do not depend on the size $n$ of the group?
\end{enumerate}
\noindent From the 60s, and following the development of abstract harmonic analysis on infinite groups and the precursor papers \cite{FH59} and \cite{Tho64}, an asymptotic representation theory has been proved for the family of \textbf{symmetric groups} $(\sym_{n})_{n \in \N}$. The elements of $\widehat{\sym_{n}}$ are the \textbf{partitions} of size $n$ (see Chapter \ref{permutation}), that is to say, the non-increasing sequences of integers $\lambda=(\lambda_{1}\geq \cdots \geq \lambda_{r})$ such that 
$$\lambda_{1}+\lambda_{2}+\cdots+\lambda_{r}=n\,\,;$$ 
consequently, any representation of $\sym_{n}$ yields a model of \textbf{random partition}, \emph{i.e.}, a probability measure on the set of partitions of size $n$. Kerov and Vershik (\cite{KV77,KV81,Ker93a}) have shown that under the Plancherel measures of the regular representations of the groups $\sym_{n}$, the partitions exhibit a gaussian concentration around a ``limit shape'' (Theorems \ref{firstasymptoticplancherel} and \ref{secondasymptoticplancherel}). }
\bigskip

\foreignlanguage{english}{The motivations of asymptotic representation theory are numerous, and during the two last decades, connections have been established between this theory and random matrix theory (\emph{cf.} \cite{BOO00,Oko00}), Ulam's problem of longest increasing subsequences in random words (\cite{BDJ99,BDJ00,Joh01}), free probability (see \cite{Bia98}), random tilings and random surfaces (see \cite{OR03}), and even Gromov-Witten theory (\cite{Oko03b}) and hydrodynamic partial differential equations (\cite{Ker93b,Ker99}). Let us put emphasis on three essential motivations:\vspace{2mm}
\begin{enumerate}
\item The Plancherel measure of a $G$-module is obviously important in the setting of \textbf{abstract harmonic analysis} on the group $G$. Thus, a good knowledge of its properties provides a new (probabilistic) point of view on the representation theory of $G$. For example, in the case of the symmetric group, this new paradigm lead Ivanov, Kerov and Olshanski to construct the algebra of partial permutations (\cite{IK99,IO02}), that can also be interpreted as an algebra of observables of integer partitions (\emph{cf.} Chapters \ref{tool} and \ref{badbeat}), and that provides simple explanations to numerous results on the structure of the (centers of the)  algebras of the symmetric groups. In this thesis, in addition to the symmetric group, we shall study the classical Lie groups over finite fields; their representation theory is substantially more complex (see Chapter \ref{general}), and one can hope that the probabilist approach will simplify its understanding, or at least bring new ideas in this domain (see for instance the results of Chapter \ref{badbeat}, and Paragraph \ref{polyobs}). \vspace{2mm}
\item The Plancherel measures of the representations of the classical groups often have combinatorial interpretations. In paticular, by using the RSK algorithm (\emph{cf.} \S\ref{ulamrsk}), Logan, Shepp, Kerov and Vershik (\cite{LS77,KV77}) have reduced the problem of \textbf{longest increasing subsequences} in a random permutation to the asymptotic study of the Plancherel measures of the regular representations of the symmetric groups, thus solving Ulam's problem. So, one can hope to solve pure combinatorial problems thanks to asymptotic representation theory; for instance, the results that we have obtained for the $q$-Plancherel measures can be stated in terms of statistics of random permutations distributed according to a potential proportional to their major index (Corollaries \ref{firstknuth} and \ref{secondknuth}).\vspace{2mm}
\item The methods and the results of asymptotic representation theory should be related to those of \textbf{random matrix theory}. These similarities will be explained in detail in Chapters \ref{determinantal} and \ref{matrix}; it is very difficult to establish a direct link between these two subjects (the best attempt is probably Okounkov's article \cite{Oko00}), but one observes the same limiting distributions, the same determinantal point proceses, the same gaussian processes, \emph{etc.} on both sides. Hence, the models of random partitions stemming from (asymptotic) representation theory should be thought of as discrete analogues of the models of random matrices, and a breakthrough in one domain should bring new ideas in order to solve problems on the other side.
\end{enumerate}}
\begin{center}
$\ast\quad\ast\quad\ast$
\end{center}
\bigskip

\foreignlanguage{english}{
During this thesis, we have established results that are analog to the classical asymptotic results regarding the Plancherel measures of the regular representations of the symmetric groups, but for other families of finite groups, and for other families of representations of the symmetric groups. According to the classification of finite simple groups, the next interesting examples after the family $(\sym_{n})_{n \in \N}$ are the families of Chevalley groups of classical type:
$$(\GL(n,\For_{q}))_{n \in \N},\,\,\,(\mathrm{U}(n,\For_{q^{2}}))_{n \in \N},\,\,\,(\Sp(2n,\For_{q}))_{n \in \N},\,\,\,(\mathrm{O}(n,\For_{q}))_{n \in \N},\,\,\,\text{etc.}$$
We have focused mainly on the case of groups of matrices $\GL(n,\For_{q})$, and we have studied families of representations of these groups that can be seen as $q$-deformations of the (left) regular representations $\sym_{n} \curvearrowright \C\sym_{n}$ (Chapters \ref{general} to \ref{qplancherelmeasure}). More precisely, we have established the gaussian concentration of the \textbf{q-Plancherel measures} (Theorems \ref{firstasymptoticqplancherel}, \ref{secondasymptoticqplancherel} and \ref{thirdasymptoticqplancherel}), which are the probability measures related to the representations
$$\GL(n,\For_{q}) \curvearrowright \C[\mathcal{F}(n,\For_{q})]\,,$$
where $\mathcal{F}(n,\For_{q})=\GL(n,\For_{q})/\GB(n,\For_{q})$ is the variety of complete flags of size $n$ over the finite field $\For_{q}$ (see \cite{FM10,Mel10a}). In fact, these phenomena of gaussian concentration occur in a much broader context, namely, for modules over finite Chevalley groups that are obtained by parabolic induction from a cuspidal character of a rational Levi subgroup (Chapter \ref{arikikoike}). For instance, we have shown that gaussian concentration holds for the analogues in type B of the $q$-Plancherel measures (Theorems \ref{asymptoticbqplancherel} and \ref{bqmix}), \emph{i.e.}, for the representations  
$$\Sp(2n,\For_{q})\curvearrowright \C[\mathcal{F}^{\mB}(n,\For_{q})]\,,$$
where $\mathcal{F}^{\mB}(n,\For_{q})=\Sp(2n,\For_{q})/\BSp(2n,\For_{q})$ is the variety of totally anisotropic complete flags of size $n$ over $\For_{q}$ (in a symplectic vector space of dimension $2n$).}\bigskip

\foreignlanguage{english}{Incidentally, this asymptotic representation theory for finite Chevalley groups can be restated for the representations of the symmetric groups, or more precisely, for the representations of the \textbf{Hecke algebras} of the symmetric groups (and the other classical Coxeter groups). We have adopted this point of view in order to adapt the arguments of Kerov and Vershik, and to prove our results of gaussian concentration. Hence, the aforementioned algebraic models become again models of random partitions. The integer partitions are objects of planar geometric nature, and they can be studied thanks to an \textbf{algebra of observables} that is built on the model of the algebra of symmetric functions (\cite{IO02}; Chapter \ref{tool}). This algebra is the analog for the set $\ym$ of all integer partitions of the algebra of rational functions on an algebraic manifold; hence, its elements are the ``polynomial'' functions of partitions. It enables to treat with the same formalism:\vspace{2mm}
\begin{itemize}
\item[-] the regular representations of the symmetric groups (Chapter \ref{plancherel}),\vspace{2mm}
\item[-] the representations of the groups $\GL(n,\For_{q})$ on the modules $\C[\mathcal{F}(n,\For_{q})]$ (Chapter \ref{qplancherelmeasure}),\vspace{2mm}
\item[-] the Schur-Weyl representations of symmetric groups on tensor products (Chapter \ref{schurweylmeasure}),\vspace{2mm}
\item[-] and the Gelfand models of symmetric groups (Chapter \ref{gelfandmeasure}).\vspace{2mm}
\end{itemize}
Regarding the Schur-Weyl representations, we were able to precise the results of P. Biane (\cite{Bia01a}), and we have shown that Kerov's central limit theorem holds also in this setting (\emph{cf.} \cite{Mel10b} and Chapter \ref{schurweylmeasure}), up to a translation along the $x$-axis. This result of gaussian concentration of the shapes of the Young diagrams holds also for the Gelfand measures of the symmetric groups (\emph{cf.} \cite{Mel10} and Chapter \ref{gelfandmeasure}), so it has a character of universality.}\bigskip

\foreignlanguage{english}{
As mentioned previously, the main tool of our asymptotic study is the algebra $\obs$ of observables, whose expectations play a similar role for random partitions as the moments for real random variables. It turns out that $\obs$ is isomorphic to a subalgebra of the \textbf{Ivanov-Kerov algebra of partial permutations}, that can be seen as a projective limit of the group algebras of the symmetric groups. This algebra can be used to show \textbf{generic identities} in the symmetric group algebras, that is to say, identities that are independent of the size $n$ of the symmetric group. We have tried to generalize this construction to the case of other algebras, \emph{e.g.} the Hecke algebras $\IH_{q}(\sym_{n})$ and the group algebras of the finite linear groups $\GL(n,\For_{q})$; in particular, we have shown an analogue of a theorem of Farahat and Higman for the centers of the Hecke algebras (\emph{cf.} \cite{Mel10c}). All these constructions can be understood in the extremely general and Bourbaki-like setting of \textbf{semilattice bundles over semigroups}. In the same context, one can formalize combinatoric problems on permutations, for instance the problem of Hurwitz numbers.
\begin{center}
$\ast\quad\ast\quad\ast$
\end{center}} 
\bigskip

\foreignlanguage{english}{
In the first part of this thesis, we recall the results that are known in asymptotic representation theory of symmetric groups. We have chosen to write a complete survey\footnote{A similar survey is proposed by A. Hora in \cite{Hora07}.} of this theory, so that the reader can understand the origin of the tools used in the asymptotic representation theory of Hecke algebras and finite linear groups. He should also have a clear idea of the kind of results that one can hope to establish in this new setting. Hence, Chapters \ref{permutation} to \ref{matrix} do not contain any new result, but they are intented to make the arguments of the ulterior chapters much more natural.}\bigskip

\foreignlanguage{english}{
The second part of this thesis is devoted to the asymptotic theory of the representations of the finite Chevalley groups on their flag varieties. The combinatoric properties of the representations of these groups are highly non trivial (especially for a probabilist reader), so we recall this theory in Chapter \ref{general}, and we treat in detail the case of the groups $\GL(n,\For_{q})$. Chapters \ref{iwahori} to \ref{arikikoike} are devoted to the heart of the problem, that is to say the asymptotic study of the $q$-Plancherel measures; in particular, we prove the asymptotic results  \ref{firstasymptoticqplancherel}, \ref{secondasymptoticqplancherel}, \ref{thirdasymptoticqplancherel}, \ref{asymptoticbqplancherel} and \ref{bqmix}, which are the analogues of the results of Kerov and Vershik for the symmetric groups. }\bigskip

\foreignlanguage{english}{In the third part, we adapt the tools used for the study of the Plancherel and $q$-Plancherel measures to the case of Schur-Weyl and Gelfand measures. Thus, Chapter \ref{schurweylmeasure} is devoted to Schur-Weyl measures, which according to the values of their parameters exhibit three distinct asymptotic behaviours; we prove notably Theorems \ref{schurweylplus} and \ref{schurweylminus}. In Chapter \ref{gelfandmeasure}, we study another family of representations of the symmetric groups, the Gelfand models; their combinatorics are related to the enumeration of square roots in the symmetric groups.}\bigskip

\foreignlanguage{english}{
Finally, the last part of this thesis presents some generalizations of the construction of Ivanov and Kerov (\cite{IK99}) and of Farahat-Higman's theorem (\cite{FH59}). In Chapter \ref{badbeat}, we show an analogue of this theorem for Geck-Rouquier elements in the centers of the Iwahori-Hecke algebras of type A (\cite{Mel10c}). In Chapter \ref{bundle}, we see that the constructions \emph{à la} Ivanov-Kerov belong to the realm of bundles over semigroups, and we apply these notions to the problem of Hurwitz numbers. It is likely that a Farahat-Higman theorem holds in the context of the centers of the algebras of the groups $\GL(n,\For_{q})$; we explain this conjecture at the end of Chapter \ref{bundle}. The scheme hereafter summarizes the plane we have just exposed, and it indicates the logical prerequisites to the reading of each chapter.}

\begin{center}\psset{unit=1mm,arrowsize=3pt 2}\pspicture(-50,-180)(100,5)
\psframe(-40,-14)(40,0)
\rput(0,-5){\textbf{Chapter 1.} Permutations, partitions,}
\rput(0,-9){representations of the symmetric group.}
\psline{->}(0,-14)(0,-25)
\psframe(-38,-35)(38,-25)
\rput(0,-30){\textbf{Chapter 2.} Observables of diagrams.}
\psline{->}(0,-35)(0,-60)(-20,-60)(-20,-65)
\psline{->}(0,-60)(24,-60)(24,-140)(36,-140)
\psline{->}(24,-140)(12,-140)
\psframe(-56,-145)(12,-135)
\psline(13,-70)(24,-70)
\psframe(-53,-75)(13,-65)
\rput(-20,-70){\textbf{Chapter 3.} Plancherel measure.}
\psframe(-58,-99)(18,-85)
\psline{->}(-20,-75)(-20,-85)
\rput(-20,-90){\textbf{Chapter 4.} Infinite symmetric group,}
\rput(-20,-94){determinantal point processes.}
\psframe(-50,-124)(10,-110)
\psline{->}(-20,-99)(-20,-110)
\rput(-20,-115){\textbf{Chapter 5.} Random matrices}
\rput(-20,-119){and ramdom permutations.}
\psframe(38,-55)(102,-45)
\rput(70,-50){\textbf{Chapter 6}. Finite linear group.}
\psline{->}(68,-55)(68,-65)
\psframe(39,-75)(97,-65)
\psline(38,-50)(2,-50)
\psline{->}(-2,-50)(-64,-50)(-64,-170)(-52,-170)
\rput(68,-70){\textbf{Chapter 7}. Hecke algebras.}
\psline(97,-70)(113,-70)
\psline{->}(68,-75)(68,-85)
\psframe(33,-95)(103,-85)
\psline{->}(24,-90)(33,-90)
\rput(68,-90){\textbf{Chapter 8}. $q$-Plancherel measure.}
\psframe(42,-119)(94,-105)
\psline{->}(68,-95)(68,-105)
\rput(68,-110){\textbf{Chapter 9}. Measures}
\rput(68,-114){of parabolic induction.}
\psframe(36,-145)(96,-135)
\rput(-22,-140){\textbf{Chapter 10}. Schur-Weyl measures.}
\rput(66,-140){\textbf{Chapter 11}. Gelfand measures.}
\psline{->}(40,-5)(113,-5)(113,-170)(92,-170)
\psline{->}(44,-170)(32,-170)
\psframe(44,-179)(92,-165)
\rput(68,-170){\textbf{Chapter 12}. Hecke-}
\rput(68,-174){-Ivanov-Kerov algebra.}
\psframe(-52,-179)(32,-165)
\rput(-10,-170){\textbf{Chapter 13}. Bundles over semigroups}
\rput(-10,-174){ and projective limits of group algebras.}
\endpspicture
\psset{arrowsize=1.5pt 2}
\end{center}
\tableofcontents
\listoffigures

\mainmatter
\pagestyle{fancy}
\fancyhead{}
\fancyfoot{}
\fancyfoot[C]{\thepage}
\renewcommand{\chaptermark}[1]{\markboth{\chaptername\ \thechapter.\ #1.}{}} 
\renewcommand{\sectionmark}[1]{\markright{\thesection.\ #1.}} 
\fancyhead[RO]{\rightmark}
\fancyhead[LE]{\leftmark}
\setlength{\headheight}{15.5pt}

\renewcommand{\cleardoublepage}{\clearpage}

\part[Représentations des groupes symétriques et mesure de Plancherel]{Représentations des groupes symétriques et asymptotique\\ de la mesure de Plancherel}
\setcounter{footnote}{0}
Cette première partie offre un panorama de la théorie des représentations du groupe symétrique, et nous insisterons en particulier sur les résultats de nature asym\-ptotique. Ainsi, nous évoquerons tour à tour :\vspace{2mm}
\begin{itemize}
\item[-] la combinatoire du groupe symétrique et de ses représentations, et le lien avec les fonctions symétriques donné par l'\textbf{isomorphisme de Frobenius-Schur} (\S\ref{frobenius}) ;\vspace{2mm}
\item[-] les \textbf{éléments de Jucys-Murphy}, et les propriétés profondes des algèbres des groupes symétriques qui se déduisent de leur étude (\S\ref{jucysmurphy}) ;\vspace{2mm}
\item[-] le problème d'Ulam et l'asymptotique de la \textbf{mesure de Plancherel} sur les partitions (chapitre \ref{plancherel}) ;\vspace{2mm}
\item[-] la théorie des représentations du groupe symétrique infini (\S\ref{syminfinite}), et l'interprétation des mesures de probabilité sur les partitions comme \textbf{processus ponctuels déterminantaux} (\S\ref{schurmeasure}) ;\vspace{2mm}
\item[-] l'analogie entre mesures sur les partitions et modèles matriciels, et en particulier l'\textbf{équi\-va\-lence de Baik-Deift-Johansson} (chapitre \ref{matrix}).\vspace{2mm}
\end{itemize}
Bien qu'aucun résultat nouveau ne soit démontré ici, nous avons jugé essentiel de consacrer une partie entière du mémoire à ces rappels. En effet, les résultats que nous avons obtenus pour les algèbres d'Hecke et les groupes $\GL(n,\For_{q})$ sont tous inspirés des résultats jusqu'alors connus pour le groupe symétrique ; et les outils que nous employerons dans les autres parties de ce mémoire sont tous issus de cette théorie. En particulier, nous présenterons dans le chapitre \ref{tool} :\vspace{2mm}
\begin{itemize}
\item[-] les \textbf{observables de diagrammes}, qui sont des outils indispensables pour l'étude asymptotique de mesures sur les partitions ;\vspace{2mm}
\item[-] les \textbf{permutations partielles}, dont des généralisations serviront à établir des identités génériques dans les algèbres d'Hecke.\vspace{2mm}
\end{itemize}
Mentionnons également le rôle important joué par l'algèbre des \textbf{fonctions symétriques}, \emph{cf.} le paragraphe \ref{symmetricfunction}.\bigskip\bigskip

Au cours de cette thèse, nous avons établi des généralisations des théorèmes \ref{farahathigman}, \ref{firstasymptoticplancherel} et \ref{secondasymptoticplancherel} ; le lecteur est invité à lire les sections consacrées à ces résultats, et à consulter à sa convenance les autres sections de cette partie. 
\begin{enumerate}
\item Dans le premier chapitre, on présente les bases de la théorie des représentations du groupe symétrique et de la théorie des fonctions symétriques. Nos principales sources sont bien sûr \cite[chapitre 1]{Mac95}, et aussi \cite{JK81} ; \cite{Zel81} pour l'étude de l'anneau de Grothendieck des représentations des groupes symétriques ; et \cite{FH59,OV04} pour l'étude des éléments de Jucys-Murphy.
\item Le second chapitre présente l'algèbre des observables de diagramme, et on y décrit les relations entre les diverses bases de cette algèbre, en suivant essentiellement \cite{IO02} et \cite{Bia03a}.
\item Le chapitre \ref{plancherel} présente la mesure de Plancherel, en partant du problème classique des plus longues sous-suites croissantes d'une permutation. On expose la loi des grands nombres de Logan-Shepp-Kerov-Vershik et le théorème central limite de Kerov ; nos principales sources sont les articles originaux \cite{LS77,KV77,Ker93a}, et de nouveau \cite{IO02} pour la preuve complète et rigoureuse de ces résultats asymptotiques.
\item Le chapitre \ref{determinantal} explique le lien entre systèmes de mesures de probabilité sur les partitions, représentations du groupe symétrique infini et processus ponctuels déterminantaux. Les mesures les plus générales que l'on sait décrire par des systèmes de particules à corrélations déterminantales sont les mesures de Schur ; on expose brièvement leur théorie et leur usage, sachant que les mesures de Plancherel des algèbres d'Hecke $\IH_{q}(\sym_{n})$ et des groupes linéaires finis $\GL(n,\For_{q})$ rentrent dans ce cadre, \emph{cf.} la seconde partie du mémoire. On s'inspire de \cite{Oko97,KOV04} pour la théorie des représentations de $\sym_{\infty}$, et de \cite{Oko03a,Oko03b} pour la théorie des mesures de Schur.
\item Enfin, le cinquième chapitre expose la correspondance de Baik-Deift-Johansson, c'est-à-dire l'équivalence asymptotique entre la mesure du \emph{Gaussian Unitary Ensemble} (GUE) sur les valeurs propres des matrices hermitiennes, et la mesure de Plancherel sur les partitions, voir \cite{BDJ99,BDJ00,BOO00,Oko00}. Ce chapitre est présent à titre purement culturel, et nous ne sommes pas parvenus à généraliser les résultats qui y sont évoqués. 
\end{enumerate}
Les derniers rappels qui seront nécessaires concernent la théorie des représentations des groupes linéaires finis $\GL(n,\For_{q})$ ; ces rappels sont reportés au chapitre \ref{general}. \bigskip\bigskip

La lecture des chapitres \ref{determinantal} et \ref{matrix} est sans doute facultative, et comme mentionné précédemment, aucun résultat nouveau n'y est présenté. Nous avons néanmoins jugé utile de compléter notre panorama par ces chapitres, car ils mettent en exergue le lien entre partitions aléatoires et matrices aléatoires. Comme nous l'avons expliqué dans l'introduction générale, les mêmes outils sont employés de part et d'autre, et les résultats qu'on obtient sont extrêmement semblables ; ainsi, le lien avec la théorie des matrices aléatoires est l'une des motivations principales de la théorie asymptotique des représentations. Dans cette correspondance (que l'on complète ici par la correspondance RSK),
\begin{align*}
&\text{matrices hermitiennes gaussiennes}\,\,\leftrightarrow\,\, \text{mesures de Plancherel}\,\,\leftrightarrow\,\, \text{permutations aléatoires}\\
&\text{matrices symétriques gaussiennes}\,\,\leftrightarrow\,\, \text{mesures de Gelfand}\,\,\leftrightarrow\,\, \text{involutions aléatoires}\\
&\text{matrices de covariance}\,\,\leftrightarrow \,\,\text{mesures de Schur-Weyl}\,\,\leftrightarrow \,\,\text{mots aléatoires}
\end{align*} 
et l'on conjecture un lien semblable entre les $\beta$-ensembles et les $\beta$-mesures de Plancherel (voir le paragraphe \ref{asymptoticgelfand}). Malheureusement, il semble qu'il n'y ait pas d'équivalent matriciel naturel des $q$-mesures de Plancherel (chapitres \ref{qplancherelmeasure} et \ref{arikikoike}); un résultat dû à Borodin laisse penser qu'un tel <<~modèle matriciel~>> serait plutôt donné par la taille des blocs de Jordan de matrices choisies aléatoirement dans un groupe de Chevalley fini, \emph{cf.} la fin de la section \ref{qgaussian}.

\chapter{Permutations, partitions et représentations du groupe symétrique}\label{permutation}

Dans ce qui suit, $n$ est un entier positif ou nul, et $\lle 1,n\rre$ désigne l'ensemble des entiers compris entre $1$ et $n$. Ce premier chapitre est consacré à l'étude des \textbf{permutations} de $\lle 1,n \rre$, c'est-à-dire les applications bijectives $\sigma : \lle 1,n\rre \to \lle 1,n\rre$.

\section{Groupe symétrique et classes de conjugaison}
Le \textbf{groupe symétrique} d'ordre $n$ est le groupe des permutations de $\lle 1,n\rre$, la loi de groupe étant la composition des applications. Ce groupe fini sera noté $\sym_{n}$, et il contient $$n!=1 \times 2 \times \cdots \times n$$ éléments. Une permutation peut être donnée par son \textbf{mot} --- c'est-à-dire la suite de ses valeurs --- ou par sa décomposition en \textbf{cycles à supports disjoints}. Ainsi, la permutation $\sigma \in \sym_{8}$ dont le mot est $38254761$ admet pour décomposition en cycles
$$\sigma=(1,3,2,8)(4,5)(6,7)\,,$$
et cette décomposition est unique à permutation des cycles près. Le \textbf{type} d'une permutation est la suite ordonnée des tailles des orbites de la permutation ; ainsi, $\sigma=38254761$ a pour type $t(\sigma)=(4,2,2)$. D'autre part, deux permutations $\sigma_{1}$ et $\sigma_{2}$  de même taille $n$ sont dites \textbf{conjuguées} s'il existe une permutation $\tau \in \sym_{n}$ telle que $\sigma_{2}=\tau \circ \sigma_{1} \circ \tau^{-1}$.

\begin{proposition}[Classes de conjugaison du groupe symétrique]~
Deux permutations de $\sym_{n}$ sont conjuguées si et seulement si elles ont même type. Ainsi, les classes de conjugaison de $\sym_{n}$ sont paramétrées par les partitions de taille $n$, c'est-à-dire les suites décroissantes d'entiers positifs $\lambda=(\lambda_{1}\geq \lambda_{2} \geq \cdots \geq \lambda_{r})$ telles que $|\lambda|=\sum_{i=1}^{r}\lambda_{i}=n$.
\end{proposition}

\noindent En effet, si $\tau$ est une permutation et si $c=(a_{1},a_{2},\ldots,a_{m})$ est un cycle, alors le conjugué $\tau c \tau^{-1}$ est le cycle 
$$ c^\tau=(\tau(a_{1}),\tau(a_{2}),\ldots,\tau(a_{m}))\,.$$
De plus, tout cycle de longueur $m$ est obtenu de la sorte en choisissant convenablement $\tau$. D'autre part, la conjugaison par $\tau$ est un morphisme de groupes $\sym_{n}\to \sym_{n}$. Par suite, si $\sigma$ est le produit de cycles disjoints de longueurs respectives $l_{1},\ldots,l_{r}$, alors les conjugués de $\sigma$ sont toutes les permutations qui sont des produits de cycles disjoints de longueurs $l_{1},\ldots,l_{r}$. Et à réindexation des cycles près, on peut supposer $l_{1}\geq l_{2} \geq \cdots \geq l_{r}$, donc les classes de conjugaison de $\sym_{n}$ sont bien indexées par les partitions de taille $n$.

\section{Partitions et fonctions symétriques}\label{symmetricfunction}
Nous noterons $\ym_{n}$ l'ensemble des \textbf{partitions} de l'entier $n$, et $\ym=\bigsqcup_{n \in \N} \ym_{n}$ l'ensemble de toutes les partitions. La \textbf{longueur} d'une partition est son nombre de parts $\ell(\lambda)$ ; son \textbf{poids} ou sa \textbf{taille} est la somme de ses parts $|\lambda|$. Nous aurons également besoin de la quantité
$$b(\lambda)=\sum_{i=1}^{\ell(\lambda)}(i-1)\,\lambda_{i}\,,$$
qui intervient dans de nombreuses formules pour les algèbres d'Hecke et les groupes linéaires finis. Par exemple, si $\lambda=(5,4,2)$, alors $|\lambda|=11$, $\ell(\lambda)=3$ et $b(\lambda)=8$. Une partition sera souvent représentée par son \textbf{diagramme de Young} (on parle aussi de diagramme de Ferrers) : c'est le tableau à $r$ lignes avec $\lambda_{1}$ cases sur la première ligne, $\lambda_{2}$ cases sur la seconde ligne, etc. Par exemple, le diagramme de la partition $(5,4,2)$ est :
\figcapt{
\yng(2,4,5)
}{Diagramme de Young de la partition $(5,4,2)$ (dessiné à la fran\c caise, c'est-à-dire du bas vers le haut).}{Diagramme de Young de la partition $(5,4,2)$}

La \textbf{partition conjuguée} de $\lambda$ est la partition de même taille $\lambda'$ obtenue en échangeant les lignes et les colonnes du diagramme ; par exemple, $(5,4,2)'=(3,3,2,2,1)$. D'autre part, le \textbf{contenu} d'une case $(x,y)$ du diagramme est la différence $c(x,y)=x-y$, et le contenu d'un diagramme est la somme $c(\lambda)$ des contenus des cases. On voit aisément que $c(\lambda)=b(\lambda')-b(\lambda)$ pour tout diagramme :
\begin{center}
{\footnotesize\psset{unit=0.7mm}
\pspicture(0,-25)(155,18)
\rput(-18,0){$c(5,4,2)\,\,=$}
\psline(0,-10)(40,-10)(40,-2)(32,-2)(32,6)(16,6)(16,14)(0,14)(0,-10)
\psline(8,14)(8,-2)(32,-2)
\psline(0,6)(16,6)(16,-10)
\psline(0,-2)(8,-2)(8,-10)
\psline(24,6)(24,-10)
\psline(32,-2)(32,-10)
\rput(4,10){$-2$}
\rput(12,10){$-1$}
\rput(4,2){$-1$}
\rput(12,2){$0$}
\rput(20,2){$1$}
\rput(28,2){$2$}
\rput(4,-6){$0$}
\rput(12,-6){$1$}
\rput(20,-6){$2$}
\rput(28,-6){$3$}
\rput(36,-6){$4$}
\rput(50,0){$=$}

\psline(115,-10)(155,-10)(155,-2)(147,-2)(147,6)(131,6)(131,14)(115,14)(115,-10)
\psline(123,14)(123,-2)(147,-2)
\psline(115,6)(131,6)(131,-10)
\psline(115,-2)(123,-2)(123,-10)
\psline(139,6)(139,-10)
\psline(147,-2)(147,-10)
\rput(119,10){$2$}
\rput(127,10){$2$}
\rput(119,2){$1$}
\rput(127,2){$1$}
\rput(135,2){$1$}
\rput(143,2){$1$}
\rput(107.5,0){$-$}

\psline(60,-10)(100,-10)(100,-2)(92,-2)(92,6)(76,6)(76,14)(60,14)(60,-10)
\psline(68,14)(68,-2)(92,-2)
\psline(60,6)(76,6)(76,-10)
\psline(60,-2)(68,-2)(68,-10)
\psline(84,6)(84,-10)
\psline(92,-2)(92,-10)
\rput(72,10){$1$}
\rput(72,2){$1$}
\rput(80,2){$2$}
\rput(88,2){$3$}
\rput(72,-6){$1$}
\rput(80,-6){$2$}
\rput(88,-6){$3$}
\rput(96,-6){$4$}
\rput(74.5,-20){$=\,\,\,b(5,4,2)'-b(5,4,2)$}
\endpspicture}
\end{center}
Enfin, une partition $\lambda$ pourra être écrite multiplicativement sous la forme 
$$\lambda=1^{m_{1}}\,2^{m_{2}}\,\cdots \,s^{m_{s}},$$
où $m_{k}=m_{k}(\lambda)$ est le nombre de parts égale à $k$ dans $\lambda$.  Cette écriture permet de déterminer la taille de la classe de conjugaison associée à $\lambda$ dans $\sym_{n}$ : en effet, le centralisateur d'une permutation de type $\lambda$ a pour cardinal 
$$z_{\lambda}=\prod_{k \geq 1} m_{k}!\,\,k^{m_{k}},$$
et la classe de conjugaison $C_{\lambda}$ a donc pour cardinal $c_{\lambda}=n!/z_{\lambda}$.\bigskip
\bigskip

En dehors du contexte précédemment évoqué, la classe combinatoire des partitions d'entiers intervient classiquement dans la théorie des fonctions symétriques. Rappelons qu'un polynôme en $m$ variables $p(x_{1}, \ldots,x_{m})$ est dit \textbf{symétrique} si
$$ p(x_{\sigma(1)},x_{\sigma(2)},\ldots,x_{\sigma(m)})=p(x_{1},x_{2},\ldots,x_{m})$$
pour toute permutation $\sigma \in \sym_{m}$. L'ensemble $\Lambda_{R}(m)$ des polynômes symétriques en $m$ variables et à coefficients dans un anneau commutatif $R$ est un anneau gradué par le degré. De plus, la spécialisation $x_{m+1}=0$ fournit un morphisme d'anneaux gradués surjectif $\Lambda_{R}(m+1) \to \Lambda_{R}(m)$, d'où une famille dirigée d'anneaux gradués $(\Lambda_{R}(m))_{m \in \N}$. On appelle \textbf{fonction symétrique} un élément de la limite projective (dans la catégorie des anneaux gradués) :
$$ \Lambda_{R} = \varprojlim_{m \to \infty} \Lambda_{R}(m)\,.$$
Un tel objet peut être spécialisé en n'importe quel alphabet fini $X=\{x_{1},x_{2},\ldots,x_{m}\}$, et aussi en des alphabets infinis dénombrables.
En faisant agir le groupe symétrique $\sym_{m}$ sur les monômes, on voit que les polynômes
$$m_{\lambda}(x_{1},\ldots,x_{m})=\sum_{i_{1} \neq i_{2} \neq \cdots \neq i_{r}} (x_{i_{1}})^{\lambda_{1}}(x_{i_{2}})^{\lambda_{2}}\cdots (x_{i_{r}})^{\lambda_{r}}$$
avec $\lambda$ partition de longueur inférieure à $m$ forment une $R$-base de $\Lambda_{R}(m)$. La limite projective ôte la restriction $\ell(\lambda)\leq m$, donc $\Lambda_{R}$ admet une base $(m_{\lambda})_{\lambda \in \ym}$ indexée par les partitions --- on dit que les $m_{\lambda}(x)$ sont les \textbf{fonctions monomiales}. Ce résultat implique en particulier l'identité $\Lambda_{R}=\Lambda_{\Z} \otimes_{\Z} R$ pour tout anneau commutatif $R$ ; la plupart des raisonnements de fonctions symétriques pourront donc être effectués dans $\Lambda=\Lambda_{\Z}$.\bigskip\bigskip

Depuis Newton, on connaît des bases algébriques de l'algèbre $\Lambda$ (ou $\Lambda_{\Q}$), à savoir, les \textbf{fonctions élémentaires}, les \textbf{fonctions homogènes} et les \textbf{fonctions sommes de puissances} :
$$e_{n}(x)=m_{1^n}(x)=\!\!\!\!\!\sum_{i_{1} \neq i_{2} \neq \cdots \neq i_{n}}\!\!\!\!\! x_{i_{1}}x_{i_{2}}\cdots x_{i_{n}}
\quad;\quad
h_{n}(x)=\sum_{\lambda \in \ym_{n}} m_{\lambda}(x)\quad;\quad p_{n}(x)=\sum_{i} (x_{i})^n\,.
$$
\begin{example}
Développées sur l'alphabet $X=\{a,b,c\}$, les fonctions symétriques $e_{3}$, $h_{3}$ et $p_{3}$ s'écrivent :
\begin{align*}&e_{3}(X)= abc\qquad;\qquad p_{3}(X)=a^{3}+b^{3}+c^{3}\qquad; \\
&h_{3}(X)=a^{3}+b^{3}+c^{3}+a^{2}b+a^{2}c+b^{2}a+b^{2}c+c^{2}a+c^{2}b+abc\,.
\end{align*}
\end{example}
\noindent On renvoie à \cite[\S1.2]{Mac95} pour une démonstration complète de la proposition suivante :
\begin{proposition}[Bases algébriques des fonctions symétriques] $\!\!\!$
 L'anneau des fonctions symétriques est librement engendré par les fonctions élémentaires ou les fonctions homogènes : 
$$\Lambda=\Z[e_{1},e_{2},\ldots,e_{n},\ldots]=\Z[h_{1},h_{2},\ldots,h_{n},\ldots]\,.$$
Il en va de même pour les sommes de puissances si l'on étend l'anneau de base au corps des nombres rationnels : $\Lambda_{\Q}=\Q[p_{1},p_{2},\ldots,p_{n},\ldots]$.
\end{proposition}
\noindent Une preuve agréable consiste à introduire les séries génératrices :
$$E(t)=\sum_{n =0}^\infty e_{n}(x)\,t^n\quad;\quad H(t)=\sum_{n =0}^\infty h_{n}(x)\,t^n\quad ; \quad P(t)=\sum_{n = 1}^{\infty} p_{n}(x)\,\frac{t^n}{n}\,.$$ 
En effet, $H(t)=E(-t)^{-1}$ et $P(t)=\log H(t)$. \bigskip
\bigskip

Dans ce qui suit, si $\lambda \in \ym$, on notera $e_{\lambda}=e_{\lambda_{1}}e_{\lambda_{2}}\cdots e_{\lambda_{r}}$, et de même pour $h_{\lambda}$ et $p_{\lambda}$. Les familles $(e_{\lambda})_{\lambda \in \ym}$, $(h_{\lambda})_{\lambda \in \ym}$ et $(p_{\lambda})_{\lambda \in \ym}$ forment des bases linéaires de $\Lambda_{\Q}$, et toutes ces bases sont homogènes : 
$$\forall \lambda,\,\,\deg m_{\lambda}=\deg e_{\lambda}=\deg h_{\lambda}=\deg p_{\lambda}=|\lambda|\,.$$
Un produit scalaire  sur l'anneau $\Lambda$ est défini par $\scal{p_{\lambda}}{p_{\mu}}=z_{\lambda} \,\delta_{\lambda\mu}$ --- c'est le \textbf{produit de Hall}. D'autre part, on peut définir un coproduit $m^{*}$ en posant
$$m^{*}(p_{n})=p_{n} \otimes 1+1\otimes p_{n}\,,$$
et en imposant que $m^{*} : \Lambda \to \Lambda \otimes \Lambda$ soit un morphisme d'anneaux. Alors, on peut montrer que $(\Lambda , m, m^{*}, \scal{\cdot}{\cdot})$ est une \textbf{algèbre de Hopf}\footnote{Les sommes de puissances ne forment pas une $\Z$-base de $\Lambda$, mais en utilisant les relations de Newton, on peut montrer que $\scal{\cdot}{\cdot}$ et $m^{*}$ sont correctement définis sur $\Z$ (et pas seulement sur $\Q$). En particulier, la base duale de $(h_{\lambda})_{\lambda \in \ym}$ est $(m_{\lambda})_{\lambda \in \ym}$, et $m^{*}(h_{n})=\sum_{p+q=n} h_{p}\otimes h_{q}$.} ; l'unité est le plongement $\eps : \Z \to \Lambda$, la counité est l'application
$$\eps^{*} : f \in \Lambda \mapsto f(0,0,\ldots ,0,\ldots) \in \Z\,,$$
et l'antipode est définie par $S(e_{\lambda})=h_{\lambda}$ et $S(h_{\lambda})=e_{\lambda}$.
De plus, l'algèbre de Hopf $\Lambda$ est \textbf{autoadjointe}, c'est-à-dire que vis-à-vis du produit scalaire $\scal{\cdot}{\cdot}$, le coproduit est l'application duale du produit, la counité est l'application duale de l'unité, et l'antipode est autoadjointe. Cette structure additionnelle jouera un rôle important dans l'isomorphisme de Frobenius-Schur, \emph{cf.} le paragraphe \S\ref{frobenius} --- on renvoie d'autre part à \cite{Car06} pour un exposé général de la théorie des algèbres de Hopf.\bigskip

\section{Tableaux et représentations des groupes symétriques}\label{tableau}
Si $G$ est un groupe (fini), on rappelle qu'une \textbf{représentation} (linéaire, complexe, de dimension finie) de $G$ est la donnée d'un morphisme de groupes $ \rho : G \to \GL(V)$, où $V$ est un espace vectoriel complexe (de dimension finie). Alternativement, une représentation de $G$ est la donnée d'un \textbf{$\C G$-module} $V$, où $\C G$ désigne l'algèbre du groupe $G$, c'est-à-dire l'espace des combinaisons linéaires formelles
$$f=\sum_{g \in G} f(g) \,g$$
d'éléments du groupe. On peut aussi voir $\C G$ comme l'algèbre des fonctions complexes sur le groupe, avec pour multiplication le produit de convolution des fonctions :
$$f_{1}*f_{2}(g)=\sum_{hk=g}f_{1}(h)\,f_{2}(k)=\sum_{h\in G}f_{1}(h)\,f_{2}(h^{-1}g)\,.$$
Toute représentation d'un groupe fini $G$ se scinde de manière unique en somme directe de \textbf{représentations irréductibles} (c'est-à-dire des $\C G$-modules \textbf{simples}), et les classes d'isomorphismes de $\C G$-modules simples forment un ensemble $\widehat{G}$ fini. D'autre part, une représentation $(\rho,V)$ d'un groupe fini $G$ est déterminée à isomorphisme près par son caractère $\varsigma^{V} : g \mapsto \tr \rho(g),$ et les caractères irréductibles forment une base orthonormale de la sous-algèbre $(\C G)^{G }$ de $\C G$ constituée des \textbf{fonctions centrales} 
$$(\C G)^{G} = \big\{f \,\,|\,\,\forall g,h,\,\,f(gh)=f(hg) \big\}= \big\{f \,\,|\,\,\forall g,h,\,\,f(hgh^{-1})=f(g) \big\},$$
le produit scalaire sur $\C G$ étant défini par 
$$\scal{f_{1}}{f_{2}}=\frac{1}{\card G}\,\sum_{g \in G} \overline{f_{1}(g)}\,f_{2}(g)\,.$$
Notons qu'alternativement, on peut voir $(\C G)^{G}$ comme le centre de l'algèbre $\C G$. Tous ces points sont exposés en détail dans \cite{JL93} ou \cite{Ser77}. Nous aurons essentiellement besoin de l'orthonormalité de la base des caractères irréductibles, qui implique en particulier le point suivant : le nombre de caractères irréductibles d'un groupe fini $G$ est toujours égal au nombre de classes de conjugaison de $G$. Un autre fait qui mérite d'être précisé ici est la décomposition de la \textbf{représentation régulière} (gauche) en somme de modules irréductibles. Ainsi, pour tout groupe fini $G$, chaque $\C G$-module irréductible $V$ intervient $\dim V$ fois dans $\C G$ : 
$$\C G \simeq_{\C G} \bigoplus_{(\rho,V) \in \widehat{G}} (\dim V)\,V\,.$$
Ce résultat peut être vu comme une conséquence du théorème de Wedderburn : en tant qu'algèbre semi-simple sur un corps algébriquement clos, $\C G$ est une somme directe d'algèbres de matrices, et la décomposition en blocs de $\C G$ est 
$$\C G \simeq \bigoplus_{(\rho,V) \in \widehat{G}} \hendo(V)\,.$$
Un isomorphisme est donné par la \textbf{transformée de Fourier abstraite}, qui à $f = \sum_{g \in G} f(g)\,g$ associe $$\widehat{f}=\bigoplus_{(\rho,V) \in \widehat{G}} \bigg(\sum_{g \in G} f(g) \,\rho(g)\bigg).$$
Cette transformée est même une isométrie d'espaces $\leb^{2}$ non commutatifs, voir la section \ref{plancherelprocess}. Ceci motivera l'introduction de la mesure de Plancherel d'un groupe fini, qui est la duale de la mesure de Haar.
\bigskip
\bigskip

Détaillons maintenant cette théorie générale dans le cas particulier où $G$ est le groupe symétrique d'ordre $n$. D'après ce qui précède, les représentations irréductibles de $\sym_{n}$ sont en même nombre que les classes de conjugaison de $\sym_{n}$, et il existe donc une indexation de ces représentations par les partitions $\lambda \in \ym_{n}$. Une description explicite des modules simples $V^{\lambda}$ est due à A. Young, et elle repose sur la combinatoire des \textbf{tableaux}. Si $\lambda$ est une partition de taille $n$, on appelle \textbf{tableau standard} de forme $\lambda$ une numérotation des cases du diagramme $\lambda$ par les entiers de $\lle 1,n \rre$, les entrées étant strictement croissantes selon les lignes et selon les colonnes. Par exemple,
\begin{center}
\vspace{-2mm}
\young(25,134)
\vspace{1mm}
\end{center}
est un tableau standard de forme $(3,2)$. D'autre part, pour toute partition $\lambda$, notons $\Delta_{\lambda}(x)$ le polynôme en $n$ variables défini par :
\begin{align*}
\Delta_{\lambda}(x)
&=\prod_{i=1}^{\ell(\lambda)}\Delta(x_{\lambda_{1}+\cdots+\lambda_{i-1}+1},\ldots,x_{\lambda_{1}+\cdots+\lambda_{i}})=\prod_{i=1}^{\ell(\lambda)}\left(\prod_{\lambda_{1}+\cdots+\lambda_{i-1}+1\leq j<k\leq \lambda_{1}+\cdots+\lambda_{i}} (x_{j}-x_{k})\right).
\end{align*}
Par exemple, $\Delta_{(3,2)}(x)=(x_{1}-x_{2})(x_{1}-x_{3})(x_{2}-x_{3})(x_{4}-x_{5})$. Plus généralement, si $T$ est un tableau standard de type $\lambda$, nous noterons $\Delta_{T}(x)$ le polynôme produit des différences de deux cases prises sur la même ligne ; par exemple, le polynôme associé au tableau de la page précédente est $\Delta_{T}(x)=(x_{1}-x_{3})(x_{1}-x_{4})(x_{3}-x_{4})(x_{2}-x_{5})$.\bigskip\bigskip

Le groupe symétrique $\sym_{n}$ agit sur $\C[x_{1},\ldots,x_{n}]$ par permutation des variables, et on appelle \textbf{module de Specht} de type $\lambda$ le $\C\sym_{n}$-module $V^{\lambda}$ engendré par le polynôme $\Delta_{\lambda}(x)$. On renvoie à \cite{JK81} ou à \cite[chapitre 7]{Ful97} pour une preuve du résultat fondamental suivant\footnote{La description des modules de Specht proposée dans les ouvrages précités met aussi en jeu les notions combinatoires de tabloïde, de sous-groupe de Young et d'idempotent de Young ; la description polynomiale évoquée ici est parfaitement équivalente, \emph{cf.} \cite[théorème IV]{Young77}.} :
\begin{proposition}[Représentations irréductibles du groupe symétrique]
Les modules $V^{\lambda}$ sont deux à deux non isomorphes, et ils forment un système complet de représentations irréductibles du groupe $\sym_{n}$. De plus, pour toute partition $\lambda$, le module $V^{\lambda}$ admet pour base les polynômes $\Delta_{T}(x)$, où $T$ parcourt l'ensemble $\mathrm{Std}(\lambda)$ des tableaux standards de forme $\lambda$.
\end{proposition}
\noindent En particulier, la dimension de $V^{\lambda}$ est le nombre de tableaux standards de forme $\lambda$, et compte tenu de la décomposition de la représentation régulière
$$\C\sym_{n}=\sum_{\lambda \in \ym_{n}} (\dim \lambda)\,V^{\lambda},$$
$n!$ est égal à la somme des carrés des cardinaux $|\mathrm{Std}(\lambda)|$. Nous donnerons plus loin une preuve combinatoire de cette \textbf{identité des carrés}. Dans ce qui suit, nous adopterons l'indexation de Specht des modules irréductibles de $\sym_{n}$, à ceci près qu'on conjuguera les diagrammes, c'est-à-dire que $V^{\lambda}=\C\sym_{n}[\Delta_{\lambda'}(x)]$. Ceci revient à lire les tableaux standards selon les colonnes, ce qui assurément ne change pas grand chose à la théorie jusqu'ici exposée.\bigskip

\section[Isomorphisme de Frobenius-Schur et combinatoire des caractères]{Isomorphisme de Frobenius-Schur et combinatoire\\ des caractères}\label{frobenius}

La description précédente des modules simples de l'algèbre $\C\sym_{n}$ ne permet pas un calcul direct des caractères irréductibles, et elle ne rend pas compte des \textbf{règles de branchement} qui régissent l'induction ou la restriction de caractères irréductibles entre $\sym_{n}$ et un sous-groupe de Young
$$\sym_{\lambda}= \sym_{\lambda_{1}} \times \sym_{\lambda_{2}} \times \cdots \times \sym_{\lambda_{r}}$$
avec $\lambda \in \ym_{n}$. Rappelons que si $H \subset G$ sont deux groupes finis, la représentation \textbf{induite} à partir d'un $\C H$-module $V$ est le $\C G$-module gauche  $$\mathrm{Ind}_{H}^{G} (V) = \C G \otimes_{\C H} V,$$ et la représentation \textbf{restreinte} $\mathrm{Res}_{H}^{G} (W)$ à partir d'un $\C G$-module gauche $W$ est le même espace vectoriel $W$, mais avec la loi extérieure restreinte à $\C H$. Les foncteurs d'induction et de restriction commutent aux sommes directes, et ils sont adjoints :
$$\hom_{G}(\mathrm{Ind}_{H}^{G}(V),W)=\hom_{H}(V,\mathrm{Res}_{H}^{G}(W))\,.$$
En revanche, l'induit d'un $\C H$-module irréductible n'est plus forcément irréductible, et de même pour la restriction d'un $\C G$-module irréductible. \bigskip
\bigskip

Examinons en particulier le cas où $G=\sym_{n+1}$ et $H=\sym_{n}$. Si $\lambda \in \ym_{n}$ et $\Lambda \in \ym_{n+1}$, on note $\lambda \nearrow \Lambda$ si le diagramme de $\Lambda$ est obtenu en rajoutant une case dans un coin du bord du diagramme de $\lambda$. Par exemple, les partitions $\Lambda$ de taille $8$ telles que $(3,2,2) \nearrow \Lambda$ sont $$(4,2,2)\,, \,(3,3,2) \text{ et } (3,2,2,1).$$
Si $T$ est un tableau standard de taille $n+1$ et de forme $\Lambda$, notons $\lambda(T)$ le diagramme obtenu à partir de $\Lambda$ en ôtant la case numérotée $n+1$. Alors, on peut montrer que le $\C\sym_{n}$-module engendré par $\Delta_{T}(x)$ est isomorphe à $V^{\lambda(T)}$, de sorte que :
\begin{proposition}[Branchement des représentations des groupes symétriques]\label{branching}
Si $\Lambda$ est une partition de taille $n+1$, alors 
$$\mathrm{Res}_{\sym_{n}}^{\sym_{n+1}}(V^{\Lambda}) \simeq_{\C\sym_{n}} \bigoplus V^{\lambda},$$
où la somme directe porte sur les partitions de taille $n$ telles que $\lambda \nearrow \Lambda$. Par propriété d'adjonction, si $\lambda$ est une partition de taille $n$, alors
$$\mathrm{Ind}_{\sym_{n}}^{\sym_{n+1}}(V^{\lambda})\simeq_{\C\sym_{n+1}} \bigoplus V^{\Lambda},$$
la somme portant cette fois sur les partitions de taille $n+1$ telles que $\lambda\nearrow \Lambda$.
\end{proposition}

\noindent On retrouve par récurrence sur $n$ l'identité $\dim V^{\lambda}=\card \mathrm{Std}(\lambda)$, car un tableau standard de forme $\lambda$ est équivalent à la donnée d'une suite de diagrammes
$$\emptyset = \lambda^{(0)} \nearrow \lambda^{(1)} \nearrow \cdots \nearrow \lambda^{(n)}=\lambda,$$
c'est-à-dire un chemin reliant l'origine à $\lambda$ dans le \textbf{graphe de Young}, voir la figure \ref{younggraph}.
\figcapt{\psset{unit=0.7mm}
\pspicture(0,-8)(80,90)
\psline[linewidth=0.2pt](5,0)(20,2)
\psline[linewidth=0.2pt](0,5)(2,20)
\psline[linewidth=0.2pt](7,15)(25,20)
\psline[linewidth=0.2pt,bordercolor=white, border=1pt](15,7)(20,25)
\psline[linewidth=0.2pt](2,25)(5,45)
\psline[linewidth=0.2pt](25,2)(45,5)
\psline[linewidth=0.2pt](30,20)(45,18.2)
\psline[linewidth=0.2pt,bordercolor=white, border=1pt](35,10)(45,22)
\psline[linewidth=0.2pt](10,35)(22,45) \psline[linewidth=0.2pt](25,25)(40,40)
\psline[linewidth=0.2pt,bordercolor=white, border=1pt](20,30)(18.2,45)
\psline[linewidth=0.2pt](5,50)(8,80)
\psline[linewidth=0.2pt](0,0)(-10,-5)
\rput(-12,-5){$\emptyset$} \psline[linewidth=0.2pt](50,5)(80,8)
\psline(0,0)(0,5)(5,5)(5,0)(0,0)
\psline(15,2)(25,2)(25,7)(15,7)(15,2) \psline(20,2)(20,7)
\psline(2,15)(2,25)(7,25)(7,15)(2,15) \psline(2,20)(7,20)
\psline(20,30)(20,20)(30,20)(30,25)(25,25)(25,30)(20,30) \psline(25,20)(25,25)(20,25)
\psline(5,35)(5,50)(10,50)(10,35)(5,35) \psline(5,40)(10,40) \psline(5,45)(10,45)
\psline(35,5)(50,5)(50,10)(35,10)(35,5) \psline(40,5)(40,10) \psline(45,5)(45,10)
\psline(45,17)(60,17)(60,22)(50,22)(50,27)(45,27)(45,17) \psline(45,22)(50,22)(50,17) \psline(55,17)(55,22)
\psline(17,45)(17,60)(22,60)(22,50)(27,50)(27,45)(17,45) \psline(22,45)(22,50)(17,50) \psline(17,55)(22,55)
\psline(40,40)(40,50)(50,50)(50,40)(40,40) \psline(45,40)(45,50) \psline(50,45)(40,45)
\psline(65,8)(85,8)(85,13)(65,13)(65,8) \psline(80,8)(80,13) \psline(70,8)(70,13) \psline(75,8)(75,13)
\psline(8,65)(8,85)(13,85)(13,65)(8,65) \psline(8,80)(13,80) \psline(8,70)(13,70) \psline(8,75)(13,75)
\endpspicture}{Les quatre premiers niveaux du graphe de Young $\ym=\bigsqcup_{n\in\N} \ym_{n}$.\label{younggraph}}{Premiers niveaux du graphe de Young}
En d'autres termes, le graphe de Young est le \textbf{diagramme de Bratteli} de la famille inductive d'algèbres complexes $(\C\sym_{n})_{n \in \N}$, voir \cite{Bra72} pour des précisions sur cette terminologie. \bigskip \bigskip

L'objet pertinent en vue d'une généralisation de ces règles de branchement est le \textbf{groupe de Grothendieck} de la catégorie des représentations des groupes symétriques.  Si $G$ est un groupe fini, notons $K(G)$ le $\Z$-module libre de base les classes d'isomorphisme de $G$-modules simples. Les éléments de $K(G)$ sont appelés \textbf{$G$-modules virtuels}, et tout module virtuel s'écrit comme différence $V \ominus W$ de deux $G$-modules, étant entendu que $V \ominus W=V' \ominus W'$ dans $K(G)$ si les deux $G$-modules $V \oplus W'$ et $V' \oplus W$ sont isomorphes. Plus loin, nous aurons également besoin des produits tensoriels $K_{R}(G)= R \otimes_{\Z} G$ pour $R$ anneau commutatif. Notons $K(\sym)$ la somme directe $\bigoplus_{n \in \N} K(\sym_{n})$, et $K_{R}(\sym)$ les versions tensorisées. Si $V$ et $W$ sont des représentations de $\sym_{p}$ et $\sym_{q}$, alors on peut construire une représentation de $\sym_{p+q}$ en considérant
$$m(V,W)=V \boxtimes W=\mathrm{Ind}_{\sym_{p} \times \sym_{q}}^{\sym_{p+q}} (V \otimes W)\,.$$
Le produit $m$ est distributif par rapport aux sommes directes et est compatible aux isomorphismes, donc peut être étendu au groupe $K(\sym)$ : on obtient ainsi une structure d'anneau gradué sur $K(\sym)$ et les versions tensorisées $K_{R}(\sym)$.\bigskip\bigskip

Les caractères déterminant à isomorphisme près les représentations irréductibles, pour tout groupe fini $G$, le groupe de Grothendieck complexifié $K_{\C}(G)$ s'identifie au centre $(\C G)^{G}$ de l'algèbre du groupe. En particulier, $K_{\C}(G)$ est muni du produit scalaire des caractères, et ce produit scalaire est l'extension sesquilinéaire de la règle
$$\scal{V}{W}=\dim \hom_{G}(V,W)\,.$$
On munit $K_{\C}(\sym)=\bigoplus_{n \in \N} K_{\C}(\sym_{n})$ du produit scalaire obtenu à partir des produits scalaires des $K_{\C}(\sym_{n})$ en imposant que la somme directe soit une somme orthogonale ; ainsi, $\scal{V^{\lambda}}{V^{\mu}}=\delta_{\lambda \mu}$ pour toutes partitions $\lambda$ et $\mu$. Finalement, en considérant l'adjoint de $m$ par rapport au produit scalaire, on définit un coproduit gradué
\begin{align*}
m^{*} : K_{\C}(\sym_{n}) &\to \bigoplus_{p+q=n} K_{\C}(\sym_{p}) \otimes K_{\C}(\sym_{q})\\
V&\mapsto \sum_{p+q=n} \mathrm{Res}^{\sym_{n}}_{\sym_{p} \times \sym_{q}}(V)
\end{align*}
sur l'algèbre $K_{\C}(\sym)$. Muni de ces opérations, $K_{\C}(\sym)$ a une structure d'algèbre de Hopf\footnote{La vérification de la structure d'algèbre de Hopf est tout à fait non triviale, et met en jeu le théorème de Mackey (voir \cite{Mack51}) ; on renvoie à \cite{Zel81} pour une preuve détaillée.} autoadjointe et graduée, l'unité et la counité se déduisant de l'identification naturelle entre $\Z$ et $K(\sym_{0})$, et l'antipode étant l'extension linéaire de la tensorisation par la représentation signature.

\begin{theorem}[Isomorphisme de Frobenius-Schur, \cite{Fro00}]\label{frobeniusschur}
Il existe une isométrie d'algèbres de Hopf autoadjointes graduées entre $K_{\C}(\sym)$ et $\Lambda_{\C}$, et si l'on impose une condition de positivité vis-à-vis des produits scalaires, alors celle-ci est unique à composition par l'antipode près. Si $C_{\mu}$ désigne la somme des permutations de type $\mu$ dans $Z(\C\sym_{n}) \simeq K_{\C}(\sym_{n})$, alors l'un des deux isomorphismes est obtenu par extension linéaire de la règle $C_{\mu} \mapsto p_{\mu}/z_{\mu}$.
\end{theorem}
\clearpage

L'isométrie donnée par le théorème \ref{frobeniusschur} est appelée \textbf{application caractéristique}, et est notée $\mathrm{ch}$ ; l'image par $\mathrm{ch}$ d'un module de Specht $V^{\lambda}$ est la \textbf{fonction de Schur} $s_{\lambda}$. Ces fonctions sont étudiées en détail dans \cite[\S1.3]{Mac95} ; la spécialisation de $s_{\lambda}$ en un alphabet fini $X=\{x_{1},\ldots,x_{n}\}$ est le quotient de fonctions antisymétriques 
$$s_{\lambda}(x_{1},\ldots,x_{n})=\frac{a_{\delta+\lambda}(x_{1},\ldots,x_{n})}{a_{\delta}(x_{1},\ldots,x_{n})}$$
avec $a_{\mu}(x_{1},\ldots,x_{n})=\sum_{\sigma \in \sym_{n}} \eps(\sigma)\,\sigma(x^\mu)=\det((x_{i}^{\mu_{j}})_{i,j})$ et $\delta=(n-1,n-2,\ldots,0)$, l'addition de partitions s'entendant ici terme à terme. Compte tenu de la propriété d'isométrie, les fonctions de Schur forment une base orthonormée de $\Lambda_{\C}$, et si $\varsigma^{\lambda}(\mu)$ désigne la valeur du caractère du module de Specht $V^{\lambda}$ en une permutation de type $\mu$, alors 
$$\varsigma^{\lambda}(\mu)=z_{\mu}\scal{\varsigma^{\lambda}}{C_{\mu}}_{K(\sym)}=\scal{s_{\lambda}}{p_{\mu}}_{\Lambda} .$$
Ce résultat constitue la \textbf{formule de Frobenius} ; nous en donnerons plus loin des analogues pour les algèbres d'Hecke et les groupes linéaires finis. Une conséquence directe de la formule de Frobenius est la \textbf{formule des équerres} pour la dimension d'un module $V^{\lambda}$ : 
$$\dim \lambda =  \scal{s_{\lambda}}{p_{1^{n}}} = \frac{n!}{\prod_{(i,j)\in \lambda}h(i,j)}\,,$$
où $h(i,j)$ est la longueur de la plus grande équerre contenue dans le diagramme $\lambda$ et de coin la case $(i,j)$. On renvoie à \cite[\S1.7]{Mac95} pour la preuve de cette identité ; le lecteur pourra également consulter la section \ref{shoji} de ce mémoire pour une preuve (indirecte) de la formule de Frobenius.
\begin{example}
La dimension de la représentation irréductible indexée par la partition $(3,2)$ est $5!/(4\times 3 \times 1 \times 2 \times 1)=5$.
\end{example}
\bigskip
\bigskip

La formule de Frobenius ramène la combinatoire des caractères irréductibles des groupes symétriques à celle des fonctions de Schur, et permet \emph{in fine} un calcul explicite de ces caractères ; nous concluons cette section en expliquant en quelques mots ces calculs. Voyons d'abord comment développer une fonction de Schur dans l'une des autres bases de l'algèbre des fonctions symétriques. Comme les fonctions sommes de puissances forment une base orthogonale, la formule de Frobenius peut être inversée, et on peut exprimer $s_{\lambda}$ en fonction des $p_{\mu}$ :
$$\forall \lambda \in \ym_{n},\,\,\, s_{\lambda}(x)=\sum_{\mu \in \ym_{n}} \varsigma^{\lambda}(\mu)\, (z_{\mu})^{-1}\,p_{\mu}(x)\,.$$
D'autre part, en interprétant les fonctions antisymétrisées $a_{\mu}(x)$ comme des déterminants, on peut exprimer les fonctions de Schur dans les bases $(h_{\lambda})_{\lambda \in \ym}$ et $(e_{\lambda})_{\lambda \in \ym}$, voir \cite[\S1.3, p. 41-42]{Mac95}. Ainsi, on dispose des formules de Jacobi-Trudi :
$$s_\lambda=\det((h_{\lambda_i-i+j})_{i,j}) \qquad;\qquad s_{\lambda'}= \det((e_{\lambda_i-i+j})_{i,j})\,.$$
\begin{example} $s_{2,1,1}=(1/8)p_{1,1,1,1}-(1/4)p_{2,1,1}-(1/8)p_{2,2}+(1/4)p_{4}$

$\qquad\qquad\,\,\,\,\,=e_{3,1}-e_{4}=h_{2,1,1}-h_{2,2}-h_{3,1}+h_{4}$.
\end{example}
\clearpage

Dans ce qui suit, nous autorisons des \textbf{diagrammes gauches}, c'est-à-dire des diagrammes $\lambda \setminus \mu$ avec $\mu \subset \lambda$. Par exemple :
\figcapt{
$$\foreignlanguage{english}{\yng(2,4,5) \,\,\,\setminus\,\,\, \yng(1,2,3)\,\,\,=\,\,\,
\young(~,:~~,::~~)}$$\vspace{-6mm}
}{Diagramme de Young de la partition gauche $(5,4,2) \setminus (3,2,1)$.}{Diagramme de Young de la partition gauche $(5,4,2) \setminus (3,2,1)$}

\noindent Un diagramme gauche est appelé \textbf{bande horizontale} si toute colonne contient au plus une case, et \textbf{bande verticale} si toute ligne contient au plus une case. Une \textbf{bande frontière}, ou \textbf{ruban} est un diagramme gauche qui ne contient pas de blocs de $2 \times 2$ cases; ainsi, les colonnes et les lignes successives se chevauchent en au plus une case. Le diagramme gauche précédent est un ruban, et ce ruban a deux composantes connexes ; les rubans connexes sont ceux dont les colonnes et les lignes successives se chevauchent en exactement une case. \bigskip\bigskip

En développant les déterminants de Jacobi-Trudi suivant une ligne ou une colonne, on peut démontrer les \textbf{règles de Pieri} : pour toute partition $\lambda$, le produit $s_{\lambda}\,h_{r}$ (respectivement, le produit $s_{\lambda}\,e_{r}$) est égal à la somme $\sum s_{\mu}$ des fonctions de Schur indexées par des partitions $\mu$ telles que $\lambda \setminus \mu$ soit une bande horizontale de poids $r$ (resp., une bande verticale de poids $r$). On retrouve en particulier le théorème de branchement \ref{branching} en prenant $r=1$, et on peut déduire des règles de Pieri et de la formule de Frobenius le résultat plus avancé suivant (voir \cite[\S1.7, exemple 5]{Mac95}, et \cite{Gre92} pour une preuve alternative):
\begin{theorem}[Formule de Murnaghan-Nakayama, \cite{Murn37, Nak40a, Nak40b}]\label{murnaghannakayama}~

\noindent Si $\mu=(\mu_1,\ldots,\mu_m)$ et $\lambda$ sont deux partitions de poids $n$, alors
 $$\varsigma^{\lambda}(\mu)=\sum_{S} \,\,(-1)^{\mathrm{ht}(S)},$$
la somme étant effectuée sur les suites de partitions $S=(\emptyset=\lambda_0\subset\lambda_1 \subset \cdots \subset \lambda_m=\lambda)$ telles que chaque $\lambda_i\setminus \lambda_{i-1}$ soit un ruban connexe de poids $\mu_i$. Ici, $\mathrm{ht}(S)$ désigne la somme des hauteurs $\mathrm{ht}(\lambda_i\setminus \lambda_{i-1})$, la hauteur d'un ruban connexe étant son nombre de lignes moins $1$.
\end{theorem}
\noindent Ainsi, on dispose d'une formule récursive pour les caractères du groupe symétrique, ce qui permet par exemple de calculer $\varsigma^{5,4,2}(3,3,1,1,1,1,1)=6$. \bigskip

\section[Théorème de Farahat-Higman et éléments de Jucys-Murphy]{Théorème de Farahat-Higman et éléments de \\ Jucys-Murphy}\label{jucysmurphy}

Pour conclure cette présentation succinte de la théorie des représentations du groupe symétrique, étudions plus en détail les centres $Z(\C\sym_{n})$ des algèbres des grou\-pes symétriques. Nous en connaissons déjà deux bases : la base des classes de conjugaison et la base des caractères irréductibles. Si $\mu$ est une partition telle que $|\mu|+\ell(\mu) \leq n$, nous noterons $\mu \rightarrow n$ la \textbf{partition complétée} obtenue en ajoutant $1$ aux $n-|\mu|$ premières parts de $\mu$ (y compris éventuellement les parts nulles). Par exemple, $(5,4,2) \rightarrow 20 = (6,5,3,1,1,1,1,1,1)$. L'opération de complétion permet d'indexer les classes de conjugaison de $\sym_{n}$ par l'ensemble de partitions
$$\{\mu \in \ym \,\,|\,\,|\mu| + \ell(\mu)  \leq n\},$$
et cette indexation est plus naturelle pour les classes de permutations qui ont beaucoup de points fixes : par exemple, $C_{(2) \rightarrow n}$ est la classe des $3$-cycles dans $\sym_{n}$.

\begin{theorem}[Farahat-Higman, \cite{FH59}]\label{farahathigman}
L'algèbre $Z(\C\sym_{n})$ est graduée par le degré $\deg C_{\mu \rightarrow n}=|\mu|$, c'est-à-dire que pour toutes partitions $\lambda$ et $\mu$, 
$$C_{\lambda \rightarrow n} \,*\, C_{\mu \rightarrow n} = \sum a_{\lambda\mu}^{\tau}(n)\, C_{\tau \rightarrow n}\,,$$
la somme étant restreinte aux partitions $\tau$ telles que $|\tau| \leq |\lambda|+|\mu|$. De plus, les coefficients $a_{\lambda\mu}^{\tau}(n)$ sont des polynômes en $n$ à valeurs entières.
\end{theorem}
\bigskip

Le caractère gradué découle de la remarque suivante : si $\sigma$ est une permutation de $C_{\mu \rightarrow n}$, alors une factorisation minimale de $\sigma$ en produit de transpositions admet $|\mu|$ termes, donc le produit $\pi$ d'une permutation de $C_{\lambda \rightarrow n}$ par une permutation de $C_{\mu \rightarrow n}$ admet une factorisation en produit de $|\lambda|+|\mu|$ transpositions. Ceci implique $t(\pi)=\tau \rightarrow n$ avec $|\tau| \leq |\lambda|+|\mu|$. Pour le caractère polynomial des constantes de structure, nous renvoyons à la section \ref{ivanovkerov}, qui expose une preuve due à Ivanov et Kerov (\emph{cf.} \cite{IK99}) et reposant sur les permutations partielles ; la généralisation de cette méthode est l'objet principal de la dernière partie du mémoire. Dans \cite{FH59}, Farahat et Higman montrent également que si $|\tau|= |\lambda|+|\mu|$, alors le coefficient 
$$a_{\lambda\mu}^{\tau}(n)=a_{\lambda\mu}^{\tau}$$ ne dépend pas de $n$. Ceci permet de construire une algèbre graduée de base indexée par toutes les partitions $\lambda \in \ym$, et de constantes de structure les coefficients $a_{\lambda\mu}^{\tau}$ pour $|\tau|=|\lambda|+|\mu|$ : 
$$FH(\sym)= \left\langle x_{\lambda},\,\,\lambda \in \ym \,\,\,\bigg |\,\,\, \forall \lambda,\mu,\,\,x_{\lambda}\,*\,x_{\mu}=\sum_{|\tau| = |\lambda|+|\mu|}a_{\lambda\mu}^{\tau}\,x_{\tau} \right\rangle.$$
Cette \textbf{algèbre de Farahat-Higman} peut être considérée intuitivement comme une limite projective des anneaux gradués $Z(\C\sym_{n})$ ; nous détaillerons plus loin une con\-struction alternative de limite projective. Certains des coefficients de structure de l'algèbre de Farahat-Higman sont calculés dans \cite{GJ92,GJ94}.\bigskip
\bigskip

Les \textbf{éléments de Jucys-Murphy} (\cite{Juc74,Mur81}) constituent un autre outil essentiel de l'étude des centres $Z(n)=Z(\C\sym_{n})$ ; ils permettent en fait de retrouver toute la théorie des représentations des groupes symétriques, voir \cite{OV04}. Si $k$ est un entier inférieur à $n$, on note
$$J_{k}=\sum_{j<k} \,(j,k) = (1,k)+(2,k)+\cdots+(k-1,k)$$
la somme des transpositions $(j<k)$ dans l'algèbre du groupe symétrique d'ordre $n$. Ainsi, $J_{1}=0$, $J_{2}=(1,2)$, $J_{3}=(1,3)+(2,3)$, etc. Si $Z(n,n-1)$ désigne le centralisateur de $\C\sym_{n-1}$ dans $\C\sym_{n}$, alors il n'est pas difficile de voir que $J_{n} \in Z(n,n-1)$ pour tout entier $n$ ; par conséquent, les éléments $J_{1},\ldots,J_{n}$ commutent et engendrent une sous-algèbre abélienne de $\C\sym_{n}$. Cette algèbre est appelée \textbf{algèbre de Gelfand-Tsetlin}, et est notée $GZ(n)$. Des manipulations combinatoires sur les transpositions permettent d'établir les faits suivants :
\begin{enumerate}
\item Pour tout $n$, le centralisateur $Z(n,n-1)$ est l'algèbre engendrée par $Z(n-1)$ et $J_{n}$ ; l'algèbre de Gelfand-Tsetlin est donc également l'algèbre engendrée par les centres $Z(1), Z(2), \ldots, Z(n)$. Elle est abélienne maximale, et joue un rôle analogue à celui d'une sous-algèbre de Cartan dans une algèbre de Lie semi-simple.
\item Plus précisément, dans 
$$\C\sym_{n}\simeq \bigoplus_{\lambda \in \ym_{n}} \hendo (V^{\lambda})\,,$$ l'algèbre de Gelfand-Tsetlin consiste en les opérateurs qui, pour tout module irréductible $V^{\lambda}$, agissent diagonalement sur la base de Young indexée par les tableaux standards de forme $\lambda$ (voir la section \ref{hurwitz} pour des précisions). En particulier,
$$\dim GZ(n)=\sum_{\lambda \in \ym_{n}}\dim \lambda\,.$$
\item Plus précisément, les valeurs propres de l'action de $J_{n}$ sur un module $V^{\lambda}$ sont les contenus des coins $(\lambda_{i},i)$ du diagramme (voir l'exemple de la figure \ref{actionjucysmurphy}), et la décomposition spectrale induite par $J_{n}$ réalise le foncteur de réduction $\mathrm{Red}_{\sym_{n-1}}^{\sym_{n}}$ : 
$$V^{\lambda}=\bigoplus_{(\lambda_{i},i) \text{ coin de }\lambda} \big[\ker (J_{n}-c(\lambda_{i},i)) \simeq_{\C\sym_{n-1}} V^{\lambda \setminus (\lambda_{i},i)}\big]\,.$$
Les éléments de Jucys-Murphy d'ordre inférieur ont donc également pour valeurs pro\-pres les contenus des cases du diagramme. Ainsi,  les éléments $J_{1},\ldots ,J_{n}$ jouent un rôle analogue aux racines simples d'une algèbre de Lie semi-simple, et les contenus jouent un rôle analogue aux poids.
\end{enumerate}
\figcapt{\psset{unit=0.7mm}\pspicture(0,0)(150,25)
\psframe*[linecolor=gray](32,0)(40,8)
\psframe*[linecolor=gray](24,8)(32,16)
\psframe*[linecolor=gray](8,16)(16,24)
\psline(0,0)(40,0)(40,8)(32,8)(32,16)(16,16)(16,24)(0,24)(0,0)
\psline(8,0)(8,24)
\psline(16,0)(16,24)
\psline(24,0)(24,16)
\psline(32,0)(32,16)
\psline(0,16)(16,16)
\psline(0,8)(32,8)
\rput(100,16){$V^{5,4,2}=\ker(J_{11}+1) \oplus \ker(J_{11}-2) \oplus \ker(J_{11}-4)$}
\rput(103,8){$=\,\,\quad V^{5,4,1}\,\,\,\quad\oplus\quad \,\,V^{5,3,2} \,\,\,\quad\oplus\quad\,\,\, V^{4,4,2}$}
\endpspicture}{L'action de $J_{n}$ sur $V^{\lambda}$ est diagonale, et elle a pour valeurs propres les contenus des coins du diagramme.\label{actionjucysmurphy}}{Action de l'élément de Jucys-Murphy $J_{n}$ sur un module de Specht $V^{\lambda}$}

\noindent On renvoie à l'article \cite{OV04} pour une preuve de toutes ces assertions. Un argument crucial de la preuve est l'utilisation du critère de Gelfand : comme $Z(n,n-1)$ est une algèbre commutative, la paire $(\sym_{n-1}\subset \sym_{n})$ est une \textbf{paire de Gelfand} forte, c'est-à-dire que la restriction d'une représentation irréductible de $\sym_{n}$ à $\sym_{n-1}$ est sans multiplicité (autrement dit, le diagramme de Bratteli des groupes symétriques est un graphe sans arête multiple).
\bigskip
\bigskip

En utilisant les décompositions spectrales successives de $V^{\lambda}$ par rapport aux éléments de Jucys-Murphy, on obtient une décomposition de $V^{\lambda}$ en somme directe de droites vectorielles indexées par les chemins $\emptyset = \lambda^{(0)} \nearrow \cdots \nearrow \lambda^{(n)}=\lambda$ reliant l'origine à $\lambda$ dans le graphe de Young, c'est-à-dire par les tableaux standards de forme $\lambda$. Une base $(v_{T})_{T \in \mathrm{Std}(\lambda)}$ associée à cette décomposition en somme directe est caractérisée par la propriété
$$\C\sym_{i}[v_{T}]=V^{\lambda^{(i)}}$$
pour tout entier $i$. On retrouve ainsi toutes les propriétés évoquées dans la section \ref{tableau}, \emph{cf.} l'article \cite{OV04}. La combinatoire des éléments de Jucys-Murphy régit donc la théorie des représentations de $\C\sym_{n}$, et de plus, ces éléments sont particulièrement utiles pour l'analyse asymptotique, voir en particulier les sections \ref{freecumulant} et \ref{bdjgeom} --- c'est en vue de ces deux sections que nous avons introduit ici les éléments de Jucys-Murphy. Finalement, les éléments de Jucys-Murphy fournissent une nouvelle description des centres $Z(n)$ des algèbres des groupes. Ainsi :
 \begin{proposition}[Jucys-Murphy, \cite{Juc74, Mur81}] Le centre $Z(n)$ est exactement l'ensemble des fonctions symétriques en les éléments de Jucys-Murphy $J_{1},\ldots,J_{n}$.
 \end{proposition} 
\noindent En effet, étant données $n$ inconnues $z_{1},\ldots,z_{n}$, on montre aisément par récurrence sur $n$ l'identité formelle
$$\prod_{i=1}^{n} (z_{i} + J_{i})= \sum_{\sigma \in \sym_{n}} \bigg(\prod_{c \text{ cycle de }\sigma} z_{\min c}\bigg)\,\sigma$$
dans $\C[\sym_{n},z_{1},\ldots,z_{n}]$. En particulier, pour $z_{1}=z_{2}=\cdots=z_{n}=z$,
$$\sum_{r=0}^{n} e_{r}(J_{1},\ldots,J_{n})\,z^{n-r}=\prod_{i=1}^{n}(z+J_{i})=\sum_{\sigma \in \sym_{n}} z^{\ell(t(\sigma))}\,\sigma= \sum_{r=0}^{n} \bigg(\sum_{|\mu|=r}C_{\mu \rightarrow n}\bigg)\,z^{n-r} ,$$
 d'où l'identité $e_{r}(J_{1},\ldots,J_{n})=\sum_{|\mu|=r}C_{\mu \rightarrow n}$. Des manipulations combinatoires permettent ensuite d'exprimer toute classe de cycles $C_{(r) \rightarrow n}$ comme fonction symétrique en les éléments de Jucys-Murphy, et finalement toute classe de conjugaison $C_{\mu \rightarrow n}$, voir \cite{Juc74}. Notons néanmoins que la fonction symétrique $f_{\mu}$ telle que $C_{\mu \rightarrow n}=f_{\mu}(J_{1},\ldots,J_{n})$ dépend en général de $n$. D'ailleurs, il y a pour chaque $n$ et chaque $\mu$ une infinité de fonctions symétriques telles que $C_{\mu \rightarrow n}=f(J_{1},\ldots,J_{n})$, puisque $Z(n)$ est de dimension finie et $\Lambda$ est de dimension infinie. On renvoie à \cite{Juc74,Mur81,LT01} pour une étude plus approfondie des éléments de Jucys-Murphy et de leurs fonctions symétriques.

\chapter{Combinatoire des observables de diagrammes}\label{tool}

En théorie asymptotique des représentations des groupes symétriques $\sym_{n}$, on considère des diagrammes $\lambda$ distribués aléatoirement suivant diverses mesures de probabilité, la taille du groupe symétrique tendant vers l'infini (\emph{cf.} les chapitres \ref{plancherel}, \ref{qplancherelmeasure}, \ref{arikikoike}, \ref{schurweylmeasure} et \ref{gelfandmeasure}). Les coordonnées $\lambda_{1},\lambda_{2},\ldots,\lambda_{r}$ sont donc des variables aléatoires, et on s'intéresse à leur loi limite modulo une éventuelle renormalisation. Mais dans ce contexte, les coordonnées standards ne sont pas forcément les plus pratiques à manipuler ; ainsi, il est utile d'introduire d'autres coordonnées de diagrammes (\emph{cf.} \cite[\S1.1]{Mac95}), ainsi que les fonctions symétriques en ces nouvelles coordonnées. Dans \cite{IO02}, V. Ivanov et S. Kerov présentent une algèbre d'observables en ces coordonnées construite sur le modèle de l'algèbre des fonctions symétriques $\Lambda$. Ce chapitre est consacré à un exposé succinct de leur article, et à la présentation de quatre bases remarquables : la base des moments en les coordonnées de Frobenius, la base des moments en les coordonnées entrelacées, la base des caractères centraux et la base des cumulants libres (\cite{Bia98,Bia03a}). \bigskip

Toutes ces observables peuvent être évaluées sur toutes les partitions, et plus généralement sur toute une classe de fonctions continues qui contient les diagrammes de Young et leurs versions renormalisées. L'espace formé par ces \textbf{diagrammes continus} est un cadre topologique agréable pour énoncer les résultats de convergence de théorie asymptotique des représentations. Ainsi, dans la dernière section du chapitre (\S\ref{markovkrein}), nous expliquerons en détail à quelle topologie sur les diagrammes continus correspond la convergence des observables, ces arguments jouant un rôle crucial dans la preuve des théorèmes limites \ref{firstasymptoticplancherel}, \ref{firstasymptoticschurweyl} et \ref{firstasymptoticgelfand} (entre autres).\bigskip
\bigskip
 
Pour ce qui suit, il est utile d'introduire les \textbf{notations de $\lambda$-anneaux}, qui sont présentées clairement par exemple dans \cite[\S1]{RRW96}. Si $$X=\{x_{1},x_{2},\ldots,x_{m}\}$$ est un alphabet fini en $m$ variables, notons $p_{k}(X)$ la somme de puissances $\sum_{i=1}^m (x_{i})^k$. Comme les fonctions sommes de puissances forment une base algébrique de $\Lambda_{\Q}$, on peut calculer à partir des $p_{k}(X)$ toutes les autres fonctions symétriques de $X$, en particulier, les fonctions de Schur $s_{\lambda}(X)$, les fonctions homogènes $h_{\lambda}(X)$ et les fonctions élémentaires $e_{\lambda}(X)$, etc. Ceci étant, si $X$ et $Y$ sont deux alphabets (finis), on note :
$$X+Y=\{x_{1},\ldots,x_{m},y_{1},\ldots,y_{l}\}\qquad;\qquad XY=\{x_{i}\,y_{j}, \,\,1\leq i \leq m,\,\,1 \leq j \leq l \}\,.$$
Les spécialisations des fonctions symétriques associées à ces alphabets sont :
$$\forall k \geq 1,\,\,p_{k}(X+Y)=p_{k}(X)+p_{k}(Y)\qquad;\qquad \forall k \geq 1,\,\,p_{k}(XY)=p_{k}(X)\,p_{k}(Y)\,.$$
Maintenant, on peut définir des alphabets <<~virtuels~>> en se contentant de décrire la spécialisation de l'algèbre $\Lambda$ (par exemple, en donnant les fonctions sommes de puissances). Ainsi, on peut définir un alphabet formel $X-Y$ dont les fonctions symétriques sont données par 
$$p_{k}(X-Y)=p_{k}(X)-p_{k}(Y)\,.$$
Un autre exemple important d'alphabet formel est l'alphabet exponentiel $E$ défini par $p_{1}(E)=1$ et $p_{k \geq 2}(E)=0$. Pour cette spécialisation de $\Lambda$, la formule de Frobenius (inversée) montre que les fonctions de Schur sont données par $|\lambda|!\,\,s_{\lambda}(E)=\dim \lambda$.
Dans ce qui suit, nous utiliserons librement ces notations ; elles seront en particulier indispensables pour la description des caractères des algèbres d'Iwahori-Hecke, \emph{cf.} la section \ref{ram}.\bigskip

\section[Coordonnées entrelacées et coordonnées de Frobenius]{Coordonnées entrelacées et coordonnées de Frobenius}
Si $\lambda$ est un diagramme de poids $n$, on peut lui associer une fonction $\omega_{\lambda}$ affine par morceaux de pente $\pm 1$ et telle que :
$$\int_{\R} \left(\omega_{\lambda}(s)-|s|\right)\,ds=2n\,.$$
En effet, il suffit de tourner la représentation du diagramme de 45 degrés et de considérer la fonction <<~bord supérieur~>>, prolongée par la valeur absolue en dehors de son support (voir la figure \ref{russiandiagram} ; on demande aussi que les cases aient pour aire $2$).
\figcapt{\psset{unit=1mm}
\pspicture(-50,0)(50,50)
\psline{->}(-50,0)(50,0)
\psline{->}(0,0)(0,50)
\rput*[0,0]{45}(0,0){\psline[linestyle=dashed, linewidth=0.25pt](0,0)(0,32) \psline[linestyle=dashed, linewidth=0.25pt](8,0)(8,32) \psline[linestyle=dashed, linewidth=0.25pt](16,0)(16,24) \psline[linestyle=dashed, linewidth=0.25pt](24,0)(24,24) \psline[linestyle=dashed, linewidth=0.25pt](32,0)(32,24) \psline[linestyle=dashed, linewidth=0.25pt](40,0)(40,8)
\psline[linestyle=dashed, linewidth=0.25pt](0,0)(40,0) \psline[linestyle=dashed, linewidth=0.25pt](0,8)(40,8) \psline[linestyle=dashed, linewidth=0.25pt](0,16)(32,16) \psline[linestyle=dashed, linewidth=0.25pt](0,24)(32,24) \psline[linestyle=dashed, linewidth=0.25pt](0,32)(8,32)
\psline[linewidth=1pt](0,56)(0,32)(8,32)(8,24)(32,24)(32,8)(40,8)(40,0)(56,0)
}
\rput(42,41.5){$y=\lambda(x)$}
\endpspicture}{Fonction $s\mapsto \lambda(s)$ associée au diagramme de Young $\lambda=(5,4,4,1)$.\label{russiandiagram}}{Fonction $s\mapsto \lambda(s)$ associée à un diagramme de Young}
Dans ce qui suit, on note $\lambda(s)=\omega_{\lambda}(s)$ et on identifie le diagramme de Young à la fonction continue associée. Cette interprétation fonctionnelle est très utile pour l'étude asymptotique de mesures de probabilité sur les partitions, en particulier parce qu'elle permet de considérer des diagrammes renormalisés, voir la section \ref{momentum}. \'Etant donné un diagramme $\lambda$, la fonction $s\mapsto \lambda(s)$ est déterminée par la suite de ses extrema locaux 
$ x_{1}<y_{1}<x_{2}<y_{2}<\cdots < x_{v-1}<y_{v-1}<x_{v},$
où les $x_{i}$ sont les minima et les $y_{i}$ sont les maxima. Par exemple, pour $\lambda=(5,4,4,1)$, les extrema locaux sont les entiers 
$$-4 < -3' < -2 < 1' < 3 < 4' < 5\,,$$
le $'$ indiquant les maxima. On dit que les $x_{i}$ et les $y_{i}$ sont les \textbf{coordonnées entrelacées} du diagramme. Ces coordonnées permettront d'exprimer relativement simplement les probabilités de transition du processus de Plancherel et du $q$-processus de Plancherel, voir les paragraphes \ref{plancherelprocess} et \ref{qplancherelprocess}.\bigskip 
\bigskip

D'autre part, si $\lambda=(\lambda_{1},\ldots,\lambda_{r},0,\ldots)$ est une partition de taille $n$, on appelle \textbf{coordonnées de Frobenius} de $\lambda$ les deux suites
$$a_{i}=\lambda_{i}-i+1/2 \qquad;\qquad b_{i}=\lambda_{i}'-i+1/2\,,$$
où $i$ varie entre $1$ et la taille $d$ de la diagonale du diagramme. Les deux suites $(a_{i})_{1\leq i \leq d}$ et $(b_{i})_{1 \leq i \leq d}$ sont strictement décroissantes et déterminent entièrement le diagramme $\lambda$. De plus, on dispose de l'identité
$$\sum_{i=1}^d a_{i}+\sum_{i=1}^d b_{i}=n=|\lambda|\,,$$
qui est évidente si l'on voit $a_{i}$ ou $b_{i}$ comme l'aire d'une ligne ou colonne partant de la diagonale, les cases de la diagonale étant coupées en deux. Dans ce qui suit, nous noterons $A=A(\lambda)=\{a_{1},\ldots,a_{d}\}$ et $-B=-B(\lambda)=\{-b_{1},\ldots,-b_{d}\}$ les deux alphabets associés aux coordonnées de Frobenius d'un diagramme ; ces deux parties de $\Z'=\Z+1/2$ ont même cardinal.
\figcapt{\psset{unit=1mm} \pspicture(-50,0)(50,47)
\psline[linewidth=0.25pt](40,40)(0,0)(-40,40)(-37.5,42.5)(-30,35)(-27.5,37.5)(-25,35)(-22.5,37.5)(-12.5,27.5)(-7.5,32.5)(-5,30)(-2.5,32.5)(2.5,27.5)(5,30)(10,25)(20,35)(22.5,32.5)(25,35)(27.5,32.5)(37.5,42.5)(40,40)
\psline[linewidth=0.25pt](-37.5,37.5)(-35,40)
\psline[linewidth=0.25pt](-35,35)(-32.5,37.5)
\psline[linewidth=0.25pt](-32.5,32.5)(-30,35)
\psline[linewidth=0.25pt](-30,30)(-25,35)
\psline[linewidth=0.25pt](-27.5,27.5)(-20,35)
\psline[linewidth=0.25pt](-25,25)(-17.5,32.5)
\psline[linewidth=0.25pt](-22.5,22.5)(-15,30)
\psline[linewidth=0.25pt](-20,20)(-12.5,27.5)
\psline[linewidth=0.25pt](-17.5,17.5)(-2.5,32.5)
\psline[linewidth=0.25pt](-15,15)(0,30)
\psline[linewidth=0.25pt](-12.5,12.5)(2.5,27.5)
\psline[linewidth=0.25pt](-10,10)(7.5,27.5)
\psline[linewidth=0.25pt](-7.5,7.5)(10,25)
\psline[linewidth=0.25pt](-5,5)(22.5,32.5)
\psline[linewidth=0.25pt](-2.5,2.5)(27.5,32.5)
\psline[linewidth=0.25pt](2.5,2.5)(-30,35)
\psline[linewidth=0.25pt](5,5)(-25,35)
\psline[linewidth=0.25pt](7.5,7.5)(-12.5,27.5)
\psline[linewidth=0.25pt](10,10)(-10,30)
\psline[linewidth=0.25pt](12.5,12.5)(-5,30)
\psline[linewidth=0.25pt](15,15)(2.5,27.5)
\psline[linewidth=0.25pt](17.5,17.5)(10,25)
\psline[linewidth=0.25pt](20,20)(12.5,27.5)
\psline[linewidth=0.25pt](22.5,22.5)(15,30)
\psline[linewidth=0.25pt](25,25)(17.5,32.5)
\psline[linewidth=0.25pt](27.5,27.5)(20,35)
\psline[linewidth=0.25pt](30,30)(27.5,32.5)
\psline[linewidth=0.25pt](32.5,32.5)(30,35)
\psline[linewidth=0.25pt](35,35)(32.5,37.5)
\psline[linewidth=0.25pt](37.5,37.5)(35,40)
\psline[border=1pt, bordercolor=white]{->}(0,0)(0,50)
\psline{->}(-50,0)(50,0)
\psline[linecolor=NavyBlue]{-*}(38.75,41.25)(38.75,0)
\psline[linecolor=NavyBlue]{-*}(26.25,33.75)(26.25,0)
\psline[linecolor=NavyBlue]{-*}(21.25,33.75)(21.25,0)
\psline[linecolor=NavyBlue]{-*}(6.25,28.75)(6.25,0)
\psline[linecolor=NavyBlue]{-*}(8.75,26.25)(8.75,0)
\psline[linecolor=NavyBlue]{-*}(1.25,28.75)(1.25,0)
\psline[linecolor=BurntOrange]{-*}(-38.75,41.25)(-38.75,0)
\psline[linecolor=BurntOrange]{-*}(-28.75,36.25)(-28.75,0)
\psline[linecolor=BurntOrange]{-*}(-23.75,36.25)(-23.75,0)
\psline[linecolor=BurntOrange]{-*}(-11.25,28.75)(-11.25,0)
\psline[linecolor=BurntOrange]{-*}(-8.75,31.25)(-8.75,0)
\psline[linecolor=BurntOrange]{-*}(-3.75,31.25)(-3.75,0)
\psline[linewidth=1pt](37.5,42.5)(40,40)
\psline[linewidth=1pt](25,35)(27.5,32.5)
\psline[linewidth=1pt](20,35)(22.5,32.5)
\psline[linewidth=1pt](5,30)(10,25)
\psline[linewidth=1pt](0,30)(2.5,27.5)
\psline[linewidth=1pt](-40,40)(-37.5,42.5)
\psline[linewidth=1pt](-30,35)(-27.5,37.5)
\psline[linewidth=1pt](-25,35)(-22.5,37.5)
\psline[linewidth=1pt](-12.5,27.5)(-7.5,32.5)
\psline[linewidth=1pt](-5,30)(-2.5,32.5)
\endpspicture}{Coordonnées de Frobenius d'un diagramme de Young.\label{frobeniuscoor}}{Coordonnées de Frobenius d'un diagramme de Young}

\section[Moments d'un diagramme et graduations sur l'algèbre d'observables]{Moments d'un diagramme et graduations sur l'algèbre\\ d'observables}\label{momentum}
\begin{definition}[Observable de diagrammes]
On appelle observable d'un diagramme $\lambda$ toute fonction symétrique en l'alphabet (virtuel) $A-(-B)$, donc, toute combinaison linéaire de produits des moments de Frobenius
$$p_{k}(\lambda)=p_{k}(A-(-B))=\sum_{i=1}^{d} \,(a_{i})^{k} -(-b_{i})^{k}.$$
\end{definition}
\noindent Par exemple, la taille $|\lambda|$ est une observable de diagrammes, car c'est le premier moment $p_{1}(\lambda)$. Dans ce qui suit, on note $\obs$ l'\textbf{algèbre des observables de diagrammes} ; l'algèbre $\obs$ est librement engendrée par les moments de Frobenius (\cite[proposition 1.5]{IO02}), et est donc isomorphe à l'algèbre $\Lambda$. En particulier, les fonctions $$\lambda \mapsto p_{\mu}(\lambda), \,\,\,\mu \in \ym$$ forment une base linéaire de $\obs$. Le \textbf{degré} d'une observable de diagrammes est obtenu en transportant par l'isomorphisme $\Lambda \simeq \obs$ le degré standard des fonctions symétriques. Autrement dit, $\deg p_{\mu}=|\mu|$. Ainsi, on dispose d'une première graduation sur l'algèbre $\obs$.\bigskip
\bigskip

Examinons maintenant une seconde base. De fa\c con analogue à ce qui précède,  on définit les \textbf{moments entrelacés} d'un diagramme $\lambda$ par :
$$\tilp_{k}(\lambda)=p_{k}(X-Y)=\sum_{i=1}^v (x_{i})^k - \sum_{i=1}^{v-1} (y_{i})^k\,.$$
En particulier, $\tilp_{1}(\lambda)=\sum_{i=1}^v x_{i} - \sum_{i=1}^{v-1} y_{i}=0$ pour tout diagramme, car cette quantité est aussi la somme des pentes de $\lambda(s)-|s|$. Soit 
\begin{align*}
G_{\lambda}(z)&=\frac{\prod_{i=1}^{v-1} (z-y_{i})}{\prod_{i=1}^{v} (z-x_{i})}=\frac{1}{z}\, \exp\left(\sum_{k=1}^{\infty}\frac{\tilp_{k}(\lambda)}{k}\,z^{-k}\right)\\
\Phi_{\lambda}(z)&=\prod_{i=1}^{d}\,\frac{z+b_{i}}{z-a_{i}}=\exp\left(\sum_{k=1}^{\infty} \frac{p_{k}(\lambda)}{k}\,z^{-k}\right)
\end{align*}
les fonctions génératrices des coordonnées entrelacées et des coordonnées de Frobenius du diagramme. On peut montrer que $z\,G_{\lambda}(z)=\Phi_{\lambda}(z-1/2)/\Phi_{\lambda}(z+1/2)$ pour tout diagramme $\lambda$, voir \cite[proposition 2.6]{IO02}. Par conséquent, les moments entrelacés $\tilp_{k}$ sont aussi des observables de diagrammes, et on dispose des formules de changement de base 
\begin{align*}
\tilp_n(\lambda)&=\sum_{k=0}^{\lfloor \frac{n-1}{2}\rfloor} \frac{1}{2^{2k}}\,\binom{n}{2k+1} \,p_{n-1-2k}(\lambda)\\
p_n(\lambda)&= \sum_{k=0}^{\lfloor \frac{n}{2} \rfloor} \frac{1}{2^{2k}}\,n^{\downarrow(2k-1)}\,C_k\,\tilp_{n+1-2k}(\lambda)
\end{align*}
où $C_k$ désigne la somme sur toutes les compositions $k=\sum_i k_i$ des inverses des produits $\prod_i -(2k_i+1)!$, et où $n^{\downarrow k}=n(n-1)\cdots (n-k+1)$, avec par convention $n^{\downarrow -1}=1/(n+1)$. On en déduit que les observables $(\tilp_{k})_{k\geq 2}$ forment une base de transcendance de $\obs$. Le \textbf{poids} des observables de diagrammes est la graduation d'algèbre définie par $\mathrm{wt}(\tilp_{k})=k$ pour $k \geq 2$. Compte tenu des formules de changement de base, $\mathrm{wt}(p_{k})=k+1$, et la composante de plus haut poids de $p_{k}$ est $\tilp_{k+1}/(k+1)$.\bigskip
\bigskip

La base des moments entrelacés est particulièrement adaptée à l'étude analytique des diagrammes de Young, et en particulier, elle permet de généraliser la définition des observables de diagrammes. Ainsi, on appellera \textbf{diagramme continu} toute fonction $s \mapsto \omega(s)$ positive, $1$-lipschitzienne et égale à $|s|$ pour $|s|$ assez grand. L'ensemble des diagrammes continus sera noté $\cym$ ; via l'interprétation fonctionnelle des diagrammes de Young, $\ym$ se plonge dans $\cym$. Si $\omega \in \cym$, on note $\sigma_{\omega}(s)=(\omega(s)-|s|)/2$ ; c'est une fonction $1$-lipschitzienne à support compacte. Les moments entrelacés d'un diagramme de Young continu sont définis par 
$$\tilp_{k}(\omega)= \int_{\R}s^{k}\, \sigma_{\omega}''(s)\,ds\,,$$
la dérivée s'entendant éventuellement au sens des distributions. On retrouve les moments entrelacés usuels pour des diagrammes de Young de partitions, et ceci permet d'étendre l'ensemble de définition d'une observable de $\ym$ à $\cym$. D'autre part, si $\omega$ est un diagramme continu et si $t$ est un nombre réel strictement positif, on note 
$$\omega^{t}(s)=\sqrt{t}\,\,\omega\big(s/\sqrt{t}\big)$$
le \textbf{diagramme renormalisé} en abscisse et en ordonnée d'un facteur $\sqrt{t}$, voir l'exemple de la figure \ref{renormalisation}.
\figcapt{\psset{unit=1mm}\pspicture(-50,0)(50,50)
\psline{->}(-50,0)(50,0)
\psline{->}(0,0)(0,50)
\parametricplot[border=2mm,bordercolor=white]{-2}{2}{t 20 mul 
                       t t 0.5 mul arcsin mul 0.017453 mul 4 t t mul neg add sqrt add 12.732395 mul t t mul t mul t mul t t mul 8 mul neg add 16 add 0.5 neg mul t mul add}
\parametricplot[border=2mm,bordercolor=white]{-2}{2}{t 20 mul 0.5 mul
                       t t 0.5 mul arcsin mul 0.017453 mul 4 t t mul neg add sqrt add 12.732395 mul t t mul t mul t mul t t mul 8 mul neg add 16 add 0.5 neg mul t mul add 0.5 mul}
\psline(-47,47)(-20,20)
\psline(47,47)(20,20)
\psline[linewidth=0.25pt](20,20)(0,0)(-20,20)
\rput(-35,39.5){$\omega$}
\rput(-13,21){$\omega^{1/4}$}
\endpspicture}{Renormalisation $\omega^t$ d'un diagramme continu $\omega$. L'aire entre la courbe et la fonction valeur absolue est multipliée par $t$. \label{renormalisation}}{Renormalisation $\omega^{t}$ d'un diagramme continu $\omega$}
Dans ce cadre, si $f$ est une observable de diagrammes homogène de poids $k$, alors $f(\omega^{t})=t^{k/2}\,f(\omega)$. Le poids est donc la graduation de $\obs$ adaptée à la renormalisation <<~isotrope~>> de diagrammes équilibrés ; dans le chapitre \ref{qplancherelmeasure}, nous verrons que le degré est une graduation adaptée à la renormalisation de certains diagrammes non équilibrés.

\section{Permutations partielles et caractères centraux}\label{centralcharacter}

Les caractères centraux forment une autre base de l'algèbre d'observables, et ils permettent d'établir le lien avec la théorie des représentations des groupes symétriques (\cite{IO02,Sni06a}). L'idée est d'associer à toute partition $\mu$ la fonction
$$\lambda \mapsto \varsigma^{\lambda}(\sigma_{\mu \sqcup 1^{n-|\mu|}})\,,$$
où $\varsigma^{\lambda}$ désigne le caractère irréductible du module de Specht $V^{\lambda}$, et avec $n=|\lambda|$. Ces applications doivent néanmoins être renormalisées, et d'autre part, il faut pouvoir traiter convenablement le cas où $|\mu| > n$. Pour ces deux raisons, il est utile de voir les caractères centraux comme éléments de l'algèbre des permutations partielles. On appelle \textbf{permutation partielle} de taille $n$ un couple $(\sigma,S)$, où $\sigma \in \sym_{n}$ et $S$ est une partie de $\lle 1,n \rre$ telle que $\sigma(x)=x$ pour tout $x \in \lle 1,n\rre \setminus S$. Alternativement, on peut voir une permutation partielle de taille $n$ comme la donnée d'un support $S \subset \lle 1,n\rre$ et d'une permutation $\sigma \in \sym(S)$. Le produit de deux permutations partielles est défini par
$$(\sigma,S)\,(\tau,T)=(\sigma\tau,S \cup T)\,,$$
et les permutations partielles de taille $n$ forment ainsi un monoïde non commutatif, de neutre 
$(\id_{\lle 1,n\rre},\emptyset)$. Nous noterons $\blg_{n}$ l'algèbre (complexe) de ce monoïde. La théorie des représentations de cette algèbre sera précisée dans la section \ref{ivanovkerov}, et l'on verra dans le chapitre \ref{bundle} que cette construction rentre dans un cadre très général. Pour l'instant, nous aurons seulement besoin de la limite projective $\blg_{\infty}$ des algèbres $\blg_{n}$. Si $N \geq n$, on peut définir une projection linéaire $\phi_{N,n} : \blg_{N} \to \blg_{n}$ en posant
$$\phi_{N,n} (\sigma,S)= \begin{cases} (\sigma,S)& \text{si }S\subset \lle 1,n\rre, \\0 &\text{sinon.} \end{cases}$$
Les applications $\phi_{N,n}$ sont des morphismes d'algèbres, et elles sont compatibles entre elles. De plus, si le degré d'une permutation partielle est défini par $\deg(\sigma,S)=\card S$, alors les applications $\phi_{N,n}$ sont compatibles avec les filtrations d'algèbres associées au degré. On note $\blg_{\infty}=\varprojlim_{n \to \infty}\blg_{n}$ la limite projective des algèbres de permutations partielles dans la catégorie des algèbres filtrées ; ses éléments sont les combinaisons linéaires formelles éventuellement infinies de permutations partielles de degrés bornés.\bigskip
\bigskip
 
Fixons maintenant une partition $\mu$, et notons $\varSigma_{\mu}$ la somme formelle infinie des permutations partielles 
$$(a_{1,1},\ldots,a_{1,\mu_{1}})(a_{2,1},\ldots,a_{2,\mu_{2}})\,\cdots\,(a_{r,1},\ldots,a_{r,\mu_{r}}) \,\,,\,\,\big\{ a_{i,j}\,\,\big|\,\,1 \leq i \leq \ell(\mu),\,\,1\leq j \leq \mu_{i} \big\}$$
où $a$ est n'importe quelle fonction injective $\{(i,j) \,\,|\,\,1 \leq i \leq \ell(\mu),\,\,1\leq j \leq \mu_{i} \} \to \N^{*} $. Si $\phi_{n}$ est la projection canonique $\blg_{\infty}\to \blg_{n}$, alors le projeté $\phi_{n}(\varSigma_{\mu})$ est donné par la même somme, mais avec l'indice de sommation $a$ restreint aux fonctions injectives à valeurs dans $\lle 1,n\rre$. Soit $\pi_{n}$ la projection de $\blg_{n}$ vers $\C\sym_{n}$ qui à une permutation partielle $(\sigma,S)$ associe la permutation $\sigma$. Le projeté $\varSigma_{\mu,n}=\pi_{n}(\phi_{n}(\varSigma_{\mu}))$ vaut $0$ lorsque $n < |\mu|$, et dans le cas contraire, c'est un multiple de la classe de conjugaison $C_{\mu \sqcup 1^{n-|\mu|}}$ :
$$\varSigma_{\mu,n}= n(n-1)\cdots(n-|\mu|+1)\,\widetilde{C}_{\mu \sqcup 1^{n-|\mu|}}=n^{\downarrow |\mu|}\,\widetilde{C}_{\mu\sqcup 1^{n-|\mu|}}\,,$$
où pour toute partition $\lambda$ de taille $n$, $\widetilde{C}_{\lambda}=C_{\lambda} /\card C_{\lambda}$. Ceci étant, le \textbf{caractère central} $\varSigma_{\mu}$ est l'observable de diagrammes définie par 
$$\varSigma_{\mu}(\lambda \in \ym_{n})= \chi^{\lambda}(\varSigma_{\mu,n})\,,$$
où $\chi^{\lambda}$ désigne le caractère irréductible renormalisé, c'est-à-dire que $(\dim \lambda)\,\chi^{\lambda}=\varsigma^{\lambda}$. L'appartenance de $\varSigma_{\mu}$ à $\obs$ est en réalité tout à fait non triviale, et elle découle des deux résultats suivants :\vspace{2mm}
\begin{enumerate}
\item Dans $\blg_{\infty}$, les éléments $\varSigma_{\mu}$ engendrent une sous-algèbre commutative $\alg_{\infty}$, et $\mathrm{wt}( \varSigma_{\mu})=|\mu|+\ell(\mu)$ est une graduation d'algèbre sur $\alg_{\infty}$. Plus précisément, 
$$\varSigma_{\mu_{1}} \,*\,\varSigma_{\mu_{2}} = \varSigma_{\mu_{1}\sqcup \mu_{2}} + 
\left(\substack{\text{combinaison linéaire de termes de poids} \\ \text{inférieur à }|\mu_{1}|+\ell(\mu_{1})+|\mu_{2}|+\ell(\mu_{2})-2}\right)$$
voir \cite[corollaire 3.8]{Sni06a}. Essentiellement, c'est parce que pour $n$ assez grand, une permutation partielle de type $\mu_{1}$ et une autre permutation partielle de type $\mu_{2}$ sont génériquement à supports disjoints, donc commutent et ont pour produit une permutation partielle de type $\mu_{1}\sqcup\mu_{2}$.\vspace{2mm}
\item D'autre part, pour tout groupe fini $G$, si $a$ et $b$ sont deux éléments du centre $(\C G)^{G}$ de l'algèbre du groupe, et si $\chi$ est un caractère irréductible normalisé de $G$, alors $\chi(ab)=\chi(a)\,\chi(b)$. En effet, dans la décomposition de $\C G$ en blocs matriciels, $(\C G)^{G}$ s'identifie à la somme directe d'espaces d'homothéties $\bigoplus_{V \in \widehat{G}} \C\id_{V}$, et la restriction de $\chi^{V}$ à $(\C G)^{G}$ est la projection sur $\C\id_{V}$ ;   c'est donc bien un morphisme d'algèbres. Par suite, pour toutes partitions $\mu_{1},\mu_{2},\lambda$,
\begin{align*}
\varSigma_{\mu_{1}}(\lambda)\times_{\C} \varSigma_{\mu_{2}}(\lambda)&= \chi^{\lambda} (\varSigma_{\mu_{1},n})\times_{\C} \chi^{\lambda}(\varSigma_{\mu_{2},n})=\chi^{\lambda}( \varSigma_{\mu_{1},n} \times_{\C\sym_{n}} \varSigma_{\mu_{2},n}) \\
&= \chi^{\lambda}(\pi_{n}\circ \phi_{n}(\varSigma_{\mu_{1}} \times_{\alg_{\infty}} \varSigma_{\mu_{2}})) = (\varSigma_{\mu_{1}} \times_{\alg_{\infty}} \varSigma_{\mu_{2}})(\lambda)
\end{align*}
c'est-à-dire que le produit d'observables $\varSigma_{\mu_{1}}$ et $\varSigma_{\mu_{2}}$ peut être effectué directement dans $\alg_{\infty}\subset \blg_{\infty}$. Par conséquent, compte tenu de la décomposition graduée des produits de caractères centraux $\varSigma_{\mu_{1}} * \varSigma_{\mu_{2}}$, il suffit de montrer que les caractères centraux des cycles $\varSigma_{k}$ sont des observables de diagrammes.\vspace{2mm}
\item En utilisant la formule de Frobenius pour $\mu=k1^{n-k}$, on peut montrer que $\varSigma_{k}(\lambda)$ est le coefficient de $z^{-1}$ dans la série de Laurent\label{wassermann}
$$-\frac{1}{k}\,\left(z-\frac{1}{2}\right)^{\downarrow k}\,\frac{\Phi_{\lambda}(z)}{\Phi_{\lambda}(z-k)}\,,$$
d'où une expression de $\varSigma_{k}$ en fonction des $p_{l \leq k}$, voir \cite[proposition 3.2]{IO02}. Par exemple, $\varSigma_{1}=p_{1}$, $\varSigma_{2}=p_{2}$, $\varSigma_{3}=p_{3}-\frac{3}{2}p_{11}+\frac{5}{4}p_{1}$ et $\varSigma_{4}=p_{4}-4p_{21}+\frac{11}{2}p_{2}$. Ainsi, les $\varSigma_{k}$ et les $\varSigma_{\mu}$ sont bien dans $\obs$.
\end{enumerate}\bigskip\bigskip

À partir de ces observations, on peut établir la proposition suivante, voir \cite[corollaire 4.3, corollaire 4.4 et proposition 4.9]{IO02} :
\begin{proposition}[Base des caractères centraux]\label{whitelight}
Les observables $(\varSigma_{\mu})_{\mu \in \ym}$ forment une base linéaire de $\obs$, et les caractères centraux des cycles $(\varSigma_{k})_{k \geq 1}$ forment une base de transcendance. La composante de plus haut degré de $\varSigma_{\mu}$ est toujours égale à $p_{\mu}$. Par suite,
$$\varSigma_{\mu_{1}} \,*\,\varSigma_{\mu_{2}}=\varSigma_{\mu_{1}\sqcup \mu_{2}}+(\text{termes de degré inférieur})\,.$$
La filtration de l'algèbre $\alg_{\infty}$ précédemment définie correspond à la filtration des poids sur l'algèbre d'observables $\obs$ ; par conséquent, on a aussi
$$\varSigma_{\mu_{1}} \,*\,\varSigma_{\mu_{2}}=\varSigma_{\mu_{1}\sqcup \mu_{2}}+(\text{termes de poids inférieur})\,.$$
\end{proposition}
\noindent Ces propriétés de factorisation des caractères centraux en plus haut degré ou en plus haut poids joueront un rôle essentiel pour l'asymptotique des caractères et des mesures sur les partitions\footnote{Nous décrirons plus en détail les termes du produit $\varSigma_{\mu_{1}}\,*\,\varSigma_{\mu_{2}}$ dans la section \ref{qgaussian} ; ainsi, le produit peut être écrit comme somme sur des appariements entre les points des cycles de type $\mu_{1}$ et les points de cycles de type $\mu_{2}$, le terme $\varSigma_{\mu_{1}\sqcup \mu_{2}}$ correspondant à l'appariement vide.}. Concluons cette section en dressant une liste des diverses descriptions isomorphiques de l'algèbre d'observables $\obs$ :\vspace{2mm}

\begin{enumerate}
\item Les moments de Frobenius $p_{\lambda}$ définissent un isomorphisme canonique $\obs \simeq \Lambda$, ce qui permet d'effectuer des calculs dans l'algèbre des fonctions symétriques.\vspace{2mm}

\item Les caractères centraux $\varSigma_{\mu}$ sont issus de l'algèbre des permutations partielles $\blg_{\infty}$, et plus précisément de l'algèbre $\alg_{\infty}$ des invariants pour l'action par conjugaison de $\sym_{\infty}$. Cette sous-algèbre commutative est l'\textbf{algèbre d'Ivanov-Kerov}, voir la section \ref{ivanovkerov}. On dispose donc d'un isomorphisme canonique $\obs\simeq \alg_{\infty}$, ce qui permet en particulier de calculer les coefficients de structure de l'algèbre $\obs$ dans la base des caractères centraux, \emph{cf.} \cite[proposition 4.5]{IO02} et \cite[proposition 6.2 et théorème 9.1]{IK99}. Nous reviendrons sur ce point dans le chapitre \ref{badbeat}.\vspace{2mm}

\item Enfin, les caractères centraux $\varSigma_{\mu}$ peuvent également être vus comme éléments de l'algèbre des fonctions symétriques décalées, voir \cite{OO98}. Un polynôme symétrique décalé en $m$ variables est un polynôme $p(x_{1},x_{2},\ldots,x_{m})$ symétrique en les nouvelles variables $y_{i}=x_{i}-i+\mathrm{const}$, la valeur de la constante étant arbitraire (on la prendra par exemple égale à $0$). Les polynômes symétriques décalés en $m$ variables forment une algèbre $\Lambda^\sharp(m)$, et on appelle \textbf{fonction symétrique décalée} un élément de la limite projective $\Lambda^\sharp= \varprojlim_{m \to \infty}\Lambda^\sharp(m)$. Une base de cette algèbre est formée des fonctions de Schur décalées :
$$s^\sharp_{\mu}(x_{1},x_{2}, \ldots,x_{m})=\frac{\det\, ((x_{i}+m-i)^{\downarrow \mu_{j}+m-j})_{i,j}}{\det\, ((x_{i}+m-i)^{\downarrow m-j})_{i,j}} \,\,\text{ si }\ell(\mu)\leq m,\text{ et }0\text{ sinon.}$$
On définit les fonctions sommes de puissances décalées en utilisant la formule de Frobenius : $$\forall \mu \in \ym_{n},\,\,\,p_{\mu}^{\sharp}=\sum_{\lambda \in \ym_{n}} \varsigma^{\lambda}(\mu)\,s_{\lambda}^{\sharp}\,.$$
Alors, on peut montrer que la fonction symétrique décalée $p_{\mu}^{\sharp}(\lambda)$ en les coordonnées standards $\lambda_{1},\ldots,\lambda_{r}$ d'un diagramme est égale à $\varSigma_{\mu}(\lambda)$, \emph{cf.} \cite[théorème 8.1]{OO98} et \cite[théorème 9.1]{IK99}. Cette identité fournit un troisième isomorphisme $\obs \simeq \Lambda^{\sharp}$.\vspace{2mm}
\end{enumerate}
Notons que l'algèbre $\Lambda^\sharp$ est également mise en jeu dans un analogue du théorème de Frobe\-nius-Schur pour les groupes linéaires complexes $\GL(m,\C)$. Ainsi, on peut montrer que $\Lambda^\sharp(m)$ est naturellement isomorphe au centre $Z(U(\mathfrak{gl}(m,\C)))$ de l'algèbre enveloppante de l'algèbre de Lie $\mathfrak{gl}(m,\C)$ --- c'est l'\textbf{isomorphisme d'Harish-Chandra}, voir par exemple \cite[\S2]{OO98}.\bigskip

\section{Cumulants libres et asymptotique des caractères}\label{freecumulant}
La dernière base importante de l'algèbre $\obs$ est la base des cumulants libres ; elle a été introduite par P. Biane dans \cite{Bia98,Bia03a}, et elle est inspirée par la théorie des probabilités libres de Voiculescu (\cite{VDN92,Bia03b}). Si $\lambda$ est un diagramme de Young, on appelle \textbf{mesure de transition} de $\lambda$ la mesure de probabilité
$$\mu_{\lambda}=\sum_{k=1}^{v} \left(  \frac{\prod_{i=1}^{v-1}x_{k}-y_{i}}{\prod_{i \neq k}x_{k}-x_{i}}\right) \,\delta_{x_{k}}\,.$$
Alors, $G_{\lambda}(z)=\int_{\R} \frac{1}{z-s}\,\mu_{\lambda}(ds)$ est la \textbf{transformée de Cauchy} de la mesure de transition de $\lambda$. De plus, 
$$G_{\lambda}(z)=\frac{1}{z}\,\exp\left(-\int_{\R} \frac{\sigma_{\lambda}'(s)}{z-s}\,ds\right),$$
et compte tenu de cette expression analytique, on peut aussi définir la mesure de transition $\mu_{\omega}$ d'un diagramme continu $\omega \in \cym$. Nous expliquerons dans la section \ref{plancherelprocess} l'origine de la terminologie <<~mesure de transition~>>.\bigskip
\bigskip

Ceci étant, pour toute mesure de probabilité $\mu$ à support compact $S\subset \R$, on appelle \textbf{R-transformée} de $\mu$ la fonction réciproque de la transformée de Cauchy $C_{\mu}$ au voisinage de l'infini :
$$C_{\mu}(z)=\int_{\R}\frac{\mu(ds)}{z-s}\sim_{\infty} \frac{1}{z} \qquad;\qquad R_{\mu}(C_{\mu}(z))=z \quad\text{et}\quad R_{\mu}(z)=\frac{1}{z}\left(1+\sum_{k=1}^{\infty} R_{k}(\mu)\,z^{k}\right)\,.$$
Les coefficients $R_{k}(\mu)$ sont appelés \textbf{cumulants libres} de la mesure $\mu$. En théorie des probabilités non commutatives, la R-transformée et les cumulants libres jouent un rôle analogue à la transformée de Fourier et aux cumulants standards des variables aléatoires. Les cumulants libres d'un diagramme sont les cumulants libres de sa mesure de transition : $R_{k}(\lambda)=R_{k}(\mu_{\lambda})$. Le premier cumulant est toujours nul, et par inversion de Lagrange, 
$$R_{k+1}(\lambda)= -\frac{1}{k}\, [z^{-1}]\,  (G_{\lambda}(z))^{-k}$$
pour tout entier $k$. On en déduit que les cumulants libres $R_{k}$ sont des observables de diagrammes, et d'autre part, la définition analytique des cumulants prouve que
$$R_{k}(\omega^{t})=t^{k/2}\,R_{k}(\omega)$$
pour tout diagramme continu $\omega$ et tout entier $k \geq 2$. Par conséquent, $R_{k+1}$ est une observable homogène de poids $k+1$, et c'est en fait la composante homogène de plus haut poids $k+1$ de $\varSigma_{k}$, car
$$\varSigma_{k}(\lambda)= - \frac{1}{k}\, [z^{-1}] \,\big(G_{\lambda}(z)\,G_{\lambda}(z-1)\,\cdots\,G_{\lambda}(z-k+1)\big)^{-1}$$
 \emph{cf.} \cite[section 5]{Bia03a} et \cite[proposition 10.1]{IO02}. Notons $R_{\mu}$ un produit $R_{\mu_{1}}R_{\mu_{2}}\cdots R_{\mu_{r}}$ de cumulants libres de diagrammes. Alors, compte tenu de la propriété de factorisation en plus haut poids des caractères centraux,
$$\varSigma_{\mu}=R_{\mu+1} + (\text{termes de poids inférieur})$$
pour toute partition $\mu$, où $\mu +1$ désigne la partition obtenue à partir de $\mu$ en ajoutant $1$ à toutes ses parts. En particulier, les $(R_{k+1})_{k \geq 1}$ forment une base de transcendance de $\obs$, et les $(R_{\mu+1})_{\mu \in \ym}$ forment une base linéaire. Dans la formule précédente, on peut montrer que les termes de poids inférieur à $|\mu|+\ell(\mu)-1$ sont en réalité de poids inférieur à $|\mu|+\ell(\mu)-2$ ; en effet, tous les poids apparaissant ont la même parité que $|\mu|+\ell(\mu)$. \bigskip
\bigskip

L'identité reliant $\varSigma_{\mu}$ à $R_{\mu+1}$ est à la base de résultats asymptotiques importants pour les caractères du groupe symétrique. \'Etant donnée une constante $A>1$, une partition $\lambda$ est dite \textbf{$A$-équilibrée} si sa longueur $\ell(\lambda)$ et sa plus grande part $\lambda_{1}$ sont inférieures à $A\sqrt{n}$, où $n=|\lambda|$. Alors, le diagramme continu $\omega=(\omega_{\lambda})^{1/|\lambda|}$ vérifie $\omega(s)=|s|$ pour tout $s \geq A$. En particulier, les diagrammes continus obtenus par renormalisation de diagrammes $A$-équilibrés forment une partie relativement compacte de $\cym$, et toute observable de diagrammes (continus) est bornée sur cette partie. On en déduit que 
\begin{align*}| \varSigma_{\mu}(\lambda) - R_{\mu+1}(\lambda)|&=\left|\sum_{k=0}^{|\mu|+\ell(\mu)-2}f_{k}(\lambda)\right|\leq \sum_{k=0}^{|\mu|+\ell(\mu)-2}|f_{k}(\lambda)|\\
&\leq n^{\frac{|\mu|+\ell(\mu)-2}{2}}\, \sum_{k=0}^{|\mu|+\ell(\mu)-2}|f_{k}(\omega)| \leq C\,n^{\frac{|\mu|+\ell(\mu)-2}{2}}
\end{align*}
où $C$ est une constante qui ne dépend que de $A$ (et de $\mu$). Dans ce qui suit, si $\mu$ est une partition et si $n$ est un entier plus grand que $|\mu|$, on note $\sigma_{\mu,n}$ une permutation de $\sym_{n}$ de type $\mu1^{n-|\mu|}$, et $|\sigma_{\mu,n}|$ le nombre minimal de transpositions dans une factorisation de $\sigma_{\mu,n}$ en cycles, c'est-à-dire $|\mu|-\ell(\mu)$. On déduit de ce qui précède :
\begin{theorem}[Expression asymptotique des caractères des groupes symétriques, \cite{Bia98}]
Pour tout $A>1$ et toute partition $\mu$, il existe une constante positive $C(A,\mu)$ telle que  $$|\chi^{\lambda}(\sigma_{\mu,n}) - n^{-|\mu|} \,R_{\mu+1}(\lambda)| \leq C(A,\mu)\,\, n^{-\frac{|\sigma_{\mu,n}|}{2}-1}$$
pour tout diagramme $A$-équilibré $\lambda$ de taille $n\geq |\mu|$. En particulier, $\chi^{\lambda}(\sigma_{\mu,n})=O(n^{-\frac{|\sigma_{\mu,n}|}{2}})$, et d'autre part, l'estimation précédente implique la factorisation asymptotique des caractères. Ainsi, si $\sigma_{\mu,n}$ et $\sigma_{\tau,n}$ commutent, alors
$$ |\chi^{\lambda}(\sigma_{\mu,n}\,\sigma_{\tau,n} )- \chi^{\lambda}(\sigma_{\mu,n})\,\chi^{\lambda}(\sigma_{\tau,n})| \leq D(A,\mu,\tau) \,\, n^{-\frac{|\sigma_{\mu,n}|+|\sigma_{\tau,n}|}{2}-1}$$
pour tout diagramme $A$-équilibré $\lambda$ de taille $n \geq |\mu|+|\tau|$.
\end{theorem}\bigskip
\bigskip

L'expression asymptotique des caractères indexés par des partitions $A$-équilibrées implique également des résultats de concentration pour la forme d'une partition apparaissant dans la décomposition d'un $\C\sym_{n}$-module construit par produit tensoriel, restriction ou induction, voir \cite[théorèmes 1.4.1, 1.5.1 et 1.6.1]{Bia98}, et \cite{Bia01a,Sni06b}. D'autre part, le lien entre caractères centraux et cumulants libres permet de démontrer très facilement des lois limites pour les grandes partitions aléatoires, voir en particulier les sections \ref{lskv} et \ref{schurplus} ; c'est l'intérêt essentiel de cette notion. Concluons ce paragraphe en évoquant deux propriétés additionnelles des cumulants libres. \vspace{2mm}
\begin{enumerate}
\item Notons $M_{\mu}(\lambda)$ les moments de la mesure de transition d'un diagramme $\lambda$ ; ce sont des observables de diagrammes, car si 
$$ G_{\lambda}(z)=\frac{1}{z}\,\exp\left(\sum_{k=1}^{\infty} \frac{\tilp_{k}(\lambda)}{k}\,z^{-k}\right)=\frac{1}{z}\,\left(1+\sum_{k=1}^\infty \widetilde{h}_{k}(\lambda)\,z^{-k}\right),$$
alors $M_{\mu}(\lambda)=\widetilde{h}_{\mu}(\lambda)$ pour toute partition $\mu$. On doit à R. Speicher une interprétation combinatoire des relations entre moments et cumulants libres d'une mesure de probabilité, \emph{cf.} \cite{Spe98,NS06}. Si $\pi$ est une partition ensembliste de l'ensemble $\lle 1,n\rre$, on note $M_{\pi}=M_{\rho}$ et $R_{\pi}=R_{\rho}$ les moments et cumulants libres indexés par la partition $\rho$ obtenue en prenant les tailles des parts de $\pi$. La partition $\pi$ est dite \textbf{non croisée} s'il n'existe pas de quadruplets $i<j<k<l$ tels que $i \sim_{\pi} k$ et $j \sim_{\pi} l$. Si l'on place les entiers $1,2,\ldots,n$ sur une droite et si l'on dessine les blocs de la partition $\pi$, alors la condition de non croisement a une signification géométrique évidente, voir la figure \ref{noncrossing}. L'ensemble $\mathfrak{N}\mathfrak{C}(n)$ des partitions non croisées de $\lle 1,n\rre$ est de cardinal le nombre de Catalan $C_{n}$, et il forme un treillis ordonné pour l'ordre de raffinement sur les partitions d'ensemble. Dans ce cadre, on peut montrer que 
$$M_{k}=\sum_{\pi \in \mathfrak{N}\mathfrak{C}(k)} R_{\pi}\qquad;\qquad R_{k}=\sum_{\pi \in \mathfrak{N}\mathfrak{C}(k)} \mu(\pi_{\mathrm{comp}})\,M_{\pi}$$
où $\mu(\pi)=\prod_{i\geq1} (-1)^{|\pi_{i}|-1}\,C_{|\pi_{i}|-1}$ est la fonction de M\"obius du treillis $\mathfrak{N}\mathfrak{C}(k)$, et $\pi_{\text{comp}}$ désigne le complémentaire de $\pi$ au sens de Kreweras, voir par exemple \cite[\S2.2.4]{Sni06a}. Ainsi, au prix d'une inversion de M\"obius, on peut ramener l'étude des cumulants libres de diagrammes $R_{\mu}(\lambda)$ à celle des moments $M_{\mu}(\lambda)$. Or, ces derniers peuvent être vus comme \textbf{moments des éléments de Jucys-Murphy} au sens des probabilités non commutatives. On renvoie à \cite{Bia98,Bia03a} et \cite[\S2.1.8 et \S4.2]{Sni06a} pour de plus amples détails au sujet de cette interprétation.
 \figcapt{\psset{unit=1mm}\pspicture(0,-5)(110,27)
 \psarc*[linecolor=MidnightBlue](25,0){25}{0}{180}
 \psarc(25,0){25}{0}{180}
 \psarc*[linecolor=MidnightBlue](85,0){25}{0}{180}
 \psarc(85,0){25}{0}{180}
 \psarc*[linecolor=white](80,0){20}{0}{180}
 \psarc(80,0){20}{0}{180}
 \psarc*[linecolor=white](105,0){5}{0}{180}
 \psarc(105,0){5}{0}{180}
 \psarc*[linecolor=white](15,0){15}{0}{180}
 \psarc(15,0){15}{0}{180}
 \psarc*[linecolor=white](35,0){5}{0}{180}
 \psarc(35,0){5}{0}{180}
 \psarc*[linecolor=white](45,0){5}{0}{180}
 \psarc(45,0){5}{0}{180}
 \psarc*[linecolor=MidnightBlue](15,0){5}{0}{180}
 \psarc(15,0){5}{0}{180}
 \psarc*[linecolor=white](15,-10){11.4}{63}{117}
 \psarc(15,-10){11.4}{63}{117}
 \psarc*[linecolor=MidnightBlue](80,0){10}{0}{180}
 \psarc(80,0){10}{0}{180}
 \psarc*[linecolor=white](80,-10){14.14}{45}{135}
 \psarc(80,-10){14.14}{45}{135}
 \pscircle*[linecolor=MidnightBlue](80,1.25){1.25}
 \pscircle(80,1.25){1.25}
 \psline(0,0)(110,0)
 \psdots(0,0)(10,0)(20,0)(30,0)(40,0)(50,0)(60,0)(70,0)(80,0)(90,0)(100,0)(110,0)
 \rput(0,-3){$1$}
 \rput(10,-3){$2$}
 \rput(20,-3){$3$}
 \rput(30,-3){$4$}
 \rput(40,-3){$5$}
 \rput(50,-3){$6$}
 \rput(60,-3){$7$}
 \rput(70,-3){$8$}
 \rput(80,-3){$9$}
 \rput(90,-3){$10$}
 \rput(100,-3){$11$}
 \rput(110,-3){$12$}
\endpspicture}{La partition $\pi=\{1,4,5,6\}\sqcup\{7,11,12\}\sqcup\{2,3\}\sqcup\{8,10\}\sqcup\{9\}$ est non croisée.\label{noncrossing}}{Partition non croisée}

\vspace{2mm}
\item La composante de plus haut poids de $\varSigma_{k}$ étant le cumulant $R_{k+1}$, on peut montrer que tout caractère central (d'un cycle) s'écrit comme polynôme à coefficients entiers en les cumulants d'ordre inférieur à son poids. Par exemple,
\begin{align*}&\varSigma_{1}=R_{2}\qquad; \qquad\varSigma_{2}=R_{3}\qquad;\qquad \varSigma_{3}=R_{4}+R_{2} \quad;\\
&  \varSigma_{4}=R_{5}+3R_{3}\qquad;\qquad\varSigma_{5}=R_{6}+15R_{4}+5R_{2}^2+8R_{2}\,.\end{align*}
Les termes suivants sont calculés à la fin de \cite{Bia03a}. On appelle \textbf{polynôme de Kerov} le polynôme universel $K_{k}$ tel que $\varSigma_{k}=K_{k}(R_{2},\ldots,R_{k+1})$. Le polynôme de Kerov associé au caractère central d'un cycle a tous ses coefficients positifs (\emph{cf.} \cite{Fer07}), et on peut donner une interprétation combinatoire des dits coefficients en termes de nombres de factorisations d'une permutation, voir \cite[théorème 1.4]{DFS08}.\vspace{2mm}
\end{enumerate}
Un dernier point que nous avons omis est la détermination de toutes les filtrations d'algèbre de $\obs$ ; ainsi, une classe très générale de telles filtrations est décrite par \cite[proposition 10.3]{IK99} et \cite[proposition 4.7]{IO02}, et nous renvoyons le lecteur à ces articles s'il souhaite pallier cette omission. Outre le degré canonique et le poids, nous aurons par la suite besoin du degré de Kerov (section \ref{sniady} et chapitre \ref{gelfandmeasure}) et du $\alpha$-degré (section \ref{schurminus}) ; nous les présenterons en temps voulu.\bigskip

\section{Correspondance de Markov-Krein et topologie des diagrammes continus}\label{markovkrein}
Pour conclure ce chapitre, nous souhaitions préciser très clairement la topologie induite sur $\cym$ par la convergence des observables de diagrammes de l'algèbre $\obs$. Ainsi, étant donnée une suite --- pour l'instant déterministe --- $(\omega_{n})_{n \in \N}$ de diagrammes continus, que peut-on dire si pour toute observable $f \in \obs$, $f(\omega_{n})$ admet une limite ? La réponse rentre dans le cadre de la \textbf{correspondance de Markov-Krein}, qui est exposée dans son cadre le plus général dans les articles \cite{Ker93b,Ker98}. Nous suivrons ici peu ou prou le second article.\bigskip\bigskip

Si $\omega$ est un diagramme de Young continu dans $\cym$, nous avons vu dans les paragraphes précédents comment lui associer :\vspace{2mm}
\begin{itemize}
\item une fonction $G_{\omega}(z)$ analytique sur le demi-plan de Poincaré $\mathbb{H}=\{z\in \C\,\,|\,\,\Im(z)>0\}$ : 
$$G_{\omega}(z)=\frac{1}{z}\exp\left(-\int_{\R}\frac{\sigma_{\omega}'(s)}{z-s}\,ds\right)\,.$$
La fonction $G_{\omega}(z)$ a sa partie imaginaire négative lorsque $z$ reste dans $\mathbb{H}$ ; de plus, $\lim_{y \to +\infty} \I y \,G_{\omega}(\I y)=1$. Enfin, elle admet une limite finie lorsque $z$ tend vers $t \in \R$ avec $t$ en dehors du support de $\sigma_{\omega}$.
\item une mesure de probabilité à support compact $\mu_{\omega}$ dont la transformée de Cauchy est la fonction $G_{\omega}(z)$ : 
$$G_{\omega}(z)=\int_{\R}\frac{\mu_{\omega}(ds)}{z-s}\,.$$ 
\end{itemize}
On obtient ainsi une bijection entre $\cym$, les mesures de probabilité sur $\R$ à support compact, et les fonctions analytiques sur $\mathbb{H}$ qui vérifient
\begin{align*}&\forall z \in \mathbb{H},\,\,\Im(N(z))\leq 0\qquad;\qquad \lim_{y \to +\infty} \I y \,N(\I y)=1\quad;\\
&\forall t \text{ réel en dehors d'un certain compact}, \,\,\lim_{z \to t,\,\,z\in \mathbb{H}} N(z) \text{ existe}.\end{align*}
Cette correspondance $\omega \leftrightarrow G_{\omega} \leftrightarrow \mu_{\omega}$ est appelée correspondance de Markov-Krein, et elle est homéomorphique pour la topologie de la convergence uniforme sur $\cym$, la topologie de la convergence en loi des mesures de probabilité, et une topologie que nous préciserons plus loin sur les fonctions analytiques à partie 
imaginaire négative. Malheureusement, ces topologies (métrisables) ne rendent pas les espaces sous-jacents complets ; par suite, il convient d'introduire des objets plus généraux pour comprendre pleinement les propriétés topologiques de cette correspondance. Définissons donc les quatres espaces métrisables complets suivants :\vspace{2mm}
\begin{itemize}
\item \textbf{diagrammes continus généralisés} : soit $\ym^{1}$ l'ensemble des fonctions $\omega : \R \to \R$ qui sont $1$-lipschitziennes et vérifient 
$$\int_{-\infty}^{-1}(1+\omega'(s))\,\frac{ds}{|s|}<\infty\qquad;\qquad \int_{1}^{+\infty}(1-\omega'(s))\,\frac{ds}{|s|}<\infty\,,$$
étant entendu qu'on identifie deux telles fonctions si elles diffèrent d'une constante. Il existe une distance sur $\ym^{1}$ qui en fait un espace métrique complet, et telle que $([\omega_{n}])_{n \in \N}$ converge vers une classe d'équivalence $[\omega]$ dans $\mathscr{Y}^{1}$ si et seulement si l'on peut choisir des représentants tels que $\omega_{n}\to \omega$ uniformément sur toute partie compacte $S \subset \R$. On a bien sûr $\ym \subset \cym \subset \ym^{1}$. \vspace{2mm}
\item \textbf{fonctions de Rayleigh} : soit $\mathscr{R}^{1}$ l'ensemble des fonctions mesurables sur $\R$ telles que 
$$ \forall s,\,\,0\leq R(s)\leq 1\qquad;\qquad\int_{-\infty}^{0}\frac{R(s)}{1+|s|}\,ds<\infty\qquad;\qquad \int_{0}^{\infty}\frac{1-R(s)}{1+|s|}\,ds<\infty\,,$$
étant entendu qu'on identifie deux telles fonctions si elles sont égales presque partout. Une suite de fonctions de Rayleigh $(R_{n})_{n \in \N}$ converge vers une fonction $R$ si
\begin{align*}\forall x \in \R,\,\,\,&\lim_{n \to \infty}\int_{-\infty}^{x}\frac{R_{n}(s)}{1+|s|}\,ds = \int_{-\infty}^{x}\frac{R(s)}{1+|s|}\,ds\quad;\\
\forall x \in \R,\,\,\,&\lim_{n \to \infty}\int_{x}^{\infty}\frac{1-R_{n}(s)}{1+|s|}\,ds=\int_{x}^{\infty}\frac{1-R(s)}{1+|s|}\,ds\,.\end{align*}
Dans ce cas, la limite $R$ est encore une fonction de Rayleigh. Il existe une distance sur $\mathscr{R}^{1}$ compatible avec cette notion de convergence, et complète ; pour cette distance, la convergence $R_{n}\to R$ dans $\mathscr{R}^{1}$ est équivalente\footnote{Attention, si une suite $(R_{n})_{n \in \N}$ de fonctions de Rayleigh converge en ce sens vers une fonction mesurable $R$, il n'y a \emph{a priori} aucune garantie pour que $R$ soit une fonction de Rayleigh dans $\mathscr{R}^{1}$. C'est la raison pour laquelle on doit attacher une distance à cette topologie ; des phénomènes du même type ont lieu pour les deux autres espaces $\ym^{1}$ et $\mathscr{N}^{1}$.} à 
$$\int_{a}^{b} R_{n}(s)\,ds \to \int_{a}^{b} R(s)\,ds$$
pour tout intervalle borné $[a,b]$.\vspace{2mm}
\item \textbf{fonctions à parties imaginaires négatives} : soit $\mathscr{N}^{1}$ l'ensemble des fonctions analytiques sur $\mathbb{H}$ qui ont leur partie imaginaire négative, et vérifient $\lim_{y \to \infty}\I y\,N(\I y)=1$. On munit $\mathscr{N}^{1}$ de la restriction de la topologie de Montel, c'est-à-dire la topologie de la convergence uniforme sur tous les compacts de $\mathbb{H}$. De nouveau, il existe une distance sur $\mathscr{N}^{1}$ compatible avec cette topologie et complète. De plus, si $(N_{n})_{n \in \N}$ suites de fonctions de $\mathscr{N}^{1}$ converge pour la topologie de Montel vers une fonction analytique (à partie imaginaire négative) $f$, alors $f$ appartient à $\mathscr{N}^{1}$ si et seulement si 
$$\lim_{Y \to \infty} \sup_{n \in \N,\,\,y\geq Y}|\I y \,f_{n}(\I y)-1|=0\,.$$
\vspace{1mm}
\item \textbf{mesures de probabilité} : on note $\mathscr{M}^{1}=\mathscr{M}^{1}(\R)=\mathscr{P}(\R)$ l'espace des mesures de probabilité (boréliennes) sur la droite réelle, et on le munit de la topologie de Skorohod de la convergence en loi, \emph{cf.} \cite[chapitre 1]{Bil69}. Il est bien connu qu'il existe une distance (naturelle !) sur $\mathscr{M}^{1}$ compatible avec cette topologie, et qui en fait un espace complet.\vspace{2mm}
\end{itemize}
Avec ces définitions, la correspondance de Markov-Krein généralisée s'énonce comme suit :
\begin{theorem}[Correspondance de Markov-Krein généralisée, \cite{Ker98}]
La correspondance de Markov-Krein entre diagrammes de Young continus, mesures à supports compacts et fonctions analytiques à partie imaginaire négative est la restriction d'une correspondance homéomorphique entre les quatre espaces métriques complets :
$$\ym^{1} \leftrightarrow \mathscr{R}^{1}\leftrightarrow \mathscr{N}^{1} \leftrightarrow \mathscr{M}^{1}\,.$$
\begin{enumerate}
\item La bijection $\mathscr{R}^{1} \leftrightarrow \ym^{1}$ est donnée par :
$$R(s)=\frac{1}{2}(1+\omega'(s))\qquad;\qquad \omega(s)=\int_{0}^{s} (2R(t)-1)\,dt\,.$$
\item L'homéomorphisme $\mathscr{N}^{1} \leftrightarrow \mathscr{R}^{1}$ est fourni par la représentation intégrale de Nevanlinna des fonctions à partie imaginaire négative : 
$$N(z)=\frac{1}{z}\,\exp\left(-\int_{-\infty}^{0}\frac{R(s)}{z-s}\,ds+\int_{0}^{\infty}\frac{1-R(s)}{z-s}\,ds\right)\,.$$
\item La correspondance $\mathscr{M}^{1} \leftrightarrow \mathscr{N}^{1}$ est la transformée de Cauchy des mesures de probabilité, qui peut être inversée par la formule de Perron-Frobenius :
$$N(z)=\int_{\R}\frac{\mu(ds)}{z-s}\qquad;\qquad\int_{a}^{b} \mu(ds)=-\frac{1}{\pi} \lim_{y \to \infty} \int_{a}^{b}\Im(N(s+\I y))\,ds\,.$$
\end{enumerate}
\end{theorem}
\bigskip

Muni de ce résultat abstrait, pla\c cons-nous dans la situation décrite au début de ce paragraphe, c'est-à-dire avec une suite $(\omega_{n})_{n \in \N}$ de diagrammes continus dont on sait que certaines observables admettent des limites. Il arrivera souvent dans la suite qu'on sache par exemple calculer les limites des cumulants libres :
$$\forall k \geq 2,\,\,\exists r_{k}=\lim_{n \to \infty} R_{k}(\omega_{n})\,.$$
Comme les cumulants libres engendrent l'algèbre $\obs$, en effectuant des changements de base, on peut alors virtuellement calculer les limites de toutes les observables $f(\omega_{n})$ ; en particulier, les moments $\widetilde{h}_{k}(\omega_{n})$ des mesures de transition $\mu_{n}$ des diagrammes ont tous des limites $m_{k}$. S'il existe une mesure de probabilité sur $\R$ qui a pour moments les nombres $m_{k}$, et qui est déterminée par ses moments, ceci implique la convergence en loi $\mu_{n} \to \mu$ dans l'espace $\mathscr{P}(\R)$. Alors, par la correspondance de Markov-Krein généralisée, il existe un diagramme continu généralisé $\omega$ tel que $\omega_{n} \to \omega$
uniformément sur tout compact de $\R$. En particulier, si les moments $m_{k}$ correspondent à une mesure à support compact $\mu$, alors cette mesure est automatiquement déterminée par ses moments, et elle correspond à un (vrai) diagramme continu $\omega$. Ainsi :
\begin{proposition}[Topologies faible et forte sur les diagrammes continus]\label{strongweak}
Soit $(\omega_{n})_{n \in \N}$ une suite de diagrammes continus qui converge pour la topologie faible vers un diagramme continu $\omega_{\infty} \in \cym$, c'est-à-dire que pour toute observable de diagrammes $f$ dans $\obs$ (ou toute observable dans une base algébrique), $f(\omega_{n})$ tend vers $f(\omega_{\infty})$. Alors, $\omega_{n}$ converge vers $\omega_{\infty}$ uniformément sur tout $\R$. 
\end{proposition} 
\noindent Autrement dit, si $(\omega_{n})_{n\in \N}$ converge vers $\omega_{\infty}$ pour la topologie faible induite sur $\cym$ par $\obs$, alors la suite converge en fait pour la topologie forte\footnote{Ce point précis et sa preuve ne sont pas vraiment mentionnés dans \cite{IO02} ; c'est essentiellement pour cela que nous avons consacré une section complète à ces subtilités topologiques.} --- notons que l'appartenance $\omega_{\infty} \in \cym$ fait partie des hypothèses de l'énoncé.
\begin{proof}
Comme un diagramme continu $\omega(s)$ est $1$-lipschitzien est égal à $|s|$ pour $s$ assez petit ou $s$ assez grand, on a forcément $\omega(s)\geq |s|$ pour tout $s$. En effet, si $\omega(s)<|s|$ pour un certain $s$, disons $s$ positif, alors $\omega(s)\leq s-\eps$ pour un certain $\eps>0$, et par caractère $1$-lispchitzien, $\omega(t) \leq t-\eps$ pour tout $t $ plus grand que $s$ ; ceci contredit l'hypothèse précédente. Maintenant, on sait déjà que sous les hypothèses de l'énoncé, $\omega_{n}$ converge vers $\omega_{\infty}$ uniformément sur toute partie compacte ; en particulier, c'est vrai sur le support $S=[a,b]$ de $\sigma_{\infty}(s)=(\omega_{\infty}(s)-|s|)/2$. Quitte à étendre ce support, on peut supposer $a \leq 0$ et $b \geq 0$. Pour tout $\eps >0$, au-delà d'un certain rang $N$, 
$$\|\omega_{n}-\omega_{\infty}\|_{\infty,S}\leq \eps\,.$$
Si $s \leq a$, alors $\omega_{n}(s)-\omega_{\infty}(s)=\omega_{n}(s)-|s|\geq 0$ par hypothèse sur les $\omega_{n}$. Mais d'autre part, $(\omega_{n}(s)-|s|)'=\omega_{n}'(s)+1 \geq 0$, donc :
$$\omega_{n}(s)-\omega_{\infty}(s)=\omega_{n}(s)-|s|\leq \omega_{n}(a)-|a| \leq \eps\,.$$
On conclut que sur $S_{-}=]-\infty,a]$, on a également $\|\omega_{n}-\omega_{\infty}\|_{\infty,S_{-}}\leq \eps$. De même, si $s \geq b$, alors $(\omega_{n}(s)-|s|)'=\omega_{n}'(s)-1 \leq 0$, 
donc :
$$0 \leq \omega_{n}(s)-\omega_{\infty}(s)=\omega_{n}(s)-|s|\leq \omega_{n}(b)-|b| \leq \eps\,.$$
et $\|\omega_{n}-\omega_{\infty}\|_{\infty,S_{+}}\leq \eps$ si $S_{+}=[b,+\infty[$. Tout ceci implique bien :
$$\|\omega_{n}-\omega_{\infty}\|_{\infty,\R}\leq \eps$$
pour $n$ assez grand, ce que l'on souhaitait. 
\end{proof}
\bigskip

Ainsi, l'algèbre des observables de diagrammes permet d'établir la convergence pour la topologie uniforme (sur tout $\R$) d'une suite $(\omega_{n})_{n \in \N}$ de diagrammes de Young continus vers un diagramme continu $\omega_{\infty}$. Néanmoins, cette notion de convergence n'est pas la plus forte que l'on puisse espérer. Ainsi, en plus de demander que $\|\omega_{n}- \omega_{\infty}\|_{\infty}$ tende vers $0$, on peut souhaiter que les supports des diagrammes $\omega_{n}$ restent bornés et convergent vers celui de $\omega_{\infty}$, c'est-à-dire que :
\begin{align*}a(\omega_{n})&=\inf\{s \in \R\,\,|\,\, \omega_{n}(s)\neq |s|\} \to a(\omega_{\infty})\,\,;\\
 b(\omega_{n})&=\sup\{s \in \R\,\,|\,\, \omega_{n}(s)\neq |s|\} \to b(\omega_{\infty})\,.
 \end{align*}
On dira dans ce cas que $(\omega_{n})_{n \in \N}$ converge \textbf{ultra-fortement}\footnote{La terminologie est locale et n'a rien à voir avec la notion de topologie ultra-forte d'algèbre d'opérateurs.} vers $\omega_{\infty}$. Dans ce qui suit, nous observerons à trois reprises ce phénomène de convergence ultra-forte, à savoir, pour les diagrammes de Young renormalisés tirés aléatoirement sous les mesures de Plancherel (chapitre \ref{plancherel}), sous les mesures de Schur-Weyl de paramètre $\alpha =1/2$ (chapitre \ref{schurweylmeasure}), et sous les mesures de Gelfand (chapitre \ref{gelfandmeasure}). Dans ces contextes, les techniques d'observables de diagrammes ne donneront rien de plus que la convergence uniforme, et nous aurons recours à d'autres types d'arguments pour établir la convergence des supports.

\chapter{Asymptotique de la mesure de Plancherel sur les partitions}\label{plancherel}

Dans ce chapitre, nous présentons les deux résultats les plus marquants de la théorie asymptotique des représentations du groupe symétrique : la \textbf{loi des grands nombres de Logan-Shepp-Kerov-Vershik} (1977) et le \textbf{théorème central limite de Kerov} (1993). Nous suivrons pour l'essentiel la preuve de \cite{IO02} ; elle repose sur les propriétés combinatoires de l'algèbre d'observables exposées dans le chapitre \ref{tool}. On emprunte également des arguments à \cite{Bia01b,Sni06b}. L'objet principal des seconde et troisième parties du mémoire sera la preuve de résultats analogues à ceux de ce chapitre, mais dans le contexte des mesures de Schur-Weyl, des mesures de Gelfand et des mesures de Plancherel des algèbres d'Hecke. De nombreux arguments seront donc réemployés plus tard, en particulier la théorie due à \'Sniady des cumulants et des cumulants disjoints d'observables de diagrammes (section \ref{sniady}).\bigskip

Dans tout ce qui suit, on utilise la notation $\esper[X]$ pour l'espérance d'une variable aléatoire $X$, et la notation $\proba[A]$ pour la probabilité  d'un événement $A$. Si une mesure $M$ est déjà fixée, on notera aussi $M[X]=\esper[X]$ l'espérance de $X$ sous $M$, et $M[A]$ la probabilité de $A$ sous $M$. \bigskip

\section{Problème d'Ulam et correspondance RSK}\label{ulamrsk}
La mesure de Plancherel est l'objet d'étude principal de ce chapitre, et cette mesure de probabilité admet une définition purement algébrique issue de la théorie des représentations et valable pour tout groupe fini $G$, \emph{cf.} la proposition \ref{defplancherel}. Ceci étant, pour le groupe symétrique $\sym_{n}$, on peut relier la mesure de Plancherel à un problème combinatoire sur les permutations. Si $\sigma$ est une permutation de mot $\sigma(1)\,\sigma(2)\cdots \sigma(n)$, on appelle \textbf{sous-mot croissant} de $\sigma$ un sous-mot
$$\sigma(i_{1})\,\sigma(i_{2})\cdots\sigma(i_{r}) \quad\text{avec}\quad \sigma(i_{1})<\sigma(i_{2})<\cdots < \sigma(i_{r})\,.$$
La plus grande longueur d'un sous-mot croissant de $\sigma$ sera notée $L(\sigma)$ ; cette quantité peut varier entre $1$ et $n$. On suppose maintenant que $\sigma=\sigma_{n}$ est tirée au hasard équiprobablement parmi les $n!$ permutations de $\sym_{n}$, et on s'intéresse à la distribution de la variable aléatoire $L_{n}=L(\sigma_{n})$. Ce problème a été soulevé dans les années 1960 par S. M. Ulam (\emph{cf.} \cite{Ulam61}), et une application ingénieuse du principe des tiroirs montre que l'espérance de $L_{n}$ vérifie
$$\esper[L_{n}]\geq \frac{1}{2}\sqrt{n-1}\,,$$
voir \cite{ES35}. J. M. Hammersley a ensuite montré que $L_{n}/\sqrt{n}$ convergeait en probabilité vers une constante finie $c$ (\emph{cf}. \cite{Hamm72}) ; finalement, il a été montré par Logan et Shepp (\cite{LS77}) et indépendamment par Kerov et Vershik (\cite{KV77}) que $c=2$, voir plus loin le théorème \ref{firstasymptoticplancherel}. L'idée est de se ramener à l'étude de partitions aléatoires en lieu et place des permutations aléatoires ; la \textbf{correspondance de Robinson-Schensted-Knuth} permet ce transport de mesure de probabilité.\bigskip
\bigskip

Compte tenu de l'identité des carrés évoquée à la fin de la section \ref{tableau}, il existe une bijection entre les $n!$ permutations de taille $n$ et les paires de tableaux standards de taille $n$ et de même forme. On doit à Robinson et Schensted (\cite{Rob38,Sch61}) une bijection explicite
$$\mathrm{RSK} : \sigma \in \sym_{n} \mapsto (P(\sigma),Q(\sigma)) \in \bigsqcup_{\lambda \in \ym_{n}} \mathrm{Std}(\lambda) \times \mathrm{Std}(\lambda)$$ 
telle que la plus grande longueur d'un sous-mot croissant de $\sigma$ corresponde à la taille de la première part de la forme commune $\lambda$ des tableaux $P(\sigma)$ et $Q(\sigma)$. L'algorithme RSK\footnote{L'algorithme de Robinson et Schensted a été généralisé dans les années 1970 par D. Knuth, d'où l'acronyme RSK.} est le suivant :
\begin{enumerate}\label{rskalgo}
\item On part des tableaux vides $P=Q=\emptyset$, et on lit le mot de $\sigma$ de gauche à droite.
\item Pour chaque lettre $j=\sigma(i)$, on insert $j$ dans $P$ de la fa\c con suivante. Si $j$ est supérieur à tous les éléments de la première ligne de $P$, on place $j$ à la fin de cette première ligne. Sinon, si $j'$ est le premier entier de la première ligne plus grand que $j$, on remplace $j'$ par $j$ dans cette ligne, et on insère $j'$ dans la seconde ligne. On réitère cette opération jusqu'à ce qu'une case ait été ajoutée au bord de $P$.
\item On ajoute la même case au bord de $Q$, et on la numérote $i$.
\end{enumerate}
Par récurrence sur la longueur $n$ du mot $\sigma$, il est clair que les deux tableaux $P(\sigma)$ et $Q(\sigma)$ sont standards et de même forme.
\figcapt{\psset{unit=1mm}\pspicture(0,5)(120,35)
\rput(0,20){$\sigma=$}\rput(12,20){\foreignlanguage{english}{\young(6,3,17,:5,:4,:29,::8)}}
\rput(50,20){$P(\sigma)=$}
\rput(64,23){\young(6,57,349,128)}
\rput(100,20){$Q(\sigma)=$}
\rput(114,23){\young(7,36,259,148)}
\endpspicture
}{Tableaux standards associés à la permutation $\sigma=631754298$ par l'algorithme RSK.\label{rsk}}{Algorithme RSK}

On peut également retrouver $P(\sigma)$ en appliquant l'algorithme du jeu de taquin au ruban standard obtenu en lisant le mot de $\sigma$, voir la figure \ref{rsk}. D'autre part, $Q(\sigma)$ est toujours égal à $P(\sigma^{-1})$ --- ce point jouera un rôle important dans le chapitre \ref{gelfandmeasure}. On renvoie à \cite[chapitres 1-4]{Ful97} pour une preuve détaillée de ces résultats et de la proposition suivante :
\begin{proposition}[Correspondance de Robinson-Schensted-Knuth, \cite{Rob38,Sch61}]
L'algorithme d'insertion de Schensted peut être inversé ; ainsi, $\mathrm{RSK}$ est une bijection entre $\sym_{n}$ et $\bigsqcup_{\lambda \in \ym_{n}} (\mathrm{Std}(\lambda))^{2}$. De plus, si $\lambda(\sigma)$ désigne la forme commune des tableaux $P(\sigma)$ et $Q(\sigma)$, alors $L(\sigma)=(\lambda(\sigma))_{1}$.
\end{proposition}
\noindent Par conséquent, la loi de $L(\sigma_{n})$ sous la mesure uniforme pour les permutations de $\sym_{n}$ est égale à la loi de $\lambda_{1}$ sous la mesure de probabilité
$$\proba[\lambda]=\frac{(\card \mathrm{Std}(\lambda))^{2}}{n!}$$
pour les partitions de taille $n$. Dans ce qui suit, nous noterons $M_{n}$ cette mesure  sur $\ym_{n}$ ; ainsi, $M_{4}(4)=1/24$, $M_{4}(31)=3/8$, $M_{4}(2^{2})=1/6$, $M_{4}(21^{2})=3/8$ et $M_{4}(1^{4})=1/24$.

\section{Mesures et processus de Plancherel}\label{plancherelprocess}
La mesure $M_{n}$ présentée dans le paragraphe précédent hérite de nombreuses propriétés algébriques, car c'est un cas particulier de \textbf{mesure de Plancherel}. Si $G$ est un groupe fini et si $V$ est un $\C G$-module, on rappelle que $V$ se scinde de manière unique en somme directe de $\C G$-modules simples :
$$V\simeq_{\C G} \bigoplus_{V^\lambda \in \widehat{G}} n_{\lambda}(V)\,V^{\lambda} ,\quad\text{avec } n_{\lambda}(V)\geq 0.$$
\begin{definition}[Mesure de Plancherel]\label{defplancherel}
On appelle mesure de Plancherel d'un $\C G$-module $V$ la mesure de probabilité
$$M_{V}(\lambda)=\frac{n_{\lambda}(V)\,\dim V^{\lambda}}{\dim V}$$
sur l'ensemble $\widehat{G}$ des classes d'isomorphismes de $\C G$-modules simples. En particulier, la mesure de Plancherel $\mathrm{Pl}_{G}$ d'un groupe fini est la mesure de Plancherel de sa représentation régulière :
$$\mathrm{Pl}_{G}(\lambda)=M_{\C G}(\lambda)=\frac{(\dim V^{\lambda})^{2}}{\card G}\,.$$
\end{definition}\bigskip
\bigskip

Dans un contexte d'analyse harmonique non commutative (\emph{cf.} \cite{Var89}), on peut interpréter la mesure de Plancherel comme la duale de la mesure de Haar du groupe ; détaillons brièvement cette interprétation. On rappelle qu'un espace $\leb^{2}$ non commutatif est une algèbre de von Neumann $\mathscr{M} \hookrightarrow \mathscr{L}(H)$ munie d'une trace, c'est-à-dire une application positive homogène $\tau : \mathscr{M}^{+} \to [0,+\infty]$ qui est traciale, fidèle et normale --- on renvoie à \cite{PX03} pour des précisions sur cette terminologie. Si la trace est finie, elle s'étend de manière unique en une forme linéaire $\tau : \mathscr{M} \to \C$, et l'espace $\leb^{2}(\mathscr{M},\tau)$ est la complétion de $\mathscr{M}$ pour la norme préhilbertienne
$$ \|a\|^{2}=\tau(a^{*}a)\,.$$
Un \textbf{espace de probabilité non commutatif} est un espace $\leb^{2}$ non commutatif dont la trace vérifie $\tau(1_{\mathscr{M}})=1$. Lorsque $\mathscr{M}$ est de dimension finie, les subtilités topologiques peuvent être omises, et à la sesquilinéarité du produit scalaire près, la notion coïncide avec celle d'\textbf{algèbre de Frobenius} (normée), \emph{cf.} \cite{Nak39}. Toute algèbre d'endomorphismes $\hendo(V)$ est un espace de probabilité non commutatif pour la trace normalisée
$$\tau(u)=\mathrm{tr}_{V}(u)=\frac{\tr(u)}{\dim V}\,,$$
et d'autre part, l'algèbre d'un groupe fini $\C G$ est un espace de probabilité non commutatif pour la trace normalisée $\mathrm{tr}_{\C G}$ obtenue par restriction de la trace normalisée de $\hendo(\C G)$ à l'image $\C G \hookrightarrow \hendo(\C G)$ de la représentation régulière gauche. Un calcul très simple montre qu'en tant qu'espace de Hilbert, $\C G$ est l'espace $\leb^{2}$ standard de la mesure de Haar du groupe :
$$\forall f = \sum_{g \in G}f(g)\,g,\quad\|f\|^{2}=\mathrm{tr}_{\C G}(f^{*}f)=\int_{G}|f(g)|^{2}\, dm_{\mathrm{Haar}}(g)=\frac{1}{\card G}\sum_{g \in G}|f(g)|^{2}\,.$$
Si $P$ est une mesure de probabilité sur $\widehat{G}$, on peut lui associer une structure d'espace de probabilité non commutatif sur $\bigoplus_{\lambda \in \widehat{G}} \hendo(V^{\lambda})$ en considérant la trace normalisée
$$\tau_{P}=\bigoplus_{\lambda \in \widehat{G}} P(\lambda)\,\,\mathrm{tr}_{V^{\lambda}}\,.$$
Alors, la transformée de Fourier abstraite $\C G \to \bigoplus_{\lambda \in\widehat{G}} \hendo(V^{\lambda})$ est une isométrie d'espaces de probabilité non commutatifs si et seulement si l'on prend $P=\mathrm{Pl}_{G}$. Ainsi, la mesure de Plancherel apparaît bien naturellement comme duale de la mesure de Haar du groupe. Nous donnerons dans le chapitre \ref{iwahori} une interprétation semblable de la $q$-mesure de Plancherel ; d'autre part, la notion d'espace de probabilité non commutatif sera utile pour comprendre les raisonnements du paragraphe \ref{sniady}.
\bigskip
\bigskip

Pour le groupe symétrique $\sym_{n}$, $\dim V^{\lambda}=\card \mathrm{Std}(\lambda)$, donc la mesure $M_{n}$ n'est nulle autre que la mesure de Plancherel $\mathrm{Pl}_{\sym_{n}}$. De plus, il existe un processus aléatoire $(\lambda^{(n)})_{n \in \N}$ à valeurs dans $\ym$ et tel que la loi marginale de $\lambda_{n}$ soit exactement $M_{n}$. En effet, considérons le processus markovien de probabilités de transition données par 
$$p(\lambda,\Lambda)= \begin{cases} 
\frac{1}{|\Lambda|}\,\frac{\dim \Lambda}{\dim \lambda} &\text{si }\lambda \nearrow \Lambda,\\
0 &\text{sinon}.\end{cases}$$
C'est le \textbf{processus de Plancherel} des groupes symétriques. Compte tenu du théorème de branchement \ref{branching}, si $n$ désigne la taille de $|\lambda|$, alors la probabilité de transition $p(\lambda,\cdot)$ est la mesure de Plancherel de $\mathrm{Ind}^{\sym_{n+1}}_{\sym_{n}} (V^{\lambda})$. Comme $$\mathrm{Ind}_{\sym_{n}}^{\sym_{n+1}}(\C\sym_{n})=\C\sym_{n+1}\otimes_{\C\sym_{n}} \C\sym_{n} = \C\sym_{n+1}\,,$$
on en déduit par récurrence sur $n$ que la loi de $\lambda^{(n)}$ est $M_{n}$. De plus, les probabilités de transition se calculent aisément à l'aide de la formule des équerres :
$$p(\lambda,\Lambda)=\frac{\prod_{(i,j) \in\lambda} h(i,j)}{\prod_{(i,j) \in \Lambda} h(i,j)}\,.$$
 On peut donc construire très facilement de grands diagrammes de Young aléatoires tirés suivant la mesure de Plancherel, voir la figure \ref{bigplancherel} --- le diagramme a été construit à l'aide du logiciel libre \texttt{sage}, \emph{cf.} \cite{sage}. Alternativement, on peut réaliser le processus de Plancherel en projetant par la correspondance RSK le \textbf{processus de Knuth}. Ce processus markovien produit une suite de permutations $(\sigma^{(n)})_{n \in \N}$ en insérant aléatoirement et de fa\c con équiprobable la lettre $n+1$ dans le mot de $\sigma^{(n)}$ pour obtenir le mot de $\sigma^{(n+1)}$. Il n'est pas difficile de montrer que si $(\sigma^{(n)})_{n\in \N}$ est un processus de Knuth, alors $(\lambda(\sigma^{(n)}))_{n \in \N}$ est un processus de Plancherel --- il suffit d'utiliser le théorème de branchement \ref{branching} et l'identité $\dim \lambda=\card \mathrm{Std}(\lambda)$.\bigskip
\bigskip

\`A partir de la formule des équerres pour les probabilités de transition du processus de Plancherel, on peut retrouver la mesure de transition des diagrammes évoquée dans la section \ref{freecumulant}. Ainsi, si $\Lambda$ est obtenu à partir de $\lambda$ en rajoutant une case au niveau du coin d'abscisse $x_{k}$, alors $p(\lambda,\Lambda)=\mu_{\lambda}(x_{k})$. Les mesures de probabilité $p(\lambda,\cdot)$ et $\mu_{\lambda}$ peuvent donc être identifiées. Partant de cette observation, Kerov a proposé un analogue différentiel déterministe du processus de Plancherel, voir par exemple \cite{Ker99}. Rappelons que si $\omega$ est un diagramme continu, sa mesure de transition $\mu_{\omega}$ est la mesure de probabilité définie par la relation
$$G_{\omega}(z)=\frac{1}{z}\,\exp\left(-\int_{\R} \frac{\sigma_{\omega}'(s)}{z-s}\,ds\right)= \int_{\R}\frac{\mu_{\omega}(ds)}{z-s}\,.$$
Un \textbf{processus de Plancherel différentiel} est la donnée d'une famille de diagrammes continus $(\omega(\cdot,t))_{t \in \R_{+}}$ qui croissent suivant leurs mesures de transition :
$$\forall s,t,\,\,\,\frac{\partial \sigma(s,t)}{\partial t}\,ds=\mu_{\omega(\cdot,t)}(ds)\,.$$
Alors, les fonctions génératrices $G_{\omega(\cdot,t)}=G(\cdot,t)$ satisfont l'\textbf{équation hydrodynamique de Burgers}
$$\frac{\partial G(z,t)}{\partial t}+G(z,t)\,\frac{\partial G(z,t)}{\partial z}=0\,,$$
et on peut montrer que les processus de Plancherel différentiels ont la même asymptotique que les processus de Plancherel <<~algébriques~>>. Le lien entre cette approche analytique et la définition algébrique des processus de Plancherel est expliqué dans les articles originaux \cite{LS77,KV77} ; l'asymptotique des processus de Plancherel y est également traduite en un problème variationnel pour une fonctionnelle sur les diagrammes continus, voir \cite{LS77}. Pour notre part, nous mènerons l'analyse asymptotique des mesures $M_{n}$ en restant dans le cadre algébrique, et en utilisant les observables introduites dans le chapitre \ref{tool}.
\figcapt{\psset{unit=1mm} \pspicture(-70,0)(70,70) \psline{<->}(-70,70)(0,0)(70,70) \parametricplot{-2}{2}{t 1.46341 mul 22.3607 mul 
t t 0.5 mul arcsin mul 0.0174533 mul t t mul neg 4 add sqrt add 0.63661977 mul 1.46341 mul 22.3607 mul }\psline[linewidth=0.25pt](0.000000,0.000000)(1.46341,1.46341)(0.000000,2.92683)(-1.46341,1.46341)(0.000000,0.000000) \psline[linewidth=0.25pt](1.46341,1.46341)(2.92683,2.92683)(1.46341,4.39024)(0.000000,2.92683)(1.46341,1.46341) \psline[linewidth=0.25pt](2.92683,2.92683)(4.39024,4.39024)(2.92683,5.85366)(1.46341,4.39024)(2.92683,2.92683) \psline[linewidth=0.25pt](4.39024,4.39024)(5.85366,5.85366)(4.39024,7.31707)(2.92683,5.85366)(4.39024,4.39024) \psline[linewidth=0.25pt](5.85366,5.85366)(7.31707,7.31707)(5.85366,8.78049)(4.39024,7.31707)(5.85366,5.85366) \psline[linewidth=0.25pt](7.31707,7.31707)(8.78049,8.78049)(7.31707,10.2439)(5.85366,8.78049)(7.31707,7.31707) \psline[linewidth=0.25pt](8.78049,8.78049)(10.2439,10.2439)(8.78049,11.7073)(7.31707,10.2439)(8.78049,8.78049) \psline[linewidth=0.25pt](10.2439,10.2439)(11.7073,11.7073)(10.2439,13.1707)(8.78049,11.7073)(10.2439,10.2439) \psline[linewidth=0.25pt](11.7073,11.7073)(13.1707,13.1707)(11.7073,14.6341)(10.2439,13.1707)(11.7073,11.7073) \psline[linewidth=0.25pt](13.1707,13.1707)(14.6341,14.6341)(13.1707,16.0976)(11.7073,14.6341)(13.1707,13.1707) \psline[linewidth=0.25pt](14.6341,14.6341)(16.0976,16.0976)(14.6341,17.5610)(13.1707,16.0976)(14.6341,14.6341) \psline[linewidth=0.25pt](16.0976,16.0976)(17.5610,17.5610)(16.0976,19.0244)(14.6341,17.5610)(16.0976,16.0976) \psline[linewidth=0.25pt](17.5610,17.5610)(19.0244,19.0244)(17.5610,20.4878)(16.0976,19.0244)(17.5610,17.5610) \psline[linewidth=0.25pt](19.0244,19.0244)(20.4878,20.4878)(19.0244,21.9512)(17.5610,20.4878)(19.0244,19.0244) \psline[linewidth=0.25pt](20.4878,20.4878)(21.9512,21.9512)(20.4878,23.4146)(19.0244,21.9512)(20.4878,20.4878) \psline[linewidth=0.25pt](21.9512,21.9512)(23.4146,23.4146)(21.9512,24.8780)(20.4878,23.4146)(21.9512,21.9512) \psline[linewidth=0.25pt](23.4146,23.4146)(24.8780,24.8780)(23.4146,26.3415)(21.9512,24.8780)(23.4146,23.4146) \psline[linewidth=0.25pt](24.8780,24.8780)(26.3415,26.3415)(24.8780,27.8049)(23.4146,26.3415)(24.8780,24.8780) \psline[linewidth=0.25pt](26.3415,26.3415)(27.8049,27.8049)(26.3415,29.2683)(24.8780,27.8049)(26.3415,26.3415) \psline[linewidth=0.25pt](27.8049,27.8049)(29.2683,29.2683)(27.8049,30.7317)(26.3415,29.2683)(27.8049,27.8049) \psline[linewidth=0.25pt](29.2683,29.2683)(30.7317,30.7317)(29.2683,32.1951)(27.8049,30.7317)(29.2683,29.2683) \psline[linewidth=0.25pt](30.7317,30.7317)(32.1951,32.1951)(30.7317,33.6585)(29.2683,32.1951)(30.7317,30.7317) \psline[linewidth=0.25pt](32.1951,32.1951)(33.6585,33.6585)(32.1951,35.1220)(30.7317,33.6585)(32.1951,32.1951) \psline[linewidth=0.25pt](33.6585,33.6585)(35.1220,35.1220)(33.6585,36.5854)(32.1951,35.1220)(33.6585,33.6585) \psline[linewidth=0.25pt](35.1220,35.1220)(36.5854,36.5854)(35.1220,38.0488)(33.6585,36.5854)(35.1220,35.1220) \psline[linewidth=0.25pt](36.5854,36.5854)(38.0488,38.0488)(36.5854,39.5122)(35.1220,38.0488)(36.5854,36.5854) \psline[linewidth=0.25pt](38.0488,38.0488)(39.5122,39.5122)(38.0488,40.9756)(36.5854,39.5122)(38.0488,38.0488) \psline[linewidth=0.25pt](39.5122,39.5122)(40.9756,40.9756)(39.5122,42.4390)(38.0488,40.9756)(39.5122,39.5122) \psline[linewidth=0.25pt](40.9756,40.9756)(42.4390,42.4390)(40.9756,43.9024)(39.5122,42.4390)(40.9756,40.9756) \psline[linewidth=0.25pt](42.4390,42.4390)(43.9024,43.9024)(42.4390,45.3659)(40.9756,43.9024)(42.4390,42.4390) \psline[linewidth=0.25pt](43.9024,43.9024)(45.3659,45.3659)(43.9024,46.8293)(42.4390,45.3659)(43.9024,43.9024) \psline[linewidth=0.25pt](45.3659,45.3659)(46.8293,46.8293)(45.3659,48.2927)(43.9024,46.8293)(45.3659,45.3659) \psline[linewidth=0.25pt](46.8293,46.8293)(48.2927,48.2927)(46.8293,49.7561)(45.3659,48.2927)(46.8293,46.8293) \psline[linewidth=0.25pt](48.2927,48.2927)(49.7561,49.7561)(48.2927,51.2195)(46.8293,49.7561)(48.2927,48.2927) \psline[linewidth=0.25pt](49.7561,49.7561)(51.2195,51.2195)(49.7561,52.6829)(48.2927,51.2195)(49.7561,49.7561) \psline[linewidth=0.25pt](51.2195,51.2195)(52.6829,52.6829)(51.2195,54.1463)(49.7561,52.6829)(51.2195,51.2195) \psline[linewidth=0.25pt](52.6829,52.6829)(54.1463,54.1463)(52.6829,55.6098)(51.2195,54.1463)(52.6829,52.6829) \psline[linewidth=0.25pt](54.1463,54.1463)(55.6098,55.6098)(54.1463,57.0732)(52.6829,55.6098)(54.1463,54.1463) \psline[linewidth=0.25pt](55.6098,55.6098)(57.0732,57.0732)(55.6098,58.5366)(54.1463,57.0732)(55.6098,55.6098) \psline[linewidth=0.25pt](57.0732,57.0732)(58.5366,58.5366)(57.0732,60.0000)(55.6098,58.5366)(57.0732,57.0732) \psline[linewidth=0.25pt](58.5366,58.5366)(60.0000,60.0000)(58.5366,61.4634)(57.0732,60.0000)(58.5366,58.5366) \psline[linewidth=0.25pt](-1.46341,1.46341)(0.000000,2.92683)(-1.46341,4.39024)(-2.92683,2.92683)(-1.46341,1.46341) \psline[linewidth=0.25pt](0.000000,2.92683)(1.46341,4.39024)(0.000000,5.85366)(-1.46341,4.39024)(0.000000,2.92683) \psline[linewidth=0.25pt](1.46341,4.39024)(2.92683,5.85366)(1.46341,7.31707)(0.000000,5.85366)(1.46341,4.39024) \psline[linewidth=0.25pt](2.92683,5.85366)(4.39024,7.31707)(2.92683,8.78049)(1.46341,7.31707)(2.92683,5.85366) \psline[linewidth=0.25pt](4.39024,7.31707)(5.85366,8.78049)(4.39024,10.2439)(2.92683,8.78049)(4.39024,7.31707) \psline[linewidth=0.25pt](5.85366,8.78049)(7.31707,10.2439)(5.85366,11.7073)(4.39024,10.2439)(5.85366,8.78049) \psline[linewidth=0.25pt](7.31707,10.2439)(8.78049,11.7073)(7.31707,13.1707)(5.85366,11.7073)(7.31707,10.2439) \psline[linewidth=0.25pt](8.78049,11.7073)(10.2439,13.1707)(8.78049,14.6341)(7.31707,13.1707)(8.78049,11.7073) \psline[linewidth=0.25pt](10.2439,13.1707)(11.7073,14.6341)(10.2439,16.0976)(8.78049,14.6341)(10.2439,13.1707) \psline[linewidth=0.25pt](11.7073,14.6341)(13.1707,16.0976)(11.7073,17.5610)(10.2439,16.0976)(11.7073,14.6341) \psline[linewidth=0.25pt](13.1707,16.0976)(14.6341,17.5610)(13.1707,19.0244)(11.7073,17.5610)(13.1707,16.0976) \psline[linewidth=0.25pt](14.6341,17.5610)(16.0976,19.0244)(14.6341,20.4878)(13.1707,19.0244)(14.6341,17.5610) \psline[linewidth=0.25pt](16.0976,19.0244)(17.5610,20.4878)(16.0976,21.9512)(14.6341,20.4878)(16.0976,19.0244) \psline[linewidth=0.25pt](17.5610,20.4878)(19.0244,21.9512)(17.5610,23.4146)(16.0976,21.9512)(17.5610,20.4878) \psline[linewidth=0.25pt](19.0244,21.9512)(20.4878,23.4146)(19.0244,24.8780)(17.5610,23.4146)(19.0244,21.9512) \psline[linewidth=0.25pt](20.4878,23.4146)(21.9512,24.8780)(20.4878,26.3415)(19.0244,24.8780)(20.4878,23.4146) \psline[linewidth=0.25pt](21.9512,24.8780)(23.4146,26.3415)(21.9512,27.8049)(20.4878,26.3415)(21.9512,24.8780) \psline[linewidth=0.25pt](23.4146,26.3415)(24.8780,27.8049)(23.4146,29.2683)(21.9512,27.8049)(23.4146,26.3415) \psline[linewidth=0.25pt](24.8780,27.8049)(26.3415,29.2683)(24.8780,30.7317)(23.4146,29.2683)(24.8780,27.8049) \psline[linewidth=0.25pt](26.3415,29.2683)(27.8049,30.7317)(26.3415,32.1951)(24.8780,30.7317)(26.3415,29.2683) \psline[linewidth=0.25pt](27.8049,30.7317)(29.2683,32.1951)(27.8049,33.6585)(26.3415,32.1951)(27.8049,30.7317) \psline[linewidth=0.25pt](29.2683,32.1951)(30.7317,33.6585)(29.2683,35.1220)(27.8049,33.6585)(29.2683,32.1951) \psline[linewidth=0.25pt](30.7317,33.6585)(32.1951,35.1220)(30.7317,36.5854)(29.2683,35.1220)(30.7317,33.6585) \psline[linewidth=0.25pt](32.1951,35.1220)(33.6585,36.5854)(32.1951,38.0488)(30.7317,36.5854)(32.1951,35.1220) \psline[linewidth=0.25pt](33.6585,36.5854)(35.1220,38.0488)(33.6585,39.5122)(32.1951,38.0488)(33.6585,36.5854) \psline[linewidth=0.25pt](35.1220,38.0488)(36.5854,39.5122)(35.1220,40.9756)(33.6585,39.5122)(35.1220,38.0488) \psline[linewidth=0.25pt](36.5854,39.5122)(38.0488,40.9756)(36.5854,42.4390)(35.1220,40.9756)(36.5854,39.5122) \psline[linewidth=0.25pt](38.0488,40.9756)(39.5122,42.4390)(38.0488,43.9024)(36.5854,42.4390)(38.0488,40.9756) \psline[linewidth=0.25pt](39.5122,42.4390)(40.9756,43.9024)(39.5122,45.3659)(38.0488,43.9024)(39.5122,42.4390) \psline[linewidth=0.25pt](40.9756,43.9024)(42.4390,45.3659)(40.9756,46.8293)(39.5122,45.3659)(40.9756,43.9024) \psline[linewidth=0.25pt](42.4390,45.3659)(43.9024,46.8293)(42.4390,48.2927)(40.9756,46.8293)(42.4390,45.3659) \psline[linewidth=0.25pt](43.9024,46.8293)(45.3659,48.2927)(43.9024,49.7561)(42.4390,48.2927)(43.9024,46.8293) \psline[linewidth=0.25pt](45.3659,48.2927)(46.8293,49.7561)(45.3659,51.2195)(43.9024,49.7561)(45.3659,48.2927) \psline[linewidth=0.25pt](46.8293,49.7561)(48.2927,51.2195)(46.8293,52.6829)(45.3659,51.2195)(46.8293,49.7561) \psline[linewidth=0.25pt](48.2927,51.2195)(49.7561,52.6829)(48.2927,54.1463)(46.8293,52.6829)(48.2927,51.2195) \psline[linewidth=0.25pt](49.7561,52.6829)(51.2195,54.1463)(49.7561,55.6098)(48.2927,54.1463)(49.7561,52.6829) \psline[linewidth=0.25pt](51.2195,54.1463)(52.6829,55.6098)(51.2195,57.0732)(49.7561,55.6098)(51.2195,54.1463) \psline[linewidth=0.25pt](52.6829,55.6098)(54.1463,57.0732)(52.6829,58.5366)(51.2195,57.0732)(52.6829,55.6098) \psline[linewidth=0.25pt](-2.92683,2.92683)(-1.46341,4.39024)(-2.92683,5.85366)(-4.39024,4.39024)(-2.92683,2.92683) \psline[linewidth=0.25pt](-1.46341,4.39024)(0.000000,5.85366)(-1.46341,7.31707)(-2.92683,5.85366)(-1.46341,4.39024) \psline[linewidth=0.25pt](0.000000,5.85366)(1.46341,7.31707)(0.000000,8.78049)(-1.46341,7.31707)(0.000000,5.85366) \psline[linewidth=0.25pt](1.46341,7.31707)(2.92683,8.78049)(1.46341,10.2439)(0.000000,8.78049)(1.46341,7.31707) \psline[linewidth=0.25pt](2.92683,8.78049)(4.39024,10.2439)(2.92683,11.7073)(1.46341,10.2439)(2.92683,8.78049) \psline[linewidth=0.25pt](4.39024,10.2439)(5.85366,11.7073)(4.39024,13.1707)(2.92683,11.7073)(4.39024,10.2439) \psline[linewidth=0.25pt](5.85366,11.7073)(7.31707,13.1707)(5.85366,14.6341)(4.39024,13.1707)(5.85366,11.7073) \psline[linewidth=0.25pt](7.31707,13.1707)(8.78049,14.6341)(7.31707,16.0976)(5.85366,14.6341)(7.31707,13.1707) \psline[linewidth=0.25pt](8.78049,14.6341)(10.2439,16.0976)(8.78049,17.5610)(7.31707,16.0976)(8.78049,14.6341) \psline[linewidth=0.25pt](10.2439,16.0976)(11.7073,17.5610)(10.2439,19.0244)(8.78049,17.5610)(10.2439,16.0976) \psline[linewidth=0.25pt](11.7073,17.5610)(13.1707,19.0244)(11.7073,20.4878)(10.2439,19.0244)(11.7073,17.5610) \psline[linewidth=0.25pt](13.1707,19.0244)(14.6341,20.4878)(13.1707,21.9512)(11.7073,20.4878)(13.1707,19.0244) \psline[linewidth=0.25pt](14.6341,20.4878)(16.0976,21.9512)(14.6341,23.4146)(13.1707,21.9512)(14.6341,20.4878) \psline[linewidth=0.25pt](16.0976,21.9512)(17.5610,23.4146)(16.0976,24.8780)(14.6341,23.4146)(16.0976,21.9512) \psline[linewidth=0.25pt](17.5610,23.4146)(19.0244,24.8780)(17.5610,26.3415)(16.0976,24.8780)(17.5610,23.4146) \psline[linewidth=0.25pt](19.0244,24.8780)(20.4878,26.3415)(19.0244,27.8049)(17.5610,26.3415)(19.0244,24.8780) \psline[linewidth=0.25pt](20.4878,26.3415)(21.9512,27.8049)(20.4878,29.2683)(19.0244,27.8049)(20.4878,26.3415) \psline[linewidth=0.25pt](21.9512,27.8049)(23.4146,29.2683)(21.9512,30.7317)(20.4878,29.2683)(21.9512,27.8049) \psline[linewidth=0.25pt](23.4146,29.2683)(24.8780,30.7317)(23.4146,32.1951)(21.9512,30.7317)(23.4146,29.2683) \psline[linewidth=0.25pt](24.8780,30.7317)(26.3415,32.1951)(24.8780,33.6585)(23.4146,32.1951)(24.8780,30.7317) \psline[linewidth=0.25pt](26.3415,32.1951)(27.8049,33.6585)(26.3415,35.1220)(24.8780,33.6585)(26.3415,32.1951) \psline[linewidth=0.25pt](27.8049,33.6585)(29.2683,35.1220)(27.8049,36.5854)(26.3415,35.1220)(27.8049,33.6585) \psline[linewidth=0.25pt](29.2683,35.1220)(30.7317,36.5854)(29.2683,38.0488)(27.8049,36.5854)(29.2683,35.1220) \psline[linewidth=0.25pt](30.7317,36.5854)(32.1951,38.0488)(30.7317,39.5122)(29.2683,38.0488)(30.7317,36.5854) \psline[linewidth=0.25pt](32.1951,38.0488)(33.6585,39.5122)(32.1951,40.9756)(30.7317,39.5122)(32.1951,38.0488) \psline[linewidth=0.25pt](33.6585,39.5122)(35.1220,40.9756)(33.6585,42.4390)(32.1951,40.9756)(33.6585,39.5122) \psline[linewidth=0.25pt](35.1220,40.9756)(36.5854,42.4390)(35.1220,43.9024)(33.6585,42.4390)(35.1220,40.9756) \psline[linewidth=0.25pt](36.5854,42.4390)(38.0488,43.9024)(36.5854,45.3659)(35.1220,43.9024)(36.5854,42.4390) \psline[linewidth=0.25pt](38.0488,43.9024)(39.5122,45.3659)(38.0488,46.8293)(36.5854,45.3659)(38.0488,43.9024) \psline[linewidth=0.25pt](39.5122,45.3659)(40.9756,46.8293)(39.5122,48.2927)(38.0488,46.8293)(39.5122,45.3659) \psline[linewidth=0.25pt](40.9756,46.8293)(42.4390,48.2927)(40.9756,49.7561)(39.5122,48.2927)(40.9756,46.8293) \psline[linewidth=0.25pt](-4.39024,4.39024)(-2.92683,5.85366)(-4.39024,7.31707)(-5.85366,5.85366)(-4.39024,4.39024) \psline[linewidth=0.25pt](-2.92683,5.85366)(-1.46341,7.31707)(-2.92683,8.78049)(-4.39024,7.31707)(-2.92683,5.85366) \psline[linewidth=0.25pt](-1.46341,7.31707)(0.000000,8.78049)(-1.46341,10.2439)(-2.92683,8.78049)(-1.46341,7.31707) \psline[linewidth=0.25pt](0.000000,8.78049)(1.46341,10.2439)(0.000000,11.7073)(-1.46341,10.2439)(0.000000,8.78049) \psline[linewidth=0.25pt](1.46341,10.2439)(2.92683,11.7073)(1.46341,13.1707)(0.000000,11.7073)(1.46341,10.2439) \psline[linewidth=0.25pt](2.92683,11.7073)(4.39024,13.1707)(2.92683,14.6341)(1.46341,13.1707)(2.92683,11.7073) \psline[linewidth=0.25pt](4.39024,13.1707)(5.85366,14.6341)(4.39024,16.0976)(2.92683,14.6341)(4.39024,13.1707) \psline[linewidth=0.25pt](5.85366,14.6341)(7.31707,16.0976)(5.85366,17.5610)(4.39024,16.0976)(5.85366,14.6341) \psline[linewidth=0.25pt](7.31707,16.0976)(8.78049,17.5610)(7.31707,19.0244)(5.85366,17.5610)(7.31707,16.0976) \psline[linewidth=0.25pt](8.78049,17.5610)(10.2439,19.0244)(8.78049,20.4878)(7.31707,19.0244)(8.78049,17.5610) \psline[linewidth=0.25pt](10.2439,19.0244)(11.7073,20.4878)(10.2439,21.9512)(8.78049,20.4878)(10.2439,19.0244) \psline[linewidth=0.25pt](11.7073,20.4878)(13.1707,21.9512)(11.7073,23.4146)(10.2439,21.9512)(11.7073,20.4878) \psline[linewidth=0.25pt](13.1707,21.9512)(14.6341,23.4146)(13.1707,24.8780)(11.7073,23.4146)(13.1707,21.9512) \psline[linewidth=0.25pt](14.6341,23.4146)(16.0976,24.8780)(14.6341,26.3415)(13.1707,24.8780)(14.6341,23.4146) \psline[linewidth=0.25pt](16.0976,24.8780)(17.5610,26.3415)(16.0976,27.8049)(14.6341,26.3415)(16.0976,24.8780) \psline[linewidth=0.25pt](17.5610,26.3415)(19.0244,27.8049)(17.5610,29.2683)(16.0976,27.8049)(17.5610,26.3415) \psline[linewidth=0.25pt](19.0244,27.8049)(20.4878,29.2683)(19.0244,30.7317)(17.5610,29.2683)(19.0244,27.8049) \psline[linewidth=0.25pt](20.4878,29.2683)(21.9512,30.7317)(20.4878,32.1951)(19.0244,30.7317)(20.4878,29.2683) \psline[linewidth=0.25pt](21.9512,30.7317)(23.4146,32.1951)(21.9512,33.6585)(20.4878,32.1951)(21.9512,30.7317) \psline[linewidth=0.25pt](23.4146,32.1951)(24.8780,33.6585)(23.4146,35.1220)(21.9512,33.6585)(23.4146,32.1951) \psline[linewidth=0.25pt](24.8780,33.6585)(26.3415,35.1220)(24.8780,36.5854)(23.4146,35.1220)(24.8780,33.6585) \psline[linewidth=0.25pt](26.3415,35.1220)(27.8049,36.5854)(26.3415,38.0488)(24.8780,36.5854)(26.3415,35.1220) \psline[linewidth=0.25pt](27.8049,36.5854)(29.2683,38.0488)(27.8049,39.5122)(26.3415,38.0488)(27.8049,36.5854) \psline[linewidth=0.25pt](29.2683,38.0488)(30.7317,39.5122)(29.2683,40.9756)(27.8049,39.5122)(29.2683,38.0488) \psline[linewidth=0.25pt](30.7317,39.5122)(32.1951,40.9756)(30.7317,42.4390)(29.2683,40.9756)(30.7317,39.5122) \psline[linewidth=0.25pt](32.1951,40.9756)(33.6585,42.4390)(32.1951,43.9024)(30.7317,42.4390)(32.1951,40.9756) \psline[linewidth=0.25pt](33.6585,42.4390)(35.1220,43.9024)(33.6585,45.3659)(32.1951,43.9024)(33.6585,42.4390) \psline[linewidth=0.25pt](35.1220,43.9024)(36.5854,45.3659)(35.1220,46.8293)(33.6585,45.3659)(35.1220,43.9024) \psline[linewidth=0.25pt](36.5854,45.3659)(38.0488,46.8293)(36.5854,48.2927)(35.1220,46.8293)(36.5854,45.3659) \psline[linewidth=0.25pt](38.0488,46.8293)(39.5122,48.2927)(38.0488,49.7561)(36.5854,48.2927)(38.0488,46.8293) \psline[linewidth=0.25pt](-5.85366,5.85366)(-4.39024,7.31707)(-5.85366,8.78049)(-7.31707,7.31707)(-5.85366,5.85366) \psline[linewidth=0.25pt](-4.39024,7.31707)(-2.92683,8.78049)(-4.39024,10.2439)(-5.85366,8.78049)(-4.39024,7.31707) \psline[linewidth=0.25pt](-2.92683,8.78049)(-1.46341,10.2439)(-2.92683,11.7073)(-4.39024,10.2439)(-2.92683,8.78049) \psline[linewidth=0.25pt](-1.46341,10.2439)(0.000000,11.7073)(-1.46341,13.1707)(-2.92683,11.7073)(-1.46341,10.2439) \psline[linewidth=0.25pt](0.000000,11.7073)(1.46341,13.1707)(0.000000,14.6341)(-1.46341,13.1707)(0.000000,11.7073) \psline[linewidth=0.25pt](1.46341,13.1707)(2.92683,14.6341)(1.46341,16.0976)(0.000000,14.6341)(1.46341,13.1707) \psline[linewidth=0.25pt](2.92683,14.6341)(4.39024,16.0976)(2.92683,17.5610)(1.46341,16.0976)(2.92683,14.6341) \psline[linewidth=0.25pt](4.39024,16.0976)(5.85366,17.5610)(4.39024,19.0244)(2.92683,17.5610)(4.39024,16.0976) \psline[linewidth=0.25pt](5.85366,17.5610)(7.31707,19.0244)(5.85366,20.4878)(4.39024,19.0244)(5.85366,17.5610) \psline[linewidth=0.25pt](7.31707,19.0244)(8.78049,20.4878)(7.31707,21.9512)(5.85366,20.4878)(7.31707,19.0244) \psline[linewidth=0.25pt](8.78049,20.4878)(10.2439,21.9512)(8.78049,23.4146)(7.31707,21.9512)(8.78049,20.4878) \psline[linewidth=0.25pt](10.2439,21.9512)(11.7073,23.4146)(10.2439,24.8780)(8.78049,23.4146)(10.2439,21.9512) \psline[linewidth=0.25pt](11.7073,23.4146)(13.1707,24.8780)(11.7073,26.3415)(10.2439,24.8780)(11.7073,23.4146) \psline[linewidth=0.25pt](13.1707,24.8780)(14.6341,26.3415)(13.1707,27.8049)(11.7073,26.3415)(13.1707,24.8780) \psline[linewidth=0.25pt](14.6341,26.3415)(16.0976,27.8049)(14.6341,29.2683)(13.1707,27.8049)(14.6341,26.3415) \psline[linewidth=0.25pt](16.0976,27.8049)(17.5610,29.2683)(16.0976,30.7317)(14.6341,29.2683)(16.0976,27.8049) \psline[linewidth=0.25pt](17.5610,29.2683)(19.0244,30.7317)(17.5610,32.1951)(16.0976,30.7317)(17.5610,29.2683) \psline[linewidth=0.25pt](19.0244,30.7317)(20.4878,32.1951)(19.0244,33.6585)(17.5610,32.1951)(19.0244,30.7317) \psline[linewidth=0.25pt](20.4878,32.1951)(21.9512,33.6585)(20.4878,35.1220)(19.0244,33.6585)(20.4878,32.1951) \psline[linewidth=0.25pt](21.9512,33.6585)(23.4146,35.1220)(21.9512,36.5854)(20.4878,35.1220)(21.9512,33.6585) \psline[linewidth=0.25pt](23.4146,35.1220)(24.8780,36.5854)(23.4146,38.0488)(21.9512,36.5854)(23.4146,35.1220) \psline[linewidth=0.25pt](24.8780,36.5854)(26.3415,38.0488)(24.8780,39.5122)(23.4146,38.0488)(24.8780,36.5854) \psline[linewidth=0.25pt](26.3415,38.0488)(27.8049,39.5122)(26.3415,40.9756)(24.8780,39.5122)(26.3415,38.0488) \psline[linewidth=0.25pt](27.8049,39.5122)(29.2683,40.9756)(27.8049,42.4390)(26.3415,40.9756)(27.8049,39.5122) \psline[linewidth=0.25pt](29.2683,40.9756)(30.7317,42.4390)(29.2683,43.9024)(27.8049,42.4390)(29.2683,40.9756) \psline[linewidth=0.25pt](30.7317,42.4390)(32.1951,43.9024)(30.7317,45.3659)(29.2683,43.9024)(30.7317,42.4390) \psline[linewidth=0.25pt](32.1951,43.9024)(33.6585,45.3659)(32.1951,46.8293)(30.7317,45.3659)(32.1951,43.9024) \psline[linewidth=0.25pt](-7.31707,7.31707)(-5.85366,8.78049)(-7.31707,10.2439)(-8.78049,8.78049)(-7.31707,7.31707) \psline[linewidth=0.25pt](-5.85366,8.78049)(-4.39024,10.2439)(-5.85366,11.7073)(-7.31707,10.2439)(-5.85366,8.78049) \psline[linewidth=0.25pt](-4.39024,10.2439)(-2.92683,11.7073)(-4.39024,13.1707)(-5.85366,11.7073)(-4.39024,10.2439) \psline[linewidth=0.25pt](-2.92683,11.7073)(-1.46341,13.1707)(-2.92683,14.6341)(-4.39024,13.1707)(-2.92683,11.7073) \psline[linewidth=0.25pt](-1.46341,13.1707)(0.000000,14.6341)(-1.46341,16.0976)(-2.92683,14.6341)(-1.46341,13.1707) \psline[linewidth=0.25pt](0.000000,14.6341)(1.46341,16.0976)(0.000000,17.5610)(-1.46341,16.0976)(0.000000,14.6341) \psline[linewidth=0.25pt](1.46341,16.0976)(2.92683,17.5610)(1.46341,19.0244)(0.000000,17.5610)(1.46341,16.0976) \psline[linewidth=0.25pt](2.92683,17.5610)(4.39024,19.0244)(2.92683,20.4878)(1.46341,19.0244)(2.92683,17.5610) \psline[linewidth=0.25pt](4.39024,19.0244)(5.85366,20.4878)(4.39024,21.9512)(2.92683,20.4878)(4.39024,19.0244) \psline[linewidth=0.25pt](5.85366,20.4878)(7.31707,21.9512)(5.85366,23.4146)(4.39024,21.9512)(5.85366,20.4878) \psline[linewidth=0.25pt](7.31707,21.9512)(8.78049,23.4146)(7.31707,24.8780)(5.85366,23.4146)(7.31707,21.9512) \psline[linewidth=0.25pt](8.78049,23.4146)(10.2439,24.8780)(8.78049,26.3415)(7.31707,24.8780)(8.78049,23.4146) \psline[linewidth=0.25pt](10.2439,24.8780)(11.7073,26.3415)(10.2439,27.8049)(8.78049,26.3415)(10.2439,24.8780) \psline[linewidth=0.25pt](11.7073,26.3415)(13.1707,27.8049)(11.7073,29.2683)(10.2439,27.8049)(11.7073,26.3415) \psline[linewidth=0.25pt](13.1707,27.8049)(14.6341,29.2683)(13.1707,30.7317)(11.7073,29.2683)(13.1707,27.8049) \psline[linewidth=0.25pt](14.6341,29.2683)(16.0976,30.7317)(14.6341,32.1951)(13.1707,30.7317)(14.6341,29.2683) \psline[linewidth=0.25pt](16.0976,30.7317)(17.5610,32.1951)(16.0976,33.6585)(14.6341,32.1951)(16.0976,30.7317) \psline[linewidth=0.25pt](17.5610,32.1951)(19.0244,33.6585)(17.5610,35.1220)(16.0976,33.6585)(17.5610,32.1951) \psline[linewidth=0.25pt](19.0244,33.6585)(20.4878,35.1220)(19.0244,36.5854)(17.5610,35.1220)(19.0244,33.6585) \psline[linewidth=0.25pt](20.4878,35.1220)(21.9512,36.5854)(20.4878,38.0488)(19.0244,36.5854)(20.4878,35.1220) \psline[linewidth=0.25pt](21.9512,36.5854)(23.4146,38.0488)(21.9512,39.5122)(20.4878,38.0488)(21.9512,36.5854) \psline[linewidth=0.25pt](23.4146,38.0488)(24.8780,39.5122)(23.4146,40.9756)(21.9512,39.5122)(23.4146,38.0488) \psline[linewidth=0.25pt](24.8780,39.5122)(26.3415,40.9756)(24.8780,42.4390)(23.4146,40.9756)(24.8780,39.5122) \psline[linewidth=0.25pt](26.3415,40.9756)(27.8049,42.4390)(26.3415,43.9024)(24.8780,42.4390)(26.3415,40.9756) \psline[linewidth=0.25pt](27.8049,42.4390)(29.2683,43.9024)(27.8049,45.3659)(26.3415,43.9024)(27.8049,42.4390) \psline[linewidth=0.25pt](29.2683,43.9024)(30.7317,45.3659)(29.2683,46.8293)(27.8049,45.3659)(29.2683,43.9024) \psline[linewidth=0.25pt](-8.78049,8.78049)(-7.31707,10.2439)(-8.78049,11.7073)(-10.2439,10.2439)(-8.78049,8.78049) \psline[linewidth=0.25pt](-7.31707,10.2439)(-5.85366,11.7073)(-7.31707,13.1707)(-8.78049,11.7073)(-7.31707,10.2439) \psline[linewidth=0.25pt](-5.85366,11.7073)(-4.39024,13.1707)(-5.85366,14.6341)(-7.31707,13.1707)(-5.85366,11.7073) \psline[linewidth=0.25pt](-4.39024,13.1707)(-2.92683,14.6341)(-4.39024,16.0976)(-5.85366,14.6341)(-4.39024,13.1707) \psline[linewidth=0.25pt](-2.92683,14.6341)(-1.46341,16.0976)(-2.92683,17.5610)(-4.39024,16.0976)(-2.92683,14.6341) \psline[linewidth=0.25pt](-1.46341,16.0976)(0.000000,17.5610)(-1.46341,19.0244)(-2.92683,17.5610)(-1.46341,16.0976) \psline[linewidth=0.25pt](0.000000,17.5610)(1.46341,19.0244)(0.000000,20.4878)(-1.46341,19.0244)(0.000000,17.5610) \psline[linewidth=0.25pt](1.46341,19.0244)(2.92683,20.4878)(1.46341,21.9512)(0.000000,20.4878)(1.46341,19.0244) \psline[linewidth=0.25pt](2.92683,20.4878)(4.39024,21.9512)(2.92683,23.4146)(1.46341,21.9512)(2.92683,20.4878) \psline[linewidth=0.25pt](4.39024,21.9512)(5.85366,23.4146)(4.39024,24.8780)(2.92683,23.4146)(4.39024,21.9512) \psline[linewidth=0.25pt](5.85366,23.4146)(7.31707,24.8780)(5.85366,26.3415)(4.39024,24.8780)(5.85366,23.4146) \psline[linewidth=0.25pt](7.31707,24.8780)(8.78049,26.3415)(7.31707,27.8049)(5.85366,26.3415)(7.31707,24.8780) \psline[linewidth=0.25pt](8.78049,26.3415)(10.2439,27.8049)(8.78049,29.2683)(7.31707,27.8049)(8.78049,26.3415) \psline[linewidth=0.25pt](10.2439,27.8049)(11.7073,29.2683)(10.2439,30.7317)(8.78049,29.2683)(10.2439,27.8049) \psline[linewidth=0.25pt](11.7073,29.2683)(13.1707,30.7317)(11.7073,32.1951)(10.2439,30.7317)(11.7073,29.2683) \psline[linewidth=0.25pt](13.1707,30.7317)(14.6341,32.1951)(13.1707,33.6585)(11.7073,32.1951)(13.1707,30.7317) \psline[linewidth=0.25pt](14.6341,32.1951)(16.0976,33.6585)(14.6341,35.1220)(13.1707,33.6585)(14.6341,32.1951) \psline[linewidth=0.25pt](16.0976,33.6585)(17.5610,35.1220)(16.0976,36.5854)(14.6341,35.1220)(16.0976,33.6585) \psline[linewidth=0.25pt](17.5610,35.1220)(19.0244,36.5854)(17.5610,38.0488)(16.0976,36.5854)(17.5610,35.1220) \psline[linewidth=0.25pt](19.0244,36.5854)(20.4878,38.0488)(19.0244,39.5122)(17.5610,38.0488)(19.0244,36.5854) \psline[linewidth=0.25pt](20.4878,38.0488)(21.9512,39.5122)(20.4878,40.9756)(19.0244,39.5122)(20.4878,38.0488) \psline[linewidth=0.25pt](21.9512,39.5122)(23.4146,40.9756)(21.9512,42.4390)(20.4878,40.9756)(21.9512,39.5122) \psline[linewidth=0.25pt](23.4146,40.9756)(24.8780,42.4390)(23.4146,43.9024)(21.9512,42.4390)(23.4146,40.9756) \psline[linewidth=0.25pt](24.8780,42.4390)(26.3415,43.9024)(24.8780,45.3659)(23.4146,43.9024)(24.8780,42.4390) \psline[linewidth=0.25pt](-10.2439,10.2439)(-8.78049,11.7073)(-10.2439,13.1707)(-11.7073,11.7073)(-10.2439,10.2439) \psline[linewidth=0.25pt](-8.78049,11.7073)(-7.31707,13.1707)(-8.78049,14.6341)(-10.2439,13.1707)(-8.78049,11.7073) \psline[linewidth=0.25pt](-7.31707,13.1707)(-5.85366,14.6341)(-7.31707,16.0976)(-8.78049,14.6341)(-7.31707,13.1707) \psline[linewidth=0.25pt](-5.85366,14.6341)(-4.39024,16.0976)(-5.85366,17.5610)(-7.31707,16.0976)(-5.85366,14.6341) \psline[linewidth=0.25pt](-4.39024,16.0976)(-2.92683,17.5610)(-4.39024,19.0244)(-5.85366,17.5610)(-4.39024,16.0976) \psline[linewidth=0.25pt](-2.92683,17.5610)(-1.46341,19.0244)(-2.92683,20.4878)(-4.39024,19.0244)(-2.92683,17.5610) \psline[linewidth=0.25pt](-1.46341,19.0244)(0.000000,20.4878)(-1.46341,21.9512)(-2.92683,20.4878)(-1.46341,19.0244) \psline[linewidth=0.25pt](0.000000,20.4878)(1.46341,21.9512)(0.000000,23.4146)(-1.46341,21.9512)(0.000000,20.4878) \psline[linewidth=0.25pt](1.46341,21.9512)(2.92683,23.4146)(1.46341,24.8780)(0.000000,23.4146)(1.46341,21.9512) \psline[linewidth=0.25pt](2.92683,23.4146)(4.39024,24.8780)(2.92683,26.3415)(1.46341,24.8780)(2.92683,23.4146) \psline[linewidth=0.25pt](4.39024,24.8780)(5.85366,26.3415)(4.39024,27.8049)(2.92683,26.3415)(4.39024,24.8780) \psline[linewidth=0.25pt](5.85366,26.3415)(7.31707,27.8049)(5.85366,29.2683)(4.39024,27.8049)(5.85366,26.3415) \psline[linewidth=0.25pt](7.31707,27.8049)(8.78049,29.2683)(7.31707,30.7317)(5.85366,29.2683)(7.31707,27.8049) \psline[linewidth=0.25pt](8.78049,29.2683)(10.2439,30.7317)(8.78049,32.1951)(7.31707,30.7317)(8.78049,29.2683) \psline[linewidth=0.25pt](10.2439,30.7317)(11.7073,32.1951)(10.2439,33.6585)(8.78049,32.1951)(10.2439,30.7317) \psline[linewidth=0.25pt](11.7073,32.1951)(13.1707,33.6585)(11.7073,35.1220)(10.2439,33.6585)(11.7073,32.1951) \psline[linewidth=0.25pt](13.1707,33.6585)(14.6341,35.1220)(13.1707,36.5854)(11.7073,35.1220)(13.1707,33.6585) \psline[linewidth=0.25pt](14.6341,35.1220)(16.0976,36.5854)(14.6341,38.0488)(13.1707,36.5854)(14.6341,35.1220) \psline[linewidth=0.25pt](16.0976,36.5854)(17.5610,38.0488)(16.0976,39.5122)(14.6341,38.0488)(16.0976,36.5854) \psline[linewidth=0.25pt](17.5610,38.0488)(19.0244,39.5122)(17.5610,40.9756)(16.0976,39.5122)(17.5610,38.0488) \psline[linewidth=0.25pt](19.0244,39.5122)(20.4878,40.9756)(19.0244,42.4390)(17.5610,40.9756)(19.0244,39.5122) \psline[linewidth=0.25pt](20.4878,40.9756)(21.9512,42.4390)(20.4878,43.9024)(19.0244,42.4390)(20.4878,40.9756) \psline[linewidth=0.25pt](21.9512,42.4390)(23.4146,43.9024)(21.9512,45.3659)(20.4878,43.9024)(21.9512,42.4390) \psline[linewidth=0.25pt](23.4146,43.9024)(24.8780,45.3659)(23.4146,46.8293)(21.9512,45.3659)(23.4146,43.9024) \psline[linewidth=0.25pt](-11.7073,11.7073)(-10.2439,13.1707)(-11.7073,14.6341)(-13.1707,13.1707)(-11.7073,11.7073) \psline[linewidth=0.25pt](-10.2439,13.1707)(-8.78049,14.6341)(-10.2439,16.0976)(-11.7073,14.6341)(-10.2439,13.1707) \psline[linewidth=0.25pt](-8.78049,14.6341)(-7.31707,16.0976)(-8.78049,17.5610)(-10.2439,16.0976)(-8.78049,14.6341) \psline[linewidth=0.25pt](-7.31707,16.0976)(-5.85366,17.5610)(-7.31707,19.0244)(-8.78049,17.5610)(-7.31707,16.0976) \psline[linewidth=0.25pt](-5.85366,17.5610)(-4.39024,19.0244)(-5.85366,20.4878)(-7.31707,19.0244)(-5.85366,17.5610) \psline[linewidth=0.25pt](-4.39024,19.0244)(-2.92683,20.4878)(-4.39024,21.9512)(-5.85366,20.4878)(-4.39024,19.0244) \psline[linewidth=0.25pt](-2.92683,20.4878)(-1.46341,21.9512)(-2.92683,23.4146)(-4.39024,21.9512)(-2.92683,20.4878) \psline[linewidth=0.25pt](-1.46341,21.9512)(0.000000,23.4146)(-1.46341,24.8780)(-2.92683,23.4146)(-1.46341,21.9512) \psline[linewidth=0.25pt](0.000000,23.4146)(1.46341,24.8780)(0.000000,26.3415)(-1.46341,24.8780)(0.000000,23.4146) \psline[linewidth=0.25pt](1.46341,24.8780)(2.92683,26.3415)(1.46341,27.8049)(0.000000,26.3415)(1.46341,24.8780) \psline[linewidth=0.25pt](2.92683,26.3415)(4.39024,27.8049)(2.92683,29.2683)(1.46341,27.8049)(2.92683,26.3415) \psline[linewidth=0.25pt](4.39024,27.8049)(5.85366,29.2683)(4.39024,30.7317)(2.92683,29.2683)(4.39024,27.8049) \psline[linewidth=0.25pt](5.85366,29.2683)(7.31707,30.7317)(5.85366,32.1951)(4.39024,30.7317)(5.85366,29.2683) \psline[linewidth=0.25pt](7.31707,30.7317)(8.78049,32.1951)(7.31707,33.6585)(5.85366,32.1951)(7.31707,30.7317) \psline[linewidth=0.25pt](8.78049,32.1951)(10.2439,33.6585)(8.78049,35.1220)(7.31707,33.6585)(8.78049,32.1951) \psline[linewidth=0.25pt](10.2439,33.6585)(11.7073,35.1220)(10.2439,36.5854)(8.78049,35.1220)(10.2439,33.6585) \psline[linewidth=0.25pt](11.7073,35.1220)(13.1707,36.5854)(11.7073,38.0488)(10.2439,36.5854)(11.7073,35.1220) \psline[linewidth=0.25pt](13.1707,36.5854)(14.6341,38.0488)(13.1707,39.5122)(11.7073,38.0488)(13.1707,36.5854) \psline[linewidth=0.25pt](14.6341,38.0488)(16.0976,39.5122)(14.6341,40.9756)(13.1707,39.5122)(14.6341,38.0488) \psline[linewidth=0.25pt](16.0976,39.5122)(17.5610,40.9756)(16.0976,42.4390)(14.6341,40.9756)(16.0976,39.5122) \psline[linewidth=0.25pt](17.5610,40.9756)(19.0244,42.4390)(17.5610,43.9024)(16.0976,42.4390)(17.5610,40.9756) \psline[linewidth=0.25pt](-13.1707,13.1707)(-11.7073,14.6341)(-13.1707,16.0976)(-14.6341,14.6341)(-13.1707,13.1707) \psline[linewidth=0.25pt](-11.7073,14.6341)(-10.2439,16.0976)(-11.7073,17.5610)(-13.1707,16.0976)(-11.7073,14.6341) \psline[linewidth=0.25pt](-10.2439,16.0976)(-8.78049,17.5610)(-10.2439,19.0244)(-11.7073,17.5610)(-10.2439,16.0976) \psline[linewidth=0.25pt](-8.78049,17.5610)(-7.31707,19.0244)(-8.78049,20.4878)(-10.2439,19.0244)(-8.78049,17.5610) \psline[linewidth=0.25pt](-7.31707,19.0244)(-5.85366,20.4878)(-7.31707,21.9512)(-8.78049,20.4878)(-7.31707,19.0244) \psline[linewidth=0.25pt](-5.85366,20.4878)(-4.39024,21.9512)(-5.85366,23.4146)(-7.31707,21.9512)(-5.85366,20.4878) \psline[linewidth=0.25pt](-4.39024,21.9512)(-2.92683,23.4146)(-4.39024,24.8780)(-5.85366,23.4146)(-4.39024,21.9512) \psline[linewidth=0.25pt](-2.92683,23.4146)(-1.46341,24.8780)(-2.92683,26.3415)(-4.39024,24.8780)(-2.92683,23.4146) \psline[linewidth=0.25pt](-1.46341,24.8780)(0.000000,26.3415)(-1.46341,27.8049)(-2.92683,26.3415)(-1.46341,24.8780) \psline[linewidth=0.25pt](0.000000,26.3415)(1.46341,27.8049)(0.000000,29.2683)(-1.46341,27.8049)(0.000000,26.3415) \psline[linewidth=0.25pt](1.46341,27.8049)(2.92683,29.2683)(1.46341,30.7317)(0.000000,29.2683)(1.46341,27.8049) \psline[linewidth=0.25pt](2.92683,29.2683)(4.39024,30.7317)(2.92683,32.1951)(1.46341,30.7317)(2.92683,29.2683) \psline[linewidth=0.25pt](4.39024,30.7317)(5.85366,32.1951)(4.39024,33.6585)(2.92683,32.1951)(4.39024,30.7317) \psline[linewidth=0.25pt](5.85366,32.1951)(7.31707,33.6585)(5.85366,35.1220)(4.39024,33.6585)(5.85366,32.1951) \psline[linewidth=0.25pt](7.31707,33.6585)(8.78049,35.1220)(7.31707,36.5854)(5.85366,35.1220)(7.31707,33.6585) \psline[linewidth=0.25pt](8.78049,35.1220)(10.2439,36.5854)(8.78049,38.0488)(7.31707,36.5854)(8.78049,35.1220) \psline[linewidth=0.25pt](10.2439,36.5854)(11.7073,38.0488)(10.2439,39.5122)(8.78049,38.0488)(10.2439,36.5854) \psline[linewidth=0.25pt](11.7073,38.0488)(13.1707,39.5122)(11.7073,40.9756)(10.2439,39.5122)(11.7073,38.0488) \psline[linewidth=0.25pt](13.1707,39.5122)(14.6341,40.9756)(13.1707,42.4390)(11.7073,40.9756)(13.1707,39.5122) \psline[linewidth=0.25pt](14.6341,40.9756)(16.0976,42.4390)(14.6341,43.9024)(13.1707,42.4390)(14.6341,40.9756) \psline[linewidth=0.25pt](-14.6341,14.6341)(-13.1707,16.0976)(-14.6341,17.5610)(-16.0976,16.0976)(-14.6341,14.6341) \psline[linewidth=0.25pt](-13.1707,16.0976)(-11.7073,17.5610)(-13.1707,19.0244)(-14.6341,17.5610)(-13.1707,16.0976) \psline[linewidth=0.25pt](-11.7073,17.5610)(-10.2439,19.0244)(-11.7073,20.4878)(-13.1707,19.0244)(-11.7073,17.5610) \psline[linewidth=0.25pt](-10.2439,19.0244)(-8.78049,20.4878)(-10.2439,21.9512)(-11.7073,20.4878)(-10.2439,19.0244) \psline[linewidth=0.25pt](-8.78049,20.4878)(-7.31707,21.9512)(-8.78049,23.4146)(-10.2439,21.9512)(-8.78049,20.4878) \psline[linewidth=0.25pt](-7.31707,21.9512)(-5.85366,23.4146)(-7.31707,24.8780)(-8.78049,23.4146)(-7.31707,21.9512) \psline[linewidth=0.25pt](-5.85366,23.4146)(-4.39024,24.8780)(-5.85366,26.3415)(-7.31707,24.8780)(-5.85366,23.4146) \psline[linewidth=0.25pt](-4.39024,24.8780)(-2.92683,26.3415)(-4.39024,27.8049)(-5.85366,26.3415)(-4.39024,24.8780) \psline[linewidth=0.25pt](-2.92683,26.3415)(-1.46341,27.8049)(-2.92683,29.2683)(-4.39024,27.8049)(-2.92683,26.3415) \psline[linewidth=0.25pt](-1.46341,27.8049)(0.000000,29.2683)(-1.46341,30.7317)(-2.92683,29.2683)(-1.46341,27.8049) \psline[linewidth=0.25pt](0.000000,29.2683)(1.46341,30.7317)(0.000000,32.1951)(-1.46341,30.7317)(0.000000,29.2683) \psline[linewidth=0.25pt](1.46341,30.7317)(2.92683,32.1951)(1.46341,33.6585)(0.000000,32.1951)(1.46341,30.7317) \psline[linewidth=0.25pt](2.92683,32.1951)(4.39024,33.6585)(2.92683,35.1220)(1.46341,33.6585)(2.92683,32.1951) \psline[linewidth=0.25pt](4.39024,33.6585)(5.85366,35.1220)(4.39024,36.5854)(2.92683,35.1220)(4.39024,33.6585) \psline[linewidth=0.25pt](5.85366,35.1220)(7.31707,36.5854)(5.85366,38.0488)(4.39024,36.5854)(5.85366,35.1220) \psline[linewidth=0.25pt](7.31707,36.5854)(8.78049,38.0488)(7.31707,39.5122)(5.85366,38.0488)(7.31707,36.5854) \psline[linewidth=0.25pt](8.78049,38.0488)(10.2439,39.5122)(8.78049,40.9756)(7.31707,39.5122)(8.78049,38.0488) \psline[linewidth=0.25pt](10.2439,39.5122)(11.7073,40.9756)(10.2439,42.4390)(8.78049,40.9756)(10.2439,39.5122) \psline[linewidth=0.25pt](11.7073,40.9756)(13.1707,42.4390)(11.7073,43.9024)(10.2439,42.4390)(11.7073,40.9756) \psline[linewidth=0.25pt](-16.0976,16.0976)(-14.6341,17.5610)(-16.0976,19.0244)(-17.5610,17.5610)(-16.0976,16.0976) \psline[linewidth=0.25pt](-14.6341,17.5610)(-13.1707,19.0244)(-14.6341,20.4878)(-16.0976,19.0244)(-14.6341,17.5610) \psline[linewidth=0.25pt](-13.1707,19.0244)(-11.7073,20.4878)(-13.1707,21.9512)(-14.6341,20.4878)(-13.1707,19.0244) \psline[linewidth=0.25pt](-11.7073,20.4878)(-10.2439,21.9512)(-11.7073,23.4146)(-13.1707,21.9512)(-11.7073,20.4878) \psline[linewidth=0.25pt](-10.2439,21.9512)(-8.78049,23.4146)(-10.2439,24.8780)(-11.7073,23.4146)(-10.2439,21.9512) \psline[linewidth=0.25pt](-8.78049,23.4146)(-7.31707,24.8780)(-8.78049,26.3415)(-10.2439,24.8780)(-8.78049,23.4146) \psline[linewidth=0.25pt](-7.31707,24.8780)(-5.85366,26.3415)(-7.31707,27.8049)(-8.78049,26.3415)(-7.31707,24.8780) \psline[linewidth=0.25pt](-5.85366,26.3415)(-4.39024,27.8049)(-5.85366,29.2683)(-7.31707,27.8049)(-5.85366,26.3415) \psline[linewidth=0.25pt](-4.39024,27.8049)(-2.92683,29.2683)(-4.39024,30.7317)(-5.85366,29.2683)(-4.39024,27.8049) \psline[linewidth=0.25pt](-2.92683,29.2683)(-1.46341,30.7317)(-2.92683,32.1951)(-4.39024,30.7317)(-2.92683,29.2683) \psline[linewidth=0.25pt](-1.46341,30.7317)(0.000000,32.1951)(-1.46341,33.6585)(-2.92683,32.1951)(-1.46341,30.7317) \psline[linewidth=0.25pt](0.000000,32.1951)(1.46341,33.6585)(0.000000,35.1220)(-1.46341,33.6585)(0.000000,32.1951) \psline[linewidth=0.25pt](1.46341,33.6585)(2.92683,35.1220)(1.46341,36.5854)(0.000000,35.1220)(1.46341,33.6585) \psline[linewidth=0.25pt](2.92683,35.1220)(4.39024,36.5854)(2.92683,38.0488)(1.46341,36.5854)(2.92683,35.1220) \psline[linewidth=0.25pt](4.39024,36.5854)(5.85366,38.0488)(4.39024,39.5122)(2.92683,38.0488)(4.39024,36.5854) \psline[linewidth=0.25pt](5.85366,38.0488)(7.31707,39.5122)(5.85366,40.9756)(4.39024,39.5122)(5.85366,38.0488) \psline[linewidth=0.25pt](7.31707,39.5122)(8.78049,40.9756)(7.31707,42.4390)(5.85366,40.9756)(7.31707,39.5122) \psline[linewidth=0.25pt](-17.5610,17.5610)(-16.0976,19.0244)(-17.5610,20.4878)(-19.0244,19.0244)(-17.5610,17.5610) \psline[linewidth=0.25pt](-16.0976,19.0244)(-14.6341,20.4878)(-16.0976,21.9512)(-17.5610,20.4878)(-16.0976,19.0244) \psline[linewidth=0.25pt](-14.6341,20.4878)(-13.1707,21.9512)(-14.6341,23.4146)(-16.0976,21.9512)(-14.6341,20.4878) \psline[linewidth=0.25pt](-13.1707,21.9512)(-11.7073,23.4146)(-13.1707,24.8780)(-14.6341,23.4146)(-13.1707,21.9512) \psline[linewidth=0.25pt](-11.7073,23.4146)(-10.2439,24.8780)(-11.7073,26.3415)(-13.1707,24.8780)(-11.7073,23.4146) \psline[linewidth=0.25pt](-10.2439,24.8780)(-8.78049,26.3415)(-10.2439,27.8049)(-11.7073,26.3415)(-10.2439,24.8780) \psline[linewidth=0.25pt](-8.78049,26.3415)(-7.31707,27.8049)(-8.78049,29.2683)(-10.2439,27.8049)(-8.78049,26.3415) \psline[linewidth=0.25pt](-7.31707,27.8049)(-5.85366,29.2683)(-7.31707,30.7317)(-8.78049,29.2683)(-7.31707,27.8049) \psline[linewidth=0.25pt](-5.85366,29.2683)(-4.39024,30.7317)(-5.85366,32.1951)(-7.31707,30.7317)(-5.85366,29.2683) \psline[linewidth=0.25pt](-4.39024,30.7317)(-2.92683,32.1951)(-4.39024,33.6585)(-5.85366,32.1951)(-4.39024,30.7317) \psline[linewidth=0.25pt](-2.92683,32.1951)(-1.46341,33.6585)(-2.92683,35.1220)(-4.39024,33.6585)(-2.92683,32.1951) \psline[linewidth=0.25pt](-1.46341,33.6585)(0.000000,35.1220)(-1.46341,36.5854)(-2.92683,35.1220)(-1.46341,33.6585) \psline[linewidth=0.25pt](0.000000,35.1220)(1.46341,36.5854)(0.000000,38.0488)(-1.46341,36.5854)(0.000000,35.1220) \psline[linewidth=0.25pt](1.46341,36.5854)(2.92683,38.0488)(1.46341,39.5122)(0.000000,38.0488)(1.46341,36.5854) \psline[linewidth=0.25pt](2.92683,38.0488)(4.39024,39.5122)(2.92683,40.9756)(1.46341,39.5122)(2.92683,38.0488) \psline[linewidth=0.25pt](4.39024,39.5122)(5.85366,40.9756)(4.39024,42.4390)(2.92683,40.9756)(4.39024,39.5122) \psline[linewidth=0.25pt](-19.0244,19.0244)(-17.5610,20.4878)(-19.0244,21.9512)(-20.4878,20.4878)(-19.0244,19.0244) \psline[linewidth=0.25pt](-17.5610,20.4878)(-16.0976,21.9512)(-17.5610,23.4146)(-19.0244,21.9512)(-17.5610,20.4878) \psline[linewidth=0.25pt](-16.0976,21.9512)(-14.6341,23.4146)(-16.0976,24.8780)(-17.5610,23.4146)(-16.0976,21.9512) \psline[linewidth=0.25pt](-14.6341,23.4146)(-13.1707,24.8780)(-14.6341,26.3415)(-16.0976,24.8780)(-14.6341,23.4146) \psline[linewidth=0.25pt](-13.1707,24.8780)(-11.7073,26.3415)(-13.1707,27.8049)(-14.6341,26.3415)(-13.1707,24.8780) \psline[linewidth=0.25pt](-11.7073,26.3415)(-10.2439,27.8049)(-11.7073,29.2683)(-13.1707,27.8049)(-11.7073,26.3415) \psline[linewidth=0.25pt](-10.2439,27.8049)(-8.78049,29.2683)(-10.2439,30.7317)(-11.7073,29.2683)(-10.2439,27.8049) \psline[linewidth=0.25pt](-8.78049,29.2683)(-7.31707,30.7317)(-8.78049,32.1951)(-10.2439,30.7317)(-8.78049,29.2683) \psline[linewidth=0.25pt](-7.31707,30.7317)(-5.85366,32.1951)(-7.31707,33.6585)(-8.78049,32.1951)(-7.31707,30.7317) \psline[linewidth=0.25pt](-5.85366,32.1951)(-4.39024,33.6585)(-5.85366,35.1220)(-7.31707,33.6585)(-5.85366,32.1951) \psline[linewidth=0.25pt](-4.39024,33.6585)(-2.92683,35.1220)(-4.39024,36.5854)(-5.85366,35.1220)(-4.39024,33.6585) \psline[linewidth=0.25pt](-2.92683,35.1220)(-1.46341,36.5854)(-2.92683,38.0488)(-4.39024,36.5854)(-2.92683,35.1220) \psline[linewidth=0.25pt](-1.46341,36.5854)(0.000000,38.0488)(-1.46341,39.5122)(-2.92683,38.0488)(-1.46341,36.5854) \psline[linewidth=0.25pt](0.000000,38.0488)(1.46341,39.5122)(0.000000,40.9756)(-1.46341,39.5122)(0.000000,38.0488) \psline[linewidth=0.25pt](1.46341,39.5122)(2.92683,40.9756)(1.46341,42.4390)(0.000000,40.9756)(1.46341,39.5122) \psline[linewidth=0.25pt](2.92683,40.9756)(4.39024,42.4390)(2.92683,43.9024)(1.46341,42.4390)(2.92683,40.9756) \psline[linewidth=0.25pt](-20.4878,20.4878)(-19.0244,21.9512)(-20.4878,23.4146)(-21.9512,21.9512)(-20.4878,20.4878) \psline[linewidth=0.25pt](-19.0244,21.9512)(-17.5610,23.4146)(-19.0244,24.8780)(-20.4878,23.4146)(-19.0244,21.9512) \psline[linewidth=0.25pt](-17.5610,23.4146)(-16.0976,24.8780)(-17.5610,26.3415)(-19.0244,24.8780)(-17.5610,23.4146) \psline[linewidth=0.25pt](-16.0976,24.8780)(-14.6341,26.3415)(-16.0976,27.8049)(-17.5610,26.3415)(-16.0976,24.8780) \psline[linewidth=0.25pt](-14.6341,26.3415)(-13.1707,27.8049)(-14.6341,29.2683)(-16.0976,27.8049)(-14.6341,26.3415) \psline[linewidth=0.25pt](-13.1707,27.8049)(-11.7073,29.2683)(-13.1707,30.7317)(-14.6341,29.2683)(-13.1707,27.8049) \psline[linewidth=0.25pt](-11.7073,29.2683)(-10.2439,30.7317)(-11.7073,32.1951)(-13.1707,30.7317)(-11.7073,29.2683) \psline[linewidth=0.25pt](-10.2439,30.7317)(-8.78049,32.1951)(-10.2439,33.6585)(-11.7073,32.1951)(-10.2439,30.7317) \psline[linewidth=0.25pt](-8.78049,32.1951)(-7.31707,33.6585)(-8.78049,35.1220)(-10.2439,33.6585)(-8.78049,32.1951) \psline[linewidth=0.25pt](-7.31707,33.6585)(-5.85366,35.1220)(-7.31707,36.5854)(-8.78049,35.1220)(-7.31707,33.6585) \psline[linewidth=0.25pt](-5.85366,35.1220)(-4.39024,36.5854)(-5.85366,38.0488)(-7.31707,36.5854)(-5.85366,35.1220) \psline[linewidth=0.25pt](-4.39024,36.5854)(-2.92683,38.0488)(-4.39024,39.5122)(-5.85366,38.0488)(-4.39024,36.5854) \psline[linewidth=0.25pt](-2.92683,38.0488)(-1.46341,39.5122)(-2.92683,40.9756)(-4.39024,39.5122)(-2.92683,38.0488) \psline[linewidth=0.25pt](-21.9512,21.9512)(-20.4878,23.4146)(-21.9512,24.8780)(-23.4146,23.4146)(-21.9512,21.9512) \psline[linewidth=0.25pt](-20.4878,23.4146)(-19.0244,24.8780)(-20.4878,26.3415)(-21.9512,24.8780)(-20.4878,23.4146) \psline[linewidth=0.25pt](-19.0244,24.8780)(-17.5610,26.3415)(-19.0244,27.8049)(-20.4878,26.3415)(-19.0244,24.8780) \psline[linewidth=0.25pt](-17.5610,26.3415)(-16.0976,27.8049)(-17.5610,29.2683)(-19.0244,27.8049)(-17.5610,26.3415) \psline[linewidth=0.25pt](-16.0976,27.8049)(-14.6341,29.2683)(-16.0976,30.7317)(-17.5610,29.2683)(-16.0976,27.8049) \psline[linewidth=0.25pt](-14.6341,29.2683)(-13.1707,30.7317)(-14.6341,32.1951)(-16.0976,30.7317)(-14.6341,29.2683) \psline[linewidth=0.25pt](-13.1707,30.7317)(-11.7073,32.1951)(-13.1707,33.6585)(-14.6341,32.1951)(-13.1707,30.7317) \psline[linewidth=0.25pt](-11.7073,32.1951)(-10.2439,33.6585)(-11.7073,35.1220)(-13.1707,33.6585)(-11.7073,32.1951) \psline[linewidth=0.25pt](-10.2439,33.6585)(-8.78049,35.1220)(-10.2439,36.5854)(-11.7073,35.1220)(-10.2439,33.6585) \psline[linewidth=0.25pt](-8.78049,35.1220)(-7.31707,36.5854)(-8.78049,38.0488)(-10.2439,36.5854)(-8.78049,35.1220) \psline[linewidth=0.25pt](-7.31707,36.5854)(-5.85366,38.0488)(-7.31707,39.5122)(-8.78049,38.0488)(-7.31707,36.5854) \psline[linewidth=0.25pt](-5.85366,38.0488)(-4.39024,39.5122)(-5.85366,40.9756)(-7.31707,39.5122)(-5.85366,38.0488) \psline[linewidth=0.25pt](-23.4146,23.4146)(-21.9512,24.8780)(-23.4146,26.3415)(-24.8780,24.8780)(-23.4146,23.4146) \psline[linewidth=0.25pt](-21.9512,24.8780)(-20.4878,26.3415)(-21.9512,27.8049)(-23.4146,26.3415)(-21.9512,24.8780) \psline[linewidth=0.25pt](-20.4878,26.3415)(-19.0244,27.8049)(-20.4878,29.2683)(-21.9512,27.8049)(-20.4878,26.3415) \psline[linewidth=0.25pt](-19.0244,27.8049)(-17.5610,29.2683)(-19.0244,30.7317)(-20.4878,29.2683)(-19.0244,27.8049) \psline[linewidth=0.25pt](-17.5610,29.2683)(-16.0976,30.7317)(-17.5610,32.1951)(-19.0244,30.7317)(-17.5610,29.2683) \psline[linewidth=0.25pt](-16.0976,30.7317)(-14.6341,32.1951)(-16.0976,33.6585)(-17.5610,32.1951)(-16.0976,30.7317) \psline[linewidth=0.25pt](-14.6341,32.1951)(-13.1707,33.6585)(-14.6341,35.1220)(-16.0976,33.6585)(-14.6341,32.1951) \psline[linewidth=0.25pt](-13.1707,33.6585)(-11.7073,35.1220)(-13.1707,36.5854)(-14.6341,35.1220)(-13.1707,33.6585) \psline[linewidth=0.25pt](-11.7073,35.1220)(-10.2439,36.5854)(-11.7073,38.0488)(-13.1707,36.5854)(-11.7073,35.1220) \psline[linewidth=0.25pt](-10.2439,36.5854)(-8.78049,38.0488)(-10.2439,39.5122)(-11.7073,38.0488)(-10.2439,36.5854) \psline[linewidth=0.25pt](-8.78049,38.0488)(-7.31707,39.5122)(-8.78049,40.9756)(-10.2439,39.5122)(-8.78049,38.0488) \psline[linewidth=0.25pt](-7.31707,39.5122)(-5.85366,40.9756)(-7.31707,42.4390)(-8.78049,40.9756)(-7.31707,39.5122) \psline[linewidth=0.25pt](-24.8780,24.8780)(-23.4146,26.3415)(-24.8780,27.8049)(-26.3415,26.3415)(-24.8780,24.8780) \psline[linewidth=0.25pt](-23.4146,26.3415)(-21.9512,27.8049)(-23.4146,29.2683)(-24.8780,27.8049)(-23.4146,26.3415) \psline[linewidth=0.25pt](-21.9512,27.8049)(-20.4878,29.2683)(-21.9512,30.7317)(-23.4146,29.2683)(-21.9512,27.8049) \psline[linewidth=0.25pt](-20.4878,29.2683)(-19.0244,30.7317)(-20.4878,32.1951)(-21.9512,30.7317)(-20.4878,29.2683) \psline[linewidth=0.25pt](-19.0244,30.7317)(-17.5610,32.1951)(-19.0244,33.6585)(-20.4878,32.1951)(-19.0244,30.7317) \psline[linewidth=0.25pt](-17.5610,32.1951)(-16.0976,33.6585)(-17.5610,35.1220)(-19.0244,33.6585)(-17.5610,32.1951) \psline[linewidth=0.25pt](-16.0976,33.6585)(-14.6341,35.1220)(-16.0976,36.5854)(-17.5610,35.1220)(-16.0976,33.6585) \psline[linewidth=0.25pt](-14.6341,35.1220)(-13.1707,36.5854)(-14.6341,38.0488)(-16.0976,36.5854)(-14.6341,35.1220) \psline[linewidth=0.25pt](-13.1707,36.5854)(-11.7073,38.0488)(-13.1707,39.5122)(-14.6341,38.0488)(-13.1707,36.5854) \psline[linewidth=0.25pt](-11.7073,38.0488)(-10.2439,39.5122)(-11.7073,40.9756)(-13.1707,39.5122)(-11.7073,38.0488) \psline[linewidth=0.25pt](-10.2439,39.5122)(-8.78049,40.9756)(-10.2439,42.4390)(-11.7073,40.9756)(-10.2439,39.5122) \psline[linewidth=0.25pt](-26.3415,26.3415)(-24.8780,27.8049)(-26.3415,29.2683)(-27.8049,27.8049)(-26.3415,26.3415) \psline[linewidth=0.25pt](-24.8780,27.8049)(-23.4146,29.2683)(-24.8780,30.7317)(-26.3415,29.2683)(-24.8780,27.8049) \psline[linewidth=0.25pt](-23.4146,29.2683)(-21.9512,30.7317)(-23.4146,32.1951)(-24.8780,30.7317)(-23.4146,29.2683) \psline[linewidth=0.25pt](-21.9512,30.7317)(-20.4878,32.1951)(-21.9512,33.6585)(-23.4146,32.1951)(-21.9512,30.7317) \psline[linewidth=0.25pt](-20.4878,32.1951)(-19.0244,33.6585)(-20.4878,35.1220)(-21.9512,33.6585)(-20.4878,32.1951) \psline[linewidth=0.25pt](-19.0244,33.6585)(-17.5610,35.1220)(-19.0244,36.5854)(-20.4878,35.1220)(-19.0244,33.6585) \psline[linewidth=0.25pt](-17.5610,35.1220)(-16.0976,36.5854)(-17.5610,38.0488)(-19.0244,36.5854)(-17.5610,35.1220) \psline[linewidth=0.25pt](-16.0976,36.5854)(-14.6341,38.0488)(-16.0976,39.5122)(-17.5610,38.0488)(-16.0976,36.5854) \psline[linewidth=0.25pt](-14.6341,38.0488)(-13.1707,39.5122)(-14.6341,40.9756)(-16.0976,39.5122)(-14.6341,38.0488) \psline[linewidth=0.25pt](-13.1707,39.5122)(-11.7073,40.9756)(-13.1707,42.4390)(-14.6341,40.9756)(-13.1707,39.5122) \psline[linewidth=0.25pt](-11.7073,40.9756)(-10.2439,42.4390)(-11.7073,43.9024)(-13.1707,42.4390)(-11.7073,40.9756) \psline[linewidth=0.25pt](-27.8049,27.8049)(-26.3415,29.2683)(-27.8049,30.7317)(-29.2683,29.2683)(-27.8049,27.8049) \psline[linewidth=0.25pt](-26.3415,29.2683)(-24.8780,30.7317)(-26.3415,32.1951)(-27.8049,30.7317)(-26.3415,29.2683) \psline[linewidth=0.25pt](-24.8780,30.7317)(-23.4146,32.1951)(-24.8780,33.6585)(-26.3415,32.1951)(-24.8780,30.7317) \psline[linewidth=0.25pt](-23.4146,32.1951)(-21.9512,33.6585)(-23.4146,35.1220)(-24.8780,33.6585)(-23.4146,32.1951) \psline[linewidth=0.25pt](-21.9512,33.6585)(-20.4878,35.1220)(-21.9512,36.5854)(-23.4146,35.1220)(-21.9512,33.6585) \psline[linewidth=0.25pt](-20.4878,35.1220)(-19.0244,36.5854)(-20.4878,38.0488)(-21.9512,36.5854)(-20.4878,35.1220) \psline[linewidth=0.25pt](-19.0244,36.5854)(-17.5610,38.0488)(-19.0244,39.5122)(-20.4878,38.0488)(-19.0244,36.5854) \psline[linewidth=0.25pt](-17.5610,38.0488)(-16.0976,39.5122)(-17.5610,40.9756)(-19.0244,39.5122)(-17.5610,38.0488) \psline[linewidth=0.25pt](-16.0976,39.5122)(-14.6341,40.9756)(-16.0976,42.4390)(-17.5610,40.9756)(-16.0976,39.5122) \psline[linewidth=0.25pt](-14.6341,40.9756)(-13.1707,42.4390)(-14.6341,43.9024)(-16.0976,42.4390)(-14.6341,40.9756) \psline[linewidth=0.25pt](-29.2683,29.2683)(-27.8049,30.7317)(-29.2683,32.1951)(-30.7317,30.7317)(-29.2683,29.2683) \psline[linewidth=0.25pt](-27.8049,30.7317)(-26.3415,32.1951)(-27.8049,33.6585)(-29.2683,32.1951)(-27.8049,30.7317) \psline[linewidth=0.25pt](-26.3415,32.1951)(-24.8780,33.6585)(-26.3415,35.1220)(-27.8049,33.6585)(-26.3415,32.1951) \psline[linewidth=0.25pt](-24.8780,33.6585)(-23.4146,35.1220)(-24.8780,36.5854)(-26.3415,35.1220)(-24.8780,33.6585) \psline[linewidth=0.25pt](-23.4146,35.1220)(-21.9512,36.5854)(-23.4146,38.0488)(-24.8780,36.5854)(-23.4146,35.1220) \psline[linewidth=0.25pt](-21.9512,36.5854)(-20.4878,38.0488)(-21.9512,39.5122)(-23.4146,38.0488)(-21.9512,36.5854) \psline[linewidth=0.25pt](-20.4878,38.0488)(-19.0244,39.5122)(-20.4878,40.9756)(-21.9512,39.5122)(-20.4878,38.0488) \psline[linewidth=0.25pt](-19.0244,39.5122)(-17.5610,40.9756)(-19.0244,42.4390)(-20.4878,40.9756)(-19.0244,39.5122) \psline[linewidth=0.25pt](-17.5610,40.9756)(-16.0976,42.4390)(-17.5610,43.9024)(-19.0244,42.4390)(-17.5610,40.9756) \psline[linewidth=0.25pt](-16.0976,42.4390)(-14.6341,43.9024)(-16.0976,45.3659)(-17.5610,43.9024)(-16.0976,42.4390) \psline[linewidth=0.25pt](-30.7317,30.7317)(-29.2683,32.1951)(-30.7317,33.6585)(-32.1951,32.1951)(-30.7317,30.7317) \psline[linewidth=0.25pt](-29.2683,32.1951)(-27.8049,33.6585)(-29.2683,35.1220)(-30.7317,33.6585)(-29.2683,32.1951) \psline[linewidth=0.25pt](-27.8049,33.6585)(-26.3415,35.1220)(-27.8049,36.5854)(-29.2683,35.1220)(-27.8049,33.6585) \psline[linewidth=0.25pt](-26.3415,35.1220)(-24.8780,36.5854)(-26.3415,38.0488)(-27.8049,36.5854)(-26.3415,35.1220) \psline[linewidth=0.25pt](-24.8780,36.5854)(-23.4146,38.0488)(-24.8780,39.5122)(-26.3415,38.0488)(-24.8780,36.5854) \psline[linewidth=0.25pt](-23.4146,38.0488)(-21.9512,39.5122)(-23.4146,40.9756)(-24.8780,39.5122)(-23.4146,38.0488) \psline[linewidth=0.25pt](-21.9512,39.5122)(-20.4878,40.9756)(-21.9512,42.4390)(-23.4146,40.9756)(-21.9512,39.5122) \psline[linewidth=0.25pt](-20.4878,40.9756)(-19.0244,42.4390)(-20.4878,43.9024)(-21.9512,42.4390)(-20.4878,40.9756) \psline[linewidth=0.25pt](-19.0244,42.4390)(-17.5610,43.9024)(-19.0244,45.3659)(-20.4878,43.9024)(-19.0244,42.4390) \psline[linewidth=0.25pt](-32.1951,32.1951)(-30.7317,33.6585)(-32.1951,35.1220)(-33.6585,33.6585)(-32.1951,32.1951) \psline[linewidth=0.25pt](-30.7317,33.6585)(-29.2683,35.1220)(-30.7317,36.5854)(-32.1951,35.1220)(-30.7317,33.6585) \psline[linewidth=0.25pt](-29.2683,35.1220)(-27.8049,36.5854)(-29.2683,38.0488)(-30.7317,36.5854)(-29.2683,35.1220) \psline[linewidth=0.25pt](-27.8049,36.5854)(-26.3415,38.0488)(-27.8049,39.5122)(-29.2683,38.0488)(-27.8049,36.5854) \psline[linewidth=0.25pt](-26.3415,38.0488)(-24.8780,39.5122)(-26.3415,40.9756)(-27.8049,39.5122)(-26.3415,38.0488) \psline[linewidth=0.25pt](-24.8780,39.5122)(-23.4146,40.9756)(-24.8780,42.4390)(-26.3415,40.9756)(-24.8780,39.5122) \psline[linewidth=0.25pt](-23.4146,40.9756)(-21.9512,42.4390)(-23.4146,43.9024)(-24.8780,42.4390)(-23.4146,40.9756) \psline[linewidth=0.25pt](-21.9512,42.4390)(-20.4878,43.9024)(-21.9512,45.3659)(-23.4146,43.9024)(-21.9512,42.4390) \psline[linewidth=0.25pt](-33.6585,33.6585)(-32.1951,35.1220)(-33.6585,36.5854)(-35.1220,35.1220)(-33.6585,33.6585) \psline[linewidth=0.25pt](-32.1951,35.1220)(-30.7317,36.5854)(-32.1951,38.0488)(-33.6585,36.5854)(-32.1951,35.1220) \psline[linewidth=0.25pt](-30.7317,36.5854)(-29.2683,38.0488)(-30.7317,39.5122)(-32.1951,38.0488)(-30.7317,36.5854) \psline[linewidth=0.25pt](-29.2683,38.0488)(-27.8049,39.5122)(-29.2683,40.9756)(-30.7317,39.5122)(-29.2683,38.0488) \psline[linewidth=0.25pt](-27.8049,39.5122)(-26.3415,40.9756)(-27.8049,42.4390)(-29.2683,40.9756)(-27.8049,39.5122) \psline[linewidth=0.25pt](-26.3415,40.9756)(-24.8780,42.4390)(-26.3415,43.9024)(-27.8049,42.4390)(-26.3415,40.9756) \psline[linewidth=0.25pt](-24.8780,42.4390)(-23.4146,43.9024)(-24.8780,45.3659)(-26.3415,43.9024)(-24.8780,42.4390) \psline[linewidth=0.25pt](-35.1220,35.1220)(-33.6585,36.5854)(-35.1220,38.0488)(-36.5854,36.5854)(-35.1220,35.1220) \psline[linewidth=0.25pt](-33.6585,36.5854)(-32.1951,38.0488)(-33.6585,39.5122)(-35.1220,38.0488)(-33.6585,36.5854) \psline[linewidth=0.25pt](-32.1951,38.0488)(-30.7317,39.5122)(-32.1951,40.9756)(-33.6585,39.5122)(-32.1951,38.0488) \psline[linewidth=0.25pt](-30.7317,39.5122)(-29.2683,40.9756)(-30.7317,42.4390)(-32.1951,40.9756)(-30.7317,39.5122) \psline[linewidth=0.25pt](-29.2683,40.9756)(-27.8049,42.4390)(-29.2683,43.9024)(-30.7317,42.4390)(-29.2683,40.9756) \psline[linewidth=0.25pt](-27.8049,42.4390)(-26.3415,43.9024)(-27.8049,45.3659)(-29.2683,43.9024)(-27.8049,42.4390) \psline[linewidth=0.25pt](-26.3415,43.9024)(-24.8780,45.3659)(-26.3415,46.8293)(-27.8049,45.3659)(-26.3415,43.9024) \psline[linewidth=0.25pt](-36.5854,36.5854)(-35.1220,38.0488)(-36.5854,39.5122)(-38.0488,38.0488)(-36.5854,36.5854) \psline[linewidth=0.25pt](-35.1220,38.0488)(-33.6585,39.5122)(-35.1220,40.9756)(-36.5854,39.5122)(-35.1220,38.0488) \psline[linewidth=0.25pt](-33.6585,39.5122)(-32.1951,40.9756)(-33.6585,42.4390)(-35.1220,40.9756)(-33.6585,39.5122) \psline[linewidth=0.25pt](-32.1951,40.9756)(-30.7317,42.4390)(-32.1951,43.9024)(-33.6585,42.4390)(-32.1951,40.9756) \psline[linewidth=0.25pt](-30.7317,42.4390)(-29.2683,43.9024)(-30.7317,45.3659)(-32.1951,43.9024)(-30.7317,42.4390) \psline[linewidth=0.25pt](-38.0488,38.0488)(-36.5854,39.5122)(-38.0488,40.9756)(-39.5122,39.5122)(-38.0488,38.0488) \psline[linewidth=0.25pt](-36.5854,39.5122)(-35.1220,40.9756)(-36.5854,42.4390)(-38.0488,40.9756)(-36.5854,39.5122) \psline[linewidth=0.25pt](-35.1220,40.9756)(-33.6585,42.4390)(-35.1220,43.9024)(-36.5854,42.4390)(-35.1220,40.9756) \psline[linewidth=0.25pt](-33.6585,42.4390)(-32.1951,43.9024)(-33.6585,45.3659)(-35.1220,43.9024)(-33.6585,42.4390) \psline[linewidth=0.25pt](-32.1951,43.9024)(-30.7317,45.3659)(-32.1951,46.8293)(-33.6585,45.3659)(-32.1951,43.9024) \psline[linewidth=0.25pt](-39.5122,39.5122)(-38.0488,40.9756)(-39.5122,42.4390)(-40.9756,40.9756)(-39.5122,39.5122) \psline[linewidth=0.25pt](-38.0488,40.9756)(-36.5854,42.4390)(-38.0488,43.9024)(-39.5122,42.4390)(-38.0488,40.9756) \psline[linewidth=0.25pt](-36.5854,42.4390)(-35.1220,43.9024)(-36.5854,45.3659)(-38.0488,43.9024)(-36.5854,42.4390) \psline[linewidth=0.25pt](-35.1220,43.9024)(-33.6585,45.3659)(-35.1220,46.8293)(-36.5854,45.3659)(-35.1220,43.9024) \psline[linewidth=0.25pt](-40.9756,40.9756)(-39.5122,42.4390)(-40.9756,43.9024)(-42.4390,42.4390)(-40.9756,40.9756) \psline[linewidth=0.25pt](-39.5122,42.4390)(-38.0488,43.9024)(-39.5122,45.3659)(-40.9756,43.9024)(-39.5122,42.4390) \psline[linewidth=0.25pt](-38.0488,43.9024)(-36.5854,45.3659)(-38.0488,46.8293)(-39.5122,45.3659)(-38.0488,43.9024) \psline[linewidth=0.25pt](-36.5854,45.3659)(-35.1220,46.8293)(-36.5854,48.2927)(-38.0488,46.8293)(-36.5854,45.3659) \psline[linewidth=0.25pt](-42.4390,42.4390)(-40.9756,43.9024)(-42.4390,45.3659)(-43.9024,43.9024)(-42.4390,42.4390) \psline[linewidth=0.25pt](-40.9756,43.9024)(-39.5122,45.3659)(-40.9756,46.8293)(-42.4390,45.3659)(-40.9756,43.9024) \psline[linewidth=0.25pt](-39.5122,45.3659)(-38.0488,46.8293)(-39.5122,48.2927)(-40.9756,46.8293)(-39.5122,45.3659) \psline[linewidth=0.25pt](-38.0488,46.8293)(-36.5854,48.2927)(-38.0488,49.7561)(-39.5122,48.2927)(-38.0488,46.8293) \psline[linewidth=0.25pt](-43.9024,43.9024)(-42.4390,45.3659)(-43.9024,46.8293)(-45.3659,45.3659)(-43.9024,43.9024) \psline[linewidth=0.25pt](-42.4390,45.3659)(-40.9756,46.8293)(-42.4390,48.2927)(-43.9024,46.8293)(-42.4390,45.3659) \psline[linewidth=0.25pt](-40.9756,46.8293)(-39.5122,48.2927)(-40.9756,49.7561)(-42.4390,48.2927)(-40.9756,46.8293) \psline[linewidth=0.25pt](-39.5122,48.2927)(-38.0488,49.7561)(-39.5122,51.2195)(-40.9756,49.7561)(-39.5122,48.2927) \psline[linewidth=0.25pt](-45.3659,45.3659)(-43.9024,46.8293)(-45.3659,48.2927)(-46.8293,46.8293)(-45.3659,45.3659) \psline[linewidth=0.25pt](-43.9024,46.8293)(-42.4390,48.2927)(-43.9024,49.7561)(-45.3659,48.2927)(-43.9024,46.8293) \psline[linewidth=0.25pt](-42.4390,48.2927)(-40.9756,49.7561)(-42.4390,51.2195)(-43.9024,49.7561)(-42.4390,48.2927) \psline[linewidth=0.25pt](-40.9756,49.7561)(-39.5122,51.2195)(-40.9756,52.6829)(-42.4390,51.2195)(-40.9756,49.7561) \psline[linewidth=0.25pt](-46.8293,46.8293)(-45.3659,48.2927)(-46.8293,49.7561)(-48.2927,48.2927)(-46.8293,46.8293) \psline[linewidth=0.25pt](-45.3659,48.2927)(-43.9024,49.7561)(-45.3659,51.2195)(-46.8293,49.7561)(-45.3659,48.2927) \psline[linewidth=0.25pt](-43.9024,49.7561)(-42.4390,51.2195)(-43.9024,52.6829)(-45.3659,51.2195)(-43.9024,49.7561) \psline[linewidth=0.25pt](-48.2927,48.2927)(-46.8293,49.7561)(-48.2927,51.2195)(-49.7561,49.7561)(-48.2927,48.2927) \psline[linewidth=0.25pt](-46.8293,49.7561)(-45.3659,51.2195)(-46.8293,52.6829)(-48.2927,51.2195)(-46.8293,49.7561) \psline[linewidth=0.25pt](-49.7561,49.7561)(-48.2927,51.2195)(-49.7561,52.6829)(-51.2195,51.2195)(-49.7561,49.7561) \psline[linewidth=0.25pt](-51.2195,51.2195)(-49.7561,52.6829)(-51.2195,54.1463)(-52.6829,52.6829)(-51.2195,51.2195) \psline[linewidth=0.25pt](-52.6829,52.6829)(-51.2195,54.1463)(-52.6829,55.6098)(-54.1463,54.1463)(-52.6829,52.6829) \psline[linewidth=0.25pt](-54.1463,54.1463)(-52.6829,55.6098)(-54.1463,57.0732)(-55.6098,55.6098)(-54.1463,54.1463) \psline[linewidth=0.25pt](-55.6098,55.6098)(-54.1463,57.0732)(-55.6098,58.5366)(-57.0732,57.0732)(-55.6098,55.6098) \endpspicture}{
Diagramme de Young aléatoire tiré suivant la mesure de Plancherel pour $n=500$. La différence $\|\lambda^*-\Omega\|_\infty$ vaut ici environ $0.07$.\label{bigplancherel}}{Diagramme de Young aléatoire tiré suivant la mesure de Plancherel}

\section{Loi des grands nombres}\label{lskv}
Les contours du diagramme aléatoire présenté dans la figure \ref{bigplancherel} épousent une forme régulière $\Omega$, et ce phénomène est d'autant plus prononcé que $n$ est grand. Le théorème fondamental suivant précise cette affirmation. Si $\lambda$ est un diagramme de Young de taille $n$, nous noterons $\lambda^{*}$ le diagramme renormalisé $\lambda^{1/n}$ ; l'aire entre le diagramme continu $\lambda^{*}$ et la fonction valeur absolue est donc égale à $2$.

\begin{theorem}[Logan-Shepp-Kerov-Vershik, \cite{LS77,KV77}]\label{firstasymptoticplancherel}
Soit $\Omega$ le diagramme continu défini par
$$\Omega(s)=\begin{cases}
            \frac{2}{\pi}\, \left(s \arcsin(\frac{s}{2})+\sqrt{4-s^2} \right) &\text{si }|s| \leq 2,\\
|s| &\text{sinon.}
            \end{cases}
$$
Ce diagramme est la forme limite universelle des diagrammes de Young tirés suivant la mesure de Plancherel. Ainsi, pour tout $\eps>0$, la probabilité $M_{n}[\|\lambda^{*}-\Omega\|_{\infty} \geq \eps]$ tend vers $0$. 
\end{theorem}\bigskip

La preuve du théorème \ref{firstasymptoticplancherel} est en réalité assez aisée si l'on utilise les cumulants libres de diagrammes. En effet, on peut facilement calculer la fonction génératrice 
$$G_{\Omega}(z)=\frac{z-\sqrt{z^{2}-4}}{2}$$
de la forme limite $\Omega$, et son inverse formel au voisinage de l'infini est $R_{\Omega}(z)=z+1/z$. Comme la topologie de la convergence uniforme sur $\cym$ peut être contrôlée par les cumulants libres (voir la section \ref{markovkrein}), il suffit de montrer la convergence en probabilité des cumulants libres :
$$\forall j\geq 2,\,\,\,R_{j}(\lambda_{n}^{*}) \to \begin{cases} 1 &\text{si }j=2,\\
0 &\text{sinon}.
\end{cases}
$$
Or, sous les mesures de Plancherel $M_{n}$, les moments des caractères centraux peuvent être rendus explicites. Ainsi,
$$M_{n}[\varSigma_{\mu}]=n^{\downarrow |\mu|}\sum_{\lambda \in \ym_{n}}\frac{(\dim \lambda)^{2}}{n!}\,\chi^{\lambda}(\mu)=n^{\downarrow |\mu|} \, \mathrm{tr}_{\C \sym_{n}}(\sigma_{\mu1^{n-|\mu|}})=n^{\downarrow k}\,\,\mathbb{1}_{\mu=1^{k}}\,,\label{startidentityplancherel}$$
car la décomposition de la trace régulière sur la base des caractères irréductibles est donnée par le mesure de Plancherel, et parce que cette trace évaluée en un élément $\sigma$ du groupe vaut $0$ sauf si $\sigma=\id_{\lle 1,n\rre}$ est le neutre. Comme $(\varSigma_{\mu})_{\mu \in \ym}$ est une base de $\obs$, on en déduit que toute espérance d'observable $f$ sous la mesure de Plancherel $M_{n}$ est un polynôme en $n$. De plus, si $\mathrm{wt}(f)=k$, alors ce polynôme est de degré inférieur à $k/2$, car ceci est vrai pour les observables $\varSigma_{\mu}$. Pour $\mu=1^{k}$, il existe une observable $f$ de poids inférieur à $2k-1$ telle que 
$$(R_{2})^{k}=\varSigma_{1^{k}} - f.$$
Par conséquent, $M_{n}[(R_{2}(\lambda))^{k}]=n^{k}+O(n^{k-1/2})$, et 
$$M_{n}[(R_{2}(\lambda^{*}))^{k}]=1+O(n^{-1/2})\,.$$ 
Le cumulant $R_{2}(\lambda_{n}^{*})$ converge donc en moments vers la constante $1$, et en particulier il con\-ver\-ge en probabilité. D'autre part, le même argument pour $j \geq 3$ et $k \geq 1$ montre que $M_{n}[(R_{j}(\lambda))^{k}]=O(n^{(jk-1)/2})$ ; ainsi,
$$M_{n}[(R_{j \geq 3}(\lambda^{*}))^{k}]=O(n^{-1/2})\,.$$
Ceci implique la convergence en probabilité $R_{j \geq 3}(\lambda_{n}^{*}) \to 0$, d'où le théorème de Logan-Shepp-Kerov-Vershik. \bigskip
\bigskip

Ainsi, en tirant parti des propriétés des diverses bases de l'algèbre d'observables $\obs$, on peut démontrer facilement des résultats asymptotiques sur les partitions tirées suivant une mesure de Plancherel. Le point important de la preuve est le calcul de l'espérance des caractères centraux, et ce calcul peut être généralisé comme suit : si $M_{V}$ est la mesure de Plancherel d'un $\C \sym_{n}$-module $V$, alors
$$M_{V}[\varSigma_{\mu}]=n^{\downarrow |\mu|} \,\chi^{V}(\sigma_{\mu1^{n-|\mu|}}),$$
$\chi^{V}$ désignant le caractère normalisé de $V$. Ensuite, on utilise les graduations de l'algèbre d'observables pour remplacer $\varSigma_{\mu}$ par sa composante de plus haut poids, ou plus généralement par sa composante en plus haut degré vis-à-vis d'une filtration d'algèbre sur $\obs$. En particulier, nous prendrons dans la seconde partie du mémoire le degré canonique, qui est adapté à la renormalisation de certains diagrammes non équilibrés.\bigskip

En utilisant les estimations \emph{a priori} sur la distribution de la longueur d'un plus long sous-mot croissant dans une permutation aléatoire (voir par exemple \cite{Hamm72}), on peut rendre le théorème \ref{firstasymptoticplancherel} plus fort : ainsi, on a en réalité convergence en probabilité au sens de la topologie ultra-forte évoquée à la fin du chapitre \ref{tool}. Comme $-2$ et $2$ sont les bornes du support de $\sigma_{\Omega}$, ceci implique bien le résultat précédemment annoncé pour la longueur du plus long sous-mot croissant d'une permutation aléatoire.\bigskip

\section{Théorème central limite}\label{sniady}
La forme limite des diagrammes de Young sous la mesure de Plancherel ayant été déterminée, intéressons-nous maintenant aux fluctuations autour de cette forme limite. Le théorème central limite de Kerov est nettement plus complexe que la loi des grands nombres, et il admet plusieurs énoncés équivalents, le plus simple étant sous doute la concentration gaussienne des caractères centraux. On doit à P. \'Sniady (\emph{cf.} \cite{Sni06b}) un cadre général permettant de traiter ces phénomènes de concentration gaussienne, et la fin de ce chapitre est consacré à l'exposé de ces méthodes, avec en vue la démonstration du théorème suivant :
\begin{theorem}[Kerov, \cite{Ker93a}]\label{secondasymptoticplancherel}
Le vecteur des caractères centraux renormalisés 
$\left(\frac{\varSigma_{k}(\lambda_{n})}{\sqrt{k}\,n^{k/2}}\right)_{k \geq 2}$
con\-ver\-ge en lois fini-dimensionnelles vers un vecteur gaussien $(\xi_{k})_{k \geq 2}$, les variables $\xi_{k}$ étant indépendantes et de variance $1$. Soit $\Delta_{n}(s)=\sqrt{n}\,(\lambda_{n}^{*}(s)-\Omega(s))$ la déviation aléatoire renormalisée d'un diagramme de Young par rapport à la forme limite. Alors, au sens des distributions, $\Delta_{n}(s)$ converge en loi vers le processus gaussien généralisé
$$\Delta(s)=\Delta(2\cos\theta)=\frac{2}{\pi}\,\sum_{k=2}^{\infty}\frac{\xi_{k}}{\sqrt{k}}\,\sin(k\theta)\,.$$
\end{theorem}
\figcapt{\vspace{-5mm}\includegraphics[scale=0.5]{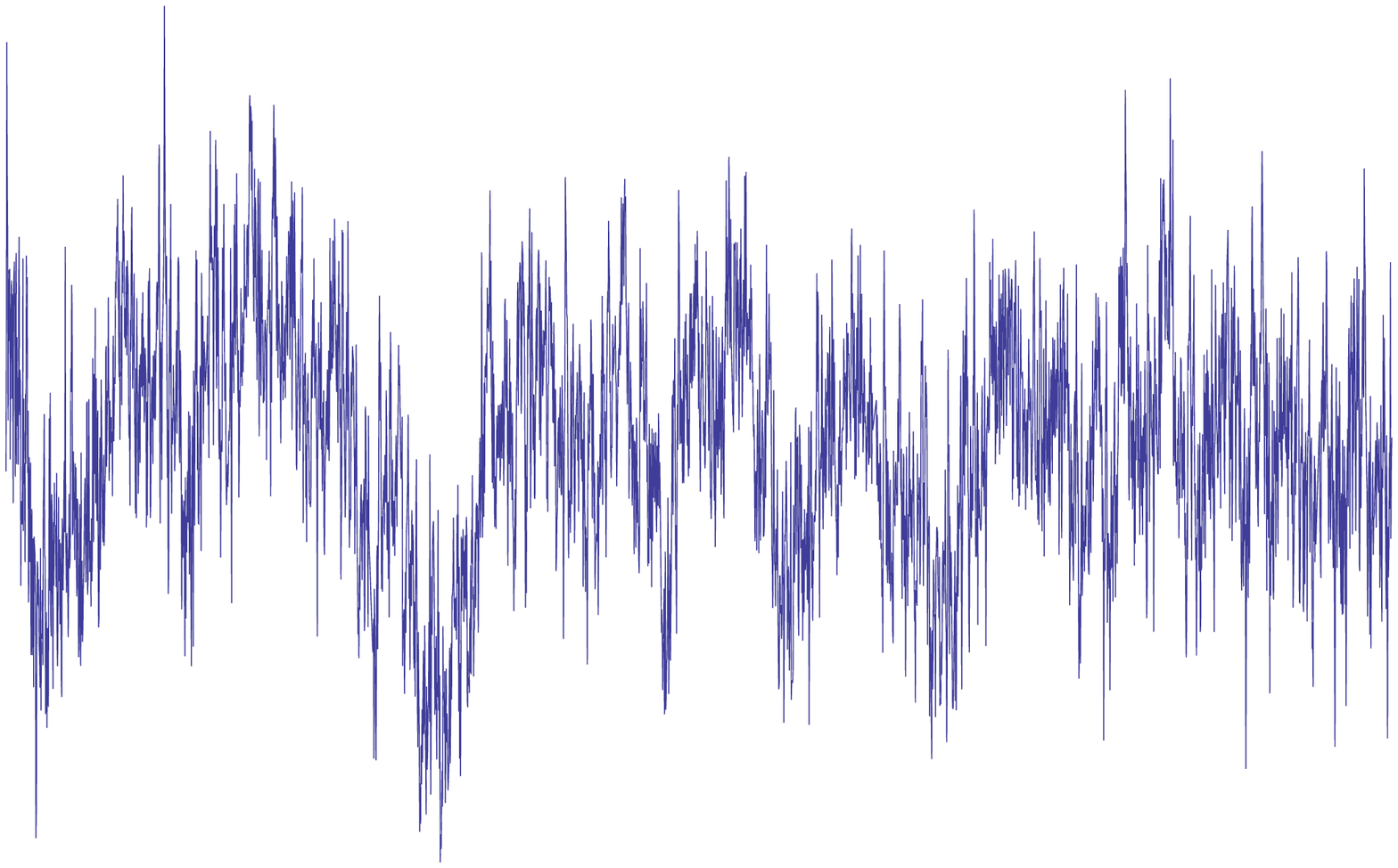}}{La déviation renormalisée d'un diagramme aléatoire $\lambda^{*}$ par rapport à la forme limite $\Omega$ converge vers un processus gaussien généralisé.\label{kerovprocess}}{Déviation renormalisée d'un diagramme de Young aléatoire}

\noindent Précisons quelque peu la seconde partie du théorème. La figure \ref{kerovprocess} représente une somme partielle (jusqu'à l'ordre $K=5000$) de la série mise en jeu dans le théorème de Kerov ; cette série est analogue à celle du \textbf{champ libre gaussien} (\cite{She07}), et elle ne converge pas ponctuellement. Néanmoins, pour toute fonction régulière $f \in \mathscr{C}^{\infty}([-2,2])$, l'intégrale 
$$\scal{\Delta}{f}=\int_{-2}^{2} f(s) \,\Delta(s)\,ds= \frac{2}{\pi}\,\sum_{k=2}^{\infty} \frac{\xi_{k}}{\sqrt{k}}\int_{-2}^{2}f(s) \sin(k\theta(s))\,ds$$
converge et est une variable gaussienne (centrée). Ces variables peuvent être explicitement calculées en décomposant $f$ sur la base des polynômes de Chebyshev de seconde espèce $u_{k}(s)=u_{k}(2\cos \theta)= \frac{\sin(k+1) \theta}{\sin \theta}$ :
$$f(s) = \sum_{k=0}^{\infty} f_{k}\,u_{k}(s)\quad \Rightarrow \quad \scal{\Delta}{f}=2\sum_{k=1}^{\infty} \frac{f_{k}\,\xi_{k+1}}{\sqrt{k+1}}\,.$$
Ainsi, la variance de $\scal{\Delta}{f}$ est la norme de Sobolev $\|f\|^{2}=4\sum_{k=1}^{\infty} |f_{k}|^{2}/(k+1)$, et la convergence en loi au sens des distributions\footnote{On peut aussi montrer que la déviation $\Delta_{n}/\sqrt{\log n}$ converge en lois fini-dimensionnelles vers un <<~vrai~>> processus gaussien, \emph{cf.} \cite{BS07}. } $\Delta_{n}\to \Delta$ signifie que $\scal{\Delta_{n}}{f}$ converge en loi vers cette variable gaussienne $\scal{\Delta}{f}$, et ce pour toute fonction $f \in \mathscr{C}^{\infty}([-2,2])$. \bigskip
\bigskip

La preuve de ce résultat se décompose en deux parties : on montre d'abord le caractère asymptotique gaussien de certaines observables de diagrammes sous la mesure de Plancherel, puis on en déduit le caractère asymptotique gaussien du processus $\Delta_{n}$ par des arguments analytiques. La seconde partie dépend essentiellement de la forme limite observée pour les diagrammes ; ainsi, on aura recours à des arguments sensiblement différents pour démontrer l'asymptotique gaussienne des $q$-mesures de Plancherel (section \ref{qgaussian}). En revanche, la première partie peut être traitée dans un cadre très général, \emph{cf.} \cite{Sni06a,Sni06b}. Nous exposons ci-après les résultats importants de cette théorie des \textbf{cumulants d'observables} due à \'Sniady. Si $X_{1},\ldots,X_{r}$ sont des variables aléatoires définies sur un même espace, leurs \textbf{cumulants joints} $k(X_{1},\ldots,X_{r})$ sont définis par la relation de récurrence :
$$\esper[X_{1}\cdots X_{r}]\,\,\,\,\,\,=\!\!\!\!\!\!\sum_{\substack{\pi \text{ partition }\\\text{d'ensemble de }\lle 1,r\rre}}k(X_{i \in \pi_{1}})\,k(X_{i \in \pi_{2}})\,\cdots\, k(X_{i \in \pi_{l}})\,.$$
Cette relation permet bien un calcul récursif des cumulants ; ainsi, $\esper[X]=k(X)$, $\esper[XY]=k(X,Y)+k(X)k(Y)$ et $\esper[XYZ]=k(X,Y,Z)+k(X,Y)k(Z)+k(X,Z)k(Y)+k(Y,Z)k(X)+k(X)k(Y)k(Z)$, donc
\begin{align*}&k(X)=\esper[X] \qquad;\qquad k(X,Y)=\esper[XY]-\esper[X]\esper[Y]=\mathrm{cov}(X,Y)\qquad;\\
&k(X,Y,Z)=\esper[XYZ]-\esper[XY]\esper[Z]-\esper[XZ]\esper[Y]-\esper[YZ]\esper[X]+2\,\esper[X]\esper[Y]\esper[Z]\,.
\end{align*}
On renvoie à \cite{Matt99,LS59,Bri69} pour des précisions sur la notion de cumulant. Alternativement, les cumulants peuvent être définis par la formule
$$k(X_{1},\ldots,X_{r})=\left.\frac{\partial^{r}}{\partial t_{1}\cdots \partial t_{r}}\right|_{t=0} \log \esper[\exp(t_{1}X_{1}+\cdots + t_{r}X_{r})]\,.$$
En particulier, les vecteurs gaussiens sont caractérisés par l'annulation des cumulants joints d'ordre supérieur à $3$. C'est ce critère que nous employerons pour démontrer les concentrations gaussiennes.\bigskip
\bigskip

Ceci étant, si $V$ est un module (réductible) sur $\C \sym_{n}$, le caractère normalisé $\chi^{V}$ fournit une structure d'espace de probabilité non commutatif sur l'algèbre $\C\sym_{n}$. On notera $\esper[X]=\chi^{V}(X)$ l'espérance des variables aléatoires non commutatives de l'algèbre. Si l'on se restreint au centre $Z(\C\sym_{n})$ de l'algèbre, les variables aléatoires considérées commutent, et les produits mis en jeu dans la relation de récurrence définissant les cumulants joints ne font plus ambiguïté. On peut donc envisager les cumulants de variables centrales de $Z(\C\sym_{n})$ --- plus généralement d'ailleurs, on peut envisager les cumulants de toute famille de variables commutantes. D'autre part, en utilisant la projection $\pi_{n}\circ \phi_{n} : \blg_{\infty} \to \blg_{n} \to \C \sym_{n}$, on peut relever $\esper$ à l'algèbre des permutations partielles $\blg_{\infty}$, qu'on considérera de même comme un espace de probabilité non commutatif. Ainsi, un module $V$ étant fixé, on peut considérer l'espérance d'un caractère central $\varSigma_{\mu}$ (vu comme élément de $\blg_{\infty}$), et également des cumulants joints
$$k(\varSigma_{\mu_{1}},\varSigma_{\mu_{2}},\ldots,\varSigma_{\mu_{r}})\,,$$
car les $\varSigma_{\mu}$ commutent dans $\alg_{\infty} \subset \blg_{\infty}$. Finalement, les caractères centraux fournissent un isomorphisme $\alg_{\infty} \simeq \obs$, donc on peut considérer les observables de diagrammes comme des variables aléatoires, et calculer leurs cumulants. En termes probabilistes, ceci revient à tirer un diagramme $\lambda$ suivant la mesure de Plancherel de $V$, et à évaluer les observables sur ces diagrammes aléatoires.\bigskip
\bigskip

Si $(\sigma,S)$ et $(\tau,T)$ sont deux permutations partielles, leur produit a été défini dans la section \ref{centralcharacter} par $(\sigma,S)(\tau,T)=(\sigma\tau,S\cup T)$. Mais on peut aussi définir dans l'algèbre $\blg_{\infty}$ un \textbf{produit disjoint}
$$(\sigma,S) \bullet (\tau,T)=\begin{cases} (\sigma \tau,S\sqcup T) &\text{si }S\cap T=\emptyset,\\
0&\text{sinon}.
\end{cases}$$
Pour ce produit, $\varSigma_{\mu_{1}}\bullet \varSigma_{\mu_{2}}=\varSigma_{\mu_{1}\sqcup \mu_{2}}$. Nous noterons $\obs^{\bullet}$ l'algèbre d'observables muni du produit disjoint, et $k^\bullet$ les cumulants d'observables correspondant à ce produit ; ainsi,
$$\esper[X_{1}\bullet X_{2}\bullet \cdots \bullet X_{r}]=\sum_{\pi \in \mathfrak{Q}(\lle 1,r\rre)} \prod_{\pi_{j} \in \pi} k^{\bullet}(X_{i \in \pi_{j}})\,,$$
où $\mathfrak{Q}(\lle 1,r\rre)$ est l'ensemble des partitions d'ensembles de $\lle 1,n\rre$.
Une technique de conditionnement (\emph{cf.} \cite{Bri69}) permet de relier les cumulants standards et les cumulants disjoints. Ainsi, si $\esper^{\mathrm{id}}$ désigne l'application identité entre $\obs$ et $\obs^{\bullet}$, alors
$$k(X_{1},\ldots,X_{r})=\sum_{\pi \in \mathfrak{Q}(\lle 1,r\rre)}k^{\bullet}(k^{\mathrm{id}}(X_{i \in \pi_{1}}),\ldots,k^{\mathrm{id}}(X_{i \in \pi_{l}}))\,.$$
En utilisant ces relations, la filtration du poids sur l'algèbre d'observables,  et une théorie <<géométrique~>> des caractères centraux (\emph{cf.} \cite{Sni06a}), P. \'Sniady a établi dans \cite{Sni06b} l'équivalence entre plusieurs propriétés de concentration gaussienne d'observables. Si $(\rho_{n},V_{n})_{n \in \N}$ est une famille de représentations des groupes symétriques $\sym_{n}$, on dit que la famille a la propriété de \textbf{factorisation asymptotique des caractères} si l'une des assertions suivantes est satisfaite :
\begin{enumerate}[(i)]
\item pour tous cycles disjoints $\sigma_{1}, \ldots,\sigma_{r}$ de longueurs $l_{1},\ldots,l_{r}$, $$k_{n}(\sigma_{1},\ldots,\sigma_{r})\,n^{\frac{l_{1}+\cdots+l_{r}+r-2}{2}}=O(1) \,\,;$$
\item pour tous entiers positifs $l_{1},\ldots,l_{r}$,
$$k_{n}(\varSigma_{l_{1}},\ldots,\varSigma_{l_{r}})\,n^{-\frac{l_{1}+\cdots+l_{r}-r+2}{2}}=O(1)\quad\text{ou}\quad k_{n}^{\bullet}(\varSigma_{l_{1}},\ldots,\varSigma_{l_{r}})\,n^{-\frac{l_{1}+\cdots+l_{r}-r+2}{2}}=O(1) \,\,;$$
\item pour tous entiers $l_{1},\ldots,l_{r} \geq 2$,
$$ k_{n}(R_{l_{1}},\ldots,R_{l_{r}})\,n^{-\frac{l_{1}+\cdots+l_{r}-2r+2}{2}}=O(1)\,.$$
\end{enumerate}
\begin{proposition}[Factorisation asymptotique et concentration gaussienne, \cite{Sni06b}]\label{factorasymptotic}
Les assertions (i)-(iii) sont équivalentes. Lorsqu'elles sont vérifiées, si pour $r\in \{1,2\}$ les $r$-cumulants joints d'un type (i)-(iii) ont une limite, alors c'est le cas pour tous les types de cumulants, et :
\begin{align*} c_{l+1}&=\lim_{n \to \infty} \esper_{n}[\sigma_{l}]\,n^{\frac{l-1}{2}}=\lim_{n\to \infty} \esper_{n}[\varSigma_{l}]\,n^{-\frac{l+1}{2}}=\lim_{n\to \infty} \esper_{n}[R_{l+1}]\,n^{-\frac{l+1}{2}}\\
v_{l_{1}+1,l_{2}+1}&=\lim_{n\to \infty} k_{n}(R_{l_{1}+1},R_{l_{2}+1})\,n^{-\frac{l_{1}+l_{2}}{2}}=
\lim_{n\to \infty} k_{n}(\varSigma_{l_{1}},\varSigma_{l_{2}})\,n^{-\frac{l_{1}+l_{2}}{2}}\\
&=\lim_{n\to \infty} k_{n}^{\bullet}(\varSigma_{l_{1}},\varSigma_{l_{2}})\,n^{-\frac{l_{1}+l_{2}}{2}}+\sum_{\substack{l_{1}=a_{1}+\cdots+a_{r}\\ l_{2}=b_{1}+\cdots + b_{r}}} \frac{l_{1}l_{2}}{r}\,c_{a_{1}+b_{1}}\cdots c_{a_{r}+b_{r}}\\
&=\lim_{n\to \infty} k_{n}(\sigma_{l_{1}},\sigma_{l_{2}})\,n^{\frac{l_{1}+l_{2}}{2}}-l_{1}l_{2}\,c_{l_{1}+1}\,c_{l_{2}+1}+\sum_{\substack{l_{1}=a_{1}+\cdots+a_{r}\\ l_{2}=b_{1}+\cdots + b_{r}}} \frac{l_{1}l_{2}}{r}\,c_{a_{1}+b_{1}}\cdots c_{a_{r}+b_{r}}\,.
\end{align*}
\end{proposition}
Dans les conditions de la proposition \ref{factorasymptotic}, si $\lambda$ est un diagramme aléatoire tiré suivant la mesure de Plancherel de $V_{n}$, alors les vecteurs 
$$r_{l}=n^{-\frac{l-1}{2}}\,(R_{l}(\lambda)-\esper_{n}[R_{l}]),\quad l\geq 2 \qquad;\qquad s_{l}=n^{\frac{l}{2}}(\chi^{\lambda}(\sigma_{l})-\esper_{n}[\sigma_{l}]),\quad l\geq 1$$
convergent en loi vers des vecteurs gaussiens, \emph{cf.} \cite[corollaire 4]{Sni06b}. En particulier, on déduit immédiatement de la proposition \ref{factorasymptotic} la première partie du théorème \ref{secondasymptoticplancherel} : en effet, on a dans ce cas
$$k_{n}(\sigma_{l_{1}},\ldots,\sigma_{l_{r}})=\begin{cases} 1 & \text{si }n=1 \text{ et } l_{1}=1,\\
0&\text{sinon},
\end{cases}$$
donc les caractères centraux renormalisés $\varSigma_{k}/n^{k/2}$ convergent vers des variables gaussiennes indépendantes. De plus, la formule de la proposition \ref{factorasymptotic} permet de calculer $$\lim_{n \to \infty} \esper[\varSigma_{k}^{2}]/n^{k}=k,$$ car dans la différence de $k_{n}(\varSigma_{k},\varSigma_{k})$ et de $k_{n}^{\bullet}(\varSigma_{k},\varSigma_{k})$, le seul terme non nul est celui correspondant aux compositions $k=1+\cdots+1$ et $k=1+\cdots+1$, et il vaut $k$. Ainsi, $\varSigma_{k}/n^{k/2} \to \mathcal{N}(0,k)$ pour $k \geq 2$, et on a du même coup l'asymptotique gaussienne des cumulants libres $R_{k+1}(\lambda)$ sous la mesure $M_{n}$,  et avec les mêmes covariances. \bigskip
\bigskip

Pour démontrer la seconde partie du théorème \ref{secondasymptoticplancherel} à partir du caractère asymptotique gaussien des cumulants libres et des caractères centraux, on exprime les moments de la déviation renormalisée en termes d'observables de diagrammes :
\begin{align*}
\sqrt{n}\int_{\R} s^{k}\, \frac{\lambda^{*}(s)-\Omega(s)}{2} \,ds &= \sqrt{n}\int_{\R}s^{k}\,\frac{\lambda^{*}(s)-|s|}{2}-\sqrt{n}\int_{\R} s^{k}\,\frac{\Omega(s)-|s|}{2}\,ds\\
&=\frac{1}{(k+1)(k+2)}\left(\frac{\tilp_{k+2}(\lambda)}{n^{\frac{k+1}{2}}} - \sqrt{n} \,\tilp_{k+2}(\Omega)\right).
\end{align*}
Dans ce contexte, il est utile d'étendre l'algèbre d'observables $\obs$ en rajoutant les puissances demi-entières (positives ou négatives) de $n=\varSigma_{1}(\lambda)$. On note $\obs^{+}$ cette nouvelle algèbre, et le \textbf{degré de Kerov} d'une observable $f=\varSigma_{\mu} \,(\varSigma_{1})^{k/2}$ avec $\mu$ sans part de taille égale à $1$ et $k \in \Z$ est $\deg_{\mathrm{K}}(f)= |\mu|+k$. En particulier, le degré de Kerov d'un caractère central $\varSigma_{\nu}$ avec $\nu$ partition quelconque est $|\nu|+m_{1}(\nu)$. On peut montrer que ceci définit bien une filtration d'algèbre sur $\obs^{+}$ ; notons que la preuve de la loi des grands nombres \ref{firstasymptoticplancherel} aurait pu être réalisée à l'aide de cette graduation (au lieu du poids). Ceci étant, si\label{bizarobs}
$$\widetilde{q}_{k}=\begin{cases} 
\frac{\tilp_{k+1}-\binom{2m}{m}(\varSigma_{1})^{m}}{(k+1)\,(\varSigma_{1})^{k/2}} &\text{si }k=2m-1\text{ est impair},\\
\frac{\tilp_{k+1}}{(k+1)\,(\varSigma_{1})^{k/2}}&\text{si }k=2m\text{ est pair}, \end{cases}$$
alors le $k$-ième moment de la déviation renormalisée est simplement $\widetilde{q}_{k+1}/(k+1)$. On écrit alors $\varSigma_{k}/(\varSigma_{1})^{k/2}$ comme combinaison linéaire d'observables $\widetilde{q}_{j}$, plus des observables dont le degré de Kerov\footnote{Incidemment, on peut aussi effectuer les raisonnements en gardant sur $\obs$ la filtration du poids ; nous expliquerons ceci dans la section \ref{schurplus}. Le degré de Kerov permet néanmoins de s'affranchir de difficultés calculatoires, et il sera particulièrement adapté à l'étude asymptotique des mesures de Gelfand, voir le chapitre \ref{gelfandmeasure}.} est strictement négatif, voir \cite[proposition 7.4]{IO02}. Ceci permet de montrer que les observables
$$\frac{\scal{\Delta_{n}}{u_{k}}}{2}=\frac{\sqrt{n}}{2}\int_{\R} u_{k}(s) \, (\lambda^{*}(s)-\Omega(s))\,ds$$
convergent pour $k\geq 1$ vers des gaussiennes indépendantes de variances $1/(k+1)$. La déviation gaussienne de la forme des diagrammes au sens du théorème \ref{firstasymptoticplancherel} s'en déduit immédiatement, car les polynômes de Chebyshev engendrent un sous-espace dense de $\mathscr{C}^{\infty}([-2,2])$. Nous donnerons plus de détail sur cette preuve dans la section \ref{schurplus} ; en effet, nous y généraliserons le théorème central limite de Kerov au cas des mesures de Schur-Weyl de paramètre $\alpha\geq 1/2$. Pour conclure, mentionnons la coïncidence suivante : les polynômes de Chebyshev de seconde espèce $u_{k}(s)$, qui fournissent une base orthogonale de l'espace gaussien $\leb^{2}(\Delta)$ sous-jacent au processus limite des déviations renormalisées, forment aussi une base orthonormale de l'espace de Hilbert
$$\leb^{2}\left(\mathbb{1}_{[-2,2]}(s)\,\frac{\sqrt{4-s^{2}}}{2\pi}\,ds\right)$$
associé à la mesure de transition du diagramme continu $\Omega$. Autrement dit, ce sont les polynômes orthogonaux pour la mesure de probabilité $\mu_{\Omega}$. Malheureusement, il semble que ce fait soit juste une coïncidence : ainsi, dans le cas des mesures de Schur-Weyl de paramètres $(\alpha=1/2,c\geq 0)$, nous verrons que les bases <<~canoniques~>> des espaces de Hilbert $\leb^{2}(\Delta_{c})$ et $\leb^{2}(\mu_{\Omega_{c}})$ sont distinctes.

\chapter{Groupe symétrique infini et processus ponctuels déterminantaux}\label{determinantal}

Ce chapitre est consacré à deux thèmes qui ne joueront pas un rôle crucial dans l'exposé des résultats que nous avons obtenus, mais qui permettront une compréhension plus subtile des dits résultats : la théorie des représentations du \textbf{groupe symétrique infini} (\S\ref{syminfinite}), et la théorie des \textbf{mesures de Schur} (\S\ref{schurmeasure}). Si $n<N$ sont deux entiers, on peut plonger le groupe symétrique $\sym_{n}$ dans $\sym_{N}$ en le faisant agir sur les $n$ premiers entiers de $\lle 1,N\rre$. On dispose ainsi d'une famille dirigée de morphismes de groupes $(i_{n, N})_{n \leq N}$, et la limite inductive 
$$\sym_{\infty}=\lim_{n \to \infty} \sym_{n}$$ 
correspondant à cette famille consiste en les permutations de $\N^{*}$ qui déplacent un nombre fini d'entiers. Les caractères irréductibles de ce groupe infini ont été déterminés dans les années 60 (\emph{cf.} \cite{Tho64}), et cette théorie a été précisée par les travaux de Kerov, Olshanski et Vershik (\cite{KV81,Ols90,KOV04}) et par la thèse d'Okounkov (\cite{Oko97}). La première partie du chapitre est consacrée à un exposé succinct de cette théorie. L'analyse harmonique sur le groupe $\sym_{\infty}$ conduit à l'étude des $z$-mesures, qui sont un particulier de mesures de Schur (\S\ref{schurmeasure}). Ces mesures de probabilité portent sur l'ensemble des diagrammes de Young, et on peut leur associer des processus ponctuels à corrélations déterminantales (\cite{Oko03a,Oko03b,OR03,Bor10}). De plus, la plupart des mesures issues de la théorie des représentations des groupes symétriques ou de leurs algèbres d'Hecke sont des mesures de Schur. Ainsi, les résultats de la section \ref{schurmeasure} s'appliqueront aux mesures de Plancherel des groupes symétriques (chapitre \ref{matrix}), aux mesures de Plancherel des groupes linéaires finis (chapitre \ref{general}), aux mesures des algèbres d'Iwahori-Hecke (chapitre \ref{iwahori}) et aux mesures de Schur-Weyl (chapitre \ref{schurweylmeasure}).\bigskip

\section{Représentations du groupe symétrique infini}\label{syminfinite}

Comme $G=\sym_{\infty}$ est un groupe infini, sa théorie des représentations doit prendre en compte les représentations unitaires $U:G \to \mathscr{U}(H)$ sur des espaces de Hilbert de dimension infinie. Fixons plus généralement un groupe $G$ infini, séparable et localement compact. L'\textbf{algèbre de von Neumann} $\mathscr{M}(G,U)$ associée à une représentation unitaire de $G$ est la clôture pour l'une des topologies faibles\footnote{La clôture d'une sous-$*$-algèbre de $\mathscr{L}(H)$ est la même pour toutes les topologies suivantes : topologie faible d'opérateur, topologie ultra-faible, topologie ultra-forte, topologie forte d'opérateur, topologie $*$-forte d'opérateur et topologie $*$-ultra-forte. On renvoie à \cite{Kir94} pour une présentation de la théorie des représentations des groupes infinis, et au monumental \cite{Tak79} pour la théorie générale des algèbres de von Neumann.} de $\mathscr{L}(H)$ de l'algèbre engendrée par les $U(g)$, $g\in G$. On dit que la représentation est irréductible si $H$ n'admet pas de sous-espace fermé $H'$ non trivial et stable par $G$ ; et que c'est une représentation facteur si $\mathscr{M}(G,U)$ est un facteur, c'est-à-dire une algèbre de von Neumann de centre trivial. Toute algèbre de von Neumann séparable s'écrit de manière unique comme intégrale directe de facteurs, et d'autre part, il est bien connu que les facteurs de von Neumann sont classifiés en trois types I, II et III. On peut montrer que :\vspace{2mm}
\begin{enumerate}
\item Une représentation facteur de type I s'écrit de manière essentiellement unique comme intégrale directe de représentations irréductibles.\vspace{2mm}
\item \`A l'inverse, une représentation facteur de type II ou III admet deux décompositions en intégrales directes de représentations irréductibles telles qu'aucune composante de la première décomposition n'apparaisse dans la seconde.\vspace{2mm}
\end{enumerate}
Le groupe $G$ est dit \textbf{modéré} s'il n'admet que des représentations facteurs de type I ; dans le cas contraire, on parle de groupe \textbf{sauvage}. Un résultat majeur de la théorie d'Harish-Chandra est le caractère modéré des groupes de Lie semi-simples, \emph{cf.} \cite{Var89}  --- ce résultat est également valable pour les groupes semi-simples $p$-adiques (Bernstein) et les groupes algébriques réels (Dixmier). Au contraire, la plupart des groupes infinis discrets sont sauvages. Ainsi, un groupe infini discret est modéré si et seulement s'il est virtuellement abélien, et si toutes ses classes de conjugaison sauf $\{e_{G}\}$ sont infinies, alors :\vspace{2mm}
\begin{enumerate}
\item La représentation régulière gauche de $G$ sur $\ell^{2}(G)$ est un facteur de type $\text{II}_{1}$.\vspace{2mm}
\item La représentation régulière bilatère de $G \times G$ sur $\ell^{2}(G)$ est irréductible.\vspace{2mm}
\end{enumerate}
En particulier, $\sym_{\infty}$ est un groupe sauvage, et la classification de toutes les représentations de $\sym_{\infty}$ n'est pas un objectif raisonnable, puisqu'on ne dispose pas de décomposition unique en représentations irréductibles. \bigskip
\bigskip

Un objectif plus raisonnable est la classification de tous les \textbf{facteurs finis} de $\sym_{\infty}$, c'est-à-dire les représentations unitaires de $\sym_{\infty}$ dont l'algèbre de von Neumann est un facteur qui admet une trace fidèle finie, \emph{i.e.}, telle que $\tau(1)<\infty$. \`A renormalisation près, on peut supposer $\tau(1)=1$, et le facteur fini est alors entièrement déterminé par l'application $\tau \circ U : \sym_{\infty} \to \C$. Cette application est ce qu'on appelle un \textbf{caractère} de $G=\sym_{\infty}$, c'est-à-dire une application $\chi : G \to \C$ qui est normalisée, centrale et définie positive :
\begin{align*}
&\chi(1)=1 \qquad;\qquad \forall g,h,\,\,\,\chi(hgh^{-1})=\chi(g)\qquad;\\
 &\forall (g_{1},\ldots,g_{n}),\,\,\,\,(\chi(g_{i}g_{j}^{-1}))_{i,j} \text{ est hermitienne et définie positive}.\end{align*}
On parle de \textbf{caractère virtuel} si l'on ôte la condition de positivité ; d'autre part, un caractère est dit \textbf{irréductible} s'il n'est pas barycentre non trivial de caractères. Si $(U,\mathscr{L}(H),\tau)$ est un facteur fini de $\sym_{\infty}$, alors $\chi=\tau\circ U$ est un caractère irréductible, et la classification des facteurs finis de $\sym_{\infty}$ se ramène donc à celle des caractères irréductibles.\bigskip
\bigskip

Ce problème a été résolu par E. Thoma (voir \cite{Tho64}), qui a montré qu'on pouvait paramétrer l'ensemble des caractères irréductibles par un simplexe de dimension infinie. S. Kerov et A. Vershik ont ensuite retrouvé la théorie de Thoma en étudiant l'asymptotique des caractères des groupes finis $\sym_{n}$, \emph{cf.} \cite{KV81}. Fixons un caractère virtuel  $\chi$  de $\sym_{\infty}$. Pour tout entier $n$, la restriction $\chi_{n}$ de $\chi$ au groupe $\sym_{n}$ est une fonction centrale vérifiant $\chi(\id)=1$. Cette restriction est donc un barycentre des caractères irréductibles $\chi^{\lambda}$ :
$$\forall n,\,\,\,\chi_{n}=\sum_{\lambda \in \ym_{n}} M_{\chi,n}(\lambda)\,\chi^{\lambda}\quad\text{avec }\sum_{\lambda \in \ym_{n}}M_{\chi,n}(\lambda)=1\,.$$
De plus, comme toutes les restrictions proviennent d'un même caractère de $\sym_{\infty}$, les $\chi_{n}$ vérifient la condition de compatibilité $\mathrm{Res}_{\sym_{n}}^{\sym_{n+1}}(\chi_{n+1})=\chi_{n}$, d'où :
$$\forall \lambda\in \ym,\,\,\, M_{\chi,n}(\lambda)=\sum_{\lambda \nearrow \Lambda} \frac{\dim \lambda}{\dim \Lambda}\,M_{\chi,n+1}(\Lambda)\,.$$
On dit que les mesures $(M_{\chi,n})_{n \in \N}$ forment un \textbf{système cohérent de mesures de probabilité} sur le graphe de Young $\ym=\bigsqcup_{n \in \N} \ym_{n}$. En d'autres termes, la fonction $\phi_{\chi}(\lambda)=M_{\chi}(\lambda)/\dim \lambda$ est \textbf{harmonique} sur le graphe de Young de la figure \ref{younggraph}, c'est-à-dire que 
$$\forall \lambda\in \ym,\,\,\,\phi_{\chi}(\lambda)=\sum_{\lambda \nearrow \Lambda}\phi_{\chi}(\Lambda)\,.$$
Réciproquement, toute fonction harmonique sur $\ym$ donne une famille $(\chi_{n})_{n \in \N}$ de caractères compatibles entre eux, donc un caractère virtuel de $\sym_{\infty}$. De plus, le caractère est positif si et seulement si les poids $M_{\chi}(\lambda)$ sont tous positifs, d'où une bijection :
$$\big\{\text{caractères (positifs) de }\sym_{\infty} \big\} \leftrightarrow \bigg\{\substack{\text{fonctions harmoniques positives sur} \\ \text{le graphe }\ym\text{, et telles que }f(\emptyset)=1}\bigg\}.$$
Notons $\mathcal{H}^{+}(\ym)$ cet ensemble de fonctions harmoniques ; il est convexe et compact pour la topologie de la convergence simple. Par le théorème de représentation intégrale de Choquet, toute fonction de $\mathcal{H}^{+}(\ym)$ s'écrit comme barycentre intégral des points extrémaux du convexe ; nous noterons $\Omega$ leur ensemble. Cet ensemble peut être vu comme la \textbf{frontière de Martin} du graphe de Young : ainsi, il existe une compactification $X=\ym \sqcup \partial\ym \hookrightarrow \mathcal{H}^{+}(\ym)$ du graphe de Young à base de noyaux de Martin (\emph{cf.} \cite{KOO97}), et cette compactification est telle que $\partial\ym=\Omega$. \bigskip\bigskip

\begin{theorem}[Simplexe de Thoma, \cite{Tho64,KV81}]\label{boundarymartin}
L'ensemble $\Omega$ s'identifie au simplexe de dimension infinie constitué des couples
$$\omega =\big(\alpha=(\alpha_{1}\geq \cdots \geq \alpha_{r} \geq  \cdots \geq 0)\,,\,\beta=(\beta_{1}\geq \cdots \geq \beta_{s} \geq  \cdots \geq 0)\big)$$
de suites positives décroissantes telles que  $\gamma=1-\sum_{i=1}^{\infty}(\alpha_{i}+\beta_{i})\geq 0$. Si $\phi$ est une fonction harmonique positive sur $\ym$, elle s'écrit sous forme intégrale $$\phi(\lambda)=\int_{\Omega} s_{\lambda}(\omega)\,m_{\phi}(d\omega)\,,$$
avec $m_{\phi}$ mesure de probabilité sur $\Omega$, et $s_{\lambda}(\omega)$ fonction de Schur généralisée définie par la spécialisation de $\Lambda$ :
$$p_{1}(\omega)=1 \qquad;\qquad p_{n \geq 2}(\omega)=p_{n}(\alpha-(-\beta))=\sum_{i \geq 1} (\alpha_{i})^{n} + (-1)^{n-1}\sum_{i\geq 1}(\beta_{i})^{n}$$
\end{theorem}
\noindent On renvoie à \cite[théorème 9.6.1]{KOV04} et aux articles originaux \cite{Tho64,KV81} pour des détails sur la preuve de ce résultat. Dans le contexte précédent, le caractère virtuel $\chi$ a donc pour restrictions :
\begin{align*}\chi_{n}(\sigma_{\mu} \in \sym_{n})&= \sum_{\lambda \in \ym_{n}} \phi_{\chi}(\lambda)\,(\dim \lambda)\,\chi^{\lambda}(\mu) = \sum_{\lambda \in \ym_{n}} \left(\int_{\Omega}s_{\lambda}(\omega)\,\varsigma^{\lambda}(\mu)\,\mu_{\chi}(d\omega)\right)\\
&= \int_{\Omega}\left(\sum_{\lambda \in \ym_{n}} s_{\lambda}(\omega)\,\varsigma^{\lambda}(\mu)\right)m_{\chi}(d\omega)=\int_{\Omega} p_{\mu}(\omega)\,m_{\chi}(d\omega)\,.\end{align*}
En particulier, les caractères irréductibles, qui correspondent aux points extrémaux du convexe compact $\mathscr{P}(\Omega)=\mathcal{H}^{+}(\ym)$, sont paramétrés par le simplexe de Thoma $\Omega$, et s'écrivent
$$\chi^{\omega}(\mu)=p_{\mu}(\omega) = \prod_{k \geq 2} \left(\sum_{i\geq 1}(\alpha_{i})^{k}+(-1)^{k-1}\sum_{i \geq 1}(\beta_{i})^{k}\right)^{m_{k}(\mu)}.$$
Ils sont donc multiplicatifs, c'est-à-dire que si $\sigma$ et $\tau$ sont deux permutations à supports disjoints, alors $\chi^{\omega}(\sigma\tau)=\chi^{\omega}(\sigma)\,\chi^{\omega}(\tau)$. Cette propriété est à rapprocher de la proposition \ref{factorasymptotic}. 
\bigskip
\bigskip

Le théorème de Thoma est en réalité très naturel si l'on souhaite réaliser $\Omega$ comme frontière <<~géodésique~>> de $\ym$. Pour réaliser cette compactification alternative, on plonge $\ym$ dans $\Omega \times [0,1]$ par 
$$\omega_{n} : \lambda\in \ym_{n} \mapsto \left\{ \bigg(\frac{a_{1}}{n},\ldots,\frac{a_{d}}{n},0,\ldots\,\,\, ;\,\,\, \frac{b_{1}}{n},\ldots,\frac{b_{d}}{n},0,\ldots \bigg),\frac{1}{n} \right\}$$
où les $a_{i}$ et les $b_{i}$ sont les coordonnées de Frobenius de $\lambda$. Alors, si $\omega$ est un point du simplexe de Thoma et si $(\lambda^{(n)})_{n \in \N}$ est une suite de partitions telle que $\omega_{n}(\lambda^{(n)})\to (\omega,0)$ dans $\Omega \times [0,1]$, 
$$\lim_{n \to \infty} \chi^{\lambda^{(n)}}(\sigma) =\chi^{\omega}(\sigma)$$
pour toute permutation $\sigma$. Ainsi :
\begin{theorem}[Kerov-Vershik, \cite{KV81,KOV04}]\label{boundarygeodesic}
Les coordonnées d'un point $\omega \in \Omega$ correspondent aux fréquences limites des coordonnées de Frobenius de partitions $\lambda^{(n)}$ telles que $\chi^{\lambda^{(n)}}$ approche $\chi^{\omega}$. De plus, si $m$ est une mesure de probabilité sur $\Omega$ correpondant au caractère (positif) $\chi$, et si $M_{\chi}=(M_{\chi,n})_{n \in \N}$ est le système cohérent de mesures de probabilité associé à $\chi$, alors 
$$\lim_{n \to \infty} (\omega_{n})_{\star}M_{\chi,n}=m$$
au sens de la convergence en loi, étant entendu que l'on identifie $\Omega$ et $\Omega \times \{0\}$. 
\end{theorem}\bigskip

En d'autres termes, compte tenu des propositions \ref{boundarymartin} et \ref{boundarygeodesic}, on voit que la frontière de Martin du graphe $\ym$ est aussi sa frontière géodésique ; ce phénomène est courant en analyse harmonique (\emph{cf.} \cite{BL06}).
\begin{examples}
Considérons le point $\omega=((0,0,\ldots),(0,0,\ldots))$. Il correspond au caractère
$$\tau(\mu)=\prod_{k \geq 2} 0^{m_{k}(\mu)}=\begin{cases} 1 & \text{si }\mu=1^{n}, \\
0 & \text{sinon.}
\end{cases}$$
Le système cohérent de mesures de probabilité correspondant est celui des mesures de Plancherel, car $\tau_{|\sym_{n}}$ est la trace normalisée usuelle de $\C\sym_{n}$. D'après le théorème \ref{boundarygeodesic}, les fréquences limites des diagrammes de Young tirés aléatoirement sous la mesure de Plancherel tendent en probabilité vers $0$. Ce résultat est corroboré par le théorème \ref{firstasymptoticplancherel} : la plus grande ligne et la plus grande colonne d'un diagramme tiré sous la mesure de Plancherel sont de l'ordre de $2\sqrt{n}=o(n)$.\bigskip

\noindent De même, considérons le point $\omega=((0,0,\ldots),(1,0,0,\ldots))$ ; il correspond au caractère signature
$$\eps(\mu)=\prod_{\substack{k \geq 2\\ k\equiv 0 \bmod 2}} (-1)^{m_{k}(\mu)}.$$ 
La restriction de $\eps$ à $\sym_{n}$ est encore la signature, donc le système cohérent de mesures de probabilité associé à $\eps$ est constitué des Dirac en les partitions $1^{n}$, $n \in \N$. Dans ce cas particulier, les fréquences des partitions sont exactement $((0,0,\ldots),(1,0,0,\ldots))$, et ce pour tout $n$.
\end{examples}
\figcapt{\psset{unit=1mm}\pspicture(0,-4)(150,87)
\psdots(75,64)(75,68)(70,72)(80,72)(75,76)(65,76)(85,76)(75,83)(67.5,81.5)(82.5,81.5)(60,80)(90,80)
\psline(75,64)(75,68)(70,72)(65,76)(60,80)
\psline(65,76)(67.5,81.5)(75,76)(80,72)(75,68)
\psline(70,72)(75,76)(75,83)
\psline(75,76)(82.5,81.5)(85,76)(80,72)
\psline(85,76)(90,80)
\psline[linestyle=dashed](60,80)(56,83)
\psline[linestyle=dashed](60,80)(62,83)
\psline[linestyle=dashed](90,80)(94,83)
\psline[linestyle=dashed](90,80)(88,83)
\psline[linestyle=dashed](75,83)(73,86)
\psline[linestyle=dashed](75,83)(77,86)
\psline[linestyle=dashed](67.5,81.5)(67.5,85.5)
\psline[linestyle=dashed](67.5,81.5)(69.5,84.5)
\psline[linestyle=dashed](67.5,81.5)(65.5,84)
\psline[linestyle=dashed](82.5,81.5)(82.5,85.5)
\psline[linestyle=dashed](82.5,81.5)(80.5,84.5)
\psline[linestyle=dashed](82.5,81.5)(84.5,84)
\rput(61,64){$\ym=\bigsqcup_{n \in \N} \ym_{n}$}
\pscircle*[linecolor=gray!30!white](25,40){25}
\psarc(25,40){25}{12}{355}
\rput(49,44){$\partial_{M}\ym$}
\rput(55,40){$=\mathrm{Extr}(\mathcal{H}^{+}(\ym))$}
\rput(20,28){$\mathcal{H}^{+}(\ym)$}
\rput(24,24){$=\mathscr{P}(\partial_{M}\ym)$}
\psline[linecolor=DarkOrchid](25,45)(25,49)(18,52)(25,54.5)
\psline[linecolor=DarkOrchid](25,49)(32,52)(25,54.5)
\psline[linecolor=DarkOrchid](32,52)(39,52)
\psline[linecolor=DarkOrchid](18,52)(11,52)
\psline[linecolor=DarkOrchid](25,54.5)(25,58)
\psline[linecolor=DarkOrchid](25,54.5)(15,56)(11,52)(5,50)
\psline[linecolor=DarkOrchid](25,54.5)(35,56)(39,52)(45,50)
\psline[linecolor=DarkOrchid,linestyle=dashed](5,50)(1.5,46.5)
\psline[linecolor=DarkOrchid,linestyle=dashed](5,50)(6,54)
\psline[linecolor=DarkOrchid,linestyle=dashed](15,56)(8,56)
\psline[linecolor=DarkOrchid,linestyle=dashed](15,56)(13,59)
\psline[linecolor=DarkOrchid,linestyle=dashed](15,56)(19,60)
\psline[linecolor=DarkOrchid,linestyle=dashed](25,58)(22,61)
\psline[linecolor=DarkOrchid,linestyle=dashed](45,50)(48.5,46.5)
\psline[linecolor=DarkOrchid,linestyle=dashed](45,50)(44,54)
\psline[linecolor=DarkOrchid,linestyle=dashed](35,56)(42,56)
\psline[linecolor=DarkOrchid,linestyle=dashed](35,56)(37,59)
\psline[linecolor=DarkOrchid,linestyle=dashed](35,56)(31,60)
\psline[linecolor=DarkOrchid,linestyle=dashed](25,58)(27,61)
\psline{->}(53,78)(25,78)(25,68)
\psline{-}(97,78)(125,78)(125,65)
\rput(20,81){compactification}
\rput(15,77){de Martin}
\rput(132,81){compactification}
\rput(136,77){géodésique}
\rput(125,62){fréquences limites des}
\rput(125,58){coordonnées de Frobenius}
\psline{->}(125,55)(125,45)
\psline*[linecolor=gray!30!white](95,18)(150,18)(140,33)(85,33)(95,18)
\psline(95,18)(150,18)(140,33)(125,33)
\psline(110,33)(85,33)(95,18)
\psline{->}(85,33)(85,50)
\psline(85,33)(84.5,33)
\psline(85,45)(84.5,45)
\rput(83,33){$0$}
\rput(83,45){$1$}
\psline[linecolor=MidnightBlue](117,39)(117,35)(122,32)(117,30)
\psline[linecolor=MidnightBlue](122,32)(128,30)
\psline[linecolor=MidnightBlue](112,32)(106,30)(111,28)(117,30)
\psline[linecolor=MidnightBlue](117,35)(112,32)(117,30)(117,28)
\psline[linecolor=MidnightBlue](117,30)(123,28)(128,30)(135,28)
\psline[linecolor=MidnightBlue](106,30)(99,28)
\psline[linecolor=MidnightBlue,linestyle=dashed](117,28)(114,26)
\psline[linecolor=MidnightBlue,linestyle=dashed](117,28)(120,26)
\psline[linecolor=MidnightBlue,linestyle=dashed](99,28)(95,27)
\psline[linecolor=MidnightBlue,linestyle=dashed](99,28)(103,27)
\psline[linecolor=MidnightBlue,linestyle=dashed](135,28)(139,27)
\psline[linecolor=MidnightBlue,linestyle=dashed](135,28)(131,27)
\psline[linecolor=MidnightBlue,linestyle=dashed](123,28)(128,27)
\psline[linecolor=MidnightBlue,linestyle=dashed](111,28)(106,27)
\rput(140,15){$\partial_{G}\ym=\Omega$}
\rput(75,5){caractères irréductibles du}
\rput(75,1){groupe symétrique infini}
\psframe(50,-2)(100,8)
\psline{->}(25,11)(25,3)(45,3)
\psline{->}(125,11)(125,3)(105,3)
\rput(135,4){fonctions}
\rput(128,0){supersymétriques}
\rput(15,4){systèmes}
\rput(22,0){cohérents de mesures}
\endpspicture}{La compactification de Martin et la compactification géodésique du graphe de Young mènent aux caractères irréductibles du groupe symétrique infini.}{Compactifications du graphe de Young}

Dans la section \ref{qbig}, nous expliquerons comment ces résultats s'adaptent au cas de l'algèbre d'Hecke infinie $\IH_{\infty,q}$, et nous en déduirons une preuve alternative de la loi des grands nombres pour la $q$-mesure de Plancherel. \'Evoquons maintenant d'autres classes de représentations du groupe symétrique infini. Dans ce qui suit, nous noterons $\sym_{\infty}(n)$ le sous-groupe de $\sym_{\infty}$ constitué des permutations qui laissent fixes les $n$ premiers entiers ; et $\overline{\sym_{\infty}}$ le groupe de toutes les permutations de $\N^{*}$. Les sous-groupes $(\sym_{\infty}(n))_{n \in \N}$ forment un système de voisinages de l'identité pour la topologie sur $\sym_{\infty}$ qui correspond à la convergence simple des permutations --- c'est la seule topologie raisonnable en dehors de la topologie discrète. On doit à A. Lieberman (voir \cite{Lie72,Ols85}) une caractérisation des représentations unitaires $U : \sym_{\infty} \to \mathscr{U}(H)$ qui sont continues pour la topologie précédemment décrite sur $\sym_{\infty}$, et pour la topologie faible d'opérateur
sur $\mathscr{U}(H)$. Ainsi, les assertions suivantes sont équivalentes :
\begin{enumerate}[(i)]
\item La représentation $U$ est continue pour la topologie de la convergence simple.
\item La représentation $U$ se prolonge par continuité en une représentation de $\overline{\sym_{\infty}}$.
\item La représentation $U$ est \textbf{modérée}, \emph{i.e.}, $\mathscr{M}(\sym_{\infty},U)$ est une algèbre de von Neumann de type $\text{I}$.
\item Tout vecteur $\xi \in H$ est limite d'une suite de vecteurs $(\xi_{n})_{n \in \N}$ avec $\xi_{n}$ laissé fixe par $\sym_{\infty}(n)$. 
\end{enumerate}
La dernière assertion revient à dire que $H=\overline{\bigcup_{n \in \N} H^{\sym_\infty(n)}}$, et  le \textbf{niveau} d'une représentation modérée est défini comme le plus petit entier $N$ tel que $H^{\sym_{\infty}(N)}\neq \{0\}$. 
\begin{proposition}[Représentations modérées du groupe symétrique infini, \cite{Oko97}]\label{tamerepresentation}
Toute représentation modérée de $\sym_{\infty}$ s'écrit de manière unique comme somme directe (hilbertienne) de représentations modérées irréductibles. Ces représentations irréductibles sont indexées par tous les diagrammes de Young de $\ym=\bigsqcup_{n \in \N} \ym_{n}$, et de plus, la $k$-ième puissance tensorielle de la représentation modérée <<~canonique~>>
$$\sym_{\infty} \to \mathscr{U}(\ell^{2}(\N^{*}))$$
par permutations des coordonnées met en jeu tous les diagrammes de Young de poids compris entre $1$ et $k$. En particulier, la représentation canonique $\sym_{\infty} \to \mathscr{U}(\ell^{2}(\N^{*}))$ correspond à la partition $(1)$ et est modérée irréductible, de type $\text{I}_{\infty}$.
\end{proposition}\bigskip\bigskip

Dans un contexte d'analyse harmonique non commutative, la proposition \ref{tamerepresentation} rapproche $K=\sym_{\infty}$ des groupes topologiques compacts,  puisqu'on dispose d'une décomposition en sommes d'irréductibles pour toute représentation modérée. Partant, il est naturel d'étudier les représentations de \textbf{paires de Gelfand} $(G,K)$ avec $\sym_{\infty}=K \subset G$. Ces paires de groupes discrets peuvent être caractérisées par la propriété suivante : si $U : G \to \mathscr{U}(H)$ est une représentation unitaire irréductible de $G$ telle que $H^{K}\neq \{0\}$ (auquel cas on parle de \textbf{représentation sphérique}), alors $H^{K}$ est de dimension $1$. Autrement dit, la multiplicité de la représentation triviale de $K$ dans la restriction d'une représentation irréductible de $G$ est toujours\footnote{On parle de \textbf{paire de Gelfand forte} si l'on a cette propriété de multiplicité simple pour toutes les représentations irréductibles de $K$, et pas seulement la représentation triviale.} $0$ ou $1$. 
\begin{example} Compte tenu du théorème de branchement \ref{branching}, $\sym_{n} \subset \sym_{n+1}$ est une paire de Gelfand forte pour tout entier $n$, \emph{cf.} \cite[\S2]{OV04}.
\end{example}\bigskip

\noindent On appelle \textbf{représentation admissible} d'une paire de Gelfand $(G,K)$ une représentation unitaire $U : G \to \mathscr{U}(H)$ telle que $\mathrm{Res}_{K}^{G}(U)$ soit modérée. Ces représentations admissibles s'écrivent de manière unique comme intégrales directes de représentations admissibles irréductibles, et pour $K=\sym_{\infty}$, elles forment la troisième classe de représentations pour lesquelles des résultats de classification sont connus. Ainsi, dans \cite{Ols90,Oko97}, diverses paires de Gelfand $\sym_{\infty} \subset G$ sont décrites. La plus naturelle est 
$$G=\sym_{\infty} \times \sym_{\infty} \qquad;\qquad K=\mathrm{diag}(G)=\sym_{\infty}\,,$$
et les représentations admissibles de cette paire sont exactement celles qui s'étendent par continuité au groupe (modéré, non localement compact) 
$$\overline{G}=\{(\sigma,\tau) \in \overline{\sym_{\infty}} \times \overline{\sym_{\infty}}\,\,|\,\,\sigma\tau^{-1} \in K=\sym_{\infty}\}\,.$$
D'autre part, ces représentations admissibles sont stables par somme directe, produit tensoriel, et restriction ou induction à partir des sous-groupes $\sym_{\infty}(n) \times \sym_{\infty}(n)$. 
\begin{proposition}[Représentations admissibles du groupe symétrique infini, \cite{Oko97}]
On identifie un point $\omega=(\alpha,\beta)$ du simplexe de Thoma à une partie $\mathcal{I}(\omega)=\alpha \cup -\beta \cup \{0\}$ de l'intervalle $[-1,1]$, avec pour chaque point $x$ de cette partie une multiplicité $n(x) \in \lle 1,\infty \rre$. Une distribution de Young $\Lambda$ adaptée à un point $\omega \in \Omega$ est la donnée de diagrammes de Young $\lambda(x)$ pour tout $x \in \mathcal{I}(\omega)$, de telle sorte que $|\Lambda|=\sum_{x \in \mathcal{I}(\omega)}\lambda(x)$ soit fini, voir la figure \ref{youngdistribution}. Alors, une représentation admissible irréductible de $(G,K)$ de niveau $N$ est entièrement déterminée par la donnée d'un point $\omega \in \Omega$ et de deux distributions de Young $\Lambda$ et $\Pi$ adaptées à ce point, et telles que 
$$|\Lambda|=|\Pi|=N\qquad;\qquad\forall x \in [-1,1],\,\,\,\ell(\lambda(x)) + \ell(\pi(x)) \leq n(x)\,.$$
\end{proposition}

\figcapt{\psset{unit=1mm}\pspicture(0,3)(80,38) \psline{[-]}(0,5)(80,5) \rput(-3,5){$-1$} \rput(83,5){$+1$} \psdots(20,5)(40,5)(55,5)(65,5)
\psline(55,5)(55,10)(52.5,12.5)(55,15)(57.5,12.5)(55,10)
\psline(20,5)(20,20)(10,30)(12.5,32.5)(15,30)(17.5,32.5)(20,30)(25,35)(27.5,32.5)(32.5,37.5)(35,35)(20,20)
\psline(12.5,27.5)(15,30)(17.5,27.5)(20,30)(22.5,27.5)(27.5,32.5)(30,30) \psline(30,35)(32.5,32.5)
\psline(15,25)(17.5,27.5)(20,25)(22.5,27.5)(25,25) \psline(27.5,27.5)(22.5,32.5)
\psline(17.5,22.5)(20,25)(22.5,22.5)
\psline(40,5)(40,15)(35,20)(40,25)(42.5,22.5)(45,25)(47.5,22.5)(40,15)
\psline(37.5,17.5)(42.5,22.5)(45,20)\psline(37.5,22.5)(42.5,17.5)
\psline(65,5)(65,15)(52.5,27.5)(60,35)(65,30)(67.5,32.5)(75,25)(65,15)
\psline(55,25)(62.5,32.5) \psline(57.5,22.5)(65,30)(72.5,22.5) \psline(57.5,32.5)(70,20) \psline(55,30)(67.5,17.5) \psline(60,20)(70,30) \psline(62.5,17.5)(72.5,27.5)
\endpspicture}{Distribution de Young $\Lambda$ de niveau $1+5+13+18=37$ adaptée à un point du simplexe de Thoma.\label{youngdistribution}}{Distribution de Young adaptée à un point du simplexe de Thoma}

\noindent Les représentations sphériques de $(G,K)$ sont précisément les représentations admissibles irréductibles de niveau $0$, c'est-à-dire que $\Lambda=\Pi=\emptyset$ ; elles sont donc en bijection avec les points du simplexe de Thoma, et étant donnée une représentation sphérique, la donnée d'un vecteur sphérique fournit effectivement un caractère irréductible de $\sym_{\infty}$, voir \cite[\S0.2, p. 6]{Oko97}.\medskip

\begin{example}
On considère la représentation régulière bilatère de $G$ sur $H=\ell^{2}(\sym_{\infty})$ ; elle est donnée par $(g,h) \cdot \sigma = g\,\sigma \,h^{-1}$, et est irréductible. La restriction de cette représentation à $K=\sym_{\infty}$ est l'action par conjugaison de $K$ sur lui-même ; en particulier, $$H^{\sym_{\infty}(n)}=\ell^{2}(\sym_{n})$$ pour tout $n$, et $H^{K}=\C1$. La représentation régulière bilatère est donc une représentation sphérique, et comme la fonction sphérique correspondante est
$$f(g,h)=\mathbb{1}_{g=h^{-1}},$$
elle correspond à la trace régulière de $\sym_{\infty}$, c'est-à-dire le caractère irréductible de paramètre de Thoma $\omega=((0,0,\ldots),(0,0,\ldots))$.
\end{example}\bigskip\bigskip

Pour conclure ce paragraphe, évoquons les \textbf{systèmes cohérents de $z$-mesures}, qui sont parmi les premiers à avoir été décrits par des systèmes de particules à corrélations déterminantales, \emph{cf.} \cite{BO98,BO00}. Un résultat classique d'analyse harmonique non commutative est le théorème de Peter-Weyl, qui assure que si $K$ est un groupe topologique compact, alors la représentation régulière bilatère de $G=K\times K$ sur $\leb^{2}(K)$ se décompose en 
$$(\leb^{2}(K),L \otimes R )=\bigoplus_{\lambda \in \widehat{G}} V^{\lambda} \otimes V^{\lambda*},$$
voir \cite[\S1]{JS91} et \cite[\S2]{Var89}. Pour le groupe symétrique infini $K=\sym_{\infty}$, ce résultat tombe malheureusement en défaut, car la représentation régulière bilatère est irréductible. Néanmoins, on peut construire des représentations régulières <<~généralisées~>> de $G=\sym_{\infty} \times \sym_{\infty}$ qui sont des intégrales directes de représentations irréductibles sphériques, et qui correspondent à des mesures non triviales sur le simplexe de Thoma. L'idée est de construire une limite projective de $G$-ensembles\footnote{L'ensemble $X$ est constitué des \textbf{permutations virtuelles}, voir \cite[\S1]{KOV04} pour une définition précise. Le défaut majeur de la construction est que $X$ n'est pas un (semi-)groupe.} $X=\varprojlim_{n \to \infty}\sym_{n}$, de munir $X$ d'une mesure limite projective $\mu_{t}$ qui est $K$-invariante et $G$-quasi-invariante, et de tordre l'action usuelle 
$$(g,h) \in G\,,\,\,\sigma \in K \mapsto (g,h)\cdot \sigma = g\sigma h^{-1}$$
par le cocycle additif qui mesure la $G$-quasi-invariance de la mesure $\mu_{t}$, \emph{cf.} \cite[\S3]{KOV04}. Ainsi, pour tout nombre complexe non nul $z$, si $t=|z|^{2}$, alors on peut construire une représentation unitaire $T_{z} : G \to \mathscr{U}(\leb^{2}(X,\mu_{t}))$, avec les propriétés suivantes :
\begin{proposition}[Systèmes cohérents de $z$-mesures, \cite{KOV04}]\label{zmeasure}
Pour tout $z \in \C\setminus \{0\}$, la représentation $T_{z}$ est l'intégrale directe de représentation sphériques :
$$T_{z}=\int_{\Omega} S^{\omega} \,P_{z}(d\omega)\,,$$
où $S^{\omega}$ désigne la représentation sphérique associé au caractère irréductible $\chi^{\omega}$. Le système cohérent de mesures de probabilité associé à $T_{z}$ et $P_{z}$ est 
$$M_{z,n}(\lambda)=\frac{(\dim \lambda)^{2}}{n!\,(t)_{n}}\,\prod_{(i,j) \in \lambda}|z+i-j|^{2}$$
où $(t)_{n}=t(t+1)(t+2)\cdots(t+n-1)$ est le symbole de Pochhammer. On retrouve le système des mesures de Plancherel lorsque $z$ tend vers l'infini, et un autre système cohérent limite lorsque $z$ tend vers $0$, avec $M_{0,n}(\lambda)=1/n$ si $\lambda$ est une partition équerre, et $0$ sinon.
\end{proposition}\bigskip\bigskip

Pour compléter la proposition \ref{zmeasure}, décrivons les \textbf{mesures spectrales} $P_{z}$ des représentations régulières généralisées. Pour tout $z \in \C\cup \{\infty\}$, les représentations $T_{z}$ et $T_{\overline{z}}$ sont équivalentes, et les représentations $T_{z}$ et $T_{-z}\otimes (\mathrm{sgn} \times \mathrm{sgn})$ sont équivalentes, \emph{cf.} \cite[propositions 3.6.1 et 4.5.1]{KOV04}. Deux mesures spectrales $P_{z}$ et $P_{z'}$ sont étrangères sauf si $z'=z$ ou $z'=\overline{z}$. On doit ensuite distinguer deux cas :\vspace{2mm}
\begin{enumerate}
\item Si $z=k$ est un entier relatif non nul, alors $T_{k}$ est somme directe de représentations $T_{pq}$ avec $k=p-q$ et $p,q \in \N$. Les mesures spectrales de ces représentations sont équivalentes aux mesures de Lebesgue des faces $(p+q-1)$-dimensionnelles 
$$\Omega(p,q)=\big\{(\alpha,\beta) \in \Omega \,\,|\,\,\alpha_{1}+\cdots +\alpha_{p}+\beta_{1}+\cdots+\beta_{q}=1\big\}$$
du simplexe de Thoma, avec des densités explicitement calculables. La mesure $M_{k}$ est portée par les diagrammes de Young ayant moins de $k$ lignes si $k>0$ (resp., par les diagrammes de Young ayant moins de $k$ colonnes si $k<0$).\vspace{2mm}
\item Si $z\in \C\setminus \Z$, alors $P_{z}$ ne charge pas les faces $\Omega(p,q)$, et est portée par la face $\Omega(0)=\{\omega \in \Omega \,\,|\,\, \gamma=0\}$. La mesure $M_{z}$ charge tous les diagrammes de Young.\vspace{2mm}
\end{enumerate}
Compte tenu de ces résultats, on peut estimer que les représentations régulières généralisées $T_{z}$ <<~engendrent~>> suffisamment de représentations sphériques, et pallient le défaut de la représentation régulière bilatère de $\sym_{\infty}$ observé précédemment. \bigskip

\section[Processus ponctuels déterminantaux et mesures de Schur]{Processus ponctuels déterminantaux et mesures de Schur}\label{schurmeasure}
Si $z$ est un paramètre complexe, les mesures $M_{z,n}$ des diagrammes de Young sont explicitement calculables, mais les valeurs du caractère normalisé $\chi_{z}$ n'ont pas d'expression simple, même en un cycle --- par exemple, les premières valeurs sont
$$\chi_{z}(1)=1 \qquad;\qquad\chi_{z}(2)=\frac{2\,\mathrm{Re}(z)}{|z|^{2}+1}\qquad;\qquad \chi_{z}(3)=\frac{|z|^{2}+4(\mathrm{Re}(z))^{2}+1}{(|z|^{2}+1)(|z|^{2}+2)}\,.$$
Par conséquent, les techniques d'observables de diagrammes présentées dans le chapitre \ref{tool} ne permettent pas l'étude asymptotique des $z$-mesures. Ceci a conduit à l'élaboration d'une seconde technique d'étude asymptotique des diagrammes aléatoires, qui repose sur la théorie des \textbf{processus ponctuels déterminantaux} (\cite{BO98,BO00,BOO00,BO05a,BO05b}). Si $X$ est un espace topologique séparable localement compact, on appelle \textbf{configuration de points} dans $X$ un multi-ensemble\footnote{Ainsi, on autorise des répétitions dans la suite $\{\xi_{1},\xi_{2},\ldots\}$, et ces points ne sont pas ordonnés.} $\xi=\{\xi_{1},\xi_{2},\ldots\}$ tel que $\card (\xi \cap A)< \infty$ pour toute partie compacte de $X$. L'ensemble $\Xi(X)$ des configurations de points est muni de la tribu engendrée par les fonctions $\xi \mapsto \card (\xi \cap A)$, et on appelle \textbf{processus ponctuel} à valeurs dans $X$ une application mesurable $\xi$ d'un espace de probabilité $(\Omega,\mathscr{F},\proba)$ vers $\Xi(X)$. On obtient ainsi une mesure de probabilité $\xi_{\star}\proba$ sur l'espace $\Xi(X)$, qu'on étudie grâce aux \textbf{mesures de corrélation} $\mu^{k}$ :
$$\mu^{k}(A_{1}\times A_{2}\times \cdots \times A_{k})=\esper\left[\sum_{i_{1}\neq i_{2}\neq \cdots \neq i_{k}}\mathbb{1}_{A_{1}}(\xi_{i_{1}}) \times \mathbb{1}_{A_{2}}(\xi_{i_{2}}) \times \cdots \times \mathbb{1}_{A_{k}}(\xi_{i_{k}})\right].$$
Chaque $\mu^{k}$ est une mesure sur $X^{k}$, et lorsque $X$ est une partie de $\R$, on appelle \textbf{fonctions de corrélation} du processus les densités $\rho^{k}(x)=d\mu^{k}(x)/dx^{k}$ par rapport aux mesures de Lebesgue (en supposant que les $\mu^{k}$ sont absolument continues par rapport aux mesures de Lebesgue).
\begin{example}
La première fonction de corrélation $\rho^{1}$ peut être vue comme la densité du processus ponctuel : ainsi, pour toute fonction $f$, $\int_{X} f(x)\,\rho^{1}(x)\,dx=\esper[\sum_{i}f(\xi)]$.
\end{example}\medskip

\noindent Lorsque $X$ est un espace discret (dénombrable), la donnée d'un processus ponctuel est équivalente à la donnée d'une probabilité $\pi$ sur l'ensemble dénombrable des parties finies de $X$, et les mesures de corrélation sont données par :
$$\mu^{k}(\{x_{1}\}\times\{x_{2}\}\times \cdots \times \{x_{k}\})=\proba[\{x_{1},x_{2},\ldots,x_{k}\} \subset \xi]\,.$$
On appellera dans ce cas fonction de corrélation l'application $\rho(A)=\proba[A \subset\xi]=\sum_{A \subset B}\pi(B)$, et on notera $\rho^{k}(x_{1},\ldots,x_{k})=\rho(\{x_{1},\ldots,x_{k}\})$.\bigskip
\bigskip

Un processus ponctuel à valeurs dans une partie de $\R$ ou dans un espace discret est dit \textbf{déterminantal} si ses fonctions de corrélation sont données par des déterminants, c'est-à-dire qu'il existe un noyau (symétrique) $K : X \times X \to X$ tel que 
$$\rho^{k}(x_{1},\ldots,x_{k})=\det((K(x_{i},x_{j}))_{1\leq i,j\leq k})\,.$$
Ces processus déterminantaux sont associés à des \textbf{opérateurs à trace} $\mathcal{K}$ positifs  sur des espaces de Hilbert. Ainsi, si $X$ est une partie de $\R$ et si $\xi$ est un processus ponctuel déterminantal de noyau $K$, l'opérateur $\mathcal{K} :\leb^{2}(X)\to\leb^{2}(X)$ défini par
$$ \mathcal{K}f(x) =\int_{X}K(x,y)f(y)\,dy$$
et ses restrictions à des sous-espaces $\leb^{2}(A \subset X)$ jouent un rôle important dans de nombreux calculs de probabilités du processus $\xi$. De même, si $X$ est un ensemble dénombrable et si $\mathcal{L}$ est un opérateur à trace sur $\ell^{2}(X)$, notant $p_{A}$ la projection sur le sous-espace $\ell^{2}(A)$ pour toute partie $A \subset X$, et $\mathcal{L}_{A}=p_{A}\mathcal{L}p_{A}$, on peut définir un processus ponctuel déterminantal en posant\footnote{Ici, $\det(\id +\mathcal{L})$ est un déterminant au sens de Fredholm.} :
$$\det(\id + \mathcal{L})=\sum_{k=0}^{\infty}\tr(\wedge^{k}\mathcal{L}) \,= \!\!\sum_{\substack{A \subset X\\ \card A<\infty}} \!\!\det(\mathcal{L}_{A}) \qquad;\qquad\pi(A)= \frac{\det(\mathcal{L}_{A})}{\det (\id+\mathcal{L})}\,.$$
Dans ce contexte, les fonctions de corrélation du processus ponctuel sont déterminantales : si $\mathcal{K}=\mathcal{L}\,(\id+\mathcal{L})^{-1}$, alors $\rho(X)=\det(\mathcal{K}_{X})$, \emph{cf.} \cite[\S2]{BO00}. On renvoie à \cite{DVJ88} pour un traitement plus complet de la théorie des processus ponctuels.
\bigskip
\bigskip

Si $\lambda$ est un diagramme de Young aléatoire, on peut lui associer un processus ponctuel à valeurs dans $\Z'=\Z+1/2$ en considérant l'ensemble $\mathcal{F}(\lambda)=A(\lambda)\cup -B(\lambda)$ des coordonnées de Frobenius. Un autre processus ponctuel est obtenu en regardant les ccordonnées de descente de la fonction $\omega_{\lambda}$ :
$$\mathcal{D}(\lambda)=\Z_{-}' \mathrel{\Delta} \mathcal{F}(\lambda)=\{\lambda_{i}-i+1/2\}_{i \geq 1}\,.$$
\begin{example}
Si $\lambda=(5,4,4,1)$, alors $\mathcal{D}(\lambda)=\big(\frac{9}{2},\frac{5}{2},\frac{3}{2},-\frac{5}{2},-\frac{9}{2},-\frac{11}{2},-\frac{13}{2},\ldots\big)$.
\end{example}
\noindent Avec ces conventions, un processus de croissance de diagrammes tel que le processus de Plancherel peut être interprété comme un système de particules s'excluant mutuellement et partant de la configuration initiale $\Z_{-}'$, voir la figure \ref{particle}.
\figcapt{\psset{unit=1mm}\pspicture(20,0)(150,62)
\psline{->}(70,0)(150,0)
\rput(154,0.5){$\Z'$}
\psline(140,60)(140,44)(135,40)(135,36)(125,28)(125,8)(115,0)
\psline(130,60)(125,56)(125,40)(120,36)(120,20)(115,16)(115,12)(110,8)(110,0)
\psline(115,60)(115,48)(110,44)(110,16)(105,12)(105,0)
\psline(110,60)(110,56)(105,52)(105,24)(100,20)(100,0)
\psline(100,60)(100,28)(95,24)(95,0)
\psline(95,60)(95,52)(90,48)(90,0)
\psline(85,60)(85,0)
\psline(80,60)(80,0)
\psdots(80,60)(85,60)(95,60)(100,60)(110,60)(115,60)(130,60)(140,60)
\psline(20,10)(60,10)(60,18)(52,18)(52,26)(36,26)(36,42)(28,42)(28,58)(20,58)(20,10)
\psline(52,10)(52,18)(44,18)(44,26)(28,26)(28,42)(20,42)
\psline(44,10)(44,18)(36,18)(36,34)(20,34)
\psline(20,50)(28,50)
\psline(20,26)(28,26)(28,18)(36,18)(36,10)
\psline(20,18)(28,18)(28,10)
\psline(117.5,0.5)(117.5,-0.5)
\rput(24,14){$1$}
\rput(32,14){$2$}
\rput(24,22){$3$}
\rput(24,30){$4$}
\rput(32,22){$5$}
\rput(24,38){$6$}
\rput(24,46){$7$}
\rput(40,14){$8$}
\rput(48,14){$9$}
\rput(40,22){$10$}
\rput(56,14){$11$}
\rput(32,30){$12$}
\rput(24,54){$13$}
\rput(32,38){$14$}
\rput(48,22){$15$}
\endpspicture}{Codage d'un tableau par l'évolution d'un système de particules s'excluant mutuellement.\label{particle}}{Codage d'un tableau par l'évolution d'un système de particules}\bigskip
\bigskip

D'autre part, si $\omega=(\alpha,\beta)$ est un point du simplexe de Thoma, on peut lui associer comme précédemment une configuration de points $\alpha \cup(-\beta)$ sur l'intervalle épointé $I^{*}=[-1,1]\setminus \{0\}$. Par conséquent, une mesure spectrale $P$ sur $\Omega$ correspond à un processus ponctuel à valeurs dans $I^{*}$, que nous noterons $\mathcal{I}_{\star}(P)$. Dans le contexte des $z$-mesures, les deux constructions de processus ponctuels précédemment décrites peuvent être reliées comme suit. Si $M_{\chi}=(M_{\chi,n})_{n\in \N}$ est un système cohérent de mesures de probabilité sur le graphe de Young $\ym$, et si $\mu$ est une mesure de probabilité sur $\N$, notons $M_{\chi,\mu}$ la mesure de probabilité sur $\ym$ définie par :
$$M_{\chi,\mu}(\lambda)=\mu(|\lambda|)\,M_{\chi,|\lambda|}(\lambda)\,.$$
\begin{example}
Pour $\theta>0$, la mesure de Plancherel poissonisée de paramètre $\theta$ est définie par
$$M_{\mathcal{P}(\theta)}(\lambda)=\frac{\theta^{|\lambda|}\,\E^{-\theta}}{|\lambda|!}\,\times\, \frac{(\dim \lambda)^{2}}{|\lambda|!}\,.$$
Les processus ponctuels correspondants sont $\mathcal{F}_{\star}(M_{\mathcal{P}(\theta)})$ et $\mathcal{D}_{\star}(M_{\mathcal{P}(\theta)})$, et l'asymptotique des mesures de Plancherel peut être retrouvée en étudiant la limite $\theta \to \infty$, voir \cite{BOO00} et la section \ref{bdjdeterminantal}.
\end{example}\bigskip\bigskip

Plus généralement, étant donnée une mesure spectrale $P_{\chi}$ sur $\Omega$, si $(M_{\chi,n})_{n \in \N}$ est le système cohérent de mesures et si $(\mu_{\theta})_{\theta \in \R_{+}}$ est une famille de mesures de probabilité sur $\N$ telle que $n \to_{\mu_{\theta}} +\infty$, alors compte tenu du théorème \ref{boundarygeodesic}, on peut s'attendre à ce que les processus ponctuels $\mathcal{F}_{\star}(M_{\chi,\mu_{\theta}})$ correctement renormalisés approchent le processus ponctuel $\mathcal{I}_{\star}(P_{\chi})$. Dans le cas des $z$-mesures, cette intuition peut être rendue rigoureuse, et tous les processus mis en jeu sont déterminantaux. Fixons un paramètre complexe $z \in \C \setminus \Z$ ; on rappelle que la mesure $P_{z}$ est portée par la face $\Omega(0)=\{\omega \in \Omega \,\,|\,\,\gamma=0\}$ du simplexe de Thoma. On note $\Omega^{*}=\Omega(0)\times\R_{+}^{*}$ ; via l'application $(\omega,s) \mapsto (s\alpha,s\beta)$, on peut identifier $\Omega^{*}$ et 
$$\bigg\{(a,b)=(a_{1}\geq a_{2} \geq \cdots \geq 0, b_{1}\geq b_{2} \geq \cdots \geq 0) \,\,\,\bigg|\,\,\,\sum_{i=1}^{\infty}a_{i}+\sum_{i=1}^{\infty}b_{i}=s \in \R_{+}^{*}\bigg\}\,. $$
Les éléments de $\Omega^{*}$ sont associés à des configurations de points dans l'intervalle épointé $\R^{*}=\R\setminus\{0\}$.
\begin{proposition}[$z$-mesures et processus ponctuels déterminantaux, \cite{BO00}]\label{asymptoticzmeasure}~
\begin{enumerate}
\item Soit $P_{z}^{*}$ la loi produit $P_{z} \otimes \{s^{t-1}\,\E^{-s}\,ds/\Gamma(t)\}$ sur $\Omega^{*}$. Le processus ponctuel $\mathcal{I}_{\star}(P^{*}_{z})$ est déterminantal, de fonctions de corrélation données par le noyau de Whittaker, \emph{cf.} \cite[\S5]{BO00}. \item Soit $\mathcal{B}(t,\zeta)$ la loi binomiale négative de paramètre $(t,\zeta)$, c'est-à-dire que 
$$\mathcal{B}(t,\zeta)(\{k\})=\frac{(t)_{k}}{k!}\,(1-\zeta)^{t}\,\zeta^{k}\,.$$
Les processus ponctuels $\mathcal{F}_{\star}(M_{z,\mathcal{B}(t,\zeta)})$ sont déterminantaux, de noyaux donnés par des fonctions hypergéométriques, \emph{cf.} \cite[\S3]{BO00}. De plus, lorsque $\zeta$ tend vers $1$, le processus renormalisé $(1-\zeta)\,\mathcal{F}_{\star}(M_{z,\mathcal{B}(t,\zeta)})$, qui est à valeurs dans $\R^{*}$, converge en loi vers $\mathcal{I}_{\star}(P_{z}^{*})$.
\end{enumerate}
\end{proposition}
\noindent On renvoie aux articles \cite{BO98,BO00,BO05b} pour des précisions sur ce résultat, qui se généralise d'ailleurs à une famille de mesures spectrales à trois paramètres $(z,z',\theta)$. L'asymptotique des mesures de Plancherel est un cas dégénéré de la proposition \ref{asymptoticzmeasure} ; nous expliquerons ceci dans la section \ref{bdjdeterminantal}.\bigskip
\bigskip

Ainsi, en s'appuyant sur l'exemple des $z$-mesures, nous avons présenté une seconde technique d'étude des diagrammes de Young aléatoires de grande taille. Notons que le caractère déterminantal des fonctions de corrélation joue un rôle essentiel dans la preuve de résultats asymptotiques tels que la proposition \ref{asymptoticzmeasure}, car les calculs sont ramenés à des fonctions de deux variables. Il convient dès lors de se demander quelles mesures $\mathcal{M}$ sur $\ym$ donnent lieu à des processus ponctuels déterminantaux. Si $\mathcal{D}_{\star}(\mathcal{M})$ est un processus ponctuel déterminantal, alors un principe de complémentation (\emph{cf.} \cite[p. 39]{BOO00}) assure que $\mathcal{F}_{\star}(\mathcal{M})=\mathcal{D}_{\star}(\mathcal{M}) \mathrel{\Delta} \Z_{-}'$ est aussi déterminantal,  les noyaux vérifiant la relation
$$K_{\mathcal{D}}=\begin{pmatrix} A & B \\ C & D\end{pmatrix}\qquad;\qquad K_{\mathcal{F}}=\begin{pmatrix} A & B \\ -C & 1-D\end{pmatrix}\,,$$
si les matrices sont découpées suivant la décomposition $\Z'=\Z_{+}' \sqcup \Z_{-}'$. Ainsi, on peut toujours se ramener aux coordonnées de descente $\mathcal{D}$, et une classe très générale de  mesures $\mathcal{S}$ telles que $\mathcal{D}_{\star}(\mathcal{S})$ soit déterminantal est décrite par A. Okounkov dans \cite{Oko01,Oko03a,Oko03b}. Cette classe contient les mesures de Plancherel et les $z$-mesures, et nous verrons dans la seconde partie qu'elle contient aussi les $q$-mesures de Plancherel et les limites des <<~lois marginales~>> des mesures de Plancherel des groupes $\GL(n,\For_{q})$ (\emph{cf.} \cite{Dud08,Ful06}). L'idée est de partir de l'\textbf{identité de Cauchy} pour les fonctions de Schur 
$$\forall a=(a_{i})_{i},b=(b_{j})_{j},\,\,\,\sum_{\lambda \in \ym} s_{\lambda}(a)\,s_{\lambda}(b)=\prod_{i,j}\frac{1}{1-a_{i}b_{j}}\,$$
cette identité étant valable formellement et pour tous paramètres $a$ et $b$ finis ou sommables. Alors, $$\mathcal{S}(\lambda)=\mathcal{S}_{a,b}(\lambda)=s_{\lambda}(a)\,s_{\lambda}(b)\,\prod_{i,j}(1-a_{i}b_{j})$$ est une mesure de probabilité sur $\ym$, qu'on appellera \textbf{mesure de Schur} de paramètres $a$ et $b$. Dans ce qui suit, on exprimera plutôt les mesures de Schur en fonction des paramètres de Miwa $$t_{k}=p_{k}(a)/k \,\,\text{ et }\,\,t_{k}'=p_{k}(b)/k\,.$$ La formule de Frobenius et l'identité de Cauchy permettent de récupérer les fonctions de Schur et la constante de normalisation $Z^{-1}=\prod_{i,j}(1-a_{i}b_{j})$ :
$$ s_{\lambda}=\sum_{\mu \in \ym_{n}} \varsigma^{\lambda}(\mu)\, \prod_{k \geq 1} \frac{(t_{k})^{m_{k}(\mu)}}{m_{k}(\mu)!}\qquad;\qquad \log Z=\sum_{k\geq 1} k\,t_{k}t_{k'}\,.$$
\begin{examples}\label{exampleschurmeasure}
Si $t=t'=(\sqrt{\theta},0,0,\ldots)$, alors $Z=\E^{\theta}$, $s_{\lambda}(a)=s_{\lambda}(b)=(\dim \lambda)\,\theta^{n/2}/n!$, et $\mathcal{S}(\lambda)=M_{\mathcal{P}(\theta)}(\lambda)$ est la mesure de Plancherel modulée selon une loi de Poisson (\cite[p. 63]{Pyn73}) de paramètre $\theta$. \bigskip

\noindent De même, si $t_{k}=\frac{\zeta^{k/2}\,z}{k}$ et si $t_{k'}=\frac{\zeta^{k/2}\,\overline{z}}{k}$, alors $Z^{-1}=(1-\zeta)^{t}$ et $s_{\lambda}(a)\,s_{\lambda}(b)=\zeta^{|\lambda|}\,s_{\lambda}(1^{z})\,s_{\lambda}(1^{\overline{z}})$, où par $s_{\lambda}(1^{z})$ on entend
$$s_{\lambda}(\underbrace{1,1,\ldots,1}_{z \text{ termes}})=\prod_{(i,j) \in \lambda}\frac{z+c(i,j)}{h(i,j)}\,,$$
voir \cite[\S1.3, exemple 4]{Mac95}. On en déduit que $\mathcal{S}(\lambda)=M_{z,\mathcal{B}(t,\zeta)}(\lambda)$.
\end{examples}
\bigskip

\begin{theorem}[Mesures de Schur, \cite{Oko01}]\label{schurmeasurearedeterminantal}
Si $\mathcal{S}=\mathcal{S}_{a,b}$ est une mesure de Schur, alors les processus ponctuels $\mathcal{D}_{\star}(\mathcal{S})$ et $\mathcal{F}_{\star}(\mathcal{S})$ sont déterminantaux, et les noyaux correspondants $K_{\mathcal{D}}$ et $K_{\mathcal{F}}$ ont des expressions explicites en termes des paramètres de Miwa. En particulier, si $T(z)=\sum_{k=1}^{\infty} t_{k}\,z^{k}-\sum_{k=1}^{\infty}t_{k}'\,z^{-k}$, alors 
$$K_{\mathcal{D}}(x,y)=\frac{1}{(2\I\pi)^{2}} \iint_{|w|<|z|} \frac{\E^{T(z)-T(w)}}{(z-w)\sqrt{zw}} \,\frac{dz\,dw}{z^{x}\,w^{-y}}\,.$$
\end{theorem}\bigskip\bigskip

Ce résultat peut être obtenu en utilisant le formalisme du \textbf{produit extérieur infini}, voir \cite{Oko01}. Soit $V$ un espace de Hilbert séparable de dimension infinie, et $(\underline{z})_{z \in \Z'}$ une base orthonormale de $V$. Le produit extérieur infini $\bigwedge^{\infty}V$ est l'espace de Hilbert de base orthonormale
$$\underline{D}=\underline{d_{1}}\wedge \underline{d_{2}} \wedge \underline{d_{3}} \wedge \cdots$$
où $D=\{d_{1} > d_{2} > d_{3} > \cdots \} \subset \Z'$ est une suite décroissante infinie telle que $D \mathrel{\Delta} \Z_{-}'$ soit fini --- autrement dit, $D$ contient une demi-droite infinie. L'espace de départ $V$ agit sur $\bigwedge^{\infty}V$ par produit extérieur : $\forall k\in \Z',\,\,\,\psi_{k}(v)=\underline{k} \wedge v $. Les opérateurs $\psi_{k}$ et leurs adjoints $\psi_{k}^{*}$ vérifient les relations d'anti-commutation des \textbf{fermions libres} :
$$\psi_k\psi_l^*+\psi_l^*\psi_k=\delta_{kl}\qquad ;\qquad \psi_k\psi_l+\psi_l\psi_k=0 \qquad;\qquad \psi_k^*\psi_l^*+\psi_l^*\psi_k^*=0 \,.$$
Autrement dit, l'algèbre engendrée par ces opérateurs est l'algèbre de Clifford de la forme bilinéaire non dégénérée définie sur $V \oplus V^{*}$ par $\scal{\psi_{k}^{*}}{\psi_{l}}=\delta_{kl}$, $\scal{\psi_{k}}{\psi_{l}}=0$ et $\scal{\psi_{k}^{*}}{\psi_{l}^{*}}=0$. Notons $:\!\psi_k\psi_l^*\!:$ le produit ordonné défini par :
$$:\!\psi_k\psi_l^*\!:\,\, = \begin{cases} \psi_k\psi_l^*&\text{if }l>0,\\
                     -\psi_l^*\psi_k&\text{if }l<0.
                     \end{cases}
$$
L'opérateur d'énergie, ou \textbf{hamiltonien}, est $H=\sum_{k \in \Z'} k\,:\! \psi_k\psi_k^*\!:$\,\,. Les opérateurs de \textbf{charge} et de \textbf{translation} sont pour leur part définis par
 $$C=\sum_{k\in \Z'} :\!\psi_k\psi_k^*\!: \qquad;\qquad R(\underline{D})=\underline{D}+1=(\underline{d_1+1})\wedge (\underline{d_2+1}) \wedge \cdots$$
 Ces divers opérateurs permettent une décomposition spectrale du produit extérieur infini. Ainsi, la base $(\underline{D})$ diagonalise la charge : 
 $$C(\underline{D})=(\card D_{+}-\card D_{-})\,\underline{D}\,,$$
où $D_{+}=D \cap \Z_{+}'$ et $D_{-}= \Z_{-}'\setminus D$. Le noyau de $C$, noté $\bigwedge_{0}^{\infty}V$, a pour base l'ensemble des coordonnées de descente\footnote{Plus généralement, un élément de base $\underline{D}$ peut être vu comme une fonction affine par morceaux égale à $|s|$ pour $s$ positif assez grand, et égale à $|s|+c$ pour $s$ négatif assez grand, $c$ désignant la charge. Les éléments de base sont donc des généralisations des partitions.} $\underline{\mathcal{D}(\lambda)}$ des partitions $\lambda \in \ym$. Dans ce qui suit, on notera $\underline{\lambda}=\underline{\mathcal{D}(\lambda)}$ ; en particulier, $\vacuum=\underline{\Z_{-}'}$. La conjugaison par $R$ translate les fermions libres : 
$$R^l\psi_kR^{-l}=\psi_{k+l} \qquad;\qquad R^l\psi_k^*R^{-l}=\psi_{k+l}^*\,.$$
Par conséquent, $R^lHR^{-l}=H+lC+l^2/2$ et $R^lCR^{-l}=C+l\,,$ de sorte que le sous-espace $\bigwedge_{c}^{\infty}V=\ker(C-c)$ est aussi $R^{c}(\bigwedge_{0}^{\infty}V)$. L'espace $\bigwedge^{\infty}V$ est la somme directe hilbertienne $\bigoplus_{c \in \Z} \bigwedge_{c}^{\infty}V$, et dans chaque sous-espace $\bigwedge_{c}^{\infty}V$, le vecteur $R^{c}(\vacuum)$ est un vecteur propre pour $H$ de valeur propre minimale.
\bigskip
\bigskip

Les \textbf{bosons libres} $\alpha_{n}=\sum_{k \in \Z'} :\!\psi_k \psi_{k+n}^*\!:$\, apportent un autre point de vue sur le produit extérieur infini. Ils satisfont les relations de commutation $$[\alpha_n,\alpha_m]=n\,\delta_{n,-m}\qquad;\qquad [\alpha_n,\psi_k]=\psi_{k-n}\qquad;\qquad [\alpha_n,\psi_k^*]=-\psi_{k+n}^*$$
 et ils engendrent l'\textbf{algèbre d'Heisenberg} (ou algèbre d'oscillation) $\mathscr{A}$. Fixons un paramètre $c \in \R$. L'algèbre $\mathscr{A}$ agit sur l'espace de Fock symétrique $S=\C[x_{1},x_{2},\ldots]=\bigotimes^{\infty}\C[x]$ :
 $$\alpha_{n>0}(P)=\frac{\partial P}{\partial x_n} \qquad;\qquad \alpha_{n<0}=(-n)x_{-n} \,P \qquad;\qquad \alpha_0(P)=c P \,.$$
 Autrement dit, les $\alpha_{n>0}$ agissent comme opérateurs d'annihilation et les $\alpha_{n<0}$ agissent comme opérateurs de création. On peut montrer que pour chaque entier $c \in \Z$, l'action de $\mathscr{A}$ sur $\bigwedge_{c}^{\infty}V$ est irréductible et équivalente à l'action $\mathscr{A} \curvearrowright S$ de paramètre $c$. L'action de $\mathscr{A}$ sur le produit extérieur infini est donc équivalente à la représentation $\mathscr{A} \curvearrowright \C[x_{1},x_{2},\ldots ; q,q^{-1}]=S\otimes \C[q,q^{-1}]$ donnée par
 $$\alpha_{n>0}(P)=\frac{\partial P}{\partial x_n}\qquad;\qquad \alpha_{n<0}=(-n)x_{-n}\,P\qquad;\qquad \alpha_0(P)=q\,\frac{\partial P}{\partial q} \,.$$
Le vecteur vide correspond au polynôme $1$ ; l'opérateur de charge correspond à $\alpha_{0}$, et l'opérateur de translation correspond à la multiplication par $q$.
\bigskip
\bigskip

Cette théorie du produit extérieur infini est complétée par la \textbf{correspondance boson-fermion} et la théorie des \textbf{opérateurs vertex}. Notons $\psi(z)$, $\psi^{*}(z)$ et $\alpha(z)$ les fonctions génératrices des fermions et des bosons :
$$\psi(z)=\sum_{n \in \Z'}z^n \,\psi_n \qquad;\qquad \psi^*(z)=\sum_{n \in \Z'} z^{-n}\,\psi_n^* \qquad;\qquad \alpha(z)=\sum_{n \in \Z} z^n\, \alpha_{-n}  \,.$$
Pour toute suite $t=(t_{1},t_{2},\ldots)$, les opérateurs vertex de paramètre $t$ sont définis par $\Gamma_+(t)=\exp ( \sum_{n\geq 1} t_n\, \alpha_n)$ et $\Gamma_-(t)=\exp ( \sum_{n\geq 1} t_n\, \alpha_{-n})=\Gamma_+(t)^*$. Avec ces conventions, si 
$$[z]=\left\{z,\frac{z^2}{2},\frac{z^3}{3},\ldots \right \},$$
alors la relation $\alpha(z)=\,:\!\psi(z)\psi^*(z)\!:\,$ peut être inversée :
$$\psi(z)=z^C\,R\,\Gamma_-([z])\,\Gamma_+([z^{-1}])^{-1}\qquad;\qquad \psi^*(z)=R^{-1}\,z^{-C}\,\Gamma_-([z])^{-1}\,\Gamma_+([z^{-1}])$$
\emph{cf.} le dernier chapitre de \cite{Kac94} pour une preuve de ces formules à l'aide d'algèbres de Lie infini-dimensionnelles. Cette correspondance boson-fermion\footnote{Tous ces résultats ont une interprétation physique naturelle. Ainsi, si $V$ désigne l'espace d'états (quantiques) d'une particule de type fermion, alors l'espace d'états d'une infinité de particules identiques est $\bigwedge^{\infty}V$, et les opérateurs $\psi_{k}$ et $\psi_{k}^{*}$ correspondent à la création et à l'annihilation d'un fermion en position $k$. Compte tenu des propriétés du produit extérieur, les configurations avec plusieurs particules dans le même état $k$ sont interdites ; c'est exactement ce que l'on attend de fermions. La correspondance boson-fermion est l'identification canonique entre $\bigwedge^{\infty}V$ et l'espace de Fock symétrique $S \otimes \C[q,q^{-1}]$ décrite précédemment ; ainsi, on peut traiter virtuellement le système de particules comme une infinité de bosons, et l'identification conserve la charge du système. De plus, à charge $c$ fixée, il existe un unique état $R^{c}(\vacuum)$ d'énergie minimale, et les autres états sont obtenus en appliquant l'opérateur vertex négatif sur cet état.} est complétée par les relations de commutation
\begin{align*}
\Gamma_{\pm}(t) \psi(z) &= H(t,z^{\pm 1})\,\psi(z) \Gamma_{\pm}(t)\\
\Gamma_{\pm}(t) \psi^*(z) &= H(t,z^{\pm 1})^{-1}\,\psi^*(z) \Gamma_{\pm}(t) \\
\Gamma_+(t)\,\Gamma_-(t')&=\exp\!\left({\sum}_{n=1}^{\infty} nt_nt_n'\right)\,\,\Gamma_-(t')\,\Gamma_+(t)
\end{align*} 
où $H(t,z)=\exp(\sum_{n=1}^\infty t_n\,z^n)$ est aussi $\sum_{n= 0}^\infty h_n(x)\,z^n$ si les $t_{n}=t_{n}(x)$ sont les paramètres de Miwa d'un ensemble de variables $x$. En supposant que c'est le cas, on en déduit l'action des opérateurs vertex sur les vecteurs de base de $\bigwedge^{\infty}V$ :
$$\Gamma_+(t)\left(\vacuum\right)=\vacuum\qquad ;\qquad \!\Gamma_-(t)\left(\vacuum\right)=\sum_{\lambda \in \ym} s_\lambda(x)\,\underline{\lambda}\qquad;\qquad\!\Gamma_-(t)\left(\underline{\mu}\right) = \sum_{\lambda \supset \mu} s_{\lambda/ \mu}(x)\,\underline{\lambda}\,.$$
Ici, $s_{\lambda/\mu}$ désigne la fonction de Schur gauche définie par la relation d'orthogonalité $$\forall \nu \in \ym,\,\,\,\langle s_\lambda \,| \,s_\mu s_\nu \rangle=\langle s_{\lambda / \mu} \,|\,s_\nu \rangle\,,$$ ou, de fa\c con équivalente, par la formule de Jacobi-Trudi $s_{\lambda / \mu}(x)=\det(h_{\lambda_i-\mu_j+j-i}(x))_{i,j}$, \emph{cf.} \cite[\S1.5]{Mac95}. Ces relations permettent de calculer les fonctions de corrélation d'une mesure de Schur de paramètre $(t,t')$. En effet, si l'on remarque que $\psi_k\psi_k^*(\underline{\lambda})$ vaut $\underline{\lambda}$ lorsque $k \in \mathcal{D}(\lambda)$ et $0$ sinon, alors :
\begin{align*}
\rho(A)&=\sum_{A \subset \mathcal{D}(\lambda)} Z^{-1}\,s_\lambda(a)\,s_\lambda(b) = Z^{-1}\,\scal{\left(\prod_{k \in A} \psi_k\psi_k^*\right) \,\Gamma_-(t)\,\vacuum}{\Gamma_-(t')\,\vacuum}\\
&=\scal{\left(Z^{-1}\,\Gamma_+(t')\Gamma_-(t)\right)\,\Gamma_-(t)^{-1}\,\left(\prod_{k \in A} \psi_k\psi_k^*\right)\,\Gamma_-(t)\,\vacuum}{\vacuum}\\
&=\scal{\left(\Gamma_-(t)\Gamma_+(t')\right)\,\Gamma_-(t)^{-1}\,\left(\prod_{k \in A} \psi_k\psi_k^*\right)\,\Gamma_-(t)\Gamma_+(t')^{-1}\Gamma_+(t')\,\vacuum}{\vacuum}\\
&=\scal{G\,\left(\prod_{k \in A}\psi_k\psi_k^*\right)\,G^{-1}\,\Gamma_+(t')\,\vacuum}{\Gamma_+(t)\,\vacuum}=\scal{\left(\prod_{k \in A} \Psi_k \Psi_k^*\right)\,\vacuum}{\vacuum}
\end{align*}
où $G=\Gamma_{+}(t')\Gamma_-(t)^{-1}$, $\Psi_k=G\psi_kG^{-1}$ et $\Psi_k^*=G\psi_k^*G^{-1}$. Finalement, en utilisant une forme de la formule de Wick, on conclut que les fonctions de corrélations sont les déterminants du noyau $K_{\mathcal{D}}(x,y)=\langle\Psi_{x}\Psi_{y}^{*}\,\vacuum|\vacuum\rangle$. De fa\c con plus explicite, si l'on introduit les fonctions génératrices
$$A(z,w)=\sum_{k,l \in \Z'}K_{\mathcal{D}}(k,l)\,z^{k}w^{-l}\quad;\quad J(x)=\frac{H(t,x)}{H(t',x^{-1})}=\exp\left(\sum_{n \geq 1} t_n \,x^n - \sum_{n \geq 1}t_n'\,x^{-n}\right)$$
alors $A(z,w)=\scal{G\psi(z)\psi^*(w)G^{-1}\,\vacuum}{\vacuum}$, et $G\psi(z)\psi^*(w)G^{-1} =(J(z)/J(w))\,\psi(z)\psi^*(w)$. Par conséquent :
$$A(z,w)=\frac{J(z)}{J(w)}\,\scal{\psi(z)\psi^*(w)\,\vacuum}{\vacuum} = \frac{J(z)}{J(w)}\, \sum_{k \in \Z_+'} (w/z)^k=\frac{\sqrt{zw}}{z-w}\,\frac{J(z)}{J(w)}$$
Comme $J(x)=\exp(T(x))$, on en déduit l'expression donnée dans l'énoncé du théorème \ref{schurmeasurearedeterminantal} par une intégrale de Cauchy sur deux contours. Alternativement, si l'on décompose $J(x)=\sum_{n \in \Z} J_{n}(t,t')\,x^{n}$ en série de Laurent, alors :
$$K_{\mathcal{D}}(x,y)=\sum_{m \in \Z_+'} J_{x+m}(t,t')\,J_{-y-m}(-t,-t')$$
et nous verrons plus loin que les fonctions $J_{n}$ sont des généralisations des fonctions de Bessel. Ainsi, les mesures de Schur donnent bien lieu à des processus ponctuels déterminantaux, et on dispose d'une expression du noyau $K_{\mathcal{D}}$ sous forme d'intégrale de Cauchy ou de série de fonctions. Plus loin, nous détaillerons cette expression dans le cas de la mesure de Plancherel (\S\ref{bdjdeterminantal}) et de la $q$-mesure de Plancherel (\S\ref{qplancherelprocess}).
\bigskip
\bigskip

Pour conclure ce chapitre, notons que les mesures de Schur ont originalement été étudiées pour comprendre l'asymptotique de surfaces aléatoires, en particulier les partitions tridimensionnelles, voir la figure \ref{planpartition}. 
\figcapt{
\psset{unit=1mm}\pspicture(0,0)(150,55)
\psframe*[linecolor=red!80!white](0,5)(8,21)
\psframe*[linecolor=red!65!white](8,5)(16,21)
\psframe*[linecolor=red!65!white](0,21)(8,29)
\psframe*[linecolor=orange!65!white](16,5)(24,13)
\psframe*[linecolor=orange!50!white](16,13)(24,21)
\psframe*[linecolor=orange!50!white](0,29)(8,37)
\psframe*[linecolor=orange!30!white](24,5)(40,13)
\psframe*[linecolor=orange!30!white](0,37)(8,45)
\psframe*[linecolor=orange!30!white](8,21)(16,37)
\psframe(0,5)(40,45)
\psline(8,5)(8,45)
\psline(16,5)(16,45)
\psline(24,5)(24,45)
\psline(32,5)(32,45)
\psline(0,13)(40,13)
\psline(0,21)(40,21)
\psline(0,29)(40,29)
\psline(0,37)(40,37)
\rput(4,9){$5$}
\rput(4,17){$5$}
\rput(4,25){$4$}
\rput(4,33){$2$}
\rput(4,41){$1$}
\rput(12,9){$4$}
\rput(12,17){$4$}
\rput(12,25){$1$}
\rput(12,33){$1$}
\rput(12,41){$0$}
\rput(20,9){$3$}
\rput(20,17){$2$}
\rput(20,25){$0$}
\rput(20,33){$0$}
\rput(20,41){$0$}
\rput(28,9){$1$}
\rput(28,17){$0$}
\rput(28,25){$0$}
\rput(28,33){$0$}
\rput(28,41){$0$}
\rput(36,9){$1$}
\rput(36,17){$0$}
\rput(36,25){$0$}
\rput(36,33){$0$}
\rput(36,41){$0$}

\psline*[linecolor=MidnightBlue!50!white](55,15)(59,12)(59,17)(55,20)(55,15)
\psline*[linecolor=MidnightBlue!50!white](63,15)(67,12)(67,17)(63,20)(63,15)
\psline*[linecolor=MidnightBlue!50!white](59,23)(63,20)(63,25)(59,28)(59,23)
\psline*[linecolor=MidnightBlue!50!white](63,31)(67,28)(67,38)(63,41)(63,31)
\psline*[linecolor=MidnightBlue!50!white](67,49)(67,44)(71,41)(71,46)(67,49)
\psline*[linecolor=MidnightBlue!50!white](75,18)(79,15)(79,25)(75,28)(75,38)(71,41)(71,26)(75,23)(75,18)
\psline*[linecolor=MidnightBlue!50!white](91,12)(91,17)(83,23)(83,18)(91,12)
\psline*[linecolor=MidnightBlue!50!white](83,33)(83,28)(79,31)(79,36)(83,33)
\psline*[linecolor=MidnightBlue!25!white](59,12)(63,15)(63,20)(59,17)(59,12)
\psline*[linecolor=MidnightBlue!25!white](67,12)(75,18)(75,23)(67,17)(67,12)
\psline*[linecolor=MidnightBlue!25!white](63,20)(71,26)(71,41)(67,38)(67,28)(63,25)(63,20)
\psline*[linecolor=MidnightBlue!25!white](71,41)(79,47)(79,52)(71,46)(71,41)
\psline*[linecolor=MidnightBlue!25!white](91,12)(95,15)(95,20)(91,17)(91,12)
\psline*[linecolor=MidnightBlue!25!white](79,15)(83,18)(83,23)(87,26)(87,36)(83,33)(83,28)(79,25)(79,15)
\psline*[linecolor=MidnightBlue!25!white](75,28)(79,31)(79,36)(83,39)(83,44)(75,38)(75,28)
\psline(55,15)(55,20)(59,23)(59,28)(63,31)(63,41)(67,44)(67,49)(75,55)(79,52)(79,47)(83,44)(83,39)(87,36)(87,26)(95,20)(95,15)(91,12)(83,18)(79,15)(75,18)(67,12)(63,15)(59,12)(55,15)
\psline(59,12)(59,17)(55,20)
\psline(63,15)(63,25)(59,28)
\psline(59,23)(67,17)(67,12)
\psline(59,17)(71,26)(71,46)(67,49)
\psline(67,44)(71,41)(67,38)(63,41)
\psline(67,38)(67,28)(63,31)
\psline(67,28)(63,25)
\psline(75,18)(75,23)(71,26)
\psline(67,17)(75,23)
\psline(79,52)(71,46)
\psline(79,47)(71,41)(75,38)(75,28)(79,25)(79,15)
\psline(91,12)(91,17)(83,23)(87,26)
\psline(91,17)(95,20)
\psline(83,23)(83,18)
\psline(87,36)(83,33)(83,28)(79,25)
\psline(83,44)(75,38)
\psline(83,39)(79,36)(83,33)
\psline(79,36)(79,31)(83,28)
\psline(79,31)(75,28)

\psline(110,15)(130,0)(150,15)(150,40)(130,55)(110,40)(110,15)
\psline(114,43)(114,23)\psline(114,17)(114,12)
\psline(118,46)(118,31)\psline(118,25)(118,15)
\psline(122,49)(122,44)\psline(122,38)(122,23)\psline(122,17)(122,12)
\psline(126,46)(126,26)\psline(126,20)(126,15)
\psline(130,49)(130,44)\psline(130,38)(130,18)
\psline(134,52)(134,47)\psline(134,41)(134,31)\psline(134,25)(134,15)
\psline(138,49)(138,39)\psline(138,33)(138,18)
\psline(142,46)(142,26)\psline(142,20)(142,15)
\psline(146,43)(146,23)\psline(146,17)(146,12)
\psline(126,3)(142,15)\psline(146,17)(150,20)
\psline(122,6)(138,18)\psline(142,20)(146,23)
\psline(118,9)(130,18)\psline(134,20)(142,26)
\psline(114,12)(118,15)\psline(122,17)(130,23)\psline(134,25)(142,31)
\psline(114,17)(126,26)\psline(130,28)(134,31)\psline(138,33)(142,36)
\psline(110,20)(114,23)\psline(118,25)(126,31)\psline(130,33)(138,39)
\psline(110,25)(118,31)\psline(122,33)(126,36)\psline(130,38)(138,44)
\psline(110,30)(118,36)\psline(122,38)(134,47)
\psline(110,35)(122,44)\psline(126,46)(134,52)
\psline(110,20)(114,17)\psline(118,15)(134,3)
\psline(114,23)(122,17)\psline(126,15)(138,6)
\psline(114,28)(118,25)\psline(122,23)(126,20)\psline(130,18)(142,9)
\psline(118,31)(122,28)\psline(126,26)(134,20)\psline(138,18)(146,12)
\psline(118,36)(122,33)\psline(126,31)(134,25)\psline(138,23)(146,17)
\psline(118,41)(122,38)\psline(126,36)(130,33)\psline(134,31)(138,28)\psline(142,26)(150,20)
\psline(122,44)(130,38)\psline(134,36)(138,33)\psline(142,31)(150,25)
\psline(122,49)(126,46)\psline(130,44)(134,41)\psline(138,39)(150,30)
\psline(126,52)(130,49)\psline(134,47)(150,35)
\endpspicture
}{Partition plane, diagramme tridimensionnel et pavage par des losanges.\label{planpartition}}{Partition plane, diagramme tridimensionnel et pavage par des losanges}

\noindent Ces diagrammes 3d sont caractérisés par l'une des assertions équivalentes suivantes :
\begin{enumerate}[(i)]
\item \'Ecrites dans un plan, les hauteurs du diagramme sont décroissantes suivantes les lignes et les colonnes --- on obtient ce qu'on appelle une \textbf{partition plane}.
\item Les suites d'entiers $\{\lambda(t)\}_{t \in \Z}$ obtenues en considérant les tranches diagonales sont des partitions, et elles sont entrelacées :
$$\cdots \prec \lambda(-2)\prec \lambda(-1) \prec \lambda(0) \succ \lambda(1) \succ \lambda(2) \succ \lambda(3) \succ \cdots$$
où $\lambda\succ \mu$ signifie que $\lambda_{1}\geq \mu_{1}\geq \lambda_{2}\geq \mu_{2}\geq \cdots$
\item Le diagramme 3d correspond au pavage d'un hexagone par des losanges.\vspace{2mm}
\end{enumerate}
On peut construire des suites aléatoires $\{\lambda(t)\}_{t \in \Z}$ de diagrammes de Young qui sont des analogues dépendant en temps des mesures de Schur, et qui correspondent à des partitions planes : il s'agit des \textbf{processus de Schur}, \emph{cf.} \cite{OR03,Bor10}. Dans ce cadre, les fonctions de corrélation des descentes des partitions $\lambda(t)$ sont déterminantales, y compris pour des descentes prises en des instants $t$ différents. Autrement dit, les corrélations des losanges du pavage sont déterminantales quelques soient les coordonnées choisies. Pour un choix adéquat des paramètres du processus de Schur, la partition plane aléatoire $\pi$ est tirée suivant une probabilité proportionnelle à $q^{\mathrm{vol}(\pi)}$ avec $q\in \,]0,1[$. En faisant tendre $q$ vers $1$ et en renormalisant correctement les partitions planes obtenues, on en déduit la forme limite des pavages par des losanges d'un grand hexagone, tous les pavages étant tirés équiprobablement, voir \cite{OR03} --- d'autre part, il est connu depuis les travaux de Kasteleyn que les corrélations des losanges d'un tel pavage sont données par des déterminants, voir \cite{Oko09} pour une généralisation de ce résultat.

\chapter{Permutations aléatoires et matrices aléatoires}\label{matrix}

Ce dernier chapitre d'introduction est consacré au lien entre la théorie asymptotique des représentations des groupes symétriques et la théorie des \textbf{grandes matrices aléatoires}. La ressemblance entre ces deux thèmes concerne en premier lieu les outils : ainsi, les espérances de polynômes traciaux jouent un rôle sensiblement identique aux observables de diagrammes du chapitre \ref{tool}, et d'autre part, les valeurs propres d'une matrice aléatoire constituent un processus ponctuel, qui est déterminantal pour le modèle du GUE (c'est-à-dire lorsque les entrées de la matrice sont des gaussiennes centrées indépendantes). Ainsi, des raisonnements semblables à ceux des paragraphes \ref{lskv} et \ref{schurmeasure} permettent de déterminer la limite des lois empiriques renormalisées des valeurs propres (loi de Wigner), ainsi que l'asymptotique <<~locale~>> au centre et au bord du spectre (noyau sinus et noyau d'Airy). Nous rappelons ceci dans la section \ref{gue}, en suivant pour l'essentiel \cite{Kon05} --- on renvoie à \cite{Meh04,AGZ09} pour un traitement plus exhaustif. Le noyau d'Airy intervient également comme limite des noyaux mis en jeu dans la description déterminantale des processus ponctuels associés aux mesures de Plancherel poissonisées. Ce fait est à l'origine de l'\textbf{équivalence de Baik-Deift-Johansson} (\cite{BDJ99,BDJ00}), qui peut être énoncée comme suit :
\begin{theorem}[Baik-Deift-Johansson, \cite{BOO00}]\label{bdj}
Soient $\lambda_{1}\geq \lambda_{2}\geq \cdots$ les plus grandes valeurs propres d'une matrice hermitienne du GUE, et 
$$Y_{i}=n^{2/3}\,\left(\frac{\lambda_{i}}{2\sqrt{n}} - 1\right)$$
leurs déviations renormalisées. On note d'autre part $\mu_{1} \geq \mu_{2} \geq \cdots$ les premières colonnes d'une partition aléatoire $\mu$ sous la mesure de Plancherel $M_{n}$, et 
$$X_{i}=n^{1/3}\,\left(\frac{\mu_{i}}{2\sqrt{n}}-1\right)$$
leurs déviations renormalisées. Par la correspondance RSK, il s'agit aussi des déviations renormalisées des tailles des plus longs sous-mots croissants d'une permutation aléatoire de taille $n$. Lorsque $n$ tend vers l'infini, les lois jointes de $(X_{1},\ldots,X_{k})$ et de $(Y_{1},\ldots,Y_{k})$ ont la même limite, et ce pour tout $k$.
\end{theorem}\bigskip

\noindent Pour $k=1$, la loi limite est la \textbf{loi de Tracy-Widom}, dont la densité est liée à l'équation différentielle de Painlevé $\text{II}$ (voir \cite{TW94}). Les sections \ref{bdjdeterminantal} et \ref{bdjgeom} sont consacrées à deux ébauches\footnote{K. Johansson a donné une troisième preuve qui utilise des polynômes orthogonaux, \emph{cf.} \cite{Joh01}.} de preuve du théorème \ref{bdj} : la première s'appuie sur les processus ponctuels déterminantaux, et la seconde est de nature géométrique --- cette approche due à Okounkov (\cite{Oko00}) a fait l'objet d'un exposé à l'Institut Henri Poincaré en décembre 2009. Comme mentionné dans l'introduction de cette partie, les résultats de ce chapitre sont présentés essentiellement à titre culturel\footnote{Le calcul du noyau associé à la mesure de Plancherel poissonisée de paramètre $\theta$ sera néanmoins repris et généralisé dans la section \ref{qplancherelprocess}.} ; néanmoins, la connection avec la théorie des matrices aléatoires est l'une des motivations principales de notre étude.\bigskip

\section[Théorie asymptotique des matrices aléatoires hermitiennes]{Théorie asymptotique des matrices aléatoires hermitiennes}\label{gue}

Soit $\Herm(n,\C)$ l'ensemble des matrices complexes de taille $n \times n$ et hermitiennes, \emph{i.e.}, $M_{ji}=\overline{M_{ij}}$ pour tout couple d'indices $(i,j)$. On note $M_{ii}=x_{ii} \in \R$, et $M_{i<j} = x_{ij}+\I y_{ij}$ avec $x_{ij} \in \R$, $y_{ij} \in \R$. La loi du \emph{Gaussian Unitary Ensemble} est la mesure de probabilité gaussienne sur $\Herm(n,\C)$ donnée par la densité
$$d\proba[M] = \prod_{i=1}^{n} \mathcal{N}_{(0,1)}(dx_{ii}) \,\prod_{i<j} \mathcal{N}_{(0,1/2)}(dx_{ij}) \,\prod_{i<j} \mathcal{N}_{(0,1/2)}(dy_{ij}) =\frac{1}{\widetilde{Z}_{n}}\,\E^{-\frac{\tr M^{2}}{2}}\,dM\, ,$$
où $\widetilde{Z}_{n}$ est une certaine constante de normalisation. Cette loi gaussienne est invariante pour l'action par conjugaison\footnote{Compte tenu de la décomposition de Cartan $\GL(n,\C) \simeq \Unit(n,\C) \times \Herm(n,\C)$, la donnée d'une mesure $\Unit(n,\C)$-invariante sur $\Herm(n,\C)$ est équivalente à la donnée d'une mesure bi-$\Unit(n,\C)$-invariante sur $\GL(n,\C)$, donc la loi du GUE fournit également une mesure gaussienne <<~symétrique~>> naturelle sur le groupe linéaire complexe.} de $\Unit(n,\C)$ sur $\Herm(n,\C)$. Notons $\lambda_{1} \geq \cdots \geq \lambda_{n}$ les valeurs propres d'une matrice du GUE ; c'est un élément de la chambre de Weyl $\overline{\Delta_{n}}$. En calculant le jacobien de l'application
\begin{align*}
\Unit(n,\C)/(\Unit(1,\C))^{n} \times \Delta_{n} &\to \Herm(n,\C) \\
[U]\,,\, (\lambda_{1}>\cdots>\lambda_{n}) &\mapsto U\,\mathrm{diag}(\lambda_{1},\ldots,\lambda_{n})\,U^{-1}
\end{align*}
d'image ouverte dense, on peut déterminer la loi jointe des valeurs propres sous la mesure du GUE :
$$d\proba[x_{1} \geq x_{2} \geq \cdots \geq x_{n}] = \frac{1}{Z_{n}}\,\Delta(x_{1}, \ldots,x_{n})^{2}\,\,\E^{-\frac{1}{2}\sum_{i=1}^{n}(x_{i})^{2}}\,dx_{1}\cdots dx_{n}$$
où $\Delta(x)=\prod_{i<j}\,(x_{j}-x_{i})$ est le Vandermonde de la matrice $((x_{i})^{j-1})_{1 \leq i,j \leq n}$, et $Z_{n}$ est une constante de normalisation.\bigskip
\bigskip

Les valeurs propres d'une matrice $M$ du GUE peuvent être regroupées en une mesure de probabilité (aléatoire) appelée \textbf{mesure empirique}, et définie par 
$$m_{M}=\frac{1}{n}\,\sum_{i=1}^{n} \delta_{\lambda_{i}}\,.$$
Si l'on note $\widetilde{m}(ds)=m(\sqrt{n}\,ds)$ la mesure empirique renormalisée, alors le résultat suivant est l'analogue pour les matrices aléatoires du théorème \ref{firstasymptoticplancherel} :
\begin{proposition}[Loi de Wigner, \cite{Wigner58}]
Au sens de la convergence en loi, $\widetilde{m} \to \mu_{\Omega}$ en probabilité, où
$$\mu_{\Omega}=\mathbb{1}_{[-2,2]}(s)\,\frac{\sqrt{4-s^{2}}}{2\pi}\,ds$$
est la loi de Wigner du demi-cercle --- c'est aussi la mesure de transition de la forme limite $\Omega$ pour les diagrammes sous la mesure de Plancherel.
\end{proposition}\bigskip

Une preuve possible de la loi de Wigner repose sur le développement topologique des espérances $\esper[\tr M^{k}]$ donné par la \textbf{formule de Wick} ; cette formule rentre dans le cadre des techniques de diagrammes de Feynman en théorie quantique des champs. Ces interprétations combinatoires des intégrales gaussiennes peuvent toutes être comprises à partir du cas trivial de la dimension $n=1$. Ainsi, pour calculer le $k$-ième moment
$$I_{k}=\frac{1}{\sqrt{2\pi}}\int_{\R} x^{k}\,\E^{-\frac{x^{2}}{2}}\,dx$$
d'une gaussienne, il est utile de perturber le champ $x^{2}/2$ par un terme linéaire $tx$ :
$$\frac{1}{\sqrt{2\pi}} \int_{\R} x^{k}\,\E^{-\frac{x^{2}}{2}+tx}\,dx=\frac{\partial^{k}}{\partial t^{k}}\left(\frac{1}{\sqrt{2\pi}} \int_{\R} \,\E^{-\frac{x^{2}}{2}+tx}\,dx\right) = \frac{\partial^{k}}{ \partial t^{k}}\left(\E^{\frac{t^{2}}{2}}\right).$$
Par conséquent, $I_{k}=\left.\frac{d^{k}}{dt^{k}}\right|_{t=0}\left(\E^{\frac{t^{2}}{2}}\right)$, et lors de ce calcul, chaque dérivation de $\E^{\frac{t^{2}}{2}}$ donne un terme $t$ qui doit être compensé par une autre dérivation, car les termes en $t^{l}$ restant à la fin du calcul sont annulés en évaluant en $t=0$. On en déduit :
$$I_{k}=\text{nombres de fa\c cons d'apparier $k$ points}=\begin{cases}(2p-1)!!&\text{si }k=2p\text{ est pair},\\0 &\text{sinon.}\end{cases}$$
\begin{example}
Détaillons le calcul de $I_{4}$ et l'interprétation combinatoire en termes d'appariements :
\psset{unit=1mm}
\pspicture(20,63)(170,64)
\psdots(90,54)(95,54)(100,54)(105,54)
\psdots(90,43)(95,43)(100,43)(105,43)
\pscurve{->}(90,43)(92.5,41)(95,40)
\psdots(90,32)(95,32)(100,32)(105,32)
\pscurve(90,32)(92.5,30.5)(95,32)
\rput(112.5,32){$+$}
\psdots(120,32)(125,32)(130,32)(135,32)
\pscurve{->}(120,32)(122.5,30)(125,29)
\pscurve{->}(125,32)(127.5,34)(130,35)

\psdots(90,20)(95,20)(100,20)(105,20)
\pscurve(90,20)(92.5,18.5)(95,20)
\pscurve{->}(100,20)(102.5,18)(105,17)
\rput(112.5,20){$+$}
\psdots(120,20)(125,20)(130,20)(135,20)
\pscurve{->}(120,20)(122.5,18)(125,17)
\pscurve(125,20)(127.5,21.5)(130,20)
\rput(142.5,20){$+$}
\psdots(150,20)(155,20)(160,20)(165,20)
\pscurve{->}(155,20)(157.5,22)(160,23)
\pscurve(150,20)(155,18)(160,20)

\psdots(90,8)(95,8)(100,8)(105,8)
\pscurve(90,8)(92.5,6.5)(95,8)
\pscurve(100,8)(102.5,6.5)(105,8)
\rput(112.5,8){$+$}
\psdots(120,8)(125,8)(130,8)(135,8)
\pscurve(120,8)(127.5,5.5)(135,8)
\pscurve(125,8)(127.5,9.5)(130,8)
\rput(142.5,8){$+$}
\psdots(150,8)(155,8)(160,8)(165,8)
\pscurve(155,8)(160,10)(165,8)
\pscurve(150,8)(155,6)(160,8)

\endpspicture
\begin{align*}
I_{4}&=\left.\frac{d^{4}}{dt^{4}}\right|_{t=0} \left(\E^{\frac{t^{2}}{2}}\right) \\
&=\left.\frac{d^{3}}{dt^{3}}\right|_{t=0} \left(t\,\E^{\frac{t^{2}}{2}}\right) \\
&=\left.\frac{d^{2}}{dt^{2}}\right|_{t=0} \left(\E^{\frac{t^{2}}{2}}+t^{2}\,\E^{\frac{t^{2}}{2}}\right) \qquad\qquad\qquad\qquad\qquad\qquad\qquad\qquad\qquad\qquad\qquad\\
&=\left.\frac{d}{dt}\right|_{t=0} \left(3t\,\E^{\frac{t^{2}}{2}}\right) \\
&=3
\end{align*}
\end{example}\bigskip
\bigskip\bigskip
\bigskip

Les mêmes techniques s'appliquent au calcul de l'espérance $\esper[\tr M^{k}]$, qui pour des raisons de parité est nulle si $k$ est impair. Ainsi :
$$\esper[\tr M^{2p}]=\sum_{i_{1},\ldots,i_{2p}} \esper[M_{i_{1}i_{2}}M_{i_{2}i_{3}}\cdots M_{i_{2p}i_{1}}] = \sum_{i_{1},\ldots,i_{2p}} \left.\frac{\partial^{2p}}{\partial M_{i_{1}i_{2}} \cdots\partial M_{i_{2p}i_{1}}} \right|_{M=0} \left(\E^{\frac{\tr M^{2}}{2}}\right).$$
Pla\c cons les entiers $1,2,\ldots,2p$ aux sommets d'un $2p$-gone régulier ; la dérivée partielle d'indices $i_{1}, \ldots,i_{2p}$ évaluée en $M=0$ est le nombre d'appariements des côtés du $2p$-gone qui sont compatibles avec l'orientation canonique, et tels que deux sommets $j_{1}$ et $j_{2}$ sont identifiés si et seulement si $i_{j_{1}}=i_{j_{2}}$. C'est donc le nombre de \textbf{cartes}\footnote{Une carte est la donnée d'une surface de Riemann et d'un graphe tracé dessus dont toutes les faces sont homéomorphes à des disques. On renvoie à \cite{LZ04} pour la théorie des cartes et ses multiples incarnations ; en particulier, le lien avec les modèles matriciels est évoqué dans le chapitre 3 de cet ouvrage.} à $p$ arêtes et $\card\{i_{1},i_{2},\ldots,i_{2p}\}$ sommets obtenues par recollement deux à deux des côtés d'un $2p$-gone. Le nombre de sommets s'écrit $p+1-2g$, où $g$ est le genre de la surface obtenue. Par conséquent :
\begin{proposition}[Formule de Wick, \cite{Hooft74}]
Si $k=2p$ est un entier pair et si $M$ est une matrice du GUE, alors
$$\esper[\tr M^{k}]=\sum_{g \geq 0} n^{p+1-2g}\,|\mathrm{Map}(k,g)|\,,$$
où $\mathrm{Map}(k,g)$ est l'ensemble des cartes planaires de genre $g$ obtenues par recollement des côtés d'un $k$-gone.
\end{proposition}
\begin{example}
À symétrie près, les cartes obtenues par recollement des côtés d'un hexagone sont toutes représentées sur la figure \ref{maphexagon}. L'espérance $\esper[\tr M^{6}]$ vaut donc $5n^{4}+10n^{2}$.

\figcapt{
\psset{unit=1mm}\pspicture(0,-2)(145,70) 
\psellipse[linecolor=violet,linewidth=0.6pt](50,15.9)(1.2,5.1)
\psframe*[linecolor=white,fillcolor=white](50,18.75)(60,19.5)
\psframe*[linecolor=white,fillcolor=white](50,16)(60,16.75)
\psframe*[linecolor=white,fillcolor=white](50,13.25)(60,14)
\psdots[linecolor=violet](49,13.9)(42,16)
\psellipse[linecolor=violet,linewidth=0.6pt](50,10)(7.5,3.9)
\psframe*[linecolor=white,fillcolor=white](40,5)(60,10)
\psellipse[linecolor=violet,linewidth=0.6pt,linestyle=dashed](50,10)(7.5,3.9)
\pscurve[linecolor=violet]{<->}(15,5)(20,10.5)(24,12.5)
\pscurve[linecolor=violet]{<->}(15,15)(20,9.5)(24,7.5)
\pscurve[linecolor=violet]{<->}(6,7.5)(7,10)(6,12.5)
\psellipse[linewidth=0.8pt](50,10)(15,11)
\psarc[linewidth=0.8pt](50,19.5){12}{225}{315}
\psarc[linewidth=0.8pt](50,-5){16}{65}{115}
\pscurve[linecolor=violet,linewidth=0.6pt](49,13.9)(45,15.5)(42,16)
\rput(10,15){$1$} \psline(12,15)(18,15) \rput(20,15){$2$} \psline(21.5,14)(26.5,11) \rput(28,10){$3$}
\psline(21.5,6)(26.5,9) \rput(20,5){$4$} \psline(12,5)(18,5) \rput(10,5){$5$} \psline(8.5,6)(3.5,9) \rput(2,10){$6$} \psline(8.5,14)(3.5,11)
\rput(30,54){$\times 3$}
\rput(30,29){$\times 2$}
\rput(30,4){$\times 6$}
\rput(105,44){$\times 3$}
\rput(105,19){$\times 1$}
\pscurve[linecolor=violet]{<->}(15,30)(20,32)(24,32.5)
\pscurve[linecolor=violet]{<->}(15,40)(20,38)(24,37.5)
\pscurve[linecolor=violet]{<->}(6,32.5)(7,35)(6,37.5)
\rput(10,40){$1$} \psline(12,40)(18,40) \rput(20,40){$2$} \psline(21.5,39)(26.5,36) \rput(28,35){$3$}
\psline(21.5,31)(26.5,34) \rput(20,30){$4$} \psline(12,30)(18,30) \rput(10,30){$5$} \psline(8.5,31)(3.5,34) \rput(2,35){$6$} \psline(8.5,39)(3.5,36)
\pscircle[linewidth=0.8pt](50,35){10}
\psdots[linecolor=violet](54,37)(49,34)(44,38)(47,29)
\pscurve[linecolor=violet,linewidth=0.6pt](47,29)(48,31)(49,34)
\pscurve[linecolor=violet,linewidth=0.6pt](54,37)(51,36)(49,34)
\pscurve[linecolor=violet,linewidth=0.6pt](44,38)(46,36)(49,34)
\psline[linecolor=violet]{<->}(15,65)(15,55)
\pscurve[linecolor=violet]{<->}(6,57.5)(7,60)(6,62.5)
\pscurve[linecolor=violet]{<->}(24,57.5)(23,60)(24,62.5)
\rput(10,65){$1$} \psline(12,65)(18,65) \rput(20,65){$2$} \psline(21.5,64)(26.5,61) \rput(28,60){$3$}
\psline(21.5,56)(26.5,59) \rput(20,55){$4$} \psline(12,55)(18,55) \rput(10,55){$5$} \psline(8.5,56)(3.5,59) \rput(2,60){$6$} \psline(8.5,64)(3.5,61)
\pscircle[linewidth=0.8pt](50,60){10}
\psdots[linecolor=violet](53,55)(48,59)(44,63)(47,54)
\pscurve[linecolor=violet,linewidth=0.6pt](47,54)(50,54.35)(53,55)
\pscurve[linecolor=violet,linewidth=0.6pt](47,54)(47.5,56.5)(48,59)
\pscurve[linecolor=violet,linewidth=0.6pt](44,63)(46.25,61)(48,59)
\psellipse[linecolor=violet,linewidth=0.6pt](125,55.9)(1.2,5.1)
\psframe*[linecolor=white,fillcolor=white](125,58.75)(135,59.5)
\psframe*[linecolor=white,fillcolor=white](125,56)(135,56.75)
\psframe*[linecolor=white,fillcolor=white](125,53.25)(135,54)
\psdots[linecolor=violet](124,53.9)(125.95,57.5)
\psellipse[linecolor=violet,linewidth=0.6pt](125,50)(7.5,3.9)
\psframe*[linecolor=white,fillcolor=white](115,45)(135,50)
\psellipse[linecolor=violet,linewidth=0.6pt,linestyle=dashed](125,50)(7.5,3.9)
\psdots[linecolor=violet](122,28.6)(131,27)
\pscurve[linecolor=violet,linewidth=0.6pt](122,28.6)(124,35.8)(128,25.8)(131,27)
\psframe*[linecolor=white,fillcolor=white](124,34.3)(128,35.3)
\psframe*[linecolor=white,fillcolor=white](124,31.55)(128,32.55)
\psframe*[linecolor=white,fillcolor=white](124,28.8)(128,29.8)
\psframe*[linecolor=white,fillcolor=white](124,26.05)(128,27.05)
\psellipse[linecolor=violet,linewidth=0.6pt](125,25)(7.5,3.9)
\psframe*[linecolor=white,fillcolor=white](115,20)(135,25)
\psellipse[linecolor=violet,linewidth=0.6pt,linestyle=dashed](125,25)(7.5,3.9)
\psline[linecolor=violet]{<->}(90,30)(90,20)
\psline[linecolor=violet]{<->}(81.5,27.5)(98.5,22.5)
\psline[linecolor=violet]{<->}(81.5,22.5)(98.5,27.5)
\rput(85,30){$1$} \psline(87,30)(93,30) \rput(95,30){$2$} \psline(96.5,29)(101.5,26) \rput(103,25){$3$}
\psline(96.5,21)(101.5,24) \rput(95,20){$4$} \psline(87,20)(93,20) \rput(85,20){$5$} \psline(83.5,21)(78.5,24) \rput(77,25){$6$} \psline(83.5,29)(78.5,26)
\psellipse[linewidth=0.8pt](125,50)(15,11)
\psellipse[linewidth=0.8pt](125,25)(15,11)
\psarc[linewidth=0.8pt](125,59.5){12}{225}{315}
\psarc[linewidth=0.8pt](125,35){16}{65}{115}
\psarc[linewidth=0.8pt](125,34.5){12}{225}{315}
\psarc[linewidth=0.8pt](125,10){16}{65}{115}
\psline[linecolor=violet]{<->}(90,55)(90,45)
\pscurve[linecolor=violet]{<->}(81.5,52.5)(90,52)(98.5,52.5)
\pscurve[linecolor=violet]{<->}(81.5,47.5)(90,48)(98.5,47.5)
\rput(85,55){$1$} \psline(87,55)(93,55) \rput(95,55){$2$} \psline(96.5,54)(101.5,51) \rput(103,50){$3$}
\psline(96.5,46)(101.5,49) \rput(95,45){$4$} \psline(87,45)(93,45) \rput(85,45){$5$} \psline(83.5,46)(78.5,49) \rput(77,50){$6$} \psline(83.5,54)(78.5,51)
\endpspicture
}{Cartes obtenues par recollement des côtés d'un hexagone. Les symétries sont indiquées en dessous de chaque carte.\label{maphexagon}}{Cartes obtenues par recollement des côtés d'un hexagone}
\end{example}

De fa\c con générale, on obtient une carte planaire ($g=0$) si et seulement si le graphe obtenu est un arbre, et il y a $C_{p}$ fa\c cons différentes de faire ceci. Or, les cartes planaires correspondent au terme dominant dans le développement topologique de $\esper[\tr M^{2p}]$, donc :
$$\esper[\tr M^{2p}] = \esper[\widetilde{m}(x^{2p})]\,n^{p+1} = C_{p}\,n^{p+1} +O(n^{p-1}).$$
Comme la mesure $\mu_{\Omega}$ a pour moments pairs les nombres de Catalan et a ses moments impairs nuls, on en déduit que $\esper[\widetilde{m}(x^{k})]$ tend vers $\mu_{\Omega}(x^{k})$ pour tout entier $k$. De plus, la formule de Wick se généralise au cas d'un produit de traces de puissances de $M$ :
$$\esper\left[\prod_{l=1}^{s} \tr M^{k_{l}}\right] = \sum_{S} \,n^{\frac{k_{1}+\cdots +k_{s}}{2}+\chi(S)-s} \,|\mathrm{Map}(k_{1},\ldots,k_{s};S)|\,,$$
la somme étant effectuée sur toutes les classes d'homéomorphismes de surfaces compactes $S$ (éventuellement non connexes), et $\mathrm{Map}(k_{1},\ldots,k_{s};S)$ désignant l'ensemble des cartes tracées sur une surface de type $S$ obtenues par recollement d'un $k_{1}$-gone, d'un $k_{2}$-gone, etc., et d'un $k_{s}$-gone. \bigskip

\begin{example}
À symétrie près, les recollements possibles des côtés de deux carrés sont tous représentés sur la figure \ref{maptwosquare}. L'espérance $\esper[(\tr M^{4})^{2}]$ vaut donc $4n^{6}+24n^{4}+77n^{2}$.

\figcapt{\psset{unit=1mm}\pspicture(0,-12)(100,47)
\psframe(0,-10)(5,-5)
\psframe(15,-10)(20,-5)
\psline[linecolor=blue]{<->}(5,-7.5)(15,-7.5)
\pscurve[linecolor=blue]{<->}(2.5,-5)(22,-5)(20,-7.5)
\pscurve[linecolor=blue]{<->}(0,-7.5)(-2,-10)(17.5,-10)
\pscurve[linecolor=blue]{<->}(2.5,-10)(9,-10)(11,-5)(17.5,-5)
\rput(35,-7.5){$\times 16,\quad\mathbb{T}^{2}$}

\psframe(0,5)(5,10)
\psframe(15,5)(20,10)
\psline[linecolor=blue]{<->}(5,7.5)(15,7.5)
\pscurve[linecolor=blue]{<->}(2.5,5)(10,3)(17.5,5)
\pscurve[linecolor=blue]{<->}(2.5,10)(22,10)(20,7.5)
\pscurve[linecolor=blue]{<->}(0,7.5)(-2,10)(17.5,10)
\rput(35,7.5){$\times 16,\quad\mathbb{T}^{2}$}

\psframe(0,20)(5,25)
\psframe(15,20)(20,25)
\psline[linecolor=blue]{<->}(0,22.5)(2.5,25)
\psline[linecolor=blue]{<->}(17.5,20)(17.5,25)
\psline[linecolor=blue]{<->}(5,22.5)(15,22.5)
\pscurve[linecolor=blue]{<->}(2.5,20)(22,20)(20,22.5)
\rput(35,22.5){$\times 32,\quad\mathbb{T}^{2}$}

\psframe(0,30)(5,35)
\psframe(10,30)(15,35)
\psline[linecolor=blue]{<->}(0,32.5)(2.5,35)
\psline[linecolor=blue]{<->}(2.5,30)(5,32.5)
\psline[linecolor=blue]{<->}(10,32.5)(12.5,35)
\psline[linecolor=blue]{<->}(12.5,30)(15,32.5)
\rput(32.5,32.5){$\times 4,\quad\mathbb{S}^{2}\mathrel{\#}\mathbb{S}^{2}$}

\psframe(0,40)(5,45)
\psframe(10,40)(15,45)
\psline[linecolor=blue]{<->}(2.5,40)(2.5,45)
\psline[linecolor=blue]{<->}(0,42.5)(5,42.5)
\psline[linecolor=blue]{<->}(12.5,40)(12.5,45)
\psline[linecolor=blue]{<->}(15,42.5)(10,42.5)
\rput(32.5,42.5){$\times 1,\quad\mathbb{T}^{2}\mathrel{\#}\mathbb{T}^{2}$}

\psframe(70,-10)(75,-5)
\psframe(85,-10)(90,-5)
\psline[linecolor=blue]{<->}(75,-7.5)(85,-7.5)
\pscurve[linecolor=blue]{<->}(70,-7.5)(68,-10)(80,-14)(92,-10)(90,-7.5)
\pscurve[linecolor=blue]{<->}(72.5,-10)(79,-10)(81,-5)(87.5,-5)
\pscurve[linecolor=blue]{<->}(72.5,-5)(79,-5)(81,-10)(87.5,-10)
\rput(105,-7.5){$\times 4,\quad\mathbb{T}^{2}$}

\psframe(70,5)(75,10)
\psframe(85,5)(90,10)
\pscurve[linecolor=blue]{<->}(72.5,5)(80,3)(87.5,5)
\pscurve[linecolor=blue]{<->}(72.5,10)(80,12)(87.5,10)
\psline[linecolor=blue]{<->}(75,7.5)(85,7.5)
\pscurve[linecolor=blue]{<->}(70,7.5)(68,5)(80,1)(92,5)(90,7.5)
\rput(105,7.5){$\times 4,\quad\mathbb{S}^{2}$}

\psframe(70,20)(75,25)
\psframe(85,20)(90,25)
\psline[linecolor=blue]{<->}(72.5,20)(72.5,25)
\psline[linecolor=blue]{<->}(87.5,20)(87.5,25)
\psline[linecolor=blue]{<->}(75,22.5)(85,22.5)
\pscurve[linecolor=blue]{<->}(70,22.5)(69,20)(80,17.5)(91,20)(90,22.5)
\rput(105,22.5){$\times 8,\quad\mathbb{T}^{2}$}

\psframe(70,30)(75,35)
\psframe(85,30)(90,35)
\psline[linecolor=blue]{<->}(70,32.5)(72.5,35)
\psline[linecolor=blue]{<->}(87.5,35)(90,32.5)
\psline[linecolor=blue]{<->}(75,32.5)(85,32.5)
\pscurve[linecolor=blue]{<->}(72.5,30)(80,28)(87.5,30)
\rput(105,32.5){$\times 16,\quad\mathbb{S}^{2}$}

\psframe(70,40)(75,45)
\psframe(80,40)(85,45)
\psline[linecolor=blue]{<->}(70,42.5)(72.5,45)
\psline[linecolor=blue]{<->}(72.5,40)(75,42.5)
\psline[linecolor=blue]{<->}(80,42.5)(85,42.5)
\psline[linecolor=blue]{<->}(82.5,40)(82.5,45)
\rput(102,42.5){$\times 4,\quad\mathbb{S}^{2}\mathrel{\#}\mathbb{T}^{2}$}

\endpspicture}{Recollements possibles des côtés de deux carrés. Les symétries sont indiquées à coté de chaque carte, ainsi que le type de la surface obtenue.\label{maptwosquare}}{Cartes obtenues par recollement des côtés de deux carrés}
\end{example}

On peut ainsi calculer $\esper[(\widetilde{m}(x^{k}))^{2}]$, et pour $k=2p$ pair, le terme dominant correspond aux cartes à deux composantes connexes planaires --- c'est ce qui maximise la caractéristique d'Euler $\chi(S)$. Ainsi, $\esper[(\widetilde{m}(x^{2p}))^{2}]\simeq (C_{p})^{2}$, donc :
$$\esper[(\widetilde{m}(x^{2p}) - \mu_{\Omega}(x^{2p}))^{2}] \to 0\,.$$
Ceci permet de conclure quant à la convergence de $\widetilde{m}$ vers la loi du demi-cercle $\mu_{\Omega}$. L'interprétation géométrique des espérances de polynômes traciaux $\esper[\prod_{l=1}^{s} \tr M^{k_{l}}] $ jouera de nouveau un rôle essentiel dans une preuve de l'équivalence de Baik-Deift-Johansson, voir la section \ref{bdjgeom}. \bigskip
\bigskip

Toujours dans le cadre du GUE, la déviation renormalisée $\sqrt{n}\,(\widetilde{m}-\mu_{\Omega})$ peut être décrite par un théorème central limite analogue à celui de la section \ref{sniady}, voir \cite{Joh98}. D'autre part, les valeurs propres $\lambda_{1},\ldots,\lambda_{n}$ peuvent être envisagées comme un processus ponctuel sur $\R$ de densité 
$$d\mathbb{Q}[x_{1},\ldots,x_{n}]=\frac{1}{n!\,Z_{n}} \,\Delta(x_{1},\ldots,x_{n})^{2}\,\E^{-\frac{1}{2}\sum_{i=1}^{n}(x_{i})^{2}}\,dx_{1}\cdots dx_{n}\,,$$
c'est-à-dire que la $k$-ième fonction de corrélation du processus s'écrit :
$$\rho^{k}(x_{1},\ldots,x_{k})=\frac{n!}{(n-k)!}\,\int_{\R^{n-k}}\frac{d\mathbb{Q}(x_{1},\ldots,x_{k},x_{k+1},\ldots,x_{n})}{dx}\,dx_{k+1}\cdots dx_{n}$$
Dans ce contexte, l'identité $$\esper\left[\prod_{i=1}^{n}(1+f(x_{i})) \right]=\sum_{k=0}^{n}\frac{1}{k!} \int_{\R^{k}} \rho^{k}(x_{1},\ldots,x_{k})\,f(x_{1})\cdots f(x_{k})\,dx_{1}\cdots dx_{k}$$
permet de nombreux calculs de probabilités ; par exemple, en prenant $f=-\mathbb{1}_{A}$ avec $A\subset \R$, on obtient 
$$\proba[\forall i,\,\,x_{i} \notin A] = \sum_{k=0}^{\infty} \frac{(-1)^{k}}{k!}\,\int_{A^{k}} \rho^{k}(x_{1},\ldots,x_{k})\,dx_{1}\cdots dx_{k}\,.$$
Envisageons plus généralement un processus ponctuel sur $\R$ dont la densité s'écrit
$$d\mathbb{Q}(x) = Z^{-1}\,\Delta(x)^{2}\,\E^{-\frac{1}{2}\sum_{i=1}^{n}R(x_{i})}\,dx\,,$$ où $R(x)$ est une fonction suffisamment grande à l'infini. Soit $(\pi_{k})_{k \in \N}$ la famille des polygônes orthogonaux associés à la mesure $\exp(-R(x)/2)\,dx$, c'est-à-dire que $\pi_{k}(x)=x^{k}+\cdots$ et 
$$\int_{\R} \pi_{i}(x) \pi_{j}(x)\,\E^{-\frac{R(x)}{2}}\,dx=c_{i}\,c_{j}\,\delta_{ij}\,.$$
Les fonctions $\phi_{k}(x)=(c_{k})^{-1}\,\exp(-R(x)/4)\,\pi_{k}(x)$ forment une base orthonormale de $\leb^{2}(\R)$, et d'autre part, $\Delta(x_{1},\ldots,x_{k})=\det((\pi_{i}(x_{j}))_{1\leq i,j\leq k})$ pour tout entier $k$. Notons $K_{n}(x,y)$ le noyau $\sum_{i=0}^{n-1}\phi_{i}(x)\,\phi_{i}(y)$. La structure particulière de ce noyau (dit reproduisant) implique le caractère déterminantal du processus ponctuel considéré :
\begin{proposition}[Polynômes orthogonaux et processus ponctuels déterminantaux, \cite{Meh04}]
Dans le cadre précédent, le processus ponctuel $\{x_{1},\ldots,x_{n}\}$ a pour densité
$$d\mathbb{Q}(x)=\frac{1}{n!}\,\det((K_{n}(x_{i},x_{j}))_{1 \leq i,j \leq n})\,dx\,,$$
et les fonctions de corrélations s'écrivent $\rho^{k}(x_{1},\ldots,x_{k})=\det((K_{n}(x_{i},x_{j}))_{1\leq i,j\leq k})$. De plus, le noyau $K_{n}(x,y)$ s'écrit à l'aide de la formule de Christoffel-Darboux :
$$K_{n}(x,y)=\frac{c_{n-1}}{c_{n}}\,\frac{\phi_{n}(x)\phi_{n-1}(y)-\phi_{n}(y)\phi_{n-1}(x)}{x-y}\,.$$
\end{proposition}\bigskip

\begin{example}
Les valeurs propres d'une matrice du GUE renormalisées par un facteur $\sqrt{n}$ correspondent au cas où $R(x)=nx^{2}$. Le polynôme $\pi_{k}$ s'écrit $\pi_{k}(x)=H_{k}(\sqrt{n}x)/n^{k/2}$, où $H_{k}$ est le $k$-ième polynôme de Hermite. L'asymptotique du processus ponctuel associé aux valeurs propres renormalisées se déduit donc dans ce cas des propriétés à l'infini des fonctions de Hermite (voir \cite{Sze39}, en particulier le chapitre 8 et les formules de Plancherel-Rotach). Précisons quelque peu ces résultats. Si $B$ est une partie (borélienne) de $\R$, alors 
$$\card \{i \,\,|\,\,x_{i} \in B\} \sim n\,\mu_{\Omega}(B)\,$$
lorsque $n$ tend vers l'infini. En particulier, si $B=B_{n,u}=[u-l/2n,u+l/2n]$ est un intervalle centré en $u$ avec $u \in \,]-2,2[$ point à l'intérieur du spectre limite, alors
$$\card \{i \,\,|\,\,x_{i} \in B_{n,u}\} \sim \rho_{\Omega}(u)=\frac{\sqrt{4-u^{2}}}{2\pi}\,.$$
Par contre, si $u=\pm 2$ est un point du bord, alors il faut choisir une renormalisation de l'intervalle $B$ en $n^{2/3}$, car :
$$\card \{i \,\,|\,\,x_{i} \in [2-\eps_{n},2]\} \sim \frac{n}{2\pi}\, \int_{2-\eps_{n}}^{2}\sqrt{4-x^{2}}\,dx \sim \frac{2}{3\pi}\,n(\eps_{n})^{3/2}\,.$$
Pour comprendre l'asymptotique <<~locale~>> du processus ponctuel des valeurs propres d'une matrice du GUE, on introduit donc deux noyaux renormalisés :
\begin{align*}\widetilde{K}_{n}^{\mathrm{bulk}}(x,y)&=\frac{1}{n\,\rho_{\Omega}(u)}\,\,K_{n}\left(u+\frac{x}{n \,\rho_{\Omega}(u)},u+\frac{y}{n \,\rho_{\Omega}(u)}\right); \\
\widetilde{K}_{n}^{\mathrm{edge}}(x,y)&=\frac{1}{n^{2/3}} \,\,K_{n}\left(2+\frac{x}{n^{2/3}},2+\frac{y}{n^{2/3}}\right).\end{align*}
\end{example}
\begin{proposition}[Asymptotique du processus ponctuel des valeurs propres, \cite{TW94}]\label{asymptoticgue}
Lorsque $n$ tend vers l'infini, le noyau renormalisé à l'intérieur du spectre converge vers le noyau sinus :
$$\forall u,\,\,\forall x,y,\,\,\,\widetilde{K}_{n}^{\mathrm{bulk}}(x,y) \to \frac{\sin \pi(x-y)}{\pi (x-y)}\,.$$
De même, le noyau renormalisé au bord du spectre converge vers le noyau d'Airy :
$$\forall x,y,\,\,\,\widetilde{K}_{n}^{\mathrm{edge}}(x,y) \to \frac{\mathrm{Ai}(x)\,\mathrm{Ai}'(y)-\mathrm{Ai}'(x)\,\mathrm{Ai}(y)}{x-y}\,,$$
où $\mathrm{Ai}(x)=\frac{1}{\pi}\,\int_{0}^{\infty} \cos(tx+t^{3}/3)\,dt$ est la fonction d'Airy, c'est-à-dire l'unique solution de $y''-xy=0$ qui ait pour asymptotique $$\mathrm{Ai}(x)\sim_{x \to \infty}\frac{\E^{-2x^{3/2}/3}}{2\sqrt{\pi}\,x^{1/4}}\,.$$
\end{proposition}\bigskip\bigskip

En combinant ces résultats et la formule donnée pour la probabilité $\proba[\forall i, \,\,x_{i} \notin A]$, C. Tracy et H. Widom sont parvenus à préciser l'asymptotique de la plus grande valeur propre $\lambda_{\mathrm{max}}$ d'une matrice du GUE. Ainsi,
\begin{align*}\proba[2Y_{1} \leq a]&=\proba\left[n^{2/3}\,\left(\frac{\lambda_{\mathrm{max}}}{\sqrt{n}}-2\right)\leq a\right]=\sum_{k=0}^{N} \frac{(-1)^{k}}{k!}\int_{]a,\infty[^{k}} \det \left(\widetilde{K}_{n}^{\mathrm{edge}}(x_{i},x_{j})\right)dx^{k}\\
&\to\,\,\, \sum_{k=0}^{\infty} \frac{(-1)^{k}}{k!}\int_{]a,\infty[^{k}} \det \left(K^{\mathrm{Airy}}(x_{i},x_{j})\right)dx^{k}=\det\left(\id-\left.\mathcal{K}^{\mathrm{Airy}}\right|_{\leb^{2}(]a,\infty[)}\right).
\end{align*}
La loi obtenue, dite \textbf{loi de Tracy-Widom}, a pour fonction de répartition $$F(x)=\exp(-\int_{x}^{\infty} (s-x)\,q^{2}(s)\,ds)\,,$$
 où $q$ est l'unique solution de l'\textbf{équation différentielle de Painlevé II} $$q''(x)=xq(x)+2q^{3}(x)$$
qui ait pour asymptotique $q(x) \sim\mathrm{Ai}(x)$. En particulier, la distribution de la plus grande valeur propre du GUE est liée à la théorie des équations de Painlevé et à la hiérarchie de Toda, voir les articles \cite{FW01,FW02}. L'objectif de ce chapitre est de montrer que les premières lignes (ou colonnes) d'une partition sous la mesure de Plancherel ont la même loi asymptotique.\bigskip

\section[\'Equivalence de Baik-Deift-Johansson : l'approche déterminantale]{\'Equivalence de Baik-Deift-Johansson : l'approche\\ déterminantale}\label{bdjdeterminantal}
Dans la section \ref{schurmeasure}, nous avons vu que si $M_{z,\mathcal{B}(t,\zeta)}$ est la mesure sur les partitions obtenue en modulant le système cohérent des $z$-mesures $(M_{z,n})_{n \in \N}$ par la loi binomiale négative $\mathcal{B}(t,\zeta)$, alors les processus ponctuels $\mathcal{F}_{\star}(M_{z,\mathcal{B}(t,\zeta)})$ et $\mathcal{D}_{\star}(M_{z,\mathcal{B}(t,\zeta)})$ sont déterminantaux. De plus, on peut donner une expression explicite des noyaux correspondants en termes de fonctions hypergéométriques (voir le théorème 3.3 de \cite{BO00}). Lorsque $z$ tend vers l'infini et $\zeta=\theta/t$ tend vers $0$, la loi $M_{z,n}$ tend vers la mesure de Plancherel $M_{n}$, et la loi binomiale négative $\mathcal{B}(t,\zeta)$ tend vers la loi de Poisson $\mathcal{P}(\theta)$. Par suite, la \textbf{mesure de Plancherel poissonisée}
$$M_{\mathcal{P}(\theta)}=\mathcal{P}(\theta)(|\lambda|)\,\,M_{|\lambda|}(\lambda)=\theta^{|\lambda|}\,\E^{-\theta}\left(\frac{\dim \lambda}{|\lambda|!}\right)^{2}$$
correspond à un processus ponctuel déterminantal $\mathcal{D}_{\star}(M_{\mathcal{P}(\theta)})$ dont le noyau s'obtient par passage à la limite du cas des $z$-mesures, \emph{cf.} \cite[théorèmes 1 et 2]{BOO00}.\bigskip\bigskip

Une autre approche possible pour le calcul du noyau de ce processus est l'utilisation du formalisme des mesures de Schur ; en effet, $M_{\mathcal{P}(\theta)}$ est la mesure de Schur de paramètres $t=t'=(\sqrt{\theta},0,0,\ldots)$, voir page \pageref{exampleschurmeasure}. Si $r$ est un paramètre réel, on rappelle que la \textbf{fonction de Bessel} de première espèce $J_{r}$ est la série de fonctions
$$J_{r}(z)=\sum_{n=0}^{\infty} \frac{(-1)^{n}}{n!\,\,\Gamma(r+n+1)}\,\left(\frac{z}{2}\right)^{r+2n}\,, $$
de sorte que $J_{m}(z)$ est le coefficient de $x^{m}$ dans le développement en série de Laurent de $\E^{zx/2}\,\E^{-z/2x}$, \emph{cf.} \cite[chapitre 2]{Wat44}. Or, pour les paramètres $t$ et $t'$ précédemment décrits, la fonction $J(x)=\exp(T(x))$ décrite dans la section \ref{schurmeasure} est exactement $\E^{zx/2}\,\E^{-z/2x}$ avec $z=2\sqrt{\theta}$. Par conséquent,
\begin{align*}
K_{\mathcal{D}_{\star}(M_{\mathcal{P}(\theta)})}(x,y)&=\sum_{n=0}^{\infty} J_{x+n+1/2}(2\sqrt{\theta})\,J_{-y-n-1/2}(-2\sqrt{\theta})\\
&=\sum_{n=0}^{\infty} J_{x+n+1/2}(2\sqrt{\theta})\,J_{y+n+1/2}(2\sqrt{\theta})\\
&=\frac{\sqrt{\theta}\,\left(J_{x-1/2}(2\sqrt{\theta})\,J_{y+1/2}(2\sqrt{\theta})-J_{x+1/2}(2\sqrt{\theta})\,J_{y-1/2}(2\sqrt{\theta})\right)}{x-y}\,,\end{align*}
la seconde égalité découlant de l'invariance de la fonction $J$ par la transformation $x\mapsto -x^{-1}$, et la troisième égalité étant une conséquence des relations de Lommel pour un produit $J_{\mu}(z)\,J_{\nu}(z)$ de fonctions de Bessel, voir \cite[\S5.4]{Wat44}. En utilisant
les formules de Debye pour l'asymptotique des fonctions de Bessel lorsque l'argument complexe $z$ tend vers l'infini (\emph{cf.} \cite[chapitre 8]{Wat44}), ainsi que des techniques de <<~dépoissonisation~>> reliant $M_{\mathcal{P}(\theta)}$ et $M_{n}$ lorsque $\theta \sim n \to \infty$, on en déduit le résultat suivant :
\begin{proposition}[Asymptotique du processus ponctuel associé à la mesure de Plancherel, \cite{BOO00}]\label{asymptoticpoisson}
Soit $(x^{n})_{n \in \N}=(x^{n}_{1}<\cdots<x^{n}_{k})_{n \in \N}$ une suite de vecteurs d'entiers telles que les limites finies ou infinies
$$ a_{i}= \lim_{n \to \infty} \frac{x^{n}_{i}}{\sqrt{n}} \qquad;\qquad d_{ij}= \lim_{n \to \infty}(a_{i}-a_{j})$$
existent pour tous indices $i$ et $j$. On note $\rho_{n}$ la fonction de corrélation discrète du processus ponctuel $\mathcal{D}_{\star}(M_{n})$, où $M_{n}$ est la mesure de Plancherel d'ordre $n$.
\begin{enumerate}
\item Si $x^{n}$ se scinde en deux parties $y^{n}$ et $z^{n}$ dont la distance tend vers l'infini, alors $\lim_{n \to \infty} \rho_{n}(x^{n})=(\lim_{n \to \infty} \rho_{n}(y^{n}))\times(\lim_{n \to \infty} \rho_{n}(z^{n}))$. Autrement dit, les comportements locaux en des endroits distincts sont asymptotiquement indépendants.
\item Si $x^{n}$ reste non scindée et si $a=a_{1}=\cdots=a_{k}$ est un point de $]-2,2[\,$, alors $\lim_{n \to \infty} \rho_{n}(x^{n})= \det(S^{a}(d_{ij}))$, où $S^{a}$ est le noyau sinus discret :
$$S^{a}(x)=\frac{\sin(\arccos(a/2)\,x)}{\pi x}\,.$$
\end{enumerate}
\end{proposition}
\noindent Au bord de la forme limite (c'est-à-dire avec $a=\pm 2$), le noyau sinus discret convenablement renormalisé dégénère en le noyau d'Airy, ce qui compte tenu de la proposition \ref{asymptoticgue} fournit une preuve de l'équivalence de Baik-Deift-Johansson, \emph{cf.} \cite[\S4]{BOO00}. Cette preuve est assurément la plus naturelle, puisqu'elle repose simplement sur l'étude asymptotique des fonctions de corrélation des processus $X=(X_{1},X_{2},\ldots)$ et $X=(Y_{1},Y_{2},\ldots)$. En contrepartie, elle utilise des propriétés fines des fonctions de Bessel et des fonctions d'Airy, et ces arguments semblent difficiles à généraliser ; nous reviendrons sur ce point dans la section \ref{qplancherelprocess}.\bigskip

\section[\'Equivalence de Baik-Deift-Johansson : l'approche géométrique]{\'Equivalence de Baik-Deift-Johansson : l'approche\\ géométrique}\label{bdjgeom}

Une seconde preuve de nature plus conceptuelle est due à A. Okounkov (\cite{Oko00}), et elle repose sur des interprétations combinatoires (géométriques) des polynômes traciaux de matrices aléatoires et des polynômes traciaux en les éléments de Jucys-Murphy dans l'algèbre du groupe symétrique. Nous concluons ce chapitre et la partie d'introduction de ce mémoire en exposant les arguments principaux de cette approche, qui laisse entrevoir des connections inattendues entre la théorie asymptotique des représentations et des problèmes géométriques tels que la théorie des invariants de Gromov-Witten et le décompte de structures combinatoires sur des surfaces (cartes et revêtements). Cette approche a motivé la recherche d'une interprétation géométrique des observables de diagrammes présentées dans le chapitre \ref{tool} ; on renvoie en particulier à \cite{Sni06a,Sni06b} pour de plus amples détails sur ce point.\bigskip
\bigskip

Compte tenu du théorème de Lévy reliant la convergence en loi d'une suite de variables aléatoires à celle des fonctions caractéristiques (voir l'exemple 5.5 dans \cite[chapitre 1]{Bil69}), le théorème \ref{bdj} est équivalent à l'assertion suivante : si $\widehat{X}(\xi)=\sum_{j=1}^{\infty} \exp(X_{j}\,\xi)$ et $\widehat{Y}(\xi)=\sum_{j=1}^{\infty} \exp(Y_{j}\,\xi)$ sont les transformées de Laplace des mesures $\sum_{j=1}^{\infty} \delta_{X_{j}}$ et $\sum_{j=1}^{\infty} \delta_{Y_{j}}$, alors pour tout entier $s\geq 1$ et tous paramètres $\xi_{1},\ldots,\xi_{s}$ strictement positifs,
$$\lim_{n \to \infty} \esper\left[\widehat{X}(\xi_{1})\cdots \widehat{X}(\xi_{s})\right]=\lim_{n \to \infty} \esper\left[\widehat{Y}(\xi_{1})\cdots \widehat{Y}(\xi_{s})\right]\,,$$
l'existence de la limite pour le terme de droite étant garantie par le théorème de Tracy-Widom.
Fixons des paramètres positifs $\xi_{1},\ldots,\xi_{s}$, et des entiers positifs pairs $k_{i} \sim \xi_{i}\,n^{2/3}$ ; on note $k=\sum_{i=1}^{s}k_{i}$. Le polynôme tracial $$M(k_{1},\ldots,k_{s})= (2\sqrt{n})^{-k}\,\,\esper\left[\prod_{i=1}^{s} \tr M^{k_{i}}\right]$$ est principalement déterminé par les puissances des plus grandes valeurs propres en valeur absolue (donc, celles au bord du spectre). Comme
$$\left(\frac{\lambda_{j}}{2\sqrt{n}}\right)^{k_{i}}\,\,\sim\,\,\left(1+\frac{Y_{j}}{n^{2/3}}\right)^{\xi_{i}\,n^{2/3}} \,\,\to \,\,\exp(Y_{j}\,\xi_{i})\,,$$
 $M(k_{1},\ldots,k_{s})$ a par conséquent la même asymptotique que $\esper[\widehat{Y}(\xi_{1})\cdots \widehat{Y}(\xi_{s})]$. De plus, la formule de Wick fournit un développement topologique de $M(k_{1},\ldots,k_{s})$ :
$$M(k_{1},\ldots,k_{s})=\frac{1}{2^{k}}\sum_{S}\, n^{\chi(S)-s} \,|\mathrm{Map}(k_{1},\ldots,k_{s};S)|\,,$$
la somme étant effectuée sur les classes d'homéomorphisme de surfaces compactes éventuellement non connexes. Ainsi, du côté des matrices aléatoires, le problème est ramené à celui de l'asymptotique des nombres de cartes $ |\mathrm{Map}(k_{1},\ldots,k_{s};S)|$.\bigskip
\bigskip

Du côté des partitions aléatoires, les analogues des fonctions $M(k_{1},\ldots,k_{s})$ sont les polynômes traciaux $C(k_{1},\ldots,k_{s})=2^{-k}\,\sqrt{n}^{s-k}\,\esper[\prod_{i=1}^{s}( L_{i})^{k_{i}}]$, où les $L_{i}$ sont les éléments de Jucys-Murphy
$$L_{1}=(1,2)+(1,3)+\cdots+(1,n)\qquad;\qquad L_{2}=(2,3)+\cdots+(2,n)\qquad;\qquad \cdots$$
et où $\esper[\cdot]$ désigne l'espérance de l'espace de probabilité non commutatif $\C\sym_{n}$, c'est-à-dire $(1/n!)\,\tr(\cdot)$. Notons que les $L_{i}$ ici définis diffèrent légèrement des $J_{i}$ du paragraphe \ref{jucysmurphy} ; ils ont néanmoins les mêmes propriétés, et en particulier, leurs actions sur un module de Specht $V^{\mu}$ ont pour valeurs propres les contenus des cases de $\mu$. Ainsi,
\begin{align*}\esper[L_{1}^{k}]&=\frac{1}{n!}\sum_{\mu \in \ym_{n}} (\dim \mu )\,\varsigma^{\mu}(L_{1}^{k})=\frac{1}{n!}\sum_{\mu \in \ym_{n}} \dim \mu \sum_{(\mu_{i},i)\text{ coin de }\mu} [\dim (\mu \setminus (\mu_{i},i))]\,(\mu_{i}-i)^{k}\\
&=\sum_{\mu \in \ym_{n}} M_{n}(\mu)\sum_{(\mu_{i},i)\text{ coin de }\mu} \delta_{i}^{*}(\mu)\,(\mu_{i}-i)^{k}\,, \end{align*}
où les $\delta_{i}^{*}(\mu)=[\dim(\mu \setminus  (\mu_{i},i))]/[\dim \mu]$ sont les probabilités de cotransition du processus de Plancherel, c'est-à-dire les probabilités de transition de l'unique processus markovien sur $\ym$ qui est décroissant, compatible avec les règles de branchement et qui conserve les mesures de Plancherel. L'existence d'une forme limite $\Omega$ pour les diagrammes sous la mesure de Plancherel implique la convergence en probabilité
$$\forall i,\,\,\,\delta_{i}^{*}(\mu)\,\sqrt{n} \to 1\,,$$
voir \cite[\S3.2.2]{Oko00} --- l'asymptotique est également valable pour les probabilités de transition du processus de Plancherel. Fixons des paramètres positifs $\xi_{1},\ldots,\xi_{s}$, et des entiers positifs $k_{i}'\sim \xi_{i}\,n^{1/3}$ ; on note comme précédemment $k'=\sum_{i=1}^{s}k_{i}'\,$. D'après ce qui précède,
$$\frac{1}{2^{k_{1}'}(\sqrt{n})^{k_{1}'-1}}\,\esper\left[(L_{1})^{k_{1}'}\right]=\sum_{\mu \in \ym_{n}} M_{n}(\mu)\sum_{(\mu_{i},i)\text{ coin de }\mu} \left(\delta_{i}^{*}(\mu)\,\sqrt{n}\right)\,\left(\frac{\mu_{i}-i}{2\sqrt{n}}\right)^{k_{1}'}$$ 
et le rapport $|(\mu_{i}-i)/2\sqrt{n}|$ est maximal pour $i=1,2,\ldots$ et $i=\ell(\mu),\ell(\mu)-1,\ldots$ Dans le premier cas,
$$\left(\frac{\mu_{i}-i}{2\sqrt{n}}\right)^{k_{1}'} \sim \left(\frac{\mu_{i}}{2 \sqrt{n}}\right)^{k_{1}'}$$
car $k_{1}' \propto n^{1/3}$ et $\mu_{i} \propto n^{1/2}$ ; dans le second cas,
$$\left(\frac{\mu_{\ell(\mu)-i}-(\ell(\mu)-i)}{2\sqrt{n}}\right)^{k_{1}'} \sim \left(-\frac{ \mu_{j}'}{2\sqrt{n}}\right)^{k_{1}'} .$$
Dans chaque cas, $\delta_{i}^{*}(\mu)\,\sqrt{n}$ tend vers $1$. En utilisant l'invariance par conjugaison des diagrammes de la mesure de Plancherel, on conclut que :
$$\frac{1}{2^{k_{1}'}(\sqrt{n})^{k_{1}'-1}}\,\esper\left[(L_{1})^{k_{1}'}\right] \sim \sum_{\mu \in \ym_{n}} M_{n}(\mu) \left\{\sum_{i=1}^{\infty} \left(\frac{\mu_{i}}{2\sqrt{n}}\right)^{k_{1}'}+\left(-\frac{\mu_{i}'}{2\sqrt{n}}\right)^{k_{1}'}\right\}\,.$$
En réitérant ce raisonnement avec les éléments $L_{2},\ldots,L_{s}$, on voit plus généralement que le polynôme tracial en les éléments de Jucys-Murphy $C(k_{1}',\ldots,k_{s}')$ a la même asymptotique que
$$\sum_{\mu \in \ym_{n}} M_{n}(\lambda)\left\{\sum_{i_{1},\ldots,i_{s}=1}^{\infty} \frac{(\mu_{i_{1}})^{k_{1}'} \cdots (\mu_{i_{s}})^{k_{s}'}}{(2\sqrt{n})^{k'}}+ \cdots \right\}$$
où les points de suspension dans la somme correspondent à $2^{s}-1$ autres termes mettant en jeu les $-\mu_{i}'$. De plus, comme dans le cas des matrices aléatoires, $C(k_{1}',\ldots,k_{s}')$ a la même asymptotique que $\esper[\widehat{X}(\xi_{1})\cdots \widehat{X}(\xi_{s})]$, car
$$\left(\frac{\mu_{j}}{2\sqrt{n}}\right)^{k_{i}'} \sim \left(1+\frac{X_{j}}{n^{1/3}}\right)^{\xi_{i}\,n^{1/3}} \,\,\to\,\,\exp(X_{j}\,\xi_{i})\,.$$
La preuve <<~combinatoire~>> de l'équivalence de Baik-Deift-Johansson est donc ramenée à l'équivalence asymptotique $M(k_{1},\ldots,k_{s}) \sim C(k_{1}',\ldots,k_{s}')$.\bigskip
\bigskip

Ceci étant, on peut donner un développement topologique de $C(k_{1}',\ldots,k_{s}')$ analogue à celui des fonctions $M$, mais en termes de revêtements ramifiés de la sphère. Pour commencer, on remplace chaque $L_{i}=\sum_{j>i}(i,j)$ par 
$$\widetilde{L_{i}}=L_{i}-\sum_{j=i+1}^{s}(i,j)=(i,s+1)+(i,s+2)+\cdots+(i,n)\,;$$
ceci ne change pas l'asymptotique des polynômes traciaux $C(k_{1}',\ldots,k_{s}')$ --- nous noterons $\widetilde{C}(k_{1}',\ldots,k_{s}')$ les polynômes modifiés. Maintenant, comme $\esper[\sigma]= \mathbb{1}_{\sigma=\id}$ pour toute permutation $\sigma$, $\widetilde{C}(k_{1}',\ldots,k_{s}')$ est le nombre de solutions dans $\sym_{n}$ de l'équation
$$(1,\tau_{1})(1,\tau_{2})\cdots(1,\tau_{k_{1}'})(2,\tau_{k_{1}'+1})\cdots(2,\tau_{k_{1}'+k_{2}'})\cdots (s,\tau_{k'})=\id$$
avec chaque $\tau_{i}$ dans $\lle s+1,n\rre$. Le groupe $\sym_{n-s}$ agit sur cet intervalle, et donc sur l'ensemble des solutions de l'équation par conjugaison ; de plus, l'orbite d'une solution $\tau$ est de cardinal $(n-s)(n-s-1)\cdots(n-s-d(\tau)-1) \sim n^{d(\tau)}$, où $d(\tau)=\card \{\tau_{1},\ldots,\tau_{k'}\}$. Ainsi, 
$$\widetilde{C}(k_{1}',\ldots,k_{s}') \sim \sum_{\{\tau\}/\sym_{n-s}} n^{d(\tau)}\,.$$
À toute solution $\tau=\{\tau_{1}',\ldots,\tau_{k}'\}$, on peut associer un revêtement ramifié de 
la sphère $\C\mathbb{P}^{1}$ en procédant comme suit\footnote{On renvoie à \cite[\S1.2]{LZ04} pour des détails sur le lien entre équations dans $\sym_{n}$ et revêtements ramifiés de $\C\mathbb{P}^{1}$ donné par la monodromie dans une fibre au-dessus d'un point base. Ce point sera également évoqué dans le chapitre \ref{bundle}, avec quelques rappels sur la notion de revêtement ramifié.} :
\begin{enumerate}
\item Sur la droite projective complexe $\C\mathbb{P}^{1}$, on place un point base $0$, et on répartit $k'$ points autour --- on les notera $1,2,\ldots,k'$. Le point $0$ sera non ramifié, et les points $1, \ldots,k'$ seront simplement ramifiés.
\item Le revêtement a $s+d(\tau)$ feuillets, dont $s$ feuillets spéciaux marqués au-dessus du point base $0$.
\item Si $\gamma_{i}$ est la boucle issue du point base et contournant le point $i$, alors la monodromie de cette boucle au dessus de $0$ est la transposition du feuillet spécial $i$ avec un feuillet non marqué $\tau_{i}$, voir la figure \ref{jmcovering}.
\end{enumerate}
\figcapt{\psset{unit=1mm}
\pspicture(0,-3)(135,80)
\psline(0,15)(10,-2)(100,-2)(90,15)
\rput(105,0){$\C\mathbb{P}^{1}$}
\pscurve[linecolor=red]{->}(45.5,33)(45,41)(41,41)(44.5,32.5)
\psframe*[linecolor=white](40,30)(43.5,42)
\psframe*[linecolor=white](40,30)(45,35)
\pscurve[linecolor=red]{->}(45.5,53)(45,61)(41,61)(44.5,52.5)
\psframe*[linecolor=white](45,50)(50,63)
\psdots(45,6)(55,7)(53,11)(43,14)(35,8)
\rput(58,7){$1$}
\rput(56,11.3){$2$}
\rput(33,8.5){$k'$}
\rput(44,12){$i$}
\rput(47,15){$\gamma_{i}$}
\psdots[linecolor=red](43,40)(43,60)
\psline[linecolor=red,linestyle=dashed](43,40)(43,60)
\pscurve{->}(45.5,7)(45,15)(41,15)(44.5,6.5)
\psline(0,35)(10,18)(100,18)(97.06,23)
\psline(0,40)(10,23)(100,23)(97.06,28)
\psline(0,45)(10,28)(100,28)(90,45)
\psline(15,28)(12.06,33)(7.06,33)
\rput(11,31){$i$}\rput(11,51){$\tau_{i}$}
\psline(15,48)(12.06,53)(7.06,53)
\psline(0,65)(10,48)(100,48)(91.18,63)
\psline(0,80)(10,63)(100,63)(90,80)
\psline(105,18)(107,20)(107,26.5)(108.5,27.5)(107,28.5)(107,35)(105,37)
\rput(125,27.5){$s$ feuillets spéciaux}
\psline(105,48)(107,50)(107,56.5)(108.5,57.5)(107,58.5)(107,65)(105,67)
\rput(121,57.5){$d(\tau)$ feuillets}
\endpspicture}{Revêtement ramifié de Jucys-Murphy associé à une équation $\tau$.\label{jmcovering}}{Revêtement de Jucys-Murphy associé à une équation $\tau$}
\noindent Nous noterons $\mathrm{Cov}_{\mathrm{JM}}(k_{1}',\ldots,k_{s}')$ l'ensemble des classes d'équivalence de \textbf{revêtements ramifiés de Jucys-Murphy} de paramètres $k_{1}',\ldots,k_{s}'$, et si $S$ est une surface orientable compacte, nous noterons $\mathrm{Cov}_{\mathrm{JM}}(k_{1}',\ldots,k_{s}';S)$ l'ensemble constitué de ceux qui sont homéomorphes à $S$. Par la formule de Riemann-Hurwitz (\emph{cf.} \cite[p. 39]{DS94}), on a alors $\chi(S)=2(s+d(\tau))-k'$, donc $$n^{d(\tau)}=n^{\frac{\chi(S)-s}{2}}\times n^{\frac{k'-s}{2}}.$$ 
D'autre part, deux équations donnent le même revêtement ramifié marqué si et seulement si l'une peut être obtenue à partir de l'autre par conjugaison par $\sym_{n-s}$.  On en déduit le développement topologique (asymptotique) 
$$C(k_{1}',\ldots,k_{s}')\sim \frac{1}{2^{k}} \sum_{S}\, n^{\frac{\chi(S)-s}{2}}\,|\mathrm{Cov}_{\mathrm{JM}}(k_{1}',\ldots,k_{s}';S)|\,,$$
la somme étant de nouveau effectuée sur les classes d'homéomorphisme de surfaces compactes. L'équivalence asymptotique $M(k_{1},\ldots,k_{s}) \sim C(k_{1}',\ldots,k_{s}')$ découle dès lors du résultat suivant :
\begin{proposition}[\'Equivalence asymptotique des nombres de cartes et des nombres de revêtements de Jucys-Murphy, \cite{Oko00}]\label{gwitten}
Si $k_{i}\sim \xi_{i}\,t$ avec $t$ tendant vers l'infini, alors pour tout genre $g$, il existe un fonction $m_{g}(\xi_{1},\ldots,\xi_{s})$ telle que $$\frac{|\mathrm{Map}(k_{1},\ldots,k_{s};\mathbb{T}_{g})|}{2^{k}} \sim t^{3g-3+3s/2}\,m_{g}(\xi_{1},\ldots,\xi_{s})\,.$$
Le résultat s'étend sans difficulté au cas d'une surface $S$ non connexe. Dans le même contexte, $$|\mathrm{Map}(k_{1},\ldots,k_{s};\mathbb{T}_{g})| \sim |\mathrm{Cov}_{\mathrm{JM}}(k_{1},\ldots,k_{s};\mathbb{T}_{g})|\,,$$ ce qui implique l'équivalence asymptotique $M(k_{1},\ldots,k_{s}) \sim C(k_{1}',\ldots,k_{s}')$ avec les normalisations des $k_{i}$ et des $k_{i}'$ en $n^{2/3}$ et $n^{1/3}$.
\end{proposition}

\begin{example}Si $g=0$ et $s=1$, alors $|\mathrm{Map}(2k;g=0)|$ est le nombre d'appariements sans croisement des côtés d'un $2k$-gone, ou encore le nombre d'arbres généraux à $k$ arêtes ; c'est donc le nombre de Catalan $C_{k}=\frac{1}{k+1}\binom{2k}{k}$. D'autre part, pour tout revêtement de Jucys-Murphy, appelons valence d'une feuille $F$ non spéciale le nombre de points $i$ de $\lle 1,k \rre$ tels que la monodromie autour de $i$ permute la feuille $F$. Alors, $\mathrm{val}(F) \geq 2$, et
$$\sum_{F}\, (\mathrm{val}(F)-2)=k-2d(\tau)=2s-\chi(S).$$
Lorsque $g=0$ et $s=1$, $k=2d(\tau)$, et la valence de tous les feuillets non spéciaux est donc $2$. Autrement dit, $|\mathrm{Cov}_{\mathrm{JM}}(2k;g=0)|$ est le nombre de solutions de l'équation
$$(1,i_{1})(1,i_{2})\cdots(1,i_{2k})=\id$$
telles que les $i_{j}$ soient réunis par paires. Il n'est pas difficile de voir que ce nombre est encore le nombre de Catalan $C_{k}$ ; dans ce cas particulier, l'équivalence asymptotique $|\mathrm{Map}|\sim |\mathrm{Cov}_{JM}|$ est donc une égalité. De plus, par la formule Stirling, la fonction $m_{0}(\xi)$ existe bien, et elle vaut $\frac{1}{\sqrt{\pi}}\left(\frac{\xi}{2}\right)^{-3/2}$.
\end{example}\bigskip\bigskip

La preuve de la proposition \ref{gwitten} donnée par Okounkov dans \cite{Oko00} repose sur l'utilisation des \textbf{graphes de rubans} comme objets combinatoires intermédiaires entre les cartes et les revêtements de Jucys-Murphy. Les graphes de rubans (métrisés) fournissent un modèle combinatoire des espaces de module des courbes $\mathscr{M}_{s,g}$ (\emph{cf.} \cite{Kon92}), et on peut interpréter l'expression $M(k_{1},\ldots,k_{s})$ comme une somme de Riemann d'une intégrale sur le modèle combinatoire de Kontsevich. D'autre part, il est bien connu que le décompte des revêtements ramifiés de la sphère est lié à la théorie de l'intersection sur les espaces de modules $\mathscr{M}_{s,g}$ --- ce qu'il est convenu d'appeler la \textbf{théorie de Gromov-Witten} --- et au calcul d'intégrales de classes de cohomologie sur le compactifié de cet espace (\emph{cf.} \cite{ELSV,Zvon05} et le paragraphe \ref{hurwitz}). Par conséquent, les deux expressions $M(k_{1},\ldots,k_{s})$ et $C(k_{1}',\ldots,k_{s}')$ sont certainement des discrétisations d'expressions intégrales sur les espaces de modules $\mathscr{M}_{s,g}$, ce qui laisse entrevoir une connection entre permutations aléatoires et espaces de modules analogue aux résultats de Kontsevich dans le contexte des matrices aléatoires, voir \cite{Kon92,OP01,OP02}. Ceci apporte une motivation supplémentaire à l'étude asymptotique des représentations, et clôt notre première partie d'introduction.

\pagestyle{empty}

\clearpage
~
\clearpage
~

\pagestyle{fancy}
\fancyhead{}
\fancyfoot{}
\fancyfoot[C]{\thepage}
\renewcommand{\chaptermark}[1]{\markboth{\chaptername\ \thechapter.\ #1.}{}} 
\renewcommand{\sectionmark}[1]{\markright{\thesection.\ #1.}} 
\fancyhead[RO]{\rightmark}
\fancyhead[LE]{\leftmark}
\setlength{\headheight}{15.5pt}

\part{Asymptotique des mesures de Plancherel des algèbres d'Hecke}

La seconde partie de ce mémoire est consacrée à l'étude asymptotique de partitions  distribuées suivant des mesures de probabilité issues de la théorie des représentations des groupes réductifs finis de type Lie, et en particulier la théorie des \textbf{algèbres d'Iwahori-Hecke}. Soit $q$ une puissance d'un nombre premier, et $\For_{q}$ le corps fini de cardinal $q$. Les groupes de matrices $\GL(n,\For_{q})$ ont leurs représentations irréductibles indexées par des familles de partitions $\bblambda=(\lambda_{1},\lambda_{2},\ldots)$, et partant, on peut s'intéresser à l'asymptotique des mesures de Plancherel des groupes $\GL(n,\For_{q})$ lorsque $n$ tend vers l'infini. J. Fulman (\cite{Ful06}) et A. Dudko (\cite{Dud08}) ont montré que si la famille $\bblambda$ est tirée suivant la mesure de Plancherel de la représentation régulière de $\GL(n,\For_{q})$, alors les partitions $\lambda_{i} \in \bblambda$ sont asymptotiquement indépendantes, et asymptotiquement réparties suivant certaines mesures de Schur ; en particulier, ces partitions restent presque sûrement de taille bornée (par contre leur nombre tend vers l'infini). Ce résultat et la théorie des représentations des groupes $\GL(n,\For_{q})$ sont rappelés dans le chapitre \ref{general} ; on s'inspire une nouvelle fois de \cite{Mac95}, et aussi de \cite{DL76,DM91,Hen03}.\bigskip\bigskip

Au cours de cette thèse, nous avons étudié les mesures de Plancherel d'autres représentations de $\GL(n,\For_{q})$, à savoir, les représentations paraboliques induites à partir d'un caractère d'un tore. Sous ces \textbf{mesures d'induction parabolique}, les partitions $\lambda_{i}$ restent en nombre fini, et ont bien leur taille qui tend vers l'infini ; une étude asymptotique semblable à celle du chapitre \ref{plancherel} devient dès lors possible. Dans le cas particulier où le tore est scindé et le caractère du tore est trivial, on obtient une mesure sur les partitions appelée \textbf{q-mesure de Plancherel}.\bigskip

 L'étude de cette mesure met en jeu la théorie des algèbres d'Hecke de type A (\emph{cf.} \cite{Mat99,Ram91,HR96,RR97}), qui sont les commutants des actions des groupes $\GL(n,\For_{q})$ sur les induits paraboliques ; elle est progressivement exposée dans les chapitres \ref{iwahori} et \ref{qplancherelmeasure}, et on y reprend les résultats des articles \cite{FM10,Mel10a}. Ainsi, la quantification $\sym_{n} \to \IH_{q}(\sym_{n})$ se traduit par une quantification de la base des caractères centraux dans l'algèbre des observables de diagrammes, et on obtient alors une loi des grands nombres (théorème \ref{firstasymptoticqplancherel}) et un théorème central limite (théorèmes \ref{secondasymptoticqplancherel} et \ref{thirdasymptoticqplancherel}) tout à fait analogues aux théorèmes \ref{firstasymptoticplancherel} et \ref{secondasymptoticplancherel} --- en revanche, le comportement asymptotique des partitions est très différent de celui observé pour les mesures de Plancherel des groupes symétriques.\bigskip\bigskip

Le chapitre \ref{arikikoike} est consacré à des généralisations de ces résultats au cas d'autres algèbres d'Hecke, en particulier les algèbres d'Hecke de type B reliées aux représentations paraboliques induites dans $\Sp(2n,\For_{q})$ à partir du caractère trivial d'un tore scindé. On dispose ainsi d'une loi des grands nombres et d'un théorème central limite pour les parts d'une bipartition aléatoire sous la B-$q$-mesure de Plancherel (théorème \ref{asymptoticbqplancherel}), et d'une conjecture pour la répartition des parts entre les deux parties de la bipartition (conjecture \ref{bqmix}). On présente également dans ce chapitre un cadre général dans lequel on peut espérer des résultats analogues aux théorèmes \ref{firstasymptoticqplancherel} et \ref{secondasymptoticqplancherel} : celui des algèbres d'Hecke généralisées décrites par un théorème de Lusztig (\emph{cf.} \cite{Lus84,HL80}).\bigskip\bigskip

\chapter{Combinatoire des groupes linéaires finis}\label{general}

Dans ce chapitre, nous exposons la théorie des groupes de matrices $\GL(n,k)$, où $k=\For_{q}$ est un corps fini. Nous rappelons la combinatoire des classes de conjugaison de telles matrices (section \ref{jordanfrobenius}), et nous présentons un résultat analogue à l'isomorphisme de Frobenius-Schur \ref{frobeniusschur} (section \ref{frobeniusschurglnfq}), ce qui permettra de calculer les caractères des groupes $\GL(n,\For_{q})$ (section \ref{delignelusztig}). On montre alors que les \textbf{caractères irréductibles} de $\GL(n,\For_{q})$ sont en bijection avec les familles de partitions $\bblambda : L /\Gal \to \ym$ de poids $n$, l'ensemble d'indexation $L/\Gal$ étant en bijection avec $$\overline{\For_{q}}^{\times}\big/ \,\mathrm{Galois}(\overline{\For_{q}}/\For_{q}),$$
c'est-à-dire l'ensemble des polynômes unitaires irréductibles sur $\For_{q}$ et différents de $X$. Enfin, dans la section \ref{dudkofulman}, on rappelle les résultats de Dudko et Fulman (\emph{cf.} \cite{Dud08,Ful06}) concernant l'asymptotique de la mesure de Plancherel de $\GL(n,\For_{q})$ sur ces familles de partitions.\bigskip\bigskip

Il existe en réalité deux familles fondamentales de caractères de $\GL(n,\For_{q})$ : les caractères irréductibles, et les \textbf{caractères de Deligne-Lusztig} obtenus par induction parabolique à partir de caractères de tores. Ces derniers fourniront des mesures de probabilité sur des familles de partitions $\bblambda$ réduites à un élément, ce qui mènera à l'étude de la $q$-mesure de Plancherel (chapitres \ref{iwahori} et \ref{qplancherelmeasure}). 
\bigskip

\section[Réduction de Jordan-Frobenius et classes de conjugaison]{Réduction de Jordan-Frobenius et classes de conjugaison}\label{jordanfrobenius}

Si $q=p^{e}$ est la puissance d'un nombre premier, on rappelle qu'il existe à isomorphisme près un unique corps fini $\For_{q}$ de cardinal $q$, et on peut le réaliser comme corps de rupture (et de décomposition) du polynôme $X^{q}-X$ sur $\For_{p}=\Z/p\Z$. Les extensions de corps de $\For_{q}$ sont les $\For_{q^{d}}$, et elles sont toutes galoisiennes ; le groupe $\mathrm{Galois}(\For_{q^{d}}/\For_{q})$ est cyclique d'ordre $d$ et est engendré par le morphisme de Frobenius $$\sigma : x \mapsto x^{q}.$$ La clôture algébrique $\overline{\For_{q}}$ a pour groupe de Galois sur $\For_{q}$ la limite profinie $\overline{\Z}=\varprojlim_{d}\Z/d\Z$, et ce groupe contient $\Z = \langle \sigma \rangle$ comme sous-groupe dense.\bigskip
\bigskip

Soit $G=\GL(n,\For_{q})$ le groupe des matrices inversibles de taille $n \times n$ et à coefficients dans $k=\For_{q}$. Par représentation matricielle, on peut voir $G$ comme le groupe des isomorphismes linéaires d'un espace vectoriel $V$ de dimension $n$ sur $\For_{q}$. Un tel espace est de cardinal $q^{n}$, et le nombre de familles libres de rang $k$ dans $V=(\For_{q})^{n}$ est 
$$(q^{n}-1)(q^{n}-q)(q^{n}-q^{2})\cdots (q^{n}-q^{r-1})=\prod_{i=1}^{r}\,(q^{n}-q^{i-1})\,.$$
En effet, pour construire une famille libre $(v_{1},\ldots,v_{r})$, on choisit un vecteur non nul $v_{1}$, puis un vecteur $v_{2}$ qui n'est pas dans l'espace vectoriel $\For_{q}[v_{1}]$ engendré par $v_{1}$, puis un vecteur $v_{3}$ qui n'est pas dans $\For_{q}[v_{1},v_{2}]$, etc. jusqu'au vecteur $v_{r}$ qui n'est pas dans $\For_{q}[v_{1},\ldots,v_{r-1}]$. Il y a $q^{n}-1$ possibilités pour le premier vecteur, $q^{n}-q$ pour le second, etc. jusqu'à $q^{n}-q^{r-1}$ possibilités pour le vecteur $v_{r}$, d'où le résultat. Comme un élément de $G$ peut être vu comme la matrice de transition entre la base canonique et une base quelconque $\mathcal{B}'$ de $V$, on en déduit :
$$\card \GL(n,\For_{q})= \text{nombre de bases de }V = (q^{n}-1)(q^{n}-q)\cdots (q^{n}-q^{n-1})\,.$$
Dans cette section, on s'intéresse aux classes de conjugaison de $\GL(n,\For_{q})$, c'est-à-dire les classes d'équivalence de matrices pour la relation
$$g_{1} \sim g_{2} \quad\iff \quad \exists h \in \GL(n,\For_{q}),\,\,\,g_{1}=h\,g_{2}\,h^{-1}\,.$$
Nous allons montrer que ces classes de conjugaison sont indexées par des familles de partitions appelées \textbf{polypartitions}. Si $u\in \mathrm{M}(n,k)$ est une matrice arbitraire de taille $n \times n$ (pas nécessairement inversible), son action sur $V=k^{n}$ fournit une structure de $k[X]$-module compatible avec la $k$-structure préexistante :
$$\forall P \in k[X],\,\,\,\forall v \in V,\,\,\,P\cdot v = [P(u)](v)\,.$$
Notons $V_{u}$ l'espace $V$ considéré comme $k[X]$-module. Comme $k[X]$ est un anneau principal et $V_{u}$ est finiment engendré sur cet anneau (par exemple par une $k$-base de $V$), on peut appliquer le théorème de structure des modules finis sur les anneaux principaux : il existe des polynômes $k$-irréductibles unitaires $\phi_{1},\ldots,\phi_{r}$ et des partitions  $\mu_{1},\ldots,\mu_{r}$ telles que
$$n=\sum_{i=1}^{r}|\mu_{i}|\,\deg \phi_{i}\qquad\text{et}\qquad V_{u} \simeq_{k[X]} \bigoplus_{i=1}^{r}\bigoplus_{j=1}^{\ell(\mu_{i})} k[X]/(\phi_{i}^{\mu_{i,j}})\,.$$
De plus, les polynômes $\phi_{i}$ et les partitions $\mu_{i}$ sont entièrement déterminés par $u$, et la décomposition de $V_{u}$ en modules primaires peut être interprétée de la fa\c con suivante : dans une $k$-base $\mathcal{B}'$ de $V$, l'action de $u$ sur $V$ est donnée par une matrice diagonale par blocs, les blocs étant les matrices de Jordan $J_{(\phi_{i})^{\mu_{i,j}}}$ --- on rappelle que
$$ J_{P}=\begin{pmatrix}
0 &   &          &   &-a_{0}\\
1 &0 &          &   &-a_{1}\\
   &1 &\ddots&   &\vdots\\
   &   &\ddots&0 &-a_{n-2}\\
   &    &         &1 &-a_{n-1}
\end{pmatrix}$$
si $P(X)=X^{n}+a_{n-1}X^{n-1}+\cdots+a_{1}X+a_{0}$. Nous noterons $\bbmu=\bbmu_{u}$ l'ensemble des paires $(\phi_{i},\mu_{i})$ ; c'est une $k$-polypartition de taille $n$. La matrice par blocs précédemment évoquée sera notée $J_{\bbmu}$, et elle est conjuguée à $u$ par la matrice de transition de $\mathcal{B}$ à $\mathcal{B}'$, $\mathcal{B}$ désignant la base canonique de $V=k^{n}$. En particulier, le polynôme minimal $m_{u}$ et le polynôme $\chi_{u}$ peuvent aisément être retrouvés à partir de la polypartition $\bbmu$ :
$$m_{u}(X)=\prod_{i=1}^{r} \,(\phi_{i}(X))^{\mu_{i,\ell(\mu_{i})}}\qquad;\qquad\chi_{u}(X)=\prod_{i=1}^{r}\prod_{j=1}^{\ell(\mu_{i})} \,(\phi_{i}(X))^{\mu_{i,j}}\,.$$
La matrice $u$ est donc inversible si et seulement si $X$ n'apparaît pas dans la polypartition $\bbmu_{u}$. Si $u_{1}$ et $u_{2}$ sont deux matrices conjuguées dans $\mathrm{M}(n,k)$ par $h \in \GL(n,k)$, alors $h$ est un isomorphisme de $k[X]$-module entre $V_{u_{1}}$ et $V_{u_{2}}$, car pour tout polynôme $P$,
$$P(u_{2})\circ h=h\circ P(u_{1})\,.$$ 
Par conséquent, $V_{u_{1}}$ et $V_{u_{2}}$ ont la même décomposition en $k[X]$-modules primaires, et $u_{1}$ et $u_{2}$ sont conjuguées à la même matrice $J_{\bbmu}$. D'autre part, comme la décomposition d'un $k[X]$-module fini est unique, deux matrices $J_{\bblambda}$ et $J_{\bbmu}$ ne peuvent pas être conjuguées. Ainsi :
\begin{proposition}[Réduction de Jordan-Frobenius]
Les classes de conjugaison de matrices dans $\mathrm{M}(n,k)$ sont indexées par les $k$-polypartitions de poids $n$, et les classes de matrices inversibles correspondent aux polypartitions dans lesquelles le polynôme $X$ n'apparaît pas. Un représentant de la classe $C_{\bbmu}$ est la matrice diagonale par blocs de Jordan $J_{\bbmu}$.
\end{proposition}\bigskip\bigskip

Notons $M$ le groupe multiplicatif de $\overline{\For_{q}}$, et $\Gal$ le groupe $\mathrm{Galois}(\overline{\For}_{q}/\For_{q})$. Si $P(X)\neq X$ est un polynôme unitaire irréductible sur $\For_{q}$, alors il se scinde dans $\overline{\For_{q}}$ sous la forme
$$P(X)=\prod_{i=1}^{\deg P} (X-\sigma^{i}(a))\,,$$
où $a$ est une racine (non nulle) de $P$ qui appartient à l'extension $\For_{q^{\deg P}}$. Par conséquent, l'ensemble des polynômes irréductibles unitaires différents de $X$ s'identifie à l'ensemble quotient $M/\Gal$. Nous noterons $\ym(\For_{q})$ l'ensemble des fonctions 
$\bbmu : M/\Gal \to \ym$ qui sont presque nulles, et $\ym_{n}(\For_{q})$ le sous-ensemble constitué des fonctions de poids $n$, c'est-à-dire que 
$$\|\bbmu\|=\sum_{P \in M/\Gal} \deg P\,\,|\bbmu(P)|=n$$
si $\bbmu \in \ym_{n}(\For_{q})$. D'après ce qui précède :
\begin{corollary}[Classes de conjugaison de $\GL(n,\For_{q})$]
Les classes de conjugaison de $\GL(n,\For_{q})$ sont en bijection avec les polypartitions de $\ym_{n}(\For_{q})$.
\end{corollary}
\begin{example}
Supposons $k=\For_{3}$ et $n=2$. Pour tout $q$ et tout $n$, le polynôme $X^{q^{n}}-X$ se scinde dans $\For_{q}[X]$ en le produit de tous les polynômes irréductibles unitaires de degré $d$ divisant $n$. Par conséquent, comme 
$$X^{9}-X=X\,(X+1)\,(X+2)\,(X^{2}+1)\,(X^{2}+X+2)\,(X^{2}+2X+2)$$
dans $\For_{3}[X]$, les classes de conjugaison de $\GL(2,\For_{3})$ sont indexées par les polypartitions de $\ym_{2}(\For_{3})$ suivantes :
\begin{align*}
&\{X^{2}+1 : 1\}\quad;\quad \{X^{2}+X+2 : 1\}\quad;\quad\{X^{2}+2X+2 : 1\}\quad;\quad\{X+1 : 2\}\\
&\{X+2 : 2\}\quad;\quad\{X+1 : 1^{2}\}\quad;\quad\{X+2 : 1^{2}\}\quad;\quad\{X+1 : 1\,\,\,;\,\,\,X+2 : 1\}\,.
\end{align*}
\end{example}
\bigskip\medskip

Pour conclure cette section, calculons le cardinal de la classe de conjugaison $C_{\bbmu}$ dans $\GL(n,\For_{q})$, où $\bbmu$ est une polypartition arbitraire de $\ym_{n}(\For_{q})$. Le centralisateur de la matrice $J_{\bbmu}$ dans $\GL(n,\For_{q})$ peut être vu comme l'ensemble des automorphismes de $k[X]$-modules de $V_{\bbmu}$ ; par la formule des classes, on a donc
$$\card C_{\bbmu}=\card \GL(n,\For_{q})\,\big/\,\card \mathrm{Aut}_{k[X]}(V_{\bbmu})\,.$$
Notons $V_{\phi_{i},\mu_{i}}$ la composante $\phi_{i}$-primaire de $V_{\bbmu}$, où $\bbmu=(\phi_{i},\mu_{i})_{i \in \lle 1,r\rre}$. Il est clair que $$\mathrm{Aut}_{k[X]}(V_{\bbmu}) = \prod_{i=1}^{r} \mathrm{Aut}_{k[X]}(V_{\phi_{i},\mu_{i}})\,,$$ donc le calcul du cardinal de $C_{\bbmu}$ se ramène à celui du cardinal du groupe d'automorphismes d'un $k[X]$-module primaire $V_{\phi,\mu}$. Dans ce qui suit, on note 
$$\|\bbmu\|= \sum_{i=1}^{r} (\deg \phi_{i})\,|\mu_{i}|\qquad;\qquad b(\bbmu)=\sum_{i=1}^{r} (\deg \phi_{i})\,b(\mu_{i})\,.$$
Si $x \neq 1$ est un nombre réel, $(x;x)_{m}$ est le symbole de Pochhammer usuel, \emph{i.e.}, $(x;x)_{m}=(1-x)(1-x^{2})\cdots(1-x^{m})$.
\begin{proposition}[Cardinal d'un classe de conjugaison dans $\GL(n,\For_{q})$]\label{cardclass}
Le nombre d'automorphismes d'un $\For_{q}[X]$-module primaire $V_{\phi,\mu}$ est $$Q^{|\mu|+2b(\mu)} \,\prod_{s\geq 1} \,(Q^{-1};Q^{-1})_{m_{s}(\mu)}\,,$$ où $Q=q^{\deg \phi}$. Par conséquent, le cardinal de la classe de conjugaison $C_{\bbmu}$ dans $\GL(n,\For_{q})$ est 
$$\frac{(q^{n}-1)(q^{n}-q)\cdots(q^{n}-q^{n-1})}{q^{\|\bbmu\|+2b(\bbmu)} \, \prod_{i=1}^{r}\prod_{k \geq 1} \,(q^{-\deg \phi_{i}};q^{-\deg \phi_{i}})_{m_{k}(\mu_{i})}}\,.$$
\end{proposition}\bigskip

\noindent L'isomorphisme de Frobenius-Schur pour l'anneau de Grothendieck de la catégorie des représentations des groupes linéaires $\GL(n,\For_{q})$ donnera une preuve indirecte de ce résultat, \emph{cf.} le paragraphe suivant. Une autre preuve repose sur la \textbf{théorie de Hall} des modules sur un anneau de valuation discrète $\mathfrak{o}$, \emph{cf.} \cite[chapitre 2, théorème 1.6]{Mac95}. En effet, si $V_{\phi,\mu}$ est $\phi$-primaire --- ce qui revient à dire que le polynôme minimal de tout élément de $V_{\phi,\mu}$ est une puissance de $\phi$ --- alors on peut considérer $V_{\phi,\mu}$ comme un $k[X]_{\phi}$-module, $k[X]_{\phi}$ désignant l'anneau local constitué des fractions $P/Q$ sans puissance de $\phi$ dans $Q$. Cet anneau $\mathfrak{o}=k[X]_{\phi}$ est un anneau de valuation discrète d'idéal maximal $\mathfrak{p}=\phi\,k[X]_{\phi}$, et de corps résiduel l'unique extension de corps de $k$ de degré $\deg \phi$. Dans ce cadre, $V_{\phi,\mu}$ est un $\mathfrak{o}$-module de torsion de type $\mu$ :
$$V_{\phi,\mu}\simeq_{\mathfrak{o}} \bigoplus_{j\geq 1}\mathfrak{o}/\mathfrak{p}^{\mu_{j}}$$
et il suffit de montrer que $\mathrm{Aut}_{k[X]}(V_{\phi,\mu})=\mathrm{Aut}_{\mathfrak{o}} (V_{\phi,\mu})$ a le cardinal énoncé dans la proposition \ref{cardclass}. Ceci peut être fait par induction sur la hauteur du module de torsion, c'est-à-dire la taille maximale d'une part de $\mu$. En effet, si $V$ est un $\mathfrak{o}$-module de torsion de hauteur $h$ et de type $\mu$, alors $\mathfrak{p} \otimes_{\mathfrak{o}}V$ est un module de torsion de hauteur $h-1$ et de type $\mu-1=(\mu_{1}-1,\mu_{2}-1,\ldots)$. De plus, tout automorphisme $g$ de $V$ fournit un automorphisme $\widetilde{g}= \id_{\mathfrak{p}} \otimes g$ de $\mathfrak{p} \otimes_{\mathfrak{o}} V$, donc il suffit de montrer que la préimage par $\sim$ d'un automorphisme de $\mathfrak{p} \otimes_{\mathfrak{o}} V$ a
$$Q^{\ell(\mu)^{2}}\,(Q^{-1};Q^{-1})_{m_{1}(\mu)}$$
éléments, ce qui n'est pas très difficile.

\section{Isomorphisme de Frobenius-Schur pour $\GL(n,\For_{q})$}\label{frobeniusschurglnfq}

Détaillons maintenant la théorie des représentations du groupe $G=\GL(n,\For_{q})$, en s'intéressant en premier lieu à l'induction à partir de sous-groupes diagonaux par blocs du type $\GL(n_{1},\For_{q})\times \GL(n_{2},\For_{q}) \times \cdots \times \GL(n_{r},\For_{q})$, où $n=n_{1}+n_{2}+\cdots+n_{r}$ est une composition. Comme dans la section \ref{frobenius}, il sera utile d'introduire la somme directe $K(\For_{q})$ des groupes de Grothendieck $K(\GL(n,\For_{q}))$,  et ses versions tensorisées $K_{R}(\For_{q})=R\otimes_{\Z}K(\For_{q})=\bigoplus_{n\in \N}K_{R}(\GL(n,\For_{q}))$. Nous aurons également besoin de la théorie des \textbf{polynômes de Hall-Littlewood} et des \textbf{chaînes de modules} finis sur un anneau de valuation discrète ; rappelons-en brièvement les points principaux (voir \cite[chapitres 2 et 3]{Mac95}). Si $\mathfrak{o}$ est un anneau de valuation discrète de corps résiduel $k=\For_{q}$ et si $\lambda,\lambda^{(1)},\ldots,\lambda^{(r)}$ est une famille de partitions telle que $|\lambda|=\sum_{i=1}^{r}|\lambda^{(i)}|$, notons $G^{\lambda}_{\lambda^{(1)},\ldots,\lambda^{(r)}}$ le nombre de chaînes
$$0=M^{(0)}\subset M^{(1)}\subset \cdots \subset M^{(r)}=M$$
de $\mathfrak{o}$-modules de torsion, où  $M$ est un module de torsion de type $\lambda$ fixé, et où chaque quotient $M^{(i)}/M^{(i-1)}$ est de type $\lambda^{(i)}$. On peut montrer que ces nombres $G^{\lambda}_{\lambda^{(1)},\ldots,\lambda^{(r)}}(q)$ ne dépendent que de $q$, et de fa\c con polynomiale. L'\textbf{algèbre de Hall} de $\mathfrak{o}$ est l'espace vectoriel complexe $H(\mathfrak{o})$ de base $(u_{\lambda})_{\lambda \in \ym}$, avec les règles de multiplication
$$u_{\lambda}\,*\,u_{\mu}=\sum_{\rho} G^{\nu}_{\lambda,\mu}(q)\,u_{\nu}\,.$$
Ces règles peuvent être encodées au moyen de fonction symétriques. Pour toute partition $\lambda$, notons $P_{\lambda}(x;t)$ la fonction symétrique de Hall-Littlewood limite inverse des polynômes
$$P_{\lambda}(x_{1},\ldots,x_{n};t)=\bigg(\prod_{m_{s}(\lambda)\geq 1} \frac{1-t}{1-t^{m_{s}(\lambda)}}\bigg)\, \sum_{\sigma \in \sym_{n}} \left( x_{1}^{\lambda_{1}}\,x_{2}^{\lambda_{2}}\cdots x_{n}^{\lambda_{n}}\,\prod_{i<j} \frac{x_{i}-tx_{j}}{x_{i}-x_{j}}\right)^{\sigma}$$
dans $\Lambda[t]=\varprojlim_{n \to \infty} \C[x_{1},x_{2},\ldots,x_{n};t]^{\sym_{n}}$. Lorsque $t=0$, on retrouve les fonctions de Schur : $P_{\lambda}(x;0)=s_{\lambda}(x)$.
\begin{proposition}[Fonction caractéristique d'une algèbre de Hall, \cite{Litt61}]\label{hallalgebra}
Il existe un unique isomorphisme d'algèbres complexes $\psi : H(\mathfrak{o}) \to \Lambda$ tel que
$$\forall n \in \N,\,\,\psi(u_{(1^{n)}})=q^{-\frac{n(n-1)}{2}}\,e_{n}(x)\,.$$ 
Cet isomorphisme est donné sur la base $(u_{\lambda})_{\lambda \in \ym}$ par $\psi(u_{\lambda})=q^{-b(\lambda)}\,P_{\lambda}(x;q^{-1})$.
\end{proposition}\bigskip
\bigskip

Ceci étant, fixons un entier $n$, une composition $n=n_{1}+n_{2}+\cdots+n_{r}$, et des fonctions centrales $u_{1} \in Z(\C\GL(n_{1},\For_{q})),\ldots ,u_{r}  \in Z(\C\GL(n_{r},\For_{q}))$. Si $L$ est le groupe produit $\GL(n_{1},\For_{q})\times \cdots \times \GL(n_{r},\For_{q})$ et si $G=\GL(n,\For_{q})$, alors par analogie avec la construction du paragraphe \ref{frobenius}, on pourrait munir $K_{\C}(\For_{q})$ d'une structure d'anneau en posant
$$u_{1}\circ \cdots \circ u_{r}=\mathrm{Ind}_{L}^{G}(u_{1}\otimes \cdots \otimes u_{r})\,.$$
Mais pour ce produit, les règles de branchement sont complexes, et il n'existe pas de structure d'algèbre de Hopf raisonnable compatible avec ce produit. Ainsi, il faut plutôt considérer le \textbf{produit d'induction parabolique d'Harish-Chandra}, qui est défini comme suit. Notons $P$ le sous-groupe parabolique de $G$ constitué des matrices triangulaires par blocs
$$g=\begin{pmatrix} g_{11} & g_{12} & \cdots & g_{1r} \\
                    0 & g_{22}  & \cdots & g_{2r} \\
                    \vdots  &  \ddots       & \ddots & \vdots \\
                    0 & \cdots & 0 & g_{rr} 
\end{pmatrix},$$
avec chaque $g_{ii}$ dans $\GL(n_{i},\For_{q})$. Si $U \subset P$ est le sous-groupe des matrices par blocs unipotentes (\emph{i.e.}, $g_{ii}=1$ pour tout $i$), alors on a une suite exacte scindée de morphismes de groupes 
$$\begin{CD} 1 @>>> U @>>> P @>{\pi}>> L @>>> 1\end{CD}$$
avec $\pi(g)=(g_{11},\ldots,g_{rr})$. Autrement dit, tout élément $p\in P$ s'écrit de manière unique $p=ul$ avec $u \in U$ et $\pi(p)=l \in L$. Par conséquent, une fonction centrale $u$ sur $L$ peut être relevée en une fonction centrale sur $P$ en posant $u_{L \to P}(p)=u(\pi(p))$ --- de fa\c con équivalente, un $L$-module $V$ peut être relevé en un $P$-module en posant $p \cdot_{P} v= \pi(p) \cdot_{L} v$. Dans ce cadre, le produit d'induction parabolique de fonctions centrales est défini par 
$$u_{1}\circ u_{2}\circ \cdots \circ u_{r}= \mathrm{Ind}_{P}^{G}((u_{1}\otimes u_{2}\otimes \cdots \otimes u_{r})_{L\to P})\,.$$
On définit de fa\c con analogue le produit d'induction parabolique de modules (éventuellement virtuels) $V_{1},\ldots,V_{r}$ au-dessus des groupes $\GL(n_{1},\For_{q}), \ldots,\GL(n_{r},\For_{q})$. Cette opération est clairement distributive vis-à-vis de la somme directe de modules, d'où une structure d'anneau gradué sur $K_{\C}(\For_{q})$, le degré d'une combinaison linéaire de $\GL(n_{i},\For_{q})$-modules étant le plus grand entier $n_{i}$ mis en jeu dans la combinaison. \bigskip\bigskip

Notons que l'opération d'induction parabolique fait sens dès que $P$ est un sous-groupe parabolique de $G=\GL(n,\For_{q})$, et $L$ est un sous-groupe de Lévi, de sorte que $P=U \rtimes L$ avec $U$ sous-groupe unipotent (\emph{cf.} \cite{DM91} et \cite{Lus84}). Dans ce contexte plus général, nous noterons $\Ind_{L}^{G}(u)$ ou $\Ind_{L}^{G}(V)$ la fonction centrale ou le module obtenu par induction parabolique de $L$ vers $G$. Le foncteur additif $\Ind_{L}^{G}$ admet toujours un adjoint ${\Ind^{*}}^{G}_{L}$, d'où un coproduit
$$\Delta\bigg(V \in K_{\C}(\GL(n,\For_{q}))\bigg)=\sum_{l+m=n} {\Ind^{*}}_{\GL(l,\For_{q}) \times \GL(m,\For_{q})}^{\GL(n,\For_{q})}(V)$$
sur $K_{\C}(\For_{q})$ et une structure d'algèbre de Hopf graduée positive autoadjointe. L'objectif de cette section est d'obtenir un isomorphisme entre cette algèbre de Hopf et un produit tensoriel infini d'algèbres de fonctions symétriques.\bigskip
\bigskip

Pour commencer, calculons plus explicitement la fonction centrale $u$ obtenue par induction parabolique à partir de fonctions centrales $u_{1},\ldots,u_{r}$. Si $H \subset G$ sont deux groupes finis et si $\chi$ est un caractère de $H$, alors le caractère induit $\theta=\mathrm{Ind}_{H}^{G}(\chi)$ s'écrit
$$\theta(g)= \sum_{xH \in G/H} \chi(x^{-1}gx)\,,$$
étant entendu qu'on étend la fonction $\chi$ par $0$ en dehors de $H$. Ainsi,
$$u(g)=(u_{1} \circ \cdots \circ u_{r})(g)=\sum_{xP \in G/P} (u_{1}\otimes\cdots \otimes u_{r})_{L \to P}(x^{-1}gx)\,.$$
Or, si $F=V_{0}\subset V_{1}\subset \cdots \subset V_{r}$ est le drapeau standard de $(\For_{q})^{n}$ associé à la composition $n=n_{1}+\cdots +n_{r}$, alors $x^{-1}gx$ est dans $P$ si et seulement si $(x^{-1}gx)(F)=F$, donc, si et seulement si $g(x(F))=x(F)$. On en déduit que si $\bbmu$ est le type de $g$, alors
$$u(g)=u(\bbmu)=\sum \,u_{1}(\bbmu_{1})\times u_{2}(\bbmu_{2})\times \cdots \times u_{r}(\bbmu_{r})$$
la somme étant effectuée sur les chaînes de $k[X]$-modules $$0 =V_{0}\subset V_{1}\subset \cdots \subset V_{r}=V_{g}=V_{\bbmu}$$ avec $\dim V_{i}/V_{i-1}=n_{i}$, et $\bbmu_{i}$ désignant le type du quotient $V_{i}/V_{i-1}$. En regroupant les chaînes de modules en fonction de leurs types, on voit donc que
$$(u_{1}\circ \cdots \circ u_{r})(\bbmu)=\sum \,F^{\bbmu}_{\bbmu_{1},\ldots,\bbmu_{r}}(q)\,u_{1}(\bbmu_{1})\cdots u_{r}(\bbmu_{r})\,,$$
où la somme est prise sur les familles de polypartitions $\bbmu_{1},\ldots,\bbmu_{r}$ de tailles $n_{1},\ldots,n_{r}$, et où $F^{\bbmu}_{\bbmu_{1},\ldots,\bbmu_{r}}(q)$ est le nombre de chaînes de $\For_{q}[X]$-modules $0=V_{0}\subset \cdots \subset V_{r}=V$, avec $V$ module fixé de type $\bbmu$, et chaque quotient $V_{i}/V_{i-1}$ étant de type $\bbmu_{i}$.
Par unicité de la décomposition d'un $\For_{q}[X]$-module en modules primaires, ces nombres s'écrivent 
$$F^{\bbmu}_{\bbmu_{1},\ldots,\bbmu_{r}}(q)=\prod_{\phi \in M/\Gal} G^{\bbmu(\phi)}_{\bbmu_{1}(\phi),\ldots,\bbmu_{r}(\phi)}(q^{\deg \phi})\,.$$
Fixons quelques notations supplémentaires : pour tout polynôme $\phi \in M/\Gal$, $(X_{\phi,i})_{i \geq 1}$ est un ensemble de variables indépendantes de degrés $\deg X_{\phi,i}=\deg \phi$ ; $\Lambda_{\phi}$ est l'algèbre des fonctions symétriques complexes en les $X_{\phi,i}$ ; $\Lambda(\For_{q})$ est le produit tensoriel infini $\bigotimes_{\phi \in M/\Gal} \Lambda_{\phi}$ ; et $\psi_{\phi}$ est l'isomorphisme d'algèbres complexes entre $H(\For_{q}[X]_{\phi})$ et $\Lambda_{\phi}$ de la proposition \ref{hallalgebra}. Un élément $f(X_{\phi})$ de $\Lambda_{\phi}$ sera noté plus simplement $f(\phi)$ ; par conséquent, $\Lambda(\For_{q})$ est l'algèbre librement engendrée par les $p_{n}(\phi)$, $n \geq 1$ et $\phi \in M/\Gal$.
\begin{theorem}[Isomorphisme de Frobenius-Schur pour les groupes $\GL(n,\For_{q})$, \cite{Zel81}]
Si $\bbmu$ est une polypartition de $\ym_{n}(\For_{q})$, considérons la classe $C_{\bbmu}$ comme un élement de $K_{\C}(\GL(n,\For_{q}))$. Alors, l'application linéaire
$$ C_{\bbmu} \mapsto \bigotimes_{\phi \in M/\Gal} u_{\bbmu(\phi)}$$
est un isomorphisme d'algèbres complexes entre $K(\For_{q})$ et $\bigotimes_{\phi \in M/\Gal} H(\For_{q}[X]_{\phi})$. Si l'on compose cette application par $\bigotimes_{\phi \in M/\Gal} \psi_{\phi}$, l'application caractéristique ainsi obtenue est un isomorphisme d'algèbres de Hopf  graduées entre $K(\For_{q})$ et $\Lambda(\For_{q})$ qui vérifie
$$\forall \bbmu \in \ym(\For_{q}),\,\,\,\mathrm{ch}(C_{\bbmu})=\prod_{\phi \in M/\Gal}(q^{-\deg \phi})^{b(\bbmu(\phi))}\, P_{\bbmu(\phi)}(\phi;q^{-\deg \phi})=P_{\bbmu}\,.$$
\end{theorem}\bigskip

L'application $\mathrm{ch}$ est même un isomorphisme d'algèbres de Hopf positives autoadjointes (\emph{cf.} les derniers chapitres de \cite{Zel81}), à condition qu'on munisse $\Lambda(\For_{q})$ d'un produit scalaire qui est le produit tensoriel infini de $q^{\deg \phi}$-produits sur les anneaux $\Lambda_{\phi}$ (voir \cite[chapitre 3, \S4]{Mac95}). Pour ce produit scalaire, la base duale de $(P_{\bbmu})_{\bbmu \in \ym(\For_{q})}$ est $(Q_{\bbmu})_{\bbmu \in \ym(\For_{q})}$, avec 
$$Q_{\bbmu}=\prod_{\phi \in M/\Gal} (q^{\deg \phi})^{|\bbmu(\phi)|+b(\bbmu(\phi))}\,Q_{\bbmu(\phi)}(\phi;q^{-\deg \phi})\,,$$
où $Q_{\mu}(X;t)$ est la fonction symétrique $\prod_{k \geq 1}(t;t)_{m_{k}(\mu)}\,P_{\mu}(X;t)$. Avec ces formules, on retrouve bien
\begin{align*}\frac{1}{\card \mathrm{Aut}(V_{\bbmu})} &= \frac{\card C_{\bbmu}}{\card \GL(n,\For_{q})}= \scal{C_{\bbmu}}{C_{\bbmu}} = \scal{P_{\bbmu}}{P_{\bbmu}} = P_{\bbmu}/Q_{\bbmu}\\
&=\frac{1}{q^{\|\bbmu\|+2b(\bbmu)} \prod_{\phi} \prod_{k \geq 1}(q^{-\deg \phi};q^{-\deg \phi})_{m_{k}(\bbmu(\phi))}}\,.\end{align*}
Ainsi, la combinatoire dans $\GL(n,\For_{q})$ de l'opération d'induction à partir de sous-groupes diagonaux peut être encodée comme dans la section \ref{frobenius} par des fonctions symétriques, et les classes de conjugaison sont cette fois-ci en correspondance avec des produits de fonctions de Hall-Littlewood.\bigskip

\section[Caractères de Deligne-Lusztig et caractères irréductibles]{Caractères de Deligne-Lusztig et caractères irréductibles}\label{delignelusztig}

Il s'agit maintenant de décrire les représentations irréductibles de $\GL(n,\For_{q})$ avec le même formalisme ; les représentations de Deligne-Lusztig joueront un rôle intermédiaire dans notre présentation. Dans le groupe multiplicatif $M$, on rappelle que $\sigma$ est l'application de Frobenius $x \mapsto x^{q}$. Si $n \geq 1$, notons $M_{n}$ le groupe multiplicatif $(\For_{q^{n}})^{\times}$, qui est aussi $M^{\sigma^{n}}$ ; le groupe $M$ est la limite directe des groupes $M_{n}$. Le groupe $M_{n}$ est cyclique de cardinal $q^{n}-1$, donc son groupe dual $L_{n}=(M_{n})^{*}$ est également cyclique de cardinal $q^{n}-1$. Notons $L$ la limite directe des groupes $L_{n}$ :
$$L=\varinjlim_{n\to \infty}L_{n}= \varinjlim_{n \to \infty}\hom(M_{n},\C^{\times})=\hom(\varprojlim_{n \to \infty} M_{n},\C^{\times})=\hom(M_{\sharp},\C^{\times})\,,$$
la limite inverse $M_{\sharp}$ étant prise par rapport aux normes $x \in M_{n}\mapsto \prod_{i=1}^{n/d}\sigma^{di}(x) \in M_{d}$ définies pour $d \mathrel{|} n$. Un élément $\xi \in L$ est la donnée d'appariements $\scal{\xi}{\cdot\,}_{n} \in L_{n}$ pour $n$ assez grand (au sens de la divisibilité), de telle sorte que si $d \mathrel{|} n$, alors
$$\forall x \in M_{n},\,\,\,\scal{\xi}{x}_{n}=\scal{\xi}{N_{d,n}(x)}_{d}\,,$$
$N_{d,n}$ désignant l'application norme précédemment évoquée\footnote{En particulier, si $x \in M_{d} \subset M_{n}$, alors $\scal{\xi}{x}_{n}=(\scal{\xi}{x}_{d})^{n/d}$ est en général différent de $\scal{\xi}{x}_{d}$.}. L'application de Frobenius agit sur $L$, et $L_{n}$ s'identifie à $L^{\sigma^{n}}$. Dans l'ensemble quotient $L/\Gal$, le degré d'une classe $\Theta=[\xi]$ est $\deg \Theta=\card \Theta$. Une \textbf{polypartition duale} est une fonction presque nulle $\bblambda : L/\Gal \to \ym$, et on dit que $\bblambda$ est de poids $n$ si 
$$\|\bblambda\|=\sum_{\Theta \in L/\Gal} \deg \Theta\,\,|\bblambda(\Theta)|=n\,.$$\bigskip

Nous noterons $\ym^{*}(\For_{q})$ l'ensemble des $\For_{q}$-polypartitions duales, et $\ym^{*}_{n}(\For_{q})$ l'ensemble des polypartitions duales de poids $n$. Ces objets apparaissent naturellement dans la théorie des représentations du groupe $\GL(n,\For_{q})$, car ils paramètrent\footnote{Notons que pour tout $n$, $\ym^{*}_{n}(\For_{q})$ et $\ym_{n}(\For_{q})$ ont le même cardinal ; par contre, il n'y a pas de bijection privilégiée.} les $\GL(n,\For_{q})$-classes de conjugaison de caractères de tores maximaux ; expliquons en détail ce point. Dans $G=\GL(n,\overline{k})$, l'\textbf{application de Frobenius} $F$ est la transformation qui agit  par $\sigma$ sur chaque coordonnée ; c'est un morphisme de groupes, et $\GL(n,k)=G^{F}$.
Un \textbf{tore} de $G$ est un sous-groupe commutatif isomorphe à $M^{r}$ pour un certain entier $r \leq n$ ; on appelle tore de $G^{F}$ l'ensemble des points fixes $T^{F}$ dans un tore $T \subset G$ qui est laissé stable par l'application de Frobenius $F$. Un tore $T^{F}$ est dit \textbf{maximal} si son rang $r$ est égal à $n$, et il est dit \textbf{scindé} sur $k$ s'il est déjà isomorphe à $k^{r}$ dans $\GL(n,k)$ ; dans le cas contraire, il existe une extension de corps minimale $K|k$ telle que $T^{F}$ se scinde lorsqu'on étend les scalaires\footnote{Il est utile de voir $G$ et $T$ comme des \textbf{$k$-schémas en groupes}, c'est-à-dire, des groupes algébriques déterminés par équations polynomiales, et dont on peut considérer les solutions dans $k$ ou dans toute extension de corps de $k$ ; \emph{cf.} \cite{SGA3}.} à $K$. Dans $G=\GL(n,\overline{k})$, un tore est toujours scindé, et deux tores maximaux sont toujours conjugués sous l'action de $G$ (\emph{cf.} \cite{Che04}). En revanche, deux tores $F$-stables $S$ et $T$ ne sont pas forcément conjugués dans $G$ sous l'action de $G^{F}$, et on peut montrer que les $G^{F}$-classes de conjugaison de tores sont paramétrées par les $F$-classes de conjugaison du \textbf{groupe de Weyl}\footnote{On renvoie au chapitre \ref{iwahori} pour un rappel de la définition du groupe de Weyl et des sous-groupes de Borel d'un groupe algébrique.} $W$ de $G$, deux éléments $\omega_{1}$ et $\omega_{2}$ étant $F$-conjugués dans $W$ s'il existe $\zeta$ tel que 
$$\omega_{1}=\zeta\,\omega_{2}\,F(\zeta)^{-1}.$$
On renvoie à \cite{DM91} pour ces questions de rationalité des groupes algébriques définis sur un corps fini. Ici, $G^{F}$ admet un tore scindé maximal (à savoir, le sous-groupe diagonal $k^{n}$), et ceci implique que $F$ agit trivialement sur $W=\sym_{n}$. Les $F$-classes de conjugaison sont donc les classes usuelles, et les $G^{F}$-classes de conjugaison de tores maximaux $F$-stables sont indexées par les partitions $\lambda \in \ym_{n}$. Un tore maximal $T$ de type $\lambda$ vérifie
$$T^{F}\simeq_{\GL(n,\For_{q})} \,\,\prod_{i\geq 1} \,M_{\lambda_{i}}\,.$$
Si l'on considère maintenant les $G^{F}$-classes de conjugaison de paires $(T^{F},\zeta)$ avec $\zeta$ caractère (unidimensionnel) de $T^{F}$, alors on peut montrer qu'elles sont indexées par les polypartitions duales de $\ym^{*}_{n}(\For_{q})$ : les tores dans la classe $\bblambda$ vérifient
$$T^{F} \simeq_{\GL(n,\For_{q})} \prod_{\Theta \in L/\Gal} \prod_{i\geq 1} \,M_{\deg \Theta\,\,\bblambda(\Theta)_{i}}\,,$$
et modulo ces isomorphismes, les caractères dans la classe $\bblambda$ s'écrivent
$$\zeta\simeq_{\GL(n,\For_{q})} \prod_{\Theta \in L/\Gal}\prod_{i \geq 1} \,\scal{\xi}{\cdot\,}_{\deg \Theta\,\,\bblambda(\Theta)_{i}},$$
\emph{cf.} \cite[\S1.3]{Hen03}. Dans ce qui suit, nous notons $(T^{\bblambda},\zeta^{\bblambda})$ un tore et un caractère dans la $\GL(n,\For_{q})$-classe de conjugaison $\bblambda$.\bigskip
\bigskip

Pour toute polypartition duale $\bblambda$, on peut construire un $\GL(n,\For_{q})$-module virtuel $R^{\bblambda}$ par induction à partir du caractère $\zeta^{\bblambda}$ d'un tore $T^{\bblambda}$. Si le tore $T$ est contenu dans un sous-groupe de Borel $B \subset G$ qui est $F$-stable, il suffit d'utiliser l'induction parabolique de $T^{F}$ à $G^{F}$ \emph{via} $B^{F}$ ; malheureusement, pour de nombreux couples $(T^{\bblambda},\zeta^{\bblambda})$, le $\overline{k}$-tore sous-jacent n'est contenu dans aucun sous-groupe de Borel $F$-stable. Il y a néanmoins une fa\c con de définir l'induction parabolique de $T^{F}$ à $G^{F}$ même lorsque $T$ n'est pas inclus dans un sous groupe de Borel $F$-stable. Cette méthode est due à Deligne et Lusztig (voir \cite{DL76}), et peut être utilisée dès que $G$ est un groupe réductif de type Lie et $T$ est un tore maximal $F$-stable. L'idée est de considérer une certaine variété $X$ sur laquelle $G^{F} $ agit à gauche et $T^{F}$ agit à droite. La somme alternée
$H_{c}^{\bullet}(X)=\sum_{n\geq 0}(-1)^{n}\,H_{c}^{n}(X)$ des groupes de cohomologie à support compact de cette variété de Deligne-Lusztig est un $(G^{F},T^{F})$-bimodule virtuel, et le produit tensoriel au-dessus de $\C[T^{F}]$ de $H_{c}^{\bullet}(X)$ avec un $T^{F}$-module à gauche produit un $G^{F}$-module à gauche ; nous le noterons
$$\Ind_{T^{F}}^{G^{F}}(V)= H_{c}^{\bullet}(X)\otimes_{\C[T^{F}]}V.$$
Cette opération constitue l'\textbf{induction de Deligne-Lusztig}, et on retrouve l'induction parabolique d'Harish-Chandra lorsque $T$ est contenu dans un sous-groupe de Borel $F$-stable. Notons qu'on peut remplacer $T$ par un sous-groupe de Lévi $L$ rationnel, auquel cas l'induction de Deligne-Lusztig permet de passer outre l'absence éventuelle de sous-groupe parabolique $F$-stable contenant $L$. Ceci étant dit, les \textbf{caractères de Deligne-Lusztig} de $\GL(n,\For_{q})$, aussi appelés caractères basiques, sont les
$$R^{\bblambda}=\Ind_{T^{\bblambda}}^{G^{F}} (\zeta^{\bblambda})\,;$$
ces caractères virtuels ne dépendent du choix de représentants $(T^{\bblambda},\zeta^{\bblambda})$ dans les classes de conjugaison $\bblambda$.\bigskip
\bigskip

La description des $R^{\bblambda}$ en termes de fonctions symétriques dans $\Lambda(\For_{q})$ est due à J. A. Green, voir \cite{Gre55} ; nous reprendrons les notations de \cite[chapitre 4]{Mac95}. Si $\phi$ est un polynôme irréductible sur $\For_{q}$ et si $f \in \Lambda$, nous avons déjà défini le symbole $f(\phi) \in \Lambda(\For_{q})$ ; voyons de même comment <<~évaluer~>> une fonction symétrique en un élément $x \in M$, un caractère $\xi \in L$ ou une orbite $\Theta \in L/\Gal$. Comme $\Lambda$ est librement engendré par les sommes de puissances, il suffit de décrire les spécialisations $p_{n}(x)$, $p_{n}(\xi)$ et $p_{n}(\Theta)$. Pour $x\in M$, on pose :
$$p_{n}(x)=\begin{cases}p_{[n/\deg m_{x}]}(m_{x}) &\text{si }x \in M_{n},\\
0&\text{sinon}.
\end{cases}$$
Ainsi, $p_{n}(x)$ est homogène de degré $n$ ou $-\infty$, et met seulement en jeu les variables $X_{\phi,i}$ avec $\phi=m_{x}$ polynôme minimal de $x$. Ensuite, si $\xi \in L$, on pose :
$$p_{n}(\xi)=\begin{cases}
(-1)^{n-1}\,\sum_{x \in M_{n}} \scal{\xi}{x}_{n}\,p_{n}(x)& \text{si }\xi \in L_{n},\\
0&\text{sinon}.
\end{cases}$$
La fonction $p_{n}(\xi)$ est encore de degré $n$ ou $-\infty$ dans $\Lambda(\For_{q})$. Finalement, si $\Theta=[\xi]$ est une orbite dans $L/\Gal$, on pose $p_{n}(\Theta)=p_{[n\,\deg\Theta]}(\xi)$ ; le choix d'un représentant $\xi$ de $\Theta$ ne change pas la définition. Dans ce contexte, le théorème 14 de \cite{Gre55} peut être interprété comme suit :

\begin{proposition}[Caractères de Deligne-Lusztig et fonctions symétriques, \cite{Gre55}] Pour toute polypartition duale $\bblambda$, si $B^{\bblambda}=\mathrm{ch}(R^{\bblambda})$, alors 
$$B^{\bblambda}=(-1)^{|\bblambda|-\sum_{\Theta \in L/\Gal}|\bblambda(\Theta)|}\prod_{\Theta \in L/\Gal} p_{\bblambda(\Theta)}(\Theta)\,.$$
\end{proposition}
\begin{theorem}[Caractères irréductibles de $\GL(n,\For_{q})$ et fonctions symétriques, \cite{Gre55}]
Les caractères irréductibles de $\GL(n,\For_{q})$ sont indexées par les polypartitions duales $\bblambda \in \ym_{n}^{*}(\For_{q})$, et si $S^{\bblambda}=\mathrm{ch}(\chi^{\bblambda})$, alors
$$S^{\bblambda}=\prod_{\Theta \in L/\Gal} s_{\bblambda(\Theta)}(\Theta)\,.$$
\end{theorem}\bigskip\bigskip

Par suite, les modules de Deligne-Lusztig $(R^{\bblambda})_{\bblambda \in \ym_{n}^{*}(\For_{q})}$ forment une base linéaire du grou\-pe de Grothendieck $K_{\C}(\GL(n,\For_{q}))$ --- plus généralement d'ailleurs, si $G^{F}$ est un groupe réductif fini de type Lie, alors tout caractère irréductible apparaît dans au moins un caractère de Deligne-Lusztig, et deux caractères irréductibles qui interviennent dans le même caractère de Deligne-Lusztig ont des valeurs identiques en les éléments unipotents, ces valeurs étant données par des polynômes de Green (\emph{cf.} \cite{Lus76}). Les matrices de transition entre la base des caractères irréductibles et la base de Deligne-Lusztig s'écrivent :
\begin{align*}\eps(\bbrho)\,B^{\bbrho}&\,\,=\sum_{\forall \Theta,\,\,|\bblambda(\Theta)|=|\bbrho(\Theta)|} \left(\prod_{\Theta \in L/\Gal} \varsigma^{\bblambda(\Theta)}(\bbrho(\Theta))\right)S^{\bblambda}\\
\eps(\bblambda)\,S^{\bblambda}&\,\,=\sum_{\forall \Theta,\,\,|\bblambda(\Theta)|=|\bbrho(\Theta)|} \left(\prod_{\Theta \in L/\Gal} \frac{\varsigma^{\bblambda(\Theta)}(\bbrho(\Theta)) }{ z_{\bbrho(\Theta)} }\right)B^{\bbrho}
\end{align*}
où $\eps(\bbrho)$ est le signe qui apparaît dans la définition de $B^{\bbrho}$, et où les $\varsigma^{\lambda}(\rho)$ sont les valeurs des caractères des groupes symétriques. Ainsi, la combinatoire des représentations irréductibles et des représentations obtenues par induction de Deligne-Lusztig à partir de tores peut de nouveau être traitée à partir de fonctions symétriques dans $\Lambda(\For_{q})=\bigotimes_{\phi \in M/\Gal} \Lambda_{\phi} = \bigotimes_{\Theta \in L/\Gal} \Lambda_{\Theta}$ : les caractère irréductibles correspondent à des produits de fonctions de Schur, et les caractères de Deligne-Lusztig correspondent à des produits de fonctions puissances.\bigskip

\begin{example}
On considère la polypartition duale $\bbrho=([\mathbb{1}] : 1^{n})$, $\mathbb{1}$ désignant le caractère trivial défini par $\scal{\mathbb{1}}{x}_{n}=1$ pour tout $x$ et tout $n$. Elle correspond au tore scindé $T^{F}=k^{n}$ et à son caractère trivial $\zeta(x_{1},\ldots,x_{n})=1$. Dans ce cas particulier, la représentation induite $R^{\bbrho}$ l'est au sens d'Harish-Chandra, car $T^{F}$ est inclus dans le sous-groupe de Borel rationnel $\GB(n,k)$ constitué des matrices triangulaires supérieures. On a donc :
$$R^{\bbrho}=\mathrm{Ind}_{\GB(n,\For_{q})}^{\GL(n,\For_{q})}(1)=\C[\GL(n,\For_{q})/\GB(n,\For_{q})]\,.$$
Autrement dit, $R^{\bbrho}$ est simplement l'espace des fonctions sur la variété des drapeaux complets $\GL(n,\For_{q})/\GB(n,\For_{q})$. Compte tenu des formules de changement de base entre les $S^{\bblambda}$ et les $B^{\bbrho}$, les polypartitions duales qui interviennent dans la décomposition de ce module en irréductibles sont celles qui sont supportées par la seule orbite $[\mathbb{1}]$ ; notons-les pour simplifier $\{[\mathbb{1}] : \lambda\}=\lambda$, et notons $U^{\lambda}$ les modules irréductibles correspondant, appelés \textbf{modules unipotents}. Alors :
$$\C[\GL(n,\For_{q})/ \GB(n,\For_{q})]= \sum_{\lambda \in \ym_{n}} \varsigma^{\lambda}(1^n) \,V^{\lambda}=\sum_{\lambda \in \ym_{n}} (\dim \lambda) \,\,U^{\lambda}\,. $$
Ainsi, la multiplicité du module unipotent $U^{\lambda}$ est égale à la dimension de la représentation irréductible de $\sym_{n}$ indexée par la partition $\lambda$. Une généralisation importante de ce résultat est due à Lusztig ; nous y reviendrons dans le chapitre \ref{arikikoike}.
\end{example}\bigskip

\begin{example}
 L'objet principal des deux chapitres suivants est l'étude de la mesure de Plancherel du module $\C[\GL(n,\For_{q})/\GB(n,\For_{q})]$ ; compte tenu ce qui précède, pour calculer cette mesure de probabilité, il suffit de connaître les \textbf{degrés génériques} $D_{\lambda}(q)$, c'est-à-dire les dimensions des modules unipotents $U^{\lambda}$. Il existe en fait une formule générale pour la dimension d'un module irréductible $V^{\bblambda}$ :
 $$\dim V^{\bblambda}=\left(\prod_{i=1}^{n}q^{i}-1\right)\times\left(\prod_{\Theta \in L/\Gal} \frac{(q^{\deg \Theta})^{b(\bblambda'(\Theta))}}{\prod_{(i,j) \in \bblambda(\Theta)} (q^{\deg \Theta})^{h(i,j)}-1}\right)$$
On renvoie à \cite[chapitre 4, \S6]{Mac95} pour une preuve de cette formule ; il s'agit de calculer dans $\Lambda(\For_{q})$ le produit scalaire $\langle S^{\bblambda}|P_{(1:1^{n})}\rangle$, et on utilise en particulier l'identité
$$s_{\lambda}(q^{-1},q^{-2},\ldots)=\frac{q^{b(\lambda')}}{\prod_{(i,j) \in \lambda} q^{h(i,j)}-1}\,,$$
voir \cite[chapitre 1, \S3, exemples 1 à 5]{Mac95}. Ainsi, la théorie des fonctions symétriques fournit une formule des équerres pour $\dim \bblambda$ analogue à celle qui existe pour les représentations irréductibles des groupes symétriques (\emph{cf.} \S\ref{frobenius}). Dans le cas particulier des modules unipotents $U^{\lambda}$, la formule se spécialise en
$$D_{\lambda}(q)=\dim U^{\lambda}=q^{b(\lambda')}\frac{\prod_{i=1}^{n}q^{i}-1}{\prod_{(i,j) \in \lambda} q^{h(i,j)}-1}\,.$$
Nous reviendrons en détail sur ce point dans la section \ref{qplancherelprocess}.
\end{example}\bigskip

\section[Asymptotique des mesures de Plancherel des groupes linéaires finis]{Asymptotique des mesures de Plancherel des groupes\\ linéaires finis}\label{dudkofulman}

Pour conclure ce chapitre, exposons les résultats de A. Dudko et J. Fulman relatifs à l'asymptotique de la mesure de Plancherel des représentations régulières des groupes linéaires finis $\GL(n,\For_{q})$. Si $\bblambda \in \ym_{n}^{*}(\For_{q})$, on vient de voir que 
$$\dim \bblambda = \left(\prod_{i=1}^{n} q^{i}-1\right)\, \prod_{ \Theta \in L/\Gal} s_{\bblambda(\Theta)}(q^{-\deg \Theta},q^{-2\deg\Theta },\ldots)\,.$$
La mesure de Plancherel de $\GL(n,\For_{q})$ s'écrit donc :
$$\mathrm{Pl}_{\GL(n,\For_{q})}[\bblambda]=\frac{(\dim \bblambda)^{2}}{\card \GL(n,\For_{q})} = q^{-\frac{n(n-1)}{2}} \left(\prod_{i=1}^{n} q^{i}-1 \right)\,\prod_{\Theta \in L/\Gal} s_{\bblambda(\Theta)}^{2}(q^{-\deg \Theta},q^{-2\deg\Theta },\ldots)\,.$$
Si $v <1$ est un paramètre réel, on rappelle l'identité d'Euler :
\begin{align*}
\frac{1}{\prod_{r=0}^{\infty} (1-vq^{-r})}&= \prod_{r=0}^{\infty} \left(\sum_{n_{r}=0}^{\infty} (vq^{-r})^{n_{r}}\right)=\sum_{n=0}^{\infty} v^{n} \left(\sum_{n_{1}+ \cdots + n_{r}+ \cdots = n} q^{-(n_{1}+2n_{2}+\cdots+rn_{r}+\cdots)}\right) \\
&=\sum_{n=0}^{\infty} v^{n} \left(\sum_{\ell(\lambda) \leq n} q^{-|\lambda|}\right)=\sum_{n=0}^{\infty} v^{n} \left(\sum_{\lambda_{1}' \leq n} q^{-|\lambda'|}\right)\\
&=\sum_{n=0}^{\infty} v^{n} \left(\sum_{m_{1},\ldots,m_{n}} q^{-(m_{1}+2m_{2}+\cdots+nm_{n})}\right)=\sum_{n=0}^{\infty} \frac{v^{n}}{\prod_{i=1}^{n}(1-q^{-i})}\,.
\end{align*}
Par conséquent, l'expression suivante est une mesure de probabilité sur l'ensemble des polydiagrammes duaux $\ym^{*}(\For_{q})=\bigsqcup_{n \in \N} \ym_{n}^{*}(\For_{q})$ :
\begin{align*}M[\bblambda]&=q^{\frac{n(n+1)}{2}}\,v^{n}\,\prod_{r=0}^{\infty}(1-vq^{-r})\,\prod_{i=1}^{n}(q^{i}-1)^{-1}\, \mathrm{Pl}_{\GL(n,\For_{q})}[\bblambda]\\
&=(qv)^{n}\,\prod_{r=0}^{\infty}(1-vq^{-r})\,\prod_{\Theta \in L/\Gal} s_{\bblambda(\Theta)}^{2}(q^{-\deg \Theta},q^{-2\deg\Theta },\ldots)\,.
\end{align*}
Pour toute partie $A$ de $\ym_{n}^{*}(\For_{q})$, $\prod_{i=1}^{n}(1-q^{-i})^{-1} \, \mathrm{Pl}_{\GL(n,\For_{q})}[A]$ est le coefficient de $v^{n}$ dans la série $\prod_{r=0}^{\infty}(1-vq^{-r})^{-1}\,M[A]$. Notons $\mathcal{L}_{v,q}$ la mesure de Schur $\mathcal{S}$ de paramètres $$a=(q^{-1},q^{-2},q^{-3},\ldots)\,\text{ et }\,b=(v,vq^{-1},vq^{-2},\ldots)\,.$$ 
Pour toute partition $\lambda$,
$$\mathcal{L}_{v,q}[\lambda]=\left(\prod_{r=1}^{\infty}(1-vq^{-r})^{r}\right)\,q^{2b(\lambda')+|\lambda|}\,v^{|\lambda|}\,s_{\lambda}^{2}(q^{-1},q^{-2},\ldots)\,.$$
Fixons des partitions $\lambda_{1},\ldots,\lambda_{k}$ arbitraires, et des orbites $\Theta_{1},\ldots,\Theta_{k} \in L/\Gal$. Si $A$ est l'événement $\{\bblambda(\Theta_{1})=\lambda_{1},\ldots,\bblambda(\Theta_{k})=\lambda_{k}\}$, alors d'après ce qui précède,
$$\prod_{i=1}^{n}(1-q^{-i})^{-1}\,\mathrm{Pl}_{\GL(n,\For_{q})}[A]=[v^{n}]\,\prod_{r=0}^{\infty}(1-vq^{-r})^{-1}\,\prod_{j=1}^{k} \mathcal{L}_{v^{\deg \Theta_{j}},q^{\deg \Theta_{j}}}[\lambda_{j}]\,.$$
De plus, si $f$ est une fonction holomorphe de rayon de convergence $1$ et de série de Taylor convergente en $1$, alors $f(1)=\lim_{n \to \infty}[v^{n}]\,f(v)\,(1-v)^{-1}$. Par conséquent :
\begin{theorem}[Dudko-Fulman, \cite{Ful06,Dud08}]
Pour toutes partitions $\lambda_{j}$ et toutes orbites $\Theta_{j} \in L/\Gal$,
$$\lim_{n \to \infty} \mathrm{Pl}_{\GL(n,\For_{q})}[\bblambda(\Theta_{1})=\lambda_{1},\ldots,\bblambda(\Theta_{k})=\lambda_{k}]=\prod_{j=1}^{k} \mathcal{L}_{1,q^{\deg \Theta_{j}}}[\lambda_{j}]\,.$$
\end{theorem}
\noindent Ainsi, sous la  $\GL(n,\For_{q})$-mesure de Plancherel, les coordonnées $\bblambda(\Theta_{1}),\ldots,\bblambda(\Theta_{k})$ sont asymptotiquement indépendantes, et asymptotiquement réparties suivant des mesures de Schur dont les paramètres ne dépendent que des $q^{\deg \Theta_{j}}$. En particulier, si une orbite $\Theta$ est fixée, alors la partition aléatoire $\bblambda(\Theta)$ avec $\bblambda$ répartie suivant la $\GL(n,\For_{q})$-mesure de Plancherel <<~reste de taille bornée~>> lorsque $n$ tend vers l'infini, et une étude asymptotique semblable à celle menée dans la première partie du mémoire n'a pas de sens.

\chapter{Algèbres d'Iwahori-Hecke et $q$-mesure de Plancherel}\label{iwahori}

Dans ce chapitre, nous présentons finalement la \textbf{q-mesure de Plancherel}, et nous expliquons son lien avec la théorie des \textbf{algèbres d'Iwahori-Hecke}. L'idée initiale est d'étudier la mesure de Plancherel d'un module de $\GL(n,\For_{q})$ obtenu par induction parabolique à partir d'un caractère $\zeta$ d'un tore maximal $T$, et on s'interesse au cas le plus simple, \emph{i.e.}, $T=k^{n}$ est le sous-groupe des matrices diagonales et $\zeta$ est le caractère trivial $t \mapsto 1$. La mesure sous-jacente a été décrite pour la première fois par Kerov dans \cite{Ker92}, et elle a été étudiée en détail par E. Strahov (\cite{Stra08}) ; dans la section \ref{qplancherelprocess}, nous donnons plusieurs expressions pour cette mesure, et nous la relions à un modèle combinatoire et à divers processus aléatoires. \bigskip

La poissonisée des $q$-mesures de Plancherel est de nouveau une mesure de Schur, et les noyaux correspondants s'expriment à l'aide de fonctions de Bessel déformées, \emph{cf.} la proposition \ref{qbessel}. Malheureusement, on ne peut pas comme dans le paragraphe \ref{bdjdeterminantal} en déduire l'asymptotique des $q$-mesures de Plancherel, car le comportement à l'infini de ces $q$-fonctions de Bessel n'est pas connu. \bigskip\bigskip

Pour cette raison, l'étude asymptotique sera menée en utilisant des techniques d'observables de diagrammes, voir le chapitre \ref{qplancherelmeasure}. En particulier, nous construirons plus loin une quantification de l'algèbre $\obs$, \emph{cf.} le paragraphe \ref{quantum}. Cette déformation s'inspire de la théorie des algèbres d'Hecke (\emph{cf.} \cite{Mat99,Iwa64}) ; nous en rappelons les points principaux dans le paragraphe \ref{hecke}, et nous rappelons également dans la section \ref{ram} la théorie des caractères de l'algèbre d'Hecke de type A, en s'inspirant pour l'essentiel de l'article \cite{Ram91}. Le point important à retenir est le suivant : la \textbf{formule de Ram} \ref{ramformula} exprime les caractères de l'algèbre $\IH_{q}(\sym_{n})$ en fonction des caractères du groupe $\sym_{n}$, et ceci permettra un changement de base dans $\obs$ et un calcul des espérances des observables $\varSigma_{\lambda}$ sous la $q$-mesure de Plancherel.\bigskip

\section{$q$-mesures et $q$-processus de Plancherel}\label{qplancherelprocess}
Soit $n$ un entier positif, et $q$ la puissance d'un nombre premier ; on s'intéresse à la mesure de Plancherel du $\GL(n,\For_{q})$-module $\C[\GL(n,\For_{q})/\GB(n,\For_{q})]=R^{([\mathbb{1}]:1^{n})}$. Le cardinal du groupe des matrices triangulaires supérieures inversibles est $(q-1)^{n}\,q^{\frac{n(n-1)}{2}}$, donc
$$\dim \C[\GL(n,\For_{q})/\GB(n,\For_{q})] = \prod_{i=1}^{n} \frac{q^{i}-1}{q-1}\,.$$
Dans ce qui suit, nous utiliserons librement les notations des $q$-entiers (voir \cite{CK02}) : $\{n\}_{q}=(q^{n}-1)/(q-1)$ pour tout entier $n \geq 1$, et $\{n!\}_{q}=\prod_{i=1}^{n} \{i\}_{q}$. On retrouve les entiers et les factorielles usuelles pour $q = 1$, et $\dim \C[\GL(n,\For_{q})/\GB(n,\For_{q})]=\{n!\}_{q}$. En un certain sens, le module que l'on étudie est donc une généralisation  de l'algèbre du groupe symétrique\footnote{On voit donc $\sym_{n}$ comme le groupe linéaire d'un espace vectoriel de dimension $n$ sur le corps à un élément.}.\bigskip
\bigskip

Nous avons vu dans le chapitre précédent que les composantes irréductibles du module $\C[\GL(n,\For_{q})/\GB(n,\For_{q})]$ étaient indexées par les partitions $\lambda \in \ym_{n}$, et que la multiplicité du module unipotent $U^{\lambda}(q)$ était toujours $\dim \lambda$. On conserve cette indexation dans ce chapitre, à une conjugaison des diagrammes près : ainsi,
$$ \dim U^{\lambda}(q)=D_{\lambda}(q)=q^{b(\lambda)}\,\frac{\prod_{i=1}^{n} q^{i}-1}{\prod_{(i,j) \in \lambda} q^{h(i,j)}-1}=q^{b(\lambda)}\,\frac{\{n!\}_{q}}{ \prod_{(i,j) \in \lambda} \{h(i,j)\}_{q}}$$
--- par rapport à la formule du chapitre précédent, on a changé $q^{b(\lambda')}$ en $q^{b(\lambda)}$. Comme $\dim \lambda=\dim \lambda'$, le changement d'indexation est compatible avec la décomposition du module $\C[\GL(n,\For_{q})/\GB(n,\For_{q})]$.
\begin{definition}[$q$-mesure de Plancherel]
On appelle $q$-mesure de Plancherel la mesure de probabilité sur les partitions associée à la décomposition du module $\C[\GL(n,\For_{q})/\GB(n,\For_{q})]$ en composantes unipotentes irréductibles $U^{\lambda}(q)$. Autrement dit,
$$M_{n,q}[\lambda]= \frac{\dim \lambda \times \dim U^{\lambda}(q)}{|\GL(n,\For_{q})/\GB(n,\For_{q})|}=\frac{\dim \lambda\times D_{\lambda}(q)}{\{n!\}_{q}}\,=\frac{n!\,\,q^{b(\lambda)}}{\prod_{(i,j) \in \lambda} h(i,j)\,\{h(i,j)\}_{q}}\,.$$
\end{definition}

\begin{example}
Pour $n=4$, la $q$-mesure de Plancherel a pour valeurs :
\begin{align*}
&M_{4,q}(4)= 1/(q^6 + 3q^5 + 5q^4 + 6q^3 + 5q^2 + 3q + 1)\\
&M_{4,q}(3,1)=3q/(q^4 + 2q^3 + 2q^2 + 2q + 1)\\
&M_{4,q}(2,2)=2q^2/(q^4 + 3q^3 + 4q^2 + 3q + 1)\\
&M_{4,q}(2,1,1)=3q^3/(q^4 + 2q^3 + 2q^2 + 2q + 1)\\
&M_{4,q}(1,1,1,1)=q^6/(q^6 + 3q^5 + 5q^4 + 6q^3 + 5q^2 + 3q + 1)
\end{align*}
On retrouve pour $q=1$ les valeurs de la mesure de Plancherel $M_{4}$ données page \pageref{plancherelprocess}.
\end{example}\bigskip

\noindent En plus de la $q$-formule des équerres, il existe une autre expression de $D_{\lambda}(q)$ due à Steinberg (voir \cite[théorème 10.5.2]{GP00}), et qui a l'avantage de se généraliser au cas du type B, \emph{cf.} le chapitre \ref{arikikoike}. Si $\lambda=(\lambda_{1},\lambda_{2},\ldots,\lambda_{r})$, notons $A$ l'ensemble des parts de $\lambda+\delta$, où $\delta$ est la partition en escalier $(r-1,r-2,\ldots,0)$. Ainsi,
$$A=\{\lambda_{i}+r-i\}_{i \in \lle 1,r\rre}\,.$$
On peut montrer que $D_{\lambda}(q)$ a pour expression :
$$D_{\lambda}(q)=\frac{\prod_{i=1}^{n} q^{i}-1}{q^{\frac{r(r-1)(r-2)}{6}} } \,\frac{\prod_{a>a' \in A}q^{a}-q^{a'}}{\prod_{a \in A}\prod_{i=1}^{a} q^{i}-1}\,.$$
Par exemple, $D_{2,2}(q)=\frac{(q^{4}-1)(q^{3}-1)(q^{2}-1)(q-1)\,\,(q^{3}-q^{2})}{q^{0}\,\,(q^{3}-1)(q^{2}-1)(q-1)(q^{2}-1)(q-1)}=q^{4}+q^{2}$.
\bigskip
\bigskip

Dans \cite{Stra08}, E. Strahov utilise une troisième expression de $D_{\lambda}(q)$ pour donner un modèle combinatoire de la $q$-mesure de Plancherel. Si $\sigma$ est une permutation de $\sym_{n}$, on appelle \textbf{descente} de $\sigma$ un entier $i \in \lle 1,n \rre$ tel que $\sigma(i)>\sigma(i+1)$. Autrement dit, la $i$-ième case du ruban standard associé à $\sigma$ est au-dessus de la $i+1$-ième case. L'\textbf{indice majeur} de $\sigma$ est la somme de ses descentes ; par exemple, si $\sigma=352614$, alors $\imaj(\sigma)=2+4=6$. D'autre part, si $T$ est un tableau standard de taille $n$, on dit que $i$ est une descente de $T$ si $i+1$ est dans une ligne strictement au-dessus de la ligne de $i$ dans $T$. Par exemple, les descentes du tableau standard
\begin{center}
\vspace{-2mm}
\young(5,36,124)
\vspace{1mm}
\end{center}
sont $2$ et $4$. Comme précédemment, l'indice majeur d'un tableau standard est la somme de ses descentes. Une récurrence permet de montrer que pour toute permutation $\sigma$, l'ensemble des descentes de $\sigma$ est aussi l'ensemble des descentes du tableau standard $Q(\sigma)=P(\sigma^{-1})$ qui lui est associé par la correspondance RSK, voir \cite[volume 2, chapitre 7]{Stan91}. D'autre part, le degré générique d'une partition $\lambda$ s'écrit
$$D_{\lambda}(q) = \sum_{T \in \mathrm{Std}(\lambda)} q^{\imaj(T)}\,.$$
\begin{example}
Si $\lambda=(2,2)$, les deux tableaux standards de forme $\lambda$ sont
$$\young(34,12)\qquad\text{et}\qquad\young(24,13)\,\,\,,$$
et ils ont pour indices majeurs $2$ et $1+3=4$, d'où l'identité $D_{2,2}(q)=q^{4}+q^{2}$.
\end{example}\bigskip\bigskip

Enfin, si $\sigma$ est une permutation de taille $n$ et si $V(\sigma)$ est l'ensemble des permutations de taille $n+1$ obtenues en insérant $n+1$ dans le mot de $\sigma$, alors la fonction $$\tau \mapsto \imaj(\tau)-\imaj(\sigma)$$
est une bijection entre $V(\sigma)$ et $\lle 0,n\rre$ --- plus précisément, si $\sigma$ a $d$ descentes, alors $\tau\in V(\sigma)$ a $d$ ou $d+1$ descentes, et la fonction précédente va de $d+1$ à $n$ pour les permutations $\tau$ avec $d+1$ descentes, et de $d$ à $0$ pour les permutations $\tau$ avec $d$ descentes, étant entendu qu'on liste les voisins de $\sigma$ suivant la position de la lettre $n+1$ dans leur mot (\emph{cf.} \cite{Gup78}). Par suite :
\begin{proposition}[$q$-processus de Knuth et $q$-mesure de Plancherel]\label{qmodel}
Pour tout entier $n$, la fonction
$$P_{n,q} : \sigma \mapsto \frac{q^{\imaj(\sigma)}}{\{n!\}_{q}}$$
est une mesure de probabilité sur $\sym_{n}$, et c'est la $n$-ième marginale du processus de Knuth déformé $(\sigma^{(n)})_{n \in \N}$ défini par les probabilités de transition
$$p(\sigma \in \sym_{k},\tau \in \sym_{k+1})=\begin{cases}
\frac{q^{\imaj(\tau)-\imaj(\sigma)}}{\{k+1\}_{q}}&\text{si }\tau \in V(\sigma),\\
0&\text{sinon}.
\end{cases}$$
De plus, la mesure image de $P_{n,q}$ par la correspondance RSK est la $q$-mesure de Plancherel $M_{n,q}$.
\end{proposition}
\begin{proof}
Comme $\{k+1\}_{q}=1+q+q^{2}+\cdots+q^{k}$, la fonction de transition $p$ est bien une mesure de probabilité :
$$\forall \sigma,\,\,\,\sum_{\tau \in V(\sigma)} p(\sigma,\tau)= \frac{\sum_{\tau \in V(\sigma)} q^{\imaj(\tau)-\imaj(\sigma)}}{1+q+q^{2}+\cdots+q^{k}}=1\,.$$
La $n$-ième marginale du processus markovien $(\sigma^{(n)})_{n \in \N}$ associé à ces probabilités de transition est évidemment la mesure $P_{n,q}$. De plus, si $\lambda$ est une partition de taille $n$, alors :
\begin{align*}\mathrm{RSK}_{\star}P_{n,q}[\lambda]&=\sum_{\sigma \,:\,Q(\sigma) \in \mathrm{Std}(\lambda)}P_{n,q}[\sigma]=\sum_{T \in \mathrm{Std}(\lambda)}\sum_{Q(\sigma)=T}\frac{q^{\imaj(\sigma)}}{\{n!\}_{q}}\\
&=\frac{1}{\{n!\}_{q}}\sum_{T \in \mathrm{Std}(\lambda)}q^{\imaj(T)}\times\card\{\sigma\,\,|\,\,Q(\sigma)=T\}\\
&=\frac{\dim \lambda}{\{n!\}_{q}}\sum_{T \in \mathrm{Std}(\lambda)}q^{\imaj(T)}=\frac{\dim\lambda \times D_{\lambda}(q)}{\{n!\}_{q}}=M_{n,q}[\lambda]\,.
\end{align*}
Ainsi, la mesure non uniforme $P_{n,q}$ sur les permutations fournit un modèle combinatoire de la $q$-mesure de Plancherel : c'est l'image par la correspondance RSK de cette mesure.\end{proof}
\bigskip

Il y a toutefois une différence avec la construction du paragraphe \ref{plancherelprocess} : la projection par RSK du $q$-processus de Knuth n'est plus un processus markovien. En effet, si tel était le cas, alors pour tout tableau standard $T$ de taille $n+1$, notant $t$ le tableau $T$ moins la case $n+1$, et $\lambda^{(0)} \nearrow \lambda^{(1)} \nearrow \cdots \nearrow \lambda^{(n+1)}$ la suite de partitions correspondant au tableau $T$, on aurait
$$p(\lambda^{(n)},\lambda^{(n+1)})=\frac{\proba[\lambda(\sigma_{0})=\lambda^{(0)},\ldots,\lambda(\sigma_{n+1})=\lambda^{(n+1)}]}{\proba[\lambda(\sigma_{0})=\lambda^{(0)},\ldots,\lambda(\sigma_{n})=\lambda^{(n)}]}=\frac{q^{\imaj(T) - \imaj(t)}}{\{n+1\}_{q}}\,\frac{\dim \lambda^{(n+1)}}{\dim \lambda^{(n)}}\,.$$
On obtient une contradiction en prenant deux tableaux standards $T$ et $T'$ associés aux mêmes formes $\lambda^{(n)}$ et $\lambda^{(n+1)}$, mais tels que $\imaj(T)-\imaj (t) \neq \imaj(T')-\imaj(t')$ ; par exemple, les tableaux
$$\young(3,25,14) \qquad\text{et}\qquad\young(4,35,12)\,\,.$$
Le modèle combinatoire de la $q$-mesure de Plancherel ne fournit donc pas par correspondance RSK une généralisation raisonnable du processus de Plancherel. Il existe toutefois une généralisation markovienne de ce processus, mais on doit la définir directement au niveau du graphe Young $\ym$. Ainsi, si l'on pose 
$$p_{q}(\lambda,\Lambda)=\begin{cases}
\frac{1}{\{|\Lambda|\}_{q}}\,\frac{D_{\Lambda}(q)}{D_{\lambda}(q)}&\text{si }\lambda \nearrow \Lambda, \\
0&\text{sinon},
\end{cases}$$
alors on a défini les probabilités de transition d'un processus markovien sur $\ym$ de lois marginales les $q$-mesures de Plancherel, et lorsque $q$ tend vers $1$, on retrouve le processus de Plancherel usuel. On appelle \textbf{q-processus de Plancherel} ces nouveaux processus. Leurs probabilités de transition $p_{q}$ s'expriment facilement à l'aide de la $q$-formule des équerres :
$$ p_{q}(\lambda,\Lambda)=q^{b(\Lambda)-b(\lambda)}\,\frac{\prod_{(i,j) \in \lambda}\{h(i,j)\}_{q}}{\prod_{(i,j) \in \Lambda}\{h(i,j)\}_{q}}=\frac{\prod_{i=1}^{v-1}\{x_{k}-y_{i}\}_{q}}{\prod_{i \neq k}\{x_{k}-x_{i}\}_{q}}$$
et ils sont donc tout aussi aisément programmables que le processus de Plancherel usuel. Des exemples numériques seront proposés au début du chapitre \ref{qplancherelmeasure} ; sur la figure \ref{qplanprocessfig}, on a indiqué les degrés génériques et les $q$-probabilités de transition pour les quatre premiers niveaux du graphe de Young.
\figcapt{\footnotesize{
\psset{unit=1mm}\pspicture(-70,-105)(70,5)
\psline(-2.5,0)(2.5,0)(2.5,5)(-2.5,5)(-2.5,0)
\rput(10,2.5){$ D=1$}
\psline[linecolor=red]{->}(0,-1.5)(20,-13.5)
\rput(16,-7){\textcolor{red}{$\frac{1}{q+1}$}}
\psline[linecolor=red]{->}(0,-1.5)(-20,-11.5)
\rput(-18,-7){\textcolor{red}{$\frac{q}{q+1}$}}
\psline(15,-20)(25,-20)(25,-15)(15,-15)(15,-20) \psline(20,-20)(20,-15)
\rput(33,-17.5){$D=1$}
\psline(-17.5,-12.5)(-22.5,-12.5)(-22.5,-22.5)(-17.5,-22.5)(-17.5,-12.5) \psline(-17.5,-17.5)(-22.5,-17.5)
\rput(-30,-17.5){$D=q$}
\psline(-5,-35)(5,-35)(5,-40)(0,-40)(0,-45)(-5,-45)(-5,-35) \psline(0,-35)(0,-40)(-5,-40)
\psline[linecolor=red]{->}(-20,-24)(-40,-31)
\rput(-42,-27){\textcolor{red}{$\frac{q^2}{q^2+q+1}$}}
\psline[linecolor=red]{->}(-20,-24)(-1,-34)
\rput(-8,-24.5){\textcolor{red}{$\frac{q+1}{q^2+q+1}$}}
\psline[linecolor=red]{->}(20,-21.5)(1,-34)
\rput(7,-24){\textcolor{red}{$\frac{q^2+q}{q^2+q+1}$}}
\psline[linecolor=red]{->}(20,-21.5)(37.5,-34)
\rput(37,-27){\textcolor{red}{$\frac{1}{q^2+q+1}$}}
\rput(10.5,-43){$D=q^2+q$}
\psline[linecolor=red]{->}(-2.5,-46.5)(-2.5,-58.5)
\rput(10,-53){\textcolor{red}{$\frac{q^3+q}{q^4+2q^3+2q^2+2q+1}$}}
\psline(-37.5,-32.5)(-37.5,-47.5)(-42.5,-47.5)(-42.5,-32.5)(-37.5,-32.5) \psline(-37.5,-37.5)(-42.5,-37.5) \psline(-37.5,-42.5)(-42.5,-42.5)
\rput(-30,-40){$D=q^3$}
\psline[linecolor=red]{->}(-43.5,-40)(-50,-59)
\rput(-55.5,-50){\textcolor{red}{$\frac{q^3}{q^3+q^2+q+1}$}}
\psline(-47.5,-60)(-52.5,-60)(-52.5,-80)(-47.5,-80)(-47.5,-60)
\psline(-47.5,-65)(-52.5,-65) \psline(-47.5,-70)(-52.5,-70) \psline(-47.5,-75)(-52.5,-75)
\rput(-59,-70){$D=q^6$}
\psline(30,-35)(45,-35)(45,-40)(30,-40)(30,-35) \psline(35,-35)(35,-40) \psline(40,-35)(40,-40)
\rput(53,-37.5){$D=1$}
\psline(-5,-60)(5,-60)(5,-70)(-5,-70)(-5,-60) \psline(-5,-65)(5,-65) \psline(0,-60)(0,-70)
\rput(15.5,-67.5){$D=q^4+q^2$}
\psline(-15,-80)(-25,-80)(-25,-95)(-20,-95)(-20,-85)(-15,-85)(-15,-80) \psline(-20,-80)(-20,-85)(-25,-85) \psline(-20,-90)(-25,-90)
\psline[linecolor=red]{->}(-2.5,-46.5)(-19,-79)
\psline[linecolor=red]{->}(-2.5,-46.5)(-10,-75)(22,-84)
\psline[linecolor=red]{->}(-43.5,-40)(-45.5,-56)(-21,-79)
\psframe*[linecolor=white,fillcolor=white](-40,-75)(-20,-69)
\rput(-30,-72){\textcolor{red}{$\frac{q^2+q+1}{q^3+q^2+q+1}$}}
\psframe*[linecolor=white,fillcolor=white](-28,-63.5)(-8.5,-57)
\rput(-18,-60){\textcolor{red}{$\frac{q^4+q^3+q^2}{q^4+2q^3+2q^2+2q+1}$}}
\psline(15,-85)(30,-85)(30,-90)(20,-90)(20,-95)(15,-95)(15,-85) \psline(20,-85)(20,-90)(15,-90) \psline(25,-85)(25,-90)
\rput(34,-93){$D=q^3+q^2+q$}
\rput(-24,-98){$D=q^5+q^4+q^3$}
\psline(35,-60)(55,-60)(55,-65)(35,-65)(35,-60) \psline(40,-60)(40,-65) \psline(45,-60)(45,-65) \psline(50,-60)(50,-65)
\psline[linecolor=red]{->}(37.5,-41)(45,-59)
\psline[linecolor=red]{->}(37.5,-41)(23.5,-84)
\rput(50,-50){\textcolor{red}{$\frac{1}{q^3+q^2+q+1}$}}
\rput(63,-62.5){$D=1$}
\psframe*[linecolor=white,fillcolor=white](25,-60)(35,-54)
\rput(29,-57){\textcolor{red}{$\frac{q^3+q^2+q}{q^3+q^2+q+1}$}}
\rput(0,-82){\textcolor{red}{$\frac{q^2+q+1}{q^4+2q^3+2q^2+2q+1}$}}
\endpspicture}}{Degrés génériques et $q$-probabilités de transition sur les quatre premiers niveaux du graphe de Young.\label{qplanprocessfig}}{Degrés génériques et $q$-probabilités de transition}\bigskip

Notons que les formules pour les probabilités de transition $p_{q}$ et les $q$-mesures de Plancherel $M_{n,q}$ font sens pour tout paramètre $q>0$, et pas seulement pour les puissances de nombres premiers. De plus, comme $$h(\lambda)=\sum_{(i,j) \in \lambda}h(i,j)=b(\lambda)+b(\lambda')+|\lambda|$$
pour toute partition $\lambda$, on voit facilement que l'image de la $q$-mesure de Plancherel (resp., du $q$-processus de Plancherel) par la conjugaison de diagrammes est la $q^{-1}$-mesure de Plancherel (resp., le $q^{-1}$-processus de Plancherel). Cette symétrie permet d'étudier indifféremment les $q$-processus de Plancherel pour $q \in \,]1,+\infty[$ ou $q \in \,]0,1[$ ; ainsi, dans le chapitre suivant, nous nous concentrerons sur le cas $q<1$, même si le cas le plus naturel d'un point de vue algébrique est $q>1$. Notons qu'au niveau des permutations, la symétrie $q\leftrightarrow q^{-1}$ correspond à la multiplication par l'élément maximal $i \mapsto n+1-i$.
\bigskip
\bigskip

Dans \cite{Stra08}, E. Strahov esquisse une étude asymptotique des $q$-mesures de Plancherel en rempla\c cant le $q$-processus de Plancherel par un analogue différentiel dont une $q$-fonction génératrice satisfait une $q$-équation de Burgers ; il utilise également des $q$-observables de diagrammes. Malheureusement, les $q$-observables employées dans \cite{Stra08} n'appartiennent pas à l'algèbre $\obs$, et d'ailleurs elles ne spécialisent pas en les observables usuelles correspondantes lorsque le paramètre $q$ tend vers $1$. D'autre part, l'étude asymptotique menée par Strahov met en jeu des $q$-formes limites $\Omega_{q}$ extrêmement peu explicites (on ne sait pas du tout les dessiner...), et qui correspondent à une renormalisation simultanée des diagrammes et du paramètre $q$ par un facteur $\sqrt{n}$ ; nous verrons dans le chapitre suivant une forme beaucoup plus simple de ces résultats. Enfin, la proposition 8.2.1 de \cite{Stra08} ne donne rien de rigoureux concernant l'asymptotique des $q$-processus de Plancherel --- ou du moins nous ne voyons pas comment l'exploiter.
\bigskip\bigskip

Pour toutes ces raisons, nous avons été amené à abandonner l'approche <<~différentielle~>> dans le problème de l'asymptotique des $q$-mesures de Plancherel. Dans le chapitre \ref{qplancherelmeasure}, nous verrons que les observables de diagrammes (usuelles) permettent de résoudre entièrement le problème ; avant cela, voyons ce qu'il en est de l'approche <<~déterminantale~>>. Fixons un paramètre $q$ strictement inférieur à $1$. On a vu dans le chapitre précédent que 
$$s_{\lambda}(q,q^{2},q^{3},\ldots)=\frac{q^{-b(\lambda')}}{\prod_{(i,j) \in \lambda}q^{-h(i,j) }-1}=\frac{q^{h(\lambda)-b(\lambda')}}{\prod_{(i,j) \in \lambda}1-q^{h(i,j)}}=\frac{q^{|\lambda|+b(\lambda)}}{\prod_{(i,j) \in \lambda}1-q^{h(i,j)}}\,.$$
Par conséquent, $s_{\lambda}(1-q,q(1-q),q^{2}(1-q),\ldots)$ vaut 
$$\frac{q^{b(\lambda)}}{\prod_{(i,j) \in \lambda}{\{h(i,j)\}_{q}}}=\frac{D_{\lambda}(q)}{\{n!\}_{q}}\,.$$
Posons $b=(\sqrt{\theta}(1-q),\sqrt{\theta}(1-q)q,\sqrt{\theta}(1-q)q^{2},\ldots)$ ; cette suite correspond aux paramètres de Miwa
$$t_{k}'=\frac{1}{k}\,p_{k}(b)=\frac{\theta^{k/2}}{k}\,\,\frac{(1-q)^{k}}{1-q^{k}}\,.$$
Alors, si $t=(\sqrt{\theta},0,0,\ldots)$, la mesure de Schur associée aux paramètres $t$ et $t'$ est
\begin{align*}\mathcal{S}(\lambda)&=\frac{1}{Z}\,\left(\theta^{\frac{|\lambda|}{2}} \,\frac{\dim \lambda}{|\lambda|!}\right)\times\left(\theta^{\frac{|\lambda|}{2}}\,s_{\lambda}(1-q,q(1-q),q^{2}(1-q),\ldots)\right)\\
&=\frac{\theta^{|\lambda|}}{\E^{\lambda}|\lambda|!}\,\,\frac{\dim \lambda\times D_{\lambda}(q)}{\{n!\}_{q}}=\mathcal{P}(\theta)(|\lambda|)\,M_{|\lambda|,q}(\lambda)\,,
\end{align*}
c'est-à-dire que c'est la poissonisée de paramètre $\theta$ des $q$-mesures de Plancherel. Nous noterons cette mesure $M_{\mathcal{P}(\theta),q}$. Lorsque $q$ tend vers $1$, les paramètres $t$ et $t'$ tendent tous les deux vers $(\sqrt{\theta},0,0,\ldots)$, qui est le paramètre de la poissonisée de la mesure de Plancherel standard ; on a ainsi généralisé la construction du paragraphe \ref{bdjdeterminantal}.
\bigskip
\bigskip

Compte tenu du théorème \ref{schurmeasurearedeterminantal}, le processus ponctuel $\mathcal{D}_{\star}(M_{\mathcal{P}(\theta),q})$ est déterminantal ; calculons son noyau. La fonction $J$ introduite dans le paragraphe \ref{schurmeasure} s'écrit ici
\begin{align*}J_{q}(x)&=\exp\left(\sum_{n\geq 1}t_{n}\,x^{n}-\sum_{n \geq 1}t_{n}'\,x^{-n}\right)=\E^{\sqrt{\theta}x}\left(\prod_{i}1-b_{i}x^{-1}\right)\\
&=\E^{\sqrt{\theta}x}\left(\prod_{i=0}^{\infty}1-q^{i}\sqrt{\theta}(1-q)x^{-1}\right)=\E^{\sqrt{\theta}x}\,(\sqrt{\theta}(1-q)x^{-1};q)_{\infty}\end{align*}
en utilisant les notations de Pochhammer. Compte tenu des résultats du paragraphe \ref{schurmeasure}, la fonction génératrice du noyau $K_{\mathcal{D}_{\star}(M_{\mathcal{P}(\theta),q})}$ est donc
 $$A_{q}(z,w)=\frac{\sqrt{zw}}{z-w}\,\E^{\sqrt{\theta}(z-w)}\,\frac{(\sqrt{\theta}(1-q)z^{-1};q)_{\infty}}{(\sqrt{\theta}(1-q)w^{-1};q)_{\infty}}$$
et le $q$-calcul permet de développer les symboles de Pochhammer infinis en séries entières 
$$(t;q)_{\infty}=\sum_{n=0}^{\infty} \frac{(-1)^{n} \,q^{\binom{n}{2}}}{(q;q)_{n}}\,t^{n}\qquad;\qquad\frac{1}{(t;q)_{\infty}}=\sum_{n=0}^{\infty} \frac{t^{n}}{(q;q)_{n}}\,,$$
voir \cite{CK02}. Ainsi, $A_{q}(z,w)$ est la somme quintuple
$$A_{q}(z,w)=\sum_{k,l \in \Z'}K_{\mathcal{D}}(k,l)\,z^{k}w^{-l}=\sum_{r,s,t,u,v} \frac{(-1)^{t+u}\,\theta^{\frac{s+t+u+v}{2}} \,q^{\binom{u}{2}}}{s!\,t!\,\{u!\}_{q}\,\{v!\}_{q}}\,z^{s-u-r-1/2}\,w^{t-v+r+1/2}\,,$$
les indices $r,s,t,u,v$ appartenant à $\N$.\bigskip
\bigskip

Fixons deux demi-entiers $x$ et $y$, et introduisons le paramètre $n=\frac{s+t+u+v}{2}-\frac{x+y+1}{2}$. Le noyau $K_{\mathcal{D}}(x,y)$ est la restriction de la somme précédente à l'ensemble d'indices $\{r,s,t,u,v\}$ tels que $x=s-u-r-1/2$ et $y=v-t-r-1/2$. Par conséquent, $n=r+t+u$, et on peut écrire $K_{\mathcal{D}}(x,y)$ sous la forme :
\begin{align*}&\phantom{=}\sum_{n=0}^{\infty}\sum_{r+t+u=n}\frac{(-1)^{t+u}\,\theta^{n+\frac{x+y+1}{2}}\,q^{\binom{u}{2}}}{s!\,t!\,\{u!\}_{q}\,\{v!\}_{q}}\\
&=\sum_{n=0}^{\infty} (-1)^{n}\,\theta^{n+\frac{x+y+1}{2}} \left(\sum_{r+t+u=n} \frac{(-1)^{r}\,q^{\binom{u}{2}}}{t!\,(n+x+1/2-t)!\,\{u!\}_{q}\,\{n+y+1/2-u!\}_{q}}\right)\\
&=\sum_{n=0}^{\infty} \frac{(-1)^{n}\,\theta^{n+\frac{x+y+1}{2}}}{(x+n+\frac{1}{2})!\,\{y+n+\frac{1}{2}!\}_{q}}\left(\sum_{r+t+u=n} \binom{n+x+\frac{1}{2}}{t}\binom{n+y+\frac{1}{2}}{u}_{\!\!\!q}\,(-1)^{r}\,q^{\binom{u}{2}}\right).\end{align*}
Lorsque $q=1$, la somme sur les indices $r,t,u$ se réduit par la formule de Vandermonde à $ \binom{x+y+2n}{n}$, et la somme 
$$\sum_{n=0}^{\infty} \frac{(-1)^{n}\,\theta^{n+\frac{x+y+1}{2}}}{(x+n+\frac{1}{2})!\,(y+n+\frac{1}{2})!}\binom{x+y+2n}{n}$$
est bien le développement en série de $\sqrt{\theta}\,\frac{J_{x-1/2}(2\sqrt{\theta})\,J_{y+1/2}(2\sqrt{\theta})-J_{x+1/2}(2\sqrt{\theta})J_{y-1/2}(2\sqrt{\theta})}{x-y}$. Nos formules généralisent donc bien celles données dans le paragraphe \ref{bdjdeterminantal}. Finalement, si l'on introduit les déformations de fonctions de Bessel
$$J^{1}_{r,q}(z)=\sum_{n=0}^{\infty} \frac{(-1)^{n}\,q^{\binom{n}{2}}}{\{n!\}_{q}\,r+n!}\,\left(\frac{z}{2}\right)^{r+2n}\qquad;\qquad J^{2}_{r,q}(z)=\sum_{n=0}^{\infty} \frac{(-1)^{n}}{n! \,\{r+n!\}_{q}}\,\left(\frac{z}{2}\right)^{r+2n}\,,$$
alors on obtient une expression semblable à celle donnée page \pageref{asymptoticpoisson} dans le cas $q=1$ :
\begin{proposition}[Noyau associé à la poissonisée des $q$-mesures de Plancherel]\label{qbessel}
\'Etant donné un paramètre $q<1$, le processus ponctuel $\mathcal{D}_{\star}(M_{\mathcal{P}(\theta),q})$ est déterminantal, de noyau
$$K_{q}(x,y)=\sum_{n=0}^{\infty} J^{1}_{x+n+1/2,q}(2\sqrt{\theta})\,J^{2}_{y+n+1/2,q}(2\sqrt{\theta})\,.$$
\end{proposition}
\begin{proof}
En développant en séries les deux fonctions de Bessel déformées, on obtient une somme triple qui est exactement la deuxième ligne de la formule donnée précédemment pour $K_{\mathcal{D}}(x,y)$.
\end{proof}\bigskip

L'analogie avec les raisonnements du chapitre \ref{matrix} s'arrête malheureusement ici, car on connaît pas l'asymptotique de $K_{q}(x,y)$ lorsque $\theta$ tend vers l'infini. Les fonctions de Bessel déformées que nous avons obtenues sont à notre connaissance inédites ; en particulier, elles diffèrent sensiblement des constructions usuelles de $q$-fonctions de Bessel, voir par exemple \cite{Jac05}. Le noyau $K_{q}(x,y)$ est \emph{a priori} non intégrable, c'est-à-dire qu'on ne sait pas comme dans le cas $q=1$ l'écrire comme combinaison linéaire finie de produits de fonctions $f(x)\,g(y)$. Et même si c'était le cas, l'asymptotique des fonctions $J_{r,q}^{1}(z)$ et $J_{r,q}^{2}(z)$ lorsque $z$ tend vers l'infini semble difficile à déterminer, essentiellement parce qu'il s'agit de sommes alternées. Nous devons donc également rejeter l'approche déterminantale ; \emph{a posteriori}, nous verrons d'ailleurs qu'elle avait peu de chances d'aboutir, puisque la simple densité $K_{q}(x,x)$ doit dégénérer après renormalisation en une somme de Dirac en les points $(1-q),q(1-q),\ldots$ voir le chapitre suivant et le théorème \ref{firstasymptoticqplancherel}.\bigskip

\section{Algèbre d'Iwahori-Hecke d'un groupe de Chevalley}\label{hecke}

Lorsqu'on a étudié la mesure de Plancherel standard à l'aide de l'algèbre d'observables $\obs$, le point de départ de la preuve était l'existence d'une base d'observables $(\varSigma_{\mu})_{\mu \in \ym}$ telle que 
$$M_{n}[\varSigma_{\mu}]=n^{\downarrow k}\,\mathbb{1}_{\mu=1^{k}}\,,$$
\emph{cf.} page \pageref{startidentityplancherel}. Dans la suite de ce chapitre, on cherche une base d'observables $(\varSigma_{\mu,q})_{\mu \in \ym}$ dont on peut calculer aussi aisément les espérances, mais cette fois-ci sous les $q$-mesures de Plancherel ; ces observables seront à la base de l'étude asymptotique des mesures $M_{n,q}$. La théorie algébrique sous-jacente est celle des \textbf{algèbres d'Hecke}, et dans ce contexte, on pourra donner plus loin une interprétation de la $q$-mesure de Plancherel analogue à celle donnée au début du paragraphe \ref{plancherelprocess}, c'est-à-dire en termes d'isométries d'espaces de probabilité non commutatifs.\bigskip
\bigskip

Les algèbres d'Iwahori-Hecke peuvent être étudiées dans le contexte général des groupes de Chevalley (\emph{cf.} \cite[\S1]{Car92}), c'est-à-dire les groupes algébriques $G^{F}$ définis sur un corps fini $k=\For_{q}$, avec $G$ groupe réductif défini sur $\overline{\For_{q}}$.  Nous adopterons ce point de vue général ; ceci nous permettra d'adapter les résultats relatifs au $\GL(n,\For_{q})$-module $\C[\GL(n,\For_{q})/\GB(n,\For_{q})]$ à d'autres cas, par exemple celui du $\Sp(2n,\For_{q})$-module $\C[\Sp(2n,\For_{q})/\BSp(2n,\For_{q})]$, voir le chapitre \ref{arikikoike}. On rappelle qu'un \textbf{sous-groupe de Borel} de $G$ est un sous-groupe algébrique $B$ qui est connexe résoluble maximal ; tout tore maximal $T \subset G$ est inclus dans un tel sous-groupe, et si $N$ est le normalisateur de $T$ dans $G$, alors $B \cap N=T$ et $W=N/T$ est un \textbf{groupe de Coxeter}, c'est-à-dire un groupe engendré par des involutions et qui admet une présentation du type
$$W=\big\langle s \in S \,\,\big|\,\,\forall s \in S,\,\,s^{2}=1 \,\,\,;\,\,\,\forall s \neq t \in S,\,\,\,\exists m_{st} \in \lle 2,+\infty\rre,\,\,\,(st)^{m_{st}}=(ts)^{m_{st}}=1\big\rangle\,,$$
voir \cite[chapitre 1]{GP00}. On dit que $W$ est le \textbf{groupe de Weyl} de $G$ ; il ne dépend pas du choix d'un tore $T$ et d'un sous-groupe de Borel $B$, car toutes ces paires sont conjuguées dans $G$. L'application de Frobenius $F$ agit sur $W$, et $W^{F}=N^{F}/T^{F}$. On parle de groupe de Chevalley \textbf{non tordu} si $F$ agit trivialement sur $W$ (ou, car c'est équivalent, trivialement sur le diagramme de Dynkin de $G$) ; alors, $W=N^{F}/T^{F}$. Dans ce qui suit, nous supposerons toujours que $G$ est un groupe de Chevalley non tordu. Dans ce contexte, le groupe $G$ est muni d'une \textbf{paire BN}, ce qui se traduit en particulier par l'existence d'une \textbf{décomposition de Bruhat}, valable sur $\overline{k}$ et sur $k$ :
$$G=\bigsqcup_{w \in W} BwB \qquad;\qquad G^{F}=\bigsqcup_{w \in W} B^{F}wB^{F}\,.$$
On renvoie à \cite[\S8.4]{GP00} et à \cite{Gar97} pour des précisions sur cette terminologie, et à \cite{DM91} et \cite[\S1]{Car92} pour des précisions sur la théorie des groupes finis réductifs de type Lie.\bigskip
\bigskip

Suivant \cite{Iwa64}, nous appellerons \textbf{algèbre d'Iwahori-Hecke} d'un groupe de Chevalley $G^{F}$ non tordu sur un corps fini $\For_{q}$ la sous-algèbre de $\C[G^{F}]$ constituée des fonctions qui sont bi-$B^{F}$-invariantes, c'est-à-dire que 
$$\IH(G^{F},B^{F})=\C[B^{F}\backslash G^{F} / B^{F}]\,,$$
le produit des éléments de cette algèbre étant bien sûr le produit de convolution des fonctions. Compte tenu de la décomposition de Bruhat, $\IH(G^{F},B^{F})$ est engendrée linéairement par les classes $B^{F}wB^{F}$, et on peut montrer qu'elle est engendrée en tant qu'algèbre par les $B^{F}sB^{F}$, où $s \in S$ parcourt les générateurs du groupe de Weyl $W$. Plus précisément, si 
$$w=s_{i_{1}}s_{i_{2}}\cdots s_{i_{r}}$$
est une \textbf{expression réduite}\footnote{On dit dans ce cas que $r$ est la \textbf{longueur} de l'élément $w$, ce qu'on note $r=\ell(w)$.} de $w$, c'est-à-dire une écriture de $w$ comme produits de générateurs $s\in S$ de longueur minimale, alors
$$B^{F}wB^{F} = (B^{F}s_{i_{1}}B^{F})\,(B^{F}s_{i_{2}}B^{F})\cdots (B^{F}s_{i_{r}}B^{F})\,,$$
voir le théorème 3.2 de \cite{Iwa64}.\bigskip

\begin{example}
Si $G^{F}=\GL(n,\For_{q})$, on peut prendre pour sous-groupe de Borel le groupe $B^{F}=\GB(n,\For_{q})$ constitué des matrices triangulaires supérieures, et pour tore maximal le sous-groupe diagonal $T^{F}=(\For_{q})^{n}$. Le normalisateur $N^{F}$ est le sous-groupe des matrices monomiales, \emph{i.e.}, celles qui ont exactement un coefficient non nul sur chaque ligne et sur chaque colonne. Le groupe de Weyl $W=N^{F}/T^{F}$ s'identifie donc au groupe symétrique $\sym_{n}$, et vu comme groupe de Coxeter, ce dernier admet pour présentation
$$\sym_{n}=\left\langle s_{1},\ldots,s_{n-1} \,\,\,\bigg|\,\,\,\substack{\forall i,\,\,(s_{i})^{2}=1\qquad\qquad\quad\,\,\,\,\,\,\\
\forall i\leq n-2,\,\,s_{i}s_{i+1}s_{i}=s_{i+1}s_{i}s_{i+1}\\
\forall i,j,\,\,|i-j|\geq 2\, \Rightarrow \,s_{i}s_{j}=s_{j}s_{i}\quad\,\, }\right\rangle,$$
où les $s_{i}$ sont les transpositions élémentaires $(i,i+1)$. Les relations autres que $(s_{i})^{2}=1$ sont appelées \textbf{relations de tresse} ; elles correspondent aux opérations sur les tresses représentées sur la figure \ref{braidrelation}.
\end{example}

\figcapt{\psset{unit=0.75mm}{\footnotesize\pspicture(0,-75)(190,0)
\pscurve(25,0)(25,-20)(24.5,-30)(15.5,-40)(14.5,-50)(5.5,-60)(5,-70)
\pscurve[border=1.5mm,bordercolor=white](15,0)(14.5,-10)(5.5,-20)(5,-35)(5.5,-50)(14.5,-60)(15,-70)
\pscurve[border=1.5mm,bordercolor=white](5,0)(5.5,-10)(14.5,-20)(15.5,-30)(24.5,-40)(25,-50)(25,-70)
\rput(5,-75){$i$}
\rput(14,-75){$i+1$}
\rput(26,-75){$i+2$}
\rput(45,-75){$i$}
\rput(54,-75){$i+1$}
\rput(66,-75){$i+2$}
\rput(35,-35){$=$}
\pscurve(65,0)(64.5,-10)(55.5,-20)(54.5,-30)(45.5,-40)(45,-50)(45,-70)
\pscurve[border=1.5mm,bordercolor=white](55,0)(55.5,-10)(64.5,-20)(65,-35)(64.5,-50)(55.5,-60)(55,-70)
\pscurve[border=1.5mm,bordercolor=white](45,0)(45,-20)(45.5,-30)(54.5,-40)(55.5,-50)(64.5,-60)(65,-70)
\pscurve(105,-10)(104.5,-20)(95.5,-30)(95,-40)(95,-60)
\pscurve[border=1.5mm,bordercolor=white](95,-10)(95.5,-20)(104.5,-30)(105,-40)(105,-60)
\rput(112.5,-35){$\cdots$}
\rput(95,-65){$i$}
\rput(105,-65){$i+1$}
\rput(120,-65){$j$}
\rput(130,-65){$j+1$}
\rput(150,-65){$i$}
\rput(160,-65){$i+1$}
\rput(175,-65){$j$}
\rput(185,-65){$j+1$}
\pscurve(130,-10)(130,-30)(129.5,-40)(120.5,-50)(120,-60)
\pscurve[border=1.5mm,bordercolor=white](120,-10)(120,-30)(120.5,-40)(129.5,-50)(130,-60)
\rput(140,-35){$=$}
\pscurve(160,-10)(160,-30)(159.5,-40)(150.5,-50)(150,-60)
\pscurve[border=1.5mm,bordercolor=white](150,-10)(150,-30)(150.5,-40)(159.5,-50)(160,-60)
\rput(167.5,-35){$\cdots$}
\pscurve(185,-10)(184.5,-20)(175.5,-30)(175,-40)(175,-60)
\pscurve[border=1.5mm,bordercolor=white](175,-10)(175.5,-20)(184.5,-30)(185,-40)(185,-60)
\endpspicture}}{Relations de tresse entre les générateurs $s_{i}$ du groupe symétrique $\sym_{n}$.\label{braidrelation}}{Relations de tresse entre les générateurs $s_{i}$ du groupe symétrique $\sym_{n}$}\bigskip
\bigskip

On doit à N. Iwahori une présentation explicite de l'algèbre $\IH(G^{F},B^{F})$ en termes des générateurs $T_{s}=B^{F}sB^{F}$, où $s$ parcourt un ensemble d'involutions $S$ engendrant le groupe de Weyl $W$, \emph{cf.} \cite{Iwa64}. Ainsi, si $W=\langle S \rangle$ admet pour présentation les relations quadratiques $s^{2}=1$ et les relations de tresse
$$\forall s,t,\,\,\,\underbrace{ststst\cdots}_{m_{st} \text{ termes}}=\underbrace{tststs\cdots}_{m_{st} \text{ termes}}\,,$$
alors les relations entre les $T_{s}$ sont données par le théorème suivant :
\begin{theorem}[Iwahori-Hecke, \cite{Iwa64}]\label{iwahorihecke}
L'algèbre d'Iwahori-Hecke $\IH(G^{F},B^{F})$ admet pour présentation :
\begin{align*}&\forall s,\,\,(T_{s})^{2}=(q-1)T_{s} +q \\ 
&\forall s,t,\,\,\underbrace{T_{s}T_{t}T_{s}T_{t}\cdots}_{m_{st} \text{ termes}}=\underbrace{T_{t}T_{s}T_{t}T_{s}\cdots}_{m_{st}\text{ termes}} 
\end{align*}
Autrement dit, $\IH(G^{F},B^{F})$ a les mêmes relations de tresse que l'algèbre de groupe $\C W$, et les relations quadratiques $(s-1)(s+1)=0$ sont déformées en $(T_{s}-q)(T_{s}+1)=0$.
\end{theorem}\bigskip\bigskip

Plus généralement, un groupe de Coxeter $(W,S)$ étant fixé, nous appellerons \textbf{algèbre d'Hecke générique} de $W$ la $\C(q)$-algèbre engendrée par des symboles $T_{s}$ vérifiant les relations de tresse du groupe et les relations quadratiques $(T_{s}-q)(T_{s}+1)=0$. Cette $\C(q)$-algèbre sera notée $\IH(W)$ ; elle admet pour $\C(q)$-base linéaire les 
$$T_{w}=T_{s_{i_{1}}}T_{s_{i_{2}}}\cdots T_{s_{i_{r}}}\,,$$
où $w$ parcourt $W$ et $w=s_{i_{1}}s_{i_{2}}\cdots s_{i_{r}}$ est une expression réduite --- le choix d'une expression réduite ne change pas $T_{w}$ en vertu du théorème de Matsumoto, voir \cite[\S1.2]{GP00}. Pour toute spécialisation du paramètre $q$ différente de $q=0$, on peut montrer que la $\C$-algèbre spécialisée $\IH_{q}(W)$ est semi-simple si et seulement si le \textbf{polynôme de Poincaré}
$$W(q)=\sum_{w \in W} q^{\ell(w)}$$ ne s'annule pas, voir \cite[chapitre 7]{GP00} et \cite{Gyo95} --- c'est le \textbf{théorème de Gyoja-Uno}. L'algèbre générique $\IH(W)$ est également semi-simple, et dans le cas semi-simple, l'algèbre d'Hecke a les mêmes représentations irréductibles que l'algèbre de groupe $\C W$, et avec les mêmes multiplicités, voir \cite[théorèmes 7.4.6 et 8.1.7]{GP00}.\bigskip

\begin{example}
L'algèbre d'Hecke de type A, c'est-à-dire l'algèbre d'Hecke du groupe symétrique, est semi-simple si et seulement si
$$W(q)=\sum_{\sigma \in \sym_{n}} q^{\ell(\sigma)} \neq 0\,.$$
Or, dans l'algèbre du groupe symétrique, la somme de toutes les permutations admet pour factorisation :
$$\sum_{\sigma \in \sym_{n}} \sigma = (1+s_{n-1}+s_{n-2}s_{n-1}+\cdots+s_{1}\cdots s_{n-2}s_{n-1})\,\cdots\,(1+s_{2}+s_{1}s_{2})\,(1+s_{1})$$
En effet, les termes du premier facteur du terme à droite sont les cycles $(n)$, $(n,n-1),\ldots$ et toute permutation de $\sym_{n}$ s'écrit de manière unique 
$$\sigma=(n,n-1,\ldots,\sigma^{-1}(n))\circ \tau\,$$
où $\tau$ est une permutation de $\sym_{n-1}$. L'identité est donc vraie par récurrence sur $n$, et de plus, cette factorisation des permutations ne donne que des expressions réduites. En spécialisant l'algèbre du groupe symétrique par la règle $s_{i} \mapsto q$, on obtient donc :
$$W(q)=\sum_{\sigma \in \sym_{n}} q^{\ell(\sigma)}=(1+q+\cdots+q^{n-1})\,\cdots\,(1+q+q^{2})\,(1+q) = \{n!\}_{q}$$
L'algèbre d'Hecke du groupe symétrique est donc semi-simple si et seulement si $q \neq 0$ n'est pas une racine de l'unité d'ordre $2 \leq i \leq n$. Et dans ce cas, les modules irréductibles $V^{\lambda}(q)$ sur $\IH_{q}(\sym_{n})$ sont indexés par les partitions $\lambda \in \ym_{n}$, et ils ont pour dimensions et multiplicités les dimensions usuelles $\dim \lambda$. En particulier, pour $q$ réel strictement positif, $\IH_{q}(\sym_{n})$ est toujours une algèbre semi-simple qui a la même théorie des représentations que $\C\sym_{n}$.
\end{example}\bigskip
\bigskip

Expliquons maintenant le lien entre ces algèbres d'Hecke et la mesure de Plancherel du module des fonctions sur la variété de drapeaux $G^{F}/B^{F}$, où $G^{F}$ est toujours un groupe de Chevalley non tordu sur $\For_{q}$. On note $V=\C[G^{F}/B^{F}]$, et on considère une application linéaire $u \in \hendo(V)$ qui commute avec l'action de $G^{F}$. Comme $G^{F}$ agit transitivement sur la variété de drapeaux, $u$ est entièrement déterminée par la valeur de $u(B^{F})$, et nous noterons cette valeur
$$u\left(\sum_{b \in B^{F}} b\right)=\sum_{g \in G^{F}} h_{g}\,g\,.$$
L'image $u(B^{F})$ est dans $\C[G^{F}/B^{F}]$, donc $h$ est $B^{F}$-invariante à droite. De plus, $bB^{F}=B^{F}$ pour tout $b$ dans $B^{F}$, donc $h$ est aussi $B^{F}$-invariante à gauche. Par conséquent, le commutant de l'action de $\C G^{F}$ sur la variété de drapeaux est
$$\hendo_{G^{F}}\left(\C[G^{F}/B^{F}]\right)=\C[B^{F}\backslash G^{F} / B^{F}]=\IH(G^{F},B^{F})=\IH_{q}(W)\,,$$
c'est-à-dire que c'est l'algèbre d'Hecke du groupe de Weyl $W$. On peut donc considérer $\C[G^{F}/B^{F}]$ comme un $(\C G^{F},\IH_{q}(W))$-bimodule, et la théorie du bicommutant de Schur-Weyl (voir le paragraphe \ref{shoji}) s'applique. Ainsi, en tant que bimodule, on a la décomposition suivante en produits tensoriels d'irréductibles :
$$_{G^{F} \curvearrowright }\big\{\C[G^{F}/B^{F}]\big\}_{\curvearrowleft\, \IH_{q}(W)} =\sum_{\lambda \in \widehat{W}} U^\lambda(q) \otimes_{\C} V^\lambda(q)\,.$$
où les $U^{\lambda}(q)$ sont les classes d'isomorphisme de modules intervenant dans la décomposition de $\C[G^{F}/B^{F}]$ en $G^{F}$-modules irréductibles (les \textbf{modules unipotents}), et où les $V^{\lambda}(q)$ sont les classes d'isomorphisme de modules irréductibles sur $\IH_{q}(W)$. En particulier, pour $G=\GL(n)$ :
$$_{\GL(n,\For_q) \curvearrowright }\big\{\C[\GL(n,\For_q)/\GB(n,\For_q)]\big\}_{\curvearrowleft \,\IH_{q}(\sym_{n})} =\sum_{\lambda \in \ym_n} U^\lambda(q) \otimes_{\C} V^\lambda(q)\,.$$\bigskip

Soit $\tau$ la restriction de la trace normalisée de $\hendo_{\C}(\C[G^{F}/B^{F}])$ à l'algèbre d'Iwa\-hori-Hecke $\hendo_{G^{F}}(\C[G^{F}/B^{F}])=\IH_{q}(W)$. Sur la base $(T_{w})_{w \in W}$ de l'algèbre d'Hecke, cette trace vaut 
$$\tau(T_{w})=\mathbb{1}_{(w=e_{W})}\,,$$ 
et la trace d'un produit $\tau(T_{w}T_{w}')$ vaut $0$ si $ww' \neq e_{W}$, et $q^{\ell(w)}$ si $ww'=e_{W}$. On peut donc considérer cette trace quelque soit la valeur du paramètre $q$, même si ce n'est pas la puissance d'un nombre premier ; et l'algèbre d'Hecke (complexe) $\IH_{q}(W)$ a donc une structure naturelle d'espace de probabilité non commutatif\footnote{Notons que $\tau$ est une quantification de la trace usuelle de $\C W$ associée à la représentation régulière gauche. Si $W$ est le groupe symétrique, il existe en réalité deux quantifications possibles de $\tau$, voir \cite[\S6]{RR97} ; mais celle que nous présentons est la seule intéressante dans le contexte des groupes linéaires finis $\GL(n,\For_{q})$.}. Pour cette structure, la base duale de la base <<~canonique~>> $(T_{w})_{w \in W}$ est $(q^{\ell(w)}T_{w^{-1}})_{w \in W}$. Ceci étant, étant donnée une mesure de probabilité $P$ sur $\widehat{W}$, on peut comme dans le chapitre \ref{plancherel} lui associer une structure d'espace de probabilité non commutatif en considérant la trace normalisée
$$\tau_{P}=\bigoplus_{\lambda \in \widehat{W}} P(\lambda)\,\,\tr_{\!V^{\lambda}(q)}$$
sur l'algèbre $\bigoplus_{\lambda \in \widehat{W}} \hendo(V^{\lambda}(q))$. Alors, la transformée de Fourier abstraite $$\IH_{q}(W) \to \bigoplus_{\lambda \in \widehat{W}} \hendo(V^{\lambda}(q))$$ est une isométrie d'espaces de probabilité non commutatifs si et seulement si l'on prend pour $P$ la $q$-mesure de Plancherel
$$P(\lambda)=\frac{(\dim U^{\lambda}(q)) \times (\dim V^{\lambda})}{\card G^{F}/B^{F}}\,\,;$$
en particulier, $P=M_{n,q}$ si $G^{F}=\GL(n,\For_{q})$. Autrement dit :
\begin{proposition}[Décomposition de la trace de l'algèbre d'Hecke dans la base des $q$-caractères]\label{lastround}
Dans la base $(\chi^{\lambda}(q))_{\lambda \in \widehat{W}}$ des caractères normalisés des modules irréductibles $(V^{\lambda}(q))_{\lambda \in \widehat{W}}$,  la trace canonique $\tau$ de l'algèbre d'Hecke $\IH_{q}(W)$ admet pour décomposition 
$$\tau=\sum_{\lambda \in \widehat{W}} \frac{D_{\lambda}(q)\times\dim\lambda}{W(q)}\,\,\chi^{\lambda}(q)\,.$$
Cette décomposition reste valable pour tout paramètre $q \neq 0$ tel que $W(q) \neq 0$.
\end{proposition} 
\noindent Compte tenu de cette proposition, il est naturel de prendre pour analogues des observables $\varSigma_{\mu}$ des versions renormalisées $\varSigma_{\mu,q}$ des $q$-caractères $\chi^{\lambda}(q)$ de l'algèbre d'Hecke $\IH_{q}(\sym_{n})$. Le prochain paragraphe est consacré à l'étude approfondie de ces $q$-caractères.

\section{Caractères des algèbres d'Hecke de type A}\label{ram}
Dans cette section, $q$ est un paramètre complexe non nul et qui n'est pas une racine non triviale de l'unité, de sorte que $\IH_{q}(\sym_{n})$ est une algèbre semi-simple abstraitement isomorphe à l'algèbre de groupe $\C\sym_{n}$. Si $\mu=(\mu_{1},\ldots,\mu_{r})$ est une partition de taille $n$, on note $\sigma_{\mu}$ la permutation standard de type $\mu$, c'est-à-dire que
\begin{align*}\sigma_{\mu}&=(1,2,\ldots,\mu_{1})\,(\mu_{1}+1,\ldots,\mu_{1}+\mu_{2})\,\cdots\,(\mu_{1}+\cdots+\mu_{r-1}+1,\cdots,n)\,\\
&=(s_{1}s_{2}\cdots s_{\mu_{1}-1})\,(s_{\mu_{1}+1}\cdots s_{\mu_{1}+\mu_{2}-1})\,\cdots\,(s_{\mu_{1}+\cdots+\mu_{r-1}+1}\cdots s_{n-1})\,.\end{align*}
C'est une permutation de longueur minimale dans la classe de conjugaison $C_{\mu}$ ; on note $T_{\mu}$ l'élément correspondant dans l'algèbre d'Hecke $\IH_{q}(\sym_{n})$. La \textbf{table des caractères} de l'algèbre d'Hecke est la matrice $(\chi^{\lambda}(q,T_{\mu})=\chi^{\lambda}(q,\mu))_{\lambda,\mu \in \ym_{n}}$. \emph{A priori}, il n'est pas clair que cette table des caractères détermine entièrement les caractères de $\IH_{q}(\sym_{n})$. En effet, si deux permutations $\sigma$ et $\tau$ sont conjuguées dans $\sym_{n}$, il n'est pas vrai en général que les éléments $T_{\sigma}$ et $T_{\tau}$ soient conjugués dans l'algèbre $\IH_{q}(\sym_{n})$, et étant donnée une fonction $\phi$ telle que
$$\forall a,b \in \IH_{q}(\sym_{n}),\,\,\,\phi(ab)=\phi(ba)\,,$$ on peut avoir $\phi(T_{\sigma})\neq \phi(T_{\tau})$. Ainsi, $T_{\sigma}$ et $T_{\tau}$ ne sont pas forcément congrus modulo le sous-espace $[\IH_{q}(\sym_{n}),\IH_{q}(\sym_{n})]$.\bigskip

\begin{example}
Pour tout entier $n\geq 2$, l'algèbre d'Hecke $\IH_{q}(\sym_{n})$ a exactement deux caractères irréductibles de dimension $1$ : le \textbf{caractère signature} $T_{\sigma} \mapsto (-1)^{\ell(\sigma)}$, et le \textbf{caractère d'indice} $T_{\sigma} \mapsto q^{\ell(\sigma)}$. Pour $n=3$, les transpositions $(1,2)$ et $(1,3)$ sont évidemment conjuguées dans $\sym_{3}$, mais la première a pour longueur $1$, et la seconde a pour longueur $3$. Par conséquent, si $\chi$ est le caractère d'indice, alors
$$\chi(T_{(1,2)})=q \neq q^{3}=\chi(T_{(1,3)})\,.$$
\end{example}\bigskip

Néanmoins, si $\sigma$ et $\tau$ sont des permutations \textbf{fortement conjuguées}, c'est-à-dire qu'il existe une suite de permutations $\sigma=\sigma_{0},\sigma_{1},\ldots,\sigma_{r}=\tau$ telle que
$$\forall i,\,\,\exists x \in \sym_{n},\,\,\,\begin{cases}&x\sigma_{i}=\sigma_{i+1}x \,\text{ et }\, \ell(x\sigma_{i})=\ell(x)+\ell(\sigma_{i})\\
\text{ou}&\sigma_{i}x=x\sigma_{i+1} \,\text{ et }\,\ell(\sigma_{i}x)=\ell(\sigma_{i})+\ell(x)\,,\end{cases}
$$
alors $T_{\sigma}$ et $T_{\tau}$ sont conjugués dans $\IH_{q}(\sym_{n})$, et ils sont donc congrus modulo $[\IH_{q}(\sym_{n}),\IH_{q}(\sym_{n})]$ dans $\IH_{q}(\sym_{n})$. Par suite, ils ont la même image par tout caractère $\chi^{\lambda}(q)$. En particulier, si $\sigma$ est un élément de longueur minimale dans une classe de conjugaison $C_{\mu}$, alors on peut montrer qu'il est fortement conjugué à $\sigma_{\mu}$ (voir \cite[théorème 3.2.9]{GP00}), et on a donc 
$$\chi^{\lambda}(q,T_{\sigma})=\chi^{\lambda}(q,\mu)$$ pour tout caractère irréductible. Maintenant, si $\sigma$ n'est pas de longueur minimale dans sa classe de conjugaison, alors il existe un algorithme qui donne une $\Z[q]$-combinaison linéaire
$$\sum_{\substack{\tau \text{ minimal dans sa} \\
\text{classe de conjugaison}}} c_{\sigma\tau}\,T_{\tau}\label{reductioncharacter}$$
congrue à $T_{\sigma}$ modulo $[\IH_{q}(\sym_{n}),\IH_{q}(\sym_{n})]$, voir \cite[\S8.2]{GP00} et \cite[théorème 5.1]{Ram91}. Les valeurs $\chi^{\lambda}(q,\mu)$ déterminent donc bien les caractères de l'algèbre d'Iwahori-Hecke $\IH_{q}(\sym_{n})$. Dans ce qui suit, nous utiliserons également les notations $\varsigma^{\lambda}(q,\mu)$, qui correspondent aux caractères non normalisés ; ainsi, $\varsigma^{\lambda}(q,\mu)=(\dim\lambda)\,\,\chi^{\lambda}(q,\mu)$. Notons que tous ces résultats restent valables pour l'algèbre d'Hecke de n'importe quel groupe de Coxeter fini : la table des caractères de l'algèbre est toujours doublement indexée par les classes de conjugaison du groupe $W$.\bigskip
\bigskip

Ceci étant, le calcul explicite de la table des caractères $(\varsigma^{\lambda}(q,\mu))_{\lambda,\mu \in \ym_{n}}$ de l'algèbre d'Hecke du groupe symétrique peut être mené comme dans la section \ref{frobeniusschur}, c'est-à-dire à l'aide d'une formule de Frobenius dans l'algèbre $\Lambda$ des fonctions symétriques. Ainsi, dans \cite{Ram91}, A. Ram utilise la dualité de Schur-Weyl\footnote{Nous donnerons un bref aper\c cu de cette théorie dans la section \ref{shoji}.} sur $(\C^{N})^{\otimes n}$ entre le groupe quantique $U_{q}(\mathfrak{gl}(N,\C))$ et l'algèbre d'Hecke $\IH_{q}(\sym_{n})$ pour calculer les caractères $\varsigma^{\lambda}(q,\mu)$, et il les relie aux polynômes de Hall-Littlewood modifiés (voir \cite{RRW96} pour les propriétés de ces polynômes) :
$$q_{k}(X,q)=\frac{h_{k}(X(q-1))}{q-1} \,.$$
Ici, $X$ est un alphabet et $q$ est une variable ; $X(q-1)$ est donc la différence d'alphabets $qX-X$. Comme $h_{k}=s_{k}$, par la formule de Frobenius,
$$h_{k}(X)=\sum_{\mu \in \ym_{k}} (z_{\mu})^{-1}\,p_{\mu}(X)\qquad;\qquad q_{k}(X,q)=\frac{1}{q-1}\,\sum_{\mu \in \ym_{k}} (q^{\mu}-1)\,(z_{\mu})^{-1}\,p_{\mu}(X)$$
où $q^{\mu}-1=\prod_{i=1}^{\ell(\mu)}q^{\mu_{i}}-1$. Si $\mu=(\mu_{1},\ldots,\mu_{r})$ est une partition, notons $q_{\mu}(X,q)=\prod_{i=1}^{r} q_{\mu_{i}}(X,q)$. D'après ce qui précède, $q_{\mu}(X,1)=p_{\mu}(X)$ pour toute partition $\mu$.
\begin{theorem}[$q$-formule de Frobenius, \cite{Ram91}]\label{qfrobeniusschur}
La valeur du $q$-caractère non normalisé $\varsigma^{\lambda}(q,\mu)$ est $\varsigma^{\lambda}(q,\mu)=\scal{s_{\lambda}}{q_{\mu}(q)}$. Ainsi, pour toute partition $\mu\in \ym$,
$$q_{\mu}(X,q)=\sum_{\lambda \in \ym_{n}} \varsigma^{\lambda}(q,\mu)\,s_{\lambda}(X)\,.$$
\end{theorem}
\figcapt{
{\normalsize
\begin{tabular}{|c|c|}
\hline $\IH_{q}(\sym_{1})$ & $(1)$ \\ 
\hline $(1)$ & $1$ \\
\hline 
\end{tabular}
\qquad 
\begin{tabular}{|c|c|c|}
\hline $\IH_q(\sym_{2})$ & $(2)$ & $(1,1)$ \\ 
\hline $(2)$ & $q$ & $1$ \\ 
\hline $(1,1)$ & $-1$ & $1$\\
\hline
\end{tabular}
\qquad 
\begin{tabular}{|c|c|c|c|}
\hline $\IH_q(\sym_{3})$ & $(3)$ & $(2,1)$ & $(1,1,1)$ \\ 
\hline $(3)$ & $q^2$ & $q$ & $1$ \\ 
\hline $(2,1)$ & $-q$ & $q-1$ & $2$ \\ 
\hline $(1,1,1)$ & $1$ & $-1$ & $1$\\
\hline
\end{tabular} \bigskip \bigskip

\begin{tabular}{|c|c|c|c|c|c|}
\hline $\IH_q(\sym_{4})$ & $(4)$ & $(3,1)$ & $(2,2)$ & $(2,1,1)$ & $(1,1,1,1)$ \\ 
\hline $(4)$ & $q^3$ & $q^2$ & $q^2$ & $q$ & $1$ \\ 
\hline $(3,1)$ & $-q^2$ & $q^2-q$ & $q^2-2q$ & $2q-1$ & $3$ \\ 
\hline $(2,2)$ & $0$ & $-q$ & $q^2+1$ & $q-1$ & $2$ \\ 
\hline $(2,1,1)$ & $q$ & $-q+1$ & $-2q+1$ & $q-2$ & $3$ \\ 
\hline $(1,1,1,1)$ & $-1$ & $1$ & $1$ & $-1$ & $1$ \\
\hline 
\end{tabular}\bigskip
}
}{Tables des caractères des algèbres d'Hecke $\IH_{q}(\sym_{n})$ pour $n \leq 4$.\label{tableofcharacterhecke}}{Tables des caractères des algèbres d'Hecke $\IH_{q}(\sym_{n})$ pour $n \leq 4$}

Cette formule a de nombreuses conséquences ; en particulier, elle implique une $q$-règle de Murnaghan-Nakayama qui généralise celle donnée par le théorème \ref{murnaghannakayama}, et qui permet un calcul récursif de tous les $q$-caractères, voir la figure \ref{tableofcharacterhecke}. Ainsi, si $\mu=(\mu_{1},\ldots,\mu_{m})$, alors $\varsigma^\lambda(q,\mu)$ est égal à
$$\sum_{\emptyset=\lambda_0 \subset \lambda_1 \subset \cdots \subset \lambda_m=\lambda}\,\, \prod_{j=1}^m\left( (q-1)^{\mathrm{cc}(\lambda_j\setminus \lambda_{j-1})-1}\prod_{\nu_{j,k} \subset \lambda_j\setminus\lambda_{j-1}} \!\!\!\!\!\!(-1)^{\mathrm{r}(\nu_{j,k}) -1}\,q^{\mathrm{c}(\nu_{j,k}) -1} \right),$$
la somme étant effectuée sur les suites de partitions telles que chaque $\lambda_{j} \setminus \lambda_{j-1}$ soit un ruban de poids $\mu_{j}$ avec $\mathrm{cc}(\lambda_{j} \setminus \lambda_{j-1})$ composantes connexes $\nu_{j,k}$, chaque composante connexe ayant $\mathrm{r}(\nu_{j,k})$ lignes et $\mathrm{c}(\nu_{j,k})$ colonnes --- par rapport à la formule usuelle, notons qu'on autorise ici des rubans éventuellement non connexes. On retrouve bien le théorème \ref{murnaghannakayama} en spécialisant la formule en $q=1$. On renvoie à \cite[\S2]{RR97} pour une preuve de cette $q$-formule ; l'article \cite{HR96} traite également le cas des algèbres d'Hecke des groupes de Coxeter de type B et D, mais nous y reviendrons dans le chapitre \ref{arikikoike}.\bigskip
\bigskip

Exposons finalement la formule de Ram (\emph{cf.} \cite[théorème 5.4]{Ram91}), qui relie les $q$-cara\-ctè\-res de l'algèbre d'Hecke $\IH_{q}(\sym_{n})$ aux caractères usuels du groupe symétrique $\sym_{n}$. D'après la $q$-formule de Frobenius, $\varsigma^{\lambda}(q,\mu)$ est égal à 
\begin{align*}
\scal{q_{\mu}(X,q)}{s_{\lambda}(X)}&=\frac{1}{(q-1)^{\ell(\mu)}}\scal{h_{\mu}(X(q-1))}{s_{\lambda}(X)}\\
&=\frac{1}{(q-1)^{\ell(\mu)}}\sum_{\nu \in \ym_{n}} \frac{\scal{h_{\mu}}{p_{\nu}}\,\scal{p_{\nu}(X(q-1))}{s_{\lambda}(X)}}{\scal{p_{\nu}}{p_{\nu}}} \\
&=\frac{1}{(q-1)^{\ell(\mu)}}\sum_{\nu \in \ym_{n}} \frac{q^{\nu}-1}{z_{v}}\,\scal{h_{\mu}}{p_{\nu}}\,\scal{p_{\nu}}{s_{\lambda}}
\end{align*}
et dans la dernière expression, les produits scalaires $\scal{p_{\nu}}{s_{\lambda}}$ sont les caractères usuels $\varsigma^{\lambda}(\nu)$. Ainsi :

\begin{proposition}[Formule de Ram, \cite{Ram91}]\label{ramformula}
Pour toute partition $\lambda$, les valeurs du $q$-caractère irréductible $\varsigma^{\lambda}(q)$ s'expriment en fonction de celles du caractère non déformé $\varsigma^{\lambda}$ : 
$$\forall \mu,\,\,\,(q-1)^{\ell(\mu)} \,\varsigma^{\lambda}(q,\mu)=\sum_{\nu \in \ym_{n}} \frac{\scal{p_{\nu}}{h_{\mu}}}{\scal{p_{\nu}}{p_{\nu}}}\,(q^{\nu}-1)\,\varsigma^{\lambda}(\nu)\,.$$
\end{proposition}

\noindent Si $\ym_{n}$ est muni de l'ordre lexicographique inverse, alors cette <<~formule de changement de base~>> entre $q$-caractères de $\IH_{q}(\sym_{n})$ et caractères de $\sym_{n}$ est triangulaire : en effet, $\scal{h_{\mu}}{p_{\nu}}=0$ si $\nu \prec \mu$. On peut donc inverser la formule, et exprimer les caractères en fonction des $q$-caractères (!). Ainsi, comme $(m_{\mu})_{\mu \in \ym}$ est la base duale de $(h_{\mu})_{\mu \in \ym}$ pour le produit scalaire canonique sur l'algèbre $\Lambda$, et comme $(p_{\mu})_{\mu \in \ym}$ est une base orthogonale, on peut écrire :
\begin{align*}
(q^{\mu}-1)\,\varsigma^{\lambda}(\mu)&=\sum_{\rho} \frac{\scal{p_{\mu}}{p_{\rho}}}{\scal{p_{\rho}}{p_{\rho}}}\,(q^{\rho}-1)\,\varsigma^{\lambda}(\rho)\\
&=\sum_{\rho,\nu}\scal{p_{\mu}}{m_{\nu}}\,\frac{\scal{h_{\nu}}{p_{\rho}}}{\scal{p_{\rho}}{p_{\rho}}}\,(q^{\rho}-1)\,\varsigma^{\lambda}(\rho)\\
&=\sum_{\nu \in \ym_{n}}\scal{p_{\mu}}{m_{\nu}}\,(q-1)^{\ell(\nu)}\,\varsigma^{\lambda}(q,\nu)\,.
\end{align*}
Cette réécriture de la formule de Ram est quelque peu étrange, puisqu'on exprime maintenant les caractères du groupe symétrique en fonction d'objets beaucoup plus complexes, à savoir les caractères de l'algèbre d'Iwahori-Hecke. C'est néanmoins cette forme de la proposition \ref{ramformula} qui va permettre le calcul des espérances $M_{n,q}[\varSigma_{\mu}]$ des observables de diagrammes sous les $q$-mesures de Plancherel.

\chapter{Asymptotique de la $q$-mesure de Plancherel}\label{qplancherelmeasure}

Dans tout ce chapitre\footnote{Sauf à la toute fin de la section \ref{qgaussian}, après la démonstration du théorème \ref{secondasymptoticqplancherel}, lorsqu'on évoquera des résultats de Borodin relatifs à la forme de Jordan des matrices triangulaires unipotentes à coefficients dans $\For_{q}$ ; dans ce court paragraphe, $q$ sera bien sûr une puissance d'un nombre premier.}, $q$ est un paramètre réel compris strictement entre $0$ et $1$, et on s'intéresse :\vspace{2mm}
\begin{itemize}
\item[-] à la forme limite des partitions aléatoires $\lambda$ tirées suivant les $q$-mesures de Plancherel $M_{n,q}$ ;\vspace{2mm}
\item[-] et à la distribution des observables de diagrammes vues comme variables aléatoires, en particulier les $q$-caractères de l'algèbre d'Hecke $\IH_{q}(\sym_{n})$.\vspace{2mm}
\end{itemize}
Comme $M_{n,q}(\lambda)=M_{n,q^{-1}}(\lambda')$, nos résultats auront des analogues immédiats pour les partitions tirées suivant les $q$-mesures de Plancherel de paramètre $q>1$. 
\figcapt{\psset{unit=1.5mm}
\pspicture(0,0)(100,10)
\psline(0,0)(101,0)(101,1)(0,1)
\psline(51,2)(0,2)
\psline(28,3)(0,3)
\psline(8,4)(0,4)
\psline(7,5)(0,5)
\psline(3,6)(0,6)
\psline(1,7)(0,7)
\psline(1,8)(0,8)
\multido{\n=52+1}{50}{\psline(\n,0)(\n,1)}
\multido{\n=29+1}{23}{\psline(\n,0)(\n,2)}
\multido{\n=9+1}{20}{\psline(\n,0)(\n,3)}
\multido{\n=8+1}{1}{\psline(\n,0)(\n,4)}
\multido{\n=4+1}{4}{\psline(\n,0)(\n,5)}
\multido{\n=2+1}{2}{\psline(\n,0)(\n,6)}
\psline(0,0)(0,8)
\psline(1,0)(1,8)
\endpspicture}{Diagramme de Young aléatoire tiré suivant la $q$-mesure de Plancherel de paramètres $q=1/2$, $n=200$. Les parts de $\lambda$ sont $(101,51,28,8,7,3,1,1)$.\label{qplanhalf}}{Diagramme aléatoire tiré suivant une $q$-mesure de Plancherel}

La figure \ref{qplanhalf} présente un diagramme de Young de taille $200$ tiré suivant la mesure de Plancherel de paramètre $q=1/2$. Sur cet exemple (obtenu à l'aide du logiciel \texttt{sage}), les premières lignes du diagramme suivent presque une progression géométrique de paramètre $1/2$ : $101 \sim 200/2$, $51 \sim 200/4$, $28 \sim 200/8$, etc. Ce résultat sera confirmé par le théorème \ref{firstasymptoticqplancherel}, et nous étudierons également la déviation des lignes par rapport à la progression géométrique attendue ; cette déviation est asymptotiquement gaussienne, voir les théorèmes \ref{secondasymptoticqplancherel} et \ref{thirdasymptoticqplancherel}. Comme dans le chapitre \ref{plancherel}, l'outil essentiel sera l'algèbre des observables de diagrammes $\obs$, et nous commencerons en présentant une quantification de cette algèbre correspondant à la déformation des algèbres de groupes $\C\sym_{n}$ en algèbres d'Hecke $\IH_{q}(\sym_{n})$ (\S\ref{quantum}). La plupart des résultats de ce chapitre ont été obtenus en collaboration avec V. Féray.

\section{Quantification de l'algèbre des observables}\label{quantum}
Dans cette section, nous utiliserons librement les formules de changement de base entre les bases de fonctions symétriques $(h_{\rho})_{\rho\in \ym}$, $(e_{\rho})_{\rho\in \ym}$ et $(p_{\rho})_{\rho\in \ym}$. Ainsi :
$$h_{n} =\sum_{\rho \in \ym_{n}} \frac{p_{\rho}}{z_{\rho}} \qquad;\qquad e_{n}=\sum_{\rho \in \ym_{n}} \frac{\eps_{\rho}\,p_{\rho}}{z_{\rho}}\,,$$
où $\eps_{\rho}=(-1)^{|\rho|-\ell(\rho)}$. La première relation est simplement la formule de Frobenius pour $s_{n}=h_{n}$, étant entendu que la représentation irréductible de $\sym_{n}$ associée à la partition ligne $(n)$ est la représentation triviale de dimension $1$. La seconde formule de changement de base s'en déduit en utilisant l'antipode de l'algèbre de Hopf $\Lambda$.\bigskip\bigskip

Compte tenu des calculs présentés dans la section \ref{ram}, les $q$-analogues naturels des observables $\varSigma_{\mu}$ sont les
$$\varSigma_{\mu,q}(\lambda)=\begin{cases}n^{\downarrow |\mu|}\,\chi^{\lambda}(q,\mu\sqcup 1^{n-|\mu|})&\text{si }n=|\lambda| \geq |\mu|,\\
0&\text{sinon}.\end{cases}$$
Notons qu'il n'est pas possible d'interpréter $\varSigma_{\mu,q}$ directement comme élément de l'algèbre des permutations partielles, ou même comme élément d'une quantification de cette algèbre (une telle quantification sera présentée dans le chapitre \ref{badbeat}). Si $n=|\lambda|=|\mu|$, alors d'après la formule de Ram \ref{ramformula},
\begin{align*}(q-1)^{\ell(\mu)}\,\varSigma_{\mu,q}(\lambda)&=n!\,(q-1)^{\ell(\mu)}\,\chi^{\lambda}(q,\mu)=n!\sum_{\nu \in \ym_{n}} \frac{\scal{h_{\mu}}{p_{\nu}}}{\scal{p_{\nu}}{p_{\nu}}}\,(q^{\nu}-1)\,\chi^{\lambda}(\nu)\\
&=\sum_{\nu \in \ym_{n}}  \frac{\scal{h_{\mu}}{p_{\nu}}}{\scal{p_{\nu}}{p_{\nu}}}\,(q^{\nu}-1)\,\varSigma_{\nu}(\lambda)\,.
\end{align*}
L'identité reste trivialement vraie si $|\lambda|<|\mu|$. Maintenant, si $k=|\mu|<|\lambda|=n$, alors les produits scalaires $\scal{h_{\mu\sqcup1^{n-k}}}{p_{\nu}}$ s'écrivent :
\begin{align*}\scal{h_{\mu\sqcup1^{n-k}}}{p_{\nu}}&=\scal{h_{\mu}\,(p_{1})^{n-k}}{p_{\nu}}=\sum_{\pi_{1} \in \ym_{\mu_{1}},\ldots,\pi_{m} \in \ym_{\mu_{m}}} \frac{1}{z_{\pi_{1}}\cdots z_{\pi_{m}} } \scal{p_{\pi}(p_{1})^{n-k}}{p_{\nu}}\\
&=\sum_{\pi_{1} \in \ym_{\mu_{1}},\ldots,\pi_{m} \in \ym_{\mu_{m}}}\frac{z_{\pi\sqcup 1^{n-k}}}{z_{\pi_{1}}\cdots z_{\pi_{m}}}\,\mathbb{1}_{(\nu=\pi \sqcup 1^{n-k})}\\
&=\frac{z_{\nu}}{z_{\nu_{k}}}\,\mathbb{1}_{(\nu=\nu_{k} \sqcup 1^{n-k})}\,\sum_{\pi_{1} \in \ym_{\mu_{1}},\ldots,\pi_{m} \in \ym_{\mu_{m}}}\frac{z_{\pi}}{z_{\pi_{1}}\cdots z_{\pi_{m}}}\,\mathbb{1}_{(\nu_{k}=\pi)}=\frac{z_{\nu}}{z_{\nu_{k}}}\,\mathbb{1}_{(\nu=\nu_{k} \sqcup 1^{n-k})}\,\scal{h_{\mu}}{p_{\nu_{k}}}\,.
\end{align*}
Donc, si $k=|\mu|<|\lambda|=n$, alors :
\begin{align*}
(q-1)^{\ell(\mu)}\,\varSigma_{\mu,q}(\lambda)&=n^{\downarrow k}\,(q-1)^{\ell(\mu)}\,\chi^{\lambda}(q,\mu\sqcup 1^{n-k}) \\
&=\frac{n^{\downarrow k}}{(q-1)^{n-k}}\sum_{\nu \in \ym_{n}} \frac{\scal{h_{\mu\sqcup1^{n-k}}}{p_{\nu}}}{z_{\nu}} \,(q^{\nu}-1)\,\chi^{\lambda}(\nu)\\
&=n^{\downarrow k}\sum_{\nu_{k} \in \ym_{k}} \frac{\scal{h_{\mu}}{p_{\nu_{k}}}}{z_{\nu_{k}}} \,(q^{\nu_{k}}-1)\,\chi^{\lambda}(\nu_{k}) = \sum_{\nu \in \ym_{k}} \frac{\scal{h_{\mu}}{p_{\nu}}}{\scal{p_{\nu}}{p_{\nu}}}\,(q^{\nu}-1)\,\varSigma_{\nu}(\lambda)\,.
\end{align*}
On déduit de ces calculs le résultat important suivant :
\begin{theorem}[Quantification de l'algèbre $\obs$, \cite{FM10}]\label{quantobs}
Pour toute partition $\mu$ et tout paramètre $q>0$, la fonction $\varSigma_{\mu,q}$ est une observable de diagrammes. La famille $(\varSigma_{\mu,q})_{\mu \in \ym}$ est une base linéaire de $\obs$, et elle constitue une quantification de la base des caractères centraux. Les formules de changement de base sont :
\begin{align*}
(q-1)^{\ell(\mu)}\,\varSigma_{\mu,q}&=\sum_{\nu \in \ym_{k}} \frac{\scal{h_{\mu}}{p_{\nu}}}{\scal{p_{\nu}}{p_{\nu}}}\,(q^{\nu}-1)\,\varSigma_{\nu}\\
(q^{\mu}-1)\,\varSigma_{\mu}&=\sum_{\nu \in \ym_{k}} \scal{m_{\nu}}{p_{\mu}}\,(q-1)^{\ell(\nu)}\,\varSigma_{\nu,q}
\end{align*}
où $k=|\mu|$. En particulier, $\deg \varSigma_{\mu,q}=|\mu|$.
\end{theorem}\bigskip

La nouvelle base d'observables $(\varSigma_{\mu,q})_{\mu \in \ym}$ est particulièrement adaptée à l'étude des $q$-mesures de Plancherel, puisqu'on peut calculer sans mal $M_{n,q}[\varSigma_{\mu,q}]$. En effet, pour toute partition $\mu$ de taille $k$ inférieure à $n$, on peut écrire compte tenu de la proposition \ref{lastround} :
$$M_{n,q}[\varSigma_{\mu,q}]=n^{\downarrow k}\sum_{\lambda \in \ym_{n}}M_{n,q}(\lambda) \,\,\chi^{\lambda}(q,\mu \sqcup 1^{n-k})=n^{\downarrow k}\,\tau(T_{\mu \sqcup 1^{n-k}})=n^{\downarrow k}\,\mathbb{1}_{(\mu=1^{k})}\,.$$
Le résultat reste vrai si $k>n$, car $n^{\downarrow k}$ est dans ce cas nul. En utilisant le théorème \ref{quantobs}, on peut ensuite calculer les espérances des caractères centraux sous la $q$-mesure de Plancherel :
\begin{proposition}[Espérance des caractères centraux sous la $q$-mesure de Plancherel]\label{melodyofoblivion}
Pour toute partition $\mu$,
$$M_{n,q}[\varSigma_{\mu}]=n^{\downarrow |\mu|} \,\frac{(1-q)^{|\mu|}}{1-q^{\mu}}\,.$$
\end{proposition}
\begin{proof}
On développe $\varSigma_{\mu}$ dans la base des $q$-observables $\varSigma_{\nu,q}$ :
\begin{align*}
M_{n,q}[\varSigma_{\mu}]&=\frac{1}{q^{\mu}-1}\sum_{\nu \in \ym_{k}} \scal{m_{\nu}}{p_{\mu}}\,(q-1)^{\ell(\mu)}\,M_{n,q}[\varSigma_{\mu,q}]\\
&=\frac{n^{\downarrow k}}{q^{\mu}-1}\,\scal{m_{1^{k}}}{p_{\mu}}\,(q-1)^{\ell(1^{k})}=n^{\downarrow |\mu|}\,\eps_{\mu}\,\frac{(q-1)^{|\mu|}}{q^{\mu}-1}
\end{align*}
car $m_{1^{k}}=e_{k}=\sum_{\mu \in \ym_{k}} \frac{\eps_{\mu}p_{\mu}}{z_{\mu}}$. Le signe $\eps_{\mu}$ permet ensuite d'écrire le numérateur et le dénominateur comme produits de termes positifs.
\end{proof}\medskip

\begin{lemma}[Ordre de grandeur des espérances d'observables de diagrammes]\label{stilldre}
L'espérance $\esper[X]$ d'une observable de diagrammes $X$ sous la $q$-mesure de Plancherel $M_{n,q}$ est un $O(n^{\deg X})$, et cette estimation est optimale, au sens suivant : pour toute observable $X$, il existe une fraction rationnelle $X(q)$ dont les pôles sont des racines de l'unité, et telle que
$$\esper[X]=X(q)\,n^{\deg X}+O(n^{\deg X -1})\,.$$
\end{lemma}
\begin{proof}
D'après la proposition \ref{melodyofoblivion}, le résultat est vrai si $X=\varSigma_{\mu}$ est un caractère central, avec $\varSigma_{\mu}(q)=(1-q)^{|\mu|}/(1-q^{\mu})$. Maintenant, si $X$ est une observable quelconque, décomposons $X$ dans la base des caractères centraux :
$$X=\sum_{|\mu| \leq \deg X} c_{\mu}\,\varSigma_{\mu}\quad\text{avec les }c_{\mu} \in \C.$$
En passant à l'espérance et en utilisant l'estimation précédente, on obtient :
$$\esper[X]=\left(\sum_{|\mu|=\deg X} c_{\mu}\frac{(1-q)^{|\mu|}}{1-q^{\mu}}\right)n^{\deg X}+O(n^{\deg X -1})=X(q)\,n^{\deg X}+O(n^{\deg X-1})\,.$$
Notons que la somme entre parenthèses peut être une fraction rationnelle nulle, même si $X$ est une observable non nulle ; en effet, les séries formelles
$$\frac{1}{1-q}\,,\,\,\frac{1}{1-q^{2}}\,,\ldots,\,\,\frac{1}{1-q^{m}}\,,\ldots$$
ne sont pas algébriquement indépendantes sur $\C$.\end{proof}\bigskip

En fait, pour toute famille $(P_{n})_{n \in \N}$ de mesures de probabilité sur les partitions, l'espérance sous $P_{n}$ d'une observable de degré $X$ est un $O(n^{\deg X})$, car ceci est vrai pour les observables $p_{\mu}$ :
\begin{align*}|P_{n}[p_{\mu}]|&\leq\sum_{\lambda \in \ym_{n}} P_{n}(\lambda)\,\left|\prod_{i=1}^{\ell(\mu)} p_{\mu_{i}}(A-(-B))\right| \leq \sum_{\lambda \in \ym_{n}} P_{n}(\lambda)\,\left(\prod_{i=1}^{\ell(\mu)} p_{\mu_{i}}(A+B) \right)\\
& \leq \sum_{\lambda \in \ym_{n}} P_{n}(\lambda)\,\left(\prod_{i=1}^{\ell(\mu)} p_{1}(A+B)^{|\mu_{i}|}\right)=\sum_{\lambda \in \ym_{n}} P_{n}(\lambda)\, n^{|\mu|}=n^{|\mu|}\,.
\end{align*}
Le lemme \ref{stilldre} assure que dans le cas des $q$-mesures de Plancherel, cette estimation est <<~génériquement~>> optimale. Dans l'étude asymptotique des observables de diagrammes, ceci permettra de remplacer une observable par une autre observable de même composante de plus haut degré.
\bigskip

\begin{theorem}[Convergence en probabilité des observables sous les $q$-mesures de Plancherel, \cite{FM10}]\label{convobservable}
Si $\rightarrow_{M_{n,q}}$ désigne la convergence en probabilité sous les $q$-mesures de Plancherel, alors 
$$\frac{p_{\mu}(\lambda)}{n^{|\mu|}}\rightarrow_{M_{n,q}} \frac{(1-q)^{|\mu|}}{1-q^{\mu}}\!\qquad;\qquad\! \frac{\varSigma_{\mu}(\lambda)}{n^{|\mu|}} \rightarrow_{M_{n,q}} \frac{(1-q)^{|\mu|}}{1-q^{\mu}}\!\qquad;\qquad\! \frac{\varSigma_{\mu,q}(\lambda)}{n^{|\mu|}} \rightarrow_{M_{n,q}} \mathbb{1}_{\mu=1^{k}}$$
pour toute partition $\mu$. 
\end{theorem}
\begin{proof}
On renvoie à \cite[\S1.3]{Bil69} pour des précisions sur la notion de convergence en probabilité ; pour des variables aléatoires réelles, cette notion est compatible avec la plupart des opérations usuelles sur les variables aléatoires, par exemple l'addition de variables. Ceci étant, comme la composante homogène de plus haut degré de $\varSigma_{\mu}$ est $p_{\mu}$, il existe une observable $X$ de degré inférieur à $|\mu|-1$ telle que $p_{\mu}=\varSigma_{\mu}+X$. Alors, en vertu du lemme \ref{stilldre},
$$M_{n,q}[p_{\mu}]=M_{n,q}[\varSigma_{\mu}]+M_{n,q}[X]=n^{\downarrow |\mu|}\, \frac{(1-q)^{|\mu|}}{1-q^{\mu}}+O(n^{|\mu|-1})\,,$$
et comme $n^{\downarrow |\mu|}=n^{|\mu|}+O(n^{|\mu|-1})$, on en déduit :
$$\frac{M_{n,q}[p_{\mu}]}{n^{|\mu|}}=\frac{(1-q)^{|\mu|}}{1-q^{\mu}}+O(n^{-1})\,.$$
Comme cette formule est multiplicative en $\mu$, l'\textbf{inégalité de Bienaymé-Chebyshev} implique alors la convergence en probabilité de $p_{\mu}(\lambda)/n^{|\mu|}$ vers  $(1-q)^{|\mu|}/(1-q^{\mu})$ :
\begin{align*}\proba\!\left[\left|\frac{p_{\mu}(\lambda)}{n^{|\mu|}} - \frac{(1-q)^{|\mu|}}{1-q^{\mu}}\right|\leq \eps\right]&\leq \frac{1}{\eps^{2}}\,\,\esper\!\left[\left(\frac{p_{\mu}(\lambda)}{n^{|\mu|}} - \frac{(1-q)^{|\mu|}}{1-q^{\mu}}\right)^{2}\right]\\
&\leq\frac{1}{\eps^{2}} \left(\frac{M_{n,q}[p_{\mu\sqcup \mu}]}{n^{2|\mu|}} - 2\,\frac{(1-q)^{|\mu|}}{1-q^{\mu}}\,\frac{M_{n,q}[p_{\mu}]}{n^{|\mu|}} + \frac{(1-q)^{2|\mu|}}{1-q^{\mu \sqcup \mu}}\right)\\
&\leq\frac{1}{\eps^{2}}\,O(n^{-1})\to 0\,.
\end{align*}
en utilisant à la dernière ligne les estimations de $M_{n,q}[p_{\mu}]$ et $M_{n,q}[p_{\mu \sqcup \mu}]$. La convergence en probabilité des moments de Frobenius est donc établie. Pour les caractères centraux, on utilise la factorisation en plus haut degré
$$\varSigma_{\mu}*\varSigma_{\mu}=\varSigma_{\mu \sqcup \mu} + \text{observable }X\text{ de degré au plus égal à }2|\mu|-1\,,$$
voir le paragraphe \ref{centralcharacter}. Alors, en utilisant de nouveau l'inégalité de Bienaymé-Che\-by\-shev :
\begin{align*}
\proba\!\left[\left|\frac{\varSigma_{\mu}(\lambda)}{n^{|\mu|}} - \frac{(1-q)^{|\mu|}}{1-q^{\mu}}\right|\leq \eps\right]&\leq \frac{1}{\eps^{2}}\,\,\esper\!\left[\left(\frac{\varSigma_{\mu}(\lambda)}{n^{|\mu|}} - \frac{(1-q)^{|\mu|}}{1-q^{\mu}}\right)^{2}\right]\\
&\!\!\!\!\!\!\!\!\leq\frac{1}{\eps^{2}}\left(\frac{M_{n,q}[(\varSigma_{\mu})^{2}]}{n^{2|\mu|}}-2\,\frac{(1-q)^{|\mu|}}{1-q^{\mu}}\,\frac{M_{n,q}[\varSigma_{\mu}]}{n^{|\mu|}} + \frac{(1-q)^{2|\mu|}}{1-q^{\mu \sqcup \mu}}\right)\\
&\!\!\!\!\!\!\!\!\leq \frac{M_{n,q}[X]}{\eps^{2}\,n^{2|\mu|}}=\frac{1}{\eps^{2}}\,O(n^{-1})\to 0\,.
\end{align*}
Enfin, pour les $q$-caractères, chaque $\varSigma_{\mu,q}$ s'exprime comme combinaison linéaire d'un nombre fini de variables aléatoires $\varSigma_{\nu}$, et ces dernières convergent toutes en probabilité après renormalisation. La variable $\varSigma_{\mu,q}(\lambda)/n^{|\mu|}$ converge donc également en probabilité, et sa limite est
\begin{align*}\frac{1}{(q-1)^{\ell(\mu)}}\sum_{\nu \in \ym_{k}} \frac{\scal{h_{\mu}}{p_{\nu}}}{\scal{p_{\nu}}{p_{\nu}}}\, (q^{\nu}-1)\,\eps_{\nu}\frac{(q-1)^{|\nu|}}{q^{\nu}-1} &=(q-1)^{|\mu|-\ell(\mu)} \scal{h_{\mu}}{\sum_{\nu \in \ym_{k}} \frac{\eps_{\nu}p_{\nu}}{z_{\nu}} }\\
=(q-1)^{|\mu|-\ell(\mu)}\scal{h_{\mu}}{e_{k}}=(q-1)^{|\mu|-\ell(\mu)}&\scal{h_{\mu}}{m_{1^{k}}}=\mathbb{1}_{\mu=1^{k}}
\end{align*}
car $(h_{\mu})_{\mu \in \ym}$ est la base duale de $(m_{\mu})_{\mu \in \ym}$.
\end{proof}\bigskip

La dernière partie du théorème \ref{convobservable} peut être interprétée comme suit : si $\chi^{\lambda}(q)$ est un $q$-caractère (normalisé) tiré aléatoirement suivant la $q$-mesure de Plancherel, et si $T_{\sigma}$ est un élément d'une algèbre d'Hecke $\IH_{q}(\sym_{n})$, alors
$$\chi^{\lambda}(q,T_{\sigma}) \rightarrow_{M_{n,q}}\tau(T_{\sigma})\,,$$ où $\tau$ désigne la trace canonique des algèbres d'Hecke --- les traces des algèbres $\IH_{q}(\sym_{n})$ et $\IH_{q}(\sym_{N \geq n})$ sont compatibles avec les plongements canoniques $\IH_{q}(\sym_{n}) \to \IH_{q}(\sym_{N})$, donc les notations $T_{\sigma}$ et $\tau$ ne font pas ambiguïté.
On observe donc un \textbf{phénomène de concentration} des $q$-caractères autour de la trace, et nous verrons dans la section \ref{qcharacter} que cette concentration est gaussienne, ce qui constitue un $q$-analogue de la première partie du théorème \ref{secondasymptoticplancherel}.

\section{Loi des grands nombres pour les premières parts}\label{qbig}
Le comportement asymptotique des observables de diagrammes ayant été précisé, nous sommes finalement prêts pour démontrer :
\begin{theorem}[Loi des grands nombres pour la $q$-mesure de Plancherel, \cite{FM10}]\label{firstasymptoticqplancherel}
Soit $q$ un paramètre réel positif strictement inférieur à $1$, et $\lambda=(\lambda_{1}\geq \lambda_{2} \geq \cdots)$ un diagramme de Young de taille $n$ tiré aléatoirement sous la $q$-mesure de Plancherel $M_{n,q}$. Pour tout $i \geq 1$,
$$\frac{\lambda_{i}}{n}\rightarrow_{M_{n,q}} (1-q)\,q^{i-1}\,.$$
\end{theorem}\bigskip

\noindent On peut donner deux preuves distinctes de ce résultat : l'une est tout à fait élémentaire et met en jeu une mesure aléatoire discrète sur $\R$, et l'autre découle de la théorie des représentations de l'algèbre d'Hecke infinie (\emph{cf.} \cite{KV07,Oko97,Oko97b}) et du théorème de Kerov \ref{boundarygeodesic}. Avant de présenter ces deux approches, donnons un corollaire important de la proposition \ref{firstasymptoticqplancherel} :
\begin{corollary}[Longueur du plus long sous-mot croissant d'une permutation obtenue par le $q$-processus de Knuth]\label{firstknuth}
Soit $\sigma_{n}$ une permutation aléatoire de taille $n$ tirée suivant la mesure de probabilité $$P_{n,q}[\sigma]=\frac{q^{\imaj(\sigma)}}{\{n!\}_{q}}\,.$$
La longueur $L(\sigma_{n})$ d'un plus long sous-mot croissant de $\sigma_{n}$ a pour asymptotique
$$\frac{L(\sigma_{n})}{n} \rightarrow_{P_{n,q}} 1-q\,.$$
\end{corollary}
\noindent En effet, ceci est une conséquence immédiate de la proposition \ref{qmodel} et des propriétés de la correspondance RSK.\bigskip

\begin{proof}[Première démonstration de la loi des grands nombres]
Pour commencer, remarquons que la limite de $p_{k}(\lambda)/n^{k}$ sous la $q$-mesure de Plancherel s'écrit 
$$\frac{(1-q)^{k}}{1-q^{k}}=\sum_{i=1}^{\infty}[(1-q)\,q^{i-1}]^{k}\,.$$
Si $\lambda$ est un diagramme de Young de coordonnées de Frobenius $(a_{i},b_{i})_{1 \leq i \leq d}$ et de taille $n$, notons $X_{\lambda}$ la mesure de probabilité :
$$X_{\lambda}=\sum_{i=1}^{d} \frac{a_{i}(\lambda)}{n}\,\delta_{a_{i}(\lambda)/n} +\sum_{i=1}^{d}\frac{b_{i}(\lambda)}{n}\,\delta_{-b_{i}(\lambda)/n}\,.$$
Le $k$-ième moment de $X_{\lambda}$ est exactement $p_{k+1}(\lambda)/n^{k+1}$ ; d'après le théorème \ref{convobservable}, tous les moments de $X_{\lambda}$ convergent donc en probabilité vers ceux de la mesure
$$X_{\infty,q}=\sum_{i=1}^{\infty} (1-q)\,q^{i-1}\,\delta_{(1-q)\,q^{i-1}}\,.$$
Dans ce qui suit, on considère $X_{\lambda}$ comme une variable aléatoire à valeur mesure ; sa loi est l'image de $M_{n,q}$ par l'application $\lambda \in \ym \mapsto X_{\lambda} \in \mathscr{M}([-1,1])$, $\mathscr{M}([-1,1])$ désignant l'ensemble des mesures de Radon signées sur l'intervalle $[-1,1]$. Ce dernier espace est le dual topologique de $\mathscr{C}([-1,1])$, et l'ensemble des mesures de probabilité $\mathscr{P}([-1,1]) \subset \mathscr{M}([-1,1])$ est compact métrisable pour la topologie $*$-faible\footnote{C'est un cas particulier du \textbf{théorème de Prohorov} (\cite[chapitre 1]{Bil69}) : si $E$ est un espace séparable et métrisable complet, alors l'espace des mesures de probabilité $\mathscr{P}(E)$ est également séparable et métrisable complet. Le caractère métrisable permet en particulier de parler de convergence en probabilité d'une famille de mesures.}. En effet, si $(f_{k})_{k\in \N}$ est une suite de fonctions dense dans $\mathscr{C}([-1,1])$, alors 
$$d(m_{1},m_{2})=\sum_{k \in \N}\frac{1}{2^{k}}\,\max(1,|m_{1}(f_{k})-m_{2}(f_{k})|)$$
est une distance compatible avec la topologie de $\mathscr{M}([-1,1])$. On choisit $\{f_{k}\}_{k \in \N}=\Q[X]$ ; cette partie est dense par le théorème de Stone-Weierstrass, et compte tenu de la convergence des moments, 
$$\forall k,\,\,\,X_{\lambda}(f_{k}) \rightarrow_{M_{n,q}} X_{\infty,q}(f_{k})\quad\iff\quad d(X_{\lambda},X_{\infty,q}) \rightarrow_{M_{n,q}} 0\,.$$
Dans l'espace des mesures muni de la topologie (métrisable !) de la convergence en loi, on a donc convergence en probabilité de $X_{\lambda}$ vers $X_{\infty,q}$. Or, pour des mesures de Radon sur l'ensemble des réels, la convergence en loi est équivalente à la convergence des fonctions de répartitions 
$$F_{X_{\lambda}}(x) \to F_{X_{\infty,q}}(x)$$
en tout point $x$ où la fonction de répartition limite $F_{X_{\infty,q}}$ est continue (c'est une partie du \textbf{théorème de Portmanteau}, voir de nouveau \cite[chapitre 1]{Bil69}). Par conséquent, pour tout nombre réel $x \neq (1-q)\,q^{i-1}$, on a
$$\sum_{\substack{y \leq x\\ ny \in A(\lambda)\sqcup -B(\lambda)}}\!\!\!\!\!\!\!\!|y| \quad\rightarrow_{M_{n,q}}\quad\sum_{(1-q)\,q^{i-1}\leq x} \!\!(1-q)\,q^{i-1}\,.$$
La fin de la preuve est maintenant purement technique. Dans tout ce qui suit, on fixe un réel $\eps$ tel que $1-q>\eps>0$ ; d'après ce qui précède, pour tout point $x$ qui n'est pas un $(1-q)\,q^{i-1}$, 
$$ F_{X_{\infty,q}}(x)-\eps \leq F_{X_{\lambda}}(x) \leq F_{X_{\infty,q}}(x)+\eps$$
pour $n$ assez grand en dehors d'un événement de probabilité arbitrairement petite. Prenons donc $x=1-q+\eta$ avec $\eta>0$ suffisamment petit. Dans un voisinage de $x$, $F_{X_{\infty,q}}=1$, donc pour $n$ assez grand, en dehors d'un événement de petite probabilité,
$$F_{X_{\lambda}}(x)=\!\!\!\!\sum_{\substack{y \leq x\\ny \in A(\lambda) \sqcup -B(\lambda)}}\!\!\!\!\!\!\!\!|y| \,\,\,\,\geq\,\, F_{X_{\infty,q}}(x)-\eps =1-\eps >q\,.$$
Par conséquent, toutes les lignes\footnote{Pour être plus précis, les lignes \emph{modifiées}, c'est-à-dire qu'on regarde les $a_{i}(\lambda)=\lambda_{i}-i+1/2$.} de $\lambda$ sont plus petites que $nx$. En effet, dans le cas contraire, il y aurait  au moins une ligne strictement plus grande que $n(1-q+\eta)$, et ceci impliquerait 
$$F_{X_{\lambda}}(1) \geq F_{X_{\lambda}(x)}+1-q+\eta \geq 1+\eta >1\,,$$
ce qui est impossible pour une fonction de répartition. Ainsi, on a montré que pour $n$ assez grand, $a_{1}(\lambda)/n$ est plus petit que $x$ avec probabilité très proche de $1$. \bigskip

\noindent De même, si l'on considère $x'=1-q-\eta$, alors $F_{X_{\infty,q}}=q$ au voisinage de $x'$, donc pour $n$ assez grand, en dehors d'un événement de petite probabilité, certaines des lignes de $\lambda$ sont plus grandes que $nx'$. En effet, dans le cas contraire, on aurait
$$F_{X_{\lambda}}(1)=F_{X_{\lambda}}(x')\leq F_{X_{\infty,q}}(x') + \eps =q+\eps <1\,.$$
On conclut que pour $n$ assez grand, en dehors d'un événement de probabilité arbitrairement petite, la première colonne renormalisée $a_{1}(\lambda)/n$ est comprise entre $x'$ et $x$. Ainsi, $a_{1}(\lambda)/n$ converge en probabilité vers $1-q$, et les mêmes arguments s'appliquent clairement aux lignes suivantes $a_{2}(\lambda)/n$, $a_{3}(\lambda)/n$, etc. (faire une récurrence sur $i$). De plus, $$a_{i}(\lambda)/n-\lambda_{i}/n=(-i+1/2)/n=O(1/n)\,,$$ donc $\lambda_{i}/n$ converge bien en probabilité vers $(1-q)\,q^{i-1}$. 
\end{proof}
\noindent Remarquons qu'en utilisant la convergence des fonctions de répartition au point $x=0$, on voit que
$$\sum_{i=1}^{d } \frac{b_{i}(\lambda)}{n} \rightarrow_{M_{n,q}} 0\,\,;$$
l'ordre de grandeur des colonnes de la partition est donc un $o(n)$.\bigskip

\begin{proof}[Seconde démonstration de la loi des grands nombres]
Une autre démonstration du théorème \ref{firstasymptoticqplancherel} utilise la théorie des représentations de l'algèbre d'Hecke infinie $\IH_{q}(\sym_{\infty})$ ; cette théorie est en tout point similaire à celle exposée dans le paragraphe \ref{syminfinite} pour le groupe symétrique infini $\sym_{\infty}$. Pour tous entiers $n \leq N$, on peut plonger $\IH_{q}(\sym_{n})$ dans $\IH_{q}(\sym_{N})$ en identifiant les générateurs $T_{i}=T_{s_{i}}$ pour $i \leq n-1$ ; ces injections sont compatibles entre elles, et la limite inductive des algèbres d'Hecke $\IH_{q}(\sym_{n})$ relativement à cette famille dirigée de morphismes d'algèbres est l'\textbf{algèbre d'Hecke infinie}
$$\IH_{q}(\sym_{\infty})=\left\langle T_{1},T_{2},T_{3},\ldots \,\,\bigg|\,\,\substack{\forall i,\,\,(T_{i}-q)(T_{i}+1)=0\quad\,\,\,\,\,\\
\forall i,\,\,T_{i}T_{i+1}T_{i}=T_{i+1}T_{i}T_{i+1}\,\,\,\\
\forall i,j,\,\,|i-j|\geq 2 \Rightarrow T_{i}T_{j}=T_{j}T_{i}\,} \right\rangle\,.$$
Les caractères irréductibles de cette algèbre infini-dimensionnelle sont de nouveau indexés par les points $\omega=(\alpha,\beta)$ du simplexe de Thoma $\Omega$, et toute mesure de probabilité $m \in \mathscr{P}(\Omega)$ correspond à un caractère (positif, normalisé) $\chi(q)$ de $\IH_{q}(\sym_{\infty})$, qui s'écrit :
$$\chi(T_{\sigma \in \sym_{n}},q)=\sum_{\lambda \in \ym_{n}} \left(\int_{\Omega} s_{\lambda}(\omega) \,\varsigma^{\lambda}(q,T_{\sigma})\,m(d\omega) \right)\,.$$
En particulier, la trace régulière\footnote{Pour une algèbre d'Hecke finie $\IH_{q}(\sym_{n})\simeq \C[\GB(n,\For_{q})\backslash \GL(n,\For_{q})/\GB(n,\For_{q})]$, la trace régulière est la restriction de la trace normalisée de $\hendo_{\C}(\C[\GL(n,\For_{q})/\GB(n,\For_{q}])$. Il est naturel de se demander si la trace $\tau$ de $\IH_{q}(\sym_{\infty})$ a une interprétation semblable (avec $q$ puissance d'un nombre premier). La réponse est oui, mais on doit prendre au lieu de $\GL(\infty,\For_{q})$ le groupe $\mathrm{GLB}(\For_{q})$ introduit par Kerov et Vershik dans \cite{KV98,KV07} ; ce dernier contient $\GL(\infty,\For_{q})$ comme sous-groupe dense. Alors, on peut montrer que $\tau$ est la restriction de la trace canonique de l'algèbre de convolution $\leb^{1}(\mathrm{GLB}(\For_{q}),m_{\mathrm{Haar}})$ à la sous-algèbre des fonctions à support compact et bi-$\GB(\infty,q)$-invariantes, cette sous-algèbre étant isomorphe à $\IH_{q}(\sym_{\infty})$.} $\tau(T_{\sigma})=\mathbb{1}_{(\sigma=\id)}$ de l'algèbre d'Hecke infinie correspond au Dirac en le point 
$$\alpha_{q},\beta_{q}=\big(1-q,(1-q)q,(1-q)q^{2},\ldots\big),\big(0,0,\ldots\big)\,.$$
En effet, compte tenu des calculs effectués dans la section \ref{qplancherelprocess}, $$M_{n,q}(\lambda)=(\dim \lambda) \, s_{\lambda}(1-q,(1-q)q,(1-q)q^{2},\ldots)=(\dim\lambda) \,s_{\lambda}(\alpha_{q},\beta_{q})\,,$$
 et d'autre part, $\tau_{|\IH_{q}(\sym_{n})}=\sum_{\lambda \in \ym_{n}} M_{n,q}(\lambda)\,\chi^{\lambda}(q)$. Si $T_{\sigma}$ est un élément de base de l'algèbre d'Hecke $\IH_{q}(\sym_{n})$, on a donc bien :
 $$\tau(T_{\sigma})=\sum_{\lambda \in \ym_{n}} s_{\lambda}(\alpha_{q},\beta_{q})\, \varsigma^{\lambda}(q,T_{\sigma})=\sum_{\lambda \in \ym_{n}} \left(\int_{\Omega}s_{\lambda}(\omega) \,\varsigma^{\lambda}(q,T_{\sigma})\,\delta_{\alpha_{q},\beta_{q}}(d\omega)\right)\,.$$
 Dit autrement, le Dirac en le point $(\alpha_{q},\beta_{q})$ correspond au système cohérent de mesures sur le graphe de Young constitué par les $q$-mesures de Plancherel --- la fonction $M_{n,q}(\lambda)/\dim \lambda$ est bien harmonique sur le graphe de Young, car les $q$-caractères ont les mêmes règles de branchement que les caractères usuels. En vertu du théorème de Kerov \ref{boundarygeodesic}, on a donc convergence en loi de $\omega_{n}(M_{n,q})$ vers le Dirac en $(\alpha_{q},\beta_{q})$, c'est-à-dire que :
 $$\left(\frac{a_{1}(\lambda)}{n},\frac{a_{2}(\lambda)}{n},\ldots ; \frac{b_{1}(\lambda)}{n},\frac{b_{2}(\lambda)}{n},\ldots \right) \,\,\to\,\, \big(1-q,(1-q)q,\ldots ; 0,0,\ldots\big)\,.$$
Sur un espace métrique (ici, le simplexe de Thoma $\Omega$), la convergence en loi vers une constante implique la convergence en probabilité ; on retrouve donc bien la loi des grands nombres \ref{firstasymptoticplancherel}. \end{proof}

\noindent La seconde démonstration n'utilise absolument pas les techniques d'observables de dia\-gram\-mes développées jusqu'ici ; en contrepartie, elle s'appuie sur le résultat difficile de Kerov et Vershik (théorème \ref{boundarygeodesic}), et surtout, elle ne permet pas d'étudier la déviation gaussienne des diagrammes sous les $q$-mesures de Plancherel. D'autre part, elle ne s'adapte pas à l'étude d'algèbres d'Hecke plus générales, par exemple celles de type B (voir le chapitre \ref{arikikoike}).\bigskip

\section{Théorème central limite pour les premières parts}\label{qgaussian}
Nous souhaitons maintenant démontrer un théorème central limite pour les lignes d'un diagramme $\lambda$ tiré aléatoirement suivant une $q$-mesure de Plancherel. Comme les lignes $\lambda_{i}$ ont toutes pour ordre de grandeur $n$, en vertu du théorème central limite <<~usuel~>>, il est naturel de supposer que les déviations $\lambda_{i}-n\,(1-q)\,q^{i-1}$ auront pour ordre de grandeur $\sqrt{n}$. On pose donc :
$$Y_{i,n,q}=\sqrt{n}\,\left(\frac{\lambda_{i}}{n}-(1-q)\,q^{i-1}\right)\,.$$
\begin{theorem}[Déviation gaussienne pour la $q$-mesure de Plancherel, \cite{FM10}]\label{secondasymptoticqplancherel}
Sous une $q$-mesure de Plancherel de paramètre $q<1$, les lois fini-dimensionnelles du processus $(Y_{i,n,q})_{i \geq 1}$ con\-ver\-gent vers celles d'un processus gaussien $(Y_{i,\infty,q})_{i \geq 1}$, avec :
\begin{align*}
&\esper[Y_{i,\infty,q}]=0\qquad;\qquad\esper[(Y_{i,\infty,q})^{2}]=(1-q)\,q^{i-1}-(1-q)^{2}\,q^{2(i-1)}\quad;\\
&\mathrm{cov}(Y_{i,\infty,q},Y_{j,\infty,q})=-(1-q)^{2}\,q^{i+j-2}\,.
\end{align*}
\end{theorem}\bigskip

\noindent En particulier, deux coordonnées distinctes $Y_{i,\infty,q}$ et $Y_{j,\infty,q}$ sont toujours négativement corrélées, ce qui n'est pas très étonnant puisque la somme des déviations doit être nulle. Ainsi, si $\eps>0$ et $x_{1},\ldots,x_{r}$ sont des nombres réels positifs, alors pour $n$ assez grand, 
$$\left|\proba\!\left[ \forall i \leq r,\,\, \left|\frac{\lambda_i}{n}- (1-q)\,q^{i-1}\,\right|\leq \frac{x_i}{\sqrt{n}} \right]-\frac{1}{\sqrt{(2\pi)^r\,\det Q}}\,\int_{-x_{1}}^{x_1}\!\!\cdots\int_{-x_{r}}^{x_r}\!\E^{-\frac{{}^tXQ^{-1}X}{2}}\,dX\right| \leq \eps\,,$$
où $Q$ est la matrice symétrique de coefficients $Q_{ij}=\delta_{ij}\,(1-q)\,q^{i-1}-(1-q)^2\,q^{i+j-2}$, qui est définie positive par la règle d'Hadamard (sur chaque ligne et chaque colonne, le coefficient diagonal est positif et strictement plus grand que la somme des valeurs absolues des autres coefficients). Comme précédemment, notre théorème central limite a un corollaire immédiat pour les permutations aléatoires :
\begin{corollary}[Déviation gaussienne pour la longueur du plus long sous-mot croissant d'une permutation obtenue par le $q$-processus de Knuth]\label{secondknuth}
Soit $\sigma_{n}$ une permutation aléatoire de taille $n$ tirée suivant la mesure de probabilité $$P_{n,q}[\sigma]=\frac{q^{\imaj(\sigma)}}{\{n!\}_{q}}\,.$$
La variable aléatoire $(L(\sigma_{n})-n(1-q))/\sqrt{n}$ converge en loi vers une gaussienne centrée de variance $(1-q)-(1-q)^{2}$.
\end{corollary}\bigskip

\begin{remark}
Si l'on applique le processus d'orthonormalisation de Gram-Schmidt à la matrice de covariance donnée dans le théorème \ref{secondasymptoticqplancherel}, alors on peut exprimer les lois limites $Y_{i,\infty,q}$ en fonctions de gaussiennes standards indépendantes. Ainsi, il existe une famille $(\xi_{i})_{i \geq 1}$ de gaussiennes indépendantes, centrées et de variance $1$ telle que pour tout $i$,
$$Y_{i,\infty,q}=\sqrt{(1-q)q^{i}}\,\left(\xi_{i}-(1-q)\sum_{j=1}^{i-1} q^{\frac{i-j-2}{2}}\,\xi_{j}\right)\,,$$
l'égalité ayant lieu presque sûrement si l'on définit les $\xi_{i}$ par :
$$\xi_{i}=\frac{1}{\sqrt{(1-q)q^{i}}}\,\left(Y_{i,\infty,q}+(1-q)\sum_{j=1}^{i-1} Y_{j,\infty,q}\right)\,.$$
Ainsi, pour simuler une déviation $Y_{i,n,q}$ avec $n$ grand, on commence par tirer une gaussienne $\sqrt{(1-q)q^{i}} \,\xi_{i}$, et on lui retranche la quantité $\sqrt{(1-q)^{3}q^{i}}\,\sum_{j=1}^{i-1}q^{\frac{i-j-2}{2}}\,\xi_{j}$ pour compenser les déviations des lignes précédentes.
\end{remark}
\bigskip
\bigskip

La preuve du théorème \ref{secondasymptoticqplancherel} met en jeu les cumulants d'observables de la section \ref{sniady}, mais nous devrons revisiter quelque peu la théorie de \'Sniady, car les espérances d'observables n'ont pas le même ordre de grandeur sous la $q$-mesure de Plancherel et sous la mesure de Plancherel standard : en effet,
$$M_{n,q}[\varSigma_{\mu}]=O\left(n^{|\mu|}\right)\qquad;\qquad M_{n}[\varSigma_{\mu}]=O\left(n^{\frac{|\mu|+\ell(\mu)}{2}}\right)\,.$$
Dans tout ce qui suit, on considère les variables $\varSigma_{\mu}$ comme éléments de l'algèbre $\blg_{\infty}$ des permutations partielles. La $q$-mesure de Plancherel $M_{n,q}$ peut être vue comme une trace sur l'algèbre $\alg_{\infty}$ engendrée par les variables commutantes $\varSigma_{\mu}$, d'où une structure d'espace de probabilité non commutatif sur $\alg_{\infty}$, avec les règles :
$$\esper[\varSigma_{\mu}]=n^{\downarrow |\mu|}\,\frac{(1-q)^{|\mu|}}{1-q^{\mu}}\,.$$
Les cumulants correspondants seront notés $k(x_{1},\ldots,x_{r})$ ; d'autre part, nous rappelons l'existence d'un second produit $\bullet$ dans $\alg_{\infty}$, et les cumulants disjoints correspondants seront notés $k^{\bullet}(x_{1},\ldots,x_{r})$. Comme $\alg_{\infty}\simeq \obs$, les mêmes notions sont disponibles pour des observables de diagrammes. Ceci étant dit, la difficulté combinatoire principale de la preuve réside dans les lemmes suivants :
\begin{lemma}[Ordre de grandeur des cumulants disjoints]\label{magnitudedj}
Sous les $q$-mesures de Plancherel, pour toutes observables de diagrammes $x_{1},\ldots,x_{r}$, $k^{\bullet}(x_{1},\ldots,x_{r})=O(n^{\deg(x_{1})+\cdots+\deg(x_{r})-r+1})$.
\end{lemma}
\begin{lemma}[Ordre de grandeur des cumulants]\label{magnitude}
Sous les $q$-mesures de Plancherel, pour toutes observables de diagrammes $x_{1},\ldots,x_{r}$, $k(x_{1},\ldots,x_{r})=O(n^{\deg(x_{1})+\cdots+\deg(x_{r})-r+1})$.
\end{lemma}\bigskip
\bigskip

Il suffit en réalité de démontrer ces lemmes avec les $x_{j}$ dans une base algébrique graduée de $\obs$ ; en effet, les cumulants de variables aléatoires sont multilinéaires et se comportent bien vis-à-vis des produits de variables, voir \cite{LS59} et \cite[théorème 4.4]{Sni06b}. Ainsi, dans tout ce qui suit, les $x_{j}$ seront des caractères centraux de cycles $\varSigma_{i_{j}}$.
\begin{proof}[Preuve du lemme \ref{magnitudedj}]
Rappelons que les cumulants disjoints sont récursivement définis par les relations :
$$\esper[X_{1}\bullet X_{2} \bullet \cdots \bullet X_{r}]=\sum_{\pi \in \mathfrak{Q}(\lle 1,r \rre)}
k^{\bullet}(X_{i \in \pi_{1}})\,k^{\bullet}(X_{i \in \pi_{2}}) \cdots k^{\bullet}(X_{i \in \pi_{l}})\,.$$
D'autre part, le produit disjoint de caractères centraux de cycles $\varSigma_{i_{1}},\ldots,\varSigma_{i_{r}}$ (avec $i_{1}\geq i_{2}\geq \cdots \geq i_{r}$) est simplement $\varSigma_{i_{1},\ldots,i_{r}}$, et l'espérance d'une telle variable est
$$\esper[\varSigma_{i_{1},\ldots,i_{r}}]=n^{\downarrow i_{1}+i_{2}+\cdots +i_{r}} \frac{(1-q)^{i_{1}+i_{2}+\cdots+i_{r}}}{\prod_{j=1}^{r}(1-q^{i_{j}})}\,.$$
\`A l'exception du facteur $n^{\downarrow i_{1}+\cdots+i_{r}}$, tous les termes de cette expression sont multiplicatifs par rapport à la partition $(i_{1},\ldots,i_{r})$. Plus formellement, il existe des nombres $\alpha_{i}$ tels que 
$$\esper[\varSigma_{i_{1},\ldots,i_{r}}]=\left(\prod_{j=1}^{r} \alpha_{i_{j}}\right)n^{\downarrow i_{1}+i_{2}+\cdots+i_{r}}\,,$$
et la puissance décroissante de $n$ peut s'écrire $\esper[\varSigma_{(1^{i_{1}+\cdots +i_{r}})}]$. Par récurrence sur $r$ et en utilisant la formule reliant espérances et cumulants disjoints, on voit donc que :
$$k^{\bullet}(\varSigma_{i_{1}},\varSigma_{i_{2}},\ldots,\varSigma_{i_{r}})=\left(\prod_{j=1}^{r}\alpha_{i_{j}}\right)k^{\bullet}(\varSigma_{(1^{i_{1}})},\varSigma_{(1^{i_{2}})},\ldots,\varSigma_{(1^{i_{r}})})\,.$$
Or, $\chi^{\lambda}(\varSigma_{(1^{i})})=n^{\downarrow i}$ pour toute partition $\lambda$, donc les cumulants correspondants ne dépendent pas de la mesure de probabilité que l'on met sur les partitions (!). On peut donc utiliser le résultat de P. \'Sniady (\emph{cf.} \cite[lemme 4.8]{Sni06b}), qui a montré que :
$$k^{\bullet}(\varSigma_{(1^{i_{1}})},\varSigma_{(1^{i_{2}})},\ldots,\varSigma_{(1^{i_{r}})})=O(n^{i_{1}+i_{2}+\cdots+i_{r}-r+1})\,.\vspace{-5mm}~$$
\end{proof}\bigskip\bigskip

Il s'agit maintenant de relier cumulants classiques et cumulants disjoints, et nous utiliserons comme dans la section \ref{sniady} les cumulants $k^{\id}$. Ces derniers correspondent à <<~l'espérance~>> $\esper^{\id} : x \in \obs \mapsto x \in \obs^{\bullet}$, et ils sont définis récursivement par les relations :
$$X_{1}X_{2} \cdots X_{r}=\sum_{\pi \in \mathfrak{Q}(\lle 1,r \rre)}
k^{\id}(X_{i \in \pi_{1}})\bullet k^{\id}(X_{i \in \pi_{2}})\bullet \cdots \bullet k^{\id}(X_{i \in \pi_{l}})\,.$$
On rappelle que dans ces conditions, cumulants classiques et cumulants disjoints sont reliés par la formule :
$$k(X_{1},\ldots,X_{r})=\sum_{\pi \in \mathfrak{Q}(\lle 1,r\rre)}k^{\bullet}(k^{\mathrm{id}}(X_{i \in \pi_{1}}),\ldots,k^{\mathrm{id}}(X_{i \in \pi_{l}}))\,,$$
voir \cite[proposition 4.1]{Sni06b}. Avant d'étudier les cumulants $k^{\id}(\varSigma_{i_{1}},\ldots,\varSigma_{i_{r}})$, nous devons préciser quelque peu les résultats évoqués dans la section \ref{centralcharacter} au sujet des produits de caractères centraux. On note $\varSigma_{\mu,n}$ le projeté d'un caractère central $\varSigma_{\mu}$ dans l'algèbre $\blg_{n}$ constituée par les permutations partielles d'ordre $n$. Alors, si $\arra(n,k)$ est l'ensemble des $k$-arrangements $(a_{1}\neq a_{2} \neq \cdots \neq a_{k})$ dans $\lle 1,n\rre$, on peut écrire :
\begin{align*}\varSigma_{\mu,n}&=\sum_{\substack{a_{(i,j)} \in \lle 1,n \rre \\ a_{(i,j)}\neq a_{(i,j)'}}} (a_{11},a_{12},\ldots,a_{1\mu_1})\,(a_{21},a_{22},\ldots,a_{2\mu_2})\,\cdots\,(a_{r1},a_{r2},\ldots,a_{r\mu_r})\\
&=\sum_{\substack{\forall i,\,\,A_{i} \in \arra(n,\mu_{i})\\
\forall i \neq j,\,\,A_{i} \cap A_{j} =\emptyset}} C(A_{1})\,C(A_{2})\,\cdots\,C(A_{r})\end{align*}
où $C(A)$ est le cycle associé à un arrangement $A=(a_{1},\ldots,a_{k})$. Par conséquent, un produit $\varSigma_{\mu,n}\,\varSigma_{\nu,n}$ s'écrit 
$$\varSigma_{\mu,n} \, \varSigma_{\nu,n} = \sum_{\substack{\forall i,\,\,A_i \in \arra(n,\,\mu_i)\\ \forall i\neq j,\,\,A_i\cap A_j=\emptyset}} \sum_{ \substack{\forall i,\,\,B_i \in \arra(n,\,\nu_i)\\ \forall i\neq j,\,\,B_i\cap B_j=\emptyset}} C(A_1)\,\cdots \,C(A_r)\, C(B_1)\,\cdots \,C(B_s)\,.$$
Ce produit n'est pas $\varSigma_{\mu\sqcup\nu,n}$, car les $A_{i}$ et les $B_{j}$ peuvent avoir des éléments en commun. Pour prendre en compte les éventuels <<~recouvrements~>>, associons à toute suite $(a_{ij},b_{kl})$ intervenant dans la somme un \textbf{appariement partiel} $M$ des ensembles $I_{A}=(i,j)_{1\leq j \leq \mu_{i}}$ et $I_{B}=(k,l)_{1\leq k \leq \nu_{l}}$ :
$$(i,j) \sim_{M} (k,l) \iff a_{ij}=b_{kl}\,.$$
La somme des produits $C(A_1)\,\cdots \,C(A_r)\, C(B_1)\,\cdots \,C(B_s)$ réalisée
sur les arrangements correspondant à un appariement partiel $M$ donné est un caractère central $\varSigma_{\rho(M)}$, car si l'on sait exactement quels $a_{ij}$ sont égaux à quels $b_{kl}$, alors on peut réécrire le produit $$C(A_1)\,\cdots \,C(A_r)\, C(B_1)\,\cdots \,C(B_s)$$ comme produit de cycles à supports disjoints, et ce de fa\c con indépendante des valeurs précises des $a_{ij}$ et des $b_{kl}$. Ainsi :\bigskip

\begin{lemma}[Produits de caractères centraux]\label{productcentralcharacter}
Le produit $\varSigma_{\mu} \,\varSigma_{\nu}$ de deux caractères centraux est la somme $\sum_{M} \varSigma_{\rho(M)}$ de caractères centraux, où $M$ parcourt tous les appariements partiels possibles entre l'ensemble $I_{A}$ des indices des $a$ et l'ensemble $I_{B}$ des indices des $b$.
\end{lemma}\bigskip

\begin{example}
\`A titre d'exemple, détaillons le calcul du produit $\varSigma_{2}\,\varSigma_{3}$. L'ensemble d'indices $I_{A}$ est $\{ 1,2\}$, et l'ensemble d'indices $I_{B}$ est $\{1',2',3'\}$. \vspace{2mm}
\begin{enumerate}
\item Si les arrangements $A=(a_{1},a_{2})$ et $B=(b_{1},b_{2},b_{3})$ sont disjoints, alors $C(A)\,C(B)$ est le produit de deux cycles disjoints de tailles respectives $2$ et $3$. L'appariement partiel vide correspond donc au caractère central $\varSigma_{3,2}$.\vspace{2mm}
\item Si $A$ et $B$ ont exactement un élément commun, on peut à symétrie près supposer que c'est $a_{1}=b_{1}$. Alors, $(a_{1},a_{2})(b_{1},b_{2},b_{3})=(a_{2},b_{1},b_{2},b_{3})$ est un $4$-cycle. En sommant sur toutes les possibilités pour $a_{2}$, $b_{1}$, $b_{2}$ et $b_{3}$ (qui doivent être distincts), on obtient le caractère central $\varSigma_{4}$. Par symétrie, on obtient aussi $\varSigma_{4}$ pour tous les autres appariements partiels de taille $1$ (il y en $6$).\vspace{2mm}
\item Si $A$ et $B$ ont deux éléments en commun, on peut à symétrie près supposer que l'appariement partiel est $a_{1}=b_{1}$ et $a_{2}=b_{2}$. Alors, $(a_{1},a_{2})(b_{1},b_{2},b_{3})=(a_{1})(a_{2},a_{3})$ est le produit d'un $1$-cycle et d'un $2$-cycle à supports disjoints. En sommant sur toutes les possibilités pour $a_{1}$, $a_{2}$ et $a_{3}$, on obtient donc le caractère central $\varSigma_{2,1}$. Et de nouveau, tous les autres appariements partiels de taille $2$ donnent ce caractère central (et il y en $6$).\vspace{2mm}
\end{enumerate}
Ainsi, $\varSigma_{2}\,\varSigma_{3}=\varSigma_{3,2}+6\,\varSigma_{4}+6\,\varSigma_{2,1}$.
\end{example}\bigskip

De fa\c con générale, la taille d'une partition $\rho(M)$ apparaissant dans la somme $\sum_{M}\varSigma_{\rho(M)}=\varSigma_{\mu}\,\varSigma_{\nu}$ est toujours $|\mu|+|\nu|-|M|$, où $M$ désigne la taille de l'appariement, c'est-à-dire le nombre d'égalités $a_{ij}=b_{kl}$. Ce qui précède apporte donc une preuve de la propriété de factorisation en plus haut degré évoquée dans la proposition \ref{whitelight}. De plus, lorsqu'on effectue des produits de caractères centraux de cycles $\varSigma_{l}$ et $\varSigma_{m}$, tous les appariements de taille $1$ donnent un caractère central $\varSigma_{l+m-1}$. Comme il y a $ml$ appariements de ce type, on obtient le développement suivant :
$$\varSigma_{l}\,\varSigma_{m}=\varSigma_{l,m} +ml\,\varSigma_{l+m-1}+(\text{termes de degré inférieur à }l+m-2)\,.$$
Ce développement jouera un rôle important dans la suite. Ceci étant, revenons au problème des cumulants $k^{\id}(\varSigma_{i_{1}},\ldots,\varSigma_{i_{r}})$. Pour $r=1$, $k^{\id}(\varSigma_{i})=\varSigma_{i}$, et pour $r=2$, on a :
$$\varSigma_{i_{1}}\varSigma_{i_{2}}=k^{\id}(\varSigma_{i_{1}},\varSigma_{i_{2}})+k^{\id}(\varSigma_{i_{1}})\bullet k^{\id}(\varSigma_{i_{1}})\qquad;\qquad k^{\id}(\varSigma_{i_{1}},\varSigma_{i_{2}})=\varSigma_{i_{1}}\varSigma_{i_{2}}-\varSigma_{i_{1},i_{2}}\,.$$
Or, le produit $\varSigma_{i_{1}}\varSigma_{i_{2}}$ projeté dans $\blg_{n}$ s'écrit
$$\sum_{A \in \arra(n,i_{1}),\,B \in \arra(n,i_{2})} C(A)\,C(B)\,,$$
et le terme $\varSigma_{i_{1},i_{2}}$ projeté dans $\blg_{n}$ s'écrit
$$\sum_{\substack{A \in \arra(n,i_{1}),\,B \in \arra(n,i_{2})\\ A \cap B=\emptyset}} C(A)\,C(B)\,.$$
Par conséquent, le cumulant $k^{\id}(\varSigma_{i_{1}},\varSigma_{i_{2}})$ est la même somme, mais restreinte aux paires d'arrangements qui ont une intersection non vide :
$$k^{\id}(\varSigma_{i_{1}},\varSigma_{i_{2}})=\sum_{\substack{A \in \arra(n,i_{1}),\,B \in \arra(n,i_{2})\\ A \cap B\neq \emptyset}} C(A)\,C(B)\,.$$
Plus généralement, étant donnée une famille d'arrangements $(A_{1},A_{2},\ldots,A_{r})$, considérons la relation $\sim$ définie par $k \sim l \iff A_{k}\cap A_{l}\neq \emptyset$. Cette relation peut être complétée en une relation d'équivalence sur $\lle 1,r \rre$, et nous noterons $\pi(A_{1},\ldots,A_{r})$ la partition d'ensemble de $\lle 1, r\rre$ correspondant à cette relation d'équivalence. Alors :
\begin{lemma}[Cumulant identité d'une famille de caractères centraux de cycles]\label{sanitarium}
Projeté dans $\blg_{n}$, le cumulant $k^{\id}(\varSigma_{i_{1}},\varSigma_{i_{2}},\ldots,\varSigma_{i_{r}})$ vaut
$$\sum_{\substack{\forall j,\,\,A_{j} \in \arra(n,i_{j})\\ \pi(A_{1},\ldots,A_{r})=\{\lle 1,r \rre\}}} C(A_{1})\,C(A_{2})\,\cdots\,C(A_{r})\,.$$
\end{lemma}
\begin{proof}
On raisonne par récurrence sur $r$ ; les cas $r=1$ et $r=2$ viennent d'être prouvés. Si le résultat est vrai jusqu'au rang $s \leq r-1$, alors par définition du produit disjoint $\bullet$, pour toute partition d'ensemble non triviale $\pi \in \mathfrak{Q}(\lle 1,r \rre)$, on a :
$$k^{\id}(\varSigma_{i \in \pi_{1}})\bullet k^{\id}(\varSigma_{i \in \pi_{2}})\bullet \cdots \bullet k^{\id}(\varSigma_{i \in \pi_{l}})=\sum_{\substack{\forall j,\,\,A_{j} \in \arra(n,i_{j})\\ \pi(A_{1},\ldots,A_{r})=\pi }} C(A_{1})\,C(A_{2})\,\cdots\,C(A_{r})\,.$$
En effet, l'hypothèse de récurrence s'applique dans chaque part $\pi_{j}$, donc la somme porte sur des familles d'arrangements telles que si $k$ et $l$ sont dans la même part $\pi_{j}$, alors $k$ et $l$ sont dans la même part de $\pi(A_{1},\ldots,A_{r})$. Réciproquement, si $k$ et $l$ sont dans deux parts $\pi_{i}$ et $\pi_{j}$ distinctes, alors le produit $\bullet$ force les arrangements $A_{k}$ et $A_{l}$ apparaissant dans la somme à être disjoints (sinon, le produit disjoint des cycles correspondants est nul) ; donc, $k$ et $l$ sont dans deux parts disjointes de $\pi(A_{1},\ldots,A_{r})$. On conclut que $\pi=\pi(A_{1},\ldots,A_{r})$ pour toutes les familles d'arrangements apparaissant dans la somme, d'où l'identité proposée ci-dessus. La définition inductive des cumulants permet maintenant de conclure :
\begin{align*}k^{\id}(\varSigma_{i_{1}},\ldots,\varSigma_{i_{r}})&= \varSigma_{i_{1}}\cdots\varSigma_{i_{r}}\,\,-\!\! \sum_{\pi \neq \{\lle 1,r\rre\}}k^{\id}(\varSigma_{i \in \pi_{1}})\bullet k^{\id}(\varSigma_{i \in \pi_{2}})\bullet \cdots \bullet k^{\id}(\varSigma_{i \in \pi_{l}}) \\
&=\!\!\sum_{\forall j,\,\,A_{j} \in \arra(n,i_{j})}\!\!\!\!C(A_{1})\cdots C(A_{r})\,\,-\!\!\sum_{\pi \neq \{\lle 1,r\rre\}}\bigg( \sum_{\substack{\forall j,\,\,A_{j} \in \arra(n,i_{j})\\ \pi(A_{1},\ldots,A_{r})=\pi }} \!\!\!\!C(A_{1})\cdots C(A_{r})\bigg)\\
&=\!\!\sum_{\substack{\forall j,\,\,A_{j} \in \arra(n,i_{j})\\ \pi(A_{1},\ldots,A_{r})=\{\lle 1,r \rre\}}} \!\!\!\!C(A_{1})\cdots C(A_{r})
\end{align*}
d'où le résultat par récurrence sur $r$.
\end{proof}\bigskip

\begin{proof}[Démonstration du lemme \ref{magnitude}]
On considère le cumulant identité de caractères centraux de cycles $\varSigma_{i_{1}},\ldots,\varSigma_{i_{r}}$. Pour toute famille d'arrangements $A_{1},\ldots,A_{r}$ apparaissant dans la somme définissant ce cumulant, on a d'après le lemme précédent $$\pi(A_{1},\ldots,A_{r})=\{\lle 1,r\rre\}\,,$$ donc le cardinal de $|A_{1}\cup A_{2}\cup \cdots \cup A_{r}|$ est plus petit que $i_{1}+\cdots+i_{r}-(r-1)$, car il y a au moins $r-1$ égalités entre les éléments de ces arrangements. Le cumulant $k^{\id}(\varSigma_{i_{1}},\ldots,\varSigma_{i_{r}})$ est donc une combinaison linéaire de caractères centraux $\varSigma_{\lambda}$ avec $|\lambda|\leq i_{1}+\cdots+i_{r}-(r-1)$, et ainsi,
$$\deg\big(k^{\id}(\varSigma_{i_{1}},\ldots,\varSigma_{i_{r}})\big)\leq i_{1}+\cdots+i_{r}-(r-1)\,.$$
Notons que ce résultat est différent du théorème 4.3 de \cite{Sni06b} ; en effet, on utilise une autre graduation sur l'algèbre d'observables $\obs$. Ceci étant, on peut maintenant contrôler le degré des observables apparaissant dans les cumulants disjoints de la somme 
$$k(\varSigma_{i_{1}},\ldots,\varSigma_{i_{r}})=\sum_{\pi \in \mathfrak{Q}(\lle 1,r\rre)} k^{\bullet}\big(k^{\id}(\varSigma_{i_{j \in \pi_{1}}}),k^{\id}(\varSigma_{i_{j \in \pi_{2}}}),\ldots,k^{\id}(\varSigma_{i_{j \in \pi_{l}}})\big)\,,$$
et leur somme est toujours inférieure à $i_{1}+i_{2}+\cdots+i_{r}-r+\ell(\pi)$. D'après le lemme \ref{magnitudedj}, les cumulants disjoints sont donc tous des $$O(n^{i_{1}+i_{2}+\cdots+i_{r}-r+\ell(\pi)-(\ell(\pi)-1)})=O(n^{i_{1}+i_{2}+\cdots+i_{r}-r+1})\,,$$
ce qui achève la preuve du lemme \ref{magnitude}.
\end{proof}\bigskip\bigskip

Nous sommes maintenant en mesure de démontrer le $q$-analogue de la première partie du théorème \ref{secondasymptoticplancherel}, c'est-à-dire la déviation gaussienne des observables $\varSigma_{k}$ correctement renormalisées. Dans ce qui suit, si $k$ est un entier plus grand que $1$, nous noterons 
$$W_{k,n,q}=\sqrt{n}\left(\frac{p_{k}(\lambda) - \esper[p_{k}]}{n^{k}}\right) \qquad;\qquad Z_{k,n,q}=\sqrt{n}\left(\frac{\varSigma_{k}(\lambda) - \esper[\varSigma_{k}]}{n^{k}}\right)$$
les déviations renormalisées des moments de Frobenius et des caractères centraux, $\lambda$ étant un diagramme aléatoire sous la $q$-mesure de Plancherel $M_{n,q}$.
\begin{proposition}[Déviation gaussienne des observables sous les $q$-mesures de Plancherel]\label{devgaussianobs}
Les processus $(W_{k,n,q})_{k \geq 1}$ et $(Z_{k,n,q})_{k \geq 1}$ convergent en lois fini-dimensionnelles vers le même processus gaussien centré $(Z_{k,\infty,q})_{k \geq 1}$, les covariances de ce processus limite étant :
$$\mathrm{cov}(Z_{l,\infty,q},Z_{m,\infty,q})=lm\,(1-q)^{l+m}\,\left(\frac{1}{1-q^{l+m-1,1}}-\frac{1}{1-q^{l,m}}\right).$$
\end{proposition}\bigskip

\begin{proof}
Si $l$ et $m$ sont deux entiers positifs, alors le cumulant disjoint $k^{\bullet}(\varSigma_{l},\varSigma_{m})$ est $\esper[\varSigma_{l,m}]-\esper[\varSigma_{l}]\,\esper[\varSigma_{m}]$, c'est-à-dire :
$$\frac{(1-q)^{l+m}}{1-q^{l,m}}\left(n^{\downarrow l+m}-n^{\downarrow l}\,n^{\downarrow m}\right) =-lm\,\frac{(1-q)^{l+m}}{1-q^{l,m}}\,n^{l+m-1}+O(n^{l+m-2})\,.$$
En effet, les factorielles décroissantes ont pour développement 
$$n^{\downarrow k}=n^{k}-\frac{k(k-1)}{2}\,n^{k-1}+O(n^{k-2}),.$$
Calculons maintenant le cumulant standard $k(\varSigma_{l},\varSigma_{m})$, en utilisant le développement des produits de caractères centraux de cycles donné plus haut :
\begin{align*}
k(\varSigma_{l},\varSigma_{m})&=\esper[\varSigma_{l}\varSigma_{m}]-\esper[\varSigma_{l}]\,\esper[\varSigma_{m}]=\esper[\varSigma_{l}\varSigma_{m}-\varSigma_{l,m}]+k^{\bullet}(\varSigma_{l},\varSigma_{m})\\
&=\esper[lm\,\varSigma_{l+m-1}]-lm\,\frac{(1-q)^{l+m}}{1-q^{l,m}}\,n^{l+m-1}+O(n^{l+m-2})\\
&=lm\,(1-q)^{l+m}\,\left(\frac{1}{1-q^{l+m-1,1}}-\frac{1}{1-q^{l,m}}\right)\,n^{l+m-1}+O(n^{l+m-2})\,.
\end{align*}
Par construction, les observables $Z_{k,n,q}$ ont leur premier cumulant nul, et par multilinéarité et invariance par translation des cumulants d'ordre $2$, 
$$k(Z_{l,n,q},Z_{m,n,q})=lm\,(1-q)^{l+m}\,\left(\frac{1}{1-q^{l+m-1,1}}-\frac{1}{1-q^{l,m}}\right)+O(n^{-1})\,.$$
Compte tenu du lemme \ref{magnitude}, les cumulants d'ordre $r\geq 3$ des variables $Z_{k,n,q}$ sont pour leur part des $O(n^{1-r/2})$ :
$$k(Z_{k_{1},n,q},\ldots,Z_{k_{r},n,q})=\frac{O(n^{k_{1}+k_{2}+\cdots+k_{r}-r+1})}{n^{k_{1}-1/2}\,\cdots\,n^{k_{r}-1/2}}=O(n^{1-r/2})\,.
$$
Ils tendent donc vers $0$, et la convergence en lois fini-dimensionnelles de $(Z_{k,n,q})_{k \geq 1}$ vers le processus gaussien $(Z_{k,\infty,q})_{k \geq 1}$ est ainsi établie. \bigskip

\noindent Pour les $W_{k,n,q}$, le lemme \ref{magnitude} montre que les cumulants d'ordre $r\geq 3$ de ces variables sont encore des $O(n^{1-r/2})$. Le processus $(W_{k,n,q})_{k \geq 1}$ converge donc en lois fini-dimensionnelles vers un processus gaussien centré, et il suffit maintenant de montrer que $k(p_{l},p_{m})$ et $k(\varSigma_{l},\varSigma_{m})$ ont le même terme de plus haut degré. Or, chaque $\varSigma_{k}$ est un polynôme rationnel en les $p_{l \leq k}$, avec pour terme de plus haut degré $p_{k}$. Par multilinéarité des cumulants, le cumulant $k(\varSigma_{l},\varSigma_{m})$ est donc 
$$k(p_{l},p_{m})+(\text{cumulants d'observables }a_{i},b_{j}\text{ telles que }\deg a_{i}+\deg b_{j}\leq l+m-1)\,.$$
Le lemme \ref{magnitude} permet de conclure que $k(\varSigma_{l},\varSigma_{m})$ et $k(p_{l},p_{m})$, qui sont tous deux d'ordre de grandeur $O(n^{l+m-1})$, sont équivalents ; ceci achève la preuve du théorème.
\end{proof}\bigskip

\begin{proof}[Démonstration du théorème \ref{secondasymptoticqplancherel}]
Comme précédemment, on note $X_{\lambda}$ la mesure de probabilité sur $[-1,1]$ associée à une partition $\lambda$, et $X_{\infty,q}$ la loi limite de ces mesures de probabilité sous les $q$-mesures de Plancherel. Si $f$ est une fonction de classe $\mathscr{C}^{1}$ sur $[-1,1]$, on note
$$D_{n,q}(f)=\sqrt{n}\,(X_{\lambda}(f)-X_{\infty,q}(f))$$
la variable aléatoire correspondant à la déviation renormalisée de $X_{\lambda}(f)$ ; ici, pour toute probabilité $\mu$, $\mu(f)=\int_{-1}^{1} f(s)\,d\mu(s)$. Calculons ces déviations lorsque $f(x)=x^{l}$ :
$$D_{n,q}(x^{l})=\sqrt{n}\,\left(\frac{p_{l+1}(\lambda)}{n^{l+1}}-\frac{(1-q)^{l+1}}{1-q^{l+1}}\right)=W_{l+1,n,q} \,.$$
Par conséquent, étant donnés deux entiers strictement positifs $l$ et $m$, les déviations $D_{n,q}(x^{l})$ et $D_{n,q}(x^{m})$ convergent vers des lois gaussiennes centrées de covariance
\begin{align*}&(l+1)(m+1)\,\left(\frac{(1-q)^{l+m+1}}{1-q^{l+m+1}}-\frac{(1-q)^{l+m+2}}{1-q^{l+1,m+1}}\right)\\
&=X_{\infty,q}\big((l+1)(m+1)\,x^{l}\,x^{m}\big)-X_{\infty,q}\big((l+1)\,x^{l}\big)\,X_{\infty,q}\big((m+1)\,x^{m}\big)\\
&=X_{\infty,q}\big((x^{l+1})'(x^{m+1})'\big)-X_{\infty,q}\big((x^{l+1})'\big)\,X_{\infty,q}\big((x^{m+1})'\big)\,.\end{align*}
Dans $\mathscr{C}^{1}([-1,1])$, les fonctions polynômes sont denses pour la topologie $\mathscr{C}^{1}$. En utilisant la multilinéarité de la covariance et cette densité, on conclut que pour toutes fonctions $f_{1},\ldots,f_{r}$ de classe $\mathscr{C}^{1}$ sur $[-1,1]$, le vecteur des déviations $(D_{n,q}(f_{i}))_{i \in \lle 1,r\rre}$ converge vers un vecteur gaussien $(D_{\infty,q}(f_{i}))_{i \in \lle 1,r\rre}$ de matrice de covariance 
$$k(D_{\infty,q}(f_{i}),D_{\infty,q}(f_{j}))=X_{\infty,q}\big((xf_{i})'(xf_{j})'\big)-X_{\infty,q}\big((xf_{i})'\big)\,X_{\infty,q}\big((xf_{j})'\big)\,.$$
Si $i \geq 1$, notons $f_{i}$ une fonction positive de classe $\mathscr{C}^{1}$ sur $[-1,1]$, telle que $f_{i}(x)=1$ sur un voisinage $V_{i}$ de $(1-q)\,q^{i-1}$, et telle que $f_{i}(x)=0$ en dehors d'un voisinage $W_{i} \subset V_{i}$, voir la figure \ref{testfunction}.
\figcapt{\psset{unit=1mm}\pspicture(-20,-8)(80,20)
\psline{->}(-20,0)(80,0)
\psline{->}(0,-2)(0,20)
\pscurve[linecolor=blue](-15,0.2)(5,0.2)(20,0.2)(34,0.2)(39,1)(46,14)(50,15)(54,14)(61,1)(67,0.2)(75,0.2)
\psline(75,-1)(75,1)
\psline(50,-1)(50,1)
\psline(37,-8)(37,1)
\psline(63,-8)(63,1)
\psline{<->}(37,-7)(63,-7)
\rput(50,-9.5){$W_i$}
\rput(50,3){$q^{i-1}(1-q)$}
\endpspicture}{Fonction test $f_{i}$ pour la déviation de la mesure $X_{\lambda}$.\label{testfunction}}{Fonction test $f_{i}$ pour la déviation de la mesure $X_{\lambda}$}

\noindent On supposera également que tous les $W_{i}$ sont disjoints ; alors, $X_{\infty,q}(f_{i})=(1-q)\,q^{i-1}$ pour tout indice $i$. Comme $f_{i}'((1-q)\,q^{i-1})=0$, on a donc :
\begin{align*}
&(xf_{i})'((1-q)\,q^{i-1})=(f_{i})((1-q)\,q^{i-1})=1\\ 
&X_{\infty,q}\big((xf_{i})'^{2}\big)-X_{\infty,q}\big((xf_{i})'\big)^{2}=(1-q)\,q^{i-1}-(1-q)^{2}\,q^{2(i-1)}\,.
\end{align*}
D'autre part, si $i\neq j$, alors comme $W_{i}$ et $W_{j}$ sont disjoints, les produits de fonctions $f_{i}f_{j}$, $f_{i}'f_{j}$, $f_{i}f_{j}'$ et $f_{i}'f_{j}'$ sont tous nuls. On a donc :
\begin{align*}X_{\infty,q}\big((xf_i)'(xf_j)'\big)-X_{\infty,q}\big((xf_{i})'\big)\,X_{\infty,q}\big((xf_j)'\big)&=-X_{\infty,q}\big((xf_i)'\big)\,X_{\infty,q}\big((xf_{j})'\big) \,;\\
&=-(1-q)^2\,q^{i+j-2}\,.\end{align*}
Finalement, pour $n$ assez grand, en dehors d'un événement de probabilité arbitrairement petite, $a_{i}(\lambda)/n$ est dans $V_{i}$, $a_{i-1}(\lambda)/n$ est dans $V_{i-1}$ et $a_{i+1}(\lambda)/n$ est dans $V_{i+1}$. Donc, avec probabilité très proche de $1$, $a_{i}(\lambda)/n$ est la seule coordonnée de Frobenius de $\lambda$ dans le support de $f_{i}$, et $X_{\lambda}(f_{i})=a_{i}(\lambda)/n$. Ainsi, en dehors d'un événement de petite probabilité,
$$\sqrt{n}\left(\frac{a_{i}(\lambda)}{n}-(1-q)\,q^{i-1}\right)=\sqrt{n}\left(X_{\lambda}(f_{i})-X_{\infty,q}(f_{i})\right)=D_{n,q}(f_{i})\,,$$
donc les fonctions de répartition de ces deux variables aléatoires ont le même comportement asymptotique. Comme $a_{i}(\lambda)=\lambda_{i}(1+o(1))$, ceci achève la preuve du théorème.
\end{proof}\bigskip

Pour conclure cette section, mentionnons une analogie frappante entre les théorèmes \ref{firstasymptoticqplancherel} et \ref{secondasymptoticqplancherel}, et des résultats dus à A. Borodin concernant la forme de Jordan\footnote{Cette analogie fournit ce qu'on peut considérer être un modèle matriciel des $q$-mesures de Plancherel (au sens de la correspondance évoquée dans la première partie du mémoire).} des matrices triangulaires unipotentes sur $\For_{q}$ (voir \cite{Bor99} ; cette analogie nous a été signalée par Borodin lui-même). Dans ce qui suit, $\mathrm{GU}(n,\For_{q})$ désigne l'ensemble des matrices triangulaires supérieures de taille $n\times n$, à coefficients dans $\For_{q}$ et avec des $1$ sur toute la diagonale. C'est un sous-groupe unipotent maximal dans $\GL(n,\For_{q})$, et son cardinal est 
$$\card \mathrm{GU}(n,\For_{q})=q^{\frac{n(n-1)}{2}}\,.$$
On munit le groupe $\mathrm{GU}(n,\For_{q})$ de la mesure de probabilité uniforme $m_{n,q}$. Toute matrice $M \in \mathrm{GU}(n,\For_{q})$ est dans une classe de conjugaison du type $\bbmu=\{X-1:\mu\}$, où $\mu$ est une partition de taille $n$ ; on notera $\mu=\mu(M)$ cette \textbf{forme de Jordan}. Alors :
\begin{theorem}[Borodin, \cite{Bor99}]\label{borodinasymptotic}
Lorsque $n$ tend vers l'infini, pour tout indice $i \geq 1$, $\mu_{i}(M)/n$ converge en probabilité sous $m_{n,q}$ vers $q^{-i}(q-1)$. De plus, la déviation des lignes de $\mu(M)$ est gaussienne, et les covariances limites des variables
$$T_{i,n,q}=\sqrt{n}\left(\frac{\mu_{i}(M)}{n} - q^{-i}(q-1)\right)$$
sont $\mathrm{cov}(T_{i,\infty,q},T_{j,\infty,q})=\delta_{ij}(q-1)\,q^{-i}-(q-1)^{2}\,q^{-i-j}$. On retrouve donc les mêmes résultats que pour les colonnes d'une partition sous la $(q>1)$-mesure de Plancherel.
\end{theorem}
\noindent La preuve du théorème de Borodin repose sur une analyse approfondie d'un analogue du $q$-processus de Plancherel. Plus précisément, si $M$ est une matrice  unipotente (aléatoire) dans $\mathrm{GU}(n,\For_{q})$, alors on peut lui associer une matrice unipotente aléatoire $M'$ dans $\mathrm{GU}(n+1,\For_{q})$ en rajoutant une colonne dont les $n$ premières coordonnées sont tirées aléatoirement dans $(\For_{q})^{n}$ :
$$M'=\begin{pmatrix}
&   & &* \\
&M& &\vdots \\
&  &  &*\\
&0&  &1
\end{pmatrix}.$$
On dispose ainsi d'un processus markovien dont les marginales sont les lois $m_{n,q}$, et l'image de ce processus par $M \mapsto \mu(M)$ est un processus markovien sur le graphe de Young $\ym$. De plus, on peut évaluer assez précisément les probabilités de transition, et montrer qu'asymptotiquement les lignes de $\mu(M)$ se comportent comme des sommes de variables indépendantes ; ceci mène au théorème \ref{borodinasymptotic} via la loi des grands nombres et le théorème central limite <<~usuels~>>.\bigskip\bigskip

Ceci étant, dans une approche naïve de la $q$-mesure de Plancherel, il peut être tentant d'utiliser des méthodes semblables pour étudier l'asymptotique du $q$-processus de Plancherel, c'est-à-dire évaluer avec précision les $q$-probabilités de transition (\emph{cf.} \S\ref{qplancherelprocess})
$$p_{q}(\lambda,\Lambda)=q^{b(\Lambda)-b(\lambda)}\,\frac{\prod_{(i,j) \in \lambda}\{h(i,j)\}_{q}}{\prod_{(i,j) \in \Lambda}\{h(i,j)\}_{q}}\,,$$
et montrer que les lignes (ou colonnes selon la valeur de $q$) d'une partition se comportent asymptotiquement comme des sommes de variables indépendantes. Incidemment, c'est la première chose que nous avons tenté de faire, et on peut effectivement traiter de cette fa\c con le cas de la première ligne (ou colonne selon la valeur de $q$) et obtenir une loi des grands nombres $\lambda_{1}/n \to (1-q)$. Malheureusement, les autres lignes ou colonnes sont beaucoup plus difficiles à étudier par cette voie. D'autre part, nous ne sommes pas parvenus à établir un théorème central limite en utilisant cette approche (ne serait-ce que pour la première ligne ou colonne).\bigskip

\section{Concentration gaussienne des $q$-caractères}\label{qcharacter}
Pour conclure notre étude asymptotique des $q$-mesures de Plancherel, nous revenons sur le phénomène de concentration gaussienne des $q$-caractères autour de la trace canonique des algèbres d'Hecke. Dans le cas du groupe symétrique, rappelons que sous la mesure de Plancherel standard, l'ordre de grandeur du caractère $\chi^{\lambda}((1,2,\ldots,k\geq 2))$ est $n^{-k/2}$, et la variable renormalisée $$n^{k/2}\,\chi^{\lambda}((1,2,\ldots,k))$$ converge en loi vers un gaussienne de variance $k$ (voir le théorème \ref{secondasymptoticplancherel}) ; de plus, les caractères des cycles sont asymptotiquement indépendants. Le théorème \ref{thirdasymptoticqplancherel} énoncé ci-dessous donne l'analogue de ce résultat pour les $q$-caractères sous la $q$-mesure de Plancherel :

\begin{theorem}[Asymptotique des $q$-caractères, \cite{Mel10a}]\label{thirdasymptoticqplancherel}
Pour $k \geq 2$, notons $C_{k}$ le <<~cycle~>> $T_{1}T_{2}\cdots T_{k-1}$ dans les algèbres d'Hecke $\IH_{q}(\sym_{n})$ avec $n\geq k$. Si $\lambda$ est tiré au hasard selon la $q$-mesure de Plancherel, alors le $q$-caractère renormalisé
$$S_{k,n,q}=\sqrt{n}\,\,\chi^{\lambda}(q,C_{k})$$
converge en loi vers une gaussienne. Les covariances de ces gaussiennes limites sont :
$$\mathrm{cov}(S_{l,n,q},S_{m,n,q})=(q-q^{2})^{l+m-3}\,(1-q^{2})\,\frac{\{l-1\}_q\,\{m-1\}_{q}}{\{l+m-1\}_{q}\,\{l+m-2\}_{q}\,\{l+m-3\}_{q}}\,.$$
\end{theorem}
\noindent Le seul point nouveau dans le théorème \ref{thirdasymptoticqplancherel} est en réalité le calcul des covariances limites des $q$-caractères. En effet, compte tenu des formules de changement de base entre les symboles $\varSigma_{\mu}$ et $\varSigma_{\mu,q}$, le degré d'une observable $\varSigma_{k,q}$ est $k$ ; donc, les cumulants d'ordre $r$ supérieur à $3$
$$k\left(n^{1/2-k_{1}}\varSigma_{k_{1},q},n^{1/2-k_{2}}\varSigma_{k_{2},q},\ldots,n^{1/2-k_{r}}\varSigma_{k_{r},q}\right)$$ tendent tous vers $0$ en vertu du lemme \ref{magnitude} (c'est la même preuve que pour les observables $\varSigma_{k}$). Ceci implique bien le caractère gaussien de 
$$n^{1/2-k}\,\varSigma_{k,q}(\lambda) \sim n^{1/2}\,\chi^{\lambda}(q,C_{k})\,,$$
et il reste donc à démontrer que $k(n^{1/2-l}\varSigma_{l,q},n^{1/2-m}\varSigma_{m,q})$ converge vers la quantité indiquée dans l'énoncé du théorème \ref{thirdasymptoticqplancherel}.\bigskip
\bigskip

Pour commencer, nous allons donner une généralisation de la proposition \ref{devgaussianobs} en examinant les déviations d'observables
$$Z_{\mu,n,q}=\sqrt{n}\left(\frac{\varSigma_{\mu}(\lambda)-\esper[\varSigma_{\mu}]}{n^{|\mu|}}\right)\,.$$
Comme pour les caractères centraux des cycles, on décompose le cumulant simple $k(\varSigma_{\mu},\varSigma_{\rho})=\esper[\varSigma_{\mu}\varSigma_{\rho}]-\esper[\varSigma_{\mu}]\,\esper[\varSigma_{\rho}]$ en
$$k(\varSigma_{\mu},\varSigma_{\rho})=\esper[\varSigma_{\mu}\varSigma_{\rho}-\varSigma_{\mu\sqcup\rho}]+k^{\bullet}(\varSigma_{\mu},\varSigma_{\rho})\,.$$
Le cumulant disjoint se calcule exactement comme pour les cycles :
\begin{align*} k^{\bullet}(\varSigma_{\mu},\varSigma_{\rho})&=\frac{(1-q)^{|\mu|+|\rho|}}{1-q^{\mu \sqcup \rho}}\left(n^{\downarrow|\mu|+|\rho|}-n^{\downarrow |\mu|}\,n^{\downarrow |\rho|}\right)\\
&=-|\mu|\,|\rho|\,\frac{(1-q)^{|\mu|+|\rho|}}{1-q^{\mu \sqcup \rho}}\,n^{|\mu|+|\rho|-1}+O(n^{|\mu|+|\rho|-2})\,.\end{align*}
D'autre part, le produit $\varSigma_{\mu}\,\varSigma_{\rho}$ a pour terme de plus haut degré $\varSigma_{\mu\sqcup \rho}$, qui s'annule avec le $-\varSigma_{\mu \sqcup \rho}$. Les termes de degré immédiatement inférieur --- c'est-à-dire $|\mu|+|\rho|-1$ --- dans le produit $\varSigma_{\mu}\,\varSigma_{\rho}$ correspondent aux appariements partiels de taille $1$ entre les ensembles d'indices $I_{A}$ et $I_{B}$ correspondants aux partitions $\mu=(\mu_{1},\ldots,\mu_{r})$ et $\rho=(\rho_{1},\ldots,\rho_{s})$. Si $a_{ij}=b_{kl}$ est un tel appariement, alors le type de la permutation correspondante est 
$$(\mu_{1},\ldots,\mu_{i-1},\mu_{i+1},\ldots,\mu_{r},\rho_{1},\ldots,\rho_{k-1},\rho_{k+1},\ldots\rho_{s},\mu_{i}+\rho_{k}-1)\,,$$
étant entendu qu'on réordonne ensuite les parts de la partition. De plus, les parts $\mu_{i}$ et $\rho_{k}$ étant fixées, il y a exactement $\mu_{i}\rho_{k}$ appariements du type $a_{ij}=b_{kl}$. Ainsi :
\begin{lemma}[Second terme d'un produit de caractères centraux]
La composante homogène de degré $|\mu|+|\rho|-1$ dans un produit $\varSigma_{\mu}\,\varSigma_{\rho}$ est
$$\sum_{\substack{\mu_{i} \in \mu\\
\rho_{k} \in \rho}}\mu_{i}\rho_{k}\,\varSigma_{(\mu_{1},\ldots,\mu_{i-1},\mu_{i+1},\ldots,\mu_{r},\rho_{1},\ldots,\rho_{k-1},\rho_{k+1},\ldots\rho_{s},\mu_{i}+\rho_{k}-1)}\,.$$
\end{lemma}
\noindent L'espérance $\esper[\varSigma_{\mu}\varSigma_{\rho}-\varSigma_{\mu\sqcup\rho}]$ a donc pour terme dominant $n^{|\mu|+|\rho|-1}$ fois
$$\sum_{\substack{\mu_{i} \in \mu\\
\rho_{k} \in \rho}}\mu_{i}\rho_{k}\,\frac{(1-q)^{|\mu|+|\rho|-1}}{1-q^{\mu\sqcup\rho}}\,\frac{(1-q^{\mu_{i}})(1-q^{\rho_{k}})}{1-q^{\mu_{i}+\rho_{k}-1}}\,,$$
et en écrivant $|\mu|\,|\rho|=\sum_{i,k}\mu_{i}\,\rho_{k}$, on a donc :
\begin{align*}
k(\varSigma_{\mu},\varSigma_{\rho})&=n^{|\mu|+|\rho|-1}\,\frac{(1-q)^{|\mu|+|\rho|}}{1-q^{\mu\sqcup \rho}}\sum_{i,k} \mu_{i}\rho_{k}\,\left(\frac{1-q^{\mu_{i},\rho_{k}} }{1-q^{\mu_{i}+\rho_{k}-1,1}}-1\right)+O(n^{|\mu|+|\rho|-2})\\
&=n^{|\mu|+|\rho|-1}\,q\,\frac{(1-q)^{|\mu|+|\rho|}}{1-q^{\mu\sqcup \rho}}\sum_{i,k} \mu_{i}\rho_{k}\,\frac{1-q^{\mu_{i}-1,\rho_{k}-1} }{1-q^{\mu_{i}+\rho_{k}-1,1}}+O(n^{|\mu|+|\rho|-2})\,\,;\\
k(Z_{\mu,n,q},Z_{\rho,n,q})&=q\,\frac{(1-q)^{|\mu|+|\rho|}}{1-q^{\mu\sqcup \rho}}\sum_{i,k} \mu_{i}\rho_{k}\,\frac{1-q^{\mu_{i}-1,\rho_{k}-1} }{1-q^{\mu_{i}+\rho_{k}-1,1}}+O(n^{-1})\,.
\end{align*}
\bigskip

Dans tout ce qui suit, nous noterons $R_{k,n,q}=n^{1/2-k}\,\varSigma_{k,q}$ ; le comportement asymptotique de ces variables est trivialement le même que celui des $S_{k,n,q}$. Compte tenu des formules de changement de base du théorème \ref{quantobs}, et en utilisant la multilinéarité des cumulants, on a :
$$k(R_{l,n,q},R_{m,n,q})=\frac{1}{(q-1)^{2}}\sum_{\substack{\mu \in \ym_{l}\\\rho \in \ym_{m}}} \frac{\scal{h_{l}}{p_{\mu}}\scal{h_{m}}{p_{\rho}} }{z_{\mu}\,z_{\rho}}\,(q^{\mu\sqcup \rho}-1)\,k(Z_{\mu,n,q},Z_{\rho,n,q})\,.$$
Comme $h_{l}=s_{l}$ et $h_{m}=s_{m}$, les produits scalaires au numérateur sont égaux à $1$ par la formule de Frobenius. En utilisant la formule donnée précédemment pour les $k(Z_{\mu,n,q},Z_{\rho,n,q})$, on obtient donc 
$$ k(R_{l,n,q},R_{m,n,q})=q(1-q)^{l+m-3}\,\sum_{\substack{\mu \in \ym_{l}\\ \rho \in \ym_{m}}} \frac{(-1)^{\ell(\mu)+\ell(\rho)}}{z_{\mu}\,z_{\rho}}\,\sum_{\substack{i \in \mu\\ k \in \rho}} ik\,\frac{1-q^{i-1,k-1}}{1-q^{i+k-1}} +O(n^{-1})\,.$$ 
Dans la somme, isolons les termes correspondant aux partitions $\mu=(l),(l-1,1)$ et $\rho=(m),(m-1,1)$. La somme de ces quatre termes est :
\begin{align*}&\frac{1-q^{l-1,m-1}}{1-q^{l+m-1}}-\frac{1-q^{l-1,m-2}}{1-q^{l+m-2}}-\frac{1-q^{l-2,m-1}}{1-q^{l+m-2}}+\frac{1-q^{l-2,m-2}}{1-q^{l+m-3}}\\
&=q^{m-2}(1-q)\left(\frac{1-q^{l-1,l}}{1-q^{l+m-2,l+m-1}}-\frac{1-q^{l-2,l-1}}{1-q^{l+m-3,l+m-2}}\right)\\
&=\frac{q^{m-2}(1-q)(1-q^{l-1})}{1-q^{l+m-2}}\left(\frac{1-q^{l}}{1-q^{l+m-1}}-\frac{1-q^{l-2}}{1-q^{l+m-3}}\right)\\
&=\frac{q^{m-2}(1-q)(1-q^{l-1})}{1-q^{l+m-2}}\,\frac{q^{l-2}(1-q^{2})(1-q^{m-1})}{1-q^{l+m-1,l+m-3}}\\
&=q^{l+m-4}(1-q^{2})\,\frac{\{l-1\}_{q}\{m-1\}_{q}}{\{l+m-1\}_{q}\{l+m-2\}_{q}\{l+m-3\}_{q}}\,.\end{align*}
En multipliant par le facteur $q(1-q)^{l+m-3}$, on retrouve exactement la formule du théorème \ref{thirdasymptoticqplancherel}. Il suffit donc de montrer que tous les termes de la somme double $\sum_{\mu,\rho}$ se simplifient, à l'exception des quatre termes correspondant aux paires $(l)(m)$, $(l)(m-1,1)$, $(l-1,1)(m)$ et $(l-1,1)(m-1,1)$. Dans ce qui suit, on note $C(\sigma)$ l'ensemble des cycles d'une permutation $\sigma$.
\begin{lemma}[Formule d'inversion de M\"obius pour les permutations et les partitions]\label{stupidlemma}
Soit $f$ une fonction des entiers, $m$ un entier plus grand que $2$. 
$$\sum_{\lambda \in \ym_{m}}\frac{(-1)^{\ell(\lambda)}}{z_{\lambda}}\left(\sum_{k \in \lambda}f(\lambda)\right)=\frac{1}{m!}\sum_{\sigma \in \sym_{m}} (-1)^{|C(\sigma)|}\left(\sum_{c \in C(\sigma)}f(|c|)\right)=\frac{f(m-1)}{m-1}-\frac{f(m)}{m}\,.$$
\end{lemma}
\begin{proof}[Preuve du lemme \ref{stupidlemma}]
La première égalité est obtenue en réunissant les permutations $\sigma$ selon leur type $\lambda=t(\sigma)$, car $\ell(\lambda)=|C(\sigma)|$ et $z_{\lambda}=n!/\card C_{\lambda}$. On démontre maintenant la seconde égalité par récurrence forte sur $m$, en notant
$$F(m)= \frac{1}{m!}\sum_{\sigma \in \sym_{m}} (-1)^{|C(\sigma)|}\left(\sum_{c \in C(\sigma)}f(|c|)\right).$$
L'identité est vraie pour $m=2$. Si elle est vraie jusqu'au rang $m$, alors au rang $m+1$, si $\sigma \in \sym_{m+1}$ :
\begin{enumerate}
\item Si $\sigma$ est un $m+1$-cycle, le terme correspondant dans la somme est $-f(m+1)$, et il y a $m!$ termes de ce type, d'où une contribution 
$$A_{m+1}=\frac{1}{m+1!}(-m!\,f(m+1))=-\frac{f(m+1)}{m+1}\,.$$
\item Si $\sigma$ est un $m$-cycle qui bouge l'entier $1$, alors le terme correspondant dans la somme est $f(m)+f(1)$, et il y a $m!$ termes de ce type, d'où une contribution 
$$A_{m}=\frac{1}{m+1!}(m!\,f(m)+m!\,f(1))=\frac{f(m)}{m+1}+\frac{f(1)}{m+1}\,.$$
\item Sinon, l'entier $1$ est dans un $k$-cycle avec $k \in \lle 1,m-1\rre$, et le reste de la permutation peut être vu comme un élément de $\sym_{m+1-k}$. Il y a $m(m-1)\cdots(m-k+2)$ cycles de longueur $k$ possibles permutant l'entier $1$, d'où une contribution 
\begin{align*}A_{k}&=\frac{m^{\downarrow k-1}}{(m+1)!}\,\sum_{\sigma\in \sym_{m+1-k}}(-1)^{1+|C(\sigma)|}\,\left(f(k)+\sum_{c \in C(\sigma)}f(c)\right)\\
&=\frac{1}{m+1}\,\left(-F(m+1-k)-\frac{f(k)}{m+1-k!}\sum_{\sigma \in \sym_{m+1-k}}(-1)^{|C(\sigma)|}\right)\\
&=\frac{1}{m+1}\,\left(-\frac{f(m-k)}{m-k}+\frac{f(m+1-k)}{m+1-k}\right)\,.
\end{align*}
En effet, l'hypothèse de récurrence s'applique avec $m+1-k \in \lle 2,m\rre$, et comme l'indice $m+1-k$ est plus grand que $2$, la somme des signes des permutations est nulle. 
\end{enumerate}
La somme des contributions $A_{k}$ pour $k$ compris entre $1$ et $m-1$ est téléscopique et égale à 
$$\frac{1}{m+1}\,\left(-f(1)+\frac{f(m)}{m}\right)\,,$$
et en ajoutant les contributions $A_{m}$ et $A_{m+1}$, on obtient bien $\frac{f(m)}{m}-\frac{f(m+1)}{m+1}$, d'où le résultat au rang $m+1$.
\end{proof}\bigskip

\begin{proof}[Fin de la démonstration du théorème \ref{thirdasymptoticqplancherel}]
Fixons un entier $i$, et considérons la fonction $f(k)=k\,(1-q^{i-1,k-1})/(1-q^{i+k-1})$. D'après le lemme \ref{stupidlemma}, 
$$\sum_{\rho \in \ym_{m}} \frac{(-1)^{\ell(\rho)}}{z_{\rho}}\,\sum_{k \in \rho}ik\,\frac{1-q^{i-1,k-1}}{1-q^{i+k-1}}=-i\left(\frac{1-q^{i-1,m-1}}{1-q^{i+m-1}} - \frac{1-q^{i-1,m-2}}{1-q^{i+m-2}}\right)\,.$$
Par conséquent, la somme double mise en jeu dans le calcul de $k(R_{l,n,q},R_{m,n,q})$ se ramène à deux sommes simples :
$$\sum_{\mu,\rho}=\sum_{\mu \in \ym_{l}}\frac{(-1)^{\ell(\mu)}}{z_{\mu}}\sum_{i \in \mu}i\,\frac{1-q^{i-1,m-2}}{1-q^{i+m-2}}-\sum_{\mu \in \ym_{l}}\frac{(-1)^{\ell(\mu)}}{z_{\mu}}\sum_{i \in \mu}i\,\frac{1-q^{i-1,m-1}}{1-q^{i+m-1}}\,.$$
En réappliquant le même lemme aux fonctions $f_{2}(i)=i\,(1-q^{i-1,m-2})/(1-q^{i+m-2})$ et $f_{1}(i)=i\,(1-q^{i-1,m-1})/(1-q^{i+m-1})$, on obtient :
$$\sum_{\mu,\rho}=\frac{1-q^{l-1,m-1}}{1-q^{l+m-1}}-\frac{1-q^{l-1,m-2}}{1-q^{l+m-2}}-\frac{1-q^{l-2,m-1}}{1-q^{l+m-2}}+\frac{1-q^{l-2,m-2}}{1-q^{l+m-3}}\,.$$
C'est ce qu'il fallait démontrer.
\end{proof}
\noindent Ces calculs achèvent notre étude asymptotique de la $q$-mesure de Plancherel des algèbres d'Hecke de type A. Notons que le théorème \ref{thirdasymptoticqplancherel} reste vrai pour $q>1$, avec exactement la même formule pour les covariances ; en effet, l'hypothèse $q<1$ des paragraphes précédents n'a pas servi dans la preuve. Notons toutefois une différence entre ces deux cas : lorsque $q$ est strictement inférieur à $1$, deux $q$-caractères $\chi^{\lambda}(q,C_{k})$ et $\chi^{\lambda}(q,C_{l})$ sont toujours asymptotiquement positivement corrélés, tandis que pour $q>1$, le signe de la covariance limite dépend de la parité de $k+l$.

\chapter{Asymptotique des mesures d'induction parabolique}\label{arikikoike}

Dans ce chapitre, nous étudions les mesures de Plancherel associées à d'autres algèbres d'Hecke, et en particulier la mesure de Plancherel de l'algèbre d'Hecke de type B. Cette algèbre est une déformation $\IH_{q}(\wsym_{n})$ de l'algèbre du \textbf{groupe hyperoctahédral}, et c'est le commutant de l'action du \textbf{groupe symplectique} $\Sp(2n,\For_{q})$ sur sa variété de drapeaux, voir la section \ref{symplectic}. Les représentations irréductibles de $\IH_{q}(\wsym_{n})$ sont indexées par les paires de partitions (ou \textbf{bipartitions}) de taille $n$, et sous la B-$q$-mesure de Plancherel pour les bipartitions, on dispose d'un résultat asymptotique tout à fait analogue aux théorèmes \ref{firstasymptoticqplancherel} et \ref{secondasymptoticqplancherel}, voir la section \ref{basymptotic}, et en particulier le théorème \ref{asymptoticbqplancherel}. La preuve du théorème \ref{asymptoticbqplancherel} est adaptée des raisonnements des chapitres précédents, et l'on considérera donc des observables de bidiagrammes, notamment des caractères et des $q$-caractères renormalisés des groupes hyperoctahédraux.\bigskip\bigskip

Dans la section \ref{heckebmix}, on présente une conjecture concernant la répartition des parts entre les deux parties d'une bipartition aléatoire $(\lambda^{(1)},\lambda^{(2)})$ tirée suivant la B-$q$-mesure de Plancherel. Nos résultats asymptotiques pour la B-$q$-mesure de Plancherel sont malheureusement partiels, car la combinatoire des caractères de l'algèbre d'Hecke de type B n'est pas entièrement comprise. Plus précisément, dans le cas du groupe symétrique $\sym_{n}$, un outil essentiel des chapitres précédents était la formule de Ram \ref{ramformula}, qui relie les caractères irréductibles du groupe aux $q$-caractères de l'algèbre d'Hecke. Dans le cas du groupe hyperoctahédral $\wsym_{n}$, et plus généralement pour un produit en couronne $(\Z/r\Z)\wr \sym_{n}$, un résultat analogue a été démontré par T. Shoji (\emph{cf.} \cite{SS99,Sho00}), voir le théorème \ref{shojiformula}. Sa formule rentre dans le cadre de la dualité de Schur-Weyl entre groupes quantiques et algèbres  d'Hecke cyclotomiques --- les \textbf{algèbres d'Ariki-Koike}. Le problème se ramène alors à comprendre les liens combinatoires entre deux présentations de l'algèbre $\IH_{q}(\wsym_{n})$ : \vspace{2mm}
\begin{itemize}
\item[-] la présentation comme algèbre d'Hecke du groupe de Weyl d'un groupe algébrique, c'est-à-dire, par la théorie d'Iwahori, comme commutant de l'action d'un groupe de Chevalley sur sa variété de drapeaux ;\vspace{2mm}
\item[-] et la présentation comme algèbre d'Ariki-Koike, c'est-à-dire le commutant de l'action sur un produit tensoriel d'espaces $(\C^{m})^{\otimes n}$ du groupe quantique associé à une sous-algèbre de Lévi
$$\mathfrak{g}=\liegl(m_{1})\oplus\liegl(m_{2})\oplus \cdots \oplus \liegl(m_{r}) \subset \liegl(m=m_{1}+\cdots+m_{r})\,.$$ \vspace{2mm}
\end{itemize}
Malheureusement, les deux présentations semblent très difficiles à relier, et en particulier, on ne dispose pas de $q$-formule de Frobenius pour la valeur d'un caractère irréductible de l'algèbre d'Hecke de type B en un élément correspondant à une permutation signée de longueur minimale dans sa classe de conjugaison. Ces problèmes sont exposés en détail dans la section \ref{shoji} ; cette discussion permettra aussi de comprendre l'origine de la $q$-formule de Frobenius-Schur \ref{qfrobeniusschur}, et préparera l'étude des mesures de Schur-Weyl dans le chapitre \ref{schurweylmeasure}. \bigskip
\bigskip

Finalement, on décrit dans la section \ref{lusztig} le cadre le plus général dans lequel il semble possible de généraliser les raisonnements des chapitres \ref{iwahori} et \ref{qplancherelmeasure} : c'est le cadre des algèbres d'Hecke généralisées, qui sont les commutants de l'action d'un groupe réductif de type Lie sur un module obtenu par induction parabolique à partir d'un caractère cuspidal d'un sous-groupe de Lévi rationnel (\emph{cf.} \cite{HL80,Lus84}). \`A titre d'exemple, nous détaillons le cas des modules (quelconques !) obtenus par induction parabolique à partir d'un tore scindé de $\GL(n,\For_{q})$.\bigskip

\section{Variété de drapeaux symplectique et algèbre d'Hecke de type B}\label{symplectic}

Si $m$ est un entier plus grand que $1$ et si $q$ est la puissance d'un nombre premier, on appelle \textbf{forme antisymétrique} sur $(\For_{q})^{m}$ une forme bilinéaire $b : (\For_{q})^{m}\otimes (\For_{q})^{m} \to \For_{q}$ telle que $b(x,x)=0$ pour tout vecteur $x$. Alors, pour tous vecteurs $x$ et $y$ dans $(\For_{q})^{m}$,
$$0=b(x+y,x+y)=b(x,x)+b(x,y)+b(y,x)+b(y,y)=b(x,y)+b(y,x)\,,$$
donc $b(x,y)=-b(y,x)$. La forme est dite \textbf{non dégénérée} si la matrice de $b$ (dans n'importe quelle base) est inversible, ou, car c'est équivalent, si pour tout vecteur $x$ non nul, la forme linéaire $b(x,\cdot)$ est non nulle. Dans ce cas, la dimension de l'espace est paire, et si $m=2n$, alors il existe une base --- dite \textbf{base hyperbolique} --- $(v_{1},w_{1},v_{2},w_{2},\ldots,v_{n},w_{n})$ telle que la matrice de $b$ s'écrive
$$J_{2n}=\begin{pmatrix}
0  & 1 &     &    &          &      & \\
-1 & 0 &     &    &          &     & \\
    &    & 0  & 1 &          &      & \\
    &    & -1 & 0 &          &     & \\
    &    &     &    &\ddots&     & \\
    &    &     &    &          & 0 & 1 \\
    &    &     &    &          & -1& 0
\end{pmatrix},$$
voir \cite[chapitre XV, \S8]{Lang93}. En particulier, toutes les formes bilinéaires antisymétriques non dégénérées sur un $\For_{q}$-espace vectoriel de dimension paire sont isométriques, et on peut sans perte de généralité supposer que $b$ est la forme donnée dans la base canonique de $(\For_{q})^{2n}$ par la matrice écrite ci-dessus. Dans ce contexte, le \textbf{groupe symplectique} $\Sp(2n,\For_{q})$ est le groupe des isométries de $b$, c'est-à-dire le groupe de matrices inversibles de taille $2n \times 2n$ telles que
$$b(gx,gy)=b(x,y)$$
pour tous vecteurs $x$ et $y$. Un isomorphisme $g \in \GL(2n,\For_{q})$ est symplectique si et seulement s'il envoie une base hyperbolique sur une autre base hyperbolique de l'espace ; comme dans le paragraphe \ref{jordanfrobenius}, le cardinal de $\Sp(2n,\For_{q})$ est donc le nombre de bases hyperboliques de $(\For_{q})^{2n}$. Pour construire une telle base, on choisit d'abord un vecteur $v_{1}$ parmi les $q^{2n}-1$ vecteurs non nuls. Le vecteur $w_{1}$ doit satisfaire l'identité
$$b(v_{1},w_{1})=1\,,$$
et comme $b(v_{1},\cdot)$ est une forme linéaire non nulle, il y a $q^{2n-1}$ solutions à cette équation. Par conséquent, il y a $q^{2n-1}(q^{2n}-1)$ choix possibles pour les deux premiers vecteurs $v_{1}$ et $w_{1}$, et les vecteurs suivants $v_{2},w_{2},v_{3},w_{3},\ldots$ doivent être choisis dans l'orthogonal de $\For_{q}[v_{1},w_{1}]$, qui est un $\For_{q}$-espace vectoriel de dimension $2(n-1)$ muni d'une forme bilinéaire alternée non dégénérée. Par récurrence sur $n$, on conclut que :
\begin{align*}\card \Sp(2n,\For_{q})&= q^{(2n-1)+(2n-3)+\cdots+1}\,(q^{2n}-1)\,(q^{2(n-1)}-1)\,\cdots\,(q^{2}-1)\\
&=q^{n^{2}}\,\prod_{i=1}^{n}\,(q^{2i}-1)\,.\end{align*}
\bigskip

Dans un $\For_{q}$-espace vectoriel $V$ de dimension paire $2n$ et muni d'une forme antisymétrique non dégénérée $b$, on appelle \textbf{sous-espace totalement anisotrope} un sous-espace $U \subset V$ tel que $b_{|U\otimes U}=0$. Si tel est le cas, alors $U \subset U^{\perp}$, et comme $$2n=\dim U + \dim U^{\perp},$$ 
la dimension de $U$ ne peut pas dépasser $n$. Cette borne est atteinte en considérant par exemple $\For_{q}[v_{1},\ldots,v_{n}]$, où $(v_{1},w_{1},\ldots,v_{n},w_{n})$ est une base hyperbolique de $V$. Un \textbf{drapeau totalement anisotrope} est une suite strictement croissante de sous-espaces totalement anisotropes de $V$, et un drapeau totalement anisotrope complet est une suite
$$\{0\}=U_{0} \subset U_{1} \subset U_{2}  \subset \cdots \subset U_{n} \subset V=(\For_{q})^{2n}$$
avec chaque $U_{i}$ sous-espace totalement anisotrope de dimension $i$. Par exemple, si $\mathcal{B}=(v_{1},w_{1},\ldots,v_{n},w_{n})$ est la base canonique de $(\For_{q})^{2n}$, et si $b$ est la forme antisymétrique de matrice $J_{2n}$ dans la base $\mathcal{B}$, alors les sous-espaces 
$$U_{i}=\For_{q}[v_{1},\ldots,v_{i}]$$ forment un drapeau totalement anisotrope complet.
Il n'est pas difficile de voir que le stabilisateur d'un tel drapeau (dans $\Sp(2n,\For_{q})$) est un groupe de cardinal $q^{n^{2}}(q-1)^{n}$. Nous noterons ce groupe $\BSp(2n,\For_{q})$ ; c'est un sous-groupe de Borel de $\Sp(2n,\For_{q})$, et il contient le tore maximal
$$\mathrm{TSp}(2n,\For_{q})=\big\{M=\mathrm{diag}(\lambda_{1},\lambda_{1}^{-1},\lambda_{2},\lambda_{2}^{-1},\ldots,\lambda_{n},\lambda_{n}^{-1}) \,\,|\,\,\forall i,\,\,\lambda_{i} \in (\For_{q})^{\times}\big\}\,.$$
La \textbf{variété de drapeaux symplectique} est le quotient $\Sp(2n,\For_{q})/\BSp(2n,\For_{q})$ ; elle paramètre les drapeaux complets totalement anisotropes dans $(\For_{q})^{2n}$, et a pour cardinal\footnote{Cette identité découle du résultat suivant: le polynôme de Poincaré $\sum_{w \in W} q^{\ell(w)}$ du groupe de Weyl (ici, $\wsym_{n}$) évalué en $q$ est toujours égal au cardinal de la variété de drapeaux (ici, $\Sp(2n,\For_{q})/\BSp(2n,\For_{q})$), et ce compte tenu de la décomposition de Bruhat et de l'identité $\card B^{F}wB^{F}=q^{\ell(w)}\,\card B^{F}$.}
$$\card \Sp(2n,\For_{q})/\BSp(2n,\For_{q})=\{2n!!\}_{q}=\prod_{i=1}^{n}\, \{2i\}_{q}\,.$$ Dans les trois premiers paragraphes de ce chapitre, on s'intéresse à la décomposition du $\Sp(2n,\For_{q})$-module $\C[\Sp(2n,\For_{q})/\BSp(2n,\For_{q})]$ en sous-modules irréductibles, et à la mesure de probabilité de Plancherel sous-jacente. 
\bigskip
\bigskip

La théorie d'Iwahori-Hecke du paragraphe \ref{hecke} s'applique, et l'action de $\Sp(2n,\For_{q})$ sur $\C[\Sp(2n,\For_{q})/\BSp(2n,\For_{q})]$ admet donc pour commutant l'algèbre d'Hecke du groupe de Weyl du groupe de Chevalley $\Sp(2n,\For_{q})$. Ce groupe de Weyl est le \textbf{groupe hyperoctahédral} $\wsym_{n}$, qu'on peut aussi voir comme groupe des permutations signées de taille $n$. En effet, le normalisateur du tore $\mathrm{TSp}(2n,\For_{q})$ est le groupe des matrices de taille $2n \times 2n$ qui s'écrivent sous la forme
$$\begin{pmatrix}                          &M_{1\sigma(1)}& &\\
                             M_{2\sigma(2)}&                         & &\\
                                                      &                         & &\\
                                                      &                         &M_{n\sigma(n)}& 
\end{pmatrix},$$
où $\sigma$ est une permutation de taille $n$ indiquant la position des matrices $M_{i\sigma(i)}$, et où chaque $M_{i\sigma(i)}$ est un bloc de taille $2\times 2$ qui est de l'une des formes suivantes :
$$M_{i\sigma(i)}=\begin{pmatrix} \lambda_{i} & 0\\ 0 & \lambda_{i}^{-1}\end{pmatrix} \quad\text{ou}\quad\begin{pmatrix} 0 & \lambda_{i} \\  \lambda_{i}^{-1} & 0\end{pmatrix}.$$
Le quotient $W=N/T$ s'identifie donc au produit en couronne $(\Z/2\Z) \wr\sym_{n}$, c'est-à-dire, le produit semi-direct $(\Z/2\Z)^{n} \rtimes \sym_{n}$, le groupe symétrique $\sym_{n}$ agissant (à droite) sur les $n$-uplets $(\eps_{1},\ldots,\eps_{n}) \in (\Z/2\Z)^{n}$ par permutation des coordonnées. Nous noterons $\wsym_{n}=(\Z/2\Z)\wr \sym_{n}$ ce groupe ; 
les éléments du groupe sont les couples $((\eps_{1},\ldots,\eps_{n}) \in \{\pm\}^{n},\sigma \in \sym_{n})$, avec la loi produit :
$$\big((\eps_{1},\ldots,\eps_{n}),\sigma\big)\times\big((\theta_{1},\ldots,\theta_{n}),\tau\big)=\big((\eps_{\tau(1)}\theta_{1},\ldots,\eps_{\tau(n)}\theta_{n}),\sigma\tau\big)\,.$$
On peut donc voir les éléments de $\wsym_{n}$ comme des \textbf{permutations signées}, c'est-à-dire des permutations de l'ensemble $\lle -n,-1 \rre\sqcup \lle 1,n\rre$ telles que $\sigma(-i)=-\sigma(i)$ pour tout $i$. En effet, si $g=((\pm_{i})_{i \in \lle 1,n\rre},\sigma)$ est un élément de $\wsym_{n}$, alors son action est donnée par :
$$g(i)=\pm_{i}\,\sigma(i)\qquad;\qquad g(-i)=\mp_{i}\,\sigma(i)\,.$$\bigskip\medskip

Le cardinal de $\wsym_{n}$ est $2^{n}\,n!$ ; plus généralement, le cardinal d'un produit en couronne $(\Z/r\Z)^{n} \wr \sym_{n}$ est $r^{n}\,n!$. Déterminons les classes de conjugaison de $\wsym_{n}$. Si $g=((\eps_{i})_{i},\sigma)$ et $h=((\theta_{i})_{i}, \tau)$ sont deux permutations signées, alors
$$\widetilde{g}=h^{-1}gh=\big((\theta_{\tau^{-1}(i)}^{-1})_{i},\tau^{-1}\big)\times\big((\eps_{i})_{i},\sigma\big)\times\big((\theta_{i})_{i},\tau\big)=\big((\theta_{\tau^{-1}\sigma\tau(i)}^{-1}\eps_{\tau(i)}\theta_{i})_{i},\tau^{-1}\sigma\tau\big)\,.$$
Appelons cycle d'une permutation signée $g=(\eps,\sigma)$ un cycle de $\sigma$. Le calcul qui précède montre que les longueurs des cycles de $g$ sont les mêmes que les longueurs des cycles de $\widetilde{g}$. Fixons un cycle $c=(x,\sigma(x),\ldots,\sigma^{k-1}(x))$ de la permutation signée $g$. Le cycle correspondant pour $\widetilde{g}$ est $\tau^{-1}(c)=(\tau^{-1}(x),\tau^{-1}(\sigma(x)),\ldots,\tau^{-1}(\sigma^{k-1}(x)))$. De plus, si $(\eps_{x},\eps_{\sigma(x)},\ldots,\eps_{\sigma^{k-1}(x)})$ est la suite des changements de signe de $g$ le long du cycle $c$, alors la suite des changements de signe de $\widetilde{g}$ le long de $\tau^{-1}(c)$ est :
$$\eps_{x}\frac{\theta_{\tau^{-1}(x)}}{\theta_{\tau^{-1}(\sigma(x))}}\,,\,\eps_{\sigma(x)}\frac{\theta_{\tau^{-1}(\sigma(x))}}{\theta_{\tau^{-1}(\sigma^{2}(x))}}\,,\,\ldots\,,\,\eps_{\sigma^{k-1}(x)}\frac{\theta_{\tau^{-1}(\sigma^{k-1}(x))}}{\theta_{\tau^{-1}(\sigma^{k}(x))}}\,.$$
Les rapports de signes peuvent prendre n'importe quelles valeurs, à ceci près qu'ils doivent former un cocycle :
$$\frac{\theta_{\tau^{-1}(x)}}{\theta_{\tau^{-1}(\sigma(x))}}\,\frac{\theta_{\tau^{-1}(\sigma(x))}}{\theta_{\tau^{-1}(\sigma^{2}(x))}}\,\cdots\,\frac{\theta_{\tau^{-1}(\sigma^{k-1}(x))}}{\theta_{\tau^{-1}(\sigma^{k}(x))}}=1\,.$$
Autrement dit, le produit des changements de signe le long du cycle $c$ pour $g$ doit rester égal au produit des changements de signe le long du cycle $\tau^{-1}(c)$ pour $\widetilde{g}$. On conclut que les classes de conjugaison de $\wsym_{n}$ sont indexées par les \textbf{bipartitions} de taille $n$, c'est-à-dire les paires de partitions $(\lambda,\mu)$ telles que $|\lambda|+|\mu|=n$. Les éléments d'une classe de conjugaison $C_{(\lambda,\mu)}$ sont les produits disjoints de cycles pairs de longueurs $\lambda_{1},\lambda_{2},\ldots,\lambda_{r}$ et de cycles impairs de longueurs $\mu_{1},\mu_{2},\ldots,\mu_{s}$. Plus généralement, les classes de conjugaison de $\wsym_{n,r}=(\Z/r\Z)\wr \sym_{n}$ sont en bijection avec les familles de partitions $\Lambda=(\lambda^{(1)},\ldots,\lambda^{(r)})$ telles que $|\lambda^{(1)}|+\cdots+|\lambda^{(r)}|=n$. Nous noterons cet ensemble $\ym_{n,r}$ ; il est également en bijection avec les modules irréductibles sur $\wsym_{n,r}$, et nous noterons plus loin $V^{\Lambda}$ le module associé à une suite de partitions $\Lambda \in \ym_{n,r}$.\bigskip
\bigskip

Vu comme groupe de Coxeter, le groupe hyperoctahédral $\wsym_{n}$ admet pour présentation $\wsym_{n}=\langle s_{0},s_{1},\ldots,s_{n-1}\rangle$, avec 
\begin{align*}&\forall i\geq0,\,\,\, (s_{i})^{2}=1\\
&\forall i \geq 1,\,\,s_{i}s_{i+1}s_{i}=s_{i+1}s_{i}s_{i+1}\qquad;\qquad s_{0}s_{1}s_{0}s_{1}=s_{1}s_{0}s_{1}s_{0}\\
&\forall i,j,\,\,|i-j|\geq 2 \Rightarrow s_{i}s_{j}=s_{j}s_{i}
 \end{align*}
où les $s_{i\geq 1}$ correspondent aux transpositions élémentaires $(i,i+1) \in \sym_{n} \subset \wsym_{n}$, et où $s_{0}$ correspond au changement de signe $((-1,1,\ldots,1),\id_{\lle1,n\rre})$. Compte tenu du théorème \ref{iwahorihecke}, le commutant de l'action de $\Sp(2n,\For_{q})$ sur le module de sa variété de drapeaux est donc l'algèbre d'Iwahori-Hecke $\IH_{q}(\wsym_{n})$, de présentation :
\begin{align*}\IH_{q}(\wsym_{n})=\langle T_{0},T_{1},\ldots,T_{n-1}\rangle,\,\,\text{avec :}\quad&\forall i\geq0,\,\,\, (T_{i}-q)\,(T_{i}+1)=0\\
&T_{0}T_{1}T_{0}T_{1}=T_{1}T_{0}T_{1}T_{0}\\
&\forall i \geq 1,\,\,T_{i}T_{i+1}T_{i}=T_{i+1}T_{i}T_{i+1}\\
&\forall i,j,\,\,|i-j|\geq 2 \Rightarrow T_{i}T_{j}=T_{j}T_{i}\,.
 \end{align*}
Cette \textbf{algèbre d'Hecke de type B}\label{algheckeb} a génériquement la même théorie des représentations que le groupe hyperoctahédral ; ainsi, pour $q \neq 0$ qui n'est pas une racine non triviale de l'unité, les modules irréductibles $V^{\Lambda}(q)$ sur $\IH_{q}(\wsym_{n})$ sont en bijection avec les bipartitions $\Lambda \in \ym_{n,2}$. De plus, $\C[\Sp(2n,\For_{q})/\BSp(2n,\For_{q})]$ se décompose en somme directe de $(\Sp(2n,\For_{q}),\IH_{q}(\wsym_{n}))$-bimodules irréductibles :
\begin{align*}&_{\Sp(2n,\For_q) \curvearrowright }\big\{\C[\Sp(2n,\For_q)/\BSp(2n,\For_q)]\big\}_{\curvearrowleft \,\IH_{q}(\wsym_{n})} \\
&=\sum_{\Lambda \in \ym_{n,2}} {}_{\Sp(2n,\For_q) \curvearrowright }\big(U^\Lambda(q)\big) \otimes_{\C} \big(V^\Lambda(q)\big)_{\curvearrowleft \,\IH_{q}(\wsym_{n})}\,.
\end{align*}
\begin{definition}[B-$q$-mesure de Plancherel]
On appelle $q$-mesure de Plancherel de type B la mesure de probabilité sur les bipartitions associée à la décomposition du module de la variété de drapeaux $\C[\Sp(2n,\For_{q})/\BSp(2n,\For_{q})]$ en composantes unipotentes irréductibles $U^{\Lambda}(q)$. Autrement dit,
$$M^{\mB}_{n,q}(\Lambda)=\frac{\dim V^{\Lambda}(q) \times \dim U^{\Lambda}(q)}{|\Sp(2n,\For_{q})/\BSp(2n,\For_{q})|}=\frac{\dim \Lambda \times D_{\Lambda}(q)}{\{2n!!\}_{q}}\,,$$
où $D_{\Lambda}(q)$ est la dimension du $\Sp(2n,\For_{q})$-module $U^{\Lambda}(q)$.
\end{definition}\clearpage

Donnons une expression plus explicite de cette nouvelle $q$-mesure de Plancherel. Si $\Lambda=(\lambda^{(1)},\lambda^{(2)})$ est une bipartition, nous la représenterons par un bidiagramme de Young dont la partie supérieure droite est le diagramme de Young de $\lambda_{1}$, et dont la partie inférieure gauche est le diagramme de Young de $\lambda_{2}$ renversé. Par exemple, le diagramme de Young de la bipartition $((4,3,1),(5,2))$ de taille $15$ est :

\figcapt{\psset{unit=1mm}\pspicture(-25,-9)(20,17)
\psdot(0,0)
\psline(-25,0)(20,0)(20,5)(0,5)(0,-5)(-25,-5)(-25,0)
\psline(0,-5)(0,-10)(-10,-10)(-10,0)
\psline(-5,-10)(-5,0)
\psline(-15,-5)(-15,0)
\psline(-20,-5)(-20,0)
\psline(15,0)(15,10)(0,10)(0,5)
\psline(5,0)(5,15)(0,15)(0,10)
\psline(10,0)(10,10)
\endpspicture}{Bidiagramme de Young de la bipartition $\Lambda=((4,3,1),(5,2))$.}{Bidiagramme de Young de la bipartition $\Lambda=((4,3,1),(5,2))$}

\noindent Pour tout entier $n \geq 1$, le groupe $\wsym_{n}$ se plonge dans $\wsym_{n+1}$ en envoyant générateurs sur générateurs. Dans ce contexte, les règles d'induction et de restriction pour les modules irréductibles sont tout à fait analogues à celles énoncées dans la proposition \ref{branching} pour le groupe symétrique. Ainsi, pour toute bipartition $M \in \ym_{n+1,2}$,
$$\mathrm{Res}_{\wsym_{n}}^{\wsym_{n+1}}(V^{M}) \simeq_{\C\wsym_{n}} \bigoplus V^{\Lambda}$$
où la somme directe porte sur les bipartitions $\Lambda$ obtenues en ôtant une case du bord inférieur gauche ou du bord supérieur droit de $M$. Par la propriété d'adjonction, on a une description analogue d'un module induit de $\wsym_{n}$ à $\wsym_{n+1}$ à partir d'un module irréductible $V^{\Lambda}$. Le diagramme de Bratteli des groupes hyperoctahédraux est donc le produit gradué de deux arbres de Young, et :
\begin{align*}\dim \left(\Lambda=(\lambda^{(1)},\lambda^{(2)})\right)&=\text{nombres de chemins reliant }(\emptyset,\emptyset)\text{ à }(\lambda^{(1)},\lambda^{(2)})\text{ dans }\ym^{2}\\
&= \binom{|\Lambda|}{|\lambda^{(1)}|} \,\dim \lambda^{1}\,\dim \lambda^{2}=\frac{|\Lambda|\,!}{\prod_{(i,j) \in \lambda^{(1)}} h(i,j)\,\,\prod_{(i,j) \in \lambda^{(2)}} h(i,j) }\,.
\end{align*}
Les règles de branchement subsistent pour les algèbres d'Hecke génériques $\IH(\wsym_{n})$, donc $\dim V^{\Lambda}(q)=\dim \Lambda$ peut comme dans le cas du groupe symétrique être calculée à l'aide d'une formule d'équerres. Concernant les degrés génériques $D_{\Lambda}(q)$, ils ont été calculés en type B par P. N. Hoefmit, voir \cite{Hoe74}. On renvoie également à \cite{Ore99,Ian01} pour une approche passant par la théorie des \textbf{traces de Markov} ; le calcul mené dans ces articles est sans doute le plus simple permettant d'obtenir la formule donnée ci-après. Si $\lambda$ et $\mu$ sont deux partitions, on note $r=\max(\ell(\lambda)-1,\ell(\mu))$, et on introduit les suites
\begin{align*}A&=\{\lambda_{1}+r,\lambda_{2}+r-1,\ldots,\lambda_{r+1}\}=\{\alpha_{1},\ldots,\alpha_{r+1}\}\\
 B&=\{\mu_{1}+r-1,\mu_{2}+r-2,\ldots,\mu_{r}\}=\{\beta_{1},\ldots,\beta_{r}\}\end{align*}
étant entendu que les parts $\lambda_{i}$ ou $\mu_{j}$ peuvent être nulles si l'indice est trop grand. Alors, la formule de Steinberg donnée dans la section \ref{qplancherelprocess} admet pour analogue en type B :
$$D_{\Lambda=(\lambda,\mu)}(q)=\frac{\prod_{i=1}^{n}q^{2i}-1}{2^{r}\,q^{\frac{r(r-1)(4r+1)}{6}}}\,\frac{\prod_{a > a' \in A}q^{a}-q^{a'}\,\,\,\prod_{b>b' \in B} q^{b}-q^{b'}\,\,\,\prod_{a\in A,\,b\in B}q^{a}+q^{b}}{\prod_{a \in A}\prod_{i=1}^{a}q^{2i}-1\,\,\,\prod_{b \in B}\prod_{j=1}^{b}q^{2j}-1}\,.$$
Ceci permet de calculer explicitement en fonction de $q$ les B-$q$-mesures de Plancherel des bipartitions.
\begin{example}
Pour $n=2$, la B-$q$-mesure de Plancherel a pour valeurs :
\begin{align*}
&M^{\mB}_{n,q}(2,\emptyset)=1/(q^4 + 2q^3 + 2q^2 + 2q + 1)\\
&M^{\mB}_{n,q}(1^{2},\emptyset)=q/(2q^2 + 4q + 2)\\
&M^{\mB}_{n,q}(1,1)=q/(q^2 + 1)\\
&M^{\mB}_{n,q}(\emptyset,2)=q/(2q^2 + 4q + 2)\\
&M^{\mB}_{n,q}(\emptyset,1^{2})=q^4/(q^4 + 2q^3 + 2q^2 + 2q + 1)
\end{align*}
\end{example}
\bigskip\bigskip

De nouveau, la $q$-mesure de Plancherel de type B fait sens pour tout paramètre réel positif $q$, et on a une propriété de symétrie analogue à l'identité $M_{n,q}(\lambda)=M_{n,q^{-1}}(\lambda')$ en type A. Ainsi :
$$M_{n,q}^{\mB}(\lambda,\mu)=M_{n,q^{-1}}^{\mB}(\mu',\lambda')\,.$$
D'autre part, on peut comme en type A construire un processus markovien $(\Lambda_{n})_{n \in \N}$ d'espace d'états $\ym^{2}$, et de lois marginales les B-$q$-mesures de Plancherel. La $q$-proba\-bi\-lité de transition entre deux bipartitions est
$$p_{q}(\Lambda,M)=\begin{cases}
\frac{1}{\{2|M|\}_{q}}\,\frac{D_{M}(q)}{D_{\Lambda}(q)} &\text{si }\Lambda \nearrow M,\\
0&\text{sinon.}\end{cases}$$
Les premières probabilités de transition du B-$q$-processus de Plancherel sont représentées sur la figure \ref{bqplancherelprocess}. Notons que lorsque $q$ tend vers $1$, on retrouve les mesures de Plancherel usuelles des groupes hyperoctahédraux, qui vérifient :
$$M_{\wsym_{n}}(\lambda,\mu)=\frac{1}{2^{n}}\binom{n}{|\lambda|,|\mu|}\,M_{\sym_{|\lambda|}}(\lambda)\,M_{\sym_{|\mu|}}(\mu)\,.$$
\figcapt{\footnotesize{
\psset{unit=1mm}\pspicture(-70,-60)(70,2.5)
\psdot(0,-0.5)
\rput(7,0){$D=1$}
\psline[linecolor=red]{->}(0,-3)(20,-13)
\psline[linecolor=red]{->}(0,-3)(-20,-13)
\psdots(-17.5,-15)(17.5,-20)
\psline(-17.5,-15)(-22.5,-15)(-22.5,-20)(-17.5,-20)(-17.5,-15)\rput(-30,-18){$D=q$}
\psline(17.5,-15)(22.5,-15)(22.5,-20)(17.5,-20)(17.5,-15)\rput(30,-18){$D=1$}
\rput(18,-7.5){\textcolor{red}{$\frac{1}{q+1}$}}
\rput(-18,-7.5){\textcolor{red}{$\frac{q}{q+1}$}}
\psline(-70,-40)(-60,-40)(-60,-35)(-70,-35)(-70,-40)\psline(-65,-35)(-65,-40)\rput(-65,-45){$D=\frac{q^{3}+q}{2}$}
\psline(70,-40)(60,-40)(60,-35)(70,-35)(70,-40)\psline(65,-35)(65,-40)\rput(65,-44){$D=1$}
\psline(-40,-40)(-35,-40)(-35,-50)(-40,-50)(-40,-40)\psline(-40,-45)(-35,-45)\rput(-37.5,-53){$D=q^{4}$}
\psline(40,-40)(35,-40)(35,-50)(40,-50)(40,-40)\psline(40,-45)(35,-45)\rput(37.5,-54){$D=\frac{q^{3}+q}{2}$}
\psdots(-60,-35)(-35,-40)(35,-50)(60,-40)(0,-45)
\psline(0,-50)(0,-40)(5,-40)(5,-45)(-5,-45)(-5,-50)(0,-50)\rput(0,-55){$D=\frac{q^{3}+2q^{2}+q}{2}$}
\psline[linecolor=red]{->}(20,-22)(65,-34)
\rput(60,-27.5){\textcolor{red}{$\frac{1}{q^{3}+q^{2}+q+1}$}}
\psline[linecolor=red]{->}(20,-22)(37.5,-39)
\psframe*[linecolor=white,fillcolor=white](27,-35)(35,-27)
\rput(30,-31.5){\textcolor{red}{$\frac{q}{2(q+1)}$}}
\psline[linecolor=red]{->}(20,-22)(1,-39)
\psline[linecolor=red]{->}(-20,-22)(-65,-34)
\rput(-57,-27.5){\textcolor{red}{$\frac{1}{2(q+1)}$}}
\psline[linecolor=red]{->}(-20,-22)(-37.5,-39)
\psframe*[linecolor=white,fillcolor=white](-35,-35)(-27,-27)
\rput(-30,-31){\textcolor{red}{$\frac{q^{3}}{q^{3}+q^{2}+q+1}$}}
\psline[linecolor=red]{->}(-20,-22)(-1,-39)
\psframe*[linecolor=white,fillcolor=white](-15,-35)(15,-27)
\rput(10,-31){\textcolor{red}{$\frac{q(q+1)}{2(q^{2}+1)}$}}
\rput(-10,-31){\textcolor{red}{$\frac{q+1}{2(q^{2}+1)}$}}
\endpspicture}
}{Degrés génériques et $q$-probabilités de transition sur les deux premiers niveaux du bigraphe de Young.\label{bqplancherelprocess}}{Degrés génériques et $q$-probabilités de transition en type B}

Les expressions des probabilités de transition $p_{q}(\Lambda,M)$ peuvent être simplifiées comme suit. On note $r$ l'entier associé à une bipartition $\Lambda$ pour le calcul des degrés génériques. Si $M=(\mu^{(1)},\mu^{(2)})$ est obtenue à partir de $\Lambda=(\lambda^{(1)},\lambda^{(2)})$ en ajoutant une case à la $l$-ième ligne de $\lambda^{(1)}$ avec $l \leq r+1$, alors
$$p_{q}(\Lambda,M)=\frac{1}{\{2(\alpha_{l}+1)\}_{q}}\,\left(\prod_{\substack{i \in \lle 1,r+1\rre\\ i \neq l}} \frac{q^{\alpha_{i}} -q^{\alpha_{l}+1}}{q^{\alpha_{i}} -q^{\alpha_{l}}} \right)\,\left( \prod_{j\in \lle 1,r\rre} \frac{q^{\alpha_{l}+1}+q^{\beta_{j}} }{q^{\alpha_{l}}+q^{\beta_{j}}} \right),$$
où $A=(\alpha_{i})_{i \leq r+1}$ et $B=(\beta_{j})_{j \leq r}$ sont les coordonnées modifiées de la bipartition $\Lambda$. Si $M$ est obtenue à partir de $\Lambda$ en rajoutant une case à la $r+2$-ième ligne de $\lambda^{(1)}$ (qui était vide), alors 
$$p_{q}(\Lambda,M)=\frac{q^{2r+1}}{2}\,\left(\prod_{i \in \lle 1,r+1\rre} \frac{q^{\alpha_{i}}-1}{q^{\alpha_{i}+1}-1} \right)\,\left(\prod_{j\in \lle1,r\rre} \frac{q^{\beta_{j}}+1}{q^{\beta_{j}+1}+1}\right).$$
De même, si $M$ est obtenue à partir de $\Lambda$ en rajoutant une case à la $l$-ième ligne de $\lambda^{(2)}$ avec $l \leq r$, alors
$$p_{q}(\Lambda,M)=\frac{1}{\{2(\beta_{l}+1)\}_{q}}\,\left( \prod_{i\in \lle 1,r+1\rre} \frac{q^{\beta_{l}+1}+q^{\alpha_{i}} }{q^{\beta_{l}}+q^{\alpha_{i}}} \right)\,\left(\prod_{\substack{j \in \lle 1,r\rre\\ j\neq l}} \frac{q^{\beta_{j}} -q^{\beta_{l}+1}}{q^{\beta_{j}} -q^{\beta_{l}}} \right),$$
et si $M$ est obtenue à partir de $\Lambda$ en rajoutant une case à la $r+1$-ième ligne de $\lambda^{(2)}$ (qui était vide), alors 
$$p_{q}(\Lambda,M)=\frac{q^{2r+1}}{2}\,\left(\prod_{i \in \lle 1,r+1\rre} \frac{q^{\alpha_{i}}+1}{q^{\alpha_{i}+1}+1} \right)\,\left(\prod_{j\in \lle1,r\rre} \frac{q^{\beta_{j}}-1}{q^{\beta_{j}+1}-1}\right).$$
Avec ces expressions simplifiées des probabilités de transition, on peut aisément programmer\footnote{En type A, pour tirer aléatoirement un diagramme sous la $q$-mesure de Plancherel $M_{n,q}$, on pouvait réaliser le $q$-processus de Knuth jusqu'au rang $n$, puis projeter par l'algorithme RSK la permutation obtenue. Il est vraisemblable qu'en type B, la B-$q$-mesure de Plancherel $M_{n,q}^{\mathrm{B}}$ soit aussi l'image par RSK d'une statistique sur les permutations signées de $\wsym_{n}$, mais ce problème reste encore ouvert.} le B-$q$-processus de Plancherel en \texttt{sage}. 
\figcapt{\psset{unit=1.4mm}
\pspicture(-47,-8)(67,10)
\psdots(0,0)
\psline(67,0)(-47,0)(-47,-1)(0,-1)(0,1)(67,1)(67,0)
\psline(38,2)(0,2)
\psline(10,3)(0,3)
\psline(8,4)(0,4)
\psline(3,5)(0,5)
\psline(2,6)(0,6)
\psline(1,0)(1,7)(0,7)(0,1)
\multido{\n=39+1}{28}{\psline(\n,0)(\n,1)}
\multido{\n=11+1}{28}{\psline(\n,0)(\n,2)}
\multido{\n=9+1}{2}{\psline(\n,0)(\n,3)}
\multido{\n=4+1}{5}{\psline(\n,0)(\n,4)}
\psline(3,0)(3,5)
\psline(2,0)(2,6)
\psline(-14,-2)(0,-2)
\psline(-6,-3)(0,-3)
\psline(-2,-4)(0,-4)
\psline(-2,-5)(0,-5)
\multido{\n=15+1}{33}{\psline(-\n,0)(-\n,-1)}
\multido{\n=7+1}{8}{\psline(-\n,0)(-\n,-2)}
\multido{\n=3+1}{4}{\psline(-\n,0)(-\n,-3)}
\psline(-2,0)(-2,-5)
\psline(-1,0)(-1,-5)
\psline(0,-1)(0,-5)
\endpspicture}{Bidiagramme de Young aléatoire tiré suivant la B-$q$-mesure de Plancherel de paramètres $q=2/3$, $n=200$. Les parts de $\Lambda$ sont $(67,38,10,8,3,2,1)$ et $(47,14,6,2,2)$.\label{qplanhalfb}}{Bidiagramme aléatoire tiré suivant une B-$q$-mesure de Plancherel}
La figure \ref{qplanhalfb} présente un bidiagramme de Young tiré suivant la B-$q$-mesure de Plancherel de paramètre $q=2/3$. Donnons un autre exemple pour $q=1/2$ :
$$|\Lambda|=200\qquad;\qquad \Lambda=(108,12,4),(47,25,3,1)\,.$$
Si l'on réordonne les parts des deux partitions, on obtient une partition dont les parts suivent approximativement une progression géométrique de paramètre $1/2$ : ainsi, $108 \sim 200/2$, $47 \sim 200/4$, $25 \sim 200/8$, etc. On devine donc que les parts des deux partitions ont un comportement asymptotique identique à celui observé en type A. Ceci sera démontré dans la section \ref{basymptotic}. \bigskip\bigskip

La répartition des parts entre les deux partitions de la bipartition --- c'est-à-dire, quelle partition re\c coit la plus grande part, quelle partition re\c coit la seconde plus grande part, etc. --- reste en revanche aléatoire lorsque $n$ tend vers l'infini. Par exem\-ple, la probabilité pour que la plus grande part de $\Lambda$ se trouve dans la première partition $\lambda^{(1)}$, c'est-à-dire la probabilité
$$M_{n,q}^{\mB}[\lambda^{(1)}_{1}\geq \lambda^{(2)}_{1}]\,,$$
est asymptotiquement strictement comprise entre $0$ et $1$ : ainsi, sur $10000$ tirages de bidiagrammes sous la B-$q$-mesures de Plancherel de paramètres $n=200$ et $q=1/2$, $7934$ bipartitions avaient leur plus grande part dans $\lambda^{(1)}$, et $2066$ bipartitions avaient leur plus grande part dans $\lambda^{(2)}$. En particulier, la taille de la première part d'une des deux partitions, disons $\lambda^{(1)}$, ne converge pas en probabilité vers une constante (contrairement à ce qui passe en type A). \bigskip
\bigskip

Essentiellement, ce phénomène est lié au fait algébrique suivant :  le plongement canonique d'un produit de groupes hyperoctahédraux $\wsym_{n_{1}}\times \wsym_{n_{2}}$ dans $\wsym_{n=n_{1}+n_{2}}$ n'est pas compatible avec les structures de groupes de Coxeter, car l'élément $s_{0}'$ du second groupe $\wsym_{n_{2}}$, originellement de longueur $1$, est envoyé sur
$$s_{n_{1}}\cdots s_{1}s_{0}s_{1}\cdots s_{n_{1}}\,,$$
qui est de longueur $2n_{1}+1$. Dit autrement, le produit $\wsym_{n_{1}} \times \wsym_{n_{2}}$ n'est pas un sous-groupe parabolique de $\wsym_{n_{1}+n_{2}}$, voir \cite[chapitre 2]{GP00} pour des précisions sur cette terminologie. Ceci conduit à la non-factorisation asymptotique des espérances de certaines observables de bidiagrammes, et donc à la non-convergence de ces observables. Nous expliquerons ceci plus précisément dans le paragraphe \ref{heckebmix}. En particulier, nous énoncerons une conjecture concernant la valeur limite de la probabilité pour que la $i$-ième part de $\lambda=\lambda^{(1)}\sqcup \lambda^{(2)}$ tombe dans $\lambda^{(1)}$ (\emph{cf.} la conjecture \ref{bqmix}). Néanmoins, ceci ne détermine que partiellement le comportement des deux partitions, car les variables de Bernoulli associées à ces événements sont hautement non indépendantes. Dans les deux prochaines sections, nous verrons qu'en utilisant des techniques d'observables de diagrammes, on peut en théorie calculer tous les moments des lois limites régissant la répartition des parts entre les deux partitions. Mais malheureusement, il semble qu'il n'y ait pas de formule simple générale pour ces lois limites. 
\bigskip

\section[Résultats asymptotiques pour la B-$q$-mesure de Plancherel]{Résultats asymptotiques pour la B-$q$-mesure de Plancherel}\label{basymptotic}

Dans ce qui suit, on appelle \textbf{observable de bidiagrammes} une fonction sur les bipartitions qui est combinaison linéaire de moments mixtes 
$$p_{\mu^{(1)},\mu^{(2)}}\big(\Lambda=(\lambda^{(1)},\lambda^{(2)})\big)=p_{\mu^{(1)}}(\lambda^{(1)})\,\times\,p_{\mu^{(2)}}(\lambda^{(2)})\,,$$
\emph{i.e.}, une combinaison linéaire de produits de moments de Frobenius de $\lambda^{(1)}$ et de moments de Frobenius de $\lambda^{(2)}$. L'ensemble des observables de bidiagrammes forme une algèbre $\obs^{\mB} \simeq \obs \otimes_{\C} \obs$, et les moments mixtes $(p_{M})_{M \in \ym \times \ym}$ constituent une base de $\obs^{\mB}$. Comme dans le chapitre \ref{tool}, on peut donner une autre base de $\obs$ dont les évalutations sont des valeurs renormalisées de caractères irréductibles, et qu'on peut interpréter comme éléments de la sous-algèbre des invariants d'une \textbf{algèbre des permutations partielles signées}. Pour commencer, donnons un bref exposé de la théorie des représentations des groupes hyperoctahédraux. Si $\Lambda=(\lambda^{(1)},\lambda^{(2)})$ et $M=(\mu^{(1)},\mu^{(2)})$ sont deux bipartitions de même taille $n$, on note $\varsigma^{\Lambda}(M)$ la valeur du caractère irréductible non normalisé de $\wsym_{n}$ de type $\Lambda$ en une permutation signée de type $M$. Alors, on peut donner pour $\varsigma^{\Lambda}(M)$ une formule de Frobenius analogue à celle du paragraphe \ref{frobenius}. Ainsi, si $X$ et $Y$ sont deux alphabets indépendants,
$$\varsigma^{\Lambda}(M)=\scal{s_{\lambda^{(1)}}(X)\,s_{\lambda^{(2)}}(Y)}{p_{\mu^{(1)}}(X+Y)\,p_{\mu^{(2)}}(X-Y)}$$
étant entendu que $\Lambda(X)\otimes \Lambda(Y)$ est muni du produit tensoriel des produits scalaires sur $\Lambda(X)$ et $\Lambda(Y)$. Plus généralement d'ailleurs, pour un produit en couronne $\wsym_{n,r}=(\Z/r\Z) \wr \sym_{n}$, on peut donner une formule de Frobenius dans l'algèbre $\bigotimes_{j=1}^{r} \Lambda(X_{j})$, voir \cite[appendice B]{Mac95}. Ainsi, étant donnés des alphabets indépendants $X_{1},\ldots,X_{r}$ et une racine primitive $r$-ième de l'unité $\zeta$, notons :
\begin{align*}
P_{k}^{(i)}(X_{1},\ldots,X_{r})&=\sum_{j=1}^{r} \zeta^{-ij}p_{k}(X_{j})\\
P_{M}(X_{1},\ldots,X_{r})&=\prod_{i=1}^{r} P_{\mu^{(i)}}^{(i)}(X_{1},\ldots,X_{r})=\prod_{i=1}^{r} \prod_{j=1}^{\ell(\mu^{(i)})}P_{\mu^{(i)}_{j}}^{(i)}(X_{1},\ldots,X_{r})\\
S_{\Lambda}(X_{1},\ldots,X_{r})&=\prod_{i=1}^{r}s_{\lambda^{(i)}}(X_{i})
\end{align*}
où $\Lambda=(\lambda^{(1)},\ldots,\lambda^{(r)})$ et $M=(\mu^{(1)},\ldots,\mu^{(r)})$ sont
des $r$-uplets de partitions dans $\ym_{n,r}$. La formule de Frobenius pour les caractères de $\wsym_{n,r}$ s'écrit alors :
$$\forall M \in \ym_{n,r},\,\,\,P_{M}(X_{1},\ldots,X_{r})=\sum_{\Lambda \in \ym_{n,r}} \varsigma^{\Lambda}(M)\,S_{\Lambda}(X_{1},\ldots,X_{r})\,.$$
On retrouve la formule pour les caractères du groupe hyperoctahédral en posant $r=2$ et $\zeta=-1$.\bigskip
\bigskip

Comme dans le paragraphe \ref{centralcharacter}, on appelle permutation partielle signée d'ordre $n$ la donnée d'un couple $(\sigma,S)$, où $S$ est une partie de $\lle 1,n\rre$ et $\sigma$ appartient à $\wsym_{S}$. Dit autrement, $\sigma$ est un élément de $\wsym_{n}$ tel que $\sigma(i)=i$ pour tout $i \notin S$. La construction ferait également sens dans un produit en couronne général $\wsym_{n,r}$, voir \cite[\S5]{Wang04}. Le produit de deux permutations partielles signées est de nouveau défini par 
$$(\sigma,S)\,(\tau,T)=(\sigma\tau,S\cup T)\,.$$
L'algèbre (complexe) $\blg_{n}^{\mB}$ du semi-groupe correspondant à ces objets et à ce produit se projette sur l'algèbre de groupe $\C\wsym_{n}$ en oubliant les supports des permutations partielles signées, et exactement comme dans le paragraphe \ref{centralcharacter}, on peut construire une limite projective $\blg_{\infty}^{\mB} =\varprojlim_{n \to \infty}\blg_{n}^{\mB}$ dans la catégorie des algèbres filtrées ; ses éléments sont les combinaisons linéaires éventuellement infinies de permutations partielles signées (sans restriction sur le support).\bigskip

Le groupe $\wsym_{\infty}=\varinjlim_{n \to \infty} \wsym_{n}$ agit par conjugaison sur les éléments de $\blg_{\infty}^{\mB}$, et la sous-algèbre des invariants $\alg_{\infty}^{\mB}=(\blg_{\infty}^{\mB})^{\wsym_{\infty}}$ est une algèbre commutative engendrée linéairement par les sommes formelles infinies $\varSigma_{M}=\varSigma_{\mu^{(1)},\mu^{(2)}}$ de permutations partielles signées du type
$$\prod_{i=1}^{r} \left(\eps_{i},\big(a_{i,1},\ldots,a_{i,\mu^{(1)}_{i}}\big)\right)\,\circ\, \prod_{j=1}^{s}\left(\theta_{j},\big(b_{j,1},\ldots,b_{j,\mu^{(2)}_{j}}\big)\right),\left\{a_{i,k},b_{j,l} \,\,\bigg|\,\,\substack{1 \leq k \leq \mu^{(1)}_{i},\,\,1 \leq i \leq \ell(\mu^{(1)}) \\ 1 \leq l \leq \mu^{(1)}_{j},\,\,1 \leq j \leq \ell(\mu^{(2)})}\right\}$$
où les $a_{i,k}$ et les $b_{j,l}$ sont des entiers tous distincts, et où
$$\eps_{i}=\eps(a_{i,1})\,\eps(a_{i,2})\,\cdots\,\eps(a_{i,\mu^{(1)}_{i}})\qquad;\qquad
\theta_{j}=\theta(b_{j,1})\,\theta(b_{j,2})\,\cdots\,\theta(b_{j,\mu^{(2)}_{j}})$$
sont des éléments de $(\Z/2\Z)^{(\infty)} \subset \wsym_{\infty}$ tels que le produit des signes d'un $\eps_{i}$ soit $+1$, et le produit des signes d'un $\theta_{j}$ soit $-1$. Ainsi, le premier produit correspond à des cycles pairs de longueurs $\mu^{(1)}_{i}$, et le second produit correspond à des cycles impairs de longueurs $\mu^{(2)}_{j}$. Si $\phi_{n}$ est la projection canonique $\blg_{\infty}^{\mB} \to \blg_{n}^{\mB}$, alors $\phi_{n}(\varSigma_{M})$ est donné par la même somme, mais avec les entiers $a_{i,j}$ et $b_{k,l}$ restreints à l'intervalle $\lle 1,n\rre$. En particulier, le nombre d'éléments de $\phi_{n}(\varSigma_{M})$ est 
$$2^{|\mu^{(1)}|-\ell(\mu^{(1)}) +|\mu^{(2)}|-\ell(\mu^{(2)})}\,n(n-1)\cdots(n-|M|+1)\,.$$
Par conséquent, si $\pi_{n}$ désigne la projection canonique $\blg_{n}^{\mB} \to \C\wsym_{n}$ correspondant à l'oubli des supports,  alors $\varSigma_{M,n}=\pi_{n}(\phi_{n}(\varSigma_{M}))$ vaut $0$ si $n<|M|=|\mu^{(1)}|+|\mu^{(2)}|$, et dans le cas contraire, c'est un multiple de la classe de conjugaison $C_{\mu^{(1)}\sqcup 1^{n-|M|}, \mu^{(2)}}$ :
$$\varSigma_{M,n}=2^{|M|-\ell(M)}\,n^{\downarrow |M|}\,\widetilde{C}_{\mu^{(1)} \sqcup 1^{n-|M|},\mu^{(2)}}$$
où pour toute bipartition $\Lambda$, $\widetilde{C_{\Lambda}}=C_{\Lambda}/\card C_{\Lambda}$. Comme en type A, ceci permet d'évaluer un élément $\varSigma_{M}$ contre une bipartition $\Lambda \in \ym_{n,2}$ :
$$\varSigma_{M}(\Lambda)=\chi^{\Lambda}(\varSigma_{M,n})\,,$$
où $\chi^{\Lambda}$ désigne le caractère irréductible normalisé, c'est-à-dire que $(\dim \Lambda)\,\chi^{\Lambda}=\varsigma^{\Lambda}$.\bigskip
\bigskip

De nouveau, il n'est pas du tout évident que ces évaluations définissent des observables de bidiagrammes dans $\obs^{\mB}=\C[(p_{M})_{M \in \ym \times \ym}]$. Ce point est une conséquence des deux lemmes suivants :
\begin{lemma}[Caractère central d'un cycle impair ou pair]\label{paidinfull}
Pour tout entier $k$ et toute bipartition $\Lambda=(\lambda^{(1)},\lambda^{(2)})$, on a :
\begin{align*}
\varSigma_{k,\emptyset}(\Lambda)&=2^{k-1}\big(\varSigma_{k}(\lambda^{(1)})+\varSigma_{k}(\lambda^{(2)})\big)\,;\\
\varSigma_{\emptyset,k}(\Lambda)&=2^{k-1}\big(\varSigma_{k}(\lambda^{(1)})-\varSigma_{k}(\lambda^{(2)})\big)\,.
\end{align*}
\end{lemma}\medskip

\begin{lemma}[Graduation sur l'algèbre des permutations partielles signées et factorisation en plus haut degré]\label{followtheleader}
Pour toute bipartition $M$, on pose $\deg \varSigma_{M}=|M|$. Ceci définit une graduation d'algèbres sur $\alg_{\infty}^{\mB}$, et pour toutes bipartitions $M_{1}$ et $M_{2}$, 
$$\varSigma_{M_{1}}\,*\,\varSigma_{M_{2}}=\varSigma_{M_{1} \sqcup M_{2}}+(\text{termes de degré inférieur})\,.$$
\end{lemma}\medskip

\begin{proof}[Preuve du lemme \ref{paidinfull}]
Supposons dans un premier temps $|\Lambda|\geq k$. La formule de Frobenius en type B assure que :
\begin{align*}\varsigma^{\Lambda}(k1^{n-k},\emptyset)&=\scal{s_{\lambda^{(1)}}(X)\,s_{\lambda^{(2)}}(Y)}{p_{k}(X+Y)\,(p_{1}(X+Y))^{n-k}}\\
&=\scal{s_{\lambda^{(1)}}(X)\,s_{\lambda^{(2)}}(Y)}{p_{k}(X)\,(p_{1}(X+Y))^{n-k}}\\
&\,\,\,+\scal{s_{\lambda^{(1)}}(X)\,s_{\lambda^{(2)}}(Y)}{p_{k}(Y)\,(p_{1}(X+Y))^{n-k}}\,.
\end{align*}
Examinons le premier terme de la dernière expression. Si $|\lambda^{(1)}|<k$, alors la partie droite du produit scalaire n'a aucun terme de degré $|\lambda^{(1)}|$ en $X$ et de degré $|\lambda^{(2)}|$ en $Y$. En effet, en développant la puissance de $p_{1}(X+Y)$, on n'obtient que des termes de degré supérieur à $k$ en $X$. Le premier terme est donc nul si $|\lambda^{(1)}|<k$, et dans le cas contraire, il vaut 
\begin{align*}&\sum_{r=0}^{n-k} \binom{n-k}{r}\scal{s_{\lambda^{(1)}}(X)\,s_{\lambda^{(2)}}(Y)}{p_{k1^{r}}(X)\,p_{1^{n-k-r}}(Y)}\\
&=\binom{n-k}{|\lambda^{(1)}|-k}\scal{s_{\lambda^{(1)}}(X)\,s_{\lambda^{(2)}}(Y)}{p_{k1^{|\lambda^{(1)}|-k}}(X)\,p_{1^{|\lambda^{(2)}|}}(Y)} \\
&=\binom{n-k}{|\lambda^{(1)}|-k}\,(\dim \lambda^{(2)})\,\,\varsigma^{\lambda^{(1)}}(k1^{|\lambda^{(1)}|-k})\,.\end{align*}
Si l'on multiplie par le facteur $\frac{n^{\downarrow k}}{\dim \Lambda}=\frac{(n-k)^{\downarrow |\lambda^{(1)}|-k}\,\,|\lambda^{(1)}|!}{(\dim \lambda^{(1)})\,(\dim \lambda^{(2)})}$, on obtient donc 
$$|\lambda^{(1)}|^{\downarrow k}\,\,\frac{\varsigma^{\lambda^{(1)}}(k1^{|\lambda^{(1)}|-k}) }{\dim \lambda^{(1)}}=|\lambda^{(1)}|^{\downarrow k}\,\,\chi^{\lambda^{(1)}}(k1^{|\lambda^{(1)}|-k})=\varSigma_{k}(\lambda^{(1)})\,,$$
et d'après ce qui précède, le premier terme renormalisé par ce facteur vaut aussi $\varSigma_{k}(\lambda^{(1)})=0$ lorsque $|\lambda^{(1)}|<k$. De fa\c con symétrique, le second terme multiplié par $\frac{n^{\downarrow k}}{\dim \Lambda}$ donne $\varSigma_{k}(\lambda^{(2)})$, et ce quelque soit la taille de $\lambda^{(2)}$. Ainsi, si $|\Lambda|\geq k$, alors
$$\varSigma_{k,\emptyset}(\Lambda)=2^{k-1}\frac{n^{\downarrow k}}{\dim \Lambda}\,\varsigma^{\Lambda}(k1^{|\Lambda|-k},\emptyset)=2^{k-1}\big(\varSigma_{k}(\lambda^{(1)})+\varSigma_{k}(\lambda^{(2)})\big)\,.$$
Si $|\Lambda|<k$, alors le terme de droite est par définition nul, et les deux termes $\varSigma_{k}(\lambda^{(1)})$ et $\varSigma_{k}(\lambda^{(2)})$ le sont également, car $|\lambda^{(1)}|\leq |\Lambda|<k$ et $|\lambda^{(2)}|\leq |\Lambda|<k$. L'identité reste donc vraie pour $|\Lambda|<k$. Finalement, le cas de $\varSigma_{\emptyset,k}(\Lambda)$ se traite exactement de la même fa\c con, en rempla\c cant $p_{k}(X+Y)$ par $p_{k}(X-Y)$.
\end{proof}
\begin{proof}[Preuve du lemme \ref{followtheleader}]
Le principe est exactement le même qu'en type A. Le degré défini dans l'énoncé du lemme \ref{followtheleader} est la restriction à $\alg_{\infty}^{\mB}$ du degré sur $\blg_{\infty}^{\mB}$ défini par $\deg(\sigma,S)=|S|$. Or, étant données deux permutations partielles signées $(\sigma,S)$ et $(\tau,T)$, on a 
$$\deg\big((\sigma,S)\,(\tau,T)\big)=\deg(\sigma\tau,S\cup T)|=|S\cup T| \leq |S| + |T|\,,$$
avec égalité si et seulement si $S$ et $T$ sont disjoints. L'inégalité assure déjà que 
$$\deg (\varSigma_{M_{1}}\,*\,\varSigma_{M_{2}}) \leq \deg \varSigma_{M_{1}} + \deg \varSigma_{M_{2}}\,.$$
Lorsqu'on développe le produit de $\varSigma_{M_{1}}$ par $\varSigma_{M_{2}}$, les produits d'une permutation partielle signée $(\sigma,S)$ de type $M_{1}$ avec une permutation partielle signée $(\tau,T)$ de type $M_{2}$ de support $T$ disjoint de $S$ donnent des permutations partielles signées de type $M_{1}\sqcup M_{2}$. Ainsi, ces produits contribuent à un terme $\varSigma_{M_{1}\sqcup M_{2}}$, et les autres produits sont de degré inférieur à $|M_{1}|+|M_{2}|-1$, car ils correspondent à des ensembles non disjoints. Ainsi,
$$\deg (\varSigma_{M_{1}}\,*\,\varSigma_{M_{2}}) = \deg \varSigma_{M_{1}} + \deg \varSigma_{M_{2}}\,,$$
et le terme de plus haut degré est bien $\varSigma_{M_{1}\sqcup M_{2}}$.
\end{proof}
\bigskip

\begin{proposition}[Algèbre des observables de bidiagrammes]\label{biobs}
Les symboles $\varSigma_{M}$ appartiennent à l'algèbre $\obs^{\mB}$ des observables de bidiagrammes, et ils en forment une base linéaire lorsque $M$ parcourt $\ym \times \ym$. De plus, $\varSigma_{M}=\varSigma_{\mu^{(1)},\mu^{(2)}}$ a même composante de plus haut degré que 
$$2^{|M|-\ell(M)}\left(\prod_{i=1}^{r} p_{\mu^{(1)}_{i},\emptyset}+p_{\emptyset,\mu^{(1)}_{i}} \right)\,\left(\prod_{j=1}^{s} p_{\mu^{(2)}_{j},\emptyset}-p_{\emptyset,\mu^{(2)}_{j}}\right).$$
\end{proposition}
\begin{proof}
Grâce au lemme \ref{followtheleader}, pour montrer que les $\varSigma_{M}$ appartiennent à $\obs^{\mB}$, il suffit de traiter le cas des symboles $\varSigma_{k,\emptyset}$ et $\varSigma_{\emptyset,k}$. En effet :
$$\varSigma_{\mu^{(1)},\mu^{(2)}}=\left(\prod_{i=1}^{r}\varSigma_{\mu^{(1)}_{i},\emptyset}\right)\,\left(\prod_{j=1}^{s}\varSigma_{\emptyset, \mu^{(2)}_{j}}\right)+(\text{termes de degré inférieur})\,.
$$
Or, on sait en type A que $\varSigma_{k}$ est une certaine combinaison linéaire d'observables $p_{\mu}$, avec comme terme de plus haut degré $p_{k}$ (\emph{cf.} la proposition \ref{whitelight}). En utilisant le lemme \ref{paidinfull}, on voit donc que :
\begin{align*}
\varSigma_{k,\emptyset}&=2^{k-1}(p_{k,\emptyset}+p_{\emptyset,k})+(\text{combinaison linéaire de }p_{\mu,\emptyset} \text{ et de }p_{\emptyset,\mu}\text{ avec }|\mu|<k)\,;\\
\varSigma_{\emptyset,k}&=2^{k-1}(p_{k,\emptyset}-p_{\emptyset,k})+(\text{combinaison linéaire de }p_{\mu,\emptyset} \text{ et de }p_{\emptyset,\mu}\text{ avec }|\mu|<k)\,.
\end{align*}
Ceci prouve presque immédiatement l'ensemble de la proposition. Ainsi, l'identification $\obs \sim \alg_{\infty}$ en type A reste vraie en type B, avec $\obs^{\mB}\sim \alg_{\infty}^{\mB}$.
\end{proof}\bigskip
\bigskip

Finalement, on peut  comme dans le paragraphe \ref{quantum} quantifier la base $(\varSigma_{M})_{M \in \ym \times \ym}$ et introduire des observables $(\varSigma_{M,q})_{M \in \ym \times \ym}$ adaptées à la déformation de $\C\wsym_{n}$ en $\IH_{q}(\wsym_{n})$. Pour toute partition $\mu$, on note $\widetilde{\mu}$ la suite $(\mu_{r},\mu_{r-1}, \ldots,\mu_{1})$. Alors, si $M=(\mu^{(1)},\mu^{(2)})$ est un bipartition de taille $n$, on note $\sigma_{M}$ la permutation signée
\begin{align*}&\prod_{i=1}^{r}\left(\overline{\widetilde{\mu}^{(2)}_{1}+\cdots+\widetilde{\mu}^{(2)}_{i-1}+1}\,,\,(\widetilde{\mu}^{(2)}_{1}+\cdots+\widetilde{\mu}^{(2)}_{i-1}+1,\ldots,\widetilde{\mu}^{(2)}_{1}+\cdots+\widetilde{\mu}^{(2)}_{i})\right)\\
&\times \prod_{j=1}^{s}\left(\overline{0}\,,\,(|\mu^{(2)}|+\mu^{(1)}_{1}+\cdots+\mu^{(1)}_{j-1}+1,\ldots,|\mu^{(2)}|+\mu^{(1)}_{1}+\cdots+\mu^{(1)}_{i})\right)
\end{align*}
où $\overline{k}=(1,\ldots,1,-1_{(k)},1,\ldots,1) \in (\Z/2\Z)^{n}$, avec pour convention $\overline{0}=(1,\ldots,1)$. La permutation $\sigma_{M}$ est de longueur minimale\footnote{Cette longueur minimale vaut $|\mu^{(1)}|+|\mu^{(2)}|+2b(\mu^{(2)})-\ell(\mu^{(1)})$.} dans la classe de conjugaison $C_{M}$. \medskip
\begin{example}
Si $M=(2,2,1),(3,2)$, alors $\sigma_{M}=\overline{1}\,(1,2)\,\overline{3}\,(3,4,5)\,(6,7)(8,9)(10)\,.$
\end{example}\bigskip

On note $T_{M}$ l'élément correspondant à $\sigma_{M}$ dans $\IH_{q}(\wsym_{n})$, c'est-à-dire l'élément qui a la même décomposition réduite. Comme en type A, les caractères irréductibles $\chi^{\Lambda}(q)$ sont entièrement déterminées par leurs valeurs en les $T_{M}$, voir \cite[\S8.2]{GP00}. Ainsi, la table des caractères de l'algèbre d'Hecke de type B $\IH_{q}(\wsym_{n})$ est la matrice $(\chi^{\Lambda}(q,T_{M})=\chi^{\Lambda}(q,M))_{M \in \ym_{n,2}}$. On pose alors :
$$\varSigma_{M,q}(\Lambda)=\begin{cases}
2^{|M|-\ell(M)}\,n^{\downarrow |M|}\,\chi^{\Lambda}(q,M)&\text{si }n=|\Lambda|\geq |M|,\\
0&\text{sinon}.
\end{cases}$$
Pour $q=1$, on retrouve les caractères centraux signés $\varSigma_{M}$. Pour $q$ générique, il n'est pas clair que les symboles $\varSigma_{M,q}$ forment une base de $\obs^{\mB}$, et d'ailleurs c'est seulement une conjecture (nous donnerons plus loin les formules de changement de base pour $n \leq 2$). En effet, on ne connaît pas en type B de formule de Frobenius pour les $q$-caractères, ou de formule de changement de base entre caractères et $q$-caractères. Pour cette raison, la technique des observables de bidiagrammes ne donnera pour l'instant que des résultats partiels (essentiellement, seulement des résultats découlant du comportement asymptotique en type A).\bigskip

\begin{lemma}[Restriction d'un module irréductible de $\wsym_{n}$ à $\sym_{n}$]\label{blackstar}
Soit $\Lambda=(\lambda^{(1)},\lambda^{(2)})$ une bipartition, et $V^{\Lambda}$ le $\C\wsym_{n}$-module irréductible correspondant. On note $c_{\lambda^{(1)}\lambda^{(2)}}^{\nu}$ les coefficients de Littlewood-Richardson du produit de fonctions de Schur $s_{\lambda^{(1)}}\,s_{\lambda^{(2)}}$, c'est-à-dire que 
$$ s_{\lambda^{(1)}}(X)\,s_{\lambda^{(2)}}(X)=\sum_{\nu} \,c_{\lambda^{(1)}\lambda^{(2)}}^{\nu}\, s_{\nu}(X)\,.$$
Le module restreint $\mathrm{Res}_{\sym_{n}}^{\wsym_{n}} (V^{\Lambda})$ admet pour décomposition en irréductibles :
$$\mathrm{Res}_{\sym_{n}}^{\wsym_{n}} (V^{\Lambda}) \simeq_{\C\sym_{n}} \bigoplus_{\nu}\,c_{\lambda^{(1)}\lambda^{(2)}}^{\nu}\,V^{\nu}\,.$$
Le résultat reste vrai pour les modules sur les algèbres d'Hecke $\IH_{q}(\wsym_{n})$ et $\IH_{q}(\sym_{n})$.
\end{lemma}
\begin{proof}
Les règles de branchement pour les modules de groupes de Coxeter et les modules de leurs algèbres d'Hecke sont identiques, voir par exemple \cite[chapitre 7]{GP00}. Il suffit donc de démontrer le résultat pour les groupes $\wsym_{n}$ et $\sym_{n}$. Le caractère irréductible associé à $\Lambda$ et évalué en une permutation (non signée) $\sigma \in \sym_{n}$ de type cyclique $\mu$ s'écrit 
$$\left(\mathrm{Res}_{\sym_{n}}^{\wsym_{n}}(\varsigma^{\Lambda})\right)(\mu)=\scal{s_{\lambda^{(1)}}(X)\,s_{\lambda^{(2)}}(Y)}{p_{\mu}(X+Y)}\,.$$
Rempla\c cons provisoirement le produit de fonctions de Schur $s_{\lambda^{(1)}}(X)\,s_{\lambda^{(2)}}(Y)$ par un produit de fonctions puissances $p_{\lambda^{(1)}}(X)\,p_{\lambda^{(2)}}(Y)$, et calculons le produit scalaire
$$S=\scal{p_{\lambda^{(1)}}(X)\,p_{\lambda^{(2)}}(Y)}{p_{\mu}(X+Y)}\,.$$
On note $\mu=(\mu_{1},\ldots,\mu_{r})$, et pour toute partie $I \subset \lle 1,r\rre$, on note $\mu_{I}$ la partition obtenue en prenant les parts $\mu_{i}$ avec $i \in I$, et $\mu_{\setminus I}$ la partition obtenue en prenant les parts $\mu_{i}$ avec $i \notin I$. Alors, en développant la fonction $p_{\mu}(X+Y)$, on voit que 
$$S=\sum_{I \subset \lle 1,r\rre} \scal{p_{\lambda^{(1)}}(X)\,p_{\lambda^{(2)}}(Y)}{p_{\mu_{I}}(X)\,p_{\mu\setminus I}(Y)}=\sum_{I \subset \lle 1,r\rre} \mathbb{1}_{\lambda^{(1)}=\mu_{I}}\mathbb{1}_{\lambda^{(2)}=\mu_{\setminus I}}\,z_{\lambda^{(1)}}\,z_{\lambda^{(2)}}\,.$$
Si la partition $\lambda^{(1)}\sqcup \lambda^{(2)}$ n'est pas égale à $\mu$, tous les termes de la somme sont nuls. Sinon, il n'est pas difficile de voir que le nombre de parties $I$ tel que $\lambda^{(1)}=\mu_{I}$ et $\lambda^{(2)}=\mu_{\setminus I}$ est $z_{\mu}/(z_{\lambda^{(1)}}\,z_{\lambda^{(2)}})$, de sorte que $S=z_{\mu}$. Dans tous les cas, on a donc 
$$S=\scal{p_{\lambda^{(1)}\sqcup\lambda^{(2)}}}{p_{\mu}}=\scal{p_{\lambda^{(1)}}\,p_{\lambda^{(2)}} }{p_{\mu}}\,.$$
Maintenant, comme les $p_{\lambda}$ forment une base de l'algèbre $\Lambda$, pour toutes fonctions symétriques $h$, $k$ et $l$, 
$$\scal{h(X)\,k(Y)}{l(X+Y)}_{\Lambda(X) \otimes \Lambda(Y)}=\scal{hk}{l}_{\Lambda}\,.$$
En prenant $h=s_{\lambda^{(1)}}$, $k=s_{\lambda^{(2)}}$ et $l=p_{\mu}$, on obtient donc :
\begin{align*}\left(\mathrm{Res}_{\sym_{n}}^{\wsym_{n}}(\varsigma^{\Lambda})\right)(\mu)&=\sum_{\nu}\,c_{\lambda^{(1)}\lambda^{(2)}}^{\nu}\scal{s_{\nu}}{p_{\mu}}\\
&=\left(\sum_{\nu}\,c_{\lambda^{(1)}\lambda^{(2)}}^{\nu}\, \varsigma^{\nu}\right)(\mu)\,.
\end{align*}
Comme les caractères caractérisent les modules sous-jacents, ceci conclut la preuve du lemme.
\end{proof}\medskip

\begin{proposition}[Quantification partielle de la base des caractères centraux signés]\label{quantobsb}
Pour toute partition $\mu$, le symbole $\varSigma_{(\mu,\emptyset),q}$ est une observable de bidiagrammes dans $\obs^{\mB}$, et plus précisément une combinaison linéaire d'observables $\varSigma_{\nu,\emptyset}$. La relation de changement de base entre ces observables de bidiagrammes est essentiellement la même qu'en type A :
$$\forall \mu \in\ym_{k},\,\,\,(q-1)^{\ell(\mu)}\frac{\varSigma_{(\mu,\emptyset),q}}{2^{|\mu| - \ell(\mu)}}=\sum_{\rho \in \ym_{k}}\frac{\scal{h_{\mu}}{p_{\rho}}}{\scal{p_{\rho}}{p_{\rho}}}\,(q^{\rho}-1)\,\frac{\varSigma_{\rho,\emptyset}}{2^{|\rho|-\ell(\rho)}}\,.$$
\end{proposition}
\begin{proof}
Il suffit de montrer que pour toute bipartition $\Lambda$, la formule donnée dans l'énoncé et évaluée en $\Lambda$ est vraie. Si $|\Lambda|=n<k$, c'est évident. Sinon, on peut écrire :
\begin{align*}
(q-1)^{\ell(\mu)}\frac{\varSigma_{(\mu,\emptyset),q}(\Lambda)}{2^{|\mu| - \ell(\mu)}}&=(q-1)^{\ell(\mu)}\,n^{\downarrow |\mu|}\,\frac{\varsigma^{\Lambda}(q,(\mu\sqcup 1^{n-|\mu|},\emptyset))}{\dim \Lambda}\\
&=\sum_{\nu} \,c_{\lambda^{(1)}\lambda^{(2)}}^{\nu}\,\frac{\dim \nu}{\dim\Lambda}\left((q-1)^{\ell(\mu)}\,n^{\downarrow |\mu|}\,\frac{\varsigma^{\nu}(q,(\mu\sqcup 1^{n-|\mu|},\emptyset))}{\dim \nu} \right).
\end{align*}
À $\nu$ fixé, le terme entre parenthèses dans la somme est $(q-1)^{\ell(\mu)}\,\varSigma_{\mu}(\nu)$, et on peut utiliser la formule de changement de base en type A donnée par la proposition \ref{quantobs}. L'expression que l'on souhaite calculer devient donc 
$$\sum_{\nu} \,c_{\lambda^{(1)}\lambda^{(2)}}^{\nu}\,\frac{\dim \nu}{\dim\Lambda}\left(\sum_{\rho \in \ym_{k}} \frac{\scal{h_{\mu}}{p_{\rho}}}{\scal{p_{\rho}}{p_{\rho}}}\,(q^{\rho}-1)\,n^{\downarrow |\rho|}\,\frac{\varsigma^{\nu}(q,\rho \sqcup1^{n-|\rho|})}{\dim \nu}\right),$$
et en intervertissant les sommes, on retrouve 
$$\sum_{\rho \in \ym_{k}}\frac{\scal{h_{\mu}}{p_{\rho}}}{\scal{p_{\rho}}{p_{\rho}}}\,(q^{\rho}-1)\,n^{\downarrow |\rho|}\,\frac{\varsigma^{\Lambda}(q,(\rho \sqcup1^{n-|\rho|},\emptyset))}{\dim \Lambda}\,,$$
c'est-à-dire le terme de droite dans la formule de l'énoncé. Notons que la forme exacte donnée par le lemme \ref{blackstar} de la restriction d'un $\wsym_{n}$-module irréductible en un $\sym_{n}$-module n'a pas réellement servi : ainsi, on aurait pu écrire la preuve sans connaître précisément les coefficients $c_{\lambda^{(1)}\lambda^{(2)}}^{\mu}$.
\end{proof}
\bigskip
\bigskip

Comme dans le chapitre \ref{qplancherelmeasure}, on envisage maintenant les observables de bi\-dia\-grammes et les symboles $\varSigma_{M,q}$ (qu'il s'agisse ou non d'observables dans $\obs^{\mB}$) comme des variables aléatoires sous la B-$q$-mesure de Plancherel $M_{n,q}^{\mB}$. Comme en type A, cette B-$q$-mesure de Plancherel peut être interprétée comme le poids des $q$-caractères irréductibles normalisés $\chi^{\Lambda}(q)$ dans la trace canonique de l'algèbre d'Hecke $\IH_{q}(\wsym_{n})$ :
$$\tau(T_{\sigma})=\mathbb{1}_{(\sigma=\id_{\lle 1,n\rre})}\,.$$
Par conséquent, l'espérance d'un symbole $\varSigma_{M,q}$ est :
$$\esper[\varSigma_{M,q}]=\begin{cases} n^{\downarrow k} &\text{si }M=(\emptyset, 1^{k}),\\
0&\text{sinon}.
\end{cases}$$
Les mêmes raisonnements qu'en type A permettent dès lors d'établir :
\begin{theorem}[Asymptotique de la B-$q$-mesure de Plancherel]\label{asymptoticbqplancherel}
On suppose $q<1$. Soit $\Lambda=(\lambda^{(1)},\lambda^{(2)})$ une bipartition aléatoire tirée sous la B-$q$-mesure de Plancherel $M_{n,q}^{\mB}$, et $\lambda=\lambda^{(1)}\sqcup \lambda^{(2)}$ la partition obtenue en réordonnant les parts de $\lambda^{(1)}$ et les parts de $\lambda^{(2)}$. Pour tout $i\geq 1$,
$$\frac{\lambda_{i}}{n} \to_{M_{n,q}^{\mB}} (1-q)\,q^{i-1}\,.$$
De plus, les déviations des parts renormalisées $\lambda_{i}/n$ sont gaussiennes, et elles suivent les mêmes lois qu'en type A (voir le théorème \ref{secondasymptoticqplancherel}).
\end{theorem}\bigskip

La preuve reprend essentiellement les arguments du chapitre \ref{qplancherelmeasure}.  Ainsi, en inversant la formule de changement de base donnée par la proposition \ref{quantobsb}, on obtient une expression exacte pour l'espérance de $\varSigma_{\mu,\emptyset}$ :
$$\esper\left[\frac{\varSigma_{\mu,\emptyset}}{2^{|\mu|-\ell(\mu)}}\right]=n^{\downarrow |\mu|}\,\frac{(1-q)^{|\mu|}}{1-q^{\mu}}\,.$$
Comme en type A\footnote{Pour démontrer ceci, il convient d'utiliser la remarque générale faite après la démonstration du lemme \ref{stilldre}. En particulier, on ne peut pas adapter la preuve du lemme \ref{stilldre}, car on ne connaît pas \emph{a priori} de formule générale pour $\esper[\varSigma_{M}]$.}, l'ordre de grandeur de l'espérance d'une observable de degré $d$ est un $O(n^{d})$. Par conséquent, on peut remplacer dans les estimations une observable de bidiagrammes par une autre observable de même composante de plus haut degré, et compte tenu de la proposition \ref{biobs} :
$$\frac{1}{n^{|\mu|}}\,\esper\left[\prod_{i=1}^{\ell(\mu)} p_{\mu_{i},\emptyset}+p_{\emptyset,\mu_{i}}\right] = \frac{(1-q)^{|\mu|}}{1-q^{\mu}}+O(n^{-1})\,.$$
La même preuve qu'en type A assure donc la convergence en probabilité
$$\frac{p_{k,\emptyset}(\Lambda)+p_{\emptyset,k}(\Lambda)}{n^{k}}=\frac{p_{k}(\lambda^{(1)})+p_{k}(\lambda^{(2)})}{n^{k}} \to_{M_{n,q}^{\mB}} \frac{(1-q)^{k}}{1-q^{k}}\,.$$ 
Or, $p_{k}(\lambda=\lambda^{(1)}\sqcup \lambda^{(2)})$ est quasiment égal à $p_{k}(\lambda^{(1)})+p_{k}(\lambda^{(2)})$ ; en effet, les coordonnées de Frobenius de $\lambda$ sont à un décalage près obtenues en prenant la réunion des coordonnées de Frobenius de $\lambda^{(1)}$ et des coordonnées de Frobenius de $\lambda^{(2)}$. Il est donc aisé d'adapter la preuve du théorème \ref{firstasymptoticqplancherel} pour montrer que 
$$\frac{\lambda_{i}}{n} \to_{M_{n,q}^{\mB}} (1-q)\,q^{i-1}$$
pour tout entier $i$ plus grand que $1$. Nous laissons les détails techniques au soin du lecteur ; il n'y a aucune subtilité supplémentaire par rapport au cas traité dans le chapitre \ref{qplancherelmeasure}. De même, la déviation gaussienne s'obtient en démontrant les exacts analogues des lemmes \ref{magnitudedj}, \ref{magnitude}, \ref{productcentralcharacter} et \ref{devgaussianobs}. Par rapport au type A, la seule légère différence est la description du produit de deux observables $\varSigma_{M_{1}}\varSigma_{M_{2}}$ à l'aide d'appariements. En effet, le type du produit d'un élément $(\{\eps_{i,j},a_{i,j}\}_{i,j})$ de type $M_{1}$ avec un élément $(\{\theta_{k,l},b_{k,l}\}_{k,l})$ de type $M_{2}$ ne dépend pas seulement des éventuelles égalités $a_{i,j}=b_{k,l}$, mais aussi des signes $\eps_{i,j}$ et $\theta_{k,l}$. Par exemple, dans le produit $\varSigma_{\emptyset,2}\,\varSigma_{\emptyset,2}$, les deux termes
$$\big(\overline{1},(1,2)\big)\times \big(\overline{1},(1,2)\big) \qquad\text{et}\qquad \big(\overline{1},(1,2)\big)\times \big(\overline{2},(1,2)\big)$$
correspondent au même appariement de taille $2$ des $a_{i,j}$ avec les $b_{k,l}$, mais le premier produit vaut
$$\big(\overline{1}\,\overline{2},(1)(2)\big)$$
et est de type cyclique $(\emptyset,1^{2})$, tandis que le second produit vaut
$$(1)(2)$$
et est de type cyclique $(1^{2},\emptyset)$. Ainsi, l'appariement partiel détermine bien la partition $\mu=\mu^{(1)}\sqcup \mu^{(2)}$ dont les parts sont les tailles des cycles pairs et impairs de la permutation signée produit (ici, $\mu=1^{2}$), mais il ne détermine pas \emph{a priori} la parité des cycles correspondants : les parts de $\mu$ peuvent tomber dans $\mu^{(1)}$ ou $\mu^{(2)}$. Par conséquent, l'analogue du lemme \ref{productcentralcharacter} en type B s'énonce comme suit : le produit $\varSigma_{M}\,\varSigma_{N}$ de deux caractères centraux est égal à la somme $\sum_{S} \,\varSigma_{R(S)}$ de caractères centraux de type B, où $S$ parcourt l'ensemble des \textbf{appariements partiels signés} entre l'ensemble $I_{A}$ des indices des $a$ (les entiers apparaissant dans les cycles des permutations partielles signées de $\varSigma_{M}$) et l'ensemble $I_{B}$ des indices des $b$ (les entiers apparaissant dans les cycles des permutations partielles signées de $\varSigma_{N}$). Ici, par appariement partiel signé, on entend un appariement partiel et la donnée d'un nombre suffisant de signes pour distinguer les cycles pairs des cycles impairs. Nous ne préciserons pas plus cette définition, mais nous allons traiter un exemple en détail.\medskip

\begin{example} Détaillons le calcul du produit $\varSigma_{(2),\emptyset}\,\varSigma_{\emptyset,(3)}$. L'ensemble d'indices $I_{A}$ est $\{1,2\}$, et l'ensemble d'indices $I_{B}$ est $\{1',2',3'\}$.\vspace{2mm}
\begin{enumerate}
\item Si les arrangements $A=(a_{1},a_{2})$ et $B=(b_{1},b_{2},b_{3})$ sont disjoints, alors $$(\eps_{A},(a_{1},a_{2}))\,(\eps_{B},(b_{1},b_{2},b_{3}))$$ est le produit disjoint d'un $2$-cycle pair par un $3$-cycle impair, donc est de type cyclique $((3),(2))$. L'appariement partiel vide correspond donc au caractère central $\varSigma_{(2),(3)}$.\vspace{2mm}
\item Si $A$ et $B$ ont un point en commun, on peut à symétrie près supposer que c'est $a_{1}=b_{1}$. Alors, 
$$\big(\eps_{a_{1}}\eps_{b_{1}},(a_{1},b_{1})\big)\times \big(\theta_{b_{1}}\theta_{b_{2}}\theta_{b_{3}},(b_{1},b_{2},b_{3})\big)=\big(\phi_{a_{1}}\phi_{b_{1}}\phi_{b_{2}}\phi_{b_{3}},(a_{1},b_{1},b_{2},b_{3})\big)$$
est un $4$-cycle, avec les signes $\phi_{a_{1}}=\eps_{a_{1}}$, $\phi_{b_{1}}=\theta_{b_{1}}$, $\phi_{b_{2}}=\theta_{b_{2}}$ et $\phi_{b_{3}}=\eps_{b_{1}}\theta_{b_{3}}$. Le produit des signes le long du cycle est donc 
$$(\eps_{a_{1}}\eps_{b_{1}})\,(\theta_{b_{1}}\theta_{b_{2}}\theta_{b_{3}})=(-1)\times(1)=-1\,,$$
\emph{i.e.}, on a forcément un $4$-cycle impair. Les $6$ appariements partiels de taille $1$ apportent donc une contribution $6\,\varSigma_{\emptyset,(4)}$.\vspace{2mm}
\item Finalement, si $A$ et $B$ ont deux points en commun, on peut à symétrie près supposer que c'est $a_{1}=b_{1}$ et $a_{2}=b_{2}$. Alors, 
$$\big(\eps_{b_{1}}\eps_{b_{2}},(b_{1},b_{2})\big)\times \big(\theta_{b_{1}}\theta_{b_{2}}\theta_{b_{3}},(b_{1},b_{2},b_{3})\big)=\big(\phi_{b_{1}}\phi_{b_{2}}\phi_{b_{3}},(b_{1})(b_{2},b_{3})\big)$$
avec les signes $\phi_{b_{1}}=\eps_{b_{2}}\theta_{b_{1}}$, $\phi_{b_{2}}=\theta_{b_{2}}$ et $\phi_{b_{3}}=\eps_{b_{1}}\theta_{b_{3}}$. Comme le produit de tous les signes est impair, le type cyclique est donc $((1),(2))$ ou $((2),(1))$. Pour déterminer quel type on obtient, il faut par exemple connaître simultanément $\eps_{b_{2}}$ et $\theta_{b_{1}}$ ; on connaît alors la parité du $1$-cycle, et donc le type cyclique global. Les $6$ appariements partiels de taille $2$ apportent donc $6 \times 2^{2}=24$ termes, qui se répartissent en $12\,\varSigma_{(1),(2)} + 12\,\varSigma_{(2),(1)}$.
\end{enumerate}\bigskip

\noindent Ainsi, $\varSigma_{(2),\emptyset}\,\varSigma_{\emptyset,(3)}=\varSigma_{(2),(3)}+6\,\varSigma_{\emptyset,(4)}+12\,\varSigma_{(2),(1)}+12\,\varSigma_{(1),(2)}$. Notons que si l'on projette l'algèbre $\blg_{\infty}^{\mB}$ sur $\blg_{\infty}$ en associant à une permutation partielle signée $((\eps,\sigma),S)$ la permutation partielle $(\sigma,S)$, alors l'image d'un caractère central de type B $\varSigma_{\mu^{(1)},\mu^{(2)}}$ par cette projection est le caractère central de type A $$2^{|\mu^{(1)}|+|\mu^{(2)}|-\ell(\mu^{(1)})-\ell(\mu^{(2)})}\,\varSigma_{\mu^{(1)}\sqcup \mu^{(2)}}.$$ L'identité que l'on vient d'obtenir est donc compatible avec la relation établie en exemple après le lemme \ref{productcentralcharacter} : $\varSigma_{2}\,\varSigma_{3}=\varSigma_{3,2}+6\,\varSigma_{4}+6\,\varSigma_{2,1}$. 
\end{example}\bigskip\bigskip

Comme en type A, dans un produit $\varSigma_{M}\,\varSigma_{N}=\sum_{S} \varSigma_{R(S)}$, le degré d'un terme $\varSigma_{R(S)}$ correspondant à un appariement partiel signé $S$ est $|M|+|N|-|S|$. De plus, lorsqu'on effectue un produit de caractères centraux de cycles pairs $\varSigma_{(l),\emptyset}$ et $\varSigma_{(m),\emptyset}$, tous les appariements partiels de taille $1$ donnent un caractère central $\varSigma_{(l+m-1),\emptyset}$. En effet, étant donnés deux cycles pairs
$$\big(\eps_{a_{1}}\cdots\eps_{a_{l}},(a_{1},\ldots,a_{l})\big)\qquad\text{et}\qquad\big(\theta_{b_{1}}\cdots\theta_{b_{m}},(b_{1},\ldots,b_{m})\big)$$
s'intersectant en un seul point $a_{l}=b_{1}$, le produit des cycles est 
$$\big(\phi_{a_{1}}\cdots\phi_{a_{l-1}}\phi_{b_{1}}\cdots\phi_{b_{m}},(a_{1},\cdots,a_{l-1},b_{1},\cdots,b_{m})\big)$$
avec les signes $\phi_{a_{i\leq l-1}}=\eps_{a_{i}}$, $\phi_{b_{j\leq m-1}}=\theta_{b_{j}}$ et $\phi_{b_{m}}=\eps_{a_{l}}\theta_{b_{m}}$. Comme les deux cycles étaient pairs, le produit des signes du cycle produit est encore pair, et le type cyclique est bien  $((l+m-1),\emptyset)$. En sommant sur les $ml$ appariements partiels de taille $1$, on obtient donc le développement :
$$\varSigma_{(l,\emptyset)}\,\varSigma_{(m),\emptyset}=\varSigma_{(l,m),\emptyset}+ml\,\varSigma_{(l+m-1),\emptyset}+(\text{termes de degré inférieur à }l+m-2)\,.$$ 
Compte tenu de la proposition \ref{biobs}, les observables intéressantes pour l'asymptotique des lignes de la bipartition sont les $\varSigma_{\mu,\emptyset}/2^{|\mu|-\ell(\mu)}$. Or, en divisant l'identité précédente par $2^{l+m-2}$, on obtient :
$$\frac{\varSigma_{(l),\emptyset}}{2^{l-1}}\,\frac{\varSigma_{(m),\emptyset}}{2^{m-1}}=\frac{\varSigma_{(l,m),\emptyset}}{2^{l+m-2}}+ml\,\frac{\varSigma_{(l+m-1),\emptyset}}{2^{l+m-2}}+(\text{termes de degré inférieur à }l+m-2)\,,$$ 
ce qui est l'exact analogue du développement obtenu dans le chapitre \ref{qplancherelmeasure} pour les produits de caractères centraux de cycles. On peut donc raisonner exactement comme dans le paragraphe \ref{qgaussian}, c'est-à-dire calculer les cumulants joints des caractères centraux des cycles pairs (on a les mêmes estimations qu'en type A), établir leur déviation gaussienne et en déduire la déviation gaussienne des lignes en associant à chaque bipartition une mesure de probabilité sur $[-1,1]$, et en montrant le caractère gaussien de ces mesures aléatoires. De nouveau, les détails sont laissés au lecteur ; le seul point qui nous semblait nouveau et non évident était le développement du produit des caractères centraux de type B démontré ci-dessus.

\section[Caractères des algèbres d'Hecke de type B et répartitions des parts]{Caractères des algèbres d'Hecke de type B et répartition\\ des parts}\label{heckebmix} 

Dans tout ce paragraphe, on suppose de nouveau $q<1$, l'autre cas s'en déduisant par une symétrie des bidiagrammes. Le théorème \ref{asymptoticbqplancherel} assure que dans l'algèbre d'observables $\obs^{\mB}$, un sous-espace d'observables converge en probabilité après renormalisation par $n^{\deg(\cdot)}$, et ces observables renormalisées sont concentrées gaussiennement autour de leurs valeurs moyennes. Ce <<~coeur gaussien~>> permet de comprendre l'asymptotique des tailles des parts d'une bipartition sous la B-$q$-mesure de Plancherel. Ceci étant, si l'on souhaite maintenant comprendre précisément quelles parts tombent dans quelle partie de la bipartition, il est nécessaire de sortir du coeur gaussien et de considérer par exemple des observables du type :
$$\frac{1}{2^{k-1}}(\varSigma_{k,\emptyset}(\Lambda)+\varSigma_{\emptyset,k}(\Lambda))=\varSigma_{k}(\lambda^{(1)})\sim p_{k}(\lambda^{(1)})$$
Il est donc tentant d'établir une formule générale pour l'espérance d'un caractère central quelconque $\varSigma_{M}$. Or, les seules <<~observables~>> dont on sait calculer les espérances sous la B-$q$-mesure de Plancherel sont pour l'instant les symboles quantifiés $\varSigma_{M,q}$ ; comme dans le paragraphe \ref{quantum}, on doit donc établir des relations de changement de base entre les $\varSigma_{M}$ et les $\varSigma_{M,q}$ qui généralisent celles données par la proposition \ref{quantobsb}. Ceci revient à décomposer dans la base des $p_{M}(X,Y)$ des fonctions multi-symétriques $q_{M}(X,Y,q)$ telles que :
$$\forall \Lambda,M,\,\,\,\varsigma^{\Lambda}(q,M)=\scal{s_{\Lambda}(X,Y)}{q_{M}(X,Y,q)}_{\Lambda(X)\otimes \Lambda(Y)}$$
Par la formule de Frobenius pour les caractères de $\wsym_{n}$, la spécialisation en $q=1$ des fonctions $q_{M}$ est $q_{M}(X,Y,q)=p_{M}(X,Y)$, et d'autre part, la proposition \ref{quantobsb} peut être réinterprétée en :
$$q_{\mu,\emptyset}(X,Y,q)=q_{\mu}(X+Y,q)\,,$$
où les fonctions $q_{\mu}$ sont celles introduites dans le paragraphe \ref{ram}.\bigskip
\bigskip

Pour calculer en pratique les fonctions $q_{M}$, on peut utiliser la formule de Murna\-ghan-Nakayama pour les caractères d'algèbres d'Hecke $\IH_{q}(\wsym_{n})$ ; une telle formule est démontrée par A. Ram et T. Halverson dans \cite{HR96}, et on renvoie également à \cite[chapitre 10]{GP00}. Ainsi, si $M=(\mu^{(1)},\mu^{(2)})$ avec $\ell(\mu^{(1)})=r$ et $\ell(\mu^{(2)})=s$, alors $\varsigma^{\Lambda}(q,M)$ est égal à 
$$\sum_{(\emptyset,\emptyset) = \Lambda_{0} \subset \Lambda_{1}\subset \cdots \subset\Lambda_{r+s}=\Lambda}\Delta(\Lambda_{1}\setminus \Lambda_{0})\cdots \Delta(\Lambda_{r}\setminus \Lambda_{r-1})\,\overline{\Delta}(\Lambda_{r+1}\setminus \Lambda_{r})\cdots \overline{\Delta}(\Lambda_{r+s}\setminus \Lambda_{r+s-1})\,,$$
la somme étant effectuée sur les suites de bipartitions telles que $\Lambda_{i}\setminus \Lambda_{i-1}$ soit composée de deux rubans (éventuellement non connexes) de poids total $\mu^{(1)}_{i}$ pour $i \leq r$, et $\widetilde{\mu}^{(2)}_{i-r}$ pour $i>r$. Les poids $\Delta$ et 
$\overline{\Delta}$ s'expriment en fonction des formes des rubans, voir \cite[théorème 2.22]{HR96} ; nous ne rentrerons pas dans les détails\footnote{On notera que dans l'article \cite{HR96}, la présentation de l'algèbre d'Hecke de type B est différente de celle que nous donnons page \pageref{algheckeb}, et l'on doit modifier en conséquence toutes les formules. D'autre part, la définition des coins \emph{sharp} et \emph{dull} donnée dans l'article ne tient pas compte des cas pathologiques, par exemple lorsque le ruban considéré est réduit à une case ou à une seule ligne ; nous avons du refaire tous les calculs compte tenu de ces imprécisions.}, et nous nous contenterons de donner les tables de caractères qui se déduisent de ces formules pour $n=1,2$, voir la figure \ref{tableofcharacterheckeb}.\bigskip

\begin{figure}[ht]
\begin{center}
\normalsize{
\begin{tabular}{|c|c|c|}
\hline $\IH_q(\wsym_{1})$ & $(1),\emptyset$ & $\emptyset,(1)$ \\ 
\hline $(1),\emptyset$ & $1$ & $q$ \\
\hline $\emptyset,(1)$ & $1$ & $-1$ \\
\hline
\end{tabular}\bigskip
\bigskip

\begin{tabular}{|c|c|c|c|c|c|}
\hline $\IH_q(\wsym_{2})$ & $(1^{2}),\emptyset$ & $(1),(1)$ &  $\emptyset,(1^{2})$ & $(2),\emptyset$  & $\emptyset,(2)$ \\ 
\hline $(1^{2}),\emptyset$ &     $1$ &      $q$ &    $q^2$  &     $-1$ &     $-q$ \\ 
\hline $(1),(1)$                  &     $2$ & $q - 1$ & $-2q^2$ &  $q - 1$ &     $0$  \\ 
\hline $\emptyset,(1^{2})$ &      $1$ &     $-1$ &        $1$ &      $-1$ &     $1$ \\
\hline $(2),\emptyset$        &      $1$ &      $q$ &    $q^4$ &      $q$ &    $q^2$\\ 
\hline $\emptyset,(2)$         & $1$ &     $-1$ &    $q^2$ &      $q$ &     $-q$ \\ 
\hline 
\end{tabular}\bigskip
}
\caption[Tables des caractères des algèbres d'Hecke $\IH_{q}(\wsym_{n})$ pour $n \leq 2$]{Tables des caractères des algèbres d'Hecke $\IH_{q}(\wsym_{n})$ pour $n \leq 2$.\label{tableofcharacterheckeb}}
\end{center}
\end{figure}

La fonction $q_{M}(X,Y,q)$ est égale à $\sum_{\Lambda \in \ym_{n,2}} \varsigma^{\Lambda}(q,M)\,s_{\lambda}(X,Y)$. En décomposant $q_{M}(X,Y,q)$ dans la base $(p_{M}(X,Y))_{M}$, on obtient une expression des $q$-caractères de $\IH_{q}(\wsym_{n})$ en fonction des caractères de $\wsym_{n}$, qu'on donne ici pour les caractères renormalisés de degré inférieur à $2$ :
\begin{align*}
&q_{(1),\emptyset}(X,Y,q)= p_{1}(X+Y)\qquad\,\,\,\,;\qquad q_{\emptyset,(1)}(X,Y,q)= p_{1}(qX-Y)\quad;\\
&q_{(1^{2}),\emptyset}(X,Y,q)= p_{1^{2}}(X+Y)\qquad;\qquad q_{(1),(1)}(X,Y,q)= p_{1}(qX-Y)\,p_{1}(X+Y)\quad; \\
&q_{\emptyset,(1^{2})}(X,Y,q)=\frac{q^{2}+1}{2}\, p_{1^{2}}(qX-Y)+q(q-1)^{2}\,p_{1}(X)\,p_{1}(Y)+\frac{q^{2}-1}{2}\,p_{2}(qX-Y)\quad;\\
&q_{(2),\emptyset}(X,Y,q)=\frac{q-1}{2}\,p_{1^{2}}(X+Y)+\frac{q+1}{2}\,p_{2}(X+Y)\quad;\\
&q_{\emptyset,(2)}(X,Y,q)=\frac{q-1}{2}\,(q\,p_{1^{2}}(X)-p_{1^{2}}(Y))+\frac{q+1}{2}\,(q\,p_{2}(X)-p_{2}(Y))\\
&\quad\Rightarrow\quad \varSigma_{((1),\emptyset),q}=\varSigma_{(1),\emptyset} \qquad;\qquad  \varSigma_{(\emptyset,(1)),q}=\frac{q-1}{2}\,\varSigma_{(1),\emptyset}+\frac{q+1}{2}\,\varSigma_{\emptyset,(1)}\quad;\\
&\phantom{\quad\Rightarrow\quad}\varSigma_{((1^{2}),\emptyset),q}=\varSigma_{(1^{2}),\emptyset}\qquad;\qquad\varSigma_{((1),(1)),q}=\frac{q-1}{2}\,\varSigma_{(1^{2}),\emptyset}+\frac{q+1}{2}\,\varSigma_{(1),(1)}\quad;\\
&\phantom{\quad\Rightarrow\quad} \varSigma_{(\emptyset,(1^{2})),q}=\frac{(q^{2}-1)^{2}}{8}\,\varSigma_{(1^{2}),\emptyset}+\frac{q^{4}-1}{4}\,\varSigma_{(1),(1)}+\frac{q^{4}+6q^{2}+1}{8}\,\varSigma_{\emptyset,(1^{2})}\\
&\qquad\qquad\qquad\qquad+\frac{(q^{2}-1)^{2}}{8}\,\varSigma_{(2),\emptyset}+\frac{q^{4}-1}{8}\,\varSigma_{\emptyset,(2)}\quad;\end{align*}
\begin{align*}
&\phantom{\quad\Rightarrow\quad} \varSigma_{((2),\emptyset),q}=(q-1)\,\varSigma_{(1^{2}),\emptyset}+\frac{q+1}{2}\,\varSigma_{(2),\emptyset}\quad;\\
&\phantom{\quad\Rightarrow\quad} \varSigma_{(\emptyset,(2)),q}=\frac{(q-1)^{2}}{4}\,\varSigma_{(1^{2}),\emptyset}+ \frac{q^{2}-1}{2}\,\varSigma_{(1),(1)}+\frac{(q-1)^{2}}{4}\,\varSigma_{\emptyset,(1^{2})}\\
&\qquad\qquad\qquad\qquad+\frac{q^{2}-1}{4}\,\varSigma_{(2),\emptyset}+\frac{(q+1)^{2}}{4}\,\varSigma_{\emptyset,(2)}\,.
\end{align*}
Contrairement à ce qui se passe en type A, les fonctions multisymétriques $q_{M}(X,Y,q)$ ne factorisent pas, c'est-à-dire que 
$$q_{M}(X,Y,q) \neq \prod_{i=1}^{\ell(\mu^{(1)})} q^{(1)}_{\mu^{(1)}_{i}}(X,Y,q) \,\prod_{j=1}^{\ell(\mu^{(2)})} q^{(2)}_{\mu^{(2)}_{j}}(X,Y,q)\,.$$
Par exemple, $q_{\emptyset,(1^{2})}(X,Y,q)$ n'est pas du tout égal à $(q_{\emptyset,1}(X,Y,q))^{2}$. Comme mentionné plus haut, ceci est lié au fait que le produit de groupes hyperoctahédraux $\wsym_{m} \times \wsym_{n}$ dans $\wsym_{m+n}$ n'est pas un sous-groupe parabolique. Ainsi, on a besoin de tous les générateurs de $\wsym_{m+n}$ pour écrire les éléments de $\wsym_{m} \times \wsym_{n}$. Néanmoins, si l'on considère un produit $\wsym_{m} \times \sym_{n}$, alors ce groupe est bien un sous-groupe parabolique de $\wsym_{m+n}$ (on enlève le générateur $s_{n}$). Ceci implique une factorisation partielle des polynômes $q_{M}(X,Y,q)$ : 
$$q_{M=(\mu^{(1)},\mu^{(2)})}(X,Y,q)=q_{\mu^{(1)}}(X+Y,q)\times q_{\emptyset,\mu^{(2)}}(X,Y,q)\,.$$
où $q_{\mu^{(1)}}$ est la déformation de $p_{\mu^{(1)}}$ introduite dans le paragraphe \ref{ram}. Le problème se ramène donc au calcul des polynômes $q_{\emptyset,\mu^{(2)}}(X,Y,q)$ ; dans le cas général, on ne sait pas faire ce calcul, mais des expériences numériques nous ont néanmoins permis d'établir une conjecture intéressante pour la répartition asymptotique des parts d'une bipartition sous la $q$-mesure de Plancherel.\bigskip
\bigskip

En effet, en petit degré, on a pu inverser les relations entre $q$-caractères et caractères, c'est-à-dire exprimer tous les symboles $\varSigma_{M}$ en fonction des symboles $\varSigma_{M,q}$ de même degré\footnote{Au moins pour $q$ petit, il est évident qu'une telle inversion est possible, car la matrice de passage dans un sens est obtenue en perturbant la matrice identité. }. Comme $\esper[\varSigma_{M,q}]=0$ pour toute bipartition $M$ qui n'est pas du type $(\emptyset,(1^{k}))$, les relations présentées ci-contre ont donc permis de calculer l'espérance sous $M_{n,q}^{\mB}$ de tous les caractères centraux de degré suffisamment petit : il s'agit des coefficients d'une des lignes de la matrice de passage. Ainsi :
\begin{align*}
|M|=1\,\,\,:\qquad&\esper\left[\frac{\varSigma_{(1),\emptyset}}{n^{\downarrow 1}}\right]=1 \qquad;\qquad \esper\left[\frac{\varSigma_{\emptyset,(1)}}{n^{\downarrow 1}}\right]=\frac{(1-q)^{2}}{1-q^{2}}\quad;\\
&\\
|M|=2\,\,\,:\qquad&\esper\left[\frac{\varSigma_{(1^{2}),\emptyset}}{n^{\downarrow 2}}\right]=1\qquad;\qquad\esper\left[\frac{\varSigma_{(1),(1)}}{n^{\downarrow 2}}\right]=\frac{(1-q)^{2}}{1-q^{2}}\quad;\\
&\esper\left[\frac{\varSigma_{\emptyset,(1^{2})}}{n^{\downarrow 2}}\right]=\frac{(1-q)^{2}(1-q^{2})}{1-q^{4}}\quad;\\
&\esper\left[\frac{\varSigma_{(2),\emptyset}}{2n^{\downarrow 2}}\right]=\frac{(1-q)^{2}}{1-q^{2}}\qquad;\qquad\esper\left[\frac{\varSigma_{\emptyset,(2)}}{2n^{\downarrow 2}}\right]=\frac{(1-q)^{3}(1-q^{3})}{(1-q^{2})(1-q^{4})}\quad;\end{align*}
\begin{align*}
|M|=3\,\,\,:\qquad&\esper\left[\frac{\varSigma_{(1^{3}),\emptyset}}{n^{\downarrow 3}}\right]=1\qquad;\qquad\esper\left[\frac{\varSigma_{(1^{2}),(1)}}{n^{\downarrow 3}}\right]=\frac{(1-q)^{2}}{1-q^{2}}\quad;\\
&\esper\left[\frac{\varSigma_{(1),(1^{2})}}{n^{\downarrow 3}}\right]=\frac{(1-q)^{2}(1-q^{2})}{1-q^{4}}\quad;\\
&\esper\left[\frac{\varSigma_{\emptyset,(1^{3})}}{n^{\downarrow 3}}\right]=\frac{(1-q)^{5}(1-q^{3})(1+q+4q^{2}+q^{3}+q^{4})}{(1-q^{2})(1-q^{4})(1-q^{6})}\quad;\\
&\esper\left[\frac{\varSigma_{(2,1),\emptyset}}{2n^{\downarrow 3}}\right]=\frac{(1-q)^{2}}{1-q^{2}}
\qquad;\qquad \esper\left[\frac{\varSigma_{(1),(2)}}{2n^{\downarrow 3}}\right]=\frac{(1-q)^{3}(1-q^{3})}{(1-q^{2})(1-q^{4})} \quad;\\
&\esper\left[\frac{\varSigma_{(2),(1)}}{2n^{\downarrow 3}}\right]=\frac{(1-q)^{4}}{(1-q^{2})^{2}}\qquad;\qquad\esper\left[\frac{\varSigma_{(3),\emptyset}}{4n^{\downarrow 3}}\right]=\frac{(1-q)^{3}}{1-q^{3}}\quad;\\
&\esper\left[\frac{\varSigma_{\emptyset,(2,1)}}{2n^{\downarrow 3}}\right]=\frac{(1-q)^{5}(1-q^{3})(1+2q+q^{2}+2q^{3}+q^{4})}{(1-q^{2})(1-q^{4})(1-q^{6})}\quad;\\
&\esper\left[\frac{\varSigma_{\emptyset,(3)}}{4n^{\downarrow 3}}\right]=\frac{(1-q)^{4}(1-q^{3})(1-q^{5})}{(1-q^{2})(1-q^{4})(1-q^{6})}\,.
\end{align*}
Dès la taille $3$, on obtient pour les espérances des caractères centraux <<~purement impairs~>> des fractions rationnelles mettant en jeu des polynômes non triviaux en $q$ tels que $1+2q+q^{2}+2q^{3}+q^{4}$. Néanmoins, l'espérance des caractères centraux de cycles impairs semble prendre une expression relativement simple. Ainsi, on a pu vérifier jusqu'en taille $k=6$ l'identité suivante :
$$\esper\left[\frac{\varSigma_{\emptyset,k}}{2^{k-1}n^{\downarrow k}}\right]=(1-q)^{k}\prod_{i=1}^{k}\frac{\{2i-1\}_{q}}{\{2i\}_{q}}=(1-q)^{k} \frac{\{2k-1!!\}_{q}}{\{2k!!\}_{q}}\,.$$
Dans ce qui suit et jusqu'à la fin du paragraphe, on suppose que cette identité est vraie pour tout entier $k\geq 1$. Alors, la conjecture énoncée ci-dessous est également vraie :
\begin{conjecture}[Répartition asymptotique des parts sous la B-$q$-me\-su\-re de Plancherel]\label{bqmix}
Soit $(\lambda^{(1)},\lambda^{(2)})$ une bipartition aléatoire sous la B-$q$-mesure de Plancherel, et $\lambda=\lambda^{(1)}\sqcup \lambda^{(2)}$. La probabilité pour que la $2k+1$-ième part de $\lambda$ tombe dans $\lambda^{(1)}$ tend vers
$$c_{k}=\frac{1}{2}\left(1+\frac{(q;q^{2})_{\infty}}{(q^{2};q^{2})_{\infty}}\,\frac{(q;q^{2})_{k}}{(q^{2};q^{2})_{k}}\,q^{k}\right).$$
La probabilité pour que la $2k+2$-ième part de $\lambda$ tombe dans $\lambda^{(1)}$ tend pour sa part vers $1/2$. En particulier, la probabilité pour que la plus grande part appartienne à $\lambda^{(1)}$ tend vers 
$$c_{0}=\frac{1}{2}\left(1+\prod_{i=0}^{\infty} \frac{1-q^{2i+1}}{1-q^{2i+2}}\right).$$
\end{conjecture}\medskip

\begin{example}
Pour $q=1/2$, $c_{0}\simeq 0.805$, ce qui est très proche de la statistique observée --- $0.793$ --- lors des $10000$ tirages de bipartitions de taille $200$ sous la mesure $M_{n,q}^{\mB}$. 
\end{example}\bigskip

\begin{proof}[Preuve de la conjecture \ref{bqmix} en admettant la formule pour les espérances des cycles impairs]~

\noindent Nous utiliserons une nouvelle identité pour les symboles de Pochhammer :
$$\frac{(tx;q)_{\infty}}{(x;q)_{\infty}}=\sum_{n=0}^{\infty} \frac{(t;q)_{n}}{(q;q)_{n}}\, x^{n}\,,$$
voir \cite{CK02}. En utilisant la proposition \ref{biobs}, on voit que 
$$\frac{p_{k}(\lambda^{(1)})}{n^{k}}= \frac{\varSigma_{k,\emptyset}(\lambda)+\varSigma_{\emptyset,k}(\lambda)}{(2n)^{k}}+O(n^{-1})\,,$$
le $O$ étant uniforme sur l'ensemble des bipartitions. Comme dans le chapitre \ref{qplancherelmeasure}, associons à $\lambda^{(1)}$ une mesure (aléatoire) positive sur le segment $[0,1]$ :
$$X_{\lambda^{(1)}}=\sum_{i=1} \frac{\lambda^{(1)}_{i}}{n}\,\delta_{\frac{\lambda^{(1)}_{i}}{n}}\,.$$
On considérera également les mesures $X_{\lambda^{(2)}}$ et $X_{\lambda}=X_{\lambda^{(1)}}+X_{\lambda^{(2)}}$. Cette dernière est une mesure de probabilité, et le théorème \ref{asymptoticbqplancherel} peut être reformulé en disant que $X_{\lambda}$ converge en probabilité dans l'espace (métrisable) des mesures de probabilité sur le segment $[0,1]$ vers 
$$X_{\infty,q}=\sum_{i=0}^{\infty} (1-q)q^{i}\,\delta_{(1-q)q^{i}}\,.$$
Notons $B_{k}$ la variable de Bernoulli définie de la fa\c con suivante : $B_{k}=1$ si la $k$-ième part de $\lambda$ tombe dans $\lambda^{(1)}$, et $B_{k}=0$ dans le cas contraire. Cette définition est ambigu\"e si $\lambda$ possède plusieurs parts $\lambda_{k+1}=\lambda_{k+2}=\cdots=\lambda_{k+r+s}$ de même taille, avec $r$ parts de ce type tombant dans $\lambda^{(1)}$ et $s$ parts tombant dans $\lambda^{(2)}$ ; dans ce cas, on convient que 
$$B_{k+1}=B_{k+2}=\cdots=B_{k+r}=1\qquad;\qquad B_{k+r+1}=B_{k+r+2}=\cdots=B_{k+r+s}=0\,,$$
c'est-à-dire que les parts égales tombent d'abord dans $\lambda^{(1)}$ (lorsque $n$ tend vers l'infini, l'ambiguïté sera résolue, car les parts suivront une progression géométrique, donc seront avec grande probabilité distinctes). Alors, les mesures aléatoires $X_{\lambda^{(1)}}$ et $X_{\lambda^{(2)}}$ s'écrivent :
$$X_{\lambda^{(1)}}=\sum_{i=1}^{n}B_{i}\,\frac{\lambda_{i}}{n} \,\delta_{\frac{\lambda_{i}}{n}}\qquad;\qquad X_{\lambda^{(2)}}=\sum_{i=1}^{n}(1-B_{i})\,\frac{\lambda_{i}}{n} \,\delta_{\frac{\lambda_{i}}{n}}\,.$$
Par construction, $\frac{p_{k}(\lambda^{(1)})}{n^{k}}=X_{\lambda^{(1)}}(x^{k-1})+O(n^{-1})$. En passant aux espérances, on obtient donc :
\begin{align*}\esper\left[X_{\lambda^{(1)}}(x^{k-1})\right]&=\esper\left[\frac{p_{k}(\lambda^{(1)})}{n^{k}}\right]+O(n^{-1})=\esper\left[\frac{\varSigma_{k,\emptyset}(\lambda)+\varSigma_{\emptyset,k}(\lambda)}{(2n)^{k}}\right]+O(n^{-1})\\
&=\frac{(1-q)^{k}}{2}\left(\frac{1}{1-q^{k}}+\frac{\{2k-1!!\}_{q}}{\{2k!!\}_{q}}\right)+O(n^{-1})\,.
\end{align*}
Le terme $((1-q)^{k})/(1-q^{k})$ est exactement $X_{\infty,q}(x^{k-1})$. D'autre part, il est connu que la loi sur $[0,1]$ de $k$-ième moment $(2k-1!!)/(2k!!)$ est la \textbf{loi de l'arcsinus}, qu'on peut aussi voir comme la loi beta de paramètre $(1/2,1/2)$. Pour prendre en compte le terme $\{2k-1!!\}_{q}/\{2k!!\}_{q}$, on peut donc chercher un $q$-analogue de la loi de l'arcsinus. La $q^{2}$-loi beta de paramètre $(1/2,1/2)$ est la loi de probabilité sur $[0,1]$ définie par 
$$\beta_{q^{2}}=\frac{(q;q^{2})_{\infty}}{(q^{2};q^{2})_{\infty}}\,\sum_{i=0}^{\infty} \frac{(q;q^{2})_{k}}{(q^{2};q^{2})_{k}}\,q^{i}\,\delta_{q^{2i}}\,.$$
D'après la formule donnée précédemment, le $k$-ième moment de cette loi est 
$$\beta_{q^{2}}(x^{k})=\frac{(q;q^{2})_{\infty}}{(q^{2};q^{2})_{\infty}}\,\sum_{i=0}^{\infty} \frac{(q;q^{2})_{k}}{(q^{2};q^{2})_{k}}\,q^{i(2k+1)}=\frac{(q;q^{2})_{\infty}}{(q^{2};q^{2})_{\infty}}\,\frac{(q^{2k+2};q^{2})_{\infty}}{(q^{2k+1};q^{2})_{\infty}}=\frac{\{2k-1!!\}_{q}}{\{2k!!\}_{q}}\,.$$
Par conséquent, si 
$$Y_{q,\infty}=\frac{(q;q^{2})_{\infty}}{(q^{2};q^{2})_{\infty}}\,\sum_{i=0}^{\infty} \frac{(q;q^{2})_{k}}{(q^{2};q^{2})_{k}}\,(1-q)q^{3i}\,\delta_{(1-q)q^{2i}}\,,$$
alors les moyennes des moments de $X_{\lambda^{(1)}}$ valent :
$$\esper\left[X_{\lambda^{(1)}}(x^{k-1})\right]=\frac{X_{\infty,q}(x^{k-1})+Y_{\infty,q}(x^{k-1})}{2}+O(n^{-1})\,.$$
Par linéarité de l'équation et par densité des polynômes dans $\mathscr{C}([0,1])$, pour toute fonction $f$ continue sur $[0,1]$, l'identité ci-dessus reste valable :
$$\esper\left[X_{\lambda^{(1)}}(f)\right]=\frac{X_{\infty,q}(f)+Y_{\infty,q}(f)}{2}+O_{f}(n^{-1})\,.$$
Or, les $\lambda_{i}/n$ tendent tous en probabilité vers $(1-q)q^{i-1}$. L'expression ci-dessus est donc équivalente à 
$$\sum_{i=0}^{\infty}\esper[B_{i+1}](1-q)q^{i}\,\delta_{(1-q)\,q^{i}}(f)\,,$$
et remarquons par ailleurs que le membre de droite de l'équation est aussi une combinaison linéaire de Dirac en les $(1-q)q^{i}$. En choisissant des fonctions tests $f_{i}$ semblables à celles de la figure \ref{testfunction}, on conclut que :
\begin{align*}\esper[B_{2k+1}]&\simeq \frac{1}{2}\left(1+\frac{(q;q^{2})_{\infty}}{(q^{2};q^{2})_{\infty}}\,\frac{(q;q^{2})_{k}}{(q^{2};q^{2})_{k}}\,q^{k}\right)\,\,;\\
\esper[B_{2k+2}]&\simeq \frac{1}{2}\,.
\end{align*}
Ceci conclut la preuve, attendu que pour $n$ assez grand, les $k$ premières parts de la bipartition sont avec grande probabilité toutes distinctes, de sorte que $B_{k}$ paramètre sans ambiguïté l'événement <<~$\lambda_{k}$ tombe dans $\lambda^{(1)}$~>>.
\end{proof}\bigskip

Ainsi, si l'on est capable d'exprimer les $q$-caractères de type B en fonction des caractères de type B et de valider l'hypothèse faite sur les espérances des caractères centraux de cycles impairs, alors on peut décrire assez précisément la répartition asymptotique des parts de la bipartition. Le paragraphe suivant est consacré à une tentative d'attaque du problème des $q$-caractères de type B ; cette tentative est un échec, mais on décrira au passage les phénomènes de dualité de Schur-Weyl, ce qui motivera la discussion du chapitre \ref{schurweylmeasure}.
\bigskip

\section[Groupes quantiques, algèbres d'Ariki-Koike et dualité de Schur-Weyl]{Groupes quantiques, algèbres d'Ariki-Koike et\\ dualité de Schur-Weyl}\label{shoji}

Pour résoudre la conjecture énoncée dans le paragraphe précédent, il suffirait de disposer d'une $q$-formule de Frobenius pour les caractères de $\IH_{q}(\wsym_{n})$. Dans \cite{Sho00}, T. Shoji a démontré une telle formule, mais celle-ci ne peut être comprise que dans le cadre de la dualité de Schur-Weyl entre groupes quantiques et algèbres d'Hecke cyclotomiques. Ce paragraphe est consacré à l'exposé de ces résultats de dualité, ce qui préparera les raisonnements du chapitre suivant. Pour commencer, rappelons ce qu'est la dualité de Schur-Weyl classique. Si $m$ et $n$ sont deux entiers positifs, alors l'espace des tenseurs $V=(\C^{m})^{\otimes n}$ est un bimodule pour l'action du groupe linéaire complexe $\GL(m,\C)$ et du groupe symétrique $\sym_{n}$. Ainsi, le groupe linéaire agit sur les tenseurs simples de $V$ par
$$g \cdot (x_{1}\otimes x_{2} \otimes \cdots \otimes x_{n})=g(x_{1})\otimes g(x_{2}) \otimes \cdots \otimes g(x_{n})\,,$$
et le groupe symétrique agit (à droite) sur les tenseurs simples de $V$ par
$$ (x_{1}\otimes x_{2} \otimes \cdots \otimes x_{n})\cdot \sigma=x_{\sigma(1)}\otimes x_{\sigma(2)} \otimes \cdots \otimes x_{\sigma(n)}\,\,;$$
il est évident que ces actions commutent. On renvoie à \cite[\S3.9]{Weyl39} pour une preuve du résultat suivant :
\begin{theorem}[Dualité de Schur-Weyl, \cite{Weyl39}]\label{schurweylduality}
Soient $A$ et $B$ les algèbres engendrées dans $\hendo(V)$ par les actions respectives de $\GL(m,\C)$ et de $\sym_{n}$. Le commutant de $A$ est $B$, et le commutant de $B$ est $A$. Par suite, $V$ se décompose en somme directe de $(\GL(m,\C),\sym_{n})$-bimodules irréductibles :
$$_{\GL(m,\C) \curvearrowright } \left\{(\C^{m})^{\otimes n} \,\right\}_{\curvearrowleft \,\sym_{n}}=\bigoplus_{\lambda} \left({}_{\GL(m,\C) \curvearrowright }M^{\lambda}\right)\otimes_{\C}\left({V^{\lambda}}_{\curvearrowleft \,\sym_{n}}\right).$$
Les indices $\lambda$ parcourent l'ensemble des partitions de taille $n$ et de longueur inférieure à $m$ ; avec cette indexation, $V^{\lambda}$ est le module de Specht de type $\lambda$, et $M^{\lambda}$ est la représentation algébrique irréductible de $\GL(m,\C)$ de plus haut poids $\lambda$.
\end{theorem}\bigskip\bigskip

Ce résultat de dualité fournit une preuve de la formule de Frobenius du paragraphe \ref{frobenius}, à condition de connaître le caractère de l'action de $\GL(m,\C)$ sur le module $M^{\lambda}$. Ce dernier est donné par la \textbf{formule des caractères de Weyl}, voir \cite{Weyl39,Var89}. Ainsi, si $t=\mathrm{diag}(t_{1},\ldots,t_{m})$ est une matrice diagonale dans $\GL(m,\C)$, alors $$\varsigma_{M^{\lambda}}(t)=s_{\lambda}(t_{1},\ldots,t_{m})\,,$$
et ceci détermine entièrement $\varsigma_{M^{\lambda}}$, car les matrices diagonalisables sont denses dans $\GL(m,\C)$. Maintenant, si $(t,\sigma)$ désigne l'endomorphisme de $V=(\C^{m})^{\otimes n}$ produit de l'action de $t$ par l'action de $\sigma$, et si $\rho=t(\sigma)$, alors il est facile de voir que la trace de $(t,\sigma)$ sur $V$ est $p_{\rho}(t_{1},t_{2},\ldots,t_{m})$, puisqu'on peut écrire explicitement la matrice de représentation. Ainsi,
$$\forall m,\,\,\,p_{\rho}(t_{1},\ldots,t_{m})=\sum_{\substack{|\lambda|=n\\ \ell(\lambda)\leq m}} s_{\lambda}(t_{1},\ldots,t_{m})\,\,\varsigma^{\lambda}(\rho)\,,$$
d'où la formule de Frobenius en faisant tendre le nombre de variables $m$ vers l'infini. \bigskip
\bigskip

Ce raisonnement peut être généralisé dans deux directions : ainsi, on peut remplacer le groupe symétrique par un produit en couronne $(\Z/r\Z)\wr \sym_{n}$, et on peut également remplacer le groupe symétrique ou le produit en couronne par son algèbre d'Hecke ou d'Ariki-Koike, à condition de remplacer le groupe linéaire $\GL(m,\C)$ par un groupe quantique. Ceci conduit \emph{in fine} à une $q$-formule de Frobenius en type B, voir le théorème \ref{shojiformula} énoncé plus loin. Avant de décrire ces généralisations, il est utile de reformuler la proposition \ref{schurweylduality} en rempla\c cant le groupe $\GL(m,\C)$ par l'algèbre enveloppante universelle $U(\liegl(m),\C)$. Cette algèbre paramètre les opérateurs différentiels invariants par translation sur le groupe $\GL(m,\C)$, et en tant qu'algèbre à unité complexe, elle est engendrée par des générateurs $(e_{i})_{i \leq m-1}$, $(f_{i})_{i \leq m-1}$, $(h_{i})_{i \leq m-1}$, avec les relations :
\begin{align*}
\text{relations de Cartan :}&\quad[h_{i},e_{j}]=\alpha_{ij}\,e_{j}\quad;\quad[h_{i},f_{j}]=\alpha_{ij}\,f_{j}\quad;\quad[e_{i},f_{j}]=\delta_{ij}\,h_{i}\quad;\\
\text{relations de Serre :}&\quad[e_{i},[e_{i\pm1},e_{i}]]=0\quad;\quad[f_{i},[f_{i\pm1},f_{i}]]=0\quad;\\
\forall i,j,\,\,\,|i-j|\geq 2 &\quad\Rightarrow\quad [e_{i},e_{j}]=0\quad;\quad [f_{i},f_{j}]=0\,.
\end{align*}
où $\alpha_{ij}$ est le coefficient $(i,j)$ de la matrice de Cartan de type A$_{m}$, c'est-à-dire que $$\alpha_{ij}=2\delta_{i,j}-\delta_{i,i-1}-\delta_{i,i+1}\,.$$
L'algèbre $U(\liegl(m),\C)$ agit sur $\C^{m}$ par
$$e_{i} \mapsto v_{i}\otimes v^{i+1}\qquad;\qquad f_{i} \mapsto v_{i+1}\otimes v^{i}\qquad;\qquad h_{i} \mapsto v_{i}\otimes v^{i}-v_{i+1}\otimes v^{i+1}$$
où $(v_{1},\ldots,v_{m})$ est une base de $V=\C^{m}$, et où l'on utilise les notations d'Einstein pour les tenseurs de $\hendo(V)=V \otimes V^{*}$. D'autre part, elle admet une structure d'algèbre de Hopf pour le coproduit
$$\Delta(e_{i})=e_{i}\otimes 1 + 1\otimes e_{i}\quad;\quad \Delta(f_{i})=f_{i} \otimes 1 + 1 \otimes f_{i} \quad;\quad \Delta(h_{i})=h_{i} \otimes 1 +1 \otimes h_{i}\,.$$
Soit $\Delta^{(n)} : U(\liegl(m),\C) \to U(\liegl(m),\C)^{\otimes n}$ les morphismes d'algèbres définis récursivement par :
$$\Delta^{(2)}=\Delta\qquad;\qquad \Delta^{(n)}=(\Delta^{(n-1)}\otimes \id)\circ \Delta$$
Ces morphismes permettent de faire agir $U(\liegl(m),\C)$ sur un produit tensoriel : si $\rho$ désigne la représentation $U(\liegl(m),\C) \curvearrowright \C^{m}$ décrite ci-dessus, alors la représentation canonique $U(\liegl(m),\C) \curvearrowright (\C^{m})^{\otimes n}$ sera $\rho^{(n)} = \rho^{\otimes n}\circ \Delta^{(n)}$. Dans ce contexte, la dualité de Schur-Weyl s'écrit :
$$_{U(\liegl(m),\C) \curvearrowright } \left\{(\C^{m})^{\otimes n} \right\}_{\curvearrowleft \,\sym_{n}}=\bigoplus_{\substack{|\lambda|=n\\ \ell(\lambda)\leq m}} \left({}_{U(\liegl(m),\C) \curvearrowright \,}M^{\lambda}\right)\otimes_{\C}\left({V^{\lambda}}_{\curvearrowleft \,\sym_{n}}\right).$$
et ceci est réellement une reformulation du théorème \ref{schurweylduality}, car les actions de $\GL(m,\C)$ et de $U(\liegl(m),\C)$ engendrent la même sous-algèbre fermée de l'algèbre des d'endomorphismes $\hendo_{\C}(V)$.
\bigskip
\bigskip

Ceci étant, dans la dualité de Schur-Weyl, on peut remplacer le groupe symétrique par un produit en couronne, et dans ce cas, le commutant de l'action sur un produit tensoriel est une sous-algèbre de Lévi de $U(\liegl(m),\C)$. Fixons un entier $r$ plus grand que $2$, une racine primitive $r$-ième de l'unité $\zeta$, une composition $m=m_{1}+\cdots+m_{r}$, et une base $(v_{i,j})_{i \leq r,j \leq m_{i}}$ de $\C^{m}$ adaptée à la décomposition en somme directe $\C^{m}=\C^{m_{1}}\oplus \cdots \oplus \C^{m_{r}}$, c'est-à-dire que le sous-espace $\C^{m_{i}}$ est engendré par les vecteurs $v_{i,1},\ldots,v_{i,m_{i}}$. On rappelle que le produit en couronne $(\Z/r\Z) \wr \sym_{n}$ est le produit semi-direct $(\Z/r\Z)^{n} \rtimes \sym_{n}$, où $\sym_{n}$ agit sur $(\Z/r\Z)^{n}$ par permutation des lettres. Ce groupe est un groupe de réflexions complexes engendré par des éléments $s_{0},s_{1},\ldots,s_{n-1}$ vérifiant les relations :
\begin{align*} &(s_{0})^{r}=1\quad;\quad (s_{i\geq 1})^{2}=1\quad;\\
&s_{0}s_{1}s_{0}s_{1}=s_{1}s_{0}s_{1}s_{0}\quad;\quad\forall i\geq1,\,\,\,s_{i}s_{i+1}s_{i}=s_{i+1}s_{i}s_{i+1}\quad;\\
&\forall i,j,\,\,\,|i-j|\geq 2 \Rightarrow s_{i}s_{j}=s_{j}s_{i}\,.
\end{align*}
On fait agir $\wsym_{n,r}$ sur les tenseurs simples de $(\C^{m})^{\otimes n}$ en suivant les règles :
\begin{align*}
s_{0}(v_{i_{1},j_{1}}\otimes v_{i_{2},j_{2}}\otimes \cdots \otimes v_{i_{n},j_{n}})&= \zeta^{i_{1}}\,(v_{i_{1},j_{1}}\otimes v_{i_{2},j_{2}}\otimes \cdots \otimes v_{i_{n},j_{n}})\,\,; \\
s_{k \geq 1}(v_{i_{1},j_{1}}\otimes v_{i_{2},j_{2}}\otimes \cdots \otimes v_{i_{n},j_{n}})&=v_{i_{1},j_{1}}\otimes \cdots \otimes v_{i_{k+1},j_{k+1}} \otimes v_{i_{k},j_{k}} \otimes \cdots\otimes v_{i_{n},j_{n}}\,.
\end{align*}
Ces règles sont clairement compatibles avec les relations définissant $\wsym_{n,r}$, d'où une représentation $\wsym_{n,r} \curvearrowright (\C^{m})^{\otimes n}$.  De plus, comme $\sym_{n}=\langle s_{1},\ldots,s_{n-1}\rangle$ est inclus dans $\wsym_{n,r}$, le commutant de cette action est inclus dans l'image de $U(\liegl(m),\C)$ sur l'algèbre $\hendo_{\C}(V)$, et on peut montrer que c'est l'image de l'algèbre enveloppante de la \textbf{sous-algèbre de Lévi}
$$\mathfrak{g}=\liegl(m_{1},\C)\oplus\liegl(m_{2},\C)\oplus\cdots \oplus\liegl(m_{r},\C)\,.$$
L'algèbre enveloppante $U(\mathfrak{g},\C)$ est la sous-algèbre de $U(\liegl(m),\C)$ obtenue en retirant les générateurs $e_{i}$ et $f_{i}$ avec $i$ indice dans
$$\{m_{1},m_{1}+m_{2},m_{1}+m_{2}+m_{3},\ldots,m_{1}+\cdots+m_{r-1}\}\,,$$
c'est-à-dire dans l'ensemble des descentes de la composition $m=m_{1}+m_{2}+\cdots+m_{r}$.
Ainsi, à partir de la dualité de Schur-Weyl classique, il n'est pas difficile de montrer que :
\begin{proposition}[Dualité de Schur-Weyl cyclotomique]
Soient $A$ et $B$ les algèbre engendrées dans $\hendo_{\C}(V)$ par les actions respectives de $U(\mathfrak{g},\C)$ et de $\C\wsym_{n,r}$. Le commutant de $A$ est $B$, et le commutant de $B$ est $A$. Par suite, $V$ se décompose en somme directe de $(U(\mathfrak{g},\C),\C\wsym_{n,r})$-bimodules irréductibles :
$$_{U(\mathfrak{g},\C) \curvearrowright } \left\{(\C^{m})^{\otimes n} \right\}_{\curvearrowleft \,\C\wsym_{n,r}}=\bigoplus_{\substack{|\lambda^{(1)}|+\cdots+|\lambda^{(r)}|=n\\ \forall i,\,\,\ell(\lambda^{(i)})\leq m_{i}}} \left({}_{U(\mathfrak{g},\C) \curvearrowright \,}M^{\Lambda}\right)\otimes_{\C}\left({V^{\Lambda}}_{\curvearrowleft \,\C\wsym_{n,r}}\right)\,.$$
\end{proposition}\medskip

\noindent Comme les caractères irréductibles de $\GL(m_{1},\C)\times \GL(m_{2},\C)\times \cdots \times \GL(m_{r},\C)$ sont donnés par des produits de fonctions de Schur, on peut comme dans le cas du groupe symétrique en déduire la formule de Frobenius pour les caractères irréductibles de $\wsym_{n,r}$ (voir \cite[appendice B]{Mac95}). 
\bigskip
\bigskip

Dans une toute autre direction, on doit à M. Jimbo une quantification de la dualité de Schur-Weyl $U(\liegl(m),\C) - \C\sym_{n}$ faisant apparaître à droite l'algèbre d'Iwahori-Hecke (générique) $\IH(\sym_{n})$, et à gauche un \textbf{groupe quantique} (voir \cite{Jim86}). Pour commencer, décrivons une représentation de l'algèbre d'Hecke sur le produit tensoriel $V=(\C(q^{1/2})^{m})^{\otimes n}$. Si $(v_{1},\ldots,v_{m})$ est une base de $\C(q^{1/2})^{m}$, les éléments de base du produit tensoriel $w=v_{i_{1}}\otimes \cdots \otimes v_{i_{n}}$ peuvent être vus comme des mots de longueur $n$ sur l'alphabet $\lle 1,m\rre$. On fait agir un générateur $T_{j}$ de l'algèbre $\IH(\sym_{n},\C(q^{1/2}))=\C(q^{1/2})\otimes_{\C(q)}\IH(\sym_{n})$ par :
$$w\cdot T_{j}=\begin{cases} q\,w&\text{si }i_{j}=i_{j+1},\\
(q-1)\,w+q^{1/2}\,\widetilde{w}^{j}&\text{si }i_{j}<i_{j+1},\\
q^{1/2}\,\widetilde{w}^{j}&\text{si }i_{j}>i_{j+1},
\end{cases}$$
où $\widetilde{w}^{j}$ désigne le mot obtenu à partir de $w$ en échangeant les lettres d'indices $j$ et $j+1$. On montre sans difficulté que ces formules sont bien compatibles avec les relations définissant $\IH(\sym_{n},\C(q^{1/2}))$. On a donc bien une représentation d'algèbres, et le calcul du commutant de cette représentation met en jeu une quantification $U(\liegl(m),\C(q^{1/2}))$ de l'algèbre enveloppante universelle $U(\liegl(m),\C)$. Ainsi, considérons la $\C(q^{1/2})$-algè\-bre à unité $U(\liegl(m),\C(q^{1/2}))$ engendrée par des générateurs $(e_{i})_{i \leq m-1}$, $(f_{i})_{i \leq m-1}$, $(q^{\pm \eps_{i}})_{i \leq m}$, avec les relations :
\begin{align*}
q\text{-relations de Cartan :}&\quad q^{\pm \eps_{i}}\,q^{\pm' \eps_{j}}=q^{\pm' \eps_{j}}\,q^{\pm \eps_{i}}\quad;\quad q^{\eps_{i}}\,q^{-\eps_{i}}=q^{-\eps_{i}}\,q^{\eps_{i}}=1\,\,;\\
&\quad q^{\eps_{i}}\,e_{j}\, q^{-\eps_{i}}=\begin{cases} q^{-1/2}\,e_{j}&\text{si }j=i-1,\\
q^{1/2}\,e_{j}&\text{si }j=i,\\
e_{j} &\text{sinon } ;
\end{cases}\\
&\quad q^{\eps_{i}}\,f_{j}\, q^{-\eps_{i}}=\begin{cases} q^{1/2}\,f_{j}&\text{si }j=i-1,\\
q^{-1/2}\,f_{j}&\text{si }j=i,\\
f_{j} &\text{sinon } ;
\end{cases}\\
&\quad [e_{i},f_{j}]=\delta_{ij}\,\,\frac{q^{\eps_{i}-\eps_{i+1}} -q^{\eps_{i+1}-\eps_{i}}}{q^{1/2}-q^{-1/2}}\,\,;\\
q\text{-relations de Serre :}&\quad[e_{i},[e_{i\pm1},e_{i}]_{q}]_{q}=0\quad;\quad[f_{i},[f_{i\pm1},f_{i}]_{q}]_{q}=0\,\,;\\
\forall i,j,\,\,\,|i-j|&\geq 2 \quad\Rightarrow\quad [e_{i},e_{j}]=0\quad;\quad [f_{i},f_{j}]=0\,.
\end{align*}
On vérifie sans mal que les relations du groupe quantique $U(\liegl(m),\C(q^{1/2}))$ sont compatibles avec les associations 
$$e_{i} \mapsto v_{i}\otimes v^{i+1}\qquad;\qquad f_{i} \mapsto v_{i+1}\otimes v^{i}\qquad;\qquad q^{\pm\eps_{i}} \mapsto q^{\pm1/2}\,v_{i}\otimes v^{i}+\sum_{j \neq i}v_{j} \otimes v^{j}$$
d'où une représentation canonique de $U(\liegl(m),\C(q^{1/2}))$ sur $\C(q^{1/2})^{m}$. Pour la relever en une représentation sur un produit tensoriel $(\C(q^{1/2})^{m})^{\otimes n}$, on utilise la $q$-structure d'algèbre de Hopf\footnote{Ainsi, par groupe quantique, V. G. Drinfeld et M. Jimbo entendent une $q$-déformation d'une algèbre de Hopf obtenue en considérant l'algèbre enveloppante universelle de l'algèbre de Lie d'un groupe (ou plus généralement une algèbre de Kac-Moody), voir \cite{Dri86,Jim86}. Nos notations de groupes quantiques --- $U(\liegl(m),\C(q^{1/2}))$ ---  différent des notations classiques --- $U_{q}(\liegl(m))$ --- mais nous souhaitions rappeler tout au long de l'exposé l'anneau des coefficients sur lequel on raisonne --- ici, $R=\C(q^{1/2})$, et plus loin $R=\C(q^{1/2},u)$.} :
$$\Delta(e_{i})=q^{\eps_{i}-\eps_{i+1}}\otimes e_{i} + e_{i} \otimes 1\quad;\quad \Delta(f_{i})=1\otimes f_{i} + f_{i } \otimes q^{\eps_{i+1}-\eps_{i}} \quad;\quad \Delta(q^{\pm\eps_{i}})=q^{\pm \eps_{i}} \otimes q^{\pm \eps_{i}}\,.$$
Par itération du coproduit, on obtient des morphismes d'algèbres $$\Delta^{(n)} : U(\liegl(m),\C(q^{1/2}))\to U(\liegl(m),\C(q^{1/2}))^{\otimes n}\,,$$ et ainsi des représentations canoniques de $U(\liegl(m),\C(q^{1/2}))$ sur $V=(\C(q^{1/2})^{m})^{\otimes n}$. Dans ce contexte, l'analogue de la proposition \ref{schurweylduality} est :
\begin{proposition}[Jimbo, \cite{Jim86}]
Soient $A$ et $B$ les algèbres engendrées dans $\hendo_{\C(q^{1/2})}(V)$ par les actions respectives du groupe quantique $U(\liegl(m),\C(q^{1/2}))$ et de l'algèbre $\IH(\sym_{n},\C(q^{1/2}))$. Le commutant de $A$ est $B$, et le commutant de $B$ est $A$. Par suite, $V$ se décompose en somme directe de $(U(\liegl(m),\C(q^{1/2})),\IH(\sym_{n},\C(q^{1/2})))$-bimodules irréductibles :
\begin{align*}&_{U(\liegl(m),\C(q^{1/2})) \curvearrowright } \left\{(\C(q^{1/2})^{m})^{\otimes n} \right\}_{\curvearrowleft \,\IH(\sym_{n},\C(q^{1/2}))}\\
&=\bigoplus_{\substack{|\lambda|=n\\ \ell(\lambda)\leq m}} \left({}_{U(\liegl(m),\C(q^{1/2})) \curvearrowright \,}M^{\lambda}(q)\right)\otimes_{\C(q^{1/2})}\left({V^{\lambda}(q)}_{\curvearrowleft \,\IH(\sym_{n},\C(q^{1/2}))}\right).\end{align*}
\end{proposition}
\noindent Pour $q$ générique, le résultat subsiste en spécialisant la valeur de $q$ et en considérant tous les objets comme $\C$-espaces vectoriels ou $\C$-algèbres. La formule de Ram \ref{ramformula} s'en déduit en regardant la trace de l'action du produit d'un élément <<~diagonal~>> du groupe quantique $U(\liegl(m),\C(q^{1/2}))$ et d'un élément de l'algèbre d'Hecke\footnote{Compte tenu des résultats évoqués dans la section \ref{hecke}, et au vu des raisonnements du présent paragraphe, il est tentant de penser que les caractères des algèbres d'Hecke peuvent aussi être calculés en regardant la bitrace d'un élément de $\IH_{q}(\sym_{n})$ et d'un élément de $\GL(n,\For_{q})$ sur le module $\C[G/B]$. C'est effectivement le cas, et ce calcul a été réalisé par T. Halverson et A. Ram dans \cite{HR99}.}, voir \cite[paragraphes 3 et 4]{Ram91}.\bigskip
\bigskip

Finalement, on peut réunir les deux généralisations de la proposition \ref{schurweylduality} et décrire une dualité de Schur-Weyl entre le groupe quantique $U(\mathfrak{g},\C(q^{1/2}))$ et une algèbre d'Hecke cyclotomique obtenue par déformation de l'algèbre de groupe $\C\wsym_{n,r}$. Ces \textbf{algèbres d'Hecke cyclotomiques} ont été introduites par S. Ariki et K. Koike dans \cite{AK94}, et elles dépendent de paramètres $q,u_{1},\ldots,u_{r}$ :
\begin{align*}\IH(\wsym_{n,r})=\langle T_{0},T_{1},\ldots,T_{n-1} \rangle_{\C(q,u_{1},\ldots,u_{r})},\,\,\,\text{avec  :}\quad&(T_{0}-u_{1})(T_{0}-u_{2})\cdots(T_{0}-u_{r})=0\\
&\forall i \geq 1,\,\,(T_{i}-q)(T_{i}+1)=0\\
&T_{0}T_{1}T_{0}T_{1}=T_{1}T_{0}T_{1}T_{0}\\
&\forall i \geq 1,\,\,T_{i}T_{i+1}T_{i}=T_{i+1}T_{i}T_{i+1}\\
&\forall i,j,\,\,|i-j|\geq 2 \Rightarrow T_{i}T_{j}=T_{j}T_{i}\,.
\end{align*}
\begin{examples} Lorsqu'on spécialise $q$ en $1$ et $u_{i}$ en $\zeta^{-i}$, on retrouve l'algèbre de groupe $\C\wsym_{n,r}$. D'autre part, si $r=1$, alors pour toute valeur du paramètre $u_{1}$, l'algèbre obtenue est simplement l'algèbre d'Iwahori-Hecke du groupe symétrique. Enfin, pour $r=2$, si l'on spécialise $u_{1}$ en $q$ et $u_{2}$ en $-1$, on obtient l'algèbre d'Iwahori-Hecke du groupe hyperoctahédral. De fa\c con générique, l'algèbre d'Ariki-Koike $\IH(\wsym_{n,r})$ est semi-simple, et a la même théorie des représentations que l'algèbre de groupe $\C\wsym_{n,r}$. Ainsi, les représentations irréductibles de $\IH(\wsym_{n,r})$ sont indexées par les $r$-uples de partitions dans $\ym_{n,r}$.
\end{examples}\bigskip\bigskip

L'action de l'algèbre d'Ariki-Koike sur le produit tensoriel $V=(\C(q^{1/2},u)^{m})^{\otimes n}$ étend celle de $\IH(\sym_{n},\C(q^{1/2})) \subset \IH(\wsym_{n,r},\C(q^{1/2},u))$, et il suffit donc de décrire l'action de $T_{0}$. Si $k \in \lle 1,n-1\rre$, on définit un endomorphisme $S_{k} \in \hendo_{\C(q^{1/2},u)}(V)$ par son action sur les tenseurs simples :
$$S_{k}(w=v_{i_{1},j_{1}}\otimes v_{i_{2},j_{2}}\otimes \cdots \otimes v_{i_{n},j_{n}})= \begin{cases} T_{k}(w) & \text{si }i_{k}=i_{k+1},\\
q^{1/2}\,\widetilde{w}^{k} &\text{si }i_{k}\neq i_{k+1},  \end{cases}$$
étant entendu que les indices $(i_{k'},j_{k'})$ sont ordonnés par l'ordre lexicographique. D'autre part, si $w=v_{i_{1},j_{1}}\otimes v_{i_{2},j_{2}}\otimes \cdots \otimes v_{i_{n},j_{n}}$ on note $\varpi_{k}(w)=u_{i_{k}}\,w$ pour $k$ entre $1$ et $n$. Alors, l'action de $T_{0}$ est donnée par :
$$T_{0}(w)=T_{1}^{-1}T_{2}^{-1}\cdots T_{n-1}^{-1}S_{n-1}S_{n-2}\cdots S_{1}\varpi_{1}(w)\,.$$
Il est vrai, mais tout à fait non trivial, que cette règle est compatible avec les relations de l'algèbre d'Ariki-Koike, voir \cite[théorème 3.2]{SS99}. Alors, la dualité de Schur-Weyl généralisée s'écrit :
\begin{proposition}[Sakamoto-Shoji, \cite{SS99}]
Soient $A$ et $B$ les algèbres engendrées par les actions respectives de $U(\mathfrak{g},\C(q^{1/2},u))$ et de $\IH(\wsym_{n,r},\C(q^{1/2},u))$. Le commutant de $A$ est $B$, et le commutant de $B$ est $A$. Par suite, $V$ se décompose en somme directe de bimodules irréductibles :
\begin{align*}&_{U(\mathfrak{g},\C(q^{1/2},u)) \curvearrowright } \left\{(\C(q^{1/2},u))^{m})^{\otimes n} \right\}_{\curvearrowleft \,\IH(\wsym_{n,r},\C(q^{1/2},u))}\\
&=\bigoplus_{\substack{|\lambda^{1}|+\cdots+|\lambda^{r}|=n\\ \forall i,\,\,\ell(\lambda^{i})\leq m_{i}}} \left({}_{U(\mathfrak{g},\C(q^{1/2},u)) \curvearrowright \,}M^{\Lambda}(q)\right)\otimes_{\C}\left({V^{\Lambda}(q,u)}_{\curvearrowleft \,\IH(\wsym_{n,r},\C(q^{1/2},u))}\right).\end{align*}
\end{proposition}\bigskip\bigskip

Le résultat de dualité de Sakamoto-Shoji mérite la précision suivante. Si les entiers $m_{i}$ sont tous plus grands que $n$, alors tous les modules irréductibles de $\IH(\wsym_{n,r})$ interviennent dans la décomposition de $(\C(q^{1/2},u)^{m})^{\otimes n}$, donc la représentation de l'algèbre d'Ariki-Koike est fidèle. Ceci permet de définir un élément de $\IH(\wsym_{n,r})$ directement par son action\footnote{En supposant bien sûr que la dite action est dans le commutant de $U(\mathfrak{g},\C(q^{1/2},u))$.} sur un produit tensoriel $n$-ième, en supposant que l'espace de base est de dimension assez grande. En particulier, l'analogue quantique et cyclotomique de la formule de Frobenius-Schur pour les caractères irréductibles des algèbres d'Ariki-Koike s'énoncera dans ce cadre, voir le théorème \ref{shojiformula}. Jusqu'à la fin de ce paragraphe, on considérera donc $\IH(\wsym_{n,r},R=\C(q^{1/2},u))$ comme une sous-algèbre de $\hendo_{R}(V)$, les paramètres $m_{i}$ étant tous supposés suffisamment grands. 
\bigskip
\bigskip

La remarque du paragraphe précédent s'applique en particulier aux opérateurs $S_{k}$ et $\varpi_{k}$ : comme ils sont dans le commutant du groupe quantique $U(\mathfrak{g},R)$, ils appartient à l'algèbre d'Ariki-Koike $\IH(\wsym_{n,r},R)$. Décrivons les relations qu'ils entretiennent avec les générateurs usuels $T_{k}$. On introduit le déterminant de Vandermonde
$$\Delta(u)=\det((u_{i})^{j-1})_{i,j \in \lle 1,r\rre}=\prod_{i<j}\, (u_{j}-u_{i})\,,$$
et les polynômes $F_{j}(X)=\Delta(u)\,\prod_{k \neq j}(X-u_{k})/(u_{j}-u_{k})$. Alors, on peut montrer que :
$$T_{j}\,\varpi_{j+1}-\varpi_{j}\,T_{j}= \varpi_{j+1}\,T_{j}-T_{j}\,\varpi_{j}=\frac{q-1}{\Delta^{2}}\,\sum_{1\leq a<b\leq r}(u_{b}-u_{a})\,F_{a}(\varpi_{j})\,F_{b}(\varpi_{j+1})\,.$$
D'autre part, si $|i-j|\geq 2$, alors $T_{i}\,\varpi_{j}-\varpi_{j}\,T_{i}=0$. Il s'en déduit que les $S_{j}$ s'expriment en fonction des $T_{j}$ et des $\varpi_{j}$ :
$$\forall j,\,\,\,S_{j}=T_{j}-\frac{q-1}{\Delta^{2}}\sum_{1 \leq a <b \leq r}F_{a}(\varpi_{j})\,F_{b}(\varpi_{j+1})\,.$$
Comme $T_{0}$ s'exprime en fonction des $S_{j}$ et de $\varpi_{1}$, ceci fournit une nouvelle présentation de l'algèbre d'Ariki-Koike. Ainsi, $\IH(\wsym_{n,r},R)=\langle T_{1},\ldots,T_{n-1},\varpi_{1},\ldots,\varpi_{n}\rangle_{R}$, avec les relations :
\begin{align*}
&\forall i\geq 1,\,\,\,(\varpi_{i}-u_{1})(\varpi_{i}-u_{2})\cdots(\varpi_{i}-u_{r})=0\\
&\forall i\geq 1,\,\,\,(T_{i}-q)(T_{i}+1)=0\\
&\forall i\geq1,\,\,\,T_{i}T_{i+1}T_{i}=T_{i+1}T_{i}T_{i+1}\\
&\forall i,j,\,\,\,|i-j| \geq 2 \Rightarrow T_{i}T_{j}=T_{j}T_{i}\\
&\forall i,j,\,\,\,|i-j| \geq 2 \Rightarrow T_{i}\,\varpi_{j}=\varpi_{j}\,T_{i}\\
&\forall i,j,\,\,\,\varpi_{i}\,\varpi_{j}=\varpi_{j}\,\varpi_{i}\\
&\forall i \geq 1,\,\,\,T_{i}\,\varpi_{i+1}-\varpi_{i}\,T_{i}= \varpi_{i+1}\,T_{i}-T_{i}\,\varpi_{i}=(q-1)\,\Delta^{-2}\!\!\!\!\sum_{1\leq a<b\leq r}\!\!(u_{b}-u_{a})\,F_{a}(\varpi_{i})\,F_{b}(\varpi_{i+1})\,.
\end{align*}
Par rapport à la présentation usuelle, on a ôté l'élément particulier $T_{0}$, et on l'a remplacé par $n$ éléments <<~symétriques~>> $\varpi_{1},\ldots,\varpi_{n}$. Lorsqu'on spécialise les valeurs de $q$ et de $u$ pour retrouver l'algèbre de groupe $\C\wsym_{n,r}$, les éléments $\varpi_{i}$ spécialisent en les éléments
$$\pi_{1}=s_{0}\quad;\quad\pi_{2}=s_{1}s_{0}s_{1}\quad;\quad\pi_{i}=s_{i-1}\cdots s_{1}s_{0}s_{1}\cdots s_{i-1}\,,$$
et vu comme permutation en couronne, chaque élément $\pi_{i}$ correspond à l'élément $$((1,\ldots,1,\zeta_{(i)},1\ldots,1),\id_{\lle 1,n\rre})\,$$ où $\zeta$ est une racine $r$-ième primitive de l'unité. Par suite, tout élément de $\IH(\wsym_{n,r},R)$ s'écrit de manière unique comme combinaison linéaire de produits 
$$\varpi_{1}^{c_{1}}\,\varpi_{2}^{c_{2}}\,\cdots\,\varpi_{n}^{c_{n}}\,T_{\sigma}\,,$$
où chaque $c_{i}$ est un entier compris entre $1$ et $r$, et où $T_{\sigma}=T_{i_{1}}\cdots T_{i_{l}}$ si $\sigma=s_{i_{1}}\cdots s_{i_{l}}$ est une expression réduite de la permutation $\sigma$.\bigskip
\bigskip

Nous sommes finalement prêts à décrire l'analogue de la formule de Frobenius-Schur dans le contexte des algèbres d'Ariki-Koike. Pour tout entier $i \in \lle 1,r \rre$, on pose :
$$a(n,i) = \varpi_{n}^{i}\,T_{n-1}T_{n-2}\cdots T_{1} \in \IH(\wsym_{n,r},R)\,.$$
Si $\mu$ est une partition d'un entier $n$, on définit de même
$$a(\mu,i)=a(\mu_{1},i)\otimes a(\mu_{2},i) \otimes \cdots \otimes a(\mu_{t},i)\,,$$
étant entendu que l'algèbre parabolique $\IH(\wsym_{\mu_{1},r},R)\otimes \cdots \otimes \IH(\wsym_{\mu_{t},r},R)$ se plonge canoniquement\footnote{C'est l'intérêt principal de la seconde présentation de l'algèbre d'Ariki-Koike ; comme les éléments $\varpi_{i}$ ont tous un rôle symétrique, on peut plonger <<~diagonalement~>> un produit d'algèbres d'Ariki-Koike de tailles respectives $n_{1},\ldots,n_{t}$ dans l'algèbre de taille $n_{1}+\cdots+n_{t}$, simplement en envoyant générateurs sur générateurs. Ce n'est pas possible avec la présentation originale, car $T_{0}$ joue un rôle particulier dans chaque algèbre $\IH(\wsym_{n_{i},r})$.} dans $\IH(\wsym_{n,r},R)$, voir \cite[\S4]{Sho00}. Enfin, si $M=(\mu^{(1)},\ldots,\mu^{(r)})$ est un $r$-uplet de partitions dans $\ym_{n,r}$, on note
$$A_{M}=a(\mu^{(1)},1)\otimes a(\mu^{(2)},2)\otimes \cdots \otimes a(\mu^{(r)},r)\,,$$
où de nouveau l'on identifie $\IH(\wsym_{|\mu^{(1)}|,r},R)\otimes \cdots \otimes \IH(\wsym_{|\mu^{(r)}|,r},R)$ à une sous-algèbre parabolique de $\IH(\wsym_{n,r},R)$. Par construction, sous la spécialisation $\IH(\wsym_{n,r},R)\rightarrow \C\wsym_{n,r}$, les éléments $A_{M}$ sont envoyés sur des permutations en couronne de type $M$. \bigskip\bigskip

Notons $\mathscr{S}_{n,r}$ l'ensemble des suites d'entiers positifs ou nuls $s=(n_{1},\ldots,n_{r})$ telles que $n_{1}+n_{2}+\cdots+n_{r}=n$. \'Etant donnée une telle suite, on note $\ell(s)$ le nombre d'entiers $k$ tels que $n_{k}\neq 0$, et $u_{s}=u_{k}$, où $k$  est le plus grand entier tel que $n_{k} \neq 0$. On introduit alors des déformations des polynômes multi-symétriques $P^{(i)}_{n}(X_{1},\ldots,X_{r})$ :
$$Q^{(i)}_{n}(X_{1},\ldots,X_{r})=\sum_{s \in \mathscr{S}_{n,r}} (u_{s})^{i}\,(q-1)^{\ell(s)-1}\,\prod_{k=1}^{r}q_{n_{k}}(X_{k},q)\,.$$
Lorsque $u_{k}=\zeta^{-k}$ et $q=1$, les polynômes $q_{n_{k}}(X_{k},q)$ spécialisent en $p_{n_{k}}(X_{k})$, et on retrouve bien les fonctions $P^{(i)}_{n}(X_{1},\ldots,X_{r})$. Finalement, si $M=(\mu^{(1)},\ldots,\mu^{(r)})$ est un $r$-uplet de partitions, on note :
$$Q_{M}(X_{1},\ldots,X_{r})=\prod_{i=1}^{r}Q_{\mu^{(i)}}^{(i)}(X_{1},\ldots,X_{r})=\prod_{i=1}^{r}\prod_{j=1}^{\ell(\mu^{(i)})} Q_{\mu^{(i)}_{j}}^{(i)}(X_{1},\ldots,X_{r})\,.$$

\begin{theorem}[Shoji, \cite{Sho00}]\label{shojiformula}
Les caractères irréductibles de l'algèbre d'Ariki-Koike $\IH(\wsym_{n,r},R)$ sont donnés par la formule :
$$\forall M \in \ym_{n,r},\,\,\,Q_{M}(X_{1},\ldots,X_{r})=\sum_{\Lambda \in \ym_{n,r}} \varsigma^{\Lambda}(q,u,A_{M})\,S_{\Lambda}(X_{1},\ldots,X_{r})\,.$$
\end{theorem}\bigskip
\bigskip

Le théorème \ref{shojiformula} peut être complété par les deux remarques suivantes. D'une part, il existe un algorithme simple qui transforme tout élément $\varpi_{1}^{c_{1}}\,\varpi_{2}^{c_{2}}\,\cdots\,\varpi_{n}^{c_{n}}\,T_{\sigma}$ en une $R$-combinaison linéaire de $A_{M}$ qui prend les mêmes valeurs contre tout caractère irréductible (voir \cite[proposition 7.5]{Sho00}) --- c'est l'analogue pour les algèbres d'Ariki-Koike du résultat présenté page \pageref{reductioncharacter} pour les algèbres d'Hecke des groupes de Coxeter. Ainsi, la formule \ref{shojiformula} donne bien toutes les valeurs des caractères des algèbres $\IH(\wsym_{n,r},R)$. D'autre part, on peut comme dans le paragraphe \ref{ram} décomposer les polynômes $Q_{M}(X_{1},\ldots,X_{r})$ dans la base de $\bigotimes_{i=1}^{r}\Lambda[X_{i}]$ constituée des $P_{M}(X_{1},\ldots,X_{r})$, et ainsi obtenir l'analogue cyclotomique de la formule de Ram. Si $\Lambda \in \ym_{n,r}$, alors on note
$$z_{\Lambda}=r^{\ell(\Lambda)}\,\prod_{i=1}^{r} z_{\lambda^{(i)}}\,,$$
et $q^{\Lambda}-1$ est bien sûr $\prod_{i=1}^{r}q^{\lambda^{(i)}}-1$. D'autre part, si $\lambda$ est une partition, on note $\mathscr{M}_{\lambda,r}=\prod_{i=1}^{\ell(\lambda)}\mathscr{S}_{\lambda_{i},r}$ l'ensemble des matrices de taille $\ell(\lambda) \times r$ et à coefficients $c_{ij}\geq 0$ tels que $\sum_{j=1}^{r} c_{ij}=\lambda_{i}$. Pour $A \in \mathscr{M}_{\lambda,r}$, on note $$u_{A}=\prod_{i=1}^{\ell(\lambda)} u_{(c_{ij})_{1 \leq j \leq r}}\,,$$ 
et $\tau(A)=(\tau^{(1)},\ldots,\tau^{(r)})$ est l'élément de $\ym_{n,r}$ dont la $j$-ième partition $\tau^{(j)}$ a pour parts les $c_{ij}$. Si $\Lambda$ est une multipartition, alors $\mathscr{M}_{\Lambda,r}=\prod_{i=1}^{r} \mathscr{M}_{\lambda^{(i)},r}$, et pour $B = (B_{i})\in \mathscr{M}_{\Lambda,r}$, on note
$$u_{B}=\prod_{i=1}^{r} (u_{B_{i}})^{i} \,,$$
et $\tau(B)=\sqcup_{i=1}^{r}\tau(B_{i})$. Finalement, pour toute multipartition $\Lambda \in \ym_{n,r}$, on introduit le groupe 
$$\wsym_{\Lambda,r}=\left((\sym_{\lambda^{(1)}}) \rtimes (\Z/r\Z)^{|\lambda^{(1)}|}\right)\times\, \cdots\, \times \left((\sym_{\lambda^{(r)}}) \rtimes (\Z/r\Z)^{|\lambda^{(r)}|}\right)\,,$$ et on note de même $\IH(\wsym_{\Lambda,r})$ la sous-algèbre <<~parabolique~>> de l'algèbre d'Hecke cyclotomique $\IH(\wsym_{n,r})$. Alors :
\begin{proposition}[Formule de Shoji, \cite{Sho00}] Pour tout $\Lambda \in \ym_{n,r}$, $\theta_{\Lambda}$ désigne le caractère complexe de dimension $1$ de $\IH(\wsym_{\Lambda,r})$ défini par $\theta_{\Lambda}(T_{i\geq 1})=1$ et $\theta_{\Lambda}(\varpi_{j})=\zeta^{i}$, où $i$ est choisi de telle sorte que $|\lambda^{(1)}|+\cdots+|\lambda^{i-1}|<j \leq|\lambda^{(1)}|+\cdots+|\lambda^{(j)}|$. On note $\Theta_{\Lambda}$ le caractère induit $\mathrm{Ind}_{\IH(\wsym_{\Lambda,r})}^{\IH(\wsym_{n,r})}(\theta_{\lambda})$. Alors :
$$\forall M,\,\,\,(q-1)^{\ell(M)}\,\varsigma^{\Lambda}(q,u,A_{M})= \sum_{N \in \ym_{n,r}} \sum_{B \in \mathscr{M}_{M,r}}\,u_{B}\,\frac{\Theta_{\tau(B)}(A_{N})}{z_{N}} \,(q^{N}-1) \,\varsigma^{\Lambda}(N) \,.$$
\end{proposition}

\begin{example}
Si $r=1$, alors on retrouve exactement la formule de Ram \ref{ramformula}. En effet, $u_{1}=1$, et d'autre part, les caractères $\Theta_{\Lambda}(M)$ interviennent dans les formules de passage entre les $P_{M}$ et les $m_{\Lambda}(X_{1},\ldots,X_{r})=\prod_{i=1}^{r}m_{\lambda^{(i)}}(X_{i})$. En effet :
$$\overline{P}_{M}(X_{1},\ldots,X_{r})=\sum_{\Lambda \in \ym_{n,r}} \Theta_{\Lambda}(M)\,m_{\Lambda}(X_{1},\ldots,X_{r})\,.$$
Ainsi, en taille $r=1$, il y a exactement un élément $B \in \mathscr{M}_{M,1}$, et le coefficient $u_{B} \,\Theta_{\tau(B)}(\Lambda)$ n'est nul autre que $\scal{h_{M}}{p_{N}}$, ce qui donne la formule \ref{ramformula}.
\end{example}\bigskip
\bigskip

Au premier abord, la formule \ref{shojiformula} semble être exactement ce dont on avait besoin, et pour $r=2$, en spécialisant $u_{1}$ en $q$ et $u_{2}$ en $-1$, on obtient effectivement une formule de changement de base entre les $\varsigma^{\Lambda}(q,A_{\mu^{(1)},\mu^{(2)}})$ et les $\varsigma^{\Lambda}(\nu^{(1)},\nu^{(2)})$. Mais malheureusement, les $q$-caractères sont maintenant exprimés en les éléments $A_{M}$, et on ne sait pas comment les relier aux $T_{M}$ ; plus généralement, on ne sait pas comment passer de la présentation de Coxeter $\langle T_{0},T_{1},\ldots,T_{n-1}\rangle$ à la présentation <<~symétrique~>> $\langle T_{1},\ldots,T_{n-1},\varpi_{1},\ldots,\varpi_{r}\rangle$ dans une algèbre d'Hecke cyclotomique.\medskip

\begin{example}
Supposons $n=2$ et $r=2$ ; $\wsym_{2,2}$ est le groupe des permutations signées de taille $2$, et il contient $2^{2}\,2!=8$ éléments. L'algèbre d'Ariki-Koike $\IH(\wsym_{2,2})$ est engendrée par les éléments $T_{0},T_{1}$, ou, alternativement, par les éléments $T_{1},\varpi_{1},\varpi_{2}$. Exprimons les éléments de base associés à la première présentation de $\IH(\wsym_{2,2})$ en fonction des générateurs pour la seconde présentation :
\begin{align*}T_{0}&=\varpi_{1}+\frac{q-1}{(u_{2}-u_{1})^{2}}\,u_{1}\,(\varpi_{2}-u_{2})\,(\varpi_{1}-u_{1})\,T_{1}\quad;\\
T_{1}T_{0}&=\varpi_{2}\,T_{1}+\frac{q-1}{(u_{2}-u_{1})^{2}}\,u_{2}\,(\varpi_{2}-u_{1})\,(\varpi_{1}-u_{2})\quad;\\
T_{0}T_{1}T_{0}&=\varpi_{1}\,\varpi_{2}\,T_{1}+\frac{q-1}{(u_{2}-u_{1})^{2}}\,u_{1}u_{2}\,\big[(\varpi_{2}-u_{2})\,(\varpi_{1}-u_{1})+(\varpi_{2}-u_{1})\,(\varpi_{1}-u_{2})\big]T_{1}^{2}\,.
\end{align*}
Ces formules permettent d'exprimer tout élément de la base $\{1$, $T_{0}$, $T_{1}T_{0}$, $T_{0}T_{1}T_{0}$, $T_{1}$, $T_{0}T_{1}$, $T_{1}T_{0}T_{1}$, $T_{0}T_{1}T_{0}T_{1}\}$ comme combinaison d'éléments de la forme $\varpi_{1}^{c_{1}}\,\varpi_{2}^{c_{2}}\,T_{1}^{n_{1}}$. On peut ensuite inverser ces relations, ce qui permet en particulier de calculer les espérances des caractères $\chi^{\Lambda}(q,A_{M})$ sous la $q$-mesure de Plancherel (on spécialise $u_{1}$ en $q$ et $u_{2}$ en $-1$). Mais malheureusement, ces expressions sont extraordinairement compliquées ; par exemple, 
$$\esper[\chi^{\Lambda}(q,A_{(2),\emptyset})]=\frac{q\,(q-1)\,(q^8 - 5q^7 + 9q^6 - 6q^5 + 7q^3 + 5q^2 + 4q + 1)}{2\,(q+1)^{3}\,(q^6 - 3q^5 + 4q^4 - 3q^3 - 2q^2 - 1)}\,.$$
\end{example}\bigskip

Ainsi, la formule de Shoji ne permet pas de résoudre le problème des $q$-caractères de type B, car elle est donnée pour une présentation de $\IH_{q}(\wsym_{n})$ qui est en pratique impossible à relier à la présentation de Coxeter. Pour autant, l'idée d'utiliser une dualité de type Schur-Weyl pour calculer des caractères d'algèbres d'Hecke est sans doute la bonne ; ainsi, il est vraisemblable qu'il soit possible d'établir une $q$-formule de Frobenius en type B et pour la présentation de Coxeter en utilisant comme dans \cite{HR99} la dualité entre $\IH_{q}(\wsym_{n})$ et le groupe symplectique $\Sp(2n,\For_{q})$ pour l'action sur le module de la variété de drapeaux symplectique.\bigskip
 
\begin{remark}
Le cadre des algèbres d'Ariki-Koike permet de définir des généralisations cyclotomiques des $q$-mesures de Plancherel de type $A$ et $B$. Ainsi, considérons la trace canonique de l'algèbre d'Ariki-Koike $\IH(\wsym_{n,r})$ définie par $\tau_{n,r}(1)=1$ et $\tau_{n,r}(T_{w})=0$ pour tout $w \in \wsym_{n,r}$ différent du neutre. Cette trace est un barycentre des caractères irréductibles $\chi^{\Lambda}(q)$:
$$\tau_{n,r}(\cdot)=\sum_{\Lambda \in \ym_{n,r}} M_{n,r,q}(\Lambda)\,\chi^{\Lambda}(q,\cdot)\,.$$
On définit ainsi la \textbf{$q$-mesure de Plancherel $r$-cyclotomique} $M_{n,r,q}$, et l'on retrouve les $q$-mesures de Plancherel de type $A$ et $B$ lorsque $r=1,2$ (on retrouve aussi les mesures de Plancherel des groupes $\wsym_{n,r}$ lorsque $q=1$). Malheureusement, pour $r\geq 3$, on ne connaît pas de formule de type équerre pour les valeurs de $M_{n,r,q}$ ; d'autre part, comme en type B, on ne connaît pas de modèle combinatoire correspondant sur les permutations $r$-colorées de $\wsym_{n,r}$. En l'état, on ne peut donc même pas proposer d'exemple de $r$-uplets de partitions aléatoires sous ces mesures cyclotomiques.
\end{remark}

\section[Mesures d'induction parabolique et algèbres d'Hecke généralisées]{Mesures d'induction parabolique et algèbres d'Hecke\\ généralisées}\label{lusztig}
Pour conclure notre étude asymptotique des $q$-mesures de Plancherel, nous souhaitions décrire le cadre le plus général dans lequel il semble possible de généraliser les résultats des trois derniers chapitres et d'établir la concentration gaussienne de mesures de Plancherel. Dans ce qui suit, on fixe un groupe de Chevalley $G^{F}$ (non tordu, à centre connexe) défini sur $\For_{q}$. Une représentation irréductible de $G^{F}$ est dite \textbf{cuspidale} si elle n'apparaît dans aucune représentation obtenue par induction d'Harish-Chandra à partir d'une représentation irréductible d'un sous-groupe de Lévi strict $L^{F} \subsetneq G^{F}$ (\emph{cf.} \cite[\S4]{Car92}). Fixons un sous-groupe parabolique rationnel $P^{F}$ et une décomposition de Lévi rationnelle $P=UL$ avec $U$ radical unipotent de $P$, et $P^{F}=U^{F}\rtimes L^{F}$. Ce cadre inclut en particulier le cas où $P^{F}=B^{F}$ est un sous-groupe de Borel rationnel et $L^{F}=T^{F}$ est un tore maximal scindé. Soit $\rho$ un caractère irréductible cuspidal de $L^{F}$ ; lorsque $L^{F}=T^{F}$ est un tore maximal scindé, ceci revient simplement à demander que $\rho$ soit un caractère de dimension $1$, car $T^{F}$ n'a pas de sous-groupe de Lévi strict. On s'intéresse au $G^{F}$-module 
$$\Ind_{L^{F}}^{G^{F}}(\rho)$$
obtenu par induction parabolique d'Harish-Chandra. La mesure de probabilité sous-jacente au sens de la définition \ref{defplancherel} est appelée \textbf{mesure d'induction parabolique} du caractère cuspidal $\rho$ ; on la note $M_{\rho}$.
\begin{examples}
Si $G^{F}=\GL(n,\For_{q})$, $L^{F}=(\For_{q}^{\times})^{n}$ et $\rho=1$ est le caractère trivial, on a déjà expliqué qu'on obtenait la $q$-mesure de Plancherel de type A. Le cas de la $q$-mesure de Plancherel de type B est tout à fait analogue avec $G^{F}=\Sp(2n,\For_{q})$.
\end{examples}
\bigskip\bigskip

Ceci étant, un théorème conjecturé par Springer et démontré par Howlett, Lehrer et Geck  permet une indexation des composantes irréductibles du module induit $\Ind_{L^{F}}^{G^{F}}(\rho)$ tout à fait analogue à celle donnée par le théorème d'Iwahori \ref{iwahorihecke} dans le cas où $\rho$ est le caractère trivial d'un tore scindé (auquel cas $\Ind_{L^{F}}^{G^{F}}(\rho)=\C[G^{F}/B^{F}]$). Introduisons le \textbf{groupe de ramification}
$$W(\rho)=\{w \in N(L^{F})/L^{F}\,\,|\,\,\rho \circ w = \rho\}\,.$$
Par exemple, si $L^{F}=T^{F}$ est un tore scindé maximal et si $\rho=1$ est le caractère trivial, alors $W(1)$ est simplement le groupe de Weyl de $G^{F}$.
\begin{theorem}[Lusztig, Howlett-Lehrer, Geck, \cite{HL80,Lus84,Geck93}]
Le groupe de ramification $W(\rho)$ est un groupe de Coxeter fini engendré par des réflexions $s\in S(\rho)$. L'algèbre 
$$\hendo_{G^{F}}(\Ind_{L^{F}}^{G^{F}}(\rho))$$
commutante de l'action de $G^{F}$ sur le module parabolique induit est abstraitement isomorphe\footnote{Initialement, R. B. Howlett et G. I. Lehrer avaient montré que l'algèbre commutante était isomorphe à $\C W(\rho)_{\mu}$, où $\mu$ est un cocycle déformant les relations dans l'algèbre $\C W(\rho)$. Mais Lusztig, puis Geck ont montré que pour un groupe algébrique réductif à centre connexe, ce cocycle $\mu$ était toujours trivial. La présentation de l'algèbre commutante donnée dans \cite{HL80} devient sous cette hypothèse celle d'une algèbre d'Hecke à plusieurs paramètres.} à l'algèbre de groupe $\C W(\rho)$, et plus précisément, c'est une algèbre d'Hecke à plusieurs paramètres de $W(\rho)$ :
\begin{align*}\hendo_{G^{F}}(\Ind_{L^{F}}^{G^{F}}(\rho))=\langle T_{s},\,\,s \in S(\rho)\rangle_{\C} \,,\text{ avec :}\quad&(T_{s})^{2}=(q^{c_{s}}-1)T_{s}+q^{c_{s}}\,\,;\\
&T_{w}T_{w'}=T_{ww'}\, \text{ si }\ell(ww')=\ell(w)+\ell(w')\,.
\end{align*}
Par suite, $\Ind_{L^{F}}^{G^{F}}(\rho)$ se décompose en somme directe de $\big(G^{F},\IH(W(\rho),S(\rho),\{q^{c_{s}}\})\big)$-bimodules irréductibles :
$$_{G^{F} \curvearrowright }\big\{\Ind_{L^{F}}^{G^{F}}(\rho)\big\}_{\curvearrowleft\, \IH(W(\rho),S(\rho),\{q^{c_{s}}\})} =\sum_{\lambda \in \widehat{W}(\rho)} \left({}_{G^{F} \curvearrowright}M^\lambda \right)\otimes_{\C} \left(V^\lambda(q)_{\curvearrowleft\, \IH(W(\rho),S(\rho),\{q^{c_{s}}\})}\right).$$
Les composantes irréductibles de $\Ind_{L^{F}}^{G^{F}}(\rho)$ sont donc en bijection avec les caractères irréductibles de $W(\rho)$, et elles ont pour multiplicités les dimensions de ces caractères irréductibles.
\end{theorem}
\bigskip
\bigskip

Ainsi, exactement comme dans les cas précédemment étudiés, on peut remplacer un module pour le groupe de Chevalley $G^{F}$ par un module pour l'\textbf{algèbre d'Hecke généralisée} $\IH(W(\rho),S(\rho),\{q^{c_{s}}\})$, et de plus, le caractère normalisé $\tau$ associé à l'action de cette algèbre sur $\Ind_{L^{F}}^{G^{F}}(\rho)$ est la trace symétrique canonique de l'algèbre d'Hecke (voir \cite[théorème 8.6]{Lus84}) :
$$\tau(T_{w})=\begin{cases}
1 & \text{si }w=e_{W(\rho)},\\
0 & \text{sinon}.
\end{cases}$$
Ceci signifie que l'espérance dans l'espace de probabilité non commutatif $$(\IH(W(\rho),S(\rho),\{q^{c_{s}}\}),\tau)$$ est très facile à calculer. On a donc une famille d'observables des irréductibles $\lambda$ dont les moments sont accessibles, et le problème devient le même que dans les chapitres \ref{tool}, \ref{plancherel} et \ref{qplancherelmeasure} : relier les valeurs des caractères
$$\chi^{\lambda}(q,T_{w})$$
avec $\lambda \in \widehat{W}(\rho)$ à des caractéristiques <<~géométriques~>> des indices $\lambda$ des composantes irréductibles. En particulier, lorsque $W$ est un groupe de Coxeter produit de groupes classiques de type A, B ou D, les indices $\lambda$ correspondent à des familles de partitions ; on peut alors tenter d'exprimer les $\chi^{\lambda}(q,T_{w})$ en fonction des moments de Frobenius des partitions de ces familles.\bigskip
\bigskip

Il est à peu près clair que les calculs ne peuvent pas être menés à bout dans un contexte aussi général ; le simple cas des algèbres d'Hecke (à un paramètre) de type B a montré à quel point l'expression des caractères $\chi^{\lambda}(q,T_{w})$ en fonction d'autres observables des indices $\lambda$ pouvait être délicate. On peut néanmoins conjecturer le caractère universel des phénomènes de concentration gaussienne dans la classe des mesures d'induction parabolique ; à titre d'exemple, concluons ce chapitre en examinant toutes les mesures d'induction parabolique obtenues à partir de caractères de tores scindés d'un groupe $G=\GL(n,\For_{q})$. Si $T=(\For_{q}^{\times})^{n}$ et si $\rho$ est un caractère irréductible de $T$, alors $\rho$ peut s'écrire sous la forme :
$$\rho=\theta_{1}\otimes \theta_{2}\otimes \cdots \otimes \theta_{n}\,,$$
où chaque $\theta_{i}$ est un caractère du groupe cyclique $\For_{q}^{\times}$. Or, le dual $L_{1}$ de $\For_{q}^{\times}$ est lui aussi cyclique de cardinal $q-1$, donc chaque $\theta_{i}$ est un $\theta^{k_{i}}$, où $\theta$ est un générateur de $L_{1}$ (avec les notations du chapitre \ref{general}), et $k_{i} \in \lle 1,q-1\rre$. Quitte à conjuguer le tore par un élément du groupe de Weyl $W=\sym_{n}$, on peut donc supposer que 
$$\rho=(\theta^{1})^{\otimes n_{1}} \otimes (\theta^{2})^{\otimes n_{2}} \otimes\cdots \otimes (\theta^{q-1})^{\otimes n_{q-1}}$$
avec $n=n_{1}+n_{2}+\cdots+n_{q-1}$. Alors, il est aisé de voir que $\sigma \in \sym_{n}$ laisse stable $\rho$ si et seulement s'il appartient au sous-groupe de Young 
$$W(\rho)=\sym_{n_{1}}\times \sym_{n_{2}} \times \cdots \times \sym_{n_{q-1}}\,.$$
Ainsi, les composantes irréductibles du module parabolique $R_{T}^{G}$ sont indexées par les éléments de $\ym_{n_{1}}\times \ym_{n_{2}}\times \cdots \times \ym_{n_{q-1}}$, et d'autre part, avec les notations du chapitre \ref{general}, si $\rho$ est de type $(n_{1},\ldots,n_{q-1})$, alors $R_{T}^{G}$ n'est nul autre que $R^{\bbmu}$, où $\bbmu$ est la polypartition duale
$$[\theta] : 1^{n_{1}}\quad;\quad[\theta^{2}] : 1^{n_{2}}\quad;\quad\cdots\quad;\quad[\theta^{q-1}] : 1^{n_{q-1}}\,.$$
Ici, $[\theta^{k}]$ désigne la classe dans $L/\Gal$ de l'élément de $L=\varinjlim_{n \to \infty} L_{n}$ donné par l'appariement $\scal{\theta^{k}}{\cdot}_{1}=\theta^{k}(\cdot)$ (cette identité est suffisante pour décrire un élément de $L$). Dans $\Lambda(\For_{q})$, $\Ind^{G}_{T}$ correspond donc à la fonction $B^{\bbmu}$, c'est-à-dire, d'après la formule de changement de base donnée dans le paragraphe \ref{delignelusztig}, à
$$\sum_{\lambda_{1}\in \ym_{n_{1}},\ldots,\lambda_{q-1}\in \ym_{n_{q-1}}} \left(\prod_{i=1}^{q-1}\varsigma^{\lambda_{i}}(1^{n_{i}})\right)\,S^{\{[\theta] : \lambda_{1},\ldots,[\theta^{q-1}]: \lambda_{q-1}\}}\,.$$
Autrement dit, en tant que $\GL(n,\For_{q})$-module, $\Ind^{G}_{T}$ admet pour décomposition en irréductibles :
$$\Ind^{G}_{T}=\sum_{\lambda^{(1)}\in \ym_{n_{1}},\ldots,\lambda^{(q-1)}\in \ym_{n_{q-1}}} \left(\prod_{i=1}^{q-1}\dim \lambda^{(i)}\right)\,V^{\{[\theta] : \lambda^{(1)},\ldots,[\theta^{q-1}]: \lambda^{(q-1)}\}}\,.$$
Les multiplicités $\prod_{i=1}^{q-1}\dim \lambda^{(i)}$ sont les dimensions des $\IH(W(\rho))$-modules irréductibles intervenant dans la décomposition de $\Ind^{G}_{T}$ en somme de $(\GL(n,\For_{q}),\IH(W(\rho)))$-bimodules irréductibles ; on peut d'autre part montrer que dans le cas qui nous intéresse, les paramètres de $\IH(W(\rho))$ sont tous égaux à $q$, de sorte que l'algèbre d'Hecke est exactement la sous-algèbre de Young $\IH_{q}(\sym_{n_{1}}\times \sym_{n_{2}} \times \cdots \times \sym_{n_{q-1}}) \subset \IH_{q}(\sym_{n})$. D'après \cite[chapitre 4, \S6]{Mac95}, la dimension de la composante irréductible $V^{\lambda^{(1)},\ldots,\lambda^{(q-1)}}$ est (à conjugaison des diagrammes près)
$$\{n!\}_{q} \,\,\prod_{i=1}^{q-1}\frac{q^{b(\lambda^{(i)})}}{\prod_{x \in \lambda^{(i)}} \{h(x)\}_{q}}\,.$$
Enfin, la dimension de $R_{T}^{G}$ est $\{n!\}_{q}$ : comme le tore est scindé, on est en effet dans la même situation que pour la $q$-mesure de Plancherel de type A. On en déduit que la mesure d'induction parabolique sur l'espace $\ym_{n_{1}}\times \cdots \times \ym_{n_{q-1}}$ paramétrant les composantes irréductibles de $\Ind^{G}_{T}$ n'est nulle autre que 
$$M_{\rho}(\lambda^{(1)},\ldots,\lambda^{(q-1)})=\prod_{i=1}^{q-1}M_{n_{i},q}(\lambda^{(i)})\,.$$
Autrement dit :
\begin{proposition}[Mesures d'induction parabolique des caractères des tores scindés des grou\-pes linéaires finis]
Soit $\rho$ un caractère irréductible d'un tore scindé de $\GL(n,\For_{q})$, et $(n_{1},\ldots,n_{q-1})$ les exposants des caractères $\theta^{i}$ dans la décomposition en produit tensoriel de $\rho$. La mesure d'induction parabolique $M_{\rho}$ est le produit de $q$-mesures de Plancherel de type A de paramètres $n_{1},\ldots,n_{q-1}$.
\end{proposition}
\noindent Dans ce nouveau cas, la concentration gaussienne est donc trivialement vraie, puisque les partitions $\lambda^{(1)},\ldots,\lambda^{(q-1)}$ sont tirées indépendamment suivant des mesures asymptotiquement gaussiennes au sens du théorème \ref{secondasymptoticqplancherel}. Ce résultat va dans le sens de notre conjecture d'universalité des phénomènes gaussiens pour les mesures d'induction parabolique.

\pagestyle{empty}
\clearpage
~

\pagestyle{fancy}
\fancyhead{}
\fancyfoot{}
\fancyfoot[C]{\thepage}
\renewcommand{\chaptermark}[1]{\markboth{\chaptername\ \thechapter.\ #1.}{}} 
\renewcommand{\sectionmark}[1]{\markright{\thesection.\ #1.}} 
\fancyhead[RO]{\rightmark}
\fancyhead[LE]{\leftmark}
\setlength{\headheight}{15.5pt}

\part{Asymptotique des mesures de Schur-Weyl et des mesures de Gelfand}

Dans cette troisième partie, nous revisitons la théorie asymptotique des représentations des groupes symétriques, en nous penchant sur deux familles de représentations :\vspace{2mm}
\begin{itemize}
\item les représentations de Schur-Weyl, qui correspondent aux actions des groupes symétriques par permutations des lettres des tenseurs simples dans des produits tensoriels.\vspace{2mm}
\item les représentations de Gelfand, qui correspondent au cas où toutes les multiplicités des composantes irréductibles sont égales à $1$.\vspace{2mm}
\end{itemize}
Dans ces deux cadres, nous établissons des lois des grands nombres et des résultats de déviation gaussienne semblables à ceux exposés dans le chapitre \ref{plancherel}, et l'outil principal mis en jeu dans les preuves de ces résultats est de nouveau l'algèbre d'observables présentée dans le chapitre \ref{tool}. \bigskip

On constate ainsi que l'étude asymptotique de partitions aléatoires peut à chaque fois être décomposée en trois étapes :\vspace{2mm}
\begin{enumerate}
\item On détermine une \textbf{renormalisation} adaptée à la forme typique des grandes partitions aléatoires, qu'on peut deviner en simulant un processus aléatoire sur le graphe de Young dont les lois marginales sont celles que l'on souhaite étudier. On choisit ensuite une graduation de l'algèbre $\obs$ adaptée à cette renormalisation : cette graduation était le degré canonique pour les mesures des algèbres d'Hecke, et sera le poids pour les mesures de Schur-Weyl de paramètre $\alpha \geq 1/2$, le $\alpha$-degré pour les mesures de Schur-Weyl de paramètre $\alpha < 1/2$, et le degré de Kerov pour les mesures de Gelfand.\vspace{2mm}
\item On trouve une \textbf{base graduée d'observables} dont les espérances ont une propriété de \textbf{factorisation asymptotique}, et on relie ces observables à des observables de nature analytique telles que les moments de Frobenius ou les cumulants libres. La propriété de factorisation asymptotique implique la convergence en probabilité des observables renormalisées, ce qui permet le plus souvent d'établir une loi des grands nombres.\vspace{2mm}
\item Pour les phénomènes de \textbf{déviation gaussienne}, on utilise la théorie des \textbf{cumulants d'observables} due à \'Sniady, et on en déduit le caractère gaussien de certaines observables de diagrammes, puis éventuellement le caractère gaussien de la forme des diagrammes aléatoires.\vspace{2mm}
\end{enumerate}
Dans un contexte de théorie des représentations, la base graduée d'observables dont on peut espérer calculer les espérances est celle des caractères centraux, et ce compte tenu de la remarque faite à la fin de la section \ref{lskv}. C'est effectivement cette base qui jouait un rôle pivot pour les mesures des algèbres d'Hecke, et il en sera de même pour tous les cas étudiés dans cette troisième partie.
\bigskip
\bigskip

Détaillons maintenant les résultats obtenus, en commen\c cant par les mesures de Schur-Weyl, qui sont liées au phénomène de dualité de Schur-Weyl entre les actions de $\GL(N,\C)$ et de $\sym_{n}$ sur $(\C^{N})^{\otimes n}$. Lorsque $N$ est de l'ordre de $n^{\alpha}$ avec $\alpha\geq 1/2$, la forme limite des diagrammes aléatoires sous ces \textbf{mesures de Schur-Weyl} a été déterminée par P. Biane (voir \cite{Bia01a}). Nous montrons dans le chapitre 10 que le théorème central limite de Kerov reste vrai dans ce cadre --- c'est le résultat de l'article \cite{Mel10b} --- et nous étudions également le cas $\alpha < 1/2$. Dans ce dernier cas, nous présentons une famille de filtrations d'algèbres sur $\obs$ qui interpolent la filtration des poids et la filtration des degrés, et sont adaptées à l'étude d'asymptotiques <<~intermédiaires~>>. Ainsi, sous une mesure de Schur-Weyl de paramètre $\alpha<1/2$, une partition a asymptotiquement des lignes de l'ordre de $n^{1-\alpha}$ et des colonnes de l'ordre de $n^{\alpha}$. Dans ce même cadre, on peut tenter d'adapter la théorie de \'Sniady aux filtrations des $\alpha$-degrés pour obtenir des résultats de concentration gaussienne. Malheureusement, ces arguments ne permettent pas de traiter l'ensemble de l'intervalle $]0,1/2[$, et l'asymptotique des fluctuations des mesures de Schur-Weyl de paramètre inférieur à $1/2$ reste en partie un problème ouvert.
\bigskip
\bigskip

Le théorème central limite de Kerov s'applique également aux \textbf{mesures de Gelfand} ; sous ces mesures, la probabilité d'une partition est proportitionnelle à $\dim \lambda$ (au lieu de $(\dim \lambda)^{2}$ pour les mesures de Plancherel). La forme limite des partitions renormalisées sous ces mesures est la même que dans le théorème \ref{firstasymptoticplancherel}, mais les fluctuations sont plus grandes ; néanmoins, la déviation est de nouveau essentiellement décrite par le processus gaussien généralisé de Kerov, voir le théorème \ref{scaledkerovgaussian}.

\chapter{Asymptotique des mesures de Schur-Weyl}\label{schurweylmeasure}

Entre la renormalisation <<~isotrope~>> en $1/\sqrt{n}$ des diagrammes de Young sous la mesure de Plancherel et la renormalisation en $1/n$ des lignes des diagrammes de Young sous une $q<1$-mesure de Plancherel, il est naturel de se demander s'il existe des familles de représentations des groupes symétriques donnant lieu à des renormalisations <<~intermédiaires~>>, c'est-à-dire qu'avec grande probabilité les lignes ont pour ordre de grandeur $n^{1/2+\eps}$ et les colonnes ont pour ordre de grandeur $n^{1/2-\eps}$. Des phénomènes asymptotiques de ce type peuvent être observés dans le contexte des \textbf{mesures de Schur-Weyl}, qui sont les mesures de probabilité associées aux représentations
 des groupes $\sym_{n}$ sur des produits tensoriels d'espaces $V^{\otimes n}$. Plus précisément, si l'on considère la mesure de Plancherel associée à la représentation de $\sym_{n}$ sur le produit tensoriel $(\C^{N})^{\otimes n}$, alors on observe trois régimes asymptotiques distincts :\vspace{2mm}
 \begin{enumerate}
 \item Si $\log N/\log n \simeq \alpha >1/2$, alors le comportement asymptotique est identique à celui observé pour la mesure de Plancherel standard. \vspace{2mm}
 \item Si $\log N/\log n \simeq 1/2$, la forme limite des diagrammes renormalisés dans les deux directions suivant un facteur $1/\sqrt{n}$ dépend du paramètre $c=\sqrt{n}/N$. La déviation des diagrammes par rapport à la forme limite $\Omega_{c}$ est décrite dans un intervalle de taille $4$ par le même processus gaussien généralisé que dans le cas $\alpha >1/2$.\vspace{2mm}
 \item Enfin, si $\log N/\log n \simeq \alpha <1/2$, alors on observe un comportement asymptotique intermédiaire, avec des lignes qui ont pour ordre de grandeur $n^{1-\alpha}$ et des colonnes qui ont pour ordre de grandeur $n^{\alpha}$.\vspace{2mm}
 \end{enumerate}
Ces trois régimes sont étudiés en détail dans les sections \ref{schurplus} et \ref{schurminus}, et on y démontre les résultats nouveaux \ref{schurweylplus} et \ref{schurweylminus} (le premier est l'objet de l'article \cite{Mel10b}, et le second théorème a été obtenu en collaboration avec V. Féray). À l'aide d'une correspondance RSK généralisée, on peut réinterpréter tous ces résultats en termes de longueurs des plus longs sous-mot croissants d'un mot aléatoire de longueur $n$ sur un alphabet de taille $N$ ; on rappelle cette correspondance dans la section \ref{combinatoricschurweyl}.\bigskip

\section{Mesures de Schur-Weyl}\label{combinatoricschurweyl}
Soit $\alpha$ et $c$ deux réels strictement positifs. La \textbf{mesure de Schur-Weyl} de paramètres $(n,\alpha,c)$ est la mesure sur les partitions de taille $n$ associée à la représentation
$$ (\C^{N})^{\otimes n}\curvearrowleft \sym_{n}\,,$$
où $ n^{\alpha}\simeq cN$ et $\sym_{n}$ agit comme dans le paragraphe \ref{shoji}, c'est-à-dire par permutation des lettres dans les tenseurs simples. Nous noterons cette mesure $SW_{n,\alpha,c}$ ; elle dépend bien sûr de la valeur précise de la dimension $N$, mais nous verrons que pour l'étude asymptotique, la donnée des paramètres $\alpha$ et $c$ suffit. Compte tenu du théorème \ref{schurweylduality}, la mesure $SW_{n,\alpha,c}$ charge seulement les partitions de longueur inférieure à $N$, et elle s'écrit :
$$SW_{n,\alpha,c}[\lambda]=\frac{\dim M^{\lambda}\times \dim \lambda}{N^{n}}\,.$$
où $M^{\lambda}$ est le $\GL(N,\C)$-module irréductible de plus haut poids $\lambda$.\bigskip
\bigskip

Notons que les mesures de Schur-Weyl rentrent elles aussi dans le cadre des mesures de Schur (\S\ref{schurmeasure}). Ainsi, la poissonisée de paramètres $(\theta,N)$ des mesures de Schur-Weyl est la mesure de probabilité sur $\ym$ définie par
$$SW_{\mathcal{P}(\theta),N}[\lambda]=\frac{\theta^{|\lambda|}\,\E^{-\theta}}{|\lambda|!}\times \frac{\dim M^{\lambda} \times \dim \lambda}{N^{|\lambda|}},$$
et c'est aussi la mesure de Schur de paramètres $t=(\theta/N,0,0,\ldots)$ et $t'=(N,N/2,N/3,\ldots)$. Cette approche est suivie dans l'article \cite{BO07} ; de nouveau, nous préférerons l'approche des observables de diagrammes.\bigskip\bigskip

Soit $A=\{1,2,\ldots,N\}$ un alphabet de taille $N$, et $m=m_{1}m_{2}\cdots m_{n}$ un mot de taille $n$ sur l'alphabet $A$. Comme dans le paragraphe \ref{ulamrsk}, on peut construire deux tableaux $P(m)$ et $Q(m)$ de taille $n$ en appliquant l'algorithme d'insertion de Schensted au mot $m$.\medskip

\begin{example}
Supposons $N=5$, $n=9$ et $m=233154243$. Si l'on applique l'algorithme décrit page \pageref{rskalgo}, alors on obtient les deux tableaux :
$$\young(5,234,12334)\qquad\text{et}\qquad\young(7,469,12358)\,\,.$$
Par construction, le premier tableau $P(m)$ est \textbf{semi-standard}, c'est-à-dire qu'il est croissant suivant les lignes et strictement croissant suivant les colonnes. Le second tableau $Q(m)$ est pour sa part un tableau standard. On peut retrouver $P(m)$ en appliquant l'algorithme du jeu de taquin au ruban semi-standard obtenu à partir du mot $m$, voir la figure \ref{taquin}.
\figcapt{\psset{unit=1mm}
\pspicture(-10,0)(125,50)
\psline(0,5)(0,15)(5,15)(5,10)(25,10)(25,0)(10,0)(10,5)(0,5)\psline(0,10)(5,10)(5,5)\psline(10,10)(10,5)(25,5)\psline(15,10)(15,0)\psline(20,10)(20,0)
\rput(2.5,12.5){$5$} \rput(2.5,7.5){$2$} \rput(7.5,7.5){$3$} \rput(12.5,7.5){$3$} \rput(17.5,7.5){$4$} \rput(22.5,7.5){$4$} \rput(12.5,2.5){$1$} \rput(17.5,2.5){$2$} \rput(22.5,2.5){$3$} 
\psline{->}(35,5)(40,5)
\psline{->}(35,30)(40,30)
\psline{->}(-15,5)(-10,5)
\psline(50,5)(50,15)(55,15)(55,10)(70,10)(70,5)(75,5)(75,0)(55,0)(55,5)(50,5)\psline(50,10)(55,10)(55,5)(65,5)(65,0)\psline(60,10)(60,0)\psline(65,10)(65,5)(70,5)(70,0)
\rput(52.5,12.5){$5$} \rput(52.5,7.5){$2$} \rput(57.5,7.5){$3$} \rput(62.5,7.5){$3$} \rput(67.5,7.5){$4$} \rput(72.5,2.5){$4$} \rput(57.5,2.5){$1$} \rput(62.5,2.5){$2$} \rput(67.5,2.5){$3$}
\psline{->}(85,5)(90,5)
\psline{->}(85,30)(90,30)
\psline(100,0)(125,0)(125,5)(115,5)(115,10)(105,10)(105,15)(100,15)(100,0)\psline(100,10)(105,10)(105,0)\psline(110,10)(110,0)\psline(100,5)(115,5)(115,0)\psline(120,5)(120,0)
\rput(102.5,12.5){$5$} \rput(102.5,7.5){$2$} \rput(107.5,7.5){$3$} \rput(112.5,7.5){$4$}  \rput(122.5,2.5){$4$} \rput(107.5,2.5){$2$} \rput(112.5,2.5){$3$} \rput(117.5,2.5){$3$} \rput(102.5,2.5){$1$}

\psline(20,25)(25,25)(25,35)(20,35)(20,45)(15,45)(15,50)(0,50)(0,45)(10,45)(10,40)(15,40)(15,30)(20,30)(20,25)\psline(5,45)(5,50)\psline(10,50)(10,45)(15,45)(15,40)(20,40)\psline(15,35)(20,35)(20,30)(25,30)
\rput(2.5,47.5){$2$}\rput(7.5,47.5){$3$}\rput(12.5,47.5){$3$}\rput(12.5,42.5){$1$}\rput(17.5,42.5){$5$}\rput(17.5,37.5){$4$}\rput(17.5,32.5){$2$}\rput(22.5,32.5){$4$}\rput(22.5,27.5){$3$}

\psline(65,25)(75,25)(75,35)(70,35)(70,40)(65,40)(65,45)(50,45)(50,40)(60,40)(60,35)(65,35)(65,25)\psline(55,45)(55,40)\psline(60,45)(60,40)(65,40)(65,35)(70,35)(70,25)\psline(65,30)(75,30)\rput(52.5,42.5){$2$}\rput(57.5,42.5){$3$}\rput(62.5,42.5){$3$}\rput(62.5,37.5){$1$}\rput(67.5,37.5){$5$}\rput(67.5,32.5){$4$}\rput(67.5,27.5){$2$}\rput(72.5,32.5){$4$}\rput(72.5,27.5){$3$}

\psline(100,35)(100,40)(120,40)(120,35)(125,35)(125,25)(115,25)(115,30)(110,30)(110,35)(100,35)\psline(105,40)(105,35)\psline(110,40)(110,35)(120,35)(120,25)\psline(115,40)(115,30)(125,30)
\rput(102.5,37.5){$2$}\rput(107.5,37.5){$3$}\rput(112.5,37.5){$3$}\rput(117.5,37.5){$5$}\rput(112.5,32.5){$1$}\rput(117.5,32.5){$4$}\rput(122.5,32.5){$4$}\rput(117.5,27.5){$2$}\rput(122.5,27.5){$3$}
\endpspicture
}{Algorithme du jeu de taquin et tableau semi-standard associé à un mot.\label{taquin}}{Algorithme du jeu de taquin et tableau semi-standard associé à un mot}
\end{example}\bigskip

Comme l'algorithme d'insertion de Schensted peut être inversé, on obtient l'identité combinatoire suivante :
$$N^{n}=\sum_{\substack{\lambda \in \ym_{n}\\ \ell(\lambda) \leq N}} \bigg|\left\{\substack{\text{tableaux semi-standards de forme $\lambda$}\\ \text{et à entrées dans $\lle 1,N\rre$}} \right\}\bigg|\times \bigg|\left\{\substack{\text{tableaux standards}\\ \text{de forme }\lambda}\right\}\bigg|\,.$$
De plus, on peut montrer que la dimension du $\GL(N,\C)$-module $M^{\lambda}$ est précisément le nombre de tableaux semi-standards de forme $\lambda$ et à entrées dans $\lle1,N\rre$ (voir \cite[appendice A, \S8]{Mac95}). Enfin, comme dans le cas des mots des permutations, la première ligne de la partition $\lambda(m)$ qui est la forme commune des tableaux $P(m)$ et $Q(m)$ associés à $m$ a pour taille la longueur d'un plus long sous-mot croissant dans $m$. Par conséquent :
\begin{proposition}[Interprétation combinatoire des mesures de Schur-Weyl]
La mesure de Schur-Weyl de paramètres $(n,\alpha,c)$ a pour expression
$$SW_{n,\alpha,c}[\lambda]=\frac{\bigg|\left\{\substack{\text{tableaux semi-standards de forme $\lambda$}\\ \text{et à entrées dans $\lle 1,N\rre$}} \right\}\bigg|\times \bigg|\left\{\substack{\text{tableaux standards}\\ \text{de forme }\lambda}\right\}\bigg|}{N^{n}}\,.$$
De plus, la loi de la longueur du plus long sous-mot croissant d'un mot choisi aléatoirement parmi les $N^{n}$ mots de taille $n$ sur l'alphabet $A=\lle 1,N\rre$ est la même que la loi de la première part d'une partition sous la mesure de Schur-Weyl $SW_{n,\alpha,c}$.
\end{proposition}\bigskip

Ainsi, comme pour la mesure de Plancherel et la $q$-mesure de Plancherel de type A, on dispose d'un modèle combinatoire des mesures de Schur-Weyl. Une autre propriété importante pour la suite est le calcul des traces des permutations pour ces représentations (voir le paragraphe \ref{shoji}) : 
$$\tr(\sigma_{\mu})=p_{\mu}(1,1,\ldots,1)=N^{\ell(\mu)}\qquad;\qquad SW_{n,\alpha,c}[\sigma_{\mu}]=\frac{\tr(\sigma_{\mu})}{N^{n}}=N^{\ell(\mu)-|\mu|}=\frac{1}{N^{|\sigma_{\mu}|}}\,.$$
En particulier, les traces ont une propriété de factorisation exacte, c'est-à-dire que si $\sigma_{1}$ et $\sigma_{2}$ sont deux permutations à cycles disjoints, alors :
 $$SW_{n,\alpha,c}[\sigma_{1}\sigma_{2}]=SW_{n,\alpha,c}[\sigma_{1}]\,SW_{n,\alpha,c}[\sigma_{2}]\,.$$
Ceci implique la nullité des cumulants joints $k_{n}(\sigma_{1},\sigma_{2},\ldots,\sigma_{r})$ pour $r\geq 2$ et des permutations à cycles disjoints.\bigskip

\section{Asymptotique pour $\alpha \geq 1/2$}\label{schurplus}
Supposons le paramètre $\alpha$ strictement supérieur à $1/2$. Alors, si $\sigma_{l}$ désigne un cycle de longueur $l$, $SW_{n,\alpha,c}[\sigma_{l}]\simeq c^{l-1} \,n^{-\alpha\,(l-1)}$, donc
$$\lim_{n \to \infty} SW_{n,\alpha,c}[\sigma_{l}]\,n^{\frac{l-1}{2}}=\begin{cases}1&\text{si }l=1,\\
0&\text{sinon}. \end{cases} $$
D'autre part, on vient de voir que les cumulants joints $k_{n}(\sigma_{l_{1}},\ldots,\sigma_{l_{r}})$ avec $r \geq 2$ et des cycles disjoints étaient tous nuls. Le théorème de \'Sniady \ref{factorasymptotic} s'applique donc, et les paramètres asymptotiques sont les mêmes que pour les mesures de Plancherel usuelles des groupes $\sym_{n}$ (voir la page \pageref{factorasymptotic}). On a donc les mêmes résultats asymptotiques que dans le chapitre \ref{plancherel}  :
\begin{proposition}[Asymptotique des mesures de Schur-Weyl de paramètre $\alpha>1/2$, \cite{Bia01a}]
Si $\alpha>1/2$, alors sous les mesures de Schur-Weyl $SW_{n,\alpha,c}$, les partitions aléatoires $\lambda \in \ym_{n}$ obéissent aux lois asymptotiques des théorèmes \ref{firstasymptoticplancherel} et \ref{secondasymptoticplancherel}.
\end{proposition}\bigskip

\figcapt{
\psset{unit=1mm}
\pspicture(-45,0)(80,80)
\psline{->}(0,0)(-40,40)\psline{->}(0,0)(90,90)
\parametricplot{-1}{3}{t 27.951 mul
t 1 add 2 div  t -1 add 2 div arcsin t -1 add mul 3.14159 mul 180 div  4 t -1 add t -1 add mul -1 mul add sqrt add  3.14159 div add 27.951 mul}
\psline[linewidth=0.25pt](0,0)(-27.500,27.500)\psline[linewidth=0.25pt](0,0)(70.000,70.000)
\psline[linewidth=0.25pt](1.2500,1.2500)(0.00000,2.5000)(-1.2500,1.2500)
\psline[linewidth=0.25pt](2.5000,2.5000)(1.2500,3.7500)(0.00000,2.5000)
\psline[linewidth=0.25pt](3.7500,3.7500)(2.5000,5.0000)(1.2500,3.7500)
\psline[linewidth=0.25pt](5.0000,5.0000)(3.7500,6.2500)(2.5000,5.0000)
\psline[linewidth=0.25pt](6.2500,6.2500)(5.0000,7.5000)(3.7500,6.2500)
\psline[linewidth=0.25pt](7.5000,7.5000)(6.2500,8.7500)(5.0000,7.5000)
\psline[linewidth=0.25pt](8.7500,8.7500)(7.5000,10.000)(6.2500,8.7500)
\psline[linewidth=0.25pt](10.000,10.000)(8.7500,11.250)(7.5000,10.000)
\psline[linewidth=0.25pt](11.250,11.250)(10.000,12.500)(8.7500,11.250)
\psline[linewidth=0.25pt](12.500,12.500)(11.250,13.750)(10.000,12.500)
\psline[linewidth=0.25pt](13.750,13.750)(12.500,15.000)(11.250,13.750)
\psline[linewidth=0.25pt](15.000,15.000)(13.750,16.250)(12.500,15.000)
\psline[linewidth=0.25pt](16.250,16.250)(15.000,17.500)(13.750,16.250)
\psline[linewidth=0.25pt](17.500,17.500)(16.250,18.750)(15.000,17.500)
\psline[linewidth=0.25pt](18.750,18.750)(17.500,20.000)(16.250,18.750)
\psline[linewidth=0.25pt](20.000,20.000)(18.750,21.250)(17.500,20.000)
\psline[linewidth=0.25pt](21.250,21.250)(20.000,22.500)(18.750,21.250)
\psline[linewidth=0.25pt](22.500,22.500)(21.250,23.750)(20.000,22.500)
\psline[linewidth=0.25pt](23.750,23.750)(22.500,25.000)(21.250,23.750)
\psline[linewidth=0.25pt](25.000,25.000)(23.750,26.250)(22.500,25.000)
\psline[linewidth=0.25pt](26.250,26.250)(25.000,27.500)(23.750,26.250)
\psline[linewidth=0.25pt](27.500,27.500)(26.250,28.750)(25.000,27.500)
\psline[linewidth=0.25pt](28.750,28.750)(27.500,30.000)(26.250,28.750)
\psline[linewidth=0.25pt](30.000,30.000)(28.750,31.250)(27.500,30.000)
\psline[linewidth=0.25pt](31.250,31.250)(30.000,32.500)(28.750,31.250)
\psline[linewidth=0.25pt](32.500,32.500)(31.250,33.750)(30.000,32.500)
\psline[linewidth=0.25pt](33.750,33.750)(32.500,35.000)(31.250,33.750)
\psline[linewidth=0.25pt](35.000,35.000)(33.750,36.250)(32.500,35.000)
\psline[linewidth=0.25pt](36.250,36.250)(35.000,37.500)(33.750,36.250)
\psline[linewidth=0.25pt](37.500,37.500)(36.250,38.750)(35.000,37.500)
\psline[linewidth=0.25pt](38.750,38.750)(37.500,40.000)(36.250,38.750)
\psline[linewidth=0.25pt](40.000,40.000)(38.750,41.250)(37.500,40.000)
\psline[linewidth=0.25pt](41.250,41.250)(40.000,42.500)(38.750,41.250)
\psline[linewidth=0.25pt](42.500,42.500)(41.250,43.750)(40.000,42.500)
\psline[linewidth=0.25pt](43.750,43.750)(42.500,45.000)(41.250,43.750)
\psline[linewidth=0.25pt](45.000,45.000)(43.750,46.250)(42.500,45.000)
\psline[linewidth=0.25pt](46.250,46.250)(45.000,47.500)(43.750,46.250)
\psline[linewidth=0.25pt](47.500,47.500)(46.250,48.750)(45.000,47.500)
\psline[linewidth=0.25pt](48.750,48.750)(47.500,50.000)(46.250,48.750)
\psline[linewidth=0.25pt](50.000,50.000)(48.750,51.250)(47.500,50.000)
\psline[linewidth=0.25pt](51.250,51.250)(50.000,52.500)(48.750,51.250)
\psline[linewidth=0.25pt](52.500,52.500)(51.250,53.750)(50.000,52.500)
\psline[linewidth=0.25pt](53.750,53.750)(52.500,55.000)(51.250,53.750)
\psline[linewidth=0.25pt](55.000,55.000)(53.750,56.250)(52.500,55.000)
\psline[linewidth=0.25pt](56.250,56.250)(55.000,57.500)(53.750,56.250)
\psline[linewidth=0.25pt](57.500,57.500)(56.250,58.750)(55.000,57.500)
\psline[linewidth=0.25pt](58.750,58.750)(57.500,60.000)(56.250,58.750)
\psline[linewidth=0.25pt](60.000,60.000)(58.750,61.250)(57.500,60.000)
\psline[linewidth=0.25pt](61.250,61.250)(60.000,62.500)(58.750,61.250)
\psline[linewidth=0.25pt](62.500,62.500)(61.250,63.750)(60.000,62.500)
\psline[linewidth=0.25pt](63.750,63.750)(62.500,65.000)(61.250,63.750)
\psline[linewidth=0.25pt](65.000,65.000)(63.750,66.250)(62.500,65.000)
\psline[linewidth=0.25pt](66.250,66.250)(65.000,67.500)(63.750,66.250)
\psline[linewidth=0.25pt](67.500,67.500)(66.250,68.750)(65.000,67.500)
\psline[linewidth=0.25pt](68.750,68.750)(67.500,70.000)(66.250,68.750)
\psline[linewidth=0.25pt](70.000,70.000)(68.750,71.250)(67.500,70.000)
\psline[linewidth=0.25pt](0.00000,2.5000)(-1.2500,3.7500)(-2.5000,2.5000)
\psline[linewidth=0.25pt](1.2500,3.7500)(0.00000,5.0000)(-1.2500,3.7500)
\psline[linewidth=0.25pt](2.5000,5.0000)(1.2500,6.2500)(0.00000,5.0000)
\psline[linewidth=0.25pt](3.7500,6.2500)(2.5000,7.5000)(1.2500,6.2500)
\psline[linewidth=0.25pt](5.0000,7.5000)(3.7500,8.7500)(2.5000,7.5000)
\psline[linewidth=0.25pt](6.2500,8.7500)(5.0000,10.000)(3.7500,8.7500)
\psline[linewidth=0.25pt](7.5000,10.000)(6.2500,11.250)(5.0000,10.000)
\psline[linewidth=0.25pt](8.7500,11.250)(7.5000,12.500)(6.2500,11.250)
\psline[linewidth=0.25pt](10.000,12.500)(8.7500,13.750)(7.5000,12.500)
\psline[linewidth=0.25pt](11.250,13.750)(10.000,15.000)(8.7500,13.750)
\psline[linewidth=0.25pt](12.500,15.000)(11.250,16.250)(10.000,15.000)
\psline[linewidth=0.25pt](13.750,16.250)(12.500,17.500)(11.250,16.250)
\psline[linewidth=0.25pt](15.000,17.500)(13.750,18.750)(12.500,17.500)
\psline[linewidth=0.25pt](16.250,18.750)(15.000,20.000)(13.750,18.750)
\psline[linewidth=0.25pt](17.500,20.000)(16.250,21.250)(15.000,20.000)
\psline[linewidth=0.25pt](18.750,21.250)(17.500,22.500)(16.250,21.250)
\psline[linewidth=0.25pt](20.000,22.500)(18.750,23.750)(17.500,22.500)
\psline[linewidth=0.25pt](21.250,23.750)(20.000,25.000)(18.750,23.750)
\psline[linewidth=0.25pt](22.500,25.000)(21.250,26.250)(20.000,25.000)
\psline[linewidth=0.25pt](23.750,26.250)(22.500,27.500)(21.250,26.250)
\psline[linewidth=0.25pt](25.000,27.500)(23.750,28.750)(22.500,27.500)
\psline[linewidth=0.25pt](26.250,28.750)(25.000,30.000)(23.750,28.750)
\psline[linewidth=0.25pt](27.500,30.000)(26.250,31.250)(25.000,30.000)
\psline[linewidth=0.25pt](28.750,31.250)(27.500,32.500)(26.250,31.250)
\psline[linewidth=0.25pt](30.000,32.500)(28.750,33.750)(27.500,32.500)
\psline[linewidth=0.25pt](31.250,33.750)(30.000,35.000)(28.750,33.750)
\psline[linewidth=0.25pt](32.500,35.000)(31.250,36.250)(30.000,35.000)
\psline[linewidth=0.25pt](33.750,36.250)(32.500,37.500)(31.250,36.250)
\psline[linewidth=0.25pt](35.000,37.500)(33.750,38.750)(32.500,37.500)
\psline[linewidth=0.25pt](36.250,38.750)(35.000,40.000)(33.750,38.750)
\psline[linewidth=0.25pt](37.500,40.000)(36.250,41.250)(35.000,40.000)
\psline[linewidth=0.25pt](38.750,41.250)(37.500,42.500)(36.250,41.250)
\psline[linewidth=0.25pt](40.000,42.500)(38.750,43.750)(37.500,42.500)
\psline[linewidth=0.25pt](41.250,43.750)(40.000,45.000)(38.750,43.750)
\psline[linewidth=0.25pt](42.500,45.000)(41.250,46.250)(40.000,45.000)
\psline[linewidth=0.25pt](43.750,46.250)(42.500,47.500)(41.250,46.250)
\psline[linewidth=0.25pt](45.000,47.500)(43.750,48.750)(42.500,47.500)
\psline[linewidth=0.25pt](46.250,48.750)(45.000,50.000)(43.750,48.750)
\psline[linewidth=0.25pt](47.500,50.000)(46.250,51.250)(45.000,50.000)
\psline[linewidth=0.25pt](48.750,51.250)(47.500,52.500)(46.250,51.250)
\psline[linewidth=0.25pt](50.000,52.500)(48.750,53.750)(47.500,52.500)
\psline[linewidth=0.25pt](51.250,53.750)(50.000,55.000)(48.750,53.750)
\psline[linewidth=0.25pt](52.500,55.000)(51.250,56.250)(50.000,55.000)
\psline[linewidth=0.25pt](53.750,56.250)(52.500,57.500)(51.250,56.250)
\psline[linewidth=0.25pt](55.000,57.500)(53.750,58.750)(52.500,57.500)
\psline[linewidth=0.25pt](56.250,58.750)(55.000,60.000)(53.750,58.750)
\psline[linewidth=0.25pt](57.500,60.000)(56.250,61.250)(55.000,60.000)
\psline[linewidth=0.25pt](58.750,61.250)(57.500,62.500)(56.250,61.250)
\psline[linewidth=0.25pt](60.000,62.500)(58.750,63.750)(57.500,62.500)
\psline[linewidth=0.25pt](61.250,63.750)(60.000,65.000)(58.750,63.750)
\psline[linewidth=0.25pt](62.500,65.000)(61.250,66.250)(60.000,65.000)
\psline[linewidth=0.25pt](63.750,66.250)(62.500,67.500)(61.250,66.250)
\psline[linewidth=0.25pt](-1.2500,3.7500)(-2.5000,5.0000)(-3.7500,3.7500)
\psline[linewidth=0.25pt](0.00000,5.0000)(-1.2500,6.2500)(-2.5000,5.0000)
\psline[linewidth=0.25pt](1.2500,6.2500)(0.00000,7.5000)(-1.2500,6.2500)
\psline[linewidth=0.25pt](2.5000,7.5000)(1.2500,8.7500)(0.00000,7.5000)
\psline[linewidth=0.25pt](3.7500,8.7500)(2.5000,10.000)(1.2500,8.7500)
\psline[linewidth=0.25pt](5.0000,10.000)(3.7500,11.250)(2.5000,10.000)
\psline[linewidth=0.25pt](6.2500,11.250)(5.0000,12.500)(3.7500,11.250)
\psline[linewidth=0.25pt](7.5000,12.500)(6.2500,13.750)(5.0000,12.500)
\psline[linewidth=0.25pt](8.7500,13.750)(7.5000,15.000)(6.2500,13.750)
\psline[linewidth=0.25pt](10.000,15.000)(8.7500,16.250)(7.5000,15.000)
\psline[linewidth=0.25pt](11.250,16.250)(10.000,17.500)(8.7500,16.250)
\psline[linewidth=0.25pt](12.500,17.500)(11.250,18.750)(10.000,17.500)
\psline[linewidth=0.25pt](13.750,18.750)(12.500,20.000)(11.250,18.750)
\psline[linewidth=0.25pt](15.000,20.000)(13.750,21.250)(12.500,20.000)
\psline[linewidth=0.25pt](16.250,21.250)(15.000,22.500)(13.750,21.250)
\psline[linewidth=0.25pt](17.500,22.500)(16.250,23.750)(15.000,22.500)
\psline[linewidth=0.25pt](18.750,23.750)(17.500,25.000)(16.250,23.750)
\psline[linewidth=0.25pt](20.000,25.000)(18.750,26.250)(17.500,25.000)
\psline[linewidth=0.25pt](21.250,26.250)(20.000,27.500)(18.750,26.250)
\psline[linewidth=0.25pt](22.500,27.500)(21.250,28.750)(20.000,27.500)
\psline[linewidth=0.25pt](23.750,28.750)(22.500,30.000)(21.250,28.750)
\psline[linewidth=0.25pt](25.000,30.000)(23.750,31.250)(22.500,30.000)
\psline[linewidth=0.25pt](26.250,31.250)(25.000,32.500)(23.750,31.250)
\psline[linewidth=0.25pt](27.500,32.500)(26.250,33.750)(25.000,32.500)
\psline[linewidth=0.25pt](28.750,33.750)(27.500,35.000)(26.250,33.750)
\psline[linewidth=0.25pt](30.000,35.000)(28.750,36.250)(27.500,35.000)
\psline[linewidth=0.25pt](31.250,36.250)(30.000,37.500)(28.750,36.250)
\psline[linewidth=0.25pt](32.500,37.500)(31.250,38.750)(30.000,37.500)
\psline[linewidth=0.25pt](33.750,38.750)(32.500,40.000)(31.250,38.750)
\psline[linewidth=0.25pt](35.000,40.000)(33.750,41.250)(32.500,40.000)
\psline[linewidth=0.25pt](36.250,41.250)(35.000,42.500)(33.750,41.250)
\psline[linewidth=0.25pt](37.500,42.500)(36.250,43.750)(35.000,42.500)
\psline[linewidth=0.25pt](38.750,43.750)(37.500,45.000)(36.250,43.750)
\psline[linewidth=0.25pt](40.000,45.000)(38.750,46.250)(37.500,45.000)
\psline[linewidth=0.25pt](41.250,46.250)(40.000,47.500)(38.750,46.250)
\psline[linewidth=0.25pt](42.500,47.500)(41.250,48.750)(40.000,47.500)
\psline[linewidth=0.25pt](43.750,48.750)(42.500,50.000)(41.250,48.750)
\psline[linewidth=0.25pt](45.000,50.000)(43.750,51.250)(42.500,50.000)
\psline[linewidth=0.25pt](46.250,51.250)(45.000,52.500)(43.750,51.250)
\psline[linewidth=0.25pt](47.500,52.500)(46.250,53.750)(45.000,52.500)
\psline[linewidth=0.25pt](48.750,53.750)(47.500,55.000)(46.250,53.750)
\psline[linewidth=0.25pt](50.000,55.000)(48.750,56.250)(47.500,55.000)
\psline[linewidth=0.25pt](51.250,56.250)(50.000,57.500)(48.750,56.250)
\psline[linewidth=0.25pt](52.500,57.500)(51.250,58.750)(50.000,57.500)
\psline[linewidth=0.25pt](53.750,58.750)(52.500,60.000)(51.250,58.750)
\psline[linewidth=0.25pt](55.000,60.000)(53.750,61.250)(52.500,60.000)
\psline[linewidth=0.25pt](-2.5000,5.0000)(-3.7500,6.2500)(-5.0000,5.0000)
\psline[linewidth=0.25pt](-1.2500,6.2500)(-2.5000,7.5000)(-3.7500,6.2500)
\psline[linewidth=0.25pt](0.00000,7.5000)(-1.2500,8.7500)(-2.5000,7.5000)
\psline[linewidth=0.25pt](1.2500,8.7500)(0.00000,10.000)(-1.2500,8.7500)
\psline[linewidth=0.25pt](2.5000,10.000)(1.2500,11.250)(0.00000,10.000)
\psline[linewidth=0.25pt](3.7500,11.250)(2.5000,12.500)(1.2500,11.250)
\psline[linewidth=0.25pt](5.0000,12.500)(3.7500,13.750)(2.5000,12.500)
\psline[linewidth=0.25pt](6.2500,13.750)(5.0000,15.000)(3.7500,13.750)
\psline[linewidth=0.25pt](7.5000,15.000)(6.2500,16.250)(5.0000,15.000)
\psline[linewidth=0.25pt](8.7500,16.250)(7.5000,17.500)(6.2500,16.250)
\psline[linewidth=0.25pt](10.000,17.500)(8.7500,18.750)(7.5000,17.500)
\psline[linewidth=0.25pt](11.250,18.750)(10.000,20.000)(8.7500,18.750)
\psline[linewidth=0.25pt](12.500,20.000)(11.250,21.250)(10.000,20.000)
\psline[linewidth=0.25pt](13.750,21.250)(12.500,22.500)(11.250,21.250)
\psline[linewidth=0.25pt](15.000,22.500)(13.750,23.750)(12.500,22.500)
\psline[linewidth=0.25pt](16.250,23.750)(15.000,25.000)(13.750,23.750)
\psline[linewidth=0.25pt](17.500,25.000)(16.250,26.250)(15.000,25.000)
\psline[linewidth=0.25pt](18.750,26.250)(17.500,27.500)(16.250,26.250)
\psline[linewidth=0.25pt](20.000,27.500)(18.750,28.750)(17.500,27.500)
\psline[linewidth=0.25pt](21.250,28.750)(20.000,30.000)(18.750,28.750)
\psline[linewidth=0.25pt](22.500,30.000)(21.250,31.250)(20.000,30.000)
\psline[linewidth=0.25pt](23.750,31.250)(22.500,32.500)(21.250,31.250)
\psline[linewidth=0.25pt](25.000,32.500)(23.750,33.750)(22.500,32.500)
\psline[linewidth=0.25pt](26.250,33.750)(25.000,35.000)(23.750,33.750)
\psline[linewidth=0.25pt](27.500,35.000)(26.250,36.250)(25.000,35.000)
\psline[linewidth=0.25pt](28.750,36.250)(27.500,37.500)(26.250,36.250)
\psline[linewidth=0.25pt](30.000,37.500)(28.750,38.750)(27.500,37.500)
\psline[linewidth=0.25pt](31.250,38.750)(30.000,40.000)(28.750,38.750)
\psline[linewidth=0.25pt](32.500,40.000)(31.250,41.250)(30.000,40.000)
\psline[linewidth=0.25pt](33.750,41.250)(32.500,42.500)(31.250,41.250)
\psline[linewidth=0.25pt](35.000,42.500)(33.750,43.750)(32.500,42.500)
\psline[linewidth=0.25pt](36.250,43.750)(35.000,45.000)(33.750,43.750)
\psline[linewidth=0.25pt](37.500,45.000)(36.250,46.250)(35.000,45.000)
\psline[linewidth=0.25pt](38.750,46.250)(37.500,47.500)(36.250,46.250)
\psline[linewidth=0.25pt](40.000,47.500)(38.750,48.750)(37.500,47.500)
\psline[linewidth=0.25pt](41.250,48.750)(40.000,50.000)(38.750,48.750)
\psline[linewidth=0.25pt](42.500,50.000)(41.250,51.250)(40.000,50.000)
\psline[linewidth=0.25pt](43.750,51.250)(42.500,52.500)(41.250,51.250)
\psline[linewidth=0.25pt](45.000,52.500)(43.750,53.750)(42.500,52.500)
\psline[linewidth=0.25pt](46.250,53.750)(45.000,55.000)(43.750,53.750)
\psline[linewidth=0.25pt](47.500,55.000)(46.250,56.250)(45.000,55.000)
\psline[linewidth=0.25pt](48.750,56.250)(47.500,57.500)(46.250,56.250)
\psline[linewidth=0.25pt](50.000,57.500)(48.750,58.750)(47.500,57.500)
\psline[linewidth=0.25pt](-3.7500,6.2500)(-5.0000,7.5000)(-6.2500,6.2500)
\psline[linewidth=0.25pt](-2.5000,7.5000)(-3.7500,8.7500)(-5.0000,7.5000)
\psline[linewidth=0.25pt](-1.2500,8.7500)(-2.5000,10.000)(-3.7500,8.7500)
\psline[linewidth=0.25pt](0.00000,10.000)(-1.2500,11.250)(-2.5000,10.000)
\psline[linewidth=0.25pt](1.2500,11.250)(0.00000,12.500)(-1.2500,11.250)
\psline[linewidth=0.25pt](2.5000,12.500)(1.2500,13.750)(0.00000,12.500)
\psline[linewidth=0.25pt](3.7500,13.750)(2.5000,15.000)(1.2500,13.750)
\psline[linewidth=0.25pt](5.0000,15.000)(3.7500,16.250)(2.5000,15.000)
\psline[linewidth=0.25pt](6.2500,16.250)(5.0000,17.500)(3.7500,16.250)
\psline[linewidth=0.25pt](7.5000,17.500)(6.2500,18.750)(5.0000,17.500)
\psline[linewidth=0.25pt](8.7500,18.750)(7.5000,20.000)(6.2500,18.750)
\psline[linewidth=0.25pt](10.000,20.000)(8.7500,21.250)(7.5000,20.000)
\psline[linewidth=0.25pt](11.250,21.250)(10.000,22.500)(8.7500,21.250)
\psline[linewidth=0.25pt](12.500,22.500)(11.250,23.750)(10.000,22.500)
\psline[linewidth=0.25pt](13.750,23.750)(12.500,25.000)(11.250,23.750)
\psline[linewidth=0.25pt](15.000,25.000)(13.750,26.250)(12.500,25.000)
\psline[linewidth=0.25pt](16.250,26.250)(15.000,27.500)(13.750,26.250)
\psline[linewidth=0.25pt](17.500,27.500)(16.250,28.750)(15.000,27.500)
\psline[linewidth=0.25pt](18.750,28.750)(17.500,30.000)(16.250,28.750)
\psline[linewidth=0.25pt](20.000,30.000)(18.750,31.250)(17.500,30.000)
\psline[linewidth=0.25pt](21.250,31.250)(20.000,32.500)(18.750,31.250)
\psline[linewidth=0.25pt](22.500,32.500)(21.250,33.750)(20.000,32.500)
\psline[linewidth=0.25pt](23.750,33.750)(22.500,35.000)(21.250,33.750)
\psline[linewidth=0.25pt](25.000,35.000)(23.750,36.250)(22.500,35.000)
\psline[linewidth=0.25pt](26.250,36.250)(25.000,37.500)(23.750,36.250)
\psline[linewidth=0.25pt](27.500,37.500)(26.250,38.750)(25.000,37.500)
\psline[linewidth=0.25pt](28.750,38.750)(27.500,40.000)(26.250,38.750)
\psline[linewidth=0.25pt](30.000,40.000)(28.750,41.250)(27.500,40.000)
\psline[linewidth=0.25pt](31.250,41.250)(30.000,42.500)(28.750,41.250)
\psline[linewidth=0.25pt](32.500,42.500)(31.250,43.750)(30.000,42.500)
\psline[linewidth=0.25pt](33.750,43.750)(32.500,45.000)(31.250,43.750)
\psline[linewidth=0.25pt](35.000,45.000)(33.750,46.250)(32.500,45.000)
\psline[linewidth=0.25pt](36.250,46.250)(35.000,47.500)(33.750,46.250)
\psline[linewidth=0.25pt](37.500,47.500)(36.250,48.750)(35.000,47.500)
\psline[linewidth=0.25pt](38.750,48.750)(37.500,50.000)(36.250,48.750)
\psline[linewidth=0.25pt](40.000,50.000)(38.750,51.250)(37.500,50.000)
\psline[linewidth=0.25pt](41.250,51.250)(40.000,52.500)(38.750,51.250)
\psline[linewidth=0.25pt](42.500,52.500)(41.250,53.750)(40.000,52.500)
\psline[linewidth=0.25pt](43.750,53.750)(42.500,55.000)(41.250,53.750)
\psline[linewidth=0.25pt](-5.0000,7.5000)(-6.2500,8.7500)(-7.5000,7.5000)
\psline[linewidth=0.25pt](-3.7500,8.7500)(-5.0000,10.000)(-6.2500,8.7500)
\psline[linewidth=0.25pt](-2.5000,10.000)(-3.7500,11.250)(-5.0000,10.000)
\psline[linewidth=0.25pt](-1.2500,11.250)(-2.5000,12.500)(-3.7500,11.250)
\psline[linewidth=0.25pt](0.00000,12.500)(-1.2500,13.750)(-2.5000,12.500)
\psline[linewidth=0.25pt](1.2500,13.750)(0.00000,15.000)(-1.2500,13.750)
\psline[linewidth=0.25pt](2.5000,15.000)(1.2500,16.250)(0.00000,15.000)
\psline[linewidth=0.25pt](3.7500,16.250)(2.5000,17.500)(1.2500,16.250)
\psline[linewidth=0.25pt](5.0000,17.500)(3.7500,18.750)(2.5000,17.500)
\psline[linewidth=0.25pt](6.2500,18.750)(5.0000,20.000)(3.7500,18.750)
\psline[linewidth=0.25pt](7.5000,20.000)(6.2500,21.250)(5.0000,20.000)
\psline[linewidth=0.25pt](8.7500,21.250)(7.5000,22.500)(6.2500,21.250)
\psline[linewidth=0.25pt](10.000,22.500)(8.7500,23.750)(7.5000,22.500)
\psline[linewidth=0.25pt](11.250,23.750)(10.000,25.000)(8.7500,23.750)
\psline[linewidth=0.25pt](12.500,25.000)(11.250,26.250)(10.000,25.000)
\psline[linewidth=0.25pt](13.750,26.250)(12.500,27.500)(11.250,26.250)
\psline[linewidth=0.25pt](15.000,27.500)(13.750,28.750)(12.500,27.500)
\psline[linewidth=0.25pt](16.250,28.750)(15.000,30.000)(13.750,28.750)
\psline[linewidth=0.25pt](17.500,30.000)(16.250,31.250)(15.000,30.000)
\psline[linewidth=0.25pt](18.750,31.250)(17.500,32.500)(16.250,31.250)
\psline[linewidth=0.25pt](20.000,32.500)(18.750,33.750)(17.500,32.500)
\psline[linewidth=0.25pt](21.250,33.750)(20.000,35.000)(18.750,33.750)
\psline[linewidth=0.25pt](22.500,35.000)(21.250,36.250)(20.000,35.000)
\psline[linewidth=0.25pt](23.750,36.250)(22.500,37.500)(21.250,36.250)
\psline[linewidth=0.25pt](25.000,37.500)(23.750,38.750)(22.500,37.500)
\psline[linewidth=0.25pt](26.250,38.750)(25.000,40.000)(23.750,38.750)
\psline[linewidth=0.25pt](27.500,40.000)(26.250,41.250)(25.000,40.000)
\psline[linewidth=0.25pt](28.750,41.250)(27.500,42.500)(26.250,41.250)
\psline[linewidth=0.25pt](30.000,42.500)(28.750,43.750)(27.500,42.500)
\psline[linewidth=0.25pt](31.250,43.750)(30.000,45.000)(28.750,43.750)
\psline[linewidth=0.25pt](32.500,45.000)(31.250,46.250)(30.000,45.000)
\psline[linewidth=0.25pt](33.750,46.250)(32.500,47.500)(31.250,46.250)
\psline[linewidth=0.25pt](35.000,47.500)(33.750,48.750)(32.500,47.500)
\psline[linewidth=0.25pt](36.250,48.750)(35.000,50.000)(33.750,48.750)
\psline[linewidth=0.25pt](37.500,50.000)(36.250,51.250)(35.000,50.000)
\psline[linewidth=0.25pt](38.750,51.250)(37.500,52.500)(36.250,51.250)
\psline[linewidth=0.25pt](-6.2500,8.7500)(-7.5000,10.000)(-8.7500,8.7500)
\psline[linewidth=0.25pt](-5.0000,10.000)(-6.2500,11.250)(-7.5000,10.000)
\psline[linewidth=0.25pt](-3.7500,11.250)(-5.0000,12.500)(-6.2500,11.250)
\psline[linewidth=0.25pt](-2.5000,12.500)(-3.7500,13.750)(-5.0000,12.500)
\psline[linewidth=0.25pt](-1.2500,13.750)(-2.5000,15.000)(-3.7500,13.750)
\psline[linewidth=0.25pt](0.00000,15.000)(-1.2500,16.250)(-2.5000,15.000)
\psline[linewidth=0.25pt](1.2500,16.250)(0.00000,17.500)(-1.2500,16.250)
\psline[linewidth=0.25pt](2.5000,17.500)(1.2500,18.750)(0.00000,17.500)
\psline[linewidth=0.25pt](3.7500,18.750)(2.5000,20.000)(1.2500,18.750)
\psline[linewidth=0.25pt](5.0000,20.000)(3.7500,21.250)(2.5000,20.000)
\psline[linewidth=0.25pt](6.2500,21.250)(5.0000,22.500)(3.7500,21.250)
\psline[linewidth=0.25pt](7.5000,22.500)(6.2500,23.750)(5.0000,22.500)
\psline[linewidth=0.25pt](8.7500,23.750)(7.5000,25.000)(6.2500,23.750)
\psline[linewidth=0.25pt](10.000,25.000)(8.7500,26.250)(7.5000,25.000)
\psline[linewidth=0.25pt](11.250,26.250)(10.000,27.500)(8.7500,26.250)
\psline[linewidth=0.25pt](12.500,27.500)(11.250,28.750)(10.000,27.500)
\psline[linewidth=0.25pt](13.750,28.750)(12.500,30.000)(11.250,28.750)
\psline[linewidth=0.25pt](15.000,30.000)(13.750,31.250)(12.500,30.000)
\psline[linewidth=0.25pt](16.250,31.250)(15.000,32.500)(13.750,31.250)
\psline[linewidth=0.25pt](17.500,32.500)(16.250,33.750)(15.000,32.500)
\psline[linewidth=0.25pt](18.750,33.750)(17.500,35.000)(16.250,33.750)
\psline[linewidth=0.25pt](20.000,35.000)(18.750,36.250)(17.500,35.000)
\psline[linewidth=0.25pt](21.250,36.250)(20.000,37.500)(18.750,36.250)
\psline[linewidth=0.25pt](22.500,37.500)(21.250,38.750)(20.000,37.500)
\psline[linewidth=0.25pt](23.750,38.750)(22.500,40.000)(21.250,38.750)
\psline[linewidth=0.25pt](25.000,40.000)(23.750,41.250)(22.500,40.000)
\psline[linewidth=0.25pt](26.250,41.250)(25.000,42.500)(23.750,41.250)
\psline[linewidth=0.25pt](27.500,42.500)(26.250,43.750)(25.000,42.500)
\psline[linewidth=0.25pt](28.750,43.750)(27.500,45.000)(26.250,43.750)
\psline[linewidth=0.25pt](30.000,45.000)(28.750,46.250)(27.500,45.000)
\psline[linewidth=0.25pt](31.250,46.250)(30.000,47.500)(28.750,46.250)
\psline[linewidth=0.25pt](32.500,47.500)(31.250,48.750)(30.000,47.500)
\psline[linewidth=0.25pt](33.750,48.750)(32.500,50.000)(31.250,48.750)
\psline[linewidth=0.25pt](-7.5000,10.000)(-8.7500,11.250)(-10.000,10.000)
\psline[linewidth=0.25pt](-6.2500,11.250)(-7.5000,12.500)(-8.7500,11.250)
\psline[linewidth=0.25pt](-5.0000,12.500)(-6.2500,13.750)(-7.5000,12.500)
\psline[linewidth=0.25pt](-3.7500,13.750)(-5.0000,15.000)(-6.2500,13.750)
\psline[linewidth=0.25pt](-2.5000,15.000)(-3.7500,16.250)(-5.0000,15.000)
\psline[linewidth=0.25pt](-1.2500,16.250)(-2.5000,17.500)(-3.7500,16.250)
\psline[linewidth=0.25pt](0.00000,17.500)(-1.2500,18.750)(-2.5000,17.500)
\psline[linewidth=0.25pt](1.2500,18.750)(0.00000,20.000)(-1.2500,18.750)
\psline[linewidth=0.25pt](2.5000,20.000)(1.2500,21.250)(0.00000,20.000)
\psline[linewidth=0.25pt](3.7500,21.250)(2.5000,22.500)(1.2500,21.250)
\psline[linewidth=0.25pt](5.0000,22.500)(3.7500,23.750)(2.5000,22.500)
\psline[linewidth=0.25pt](6.2500,23.750)(5.0000,25.000)(3.7500,23.750)
\psline[linewidth=0.25pt](7.5000,25.000)(6.2500,26.250)(5.0000,25.000)
\psline[linewidth=0.25pt](8.7500,26.250)(7.5000,27.500)(6.2500,26.250)
\psline[linewidth=0.25pt](10.000,27.500)(8.7500,28.750)(7.5000,27.500)
\psline[linewidth=0.25pt](11.250,28.750)(10.000,30.000)(8.7500,28.750)
\psline[linewidth=0.25pt](12.500,30.000)(11.250,31.250)(10.000,30.000)
\psline[linewidth=0.25pt](13.750,31.250)(12.500,32.500)(11.250,31.250)
\psline[linewidth=0.25pt](15.000,32.500)(13.750,33.750)(12.500,32.500)
\psline[linewidth=0.25pt](16.250,33.750)(15.000,35.000)(13.750,33.750)
\psline[linewidth=0.25pt](17.500,35.000)(16.250,36.250)(15.000,35.000)
\psline[linewidth=0.25pt](18.750,36.250)(17.500,37.500)(16.250,36.250)
\psline[linewidth=0.25pt](20.000,37.500)(18.750,38.750)(17.500,37.500)
\psline[linewidth=0.25pt](21.250,38.750)(20.000,40.000)(18.750,38.750)
\psline[linewidth=0.25pt](22.500,40.000)(21.250,41.250)(20.000,40.000)
\psline[linewidth=0.25pt](23.750,41.250)(22.500,42.500)(21.250,41.250)
\psline[linewidth=0.25pt](25.000,42.500)(23.750,43.750)(22.500,42.500)
\psline[linewidth=0.25pt](26.250,43.750)(25.000,45.000)(23.750,43.750)
\psline[linewidth=0.25pt](27.500,45.000)(26.250,46.250)(25.000,45.000)
\psline[linewidth=0.25pt](-8.7500,11.250)(-10.000,12.500)(-11.250,11.250)
\psline[linewidth=0.25pt](-7.5000,12.500)(-8.7500,13.750)(-10.000,12.500)
\psline[linewidth=0.25pt](-6.2500,13.750)(-7.5000,15.000)(-8.7500,13.750)
\psline[linewidth=0.25pt](-5.0000,15.000)(-6.2500,16.250)(-7.5000,15.000)
\psline[linewidth=0.25pt](-3.7500,16.250)(-5.0000,17.500)(-6.2500,16.250)
\psline[linewidth=0.25pt](-2.5000,17.500)(-3.7500,18.750)(-5.0000,17.500)
\psline[linewidth=0.25pt](-1.2500,18.750)(-2.5000,20.000)(-3.7500,18.750)
\psline[linewidth=0.25pt](0.00000,20.000)(-1.2500,21.250)(-2.5000,20.000)
\psline[linewidth=0.25pt](1.2500,21.250)(0.00000,22.500)(-1.2500,21.250)
\psline[linewidth=0.25pt](2.5000,22.500)(1.2500,23.750)(0.00000,22.500)
\psline[linewidth=0.25pt](3.7500,23.750)(2.5000,25.000)(1.2500,23.750)
\psline[linewidth=0.25pt](5.0000,25.000)(3.7500,26.250)(2.5000,25.000)
\psline[linewidth=0.25pt](6.2500,26.250)(5.0000,27.500)(3.7500,26.250)
\psline[linewidth=0.25pt](7.5000,27.500)(6.2500,28.750)(5.0000,27.500)
\psline[linewidth=0.25pt](8.7500,28.750)(7.5000,30.000)(6.2500,28.750)
\psline[linewidth=0.25pt](10.000,30.000)(8.7500,31.250)(7.5000,30.000)
\psline[linewidth=0.25pt](11.250,31.250)(10.000,32.500)(8.7500,31.250)
\psline[linewidth=0.25pt](12.500,32.500)(11.250,33.750)(10.000,32.500)
\psline[linewidth=0.25pt](13.750,33.750)(12.500,35.000)(11.250,33.750)
\psline[linewidth=0.25pt](15.000,35.000)(13.750,36.250)(12.500,35.000)
\psline[linewidth=0.25pt](16.250,36.250)(15.000,37.500)(13.750,36.250)
\psline[linewidth=0.25pt](17.500,37.500)(16.250,38.750)(15.000,37.500)
\psline[linewidth=0.25pt](18.750,38.750)(17.500,40.000)(16.250,38.750)
\psline[linewidth=0.25pt](20.000,40.000)(18.750,41.250)(17.500,40.000)
\psline[linewidth=0.25pt](21.250,41.250)(20.000,42.500)(18.750,41.250)
\psline[linewidth=0.25pt](-10.000,12.500)(-11.250,13.750)(-12.500,12.500)
\psline[linewidth=0.25pt](-8.7500,13.750)(-10.000,15.000)(-11.250,13.750)
\psline[linewidth=0.25pt](-7.5000,15.000)(-8.7500,16.250)(-10.000,15.000)
\psline[linewidth=0.25pt](-6.2500,16.250)(-7.5000,17.500)(-8.7500,16.250)
\psline[linewidth=0.25pt](-5.0000,17.500)(-6.2500,18.750)(-7.5000,17.500)
\psline[linewidth=0.25pt](-3.7500,18.750)(-5.0000,20.000)(-6.2500,18.750)
\psline[linewidth=0.25pt](-2.5000,20.000)(-3.7500,21.250)(-5.0000,20.000)
\psline[linewidth=0.25pt](-1.2500,21.250)(-2.5000,22.500)(-3.7500,21.250)
\psline[linewidth=0.25pt](0.00000,22.500)(-1.2500,23.750)(-2.5000,22.500)
\psline[linewidth=0.25pt](1.2500,23.750)(0.00000,25.000)(-1.2500,23.750)
\psline[linewidth=0.25pt](2.5000,25.000)(1.2500,26.250)(0.00000,25.000)
\psline[linewidth=0.25pt](3.7500,26.250)(2.5000,27.500)(1.2500,26.250)
\psline[linewidth=0.25pt](5.0000,27.500)(3.7500,28.750)(2.5000,27.500)
\psline[linewidth=0.25pt](6.2500,28.750)(5.0000,30.000)(3.7500,28.750)
\psline[linewidth=0.25pt](7.5000,30.000)(6.2500,31.250)(5.0000,30.000)
\psline[linewidth=0.25pt](8.7500,31.250)(7.5000,32.500)(6.2500,31.250)
\psline[linewidth=0.25pt](10.000,32.500)(8.7500,33.750)(7.5000,32.500)
\psline[linewidth=0.25pt](11.250,33.750)(10.000,35.000)(8.7500,33.750)
\psline[linewidth=0.25pt](12.500,35.000)(11.250,36.250)(10.000,35.000)
\psline[linewidth=0.25pt](13.750,36.250)(12.500,37.500)(11.250,36.250)
\psline[linewidth=0.25pt](15.000,37.500)(13.750,38.750)(12.500,37.500)
\psline[linewidth=0.25pt](16.250,38.750)(15.000,40.000)(13.750,38.750)
\psline[linewidth=0.25pt](-11.250,13.750)(-12.500,15.000)(-13.750,13.750)
\psline[linewidth=0.25pt](-10.000,15.000)(-11.250,16.250)(-12.500,15.000)
\psline[linewidth=0.25pt](-8.7500,16.250)(-10.000,17.500)(-11.250,16.250)
\psline[linewidth=0.25pt](-7.5000,17.500)(-8.7500,18.750)(-10.000,17.500)
\psline[linewidth=0.25pt](-6.2500,18.750)(-7.5000,20.000)(-8.7500,18.750)
\psline[linewidth=0.25pt](-5.0000,20.000)(-6.2500,21.250)(-7.5000,20.000)
\psline[linewidth=0.25pt](-3.7500,21.250)(-5.0000,22.500)(-6.2500,21.250)
\psline[linewidth=0.25pt](-2.5000,22.500)(-3.7500,23.750)(-5.0000,22.500)
\psline[linewidth=0.25pt](-1.2500,23.750)(-2.5000,25.000)(-3.7500,23.750)
\psline[linewidth=0.25pt](0.00000,25.000)(-1.2500,26.250)(-2.5000,25.000)
\psline[linewidth=0.25pt](1.2500,26.250)(0.00000,27.500)(-1.2500,26.250)
\psline[linewidth=0.25pt](2.5000,27.500)(1.2500,28.750)(0.00000,27.500)
\psline[linewidth=0.25pt](3.7500,28.750)(2.5000,30.000)(1.2500,28.750)
\psline[linewidth=0.25pt](5.0000,30.000)(3.7500,31.250)(2.5000,30.000)
\psline[linewidth=0.25pt](6.2500,31.250)(5.0000,32.500)(3.7500,31.250)
\psline[linewidth=0.25pt](7.5000,32.500)(6.2500,33.750)(5.0000,32.500)
\psline[linewidth=0.25pt](8.7500,33.750)(7.5000,35.000)(6.2500,33.750)
\psline[linewidth=0.25pt](10.000,35.000)(8.7500,36.250)(7.5000,35.000)
\psline[linewidth=0.25pt](11.250,36.250)(10.000,37.500)(8.7500,36.250)
\psline[linewidth=0.25pt](12.500,37.500)(11.250,38.750)(10.000,37.500)
\psline[linewidth=0.25pt](13.750,38.750)(12.500,40.000)(11.250,38.750)
\psline[linewidth=0.25pt](-12.500,15.000)(-13.750,16.250)(-15.000,15.000)
\psline[linewidth=0.25pt](-11.250,16.250)(-12.500,17.500)(-13.750,16.250)
\psline[linewidth=0.25pt](-10.000,17.500)(-11.250,18.750)(-12.500,17.500)
\psline[linewidth=0.25pt](-8.7500,18.750)(-10.000,20.000)(-11.250,18.750)
\psline[linewidth=0.25pt](-7.5000,20.000)(-8.7500,21.250)(-10.000,20.000)
\psline[linewidth=0.25pt](-6.2500,21.250)(-7.5000,22.500)(-8.7500,21.250)
\psline[linewidth=0.25pt](-5.0000,22.500)(-6.2500,23.750)(-7.5000,22.500)
\psline[linewidth=0.25pt](-3.7500,23.750)(-5.0000,25.000)(-6.2500,23.750)
\psline[linewidth=0.25pt](-2.5000,25.000)(-3.7500,26.250)(-5.0000,25.000)
\psline[linewidth=0.25pt](-1.2500,26.250)(-2.5000,27.500)(-3.7500,26.250)
\psline[linewidth=0.25pt](0.00000,27.500)(-1.2500,28.750)(-2.5000,27.500)
\psline[linewidth=0.25pt](1.2500,28.750)(0.00000,30.000)(-1.2500,28.750)
\psline[linewidth=0.25pt](2.5000,30.000)(1.2500,31.250)(0.00000,30.000)
\psline[linewidth=0.25pt](3.7500,31.250)(2.5000,32.500)(1.2500,31.250)
\psline[linewidth=0.25pt](5.0000,32.500)(3.7500,33.750)(2.5000,32.500)
\psline[linewidth=0.25pt](6.2500,33.750)(5.0000,35.000)(3.7500,33.750)
\psline[linewidth=0.25pt](7.5000,35.000)(6.2500,36.250)(5.0000,35.000)
\psline[linewidth=0.25pt](8.7500,36.250)(7.5000,37.500)(6.2500,36.250)
\psline[linewidth=0.25pt](10.000,37.500)(8.7500,38.750)(7.5000,37.500)
\psline[linewidth=0.25pt](11.250,38.750)(10.000,40.000)(8.7500,38.750)
\psline[linewidth=0.25pt](-13.750,16.250)(-15.000,17.500)(-16.250,16.250)
\psline[linewidth=0.25pt](-12.500,17.500)(-13.750,18.750)(-15.000,17.500)
\psline[linewidth=0.25pt](-11.250,18.750)(-12.500,20.000)(-13.750,18.750)
\psline[linewidth=0.25pt](-10.000,20.000)(-11.250,21.250)(-12.500,20.000)
\psline[linewidth=0.25pt](-8.7500,21.250)(-10.000,22.500)(-11.250,21.250)
\psline[linewidth=0.25pt](-7.5000,22.500)(-8.7500,23.750)(-10.000,22.500)
\psline[linewidth=0.25pt](-6.2500,23.750)(-7.5000,25.000)(-8.7500,23.750)
\psline[linewidth=0.25pt](-5.0000,25.000)(-6.2500,26.250)(-7.5000,25.000)
\psline[linewidth=0.25pt](-3.7500,26.250)(-5.0000,27.500)(-6.2500,26.250)
\psline[linewidth=0.25pt](-2.5000,27.500)(-3.7500,28.750)(-5.0000,27.500)
\psline[linewidth=0.25pt](-1.2500,28.750)(-2.5000,30.000)(-3.7500,28.750)
\psline[linewidth=0.25pt](0.00000,30.000)(-1.2500,31.250)(-2.5000,30.000)
\psline[linewidth=0.25pt](1.2500,31.250)(0.00000,32.500)(-1.2500,31.250)
\psline[linewidth=0.25pt](2.5000,32.500)(1.2500,33.750)(0.00000,32.500)
\psline[linewidth=0.25pt](3.7500,33.750)(2.5000,35.000)(1.2500,33.750)
\psline[linewidth=0.25pt](5.0000,35.000)(3.7500,36.250)(2.5000,35.000)
\psline[linewidth=0.25pt](6.2500,36.250)(5.0000,37.500)(3.7500,36.250)
\psline[linewidth=0.25pt](-15.000,17.500)(-16.250,18.750)(-17.500,17.500)
\psline[linewidth=0.25pt](-13.750,18.750)(-15.000,20.000)(-16.250,18.750)
\psline[linewidth=0.25pt](-12.500,20.000)(-13.750,21.250)(-15.000,20.000)
\psline[linewidth=0.25pt](-11.250,21.250)(-12.500,22.500)(-13.750,21.250)
\psline[linewidth=0.25pt](-10.000,22.500)(-11.250,23.750)(-12.500,22.500)
\psline[linewidth=0.25pt](-8.7500,23.750)(-10.000,25.000)(-11.250,23.750)
\psline[linewidth=0.25pt](-7.5000,25.000)(-8.7500,26.250)(-10.000,25.000)
\psline[linewidth=0.25pt](-6.2500,26.250)(-7.5000,27.500)(-8.7500,26.250)
\psline[linewidth=0.25pt](-5.0000,27.500)(-6.2500,28.750)(-7.5000,27.500)
\psline[linewidth=0.25pt](-3.7500,28.750)(-5.0000,30.000)(-6.2500,28.750)
\psline[linewidth=0.25pt](-2.5000,30.000)(-3.7500,31.250)(-5.0000,30.000)
\psline[linewidth=0.25pt](-1.2500,31.250)(-2.5000,32.500)(-3.7500,31.250)
\psline[linewidth=0.25pt](0.00000,32.500)(-1.2500,33.750)(-2.5000,32.500)
\psline[linewidth=0.25pt](1.2500,33.750)(0.00000,35.000)(-1.2500,33.750)
\psline[linewidth=0.25pt](-16.250,18.750)(-17.500,20.000)(-18.750,18.750)
\psline[linewidth=0.25pt](-15.000,20.000)(-16.250,21.250)(-17.500,20.000)
\psline[linewidth=0.25pt](-13.750,21.250)(-15.000,22.500)(-16.250,21.250)
\psline[linewidth=0.25pt](-12.500,22.500)(-13.750,23.750)(-15.000,22.500)
\psline[linewidth=0.25pt](-11.250,23.750)(-12.500,25.000)(-13.750,23.750)
\psline[linewidth=0.25pt](-10.000,25.000)(-11.250,26.250)(-12.500,25.000)
\psline[linewidth=0.25pt](-8.7500,26.250)(-10.000,27.500)(-11.250,26.250)
\psline[linewidth=0.25pt](-7.5000,27.500)(-8.7500,28.750)(-10.000,27.500)
\psline[linewidth=0.25pt](-6.2500,28.750)(-7.5000,30.000)(-8.7500,28.750)
\psline[linewidth=0.25pt](-5.0000,30.000)(-6.2500,31.250)(-7.5000,30.000)
\psline[linewidth=0.25pt](-3.7500,31.250)(-5.0000,32.500)(-6.2500,31.250)
\psline[linewidth=0.25pt](-2.5000,32.500)(-3.7500,33.750)(-5.0000,32.500)
\psline[linewidth=0.25pt](-17.500,20.000)(-18.750,21.250)(-20.000,20.000)
\psline[linewidth=0.25pt](-16.250,21.250)(-17.500,22.500)(-18.750,21.250)
\psline[linewidth=0.25pt](-15.000,22.500)(-16.250,23.750)(-17.500,22.500)
\psline[linewidth=0.25pt](-13.750,23.750)(-15.000,25.000)(-16.250,23.750)
\psline[linewidth=0.25pt](-12.500,25.000)(-13.750,26.250)(-15.000,25.000)
\psline[linewidth=0.25pt](-11.250,26.250)(-12.500,27.500)(-13.750,26.250)
\psline[linewidth=0.25pt](-10.000,27.500)(-11.250,28.750)(-12.500,27.500)
\psline[linewidth=0.25pt](-8.7500,28.750)(-10.000,30.000)(-11.250,28.750)
\psline[linewidth=0.25pt](-18.750,21.250)(-20.000,22.500)(-21.250,21.250)
\psline[linewidth=0.25pt](-17.500,22.500)(-18.750,23.750)(-20.000,22.500)
\psline[linewidth=0.25pt](-16.250,23.750)(-17.500,25.000)(-18.750,23.750)
\psline[linewidth=0.25pt](-15.000,25.000)(-16.250,26.250)(-17.500,25.000)
\psline[linewidth=0.25pt](-13.750,26.250)(-15.000,27.500)(-16.250,26.250)
\psline[linewidth=0.25pt](-12.500,27.500)(-13.750,28.750)(-15.000,27.500)
\psline[linewidth=0.25pt](-11.250,28.750)(-12.500,30.000)(-13.750,28.750)
\psline[linewidth=0.25pt](-10.000,30.000)(-11.250,31.250)(-12.500,30.000)
\psline[linewidth=0.25pt](-20.000,22.500)(-21.250,23.750)(-22.500,22.500)
\psline[linewidth=0.25pt](-18.750,23.750)(-20.000,25.000)(-21.250,23.750)
\psline[linewidth=0.25pt](-17.500,25.000)(-18.750,26.250)(-20.000,25.000)
\psline[linewidth=0.25pt](-16.250,26.250)(-17.500,27.500)(-18.750,26.250)
\psline[linewidth=0.25pt](-15.000,27.500)(-16.250,28.750)(-17.500,27.500)
\psline[linewidth=0.25pt](-13.750,28.750)(-15.000,30.000)(-16.250,28.750)
\psline[linewidth=0.25pt](-21.250,23.750)(-22.500,25.000)(-23.750,23.750)
\psline[linewidth=0.25pt](-20.000,25.000)(-21.250,26.250)(-22.500,25.000)
\psline[linewidth=0.25pt](-18.750,26.250)(-20.000,27.500)(-21.250,26.250)
\psline[linewidth=0.25pt](-17.500,27.500)(-18.750,28.750)(-20.000,27.500)
\psline[linewidth=0.25pt](-22.500,25.000)(-23.750,26.250)(-25.000,25.000)
\psline[linewidth=0.25pt](-21.250,26.250)(-22.500,27.500)(-23.750,26.250)
\psline[linewidth=0.25pt](-20.000,27.500)(-21.250,28.750)(-22.500,27.500)
\psline[linewidth=0.25pt](-18.750,28.750)(-20.000,30.000)(-21.250,28.750)
\psline[linewidth=0.25pt](-23.750,26.250)(-25.000,27.500)(-26.250,26.250)
\psline[linewidth=0.25pt](-22.500,27.500)(-23.750,28.750)(-25.000,27.500)
\psline[linewidth=0.25pt](-21.250,28.750)(-22.500,30.000)(-23.750,28.750)
\psline[linewidth=0.25pt](-25.000,27.500)(-26.250,28.750)(-27.500,27.500)
\psline[linewidth=0.25pt](-23.750,28.750)(-25.000,30.000)(-26.250,28.750)
\endpspicture}{Diagramme de Young aléatoire tiré suivant la mesure de Schur-Weyl de paramètres $n=500$, $\alpha=1/2$ et $c=1$. \label{largeschurweyl}}{Diagramme de Young aléatoire tiré suivant une mesure de Schur-Weyl}
Supposons maintenant $\alpha = 1/2$ (voir la figure \ref{largeschurweyl}). Les cumulants joints $k_{n}(\sigma_{l_{1}},\ldots,\sigma_{l_{r}})$  avec $r\geq 2$ restent nuls, et pour $r=1$, les limites des espérances renormalisées des cycles s'écrivent :
$$\lim_{n \to \infty}SW_{n,\alpha,c}[\sigma_{l}]\,n^{\frac{l-1}{2}}=\lim_{n\to \infty}SW_{n,\alpha,c}[\varSigma_{l}]\,n^{-\frac{l+1}{2}}=c^{l-1}\,.$$
On peut remplacer les caractères centraux $\varSigma_{l}$ par leurs composantes de plus haut poids, c'est-à-dire les cumulants libres $R_{l+1}$ ; par conséquent, notant $\lambda^{*}=\lambda^{\frac{1}{n}}$ un diagramme de taille $n$ renormalisé en abscisse et en ordonnée par un facteur $1/\sqrt{n}$, on voit que 
$$\lim_{n \to \infty} SW_{n,\alpha,c}[R_{l+1}(\lambda^{*})]=c^{l-1}\,.$$
Alors, la propriété de factorisation asymptotique permet d'établir la convergence en probabilité (et en moments) :
$$\forall l \geq 2,\,\,R_{l}(\lambda^{*}) \rightarrow_{SW_{n,\alpha,c}} c^{l-2}\,.$$
Le calcul de la forme limite correspondant à ces cumulants libres est dû à P. Biane, \emph{cf.} \cite{Bia01a}. Un paramètre $c$ étant fixé, notons $\mu_{c}$ la mesure de transition de la forme limite, $G_{c}$ sa transformée de Cauchy et $R_{c}$ sa $R$-transformée. Alors, par définition des cumulants libres :
$$R_{c}(z)=\frac{1}{z}\left(1+\sum_{l =1}^{\infty}R_{l}(\mu_{c})\,z^{l}\right)=\frac{1}{z}+\frac{z}{1-cz}\,.$$
En évaluant cette relation en $z=G_{c}(y)$, on obtient l'équation 
$$ (1+cy)\,[G_{c}(y)]^{2}-(c+y)\,G_{c}(y)+1=0\,,$$
et comme $G_{c}(y)\simeq_{y \to \infty} 1/y$, ceci conduit à 
$$G_{c}(y)=\frac{c+y- \sqrt{(c-y)^{2}-4}}{2(1+cy)}=\frac{2}{y+c+\sqrt{(y-c)^{2}-4}}\,,$$
où $y$ est une variable complexe, et la racine carrée holomorphe est définie sur $\C\setminus \R_{-}$ et choisie de telle sorte que $\sqrt{1}=1$. Des techniques d'inversion de Stieltjes permettent finalement de retrouver les formes limites correspondantes $\Omega_{c}$. Ainsi, notons $\Omega_{c}$ le diagramme continu défini par :
\begin{align*}&\Omega_{0}(s)=\Omega(s)=\begin{cases}
            \frac{2}{\pi} \left(s \arcsin(\frac{s}{2})+\sqrt{4-s^2} \right) &\text{si }|s| \leq 2,\\
|s| &\text{sinon ;}
            \end{cases}\\
&\Omega_{c\in \,]0,1[}(s)=\begin{cases}
            \frac{2}{\pi} \left(s \arcsin(\frac{s+c}{2\sqrt{1+sc}})+\frac{1}{c}\arccos(\frac{2+sc-c^{2}}{2\sqrt{1+sc}})+\frac{\sqrt{4-(s-c)^2}}{2} \right) &\text{si }|s-c| \leq 2,\\
|s| &\text{sinon ;}
            \end{cases}\\
                        &\Omega_{1}(s)=\begin{cases}
            \frac{s+1}{2}+\frac{1}{\pi}\left((s-1)\arcsin(\frac{s-1}{2})+\sqrt{4-(s-1)^{2}}\right)&\text{si }|s-1|\leq 2,\\
            |s|&\text{sinon ;}
            \end{cases}
\end{align*}

$$\Omega_{c >1}(s)=\begin{cases}s+\frac{2}{c}&\text{si } s\in\,]\frac{-1}{c},c-2[\,,\\
            \frac{2}{\pi} \left(s \arcsin(\frac{s+c}{2\sqrt{1+sc}})+\frac{1}{c}\arccos(\frac{2+sc-c^{2}}{2\sqrt{1+sc}})+\frac{\sqrt{4-(s-c)^2}}{2} \right) &\text{si }|s-c| \leq 2,\\
|s| &\text{sinon.}
            \end{cases}
$$
Notons que $\Omega_{1}(s \in [-1,3])$ est bien la limite des $\Omega_{c}(s)$ avec $c<1$ ou $c>1$ ; de même, $\Omega_{0}(s \in [-2,2])$ est bien la limite des $\Omega_{c}(s)$ avec $0<c<1$. Les diagrammes continus qui correspondent aux formules ci-dessus sont représentés sur la figure \ref{omegac} --- on a repris exactement la figure 2 de \cite{Bia01a}.
\begin{proposition}[Formes limites des diagrammes sous les mesures de Schur-Weyl de paramètre $\alpha=1/2$, \cite{Bia01a}]\label{firstasymptoticschurweyl}
Pour tout paramètre $c>0$ et tout $\eps>0$, $SW_{n,1/2,c}[\|\lambda^{*}-\Omega_{c}\|_{\infty}\geq \eps]$ tend vers\footnote{En réalité, on a même convergence en probabilité au sens de la topologie ultra-forte sur $\cym$, car il n'est pas très difficile d'évaluer la distribution des longueurs des plus longs sous-mots croissants et décroissants d'un mot aléatoire de $\lle 1,N\rre^n$.} $0$. Le cas $c=0$ correspond aux paramètres $\alpha>1/2$, et on a alors le même résultat avec $\Omega_{0}$ pour forme limite.
\end{proposition}

\figcapt{\psset{unit=0.5mm}
\pspicture(-50,-73)(170,50)
\psline{->}(-50,0)(50,0)
\psline{->}(-50,-70)(50,-70)
\psline{->}(70,0)(170,0)
\psline{->}(70,-70)(170,-70)
\psline{->}(0,0)(0,50)
\psline{->}(0,-70)(0,-20)
\psline{->}(120,0)(120,50)
\psline{->}(120,-70)(120,-20)
\psline(-45,45)(0,0)(45,45)
\psline(-45,-25)(0,-70)(45,-25)
\psline(75,45)(120,0)(165,45)
\psline(75,-25)(120,-70)(180,-10)
\rput(-40,5){\textcolor{BurntOrange}{$c=0$}}
\rput(-40,-65){\textcolor{violet}{$c=1$}}
\rput(80,5){\textcolor{red}{$c=1/2$}}
\rput(80,-65){\textcolor{blue}{$c=2$}}
\parametricplot[linecolor=BurntOrange]{-2}{2}{t 15 mul
t 2 div arcsin t mul 3.14159 mul 180 div t t mul -1 mul 4 add sqrt add 30 3.14159 div mul 0.25 add}
\psline[linecolor=BurntOrange](-45,45.25)(-30,30.25)
\psline[linecolor=BurntOrange](45,45.25)(30,30.25)
\parametricplot[linecolor=red]{-1.5}{2.5}{t 15 mul 120 add
t 0.5 add t 2 div 1 add sqrt 2 mul div arcsin t mul 1.75 t 2 div add t 2 div 1 add sqrt 2 mul div arccos 2 mul add 3.14159 mul 180 div  t 0.5 -1 mul add t 0.5 -1 mul add mul -1 mul 4 add sqrt 2 div add 30 3.14159 div mul 0.25 add}
\psline[linecolor=red](97.5,22.75)(75,45.25)
\psline[linecolor=red](165,45.25)(157.5,37.75)
\parametricplot[linecolor=violet]{-1}{3}{t 15 mul
t 1 add sqrt 2 div arcsin t mul t 1 add sqrt 2 div arccos add 3.14159 mul 180 div  t -1 add t -1 add mul -1 mul 4 add sqrt 2 div add 30 3.14159 div mul -69.75 add}
\psline[linecolor=violet](-15,-54.75)(-45,-24.75)
\parametricplot[linecolor=blue]{0}{4}{t 15 mul 120 add
t 2 add t 2 mul 1 add sqrt 2 mul div arcsin t mul t 2 mul -2 add t 2 mul 1 add sqrt 2 mul div arccos 0.5 mul add 3.14159 mul 180 div  t t mul -1 mul 4 t mul add sqrt 2 div add 30 3.14159 div mul -69.75 add}
\psline[linecolor=blue](75,-24.75)(112.5,-62.25)(120,-54.75)
\endpspicture}{Formes limites $\Omega_{c}$ des diagrammes de Young tirés suivant des mesures de Schur-Weyl de paramètre $\alpha=1/2$.\label{omegac}}{Formes limites pour les mesures de Schur-Weyl de paramètre $\alpha=1/2$}

En particulier, on observe une transition de phase lorsque le paramètre $c$ s'approche de la valeur $1$ : par exemple, la tangente à la courbe limite en $(c-2)^{+}$ vaut $-1$ si $c<1$, $0$ si $c=1$ et $+1$ si $c>1$. Le comportement asymptotique local au voisinage de $s=-1$ et pour $c=1$ peut être décrit par des processus déterminantaux reliés au noyau d'Hermite discret, voir \cite{BO07}. Dans la suite de cette section, on s'intéresse aux déviations des diagrammes de Young par rapport à leur forme limite $\Omega_{c}$ ; ce problème était jusqu'alors resté non résolu, et nous avons découvert que le théorème central limite de Kerov \ref{secondasymptoticplancherel} se généralisait dans ce cadre (\emph{cf.} \cite{Mel10b}). Pour commencer, calculons la mesure de transition $\mu_{c}$ dont la transformée de Cauchy est la fonction $G_{c}$ précédemment décrite. D'après le chapitre \ref{plancherel}, si $c=0$, alors 
$$d\mu_{0}(s)=\mathbb{1}_{s \in [-2,2]}\,\frac{\sqrt{4-s^{2}}}{2\pi}\,ds$$
est la loi de Wigner du demi-cercle.

\begin{lemma}[Déformations de la loi de Wigner]
La mesure de transition $\mu_{c}$ a pour expression :
\begin{align*}d\mu_{c\leq 1}(s)&=\mathbb{1}_{s \in [c-2,c+2]}\,\frac{\sqrt{4-(s-c)^{2}}}{2\pi(1+cs)}\,ds\,\,;\\
d\mu_{c> 1}(s)&=\mathbb{1}_{s \in [c-2,c+2]}\,\frac{\sqrt{4-(s-c)^{2}}}{2\pi(1+cs)}\,ds+\left(1-\frac{1}{c^{2}}\right)\delta_{-\frac{1}{c}}(s)\,.\end{align*}
\end{lemma}

\figcapt{\psset{unit=0.5mm}
\pspicture(-50,-73)(170,53)
\rput(-40,5){\textcolor{BurntOrange}{$c=0$}}
\rput(-40,-65){\textcolor{violet}{$c=1$}}
\rput(80,5){\textcolor{red}{$c=1/2$}}
\rput(80,-65){\textcolor{blue}{$c=2$}}
\parametricplot*[linecolor=BurntOrange!50!white]{-2}{2}{t 15 mul
4 t t mul -1 mul add sqrt 15 mul}
\parametricplot[linecolor=BurntOrange]{-2}{2}{t 15 mul
4 t t mul -1 mul add sqrt 15 mul}
\parametricplot*[linecolor=red!50!white]{-1.5}{2.5}{t 15 mul 120 add
4 t -0.5 add t -0.5 add mul -1 mul add sqrt 15 mul 1 0.5 t mul add div}
\parametricplot[linecolor=red]{-1.5}{2.5}{t 15 mul 120 add
4 t -0.5 add t -0.5 add mul -1 mul add sqrt 15 mul 1 0.5 t mul add div}
\pscustom[linecolor=violet!50!white,fillstyle=solid,fillcolor=violet!50!white]{
\parametricplot{-0.80}{3}{t 15 mul
4 t -1 add t -1 add mul -1 mul add sqrt 15 mul 1 t  add div -70 add}
\psline(45,-70)(-15,-70)}
\psframe*[linecolor=violet!50!white](-15,-70)(-12,-5)
\parametricplot[linecolor=violet]{-0.80}{3}{t 15 mul
4 t -1 add t -1 add mul -1 mul add sqrt 15 mul 1 t  add div -70 add}
\psline[linecolor=violet](-15,-70)(-15,-5)
\parametricplot*[linecolor=blue!50!white]{0}{4}{t 15 mul 120 add
4 t -2 add t -2 add mul -1 mul add sqrt 15 mul 1 t 2 mul add div -70 add}
\parametricplot[linecolor=blue]{0}{4}{t 15 mul 120 add
4 t -2 add t -2 add mul -1 mul add sqrt 15 mul 1 t 2 mul add div -70 add}
\psline[linecolor=blue,linewidth=2pt](112.5,-70)(112.5,-5)
\psline{->}(-50,0)(50,0)
\psline{->}(-50,-70)(50,-70)
\psline{->}(70,0)(170,0)
\psline{->}(70,-70)(200,-70)
\psline{->}(0,0)(0,50)
\psline{->}(0,-70)(0,-20)
\psline{->}(120,0)(120,50)
\psline{->}(120,-70)(120,-20)
\endpspicture
}{D\'eformations de la mesure de Wigner correspondant aux mesures de Schur-Weyl de paramètre $\alpha=1/2$.}{Déformations de la mesure de Wigner}

\begin{proof}
Supposons dans un premier temps $c<1$ ; alors, la fonction $G_{c}(z)$ est holomorphe sur $\C\setminus \R_{-}$ et reste bornée au voisinage de la demi-droite $\R_{-}$ (en particulier, sa singularité en $z=-1/c$ est effa\c cable). La loi $\mu_{c}$ dont la transformée de Cauchy est $G_{c}$ peut dans ce cas être retrouvée à l'aide de la \textbf{formule d'inversion de Perron-Stieltjes} :
$$\frac{d\mu_{c}(s)}{ds}=\lim_{\eps \to 0^{+}} \frac{G_{c}(s-\I\eps)-G_{c}(s+\I\eps)}{2\I \pi}\,;$$
cette formule classique découle d'un calcul de résidu. Si $s$ est strictement compris entre $c-2$ et $c+2$, alors $(s-c)^{2}-4$ est un réel strictement négatif, donc 
$$\lim_{\eps \to 0^{\pm}}\sqrt{(s+\I\eps-c)^{2}-4}=\pm\I\sqrt{4-(s-c)^{2}}\,.$$
D'autre part, tous les autres termes dans l'expression de $G_{c}$ sont continus au voisinage $s$ ; par conséquent, 
$$\frac{d\mu_{c}(s \in \,]c-2,c+2[)}{ds}=\frac{(\I\sqrt{4-(s-c)^{2}}) - (-\I\sqrt{4-(s-c)^{2}})}{4\I\pi(1+cs)}=\frac{\sqrt{4-(s-c)^{2}}}{2\pi(1+cs)}\,.$$
Si $|s-c|\geq 2$, alors $(s-c)^{2}-4$ est un réel positif, donc la racine carrée reste continue au voisinage de ce point. Par conséquent,
$$\frac{d\mu_{c}(s \notin \,]c-2,c+2[)}{ds}=0\,.$$
L'expression de $\mu_{c< 1}(ds)$ est donc établie, et notons qu'un calcul permet effectivement de retrouver $$\int_{s=c-2}^{s=c+2}\frac{\sqrt{4-(s-c)^{2}}}{(1+cs)(z-s)}\,ds=G_{c}(z)\,\, ;$$ il suffit d'effectuer le changement de variable $s-c=2\sin\theta$, de développer en série entière les expressions ainsi obtenues et d'utiliser les formules de Wallis pour les intégrales de puissances de la fonction sinus. Lorsque $c\geq 1$, les mêmes formules sont valables pour tout point $s\neq -1/c$, mais il faut rajouter un Dirac en $-1/c$ pour prendre en compte le pôle de la fonction $G_{c}(z)$ en $z=-1/c$. Comme le résidu de $G_{c}(z)$ en $-1/c$ est 
$$\frac{-1/c+c-\sqrt{(c+1/c)^{2}-4}}{2c}=\frac{-(c-1/c)-\sqrt{(c-1/c)^{2}}}{2c}=-\left(1-\frac{1}{c^{2}}\right)\,,$$
on obtient bien la formule\footnote{Notons que ces mesures déformées sont en réalité bien connues : il s'agit des lois images par des homothéties des distributions de Mar\v{c}enko-Pastur reliées aux valeurs propres de matrices aléatoires de Wishart, \emph{cf.} \cite{Bia01a}. On renvoie d'autre part à \cite{LP09} pour l'analogue dans le contexte des matrices aléatoires et des distributions de Mar\v{c}enko-Pastur de nos résultats de déviation gaussienne.} de l'énoncé. 
\end{proof}
\bigskip
\bigskip

Pour la mesure de Wigner $\mu_{0}$, les moments pairs sont donnés par les nombres de Catalan : $\widetilde{h}_{2n}(\Omega_{0})=\mu_{0}(s^{2n})=C_{n}$. Déterminons de même les moments des lois déformées $\mu_{c}$. 
\begin{lemma}[Deux identités hypergéométriques]\label{whatever}
Pour tous entiers positifs $m$ et $k$,
\begin{align*}&\sum_{l=0}^{k}\sum_{u=0}^{m} \frac{(2k+2+u-2l)\,\,m+2l+2-u!\,2k-2l+u!}{(m-u+l+2)(m-u+l+1)\,m-u!\,u!\,l!\,l+1!\,k-l!\,k-l+1!}\\
&= \frac{m+2k+4}{m+k+2}\,\times\,\frac{m+2k+2!}{m!\, k!\,k+2!}\,.\end{align*}
D'autre part, pour tous entiers positifs $Z$ et $\alpha<\beta$, 
$$\sum_{z=0}^{Z} \frac{\alpha+z!}{\beta+z!}\binom{Z}{z} (-1)^{z}= \frac{\alpha!\,\beta-\alpha+Z-1!}{\beta-\alpha-1!\,\beta+Z!}\,.$$
\end{lemma}
\begin{proof}
Ceci rentre dans le cadre des identités hypergéométriques (multivariées) que l'on peut démontrer à l'aide de la \textbf{théorie de Wilf-Zeilberger}, voir \cite{PWZ97,WZ92a,WZ92b}. Ainsi, dans chaque cas, il existe un algorithme qui exhibe des relations de récurrence vérifiées par les deux membres de l'équation, et d'autre part, il est aisé de vérifier que les identités sont vraies pour suffisamment de valeurs des paramètres ($m$ et $k$ pour la première identité, et $Z$ pour la seconde) ; partant, les identités sont effectivement vraies pour toutes valeurs des paramètres. Nous laissons les détails du calcul des relations de récurrence au lecteur, ou à tout logiciel de calcul formel (par exemple \texttt{Mathematica}).
\end{proof}\bigskip

\begin{lemma}[Moments des lois de Wigner déformées]\label{momentwignerdeform}
Le $n$-ième moment de la loi $\mu_{c}$ est égal à :
$$\widetilde{h}_{n}(\Omega_{c})=\mu_{c}(s^{n})=\sum_{k=1}^{\lfloor \frac{n}{2}\rfloor}\frac{n^{\downarrow 2k}}{(n-k+1)(n-k)\,k!\,k-1!}\,c^{n-2k}\,.$$
Le $n$-ième moment entrelacé correspondant est égal à :
$$\tilp_{n}(\Omega_{c})=\sum_{k=1}^{\lfloor \frac{n}{2}\rfloor}\frac{n^{\downarrow 2k}}{(n-k)\,k!\,k-1!}\,c^{n-2k}\,.$$
\end{lemma}
\begin{proof}
On a vu dans le chapitre \ref{tool} que les moments des diagrammes continus étaient reliés aux fonctions génératrices des diagrammes par la formule :
$$z^{-1}\,G_{\omega}(z^{-1})=1+\sum_{k=1}^{\infty}\widetilde{h}_{k}(\omega)\,z^{k}=\exp\left(\sum_{k=1}^{\infty} \frac{\tilp_{k}(\omega)}{k}\,z^{k}\right)\,.$$
Il suffit donc \emph{a priori} d'effectuer deux développements en série ; c'est ce que nous faisons dans la suite, mais les calculs sont assez apocalyptiques. Notons $H_{c}(z)=z^{-1}\,G_{\omega}(z^{-1})=\frac{1}{2(z+c)}\big(\frac{1}{z}+c-\sqrt{(\frac{1}{z}-c)^{2}-4}\big)=\frac{2}{1+cz+\sqrt{(1-cz)^{2}-4z^{2}}}$. D'après ce qui précède,
$$H_{c}(z) = 1+\sum_{n=1}^{\infty} \widetilde{h}_{n}(\Omega_{c})\,z^{n}\qquad;\qquad \frac{z\,H_{c}'(z)}{H_{c}(z)}=P_{c}(z)=\sum_{n=1}^{\infty} \tilp_{n}(\Omega_{c})\,z^{n}\,.$$
Le développement en série de $2/H_{c}(z)$ donne 
\begin{align*}
\frac{2}{H_{c}(z)}&=1+cz+(1-cz)\sqrt{1-4\left(\frac{z}{1-cz}\right)^{2}}\\
&=1+cz+(1-cz)\left(1-2\sum_{m=0}^{\infty}C_{m}\left(\frac{z}{1-cz}\right)^{2(m+1)}\right)\\
&=2\left(1-\sum_{m=0}^{\infty}C_{m}\,\frac{z^{2m+2}}{(1-cz)^{2m+1}}\right)\,,
\end{align*}
étant entendu que $\sum_{m=0}^{\infty}C_{m}\,x^{m}=\frac{1-\sqrt{1-4x}}{2x}$. Par conséquent,
$$H_{c}(z)=1+\sum_{m=0}^{\infty}\sum_{r=1}^{\infty}\left\{\sum_{m_{1}+m_{2}+\cdots+m_{r}=m}\!\!\!\!\! C_{m_{1}}C_{m_{2}}\cdots C_{m_{r}}\right\}\frac{z^{2m+2r}}{(1-cz)^{2m+r}}\,.
$$
Soit $f(m,r)$ la somme entre crochets. On a $f(m,0)=0$, $f(m,1)=C_{m}$, et compte tenu de la formule de récurrence pour les nombres de Catalan, $f(m,r)=f(m+1,r-1)-f(m+1,r-2)$. Par conséquent, on peut démontrer par récurrence que
$$f(m,r)=r\,\frac{2m+r-1!}{m!\,m+r!}\,.$$
Finalement, on utilise le développement $\frac{1}{(1-cz)^{2m+r}}=\sum_{l=0}^{\infty} \binom{2m+r+l-1}{l}\,c^{l}z^{l}$ pour obtenir:
\begin{align*}H_{c}(z)&=1+\sum_{m=0}^{\infty}\sum_{r=1}^{\infty}\sum_{l=0}^{\infty} r\frac{2m+r+l-1!}{l!\,2m+r-1!}\,\frac{2m+r-1!}{m!\,m+r!}\,c^{l}\,z^{2m+2r+l}\\
&=1+\sum_{n=2}^{\infty}z^{n}\left\{\sum_{k=1}^{\lfloor \frac{n}{2}\rfloor }\left(\sum_{r=1}^{k}r\, \frac{(n-r-1)^{\downarrow 2k-r-1}}{k!\,k-r!}\right)c^{n-2k}\right\}\,.
\end{align*}
Pour tous entiers positifs $a$, $b$ et $c$, $\binom{a+1}{b+c+1}=\sum_{r\geq 1} \binom{r}{b}\binom{a-r}{c}$ --- cette identité peut être démontrée bijectivement en regroupant les parties de taille $b+c+1$ dans $\lle 1,a+1\rre$ suivant la valeur de leur $(b+1)$-ième élément. En particulier, 
$$\sum_{r\geq 1}\,r\,\binom{n-r-1}{n-k-1}=\binom{n}{n-k+1}$$
si l'on prend $a=n-1$, $b=1$ et $c=n-k-1$. Ainsi, le terme entre parenthèses dans l'expression précédente de $H_{c}(z)$ est simplement
$$\frac{n^{\downarrow 2k}}{(n-k+1)(n-k)\,k!\,k-1!}\,,$$
et on obtient bien l'expression attendue pour $\widetilde{h}_{n}(\Omega_{c})$.
\bigskip
\bigskip

La preuve de l'identité pour les moments entrelacés est tout aussi laborieuse\footnote{On pourrait aussi utiliser les relations de Newton pour calculer récursivement les $\tilp_{n}(\Omega_{c})$ à partir des $\widetilde{h}_{n}(\Omega_{c})$, mais ces calculs reviennent essentiellement à ceux présentés ci-après.}. Pour commencer, on remarque que :
\begin{align*}P_{c}(z)&=-\frac{z}{1+cz+\sqrt{(1-cz)^{2}-4z^{2}}}\left(c - \frac{c(1-cz)+4z}{\sqrt{(1-cz)^{2}-4z^{2}}}\right)\\
&=\frac{z\,H_{c}(z)}{2}\left(\frac{c(1-cz)+4z}{\sqrt{(1-cz)^{2}-4z^{2}}}-c\right)\\
&=\frac{z\,H_{c}(z)}{2}\left(c\sum_{n=1}^{\infty}\binom{2n}{n}\left(\frac{z}{1-cz}\right)^{2n}+4\sum_{n=0}^{\infty}\binom{2n}{n}\left(\frac{z}{1-cz}\right)^{2n+1}\right)\\
&=\frac{z\,H_{c}(z)}{2}\left(\sum_{n=0}^{\infty}\sum_{r=0}^{\infty}\left\{cz\binom{2n+1+r}{r}\binom{2n+2}{n+1}+ 4\binom{2n+r}{r}\binom{2n}{n}\right\}c^{r}z^{2n+1+r}\right)\,.
\end{align*}
Si $n_{1}\geq 2$, le coefficient de $z^{n_{1}}$ dans la série $z\,H_{c}(z)/2$ est d'après ce qui précède égal à 
$$ \frac{1}{2}\sum_{k=1}^{\lfloor \frac{n_{1}-1}{2}\rfloor} \frac{n_{1}-1!}{(n_{1}-k)(n_{1}-k-1)\,k!\,k-1!\,n_{1}-1-2k!}\,c^{n_{1}-2k-1}\,\,;$$
et si $n_{2} \geq 1$, le coefficient de $z^{n_{2}}$ dans le second terme $T_{c}(z)$ du produit $P_{c}(z)$ est égal à 
\begin{align*}&2\sum_{l=0}^{\lfloor \frac{n_{2}}{2}\rfloor-1}\frac{n_{2}-1!}{n_{2}-2l-2!\,l!\,l+1!}\,c^{n_{2}-2l-1}+4\sum_{l=0}^{\lfloor\frac{n_{2}-1}{2} \rfloor}\frac{n_{2}-1!}{n_{2}-2l-1!\,l!\,l!}\,c^{n_{2}-2l-1}\\
&=\frac{2}{n_{2}}\sum_{l=0}^{\lfloor\frac{n_{2}-1}{2}\rfloor} \frac{n_{2}+1!}{n_{2}-2l-1!\,l!\,l+1!}\,c^{n_{2}-2l-1}\,.\end{align*}
De plus, le coefficient de $z$ dans $z\,H_{c}(z)/2$ vaut $1/2$. Par suite, comme $$[z^{n}](P_{c}(z))=\sum_{n_{1}+n_{2}=n}[z^{n_{1}}](zH_{c}(z)/2)\,\times\,[z^{n_{2}}](T_{c}(z))\,,$$ on voit que pour $n\geq 3$, $\tilp_{n}(\Omega_{c})$ est la somme des expressions suivantes :
\begin{align*}A&=\frac{1}{n-1}\sum_{k=1}^{\lfloor\frac{n}{2} \rfloor} \frac{n!}{n-2k!\,k!\,k-1!}\,c^{n-2k}\,\,;\\
B&=\!\!\!\!\!\!\sum_{\substack{n_{1}+n_{2}=n\\n_{1}\geq 2,\,\,n_{2}\geq 1\\2\leq 2k_{1}\leq n_{1}-1\\ 0 \leq 2k_{2} \leq n_{2}-1}}\frac{(n_{2}+1)\,n_{1}-1!\,n_{2}-1!\,n_{1}-k_{1}-2!}{n_{1}-k_{1}!\,k_{1}!\,k_{1}-1!\, n_{1}-2k_{1}-1!\,n_{2}-2k_{2}-1!\,k_{2}!\,k_{2}+1!}\,c^{n-2k_{1}-2k_{2}-2}\\
&=\!\!\!\!\!\!\sum_{\substack{n_{1}+n_{2}=n-3\\n_{1}\geq 0,\,\,n_{2}\geq 0\\0\leq 2k_{1}\leq n_{1}-1\\ 0 \leq 2k_{2} \leq n_{2}}}\frac{(n_{2}+2)\,n_{1}+1!\,n_{2}!\,n_{1}-k_{1}-1!}{n_{1}-k_{1}+1!\,k_{1}!\,k_{1}+1!\, n_{1}-2k_{1}-1!\,n_{2}-2k_{2}!\,k_{2}!\,k_{2}+1!}\,c^{n-2k_{1}-2k_{2}-4}\\
&=\!\!\!\!\!\!\sum_{\substack{0\leq u \leq n-4-2k\\0\leq k = k_{1}+k_{2}}} \frac{f(k_{1},k_{2},u,n)\,n-2k_{2}-2-u!\,2k_{2}+u!}{n-u-2k-4!\,u!\,k_{1}!\,k_{1}+1!\,k_{2}!\,k_{2}+1!}\,c^{n-2k-4}
\end{align*}
où $f(k_{1},k_{2},u,n)=\frac{2k_{2}+2+u}{(n-u-k_{1}-2k_{2}-2)(n-u-k_{1}-2k_{2}-3)}$. Les moments $\tilp_{n}(\Omega_{c})$ sont donc bien des polynômes en $c$ avec tous les termes de degré pair ou tous les termes de degré impair. Le terme de degré ${n-2}$ en $c$ provient seulement de $A$, et est égal à 
$$\frac{1}{n-1}\,\frac{n!}{n-2!}=n=\frac{n^{\downarrow 2} }{n-1\,1!\,0!}\,.$$
Ainsi, la formule de l'énoncé est vraie pour le coefficient de $c^{n-2}$. Pour les autres coefficients, il s'agit de montrer que si $k \geq 0$ et $n-2k-4\geq 0$, alors :
$$\sum_{\substack{0\leq u \leq n-4-2k\\ k = k_{1}+k_{2}}} \frac{f(k_{1},k_{2},u,n)\,n-2k_{2}-2-u!\,2k_{2}+u!}{n-u-2k-4!\,u!\,k_{1}!\,k_{1}+1!\,k_{2}!\,k_{2}+1!}=\frac{n^{\downarrow 2k+4}}{(n-k-2)(n-1)\,k+2!\,k!}\,.$$
En effet, le terme de droite est la différence des coefficients de $c^{n-4-2k}$ dans la formule de l'énoncé et dans le terme $A$. Cette dernière identité est équivalente à la première partie du lemme \ref{whatever}. Finalement, on laisse au lecteur le soin de vérifier que la formule pour $\tilp_{n}(\Omega_{c})$ reste vraie pour $n=1,2$.
\end{proof}
\bigskip
\bigskip

\noindent Dans l'algèbre des observables étendue $\obs^{+}$, introduisons les éléments
$$\widetilde{q}_{k,c}=\frac{\tilp_{k+1}}{(k+1)\,(\varSigma_{1})^{k/2}}-\sum_{l=1}^{\lfloor\frac{k+1}{2}\rfloor} \frac{k^{\downarrow 2l-1}}{(k+1-l)\,l!\,l-1!}\,c^{k+1-2l}\,(\varSigma_{1})^{1/2}\,.$$
Lorsque $c=0$, on retrouve les observables $\widetilde{q}_{k}$ présentées page \pageref{bizarobs}.
\begin{lemma}[Moments de la déviation d'un diagramme]\label{momentdeviation}
Soit $\lambda$ un diagramme de Young, $\lambda^{*}$ le diagramme continu renormalisé en abscisse et en ordonnée par un facteur $1/\sqrt{n}$, et $\Delta_{\lambda,c}(s)=\lambda^{*}(s)-\Omega_{c}(s)$. Pour tout entier $k$,
$$\frac{\sqrt{n}}{2}\int_{\R} s^{k}\,\Delta_{\lambda,c}(s)\,ds=\frac{\widetilde{q}_{k+1,c}(\lambda)}{k+1}\,.$$
\end{lemma}
\begin{proof}
D'après la seconde partie de la proposition \ref{momentwignerdeform}, l'observable $\widetilde{q}_{k,c}(\lambda)$ est simplement égale à 
$$\frac{(\varSigma_{1}(\lambda))^{1/2}}{k+1}\left(\frac{\tilp_{k+1}(\lambda)}{(\varSigma_{1}(\lambda))^{k+1/2}}-\tilp_{k+1}(\Omega_{c})\right)=\frac{\sqrt{n}}{k+1}\big(\tilp_{k+1}(\lambda^{*})-\tilp_{k+1}(\Omega_{c})\big)\,.$$
De plus, $\tilp_{k+2}(\lambda)=\int_{\R}s^{k+2}\,\sigma_{\lambda}''(s)\,ds=\frac{1}{(k+1)(k+2)}\int_{\R}s^{k}\,(\frac{\lambda(s)-|s|}{2})\,ds$ pour tout dia\-gram\-me continu $\lambda \in \cym$. Ainsi,
\begin{align*}
\frac{\widetilde{q}_{k+1,c}(\lambda)}{k+1}&=\frac{\sqrt{n}}{(k+2)(k+1)}\big(\tilp_{k+2}(\lambda^{*})-\tilp_{k+2}(\Omega_{c})\big)\\
&=\frac{\sqrt{n}}{2}\int_{\R}s^{k}\big((\lambda^{*}(s)-|s|)-(\Omega_{c}(s)-|s|)\big)\,ds=\frac{\sqrt{n}}{2}\int_{\R} s^{k}\,\Delta_{\lambda,c}(s)\,ds\,.
\end{align*}
\end{proof}

\begin{lemma}[Moments décalés de la déviation d'un diagramme]\label{relapse}
Avec les mêmes notations, pour tout entier $k$ et tout diagramme de Young $\lambda$,
\begin{align*}&\frac{\sqrt{n}}{2}\int_{\R}(s-c)^{k}\,\Delta_{\lambda,c}(s)\,ds=\\
&\frac{\sqrt{n}}{(k+1)(k+2)}\left(\sum_{2 \leq l \leq k+2}\!\! \binom{k+2}{l} (-c)^{k+2-l}\,\frac{\tilp_{l}(\lambda)}{(\varSigma_{1}(\lambda))^{l/2}}\,\,-\!\!\!\!\sum_{2 \leq 2m \leq k+2}\!\! \binom{k+2}{m} (-c)^{k+2-2m}\right).
\end{align*}
\end{lemma}
\begin{proof}
On développe $(s-c)^{k}$ par la formule du binôme de Newton :
\begin{align*}
\frac{\sqrt{n}}{2}\scal{(s-c)^{k}}{\Delta_{\lambda,c}}&=\sum_{l=0}^{k}\binom{k}{l}(-c)^{k-l}\,\frac{\sqrt{n}}{2}\int_{\R}s^{l}\,\Delta_{\lambda,c}(s)\,ds\\
&=\frac{1}{k+1}\sum_{l=0}^{k}\binom{k+1}{l+1}(-c)^{k-l}\,\widetilde{q}_{l+1,c}(\lambda)\\
&=\frac{1}{k+1}\sum_{l=1}^{k+1}\binom{k+1}{l}(-c)^{k+1-l}\,\widetilde{q}_{l,c}(\lambda)\\
&=\frac{\sqrt{n}}{(k+1)(k+2)}\sum_{l=2}^{k+2}\binom{k+2}{l}(-c)^{k+2-l}\,\frac{\tilp_{l}(\lambda)}{(\varSigma_{1}(\lambda))^{l/2}}-\sqrt{n}\,A(c)\,,
\end{align*}
où $A(c)$ est l'expression suivante :
\begin{align*}
A(c)&=\sum_{m=1}^{\lfloor \frac{k}{2}\rfloor +1} \sum_{l=2m-1}^{k+1} \frac{k!}{k+1-l!\,l-2m+1!\,m!\,m-1!\, (l+1-m)}\,(-1)^{k+1-l}c^{k+2-2m} \\
&=\sum_{m=0}^{\lfloor \frac{k}{2}\rfloor} \frac{k!}{k-2m!\,m!\,m+1!}\left\{\sum_{u=0}^{k-2m} \binom{k-2m}{u}\,\frac{(-1)^{k-u-2m}}{u+m+1}\right\}c^{k-2m}\,.
\end{align*}
Le terme entre crochets peut être calculé en écrivant $\frac{1}{u+m+1}=\int_{0}^{1} x^{u+m}\,dx$. Ainsi,
\begin{align*}
\sum_{u=0}^{k-2m} \binom{k-2m}{u}\,\frac{(-1)^{k-u-2m}}{u+m+1}&=\int_{0}^{1} \sum_{u=0}^{k-2m} \binom{k-2m}{u}\,(-1)^{k-u-2m}\,x^{u+m}\,dx\\
&=\int_{0}^{1}x^{m}\,(x-1)^{k-2m}\,dx=(-1)^{k-2m}\,\frac{k-2m!\,m!}{k-m+1!}\,.
\end{align*}
Par conséquent, $A(c)=\sum_{m=0}^{\lfloor \frac{k}{2}\rfloor} \frac{k!}{k-m+1!\,m+1!}\,(-c)^{k-2m}=\frac{1}{(k+1)(k+2)}\sum_{m=0}^{\lfloor \frac{k}{2}\rfloor} \binom{k+2}{m+1}\,(-c)^{k-2m}$, et on en déduit la formule de l'énoncé.
\end{proof}\bigskip\bigskip

Ceci étant, sous une mesure de Schur-Weyl de paramètres $\alpha=1/2$ et $c>0$, le théorème de \'Sniady \ref{factorasymptotic} s'applique de nouveau, avec pour paramètres :
\begin{align*}c_{l+1}&=\lim_{n \to \infty}k_{n}(\varSigma_{l})\,n^{-\frac{l+1}{2}}=\lim_{n \to \infty}k_{n}(R_{l+1})\,n^{-\frac{l+1}{2}}=c^{l-1}\,\,;\\
v_{l+1,m+1}&=\lim_{n \to \infty} k_{n}(\varSigma_{l},\varSigma_{m})\,n^{-\frac{l+m}{2}}=\lim_{n \to \infty} k_{n}(R_{l+1},R_{m+1})\,n^{-\frac{l+m}{2}}\\
&=\lim_{n \to \infty} k_{n}(\sigma_{l},\sigma_{m})\,n^{\frac{l+m}{2}}-lm\,c^{l+m-2}+\sum_{\substack{l=a_{1}+\cdots+a_{r}\\m=b_{1}+\cdots+b_{r}}}\!\! \frac{lm}{r}\,c_{a_{1}+b_{1}}\cdots c_{a_{r}+b_{r}}\\
&=-lm\,c^{l+m-2}+\sum_{\substack{l=a_{1}+\cdots+a_{r}\\m=b_{1}+\cdots+b_{r}}}\!\! \frac{lm}{r}\,c^{l+m-2r}=\sum_{r \geq 2} \binom{l}{r}\binom{m}{r} r\,\,c^{l+m-2r}\,.
\end{align*}
Ainsi, les observables (généralisées) 
$$W_{l \geq 1}=(\varSigma_{1})^{1/2}\left(\frac{\varSigma_{l}}{(\varSigma_{1})^{\frac{l+1}{2}}}-c^{l-1}\right)\qquad;\qquad X_{l\geq 2}=(\varSigma_{1})^{1/2}\left(\frac{R_{l}}{(\varSigma_{1})^{\frac{l}{2}}}-c^{l-2}\right)$$ sont asymptotiquement gaussiennes centrées, avec les covariances limites indiquées par la formule ci-dessus. Dans ce qui suit, nous développons les moments décalés de la déviation d'un diagramme en fonction de ces nouvelles observables.\bigskip

\begin{proposition}[Comportement asymptotique de fonctionnelles linéaires de la déviation]\label{linearfunctional}
Soit $p(s)$ un polynôme en $s$, et $\lambda \in \ym_{n}$ un diagramme de Young aléatoire tiré suivant une mesure de Schur-Weyl de paramètres $\alpha=1/2$ et $c>0$. La variable aléatoire
$$\sqrt{n}\scal{p(s)}{\Delta_{\lambda,c}(s)}=\sqrt{n}\int_{\R} p(s)\,\Delta_{\lambda,c}(s)\,ds$$
converge en loi vers une variable gaussienne centrée.
\end{proposition}
\begin{proof}
La partie importante de la proposition est en fait l'existence d'une limite et son caractère centré, voir la preuve de la proposition \ref{recovery} ci-après. D'après ce qui précède, pour tout $l$ plus grand que $2$, 
$$\sqrt{n}\left( \frac{R_{l}(\lambda)}{n^{\frac{l}{2}}} - c^{l-2}\right)=\sqrt{n}\left( R_{l}(\lambda^{*}) - R_{l}(\Omega_{c})\right)$$
converge vers une variable gaussienne centrée, et on a même convergence en loi jointe vers un vecteur gaussien (centré). De plus, le résultat reste vrai pour $l=1$ (dans ce cas la variable aléatoire est nulle pour tout $n$). Soit $\mu=(\mu_{1},\ldots,\mu_{r})\in \ym$. La déviation $\sqrt{n}\,(R_{\mu}(\lambda^{*})-R_{\mu}(\Omega_{c}))$ peut être développée en 
$$\sum_{i=0}^{r-1}R_{\mu_{1}\ldots\mu_{r-i-1}}(\lambda^{*})\,\left[\sqrt{n}\,(R_{\mu_{r-i}}(\lambda^{*})-R_{\mu_{r-i}}(\Omega_{c}))\right]\,R_{\mu_{r-i+1}\ldots\mu_{r}}(\Omega_{c})\,.$$
Chaque $R_{\mu_{1}\ldots\mu_{r-i-1}}(\lambda^{*})$ converge vers la constante $R_{\mu_{1}\ldots\mu_{r-i-1}}(\Omega_{c})$, et d'autre part, cha\-que déviation $\left[\sqrt{n}\,(R_{\mu_{r-i}}(\lambda^{*})-R_{\mu_{r-i}}(\Omega_{c}))\right]$ converge vers une variable gaussienne centrée ; de plus, toutes ses convergences sont jointes. Par conséquent, la limite des $\sqrt{n}(R_{\mu}(\lambda^{*})-R_{\mu}(\Omega_{c}))$ est un élément d'un espace gaussien dont toutes les variables sont centrées :
$$\forall \mu,\,\,\,\Delta_{n,c} (R_{\mu})=\sqrt{n}\,(R_{\mu}(\lambda^{*})-R_{\mu}(\Omega_{c})) \quad\longrightarrow \quad\mathcal{N}(0,(\sigma_{\mu,c})^{2})\,.$$
Comme précédemment, on peut même raffiner ce résultat et montrer que pour tout vecteur
$(\Delta_{n,c}(R_{\mu^{(1)}}), \ldots, \Delta_{n,c}(R_{\mu^{(s)}}))$ de déviations, on a convergence en loi jointe vers un vecteur gaussien centré. Comme les $R_{\mu}$ forment une base de l'algèbre $\obs$, on a le même résultat pour toute famille d'observables $f^{(1)},\ldots,f^{(s)}$ : le vecteur des fluctuations renormalisées $(\Delta_{n,c}(f^{(1)}),\ldots,\Delta_{n,c}(f^{(s)}))$ converge en loi jointe vers un vecteur gaussien centré. Or, le lemme \ref{momentdeviation} assure que
$$\frac{\sqrt{n}}{2}\scal{s^{k}}{\Delta_{\lambda,c}(s)}=\frac{\Delta_{n,c}(\tilp_{k+2})}{(k+1)(k+2)}$$
est une fluctuation renormalisée d'observable ; la proposition est donc démontrée.
\end{proof}\bigskip

\begin{lemma}[Formule de changement de base entre les moments entrelacés et les caractères centraux]\label{otherchangeofbasis}
Pour $k \geq 2$, $\tilp_{k}$ est la composante de plus haut poids de l'observable
$$\sum_{\substack{\mu=1^{m_{1}}2^{m_{2}}\cdots s^{m_{s}}\\|\mu|+\ell(\mu)=k }} \frac{k^{\downarrow \ell(\mu)}}{\prod_{i\geq 1} m_{i}!}\prod_{i\geq 1}(\varSigma_{i})^{m_{i}}\,.$$
\end{lemma}
\begin{proof}
On renvoie à \cite[proposition 3.7]{IO02} ; il s'agit essentiellement de combiner la formule de changement de base entre les moments de Frobenius et les moments entrelacés et la formule de changement de base entre les caractères centraux et les moments de Frobenius, voir page \pageref{wassermann}. Notons que cette formule permettra de se passer du degré de Kerov\footnote{La filtration de Kerov est néanmoins intéressante, et nous verrons dans le chapitre suivant qu'elle est tout à fait adaptée à l'étude asymptotique des mesures de Gelfand.} introduit à la fin du chapitre \ref{plancherel} ; ainsi, tous nos raisonnements utiliseront la filtration du poids sur l'algèbre d'observables $\obs$.
\end{proof}
\begin{examples}
Les premiers moments entrelacés s'écrivent dans la base des caractères centraux :
\begin{align*}
\tilp_{1}&=0\quad\qquad\,;\qquad\tilp_{2}=2\varSigma_{1}\\
\tilp_{3}&=3\varSigma_{2}\qquad;\qquad\tilp_{4}=4\varSigma_{3}+6\varSigma_{1,1}+2\varSigma_{1}\\
\tilp_{5}&=5\varSigma_{4}+20\varSigma_{2,1}+15\varSigma_{2}\\
\tilp_{6}&=6\varSigma_{5}+30\varSigma_{3,1}+15\varSigma_{2,2}+20\varSigma_{1,1,1}+30\varSigma_{1,1}+60\varSigma_{3}+2\varSigma_{1}
\end{align*}
Ainsi, on voit par exemple que $\tilp_{6}$ est la composante homogène de plus haut poids de $6\varSigma_{5}+30\varSigma_{3,1}+15\varSigma_{2,2}+20\varSigma_{1,1,1}$, et compte tenu de la factorisation en plus haut poids des caractères centraux, ceci est en accord avec le lemme énoncé ci-dessus.
\end{examples}\bigskip

\begin{lemma}[Deux identités sur les partitions]\label{twoidenpart}
Pour tous entiers positifs $s\geq l$,
$$\sum_{\substack{|\mu|=s\\ \ell(\mu)=l}}\frac{1}{\prod_{i \geq 1}m_{i}(\mu)!}=\frac{1}{l!}\,\binom{s-1}{l-1}\,.$$
Plus généralement, étant données des variables (aléatoires) $Y_{p\geq 1}$, on a
$$\sum_{\substack{|\mu|=s\\ \ell(\mu)=l}}\frac{1}{\prod_{i \geq 1}m_{i}(\mu)!}\left(\sum_{p \in \mu}Y_{p}\right)=\frac{1}{l-1!}\,\sum_{u=0}^{s-l} \binom{l-2+u}{u}\,Y_{s-l+1-u}\,.$$
\end{lemma}
\begin{proof}
La première formule est équivalente à 
$$\sum_{\substack{|\mu|=s\\ \ell(\mu)=l}}\binom{l}{m_{1}(\mu),\ldots,m_{s}(\mu)}=\binom{s-1}{l-1}\,,$$
et les deux termes décomptent le nombre de $l$-uplets d'entiers strictement positifs $(s_{1},\ldots,s_{l})$ tels que $s_{1}+\cdots+s_{l}=s$. La première identité est donc vraie, et pour la seconde formule, on peut effectuer la manipulation suivante :
\begin{align*}S_{2}&=\sum_{\substack{|\mu|=s\\ \ell(\mu)=l}}\frac{1}{\prod_{i \geq 1}m_{i}(\mu)!}\left(\sum_{p \in \mu}Y_{p}\right)=\sum_{\substack{|\mu|=s\\ \ell(\mu)=l}}\frac{1}{\prod_{i \geq 1}m_{i}(\mu)!}\left(\sum_{i=1}^{s-l+1}m_{i}(\mu)\,Y_{i}\right)\\
&=\sum_{\substack{|\mu|=s\\ \ell(\mu)=l}}\sum_{i=1}^{s-l+1}\frac{1}{\left(\prod_{\substack{j \geq 1\\ j \neq i}}m_{j}(\mu)!\right)\,(\max(0,m_{i}(\mu)-1))!}\,Y_{i}\,.
\end{align*}
En effet, la plus grande part possible d'une partition de longueur $l$ et de taille $s$ est $s-l+1$. Alors, en intervertissant l'ordre de sommation, on obtient :
$$S_{2}=\sum_{i=1}^{s-l+1}Y_{i}\left(\sum_{\substack{|\mu|=s-i\\\ell(\mu)=l-1}} \frac{1}{\prod_{i \geq 1}m_{i}(\mu)!}\right)=\frac{1}{l-1!}\sum_{i=1}^{s-l+1}\binom{s-i-1}{l-2}\,Y_{i}\,,$$
d'où la seconde formule en effectuant le changement d'indices $u=s-l+1-i$.
\end{proof}\bigskip

\begin{lemma}[Développement gaussien des moments entrelacés sous les mesures de Schur-Weyl]\label{nasty}
Pour tout entier $k \geq 1$, l'observable renormalisée $\tilp_{k}(\lambda)/n^{\frac{k-1}{2}}$ est égale à
$$\sqrt{n}\left(\sum_{l=1}^{\lfloor \frac{k}{2}\rfloor} \frac{k!}{(k-l)\,k-2l!\,l-1!\,l!}\,c^{k-2l }\right)+\sum_{l=1}^{\lfloor \frac{k}{2}\rfloor}\frac{k!}{k-l!\,l-1!} \left( \sum_{u=0}^{k-2l} \binom{l-2+u}{u}c^{u}\,W_{1+k-2l-u} \right)$$
plus une variable aléatoire $V_{k}(c)$ qui sous les mesures de Schur-Weyl $SW_{n,1/2,c}$ converge en probabilité vers une constante $L_{k}(c)$.
\end{lemma}
\begin{proof}
La proposition \ref{otherchangeofbasis} permet d'écrire :
$$\tilp_{k}=\left(\sum_{|\mu|+\ell(\mu)=k}\frac{k^{\downarrow \ell(\mu)}}{\prod_{i \geq 1}m_{i}(\mu)!} \prod_{i \geq 1}(\varSigma_{i})^{m_{i}}\right)+D_{k-1}\,,$$
où $D_{k-1}$ est une observable de poids $k-1$. Décomposons $D_{k-1}$ dans la base (graduée) des cumulants libres :
$$D_{k-1}=\sum_{\substack{|\mu| \leq k-1\\ m_{1}(\mu)=0}}d_{\mu,k}\,R_{\mu}\,.$$
Comme $R_{\mu}/n^{|\mu|/2}$ converge en probabilité et en moments vers $c^{|\mu|-2\ell(\mu)}$, la variable $V_{k}(c)=D_{k-1}/n^{(k-1)/2}$ converge en probabilité vers
$$L_{k}(c)=\sum_{\substack{|\mu|=k-1  \\ m_{1}(\mu)=0}} d_{\mu,k}\,c^{|\mu|-2\ell(\mu)}\,.$$
Concentrons-nous maintenant sur la partie restante $\tilp_{k}-D_{k-1}$. On utilise l'astuce de calcul suivante :
\begin{align*}\frac{\tilp_{k}-D_{k-1}}{n^{\frac{k-1}{2}}}&=\sqrt{n}\left(\sum_{|\mu|+\ell(\mu)=k}\frac{k^{\downarrow \ell(\mu)}}{\prod_{i \geq 1}m_{i}(\mu)!} \prod_{i \geq 1}\left(\frac{\varSigma_{i}}{n^{\frac{i+1}{2}}}\right)^{m_{i}(\mu)}\right)\\
&=\sqrt{n}\left(\sum_{|\mu|+\ell(\mu)=k}\frac{k^{\downarrow \ell(\mu)}}{\prod_{i \geq 1}m_{i}(\mu)!} \prod_{i \geq 1}\left(\frac{W_{i}}{\sqrt{n}}+c^{i-1}\right)^{m_{i}(\mu)}\right).
\end{align*}
Une partition $\mu$ étant fixée, développons le produit $\prod_{i\geq 1}\left(\frac{W_{i}}{\sqrt{n}}+c^{i-1}\right)^{m_{i}(\mu)}$. Les deux premiers termes du développement sont
$$A(\mu)+\frac{B(\mu)}{\sqrt{n}}=c^{|\mu|-\ell(\mu)}+\frac{1}{\sqrt{n}}\sum_{p \in \mu} W_{p} \,c^{|\mu|-\ell(\mu)-p+1}\,,$$
et le reste est un $O(n^{-1})$, de sorte que :
$$\frac{\tilp_{k}-D_{k-1}}{n^{\frac{k-1}{2}}}=\sqrt{n}\left(\sum_{|\mu|+\ell(\mu)=k}\frac{k^{\downarrow \ell(\mu)}\,A(\mu)}{\prod_{i \geq 1}m_{i}(\mu)!}\right) + \left(\sum_{|\mu|+\ell(\mu)=k}\frac{k^{\downarrow \ell(\mu)}\,B(\mu) }{\prod_{i \geq 1}m_{i}(\mu)!}\right)+ O(n^{-1/2})\,.$$
On utilise alors le lemme \ref{twoidenpart} pour simplifier les deux sommes. La première somme est simplement égale à
$$\sqrt{n}\left(\sum_{l=1}^{\lfloor \frac{k}{2}\rfloor}\sum_{\substack{\mu \in \ym_{k-l}\\ |\mu|=l}} \frac{k!}{k-l!}
\,\frac{c^{k-2l }}{\prod_{i \geq 1}m_{i}(\mu)!}\right)=\sqrt{n}\left(\sum_{l=1}^{\lfloor \frac{k}{2}\rfloor} \frac{k!}{(k-l)\,k-2l!\,l-1!\,l!}\,c^{k-2l }\right),$$
et pour la seconde somme, on pose $Y_{p}=W_{p}/c^{p-1}$, de fa\c con à obtenir
$$\sum_{l=1}^{\lfloor \frac{k}{2}\rfloor} \sum_{\substack{\mu \in \ym_{k-l}\\ \ell(\mu)=l}}\sum_{p \in \mu}\frac{k!}{k-l!}\,\frac{c^{k-2l}\,Y_{p}}{{\prod_{i \geq 1}m_{i}(\mu)!}}=\sum_{l=1}^{\lfloor \frac{k}{2}\rfloor}\frac{k!}{k-l!\,l-1!} \left( \sum_{u=0}^{k-2l} \binom{l-2+u}{u}Y_{1+k-2l-u} \right) c^{k-2l}\,.$$
Comme $Y_{1+k-2l-u}\,c^{k-2l}=W_{1+k-2l-u}\,c^{u}$, le lemme est démontré.
\end{proof}\bigskip

\begin{proposition}[Développement gaussien des moments décalés de la déviation d'un diagramme]\label{recovery}
Le $k$-ième moment translaté de la déviation $\Delta_{\lambda,c}$ examiné dans le lemme \ref{relapse} est égal à 
$$\frac{1}{k+1}\sum_{l=0}^{k-1}\left(\sum_{m=0}^{\lfloor \frac{l}{2} \rfloor} \binom{k+1}{m}\binom{k+1-2m}{l-2m}\,(-c)^{l-2m}\right)W_{k+1-l}(\lambda)\,,$$
plus une variable aléatoire qui converge vers $0$ sous les mesures de Schur-Weyl de paramètres $\alpha=1/2$ et $c$.
\end{proposition}
\begin{proof}
Compte tenu des lemmes \ref{relapse} et \ref{nasty}, $\frac{\sqrt{n}}{2}\int_{\R}(s-c)^{k}\,\Delta_{\lambda,c}(s)\,ds$ est la somme des quatres termes suivants :
\begin{align*}
A&=\sqrt{n}\,\left(\sum_{2\leq l \leq k+2}\sum_{m=1}^{\lfloor \frac{l}{2}\rfloor } \frac{k!}{(l-m)\,k+2-l!\,l-2m!\,m-1!\,m!}\,(-1)^{k+2-l}\, c^{k+2-2m}\right)\,;\\
B&=\sum_{2\leq l \leq k+2}\sum_{m=1}^{\lfloor \frac{l}{2}\rfloor }\sum_{u=0}^{l-2m} \frac{k!}{k+2-l!\,l-m!\,m-1!}\binom{m-2+u}{u}\,(-1)^{k+2-l}\,c^{k+2-l+u}\,W_{1+l-2m-u}\,\,;\\
C&=\sum_{2\leq l \leq k+2} \frac{k!}{l!\,k+2-l!}\,(-c)^{k+2-l}\,V_{l}(c)\,\,;\\
D&=-\sqrt{n}\,\left(\sum_{2\leq 2m \leq k+2}\frac{k!}{m!\,k+2-m!}\,(-c)^{k+2-2m}\right).
\end{align*}
Dans $A$, on intervertit l'ordre de sommation pour obtenir
$$\frac{A}{\sqrt{n}}=\sum_{2\leq 2m \leq k+2}\frac{k!\,(-c)^{k+2-2m}}{m!\,k+2-2m!\,m-1!}\left\{\sum_{l=2m}^{k+2}(-1)^{l}\frac{k+2-2m!}{(l-m)\,k+2-l!\,l-2m!}\right\},$$
et on calcule l'expression entre crochets en écrivant :
\begin{align*}\sum_{u=0}^{k+2-2m}\binom{k+2-2m}{u}\frac{(-1)^{u}}{u+m}&=\int_{0}^{1}\left(\sum_{u=0}^{k+2-2m}\binom{k+2-2m}{u}(-x)^{u}\right)x^{m-1}\,dx\\
&=\int_{0}^{1}(1-x)^{k+2-2m}\,x^{m-1}\,dx=\frac{k+2-2m!\,m-1!}{k+2-m!}\,.
\end{align*}
Ainsi, $A+D=0$. D'autre part, $B$ est une combinaison linéaire de variables asymptotiquement gaussiennes centrées, et $C$ tend vers  $\sum_{2\leq l \leq k+2} \frac{k!}{l!\,k+2-l!}\,(-c)^{k+2-l}\,L_{l}(c)$, qui est une constante. Mais compte tenu du lemme \ref{linearfunctional}, le $k$-ième moment translaté de la déviation converge vers une variable gaussienne centrée, donc
$$\sum_{2\leq l \leq k+2} \frac{k!}{l!\,k+2-l!}\,(-c)^{k+2-l}\,L_{l}(c)=0\,.$$
Finalement, on peut simplifier $B$ comme suit. D'une part, comme $W_{1}=0$, on peut supprimer les termes tels que $1+l-2m-u=1$. D'autre part, si l'on effectue le changement d'indices $1+l-2m-u=k+1-x$ avec $x \in \lle 0,k-1\rre$, alors on obtient :
 \begin{align*}
 B&=\sum_{x=0}^{k-1}W_{k+1-x}\left(\sum_{\substack{k+2\geq l \geq 2m \geq 2\\l-2m\geq k-x}} \frac{(-1)^{k+2-l}\,c^{x-2m+2}\,k!\,l-m-k+x-2!}{k+2-l!\,l-m!\,m-1!\,m-2!\,l-2m-k+x!}\right)\\
 &=\sum_{x=0}^{k-1}W_{k+1-x}\left(\sum_{y=0}^{\lfloor\frac{x}{2}\rfloor} \sum_{z=0}^{x-2y} \frac{(-1)^{x-2y-z}\,c^{x-2y}\,k!\,y+z-1!}{y!\,y-1!\,z!\,k+1-x+y+z!\,x-2y-z!}\right)
 \end{align*}
avec $y=m-1$ et $z=l-2m-k+x$. Enfin, compte tenu de la seconde partie du lemme \ref{whatever}, posant $Z=x-2y$, $\alpha=y-1$ et $\beta=y+k+1-x$, on peut clore la troisième somme :
\begin{align*}
B&=\sum_{x=0}^{k-1}W_{k+1-x}\left(\sum_{y=0}^{\lfloor \frac{x}{2} \rfloor}\frac{(-c)^{x-2y}\,k!}{y!\,y-1!\,x-2y!}\sum_{z=0}^{Z=x-2y} \frac{y-1+z!}{y+k+1-x+z!}\binom{Z}{z}(-1)^{z}\right)\\
&=\sum_{x=0}^{k-1}W_{k+1-x}\left(\sum_{y=0}^{\lfloor \frac{x}{2} \rfloor}\frac{(-c)^{x-2y}\,k!\,k+1-2y!}{y!\,x-2y!\,k+1-x!\,k+1-y!}\right)\\
&=\frac{1}{k+1}\sum_{x=0}^{k-1}W_{k+1-x}\left(\sum_{y=0}^{\lfloor \frac{x}{2} \rfloor} \binom{k+1}{y}\binom{k+1-2y}{x-2y}\,(-c)^{x-2y}\right)\,.
\end{align*}
On laisse au lecteur le soin de vérifier que l'abus de notation qui est fait en considérant des expressions telles que $y-1!$ avec $y=0$ ne change pas le résultat.
\end{proof}\bigskip

Soit $u_{k}(s)$ le polynôme de Chebyshev de seconde espèce, qu'on normalise de telle sorte que $u_{k}(2 \cos\theta)=\frac{\sin(k+1)\theta}{\sin \theta}$. Ces polynômes vérifient la relation de récurrence 
$$u_{k+2}(X)=X\,u_{k+1}(X)-u_{k}(X)\,.$$
Les premières valeurs sont $u_{0}(X)=1$, $u_{1}(X)=X$, $u_{2}(X)=X^{2}-1$, $u_{3}(X)=X^{3}-2X$ et $u_{4}(X)=X^{4}-3X^{2}+1$. 
\begin{theorem}[Polynômes de Chebyshev de seconde espèce et déviation d'un diagramme]\label{almosttheorem}
Pour tout entier $k$, l'observable 
$\frac{\sqrt{n}}{2}\scal{u_{k}(s-c)}{\Delta_{\lambda,c}}=\frac{\sqrt{n}}{2}\int u_{k}(s-c)\,\Delta_{\lambda,c}(s)\,ds$ est égale à 
$$\frac{1}{k+1}\sum_{l=0}^{k-1}\binom{k+1}{l}(-c)^{l}\,W_{k+1-l}(\lambda)\,,$$
plus des observables qui sous les mesures de Schur-Weyl de paramètre $1/2$ tendent en probabilité (et en moments) vers $0$ lorsque $n$ tend vers l'infini. 
\end{theorem}
\begin{proof}
On suit essentiellement le raisonnement de \cite[\S7]{IO02}. Le développement explicite du polynôme $u_{k}(X)$ est :
$$u_{k}(X)=\sum_{m=0}^{\lfloor \frac{k}{2}\rfloor} (-1)^{m}\binom{k-m}{m}\,X^{k-2m}$$
--- ceci peut par exemple être démontré à l'aide de la formule de récurrence sur les polynômes de Chebyshev. \'Etant donnés deux ensembles de variables $a_{0},a_{1},\ldots$ et $b_{0},b_{1},\ldots$, les relations suivantes sont équivalentes :
$$\left\{b_{k}=\sum_{m=0}^{\lfloor\frac{k}{2} \rfloor}\binom{k}{m} \,a_{k-2m}\right\}_{\!k \geq 0}\quad\iff\quad\left\{a_{k}=\sum_{m=0}^{\lfloor\frac{k}{2} \rfloor}(-1)^{m}\,\frac{k}{k-m}\binom{k-m}{m} \,b_{k-2m}\right\}_{\!k\geq 0}$$
voir \cite[p. 36]{IO02}. Notons $B_{k+1}/(k+1)$ la variable aléatoire calculée dans la proposition \ref{recovery}. Alors, la fonctionnelle linéaire de la déviation associé au polynôme de Chebyshev translaté $u_{k}(s-c)$ s'écrit :
\begin{align*}\frac{\sqrt{n}}{2}\scal{u_{k}(s-c)}{\Delta_{\lambda,c}(s)}&=\sum_{m=0}^{\lfloor \frac{k}{2} \rfloor} (-1)^{m} \binom{k-m}{m}\,\frac{B_{k+1-2m}}{k+1-2m}\\
&=\sum_{m=0}^{\lfloor \frac{k}{2} \rfloor} (-1)^{m} \binom{k+1-m}{m}\,\frac{B_{k+1-2m}}{k+1-m}\,.\end{align*}
Comme $B_{0}=B_{1}=0$, l'indice $m$ peut en réalité être pris jusqu'à $\lfloor \frac{k+1}{2}\rfloor$, donc l'expression ci-dessus est exactement $A_{k+1}/(k+1)$, où les $A$ et les $B$ vérifient les relations présentées précédemment. Il suffit donc d'identifier les variables aléatoires $A_{k}$ telles que 
$$\forall k,\,\,\,\sum_{l=0}^{k-1}\left(\sum_{m=0}^{\lfloor \frac{l}{2} \rfloor} \binom{k+1}{m}\binom{k+1-2m}{l-2m}\,(-c)^{l-2m}\right)W_{k+1-l}=\sum_{m=0}^{\lfloor \frac{k+1}{2}\rfloor} \binom{k+1}{m}\,A_{k+1-2m}\,.$$
Dans le membre de gauche, inversons l'ordre de sommation :
\begin{align*}B_{k+1}&=\sum_{m=0}^{\lfloor \frac{k-1}{2}\rfloor}\binom{k+1}{m}\left(\sum_{l=2m}^{k-1}\binom{k+1-2m}{l-2m}\,(-c)^{l-2m}\,W_{k+1-l}\right)\\ 
&=\sum_{m=0}^{\lfloor \frac{k-1}{2}\rfloor}\binom{k+1}{m}\left(\sum_{l=0}^{k-1-2m}\binom{k+1-2m}{l}\,(-c)^{l}\,W_{k+1-2m-l}\right)\\
&=\sum_{m=0}^{\lfloor \frac{k+1}{2}\rfloor}\binom{k+1}{m}\left(\sum_{l=0}^{k-1-2m}\binom{k+1-2m}{l}\,(-c)^{l}\,W_{k+1-2m-l}\right)\quad\text{car }W_{0}=W_{1}=0\,.
\end{align*}
On a donc $A_{k+1}=\sum_{l=0}^{k-1}\binom{k+1}{l}\,(-c)^{l}\,W_{k+1-l}$, ce qu'il fallait démontrer.
\end{proof}\bigskip

Finalement, la proposition \ref{almosttheorem} conduit à un théorème central limite pour les mesures de Schur-Weyl tout à fait analogue au théorème \ref{secondasymptoticplancherel} :
\begin{theorem}[Déviation des diagrammes par rapport à leur forme limite sous les mesures de Schur-Weyl de paramètre $\alpha=1/2$, \cite{Mel10b}]\label{schurweylplus}
Soit $(\xi_{k})_{k \geq 2}$ une famille de variables gaussiennes indépendantes, centrées et toutes de variance $1$. Sous une mesure de Schur-Weyl de paramètres $\alpha=1/2$ et $c>0$, pour tout $k \geq 1$, la variable
$$\frac{\sqrt{n}}{2}\scal{u_{k}(s-c)}{\Delta_{\lambda,c}}=\frac{\sqrt{n}}{2}\int_{\R}u_{k}(s-c)\,\Delta_{\lambda,c}(s)\,ds$$
converge en loi vers $\xi_{k+1}/\sqrt{k+1}$. Par conséquent, au sens des distributions tempérées sur l'intervalle $[c-2,c+2]$, le processus $\sqrt{n}\,\Delta_{\lambda,c}(s)$ converge en loi vers le processus gaussien généralisé
$$\Delta_{c}(s)=\Delta(c+2\cos\theta)=\frac{2}{\pi}\sum_{k=2}^{\infty}\frac{\xi_{k}}{\sqrt{k}}\,\sin(k\theta)\,,$$
c'est-à-dire le processus du théorème \ref{secondasymptoticplancherel} décalé en abscisse de $c$. 
\end{theorem}
\noindent Ainsi, le processus de théorème central limite de Kerov est universel pour les mesures de Schur-Weyl de paramètre $\alpha=1/2$, $c>0$ ou $\alpha>1/2$ (ce qui correspond au cas $c=0$). 
\begin{proof}
Par rapport au cas $c=0$, une différence essentielle est la non indépen\-dan\-ce asymptotique des variables $W_{l}$ : ces variables convergent en lois fini-dimension\-nelles vers un vecteur gaussien  $(W_{l,\infty})_{l \geq 2}$ dont la matrice de covariance n'est pas diagonale. Néanmoins, si l'on applique le processus d'orthonormalisation de Gram-Schmidt à la base $(W_{l,\infty})_{l \geq 2}$ de l'espace gaussien limite, on peut retrouver des variables gaussiennes indépendantes. Ainsi, étant donné un vecteur gaussien $(W_{l,\infty})_{l \geq 2}$ dont les covariances sont données par 
$$\mathrm{cov}(W_{l,\infty},W_{m,\infty})=\sum_{r \geq 2}\binom{l}{r}\binom{m}{r}\,r\,c^{l+m-2r}\,,$$
notons $(T_{l,\infty})_{l \geq 2}$ la famille de variables aléatoires définies par :
$$T_{l,\infty}=\sum_{j=0}^{l-2}\binom{l}{j}(-c)^{j}\,W_{l-j,\infty}\,.$$
Les relations de changement de base inverses sont bien sûr $W_{l,\infty}=\sum_{j=0}^{l-2}\binom{l}{j}\,c^{j}\,T_{l-j,\infty}$\,, et on déduit de ceci les covariances des variables gaussiennes $T_{l,\infty}$ :
$$\mathrm{cov}(T_{l,\infty},T_{m,\infty})=\mathbb{1}_{l=m}\,l\,.$$
Autrement dit, $(\xi_{l})_{l \geq 2}=(T_{l,\infty}/\sqrt{l})_{l \geq 2}$ est l'orthonormalisée de Gram-Schmidt de $(W_{l,\infty})_{l \geq 2}$ dans l'espace gaussien. Or, d'après le théorème \ref{almosttheorem}, 
$$\frac{\sqrt{n}}{2}\scal{u_{k}(s-c)}{\Delta_{\lambda,c}}\simeq_{n \to \infty} \frac{1}{k+1}\sum_{l=0}^{k-1}\binom{k+1}{l}\,(-c)^{l}\,W_{k+1-l,\infty}=\frac{T_{k+1,\infty}}{k+1} = \frac{\xi_{k+1}}{\sqrt{k+1}}$$
lorsque $n$ tend vers l'infini. La première partie du théorème est donc établie, et la seconde partie s'en déduit exactement comme dans la section \ref{sniady}.
\end{proof}\bigskip

Lorsque $c<1$, le théorème \ref{schurweylplus} décrit entièrement la déviation asymptotique, puisque celle-ci est entièrement incluse dans l'intervalle $[c-2,c+2]$. Au contraire, lorsque $c\geq 1$, la zone $s \in [-1/c,c-2]$ est imparfaitement décrite par le théorème, car les polynômes de Chebyshev <<~décalés~>> $u_{k}(s-c)$ forment une base orthogonale de l'espace $\leb^{2}(d\mu_{0}(s-c))$, et cette mesure --- c'est-à-dire la loi de Wigner décalée d'un facteur $c$ --- ne charge pas l'intervalle $[-1/c,c-2]$. Un fait trivial mérite néanmoins d'être mentionné : dans cette zone, le diagramme reste toujours sous sa forme limite, puisque la mesure de Schur-Weyl ne charge que les diagrammes de longueur $\ell(\lambda) \leq N=n^{1/2}/c$.\bigskip
\bigskip

Concluons cette section par une remarque. Il faut admettre que la preuve que nous donnons du théorème \ref{almosttheorem} repose sur des simplifications miraculeuses dans les calculs, et l'universalité du processus de Kerov ne paraît pas en l'état être un résultat profond sur les partitions aléatoires. Dans le cas de la mesure de Plancherel, nous avions noté la coïncidence suivante : les polynômes orthogonaux de la mesure de transition limite $\mu_{0}$ (la loi de Wigner), c'est-à-dire les polynômes de Chebyshev de seconde espèce, formaient également une base orthogonale de l'espace $\leb^{2}(\Delta)$, c'est-à-dire que les
$$ \lim_{n \to \infty} \frac{\sqrt{n}}{2}\int_{\R} u_{k}(s) \Delta_{\lambda}(s)\,ds = \frac{1}{2}\int_{\R} u_{k}(s) \Delta(s)\,ds $$ 
sont des variables gaussiennes indépendantes engendrant l'espace gaussien limite. Mais pour les mesures de Schur-Weyl, les polynômes de Chebyshev décales $u_{k}(s-c)$, qui forment une base orthogonale de l'espace $\leb^{2}(\Delta_{c})$, ne sont pas les polynômes orthogonaux associés à la mesure $\mu_{c}$ --- ce sont les polynômes orthogonaux associés à la mesure
$$ \mathbb{1}_{s \in [c-2,c+2]}\,\frac{\sqrt{4-(s-c)^{2}}}{2\pi}\,ds\,.$$
Il semble donc qu'on ne puisse pas employer dans un contexte général des méthodes de polynômes orthogonaux, et l'apparition des polynômes de Chebyshev de seconde espèce dans l'étude des déviations semble plus liée à la combinatoire des observables de diagrammes qu'à l'expression de la forme limite.
\bigskip

\section{Asymptotique pour $\alpha < 1/2$}\label{schurminus}
Finalement, étudions le troisième régime asymptotique des mesures de Schur-Weyl, c'est-à-dire avec $\alpha<1/2$ et $c>0$. Comme précédemment, on peut calculer explicitement les espérances $SW_{n,\alpha,c}[\varSigma_{\mu}]$ pour toute partition $\mu$ :
$$SW_{n,\alpha,c}[\varSigma_{\mu}]=n^{\downarrow |\mu|}\,N^{\ell(\mu)-|\mu|}\simeq c^{|\mu|-\ell(\mu)} \,n^{|\mu|+\alpha\ell(\mu)-\alpha|\mu|}\,.$$
Dans ce contexte, il est utile d'introduire une nouvelle filtration sur l'algèbre $\obs$. Ainsi, pour $\alpha<1/2$, définissons le \textbf{$\alpha$-degré} d'une observable $p_{\mu}$ par 
\begin{align*}\deg_{\alpha}(p_{\mu})&=(1-\alpha)|\mu|+\alpha\ell(\mu)\\
&=(1-2\alpha)|\mu|+\alpha(|\mu|+\ell(\mu))=(1-2\alpha)\deg (p_{\mu})+\alpha \,\mathrm{wt}(p_{\mu})\,.\end{align*}
Comme $\deg (p_{\mu})=\sum_{i=1}^{\ell(\mu)} \deg (p_{\mu_{i}})$ et $\mathrm{wt}(p_{\mu})=\sum_{i=1}^{\ell(\mu)} \mathrm{wt}(p_{\mu_{i}})$, le $\alpha$-degré est multiplicatif, et peut donc être étendu en une filtration d'algèbres sur $\obs$. D'autre part, on constate que le $\alpha$-degré interpole le degré des observables et le poids des observables ; ainsi, $\deg_{0}(\cdot)=\deg(\cdot)$ et $\deg_{1/2}(\cdot)=\frac{\mathrm{wt}(\cdot)}{2}$.
\begin{lemma}[Factorisation asymptotique pour le $\alpha$-degré]
Pour toute partition $\mu$, $\varSigma_{\mu}$ est d'$\alpha$-degré $(1-\alpha)|\mu|+\alpha\ell(\mu)$, et sa composante homogène de plus haut $\alpha$-degré est $p_{\mu}$. Par conséquent, les caractères centraux ont la propriété de factorisation asymptotique en plus haut $\alpha$-degré :
$$\varSigma_{\mu_{1}} \,*\,\varSigma_{\mu_{2}}=\varSigma_{\mu_{1}\sqcup \mu_{2}}+(\text{termes d'$\alpha$-degré inférieur})\,.$$
\end{lemma}
\begin{proof}
On décompose $\varSigma_{\mu}$ dans la base des $p_{\nu}$, et on utilise l'expression de $\deg_{\alpha}$ comme interpolation entre le degré et le poids des observables. Ainsi, si $\varSigma_{\mu}=\sum_{\mu} c_{\mu}^{\nu}\,p_{\nu}$, alors pour toute partition $\nu$ telle que $c_{\mu}^{\nu}\neq 0$, on a :\vspace{1mm}
\begin{itemize}
\item $|\nu|\leq |\mu|$ avec égalité si et seulement si $\nu=\mu$.\vspace{1mm}
\item $|\nu|+\ell(\nu)\leq |\mu|+\ell(\mu)$.\vspace{1mm}
\end{itemize}
Comme $1-2\alpha>0$, on en déduit que le $\alpha$-degré de $\varSigma_{\mu}$ est égal à $$(1-2\alpha)|\mu|+\alpha(|\mu|+\ell(\mu))=(1-\alpha)|\mu|+\alpha\ell(\mu)\,,$$ et l'unique moment $p_{\nu}$ apparaissant dans $\varSigma_{\mu}$ et d'$\alpha$-degré $(1-\alpha)|\mu|+\alpha\ell(\mu)$ est $p_{\mu}$. La propriété de factorisation asymptotique s'en déduit immédiatement. \end{proof}
\begin{theorem}[Convergence en probabilité des observables sous les mesures de Schur-Weyl de paramètre $\alpha<1/2$, \cite{FM10}]\label{convobsschurminus}
Sous les mesures de Schur-Weyl $SW_{n,\alpha,c}$ avec $\alpha<1/2$, on a :
$$\frac{p_{\mu}(\lambda)}{n^{(1-\alpha)|\mu|+\alpha\ell(\mu)}}\to_{SW_{n,\alpha,c}} c^{|\mu|-\ell(\mu)}\qquad;\qquad\frac{\varSigma_{\mu}(\lambda)}{n^{(1-\alpha)|\mu|+\alpha\ell(\mu)}}\to_{SW_{n,\alpha,c}} c^{|\mu|-\ell(\mu)}$$
pour toute partition $\mu$.
\end{theorem}
\begin{proof}
Pour toute observable de diagrammes $f$, $SW_{n,\alpha,c}[f]=\esper[f]$ est un $O(n^{\deg_{\alpha}(f)})$ ; en effet, ceci est vrai sur la base des caractères centraux. On sait déjà que 
$$\esper\left[\frac{\varSigma_{\mu}}{n^{(1-\alpha)|\mu|+\alpha\ell(\mu)}}\right]\to c^{|\mu|-\ell(\mu)}$$
lorsque $n$ tend vers l'infini. Mais de plus, $\esper\left[\left(\frac{\varSigma_{\mu}}{n^{(1-\alpha)|\mu|+\alpha\ell(\mu)}}\right)^{2}\right]$ est égal à 
\begin{align*}&\esper\left[\frac{\varSigma_{\mu\sqcup\mu}}{n^{2(1-\alpha)|\mu|+2\alpha\ell(\mu)}}\right]+\frac{\esper[\text{observable d'$\alpha$-degré inférieur à $2(1-\alpha)|\mu|+2\alpha\ell(\mu)$}]}{n^{2(1-\alpha)|\mu|+2\alpha\ell(\mu)}}\\
&=(c^{|\mu|-\ell(\mu)})^{2}+o(1)\,.
\end{align*}
Par Bienaymé-Chebyshev, on en déduit la convergence en probabilité de $\varSigma_{\mu}/n^{\deg_{\alpha}(\varSigma_{\mu})}$ vers $c^{|\mu|-\ell(\mu)}$. Et comme la composante de plus haut $\alpha$-degré de $\varSigma_{\mu}$ est $p_{\mu}$, on a la même propriété pour les moments de Frobenius renormalisés.
\end{proof}
\noindent Ce théorème impliquera l'existence d'une forme limite pour les diagrammes renormalisés de fa\c con non isotrope ; \emph{cf.} le théorème \ref{schurweylminus}.
\bigskip
\bigskip

Compte tenu des résultats connus pour les mesures de Plancherel, les $q$-mesures de Plancherel et les mesures de Schur-Weyl de paramètre $\alpha \geq 1/2$, il est naturel de chercher à préciser le théorème \ref{convobsschurminus}  par des résultats de déviation gaussienne. Dans ce contexte, les techniques de cumulants d'observables fonctionnent à nouveau, mais elles donnent seulement un résultat partiel. Pour commencer, démontrons l'analogue pour le $\alpha$-degré des lemmes \ref{magnitudedj} et \ref{magnitude}. Dans ce qui suit, on note $k_{\alpha}$ les cumulants d'observables sous les mesures de Schur-Weyl $SW_{n,\alpha,c}$ avec $\alpha<1/2$, et $k_{\alpha}^{\bullet}$ les cumulants disjoints sous ces mêmes mesures.
\begin{lemma}[Ordre de grandeur des $(\alpha<1/2)$-cumulants disjoints]\label{alphamagnitudedj}
Pour toutes observables de diagrammes $x_{1},\ldots,x_{r}$, $k_{\alpha}^{\bullet}(x_{1},\ldots,x_{r})=O(n^{\deg_{\alpha}(x_{1}) + \cdots+\deg_{\alpha}(x_{r})-r+1})$.
\end{lemma}
\begin{proof}
De nouveau, il suffit de montrer le résultat pour une base algébrique de $\obs$, par exemple les caractères centraux des cycles. Ainsi, dans tous ce qui suit, les $x_{j}$ seront des observables $\varSigma_{i_{j}}$. Le produit disjoint de telles observables $\varSigma_{i_{1}},\ldots,\varSigma_{i_{r}}$ est simplement $\varSigma_{i_{1},i_{2},\ldots,i_{r}}$, et l'espérance d'une telle variable est 
$$\esper[\varSigma_{i_{1},\ldots,i_{r}}]=n^{\downarrow i_{1}+i_{2}+\cdots+i_{r}}\,\frac{1}{N^{(i_{1}-1)+\cdots+(i_{r}-1)}}\,.$$
Exactement comme dans la preuve du lemme \ref{magnitudedj}, il existe donc des nombres $\alpha_{i}$ tels que 
$$\esper[\varSigma_{i_{1},\ldots,i_{r}}]=\left(\prod_{j=1}^{r}\alpha_{i_{j}}\right)n^{\downarrow i_{1}+i_{2}+\cdots+i_{r}}\,,$$
et la seule différence par rapport au cas précédemment traité est que les $\alpha_{i}=N^{-(i-1)}$ dépendent maintenant de $n$. Néanmoins, on peut suivre exactement le même raisonnement que pour le lemme \ref{magnitudedj}, et ainsi,
$$k_{\alpha}^{\bullet}(\varSigma_{i_{1}},\ldots,\varSigma_{i_{r}})= \left(\prod_{j=1}^{r}\alpha_{i_{j}}\right)k_{\alpha}^{\bullet}(\varSigma_{(1^{i_{1}})},\ldots,\varSigma_{(1^{i_{r}})})=\left(\prod_{j=1}^{r}\alpha_{i_{j}}\right)O(n^{i_{1}+i_{2}+\cdots+i_{r}-r+1})\,.$$
Maintenant, $\prod_{j=1}^{r}\alpha_{i_{j}}=c^{\alpha(i_{1}+\cdots+i_{r}-r)} n^{\alpha r-\alpha(i_{1}+\cdots+i_{r})}$, donc 
$$k_{\alpha}^{\bullet}(\varSigma_{i_{1}},\ldots,\varSigma_{i_{r}})=O(n^{(1-\alpha)(i_{1}+\cdots+i_{r}) +\alpha r - r+1})=O(n^{\deg_{\alpha}(\varSigma_{i_{1}}) + \cdots+\deg_{\alpha}(\varSigma_{i_{r}})-r+1 })\,.\vspace{-5mm}~$$
\end{proof}\bigskip

\begin{lemma}[Ordre de grandeur des $(\alpha<1/2)$-cumulants]\label{alphamagnitude}
Pour toutes observables de diagrammes $x_{1},\ldots,x_{r}$, $k_{\alpha}(x_{1},\ldots,x_{r})=O(n^{\deg_{\alpha}(x_{1}) + \cdots+\deg_{\alpha}(x_{r})-r+1})$.
\end{lemma}
\begin{proof}
De nouveau, on est amené à estimer les cumulants identité $k^{\id}$ d'observables, cette fois-ci vis-à-vis de la filtration du $\alpha$-degré. Dans la preuve du lemme \ref{magnitude}, on a montré que :
$$\deg (k^{\id}(x_{1},\ldots,x_{r}))\leq \deg( x_{1})+\cdots+\deg( x_{r})-r+1\,.$$
D'autre part, il est démontré dans \cite{Sni06b} que :
$$\mathrm{wt}(k^{\id}(x_{1},\ldots,x_{r}))\leq \mathrm{wt}(x_{1})+\cdots+\mathrm{wt}(x_{r})-2r+2\,.$$
Comme $\deg_{\alpha}(\cdot)=(1-2\alpha)\deg(\cdot)+\alpha\mathrm{wt}(\cdot)$, on en déduit immédiatement :
$$\deg_{\alpha} (k^{\id}(x_{1},\ldots,x_{r}))\leq \deg_{\alpha} (x_{1})+\cdots+\deg_{\alpha} (x_{r})-r+1\,.$$
Ainsi, un cumulant identité de caractères centraux de cycles $k^{\id}(\varSigma_{i_{1}},\ldots,\varSigma_{i_{r}})$ est toujours d'$\alpha$-degré inférieur à $(1-\alpha)(i_{1}+\cdots+i_{r})+\alpha r-(r-1)$. On écrit alors :
$$k_{\alpha}(\varSigma_{i_{1}},\ldots,\varSigma_{i_{r}})=\sum_{\pi \in \mathfrak{Q}(\lle 1,r\rre)} k^{\bullet}\big(k_{\alpha}^{\id}(\varSigma_{i_{j \in \pi_{1}}}),k^{\id}(\varSigma_{i_{j \in \pi_{2}}}),\ldots,k^{\id}(\varSigma_{i_{j \in \pi_{l}}})\big)\,,$$
Le $\alpha$-degré d'un $k^{\id}(\varSigma_{i_{j \in \pi_{k}}})$ est inférieur à $(1-\alpha)\sum_{j \in \pi_{k}} i_{j}+\alpha |\pi_{k}| - (| \pi_{k}|-1)$ d'après ce qui précède. Une partition $\pi$ étant fixée, la somme de ces $\alpha$-degrés est donc toujours inférieure à $(1-\alpha)(i_{1}+\cdots+i_{r}-r)+\ell(\pi)$, et en utilisant le lemme \ref{alphamagnitudedj}, on conclut que le cumulant disjoint correspondant est un
$$O(n^{(1-\alpha)(i_{1}+\cdots+i_{r}-r)+1})=O(n^{\deg_{\alpha}(\varSigma_{i_{1}})+\cdots+\deg_{\alpha}(\varSigma_{i_{r}})-r+1})\,.\vspace{-5mm}~$$
\end{proof}\bigskip

Soit $\beta>0$ un paramètre réel. On introduit exactement comme dans le chapitre \ref{qplancherelmeasure} des déviations renormalisées 
$$ Z_{k,n,\alpha,c}=n^{\beta}\left(\frac{\varSigma_{k}(\lambda)-\esper[\varSigma_{k}]}{n^{(1-\alpha)k+\alpha}}\right).$$ 
Calculons les covariances de ces variables centrées. Pour commencer, remarquons que dans un produit $\varSigma_{l}\,\varSigma_{m}$ :\vspace{2mm}
\begin{enumerate}
\item Vis-à-vis de la filtration du degré, les deux principaux termes sont $\varSigma_{l,m}$ (de degré $l+m$) et $ml\,\varSigma_{l+m-1}$ (de degré $l+m-1$) ; tous les autres termes sont de degré au plus égal à $l+m-2$.\vspace{2mm}
\item Vis-à-vis de la filtration du poids, le principal terme est $\varSigma_{l,m}$ (de poids $l+m+2$), et tous les autres termes sont de poids inférieur à $l+m$ (et non $l+m+1$), voir par exemple \cite[corollaire 3.8]{Sni06a}. En particulier, $ml\,\varSigma_{l+m-1}$ est exactement de poids $l+m$.\vspace{2mm}
\end{enumerate}
En utilisant l'expression du $\alpha$-degré en fonction du degré et du poids, on conclut que :\vspace{2mm}
\begin{enumerate}
\item Le terme $\varSigma_{l,m}$ est d'$\alpha$-degré $(1-\alpha)(l+m)+2\alpha$, et le terme $ml\,\varSigma_{l+m-1}$ est d'$\alpha$-degré $(1-\alpha)(l+m-1)+\alpha=(1-\alpha)(l+m)+2\alpha-1$.\vspace{2mm}
\item Tous les autres termes sont d'$\alpha$-degré au plus égal à $(1-2\alpha)(l+m-2)+\alpha(l+m)=(1-\alpha)(l+m)+4\alpha-2$.\vspace{2mm}
\end{enumerate}
Comme dans la preuve du théorème \ref{devgaussianobs}, on décompose $k_{\alpha}(\varSigma_{l},\varSigma_{m})$ en $k_{\alpha}^{\bullet}(\varSigma_{l},\varSigma_{m})+\esper[\varSigma_{l,m}]-\esper[\varSigma_{l}]\,\esper[\varSigma_{m}]$. Le premier terme donne
\begin{align*}\esper[\varSigma_{l}\,\varSigma_{m}-\varSigma_{l,m}]&=\esper[ml\,\varSigma_{l+m-1}]+O(n^{(1-\alpha)(l+m)+4\alpha-2})\,,\\
&=ml\,c^{l+m-2}\,n^{(1-\alpha)(l+m)+2\alpha-1}+O(n^{(1-\alpha)(l+m)+4\alpha-2})\end{align*}
et comme dans la preuve du théorème \ref{devgaussianobs}, le second terme donne
$$N^{2-l-m}(n^{\downarrow l+m}-n^{\downarrow l}\,n^{\downarrow m})=-ml\,c^{l+m-2}\,n^{(1-\alpha)(l+m)+2\alpha-1}+O(n^{(1-\alpha)(l+m)+2\alpha-2})\,.$$
Ainsi, les deux termes de plus haut degré se simplifient, et 
$$k_{\alpha}(Z_{l,n,\alpha,c},Z_{m,n,\alpha,c})=O(n^{2\alpha+2\beta-2})+O(n^{2\beta-2})=O(n^{2\alpha+2\beta-2})\,.$$
Pour tout $\beta<1-\alpha$, les covariances tendent donc vers $0$, et :
\begin{proposition}[Déviation des observables sous les mesures de Schur-Weyl de paramètre $\alpha<1/2$]\label{notprecisenough}
Si $\alpha<1/2$, lorsque $n$ tend vers l'infini,
$$\forall \eps>0,\,\,\,\frac{\varSigma_{k}(\lambda)-\esper[\varSigma_{k}]}{n^{(1-\alpha)k+\alpha}}=o\left(\frac{1}{n^{1-\alpha-\eps}}\right)$$
la convergence s'entendant par exemple en probabilité.
\end{proposition}\bigskip\bigskip

Partant, il est naturel de s'intéresser aux déviations renormalisées de paramètre $\beta=1-\alpha$, et on fixe donc cette valeur de $\beta$ dans ce qui suit.
Compte tenu du lemme \ref{alphamagnitude}, les cumulants d'ordre $r\geq 2$ des variables $Z_{k,n,\alpha,c}$ sont des $O(n^{1-\alpha r})$ :
$$k_{\alpha}(Z_{k_{1},n,\alpha,c},\ldots,Z_{k_{r},n,\alpha,c})=\frac{n^{(1-\alpha)r}\,O(n^{(1-\alpha)(k_{1}+\cdots+k_{r})+\alpha r -r+1})}{n^{(1-\alpha)(k_{1}+\cdots+k_{r})+\alpha r}}=O(n^{1-\alpha r})\,.$$
Malheureusement, cette estimation n'est pas suffisante pour démontrer la convergence gaussienne : par exemple, si $\alpha=1/3$, alors on sait juste que les troisièmes cumulants sont des $O(1)$. Le calcul qui précède pour $r=2$ montre que l'estimation n'est pas optimale : en effet, pour toute valeur de $\alpha$, les cumulants d'ordre $2$ sont bien des $O(1)$, et non des $O(n^{1-2\alpha})$. Nous reviendrons plus loin sur ce problème. \bigskip\bigskip

Supposons pour l'instant $\alpha>1/3$. Alors, $1-\alpha r$ est bien strictement inférieur à $0$ pour $r\geq 3$, donc $(Z_{k,n,\alpha,c})_{k \geq 2}$ converge en lois fini-dimensionnelles vers un processus gaussien $(Z_{k,\infty,\alpha,c})_{k \geq 2}$. Les covariances de ce processus limite correspondent aux termes de $\alpha$-degré exactement égal à $(1-\alpha)(l+m)+4\alpha-2$ dans le produit $\varSigma_{l}\,\varSigma_{m}$. Si l'on écrit ce produit comme somme sur des appariements $M$ de termes $\varSigma_{\rho(M)}$, alors les termes $\varSigma_{\rho(M)}$ d'$\alpha$-degré égal à $(1-\alpha)(l+m)+4\alpha-2$ correspondent à des produits d'un $l$-cycle et d'un $m$-cycle s'intersectant en deux points, et donnant deux cycles disjoints ; en effet, il faut à la fois $|\rho(M)|=l+m-2$ et $|\rho(M)|+\ell(\rho(M))=l+m$. Or, un produit de deux cycles de longueurs $l$ et $m$ s'intersectant en deux points d'écarts $(i,l-i)$ dans le $l$-cycle et $(k,m-k)$ dans le $m$-cycle  est toujours un produit de deux cycles disjoints de longueurs $i+k-1$ et $l-i+m-k-1$ ; \emph{cf.} la figure \ref{prodtwocycles}. 
\figcapt{\psset{unit=0.9mm}\pspicture(0,-10)(140,53)
\rput(0,50){$\textcolor{blue}{1}$}\psline[linecolor=blue]{->}(5,50)(15,50)\rput(20,50){$\textcolor{blue}{\cdots}$}\psline[linecolor=blue]{->}(25,50)(35,50)
\rput(40,50){$\textcolor{blue}{i-1}$}\psline[linecolor=blue]{->}(45,50)(55,50)\rput(60,50){$\textcolor{blue}{i}$}
\psline[linecolor=blue]{->}(60,45)(60,40)\rput(60,35){$\textcolor{blue}{i+1}$}\psline[linecolor=blue]{->}(60,30)(60,25)
\rput(60,20){$\textcolor{blue}{\vdots}$}\psline[linecolor=blue]{->}(60,15)(60,10) \rput(60,5){$\textcolor{blue}{l-1}$} \psline[linecolor=blue]{->}(60,0)(60,-5)  \rput(60,-10){$\textcolor{blue}{l}$} \psline[linecolor=blue]{->}(55,-10)(5,45)
\psline[doubleline=true](65,50)(75,50) 
\psline[doubleline=true](65,-10)(75,-10) 
\rput(80,50){$\textcolor{violet}{i}$}
\psline[linecolor=violet]{->}(80,45)(80,40)\rput(80,35){$\textcolor{violet}{l+1}$}\psline[linecolor=violet]{->}(80,30)(80,25)
\rput(80,20){$\textcolor{violet}{\vdots}$}\psline[linecolor=violet]{->}(80,15)(80,10) \rput(80,5){$\textcolor{violet}{l+k-1}$} \psline[linecolor=violet]{->}(80,0)(80,-5)  \rput(80,-10){$\textcolor{violet}{l}$}
\psline[linecolor=violet]{->}(85,-10)(95,-10) \rput(100,-10){$\textcolor{violet}{l+k}$} \psline[linecolor=violet]{->}(105,-10)(115,-10) \rput(120,-10){$\textcolor{violet}{\cdots}$} \psline[linecolor=violet]{->}(125,-10)(130,-10) \rput(140,-10){$\textcolor{violet}{l+m-2}$} \psline[linecolor=violet]{->}(135,-5)(85,50)
\endpspicture
\bigskip\bigskip
\bigskip

\pspicture(0,-10)(140,50)
\rput(0,50){$\textcolor{red}{1}$}\psline[linecolor=red]{->}(5,50)(15,50)\rput(20,50){$\textcolor{red}{\cdots}$}\psline[linecolor=red]{->}(25,50)(35,50)
\rput(40,50){$\textcolor{red}{i-1}$}\psline[linecolor=red]{->}(45,50)(75,35)\rput(80,35){$\textcolor{red}{l+1}$}\psline[linecolor=red]{->}(80,30)(80,25)
\rput(80,20){$\textcolor{red}{\vdots}$}\psline[linecolor=red]{->}(80,15)(80,10) \rput(80,5){$\textcolor{red}{l+k-1}$} \psline[linecolor=red]{->}(80,0)(65,-10)  \rput(60,-10){$\textcolor{red}{l}$}
\psline[linecolor=red]{->}(55,-10)(5,45)
\rput(80,50){$\textcolor{BurntOrange}{i}$}
\psline[linecolor=BurntOrange,border=1mm,bordercolor=white]{->}(75,50)(60,40)\rput(60,35){$\textcolor{BurntOrange}{i+1}$}\psline[linecolor=BurntOrange]{->}(60,30)(60,25)
\rput(60,20){$\textcolor{BurntOrange}{\vdots}$}\psline[linecolor=BurntOrange]{->}(60,15)(60,10) \rput(60,5){$\textcolor{BurntOrange}{l-1}$} \psline[linecolor=BurntOrange,border=1mm,bordercolor=white]{->}(60,0)(95,-10)  
 \rput(100,-10){$\textcolor{BurntOrange}{l+k}$} \psline[linecolor=BurntOrange]{->}(105,-10)(115,-10) \rput(120,-10){$\textcolor{BurntOrange}{\cdots}$} \psline[linecolor=BurntOrange]{->}(125,-10)(130,-10) \rput(140,-10){$\textcolor{BurntOrange}{l+m-2}$} \psline[linecolor=BurntOrange]{->}(135,-5)(85,50)
\endpspicture
}{Le produit de deux cycles $C_{l}$ (en bleu) et $C_{m}$ (en violet) s'intersectant en deux points d'écarts $(i,l-i)$ dans $C_{l}$ et $(k,m-k)$ dans $C_{m}$ est égal au produit de deux cycles disjoints de longueurs respectives $i+k-1$ (en rouge) et $(l-i)+(m-k)-1$ (en orange).\label{prodtwocycles}}{Produit de deux cycles $C_{l}$ et $C_{m}$ s'intersectant en deux points}

De plus, le nombre de termes de ce type est le nombre d'appariements de taille $2$ entre $\lle 1,l\rre$ et $\lle 1',m'\rre$, c'est-à-dire
$$2\binom{l}{2}\binom{m}{2}=\frac{l(l-1)m(m-1)}{2}\,,$$
et chacun apporte une contribution $n^{(1-\alpha)(l+m)+4\alpha-2}\,c^{l+m-4}$ dans le calcul du cumulant disjoint $k_{\alpha}^{\bullet}(\varSigma_{l},\varSigma_{m})$ (quelque soit la valeur de $i$ et de $k$). On conclut que 
$$\lim_{n \to \infty} k_{\alpha}(Z_{l,n,\alpha,c},Z_{m,n,\alpha,c})=\frac{lm(l-1)(m-1)\,c^{l+m-4}}{2}$$
si $l\geq 2$ et $m\geq 2$. Ainsi, sous les mesures de Schur-Weyl de paramètre $\alpha$ strictement compris entre $1/3$ et $1/2$, les observables 
$$Z_{l,n,\alpha,c}(\lambda)=\frac{\varSigma_{l}(\lambda)-\esper[\varSigma_{l}]}{n^{(1-\alpha)(l-1)+\alpha}}$$
convergent conjointement vers un processus gaussien $(Z_{l,\infty,\alpha,c})_{l\geq 2}$ de matrice de covariance 
$$k(Z_{l,\infty,\alpha,c},Z_{m,\infty,\alpha,c})=2\,\binom{l}{2}\binom{m}{2}\,c^{l+m-4}\,.$$
Comme cette expression est multiplicative en $l$ et en $m$, l'espace gaussien correspondant est de rang $1$, et 
$$\forall l \geq 2,\,\,Z_{l,\infty,\alpha} = \sqrt{2}\,\binom{l}{2}\,c^{l-2}\,X$$
avec $X$ variable gaussienne centrée de variance $1$. Ainsi :

\begin{theorem}[Déviation gaussienne des observables sous les mesures de Schur-Weyl de paramètre $1/3<\alpha<1/2$]\label{devgaussalpha}
Sous les mesures de Schur-Weyl de paramètre $\alpha$ strictement compris entre $1/3$ et $1/2$, les observables 
$$Z_{l,n,\alpha,c}(\lambda)=\frac{\varSigma_{l}(\lambda)-\esper[\varSigma_{l}]}{n^{(1-\alpha)(l-1)+\alpha}}$$
convergent conjointement vers un vecteur gaussien $\mathcal{N}(0,1)\times (\sqrt{2}\,\binom{l}{2}\,c^{l-2})_{l \geq 2}$.
\end{theorem}
\noindent En réalité, on conjecture que le théorème \ref{devgaussalpha} reste vrai pour toute valeur de $\alpha$ comprise entre $0$ et $1/2$ ; en effet, il semble que les cumulants d'ordre $r\geq 2$ des variables $Z_{k,n,\alpha,c}$ soient des $O(n^{(2-r)\alpha}) \ll O(n^{1-\alpha r})$, ce qui impliquerait le résultat directement. Malheureusement, l'outil du $\alpha$-degré n'est pas suffisant pour ce type d'estimation. Notons néanmoins que le calcul fait jusqu'ici montre que les $Z_{k,n,\alpha,c}$ sont asymptotiquement de deuxième moment fini ; ainsi, les estimations
$$ \frac{\varSigma_{k}(\lambda)-\esper[\varSigma_{k}]}{n^{(1-\alpha)k+\alpha}}=O\left(\frac{1}{n^{1-\alpha}}\right)$$
sont exactes, étant entendu que la constante du $O$ est une variable (!) de carré intégrable. On a donc un résultat un peu plus fort que la proposition \ref{notprecisenough}.\bigskip

\begin{remark}
Le théorème \ref{devgaussalpha} ne se traduit pas directement par un résultat sur les déviations des $p_{l}(\lambda)$. En effet, on sait que $p_{l}$ et $\varSigma_{l}$ ont même composante de plus haut $\alpha$-degré (à savoir, $\alpha+(1-\alpha)l$), et la différence est en fait d'$\alpha$-degré $2\alpha+(1-\alpha)(l-1)$. Par exemple, 
$$\varSigma_{6}=p_{6}-6\,p_{4,1}+\frac{145}{4}\,p_{4}-6\,p_{3,2}+18\,p_{2,1,1}-\frac{219}{2}\,p_{2,1}+\frac{1627}{16}\,p_{2}\,,$$
et les termes de plus haut $\alpha$-degré dans $\varSigma_{6}-p_{6}$ sont donc $-6\,p_{4,1}-6\,p_{3,2}$, et ils sont d'$\alpha$-degré $2\alpha+5(1-\alpha)$. L'observable $p_{l}-\varSigma_{l}$ est donc d'ordre de grandeur $n^{2\alpha+(l-1)(1-\alpha)}$, et d'autre part, l'ordre de grandeur de la déviation gaussienne de $\varSigma_{l}$ est $n^{\alpha+l(1-\alpha)}$, qui est strictement inférieur. 
\end{remark}
\bigskip\bigskip
\addtocontents{lof}{\protect\clearpage}

Donnons finalement une interprétation géométrique des résultats de convergence des observables de diagrammes. On dessine le diagramme de Young en redimensionnant les lignes par un facteur $1/n^{1-\alpha}$, et les colonnes par un facteur $1/n^{\alpha}$. Notons $\lambda^{*\alpha}$ cette région du plan --- elle a pour aire $1$. Alors :
\begin{theorem}[Forme limite des diagrammes sous les mesures de Schur-Weyl de paramètre $\alpha<1/2$, \cite{FM10}]\label{schurweylminus}
Sous les mesures de Schur-Weyl de paramètres $\alpha<1/2$ et $c>0$, lorsque $n$ tend vers l'infini, la région $\lambda^{*\alpha}$ converge en probabilité vers le rectangle de longueur $c$ et de hauteur $c^{-1}$. Plus précisément, étant fixées des constantes positives $\eta$ et $\eps$, pour $n$ assez grand, la probabilité pour que le bord de la région $\lambda^{*}$ appartienne à la région hachurée de la figure \ref{zebra} est plus grande que $1-\eps$.
\end{theorem}
\figcapt{\psset{unit=1mm}\pspicture(-2,-2)(100,45)
\pspolygon[fillstyle=hlines](0,25)(0,22)(58,22)(58,0)(90,0)(90,3)(62,3)(62,25)
\psline[linecolor=white,linewidth=2pt](90,0)(90,3)
\psline{->}(-2,0)(95,0)
\psline{->}(0,-2)(0,40)
\rput(105,0){$\times n^{1-\alpha}$}
\rput(0,45){$\times n^{\alpha}$}
\rput(60,-3){$c$}
\rput(-4,25){$c^{-1}$}
\psline{[-]}(30,25)(30,22)
\rput(30,28){$c^{-1}\eta$}
\psline{[-]}(58,13)(62,13)
\rput(66,13){$c\eta$}
\psline(-1,25)(0,25)
\psline(60,-1)(60,1)
\endpspicture}{Forme quasi-rectangle typique des diagrammes de Young renormalisés sous les mesures de Schur-Weyl de paramètres $\alpha<1/2$ et $c>0$.\label{zebra}}{Formes limites pour les mesures de Schur-Weyl de paramètre $\alpha<1/2$, I}
\begin{proof}
Comme dans le chapitre \ref{qplancherelmeasure}, on associe à toute partition $\lambda$ une mesure de probabilité sur $\R$ :
$$X_{\lambda,\alpha}=\sum_{i=1}^{d} \frac{a_{i}(\lambda)}{n}\,\delta_{a_{i}(\lambda)/n^{1-\alpha}}+\sum_{i=1}^{d}\frac{b_{i}(\lambda)}{n}\,\delta_{-b_{i}(\lambda)/n^{1-\alpha}}\,.$$
Par construction, le $k$-ième moment de $X_{\lambda,\alpha}$ est $p_{k+1}(\lambda)/n^{(1-\alpha)(k+1)+\alpha}$, donc converge vers $c^{k}$. Les moments de $X_{\lambda,\alpha}$ convergent donc vers ceux de $\delta_{c}$, qui est une mesure de probabilité sur $\R$ caractérisée par ces moments ; par conséquent, $X_{\lambda,\alpha}$ converge pour la topologie faible vers $\delta_{c}$ (en probabilité). Notons que l'on sait déjà que la hauteur $\ell(\lambda)$ est majorée avec probabilité $1$ par $N=n^{\alpha}/c$. Par conséquent, le poids de la partie négative 
$$\sum_{i=1}^{d} \frac{b_{i}(\lambda)}{n}\,\delta_{-b_{i}(\lambda)/n^{1-\alpha}}$$
de la mesure (aléatoire) $X_{\lambda,\alpha}$ est majoré avec probabilité $1$ par 
$$M=\sum_{i=1}^{N} \frac{N}{n}=\frac{N^{2}}{n}\simeq \frac{n^{2\alpha-1}}{c^{2}} \to 0\,.$$
On peut donc remplacer $X_{\lambda,\alpha}$ par la mesure $Y_{\lambda,\alpha}=\sum_{i=1}^{d} \frac{a_{i}(\lambda)}{n}\,\delta_{a_{i}(\lambda)/n^{1-\alpha}}$ : son poids converge vers $1$ et elle converge en probabilité vers $\delta_{c}$ (au sens de la convergence en loi des mesures de Radon). Finalement, on peut remplacer les $a_{i}(\lambda)$ par les $\lambda_{i}$. En effet, fixons une fonction $f$ bornée et uniformément continue sur $\R$. La différence $\frac{a_{i}(\lambda)-\lambda_{i}}{n^{1-\alpha}}$ est bornée par $\frac{N}{n^{1-\alpha}}=\frac{n^{2\alpha-1}}{c}$, donc pour $n$ assez grand, 
$$\left|f\left(\frac{\lambda_{i}}{n^{1-\alpha}}\right)-f\left(\frac{a_{i}(\lambda)}{n^{1-\alpha}}\right)\right|\leq \eps$$
pour tout $i$, et ce par uniforme continuité de $f$ (l'inégalité est vraie presque sûrement, et pas seulement avec grande probabilité). Alors, notant $Z_{\lambda,\alpha}=\sum_{i=1}^{N}\frac{\lambda_{i}}{n}\,\delta_{\lambda_{i}/n^{1-\alpha}}$, on a :
\begin{align*} \left|Y_{\lambda,\alpha}(f)-Z_{\lambda,\alpha}(f)\right|&\leq \left|\sum_{i=1}^{d}\left(\frac{a_{i}(\lambda)}{n}-\frac{\lambda_{i}}{n}\right)f\left(\frac{a_{i}(\lambda)}{n^{1-\alpha}}\right)\right|+\left|\sum_{i=1}^{d} \frac{\lambda_{i}}{n}\left(f\left(\frac{a_{i}(\lambda)}{n^{1-\alpha}}\right)-f\left(\frac{\lambda_{i}}{n^{1-\alpha}}\right)\right)\right|\\
&\quad+\left|\sum_{i=d+1}^{N}\frac{\lambda_{i}}{n}f\left(\frac{\lambda_{i}}{n^{1-\alpha}}\right)\right|\\
&\leq \frac{N^{2}}{n}\|f\|_{\infty}+\eps+\left(\sum_{i=d+1}^{N}\frac{\lambda_{i}}{n}\right)\|f\|_{\infty}\,.
\end{align*}
De plus, $\sum_{i=d+1}^{N}\frac{\lambda_{i}}{n}\leq |Y_{\lambda,\alpha}(1)-Z_{\lambda,\alpha}(1)|+\frac{N^{2}}{n}$, qui tend vers $0$ car $Z_{\lambda,\alpha}$ est une mesure de probabilité, et le poids de $Y_{\lambda,\alpha}$ tend vers $1$. Ainsi, on a montré que pour toute fonction $f$ bornée uniformément continue sur $\R$, $|Y_{\lambda,\alpha}(f)-Z_{\lambda,\alpha}(f)|$ tend vers $0$ (avec probabilité $1$), et par suite, $Z_{\lambda,\alpha}(f)$ converge en probabilité vers $\delta_{c}(f)=f(c)$. Or, les fonctions uniformément continues suffisent pour tester la convergence en loi des mesures de probabilité, voir \cite[chapitre 1, théorème 2.1]{Bil69} ; ainsi, la mesure $Z_{\lambda,\alpha}$ converge en probabilité vers $\delta_{c}$ pour la topologie faible sur les mesures de probabilité.\bigskip

Fixons maintenant des réels positifs $\theta$, $\eta$ et $\eps$, et considérons une fonction $f$ continue, qui vaut $0$ en dehors de $[c(1-\eta/2),c(1+\eta/2)]$ et $1$ en $c$. Comme $f(c)=1$, pour $n$ assez grand, avec probabilité plus grande que $1-\eps$,
$$1-\theta \leq \sum_{i=1}^{N} \frac{\lambda_{i}}{n}f\left(\frac{\lambda_{i}}{n^{1-\alpha}}\right) \leq 1+\theta\,.$$
La somme au milieu est majorée par $\sum_{\left\{i\,\,|\,\,\lambda_{i}/n^{1-\alpha} \in [c(1-\eta/2),c(1+\eta/2)]\right\}}\frac{\lambda_{i}}{n}$, et donc par $$\frac{c(1+\eta/2)}{n^{\alpha}}\,\,\card\left\{i\,\,\bigg|\,\,\frac{\lambda_{i}}{n^{1-\alpha}} \in [c(1-\eta/2),c(1+\eta/2)]\right\}\,.$$ Ainsi, le nombre de lignes $\lambda_{i}$ telles que $c(1-\eta/2)\leq \lambda_{i}/n^{1-\eta} \leq c(1+\eta/2)$ est plus grand que
$$\frac{1-\theta}{c(1+\eta/2)} \,n^{\alpha}\geq \frac{(1-\theta)(1-\eta/2)}{c}\,n^{\alpha}\geq \frac{1-\eta}{c}\,n^{\alpha}$$
pour $\theta$ assez petit, et ce avec probabilité plus grande que $1-\eps$. C'est exactement ce qui est représenté sur la figure \ref{zebra}.
\end{proof}
\bigskip

En réalité, on peut raffiner le résultat donné par le théorème \ref{schurweylminus} et ôter la partie inférieure droite de la zone hachurée. Ainsi, pour tous $\eps,\eta>0$ et tout $n$ assez grand, les diagrammes choisis sous les mesures de Schur-Weyl $SW_{n,\alpha,c}$ sont après renormalisation compris dans la zone hachurée de la figure suivante ave probabilité plus grande que $1-\eps$ :

\figcapt{\psset{unit=1mm}\pspicture(-2,-2)(100,45)
\pspolygon[fillstyle=hlines](0,25)(0,22)(58,22)(58,0)(62,0)(62,3)(62,25)
\psline{->}(-2,0)(95,0)
\psline{->}(0,-2)(0,40)
\rput(105,0){$\times n^{1-\alpha}$}
\rput(0,45){$\times n^{\alpha}$}
\rput(60,-3){$c$}
\rput(-4,25){$c^{-1}$}
\psline{[-]}(30,25)(30,22)
\rput(30,28){$c^{-1}\eta$}
\psline{[-]}(58,13)(62,13)
\rput(66,13){$c\eta$}
\psline(-1,25)(0,25)
\psline(60,-1)(60,1)
\endpspicture}{Forme rectangle typique des diagrammes de Young renormalisés sous les mesures de Schur-Weyl de paramètre $\alpha<1/2$ et $c>0$.\label{zebra2}}{Formes limites pour les mesures de Schur-Weyl de paramètre $\alpha<1/2$, II}
\noindent En effet, K. Johansson a montré dans \cite[théorème 1.7]{Joh01} que pour toutes les mesures de Schur-Weyl $SW_{n,\alpha,c}$, la longueur de la première ligne d'une partition pouvait s'écrire sous la forme 
$$\lambda_{1}=c\,n^{1-\alpha}+2\sqrt{n}+(1+c\,n^{1/2-\alpha})^{2/3}\,n^{1/6}\,X_{1}$$
où $X_{1}$ suit asymptotiquement une loi de Tracy-Widom (la même que dans le chapitre \ref{matrix}). Ainsi, lorsque $\alpha<1/2$, $\lambda_{1}=c\,n^{1-\alpha}+2\sqrt{n}+c\,\sqrt{n}\,n^{-2\alpha/3}\,X_{1}\simeq c\,n^{1-\alpha}$, d'où le résultat raffiné de la figure \ref{zebra2}.

\chapter{Asymptotique des mesures de Gelfand}\label{gelfandmeasure}

L'un des résultats marquants du chapitre précédent était la validité (à une translation près) du théorème central limite de Kerov dans le cadre des mesures de Schur-Weyl de paramètre $\alpha=1/2$ (\emph{cf.} le théorème \ref{schurweylplus}). Partant, il est naturel de se demander s'il existe d'autres mesures sur les partitions qui présentent ce même comportement asymptotique. Dans ce chapitre, nous reprenons les résultats de l'article \cite{Mel10} et nous étudions une dernière famille intéressante de mesures sur les partitions associées à des représentations des groupes symétriques : les \textbf{mesures de Gelfand}. Il s'agit des mesures de probabilité définies par :
$$\Gel_{n}[\lambda \in \ym_{n}]=\frac{\dim \lambda}{\sum_{\mu \in \ym_{n}} \dim \mu}\,.$$
La fonction de partition $\sum_{\mu \in \ym_{n}} \dim \mu$ peut être interprétée combinatoirement comme le nom\-bre d'involutions de taille $n$ ; nous expliquerons ceci dans la section \ref{gelfandmodel}, et nous présenterons des représentations explicites des groupes $\sym_{n}$ dont les mesures de Plancherel (au sens de la définition \ref{defplancherel}) sont ces mesures de Gelfand.  La figure ci-dessous représente un diagramme de Young tiré sous la mesure de Plancherel de paramètre $n=400$ (à gauche), et un diagramme de Young tiré sous la mesure de Gelfand de paramètre $n=400$ (à droite) :
\figcapt{
\psset{unit=0.5mm}
\pspicture(-75,0)(230,85)
\parametricplot{-2}{2}{t 40 mul 
t t 0.5 mul arcsin mul 0.0174533 mul t t mul neg 4 add sqrt add 0.63661977 mul 40 mul }
\psline{->}(0,0)(-82,82)\psline{->}(0,0)(82,82)
\psline(0,0)(-70.000,70.000)\psline(0,0)(66.000,66.000)
\psline(2.0000,2.0000)(0.00000,4.0000)(-2.0000,2.0000)
\psline(4.0000,4.0000)(2.0000,6.0000)(0.00000,4.0000)
\psline(6.0000,6.0000)(4.0000,8.0000)(2.0000,6.0000)
\psline(8.0000,8.0000)(6.0000,10.000)(4.0000,8.0000)
\psline(10.000,10.000)(8.0000,12.000)(6.0000,10.000)
\psline(12.000,12.000)(10.000,14.000)(8.0000,12.000)
\psline(14.000,14.000)(12.000,16.000)(10.000,14.000)
\psline(16.000,16.000)(14.000,18.000)(12.000,16.000)
\psline(18.000,18.000)(16.000,20.000)(14.000,18.000)
\psline(20.000,20.000)(18.000,22.000)(16.000,20.000)
\psline(22.000,22.000)(20.000,24.000)(18.000,22.000)
\psline(24.000,24.000)(22.000,26.000)(20.000,24.000)
\psline(26.000,26.000)(24.000,28.000)(22.000,26.000)
\psline(28.000,28.000)(26.000,30.000)(24.000,28.000)
\psline(30.000,30.000)(28.000,32.000)(26.000,30.000)
\psline(32.000,32.000)(30.000,34.000)(28.000,32.000)
\psline(34.000,34.000)(32.000,36.000)(30.000,34.000)
\psline(36.000,36.000)(34.000,38.000)(32.000,36.000)
\psline(38.000,38.000)(36.000,40.000)(34.000,38.000)
\psline(40.000,40.000)(38.000,42.000)(36.000,40.000)
\psline(42.000,42.000)(40.000,44.000)(38.000,42.000)
\psline(44.000,44.000)(42.000,46.000)(40.000,44.000)
\psline(46.000,46.000)(44.000,48.000)(42.000,46.000)
\psline(48.000,48.000)(46.000,50.000)(44.000,48.000)
\psline(50.000,50.000)(48.000,52.000)(46.000,50.000)
\psline(52.000,52.000)(50.000,54.000)(48.000,52.000)
\psline(54.000,54.000)(52.000,56.000)(50.000,54.000)
\psline(56.000,56.000)(54.000,58.000)(52.000,56.000)
\psline(58.000,58.000)(56.000,60.000)(54.000,58.000)
\psline(60.000,60.000)(58.000,62.000)(56.000,60.000)
\psline(62.000,62.000)(60.000,64.000)(58.000,62.000)
\psline(64.000,64.000)(62.000,66.000)(60.000,64.000)
\psline(66.000,66.000)(64.000,68.000)(62.000,66.000)
\psline(0.00000,4.0000)(-2.0000,6.0000)(-4.0000,4.0000)
\psline(2.0000,6.0000)(0.00000,8.0000)(-2.0000,6.0000)
\psline(4.0000,8.0000)(2.0000,10.000)(0.00000,8.0000)
\psline(6.0000,10.000)(4.0000,12.000)(2.0000,10.000)
\psline(8.0000,12.000)(6.0000,14.000)(4.0000,12.000)
\psline(10.000,14.000)(8.0000,16.000)(6.0000,14.000)
\psline(12.000,16.000)(10.000,18.000)(8.0000,16.000)
\psline(14.000,18.000)(12.000,20.000)(10.000,18.000)
\psline(16.000,20.000)(14.000,22.000)(12.000,20.000)
\psline(18.000,22.000)(16.000,24.000)(14.000,22.000)
\psline(20.000,24.000)(18.000,26.000)(16.000,24.000)
\psline(22.000,26.000)(20.000,28.000)(18.000,26.000)
\psline(24.000,28.000)(22.000,30.000)(20.000,28.000)
\psline(26.000,30.000)(24.000,32.000)(22.000,30.000)
\psline(28.000,32.000)(26.000,34.000)(24.000,32.000)
\psline(30.000,34.000)(28.000,36.000)(26.000,34.000)
\psline(32.000,36.000)(30.000,38.000)(28.000,36.000)
\psline(34.000,38.000)(32.000,40.000)(30.000,38.000)
\psline(36.000,40.000)(34.000,42.000)(32.000,40.000)
\psline(38.000,42.000)(36.000,44.000)(34.000,42.000)
\psline(40.000,44.000)(38.000,46.000)(36.000,44.000)
\psline(42.000,46.000)(40.000,48.000)(38.000,46.000)
\psline(44.000,48.000)(42.000,50.000)(40.000,48.000)
\psline(46.000,50.000)(44.000,52.000)(42.000,50.000)
\psline(48.000,52.000)(46.000,54.000)(44.000,52.000)
\psline(50.000,54.000)(48.000,56.000)(46.000,54.000)
\psline(52.000,56.000)(50.000,58.000)(48.000,56.000)
\psline(54.000,58.000)(52.000,60.000)(50.000,58.000)
\psline(56.000,60.000)(54.000,62.000)(52.000,60.000)
\psline(58.000,62.000)(56.000,64.000)(54.000,62.000)
\psline(60.000,64.000)(58.000,66.000)(56.000,64.000)
\psline(62.000,66.000)(60.000,68.000)(58.000,66.000)
\psline(-2.0000,6.0000)(-4.0000,8.0000)(-6.0000,6.0000)
\psline(0.00000,8.0000)(-2.0000,10.000)(-4.0000,8.0000)
\psline(2.0000,10.000)(0.00000,12.000)(-2.0000,10.000)
\psline(4.0000,12.000)(2.0000,14.000)(0.00000,12.000)
\psline(6.0000,14.000)(4.0000,16.000)(2.0000,14.000)
\psline(8.0000,16.000)(6.0000,18.000)(4.0000,16.000)
\psline(10.000,18.000)(8.0000,20.000)(6.0000,18.000)
\psline(12.000,20.000)(10.000,22.000)(8.0000,20.000)
\psline(14.000,22.000)(12.000,24.000)(10.000,22.000)
\psline(16.000,24.000)(14.000,26.000)(12.000,24.000)
\psline(18.000,26.000)(16.000,28.000)(14.000,26.000)
\psline(20.000,28.000)(18.000,30.000)(16.000,28.000)
\psline(22.000,30.000)(20.000,32.000)(18.000,30.000)
\psline(24.000,32.000)(22.000,34.000)(20.000,32.000)
\psline(26.000,34.000)(24.000,36.000)(22.000,34.000)
\psline(28.000,36.000)(26.000,38.000)(24.000,36.000)
\psline(30.000,38.000)(28.000,40.000)(26.000,38.000)
\psline(32.000,40.000)(30.000,42.000)(28.000,40.000)
\psline(34.000,42.000)(32.000,44.000)(30.000,42.000)
\psline(36.000,44.000)(34.000,46.000)(32.000,44.000)
\psline(38.000,46.000)(36.000,48.000)(34.000,46.000)
\psline(40.000,48.000)(38.000,50.000)(36.000,48.000)
\psline(42.000,50.000)(40.000,52.000)(38.000,50.000)
\psline(44.000,52.000)(42.000,54.000)(40.000,52.000)
\psline(46.000,54.000)(44.000,56.000)(42.000,54.000)
\psline(48.000,56.000)(46.000,58.000)(44.000,56.000)
\psline(50.000,58.000)(48.000,60.000)(46.000,58.000)
\psline(-4.0000,8.0000)(-6.0000,10.000)(-8.0000,8.0000)
\psline(-2.0000,10.000)(-4.0000,12.000)(-6.0000,10.000)
\psline(0.00000,12.000)(-2.0000,14.000)(-4.0000,12.000)
\psline(2.0000,14.000)(0.00000,16.000)(-2.0000,14.000)
\psline(4.0000,16.000)(2.0000,18.000)(0.00000,16.000)
\psline(6.0000,18.000)(4.0000,20.000)(2.0000,18.000)
\psline(8.0000,20.000)(6.0000,22.000)(4.0000,20.000)
\psline(10.000,22.000)(8.0000,24.000)(6.0000,22.000)
\psline(12.000,24.000)(10.000,26.000)(8.0000,24.000)
\psline(14.000,26.000)(12.000,28.000)(10.000,26.000)
\psline(16.000,28.000)(14.000,30.000)(12.000,28.000)
\psline(18.000,30.000)(16.000,32.000)(14.000,30.000)
\psline(20.000,32.000)(18.000,34.000)(16.000,32.000)
\psline(22.000,34.000)(20.000,36.000)(18.000,34.000)
\psline(24.000,36.000)(22.000,38.000)(20.000,36.000)
\psline(26.000,38.000)(24.000,40.000)(22.000,38.000)
\psline(28.000,40.000)(26.000,42.000)(24.000,40.000)
\psline(30.000,42.000)(28.000,44.000)(26.000,42.000)
\psline(32.000,44.000)(30.000,46.000)(28.000,44.000)
\psline(34.000,46.000)(32.000,48.000)(30.000,46.000)
\psline(36.000,48.000)(34.000,50.000)(32.000,48.000)
\psline(38.000,50.000)(36.000,52.000)(34.000,50.000)
\psline(40.000,52.000)(38.000,54.000)(36.000,52.000)
\psline(42.000,54.000)(40.000,56.000)(38.000,54.000)
\psline(44.000,56.000)(42.000,58.000)(40.000,56.000)
\psline(-6.0000,10.000)(-8.0000,12.000)(-10.000,10.000)
\psline(-4.0000,12.000)(-6.0000,14.000)(-8.0000,12.000)
\psline(-2.0000,14.000)(-4.0000,16.000)(-6.0000,14.000)
\psline(0.00000,16.000)(-2.0000,18.000)(-4.0000,16.000)
\psline(2.0000,18.000)(0.00000,20.000)(-2.0000,18.000)
\psline(4.0000,20.000)(2.0000,22.000)(0.00000,20.000)
\psline(6.0000,22.000)(4.0000,24.000)(2.0000,22.000)
\psline(8.0000,24.000)(6.0000,26.000)(4.0000,24.000)
\psline(10.000,26.000)(8.0000,28.000)(6.0000,26.000)
\psline(12.000,28.000)(10.000,30.000)(8.0000,28.000)
\psline(14.000,30.000)(12.000,32.000)(10.000,30.000)
\psline(16.000,32.000)(14.000,34.000)(12.000,32.000)
\psline(18.000,34.000)(16.000,36.000)(14.000,34.000)
\psline(20.000,36.000)(18.000,38.000)(16.000,36.000)
\psline(22.000,38.000)(20.000,40.000)(18.000,38.000)
\psline(24.000,40.000)(22.000,42.000)(20.000,40.000)
\psline(26.000,42.000)(24.000,44.000)(22.000,42.000)
\psline(28.000,44.000)(26.000,46.000)(24.000,44.000)
\psline(30.000,46.000)(28.000,48.000)(26.000,46.000)
\psline(32.000,48.000)(30.000,50.000)(28.000,48.000)
\psline(34.000,50.000)(32.000,52.000)(30.000,50.000)
\psline(36.000,52.000)(34.000,54.000)(32.000,52.000)
\psline(38.000,54.000)(36.000,56.000)(34.000,54.000)
\psline(40.000,56.000)(38.000,58.000)(36.000,56.000)
\psline(-8.0000,12.000)(-10.000,14.000)(-12.000,12.000)
\psline(-6.0000,14.000)(-8.0000,16.000)(-10.000,14.000)
\psline(-4.0000,16.000)(-6.0000,18.000)(-8.0000,16.000)
\psline(-2.0000,18.000)(-4.0000,20.000)(-6.0000,18.000)
\psline(0.00000,20.000)(-2.0000,22.000)(-4.0000,20.000)
\psline(2.0000,22.000)(0.00000,24.000)(-2.0000,22.000)
\psline(4.0000,24.000)(2.0000,26.000)(0.00000,24.000)
\psline(6.0000,26.000)(4.0000,28.000)(2.0000,26.000)
\psline(8.0000,28.000)(6.0000,30.000)(4.0000,28.000)
\psline(10.000,30.000)(8.0000,32.000)(6.0000,30.000)
\psline(12.000,32.000)(10.000,34.000)(8.0000,32.000)
\psline(14.000,34.000)(12.000,36.000)(10.000,34.000)
\psline(16.000,36.000)(14.000,38.000)(12.000,36.000)
\psline(18.000,38.000)(16.000,40.000)(14.000,38.000)
\psline(20.000,40.000)(18.000,42.000)(16.000,40.000)
\psline(22.000,42.000)(20.000,44.000)(18.000,42.000)
\psline(24.000,44.000)(22.000,46.000)(20.000,44.000)
\psline(26.000,46.000)(24.000,48.000)(22.000,46.000)
\psline(28.000,48.000)(26.000,50.000)(24.000,48.000)
\psline(30.000,50.000)(28.000,52.000)(26.000,50.000)
\psline(32.000,52.000)(30.000,54.000)(28.000,52.000)
\psline(34.000,54.000)(32.000,56.000)(30.000,54.000)
\psline(36.000,56.000)(34.000,58.000)(32.000,56.000)
\psline(-10.000,14.000)(-12.000,16.000)(-14.000,14.000)
\psline(-8.0000,16.000)(-10.000,18.000)(-12.000,16.000)
\psline(-6.0000,18.000)(-8.0000,20.000)(-10.000,18.000)
\psline(-4.0000,20.000)(-6.0000,22.000)(-8.0000,20.000)
\psline(-2.0000,22.000)(-4.0000,24.000)(-6.0000,22.000)
\psline(0.00000,24.000)(-2.0000,26.000)(-4.0000,24.000)
\psline(2.0000,26.000)(0.00000,28.000)(-2.0000,26.000)
\psline(4.0000,28.000)(2.0000,30.000)(0.00000,28.000)
\psline(6.0000,30.000)(4.0000,32.000)(2.0000,30.000)
\psline(8.0000,32.000)(6.0000,34.000)(4.0000,32.000)
\psline(10.000,34.000)(8.0000,36.000)(6.0000,34.000)
\psline(12.000,36.000)(10.000,38.000)(8.0000,36.000)
\psline(14.000,38.000)(12.000,40.000)(10.000,38.000)
\psline(16.000,40.000)(14.000,42.000)(12.000,40.000)
\psline(18.000,42.000)(16.000,44.000)(14.000,42.000)
\psline(20.000,44.000)(18.000,46.000)(16.000,44.000)
\psline(22.000,46.000)(20.000,48.000)(18.000,46.000)
\psline(24.000,48.000)(22.000,50.000)(20.000,48.000)
\psline(26.000,50.000)(24.000,52.000)(22.000,50.000)
\psline(28.000,52.000)(26.000,54.000)(24.000,52.000)
\psline(30.000,54.000)(28.000,56.000)(26.000,54.000)
\psline(32.000,56.000)(30.000,58.000)(28.000,56.000)
\psline(-12.000,16.000)(-14.000,18.000)(-16.000,16.000)
\psline(-10.000,18.000)(-12.000,20.000)(-14.000,18.000)
\psline(-8.0000,20.000)(-10.000,22.000)(-12.000,20.000)
\psline(-6.0000,22.000)(-8.0000,24.000)(-10.000,22.000)
\psline(-4.0000,24.000)(-6.0000,26.000)(-8.0000,24.000)
\psline(-2.0000,26.000)(-4.0000,28.000)(-6.0000,26.000)
\psline(0.00000,28.000)(-2.0000,30.000)(-4.0000,28.000)
\psline(2.0000,30.000)(0.00000,32.000)(-2.0000,30.000)
\psline(4.0000,32.000)(2.0000,34.000)(0.00000,32.000)
\psline(6.0000,34.000)(4.0000,36.000)(2.0000,34.000)
\psline(8.0000,36.000)(6.0000,38.000)(4.0000,36.000)
\psline(10.000,38.000)(8.0000,40.000)(6.0000,38.000)
\psline(12.000,40.000)(10.000,42.000)(8.0000,40.000)
\psline(14.000,42.000)(12.000,44.000)(10.000,42.000)
\psline(16.000,44.000)(14.000,46.000)(12.000,44.000)
\psline(18.000,46.000)(16.000,48.000)(14.000,46.000)
\psline(20.000,48.000)(18.000,50.000)(16.000,48.000)
\psline(22.000,50.000)(20.000,52.000)(18.000,50.000)
\psline(24.000,52.000)(22.000,54.000)(20.000,52.000)
\psline(26.000,54.000)(24.000,56.000)(22.000,54.000)
\psline(-14.000,18.000)(-16.000,20.000)(-18.000,18.000)
\psline(-12.000,20.000)(-14.000,22.000)(-16.000,20.000)
\psline(-10.000,22.000)(-12.000,24.000)(-14.000,22.000)
\psline(-8.0000,24.000)(-10.000,26.000)(-12.000,24.000)
\psline(-6.0000,26.000)(-8.0000,28.000)(-10.000,26.000)
\psline(-4.0000,28.000)(-6.0000,30.000)(-8.0000,28.000)
\psline(-2.0000,30.000)(-4.0000,32.000)(-6.0000,30.000)
\psline(0.00000,32.000)(-2.0000,34.000)(-4.0000,32.000)
\psline(2.0000,34.000)(0.00000,36.000)(-2.0000,34.000)
\psline(4.0000,36.000)(2.0000,38.000)(0.00000,36.000)
\psline(6.0000,38.000)(4.0000,40.000)(2.0000,38.000)
\psline(8.0000,40.000)(6.0000,42.000)(4.0000,40.000)
\psline(10.000,42.000)(8.0000,44.000)(6.0000,42.000)
\psline(12.000,44.000)(10.000,46.000)(8.0000,44.000)
\psline(14.000,46.000)(12.000,48.000)(10.000,46.000)
\psline(16.000,48.000)(14.000,50.000)(12.000,48.000)
\psline(18.000,50.000)(16.000,52.000)(14.000,50.000)
\psline(20.000,52.000)(18.000,54.000)(16.000,52.000)
\psline(-16.000,20.000)(-18.000,22.000)(-20.000,20.000)
\psline(-14.000,22.000)(-16.000,24.000)(-18.000,22.000)
\psline(-12.000,24.000)(-14.000,26.000)(-16.000,24.000)
\psline(-10.000,26.000)(-12.000,28.000)(-14.000,26.000)
\psline(-8.0000,28.000)(-10.000,30.000)(-12.000,28.000)
\psline(-6.0000,30.000)(-8.0000,32.000)(-10.000,30.000)
\psline(-4.0000,32.000)(-6.0000,34.000)(-8.0000,32.000)
\psline(-2.0000,34.000)(-4.0000,36.000)(-6.0000,34.000)
\psline(0.00000,36.000)(-2.0000,38.000)(-4.0000,36.000)
\psline(2.0000,38.000)(0.00000,40.000)(-2.0000,38.000)
\psline(4.0000,40.000)(2.0000,42.000)(0.00000,40.000)
\psline(6.0000,42.000)(4.0000,44.000)(2.0000,42.000)
\psline(8.0000,44.000)(6.0000,46.000)(4.0000,44.000)
\psline(10.000,46.000)(8.0000,48.000)(6.0000,46.000)
\psline(12.000,48.000)(10.000,50.000)(8.0000,48.000)
\psline(14.000,50.000)(12.000,52.000)(10.000,50.000)
\psline(-18.000,22.000)(-20.000,24.000)(-22.000,22.000)
\psline(-16.000,24.000)(-18.000,26.000)(-20.000,24.000)
\psline(-14.000,26.000)(-16.000,28.000)(-18.000,26.000)
\psline(-12.000,28.000)(-14.000,30.000)(-16.000,28.000)
\psline(-10.000,30.000)(-12.000,32.000)(-14.000,30.000)
\psline(-8.0000,32.000)(-10.000,34.000)(-12.000,32.000)
\psline(-6.0000,34.000)(-8.0000,36.000)(-10.000,34.000)
\psline(-4.0000,36.000)(-6.0000,38.000)(-8.0000,36.000)
\psline(-2.0000,38.000)(-4.0000,40.000)(-6.0000,38.000)
\psline(0.00000,40.000)(-2.0000,42.000)(-4.0000,40.000)
\psline(2.0000,42.000)(0.00000,44.000)(-2.0000,42.000)
\psline(4.0000,44.000)(2.0000,46.000)(0.00000,44.000)
\psline(6.0000,46.000)(4.0000,48.000)(2.0000,46.000)
\psline(8.0000,48.000)(6.0000,50.000)(4.0000,48.000)
\psline(10.000,50.000)(8.0000,52.000)(6.0000,50.000)
\psline(-20.000,24.000)(-22.000,26.000)(-24.000,24.000)
\psline(-18.000,26.000)(-20.000,28.000)(-22.000,26.000)
\psline(-16.000,28.000)(-18.000,30.000)(-20.000,28.000)
\psline(-14.000,30.000)(-16.000,32.000)(-18.000,30.000)
\psline(-12.000,32.000)(-14.000,34.000)(-16.000,32.000)
\psline(-10.000,34.000)(-12.000,36.000)(-14.000,34.000)
\psline(-8.0000,36.000)(-10.000,38.000)(-12.000,36.000)
\psline(-6.0000,38.000)(-8.0000,40.000)(-10.000,38.000)
\psline(-4.0000,40.000)(-6.0000,42.000)(-8.0000,40.000)
\psline(-2.0000,42.000)(-4.0000,44.000)(-6.0000,42.000)
\psline(0.00000,44.000)(-2.0000,46.000)(-4.0000,44.000)
\psline(2.0000,46.000)(0.00000,48.000)(-2.0000,46.000)
\psline(4.0000,48.000)(2.0000,50.000)(0.00000,48.000)
\psline(6.0000,50.000)(4.0000,52.000)(2.0000,50.000)
\psline(8.0000,52.000)(6.0000,54.000)(4.0000,52.000)
\psline(-22.000,26.000)(-24.000,28.000)(-26.000,26.000)
\psline(-20.000,28.000)(-22.000,30.000)(-24.000,28.000)
\psline(-18.000,30.000)(-20.000,32.000)(-22.000,30.000)
\psline(-16.000,32.000)(-18.000,34.000)(-20.000,32.000)
\psline(-14.000,34.000)(-16.000,36.000)(-18.000,34.000)
\psline(-12.000,36.000)(-14.000,38.000)(-16.000,36.000)
\psline(-10.000,38.000)(-12.000,40.000)(-14.000,38.000)
\psline(-8.0000,40.000)(-10.000,42.000)(-12.000,40.000)
\psline(-6.0000,42.000)(-8.0000,44.000)(-10.000,42.000)
\psline(-4.0000,44.000)(-6.0000,46.000)(-8.0000,44.000)
\psline(-2.0000,46.000)(-4.0000,48.000)(-6.0000,46.000)
\psline(0.00000,48.000)(-2.0000,50.000)(-4.0000,48.000)
\psline(2.0000,50.000)(0.00000,52.000)(-2.0000,50.000)
\psline(4.0000,52.000)(2.0000,54.000)(0.00000,52.000)
\psline(-24.000,28.000)(-26.000,30.000)(-28.000,28.000)
\psline(-22.000,30.000)(-24.000,32.000)(-26.000,30.000)
\psline(-20.000,32.000)(-22.000,34.000)(-24.000,32.000)
\psline(-18.000,34.000)(-20.000,36.000)(-22.000,34.000)
\psline(-16.000,36.000)(-18.000,38.000)(-20.000,36.000)
\psline(-14.000,38.000)(-16.000,40.000)(-18.000,38.000)
\psline(-12.000,40.000)(-14.000,42.000)(-16.000,40.000)
\psline(-10.000,42.000)(-12.000,44.000)(-14.000,42.000)
\psline(-8.0000,44.000)(-10.000,46.000)(-12.000,44.000)
\psline(-6.0000,46.000)(-8.0000,48.000)(-10.000,46.000)
\psline(-4.0000,48.000)(-6.0000,50.000)(-8.0000,48.000)
\psline(-2.0000,50.000)(-4.0000,52.000)(-6.0000,50.000)
\psline(-26.000,30.000)(-28.000,32.000)(-30.000,30.000)
\psline(-24.000,32.000)(-26.000,34.000)(-28.000,32.000)
\psline(-22.000,34.000)(-24.000,36.000)(-26.000,34.000)
\psline(-20.000,36.000)(-22.000,38.000)(-24.000,36.000)
\psline(-18.000,38.000)(-20.000,40.000)(-22.000,38.000)
\psline(-16.000,40.000)(-18.000,42.000)(-20.000,40.000)
\psline(-14.000,42.000)(-16.000,44.000)(-18.000,42.000)
\psline(-12.000,44.000)(-14.000,46.000)(-16.000,44.000)
\psline(-10.000,46.000)(-12.000,48.000)(-14.000,46.000)
\psline(-8.0000,48.000)(-10.000,50.000)(-12.000,48.000)
\psline(-6.0000,50.000)(-8.0000,52.000)(-10.000,50.000)
\psline(-28.000,32.000)(-30.000,34.000)(-32.000,32.000)
\psline(-26.000,34.000)(-28.000,36.000)(-30.000,34.000)
\psline(-24.000,36.000)(-26.000,38.000)(-28.000,36.000)
\psline(-22.000,38.000)(-24.000,40.000)(-26.000,38.000)
\psline(-20.000,40.000)(-22.000,42.000)(-24.000,40.000)
\psline(-18.000,42.000)(-20.000,44.000)(-22.000,42.000)
\psline(-16.000,44.000)(-18.000,46.000)(-20.000,44.000)
\psline(-14.000,46.000)(-16.000,48.000)(-18.000,46.000)
\psline(-12.000,48.000)(-14.000,50.000)(-16.000,48.000)
\psline(-10.000,50.000)(-12.000,52.000)(-14.000,50.000)
\psline(-30.000,34.000)(-32.000,36.000)(-34.000,34.000)
\psline(-28.000,36.000)(-30.000,38.000)(-32.000,36.000)
\psline(-26.000,38.000)(-28.000,40.000)(-30.000,38.000)
\psline(-24.000,40.000)(-26.000,42.000)(-28.000,40.000)
\psline(-22.000,42.000)(-24.000,44.000)(-26.000,42.000)
\psline(-20.000,44.000)(-22.000,46.000)(-24.000,44.000)
\psline(-18.000,46.000)(-20.000,48.000)(-22.000,46.000)
\psline(-16.000,48.000)(-18.000,50.000)(-20.000,48.000)
\psline(-14.000,50.000)(-16.000,52.000)(-18.000,50.000)
\psline(-32.000,36.000)(-34.000,38.000)(-36.000,36.000)
\psline(-30.000,38.000)(-32.000,40.000)(-34.000,38.000)
\psline(-28.000,40.000)(-30.000,42.000)(-32.000,40.000)
\psline(-26.000,42.000)(-28.000,44.000)(-30.000,42.000)
\psline(-24.000,44.000)(-26.000,46.000)(-28.000,44.000)
\psline(-22.000,46.000)(-24.000,48.000)(-26.000,46.000)
\psline(-20.000,48.000)(-22.000,50.000)(-24.000,48.000)
\psline(-18.000,50.000)(-20.000,52.000)(-22.000,50.000)
\psline(-16.000,52.000)(-18.000,54.000)(-20.000,52.000)
\psline(-34.000,38.000)(-36.000,40.000)(-38.000,38.000)
\psline(-32.000,40.000)(-34.000,42.000)(-36.000,40.000)
\psline(-30.000,42.000)(-32.000,44.000)(-34.000,42.000)
\psline(-28.000,44.000)(-30.000,46.000)(-32.000,44.000)
\psline(-26.000,46.000)(-28.000,48.000)(-30.000,46.000)
\psline(-24.000,48.000)(-26.000,50.000)(-28.000,48.000)
\psline(-22.000,50.000)(-24.000,52.000)(-26.000,50.000)
\psline(-20.000,52.000)(-22.000,54.000)(-24.000,52.000)
\psline(-36.000,40.000)(-38.000,42.000)(-40.000,40.000)
\psline(-34.000,42.000)(-36.000,44.000)(-38.000,42.000)
\psline(-32.000,44.000)(-34.000,46.000)(-36.000,44.000)
\psline(-30.000,46.000)(-32.000,48.000)(-34.000,46.000)
\psline(-28.000,48.000)(-30.000,50.000)(-32.000,48.000)
\psline(-26.000,50.000)(-28.000,52.000)(-30.000,50.000)
\psline(-24.000,52.000)(-26.000,54.000)(-28.000,52.000)
\psline(-38.000,42.000)(-40.000,44.000)(-42.000,42.000)
\psline(-36.000,44.000)(-38.000,46.000)(-40.000,44.000)
\psline(-34.000,46.000)(-36.000,48.000)(-38.000,46.000)
\psline(-32.000,48.000)(-34.000,50.000)(-36.000,48.000)
\psline(-30.000,50.000)(-32.000,52.000)(-34.000,50.000)
\psline(-28.000,52.000)(-30.000,54.000)(-32.000,52.000)
\psline(-26.000,54.000)(-28.000,56.000)(-30.000,54.000)
\psline(-40.000,44.000)(-42.000,46.000)(-44.000,44.000)
\psline(-38.000,46.000)(-40.000,48.000)(-42.000,46.000)
\psline(-36.000,48.000)(-38.000,50.000)(-40.000,48.000)
\psline(-34.000,50.000)(-36.000,52.000)(-38.000,50.000)
\psline(-32.000,52.000)(-34.000,54.000)(-36.000,52.000)
\psline(-30.000,54.000)(-32.000,56.000)(-34.000,54.000)
\psline(-42.000,46.000)(-44.000,48.000)(-46.000,46.000)
\psline(-40.000,48.000)(-42.000,50.000)(-44.000,48.000)
\psline(-38.000,50.000)(-40.000,52.000)(-42.000,50.000)
\psline(-36.000,52.000)(-38.000,54.000)(-40.000,52.000)
\psline(-34.000,54.000)(-36.000,56.000)(-38.000,54.000)
\psline(-32.000,56.000)(-34.000,58.000)(-36.000,56.000)
\psline(-44.000,48.000)(-46.000,50.000)(-48.000,48.000)
\psline(-42.000,50.000)(-44.000,52.000)(-46.000,50.000)
\psline(-40.000,52.000)(-42.000,54.000)(-44.000,52.000)
\psline(-38.000,54.000)(-40.000,56.000)(-42.000,54.000)
\psline(-36.000,56.000)(-38.000,58.000)(-40.000,56.000)
\psline(-34.000,58.000)(-36.000,60.000)(-38.000,58.000)
\psline(-46.000,50.000)(-48.000,52.000)(-50.000,50.000)
\psline(-44.000,52.000)(-46.000,54.000)(-48.000,52.000)
\psline(-42.000,54.000)(-44.000,56.000)(-46.000,54.000)
\psline(-40.000,56.000)(-42.000,58.000)(-44.000,56.000)
\psline(-38.000,58.000)(-40.000,60.000)(-42.000,58.000)
\psline(-48.000,52.000)(-50.000,54.000)(-52.000,52.000)
\psline(-46.000,54.000)(-48.000,56.000)(-50.000,54.000)
\psline(-44.000,56.000)(-46.000,58.000)(-48.000,56.000)
\psline(-42.000,58.000)(-44.000,60.000)(-46.000,58.000)
\psline(-50.000,54.000)(-52.000,56.000)(-54.000,54.000)
\psline(-48.000,56.000)(-50.000,58.000)(-52.000,56.000)
\psline(-46.000,58.000)(-48.000,60.000)(-50.000,58.000)
\psline(-52.000,56.000)(-54.000,58.000)(-56.000,56.000)
\psline(-50.000,58.000)(-52.000,60.000)(-54.000,58.000)
\psline(-48.000,60.000)(-50.000,62.000)(-52.000,60.000)
\psline(-54.000,58.000)(-56.000,60.000)(-58.000,58.000)
\psline(-52.000,60.000)(-54.000,62.000)(-56.000,60.000)
\psline(-56.000,60.000)(-58.000,62.000)(-60.000,60.000)
\psline(-54.000,62.000)(-56.000,64.000)(-58.000,62.000)
\psline(-58.000,62.000)(-60.000,64.000)(-62.000,62.000)
\psline(-56.000,64.000)(-58.000,66.000)(-60.000,64.000)
\psline(-60.000,64.000)(-62.000,66.000)(-64.000,64.000)
\psline(-62.000,66.000)(-64.000,68.000)(-66.000,66.000)
\psline(-64.000,68.000)(-66.000,70.000)(-68.000,68.000)
\psline(-66.000,70.000)(-68.000,72.000)(-70.000,70.000)
\endpspicture

\vspace{-4.65cm}
\psset{unit=0.5429mm}
\pspicture(-230,0)(75,85)
\parametricplot{-2}{2}{t 36.842 mul 
t t 0.5 mul arcsin mul 0.0174533 mul t t mul neg 4 add sqrt add 0.63661977 mul 36.842 mul }
\psline{->}(0,0)(-75,75)\psline{->}(0,0)(75,75)
\psline(0,0)(-70.000,70.000)\psline(0,0)(64.474,64.474)
\psline(1.8421,1.8421)(0.00000,3.6842)(-1.8421,1.8421)
\psline(3.6842,3.6842)(1.8421,5.5263)(0.00000,3.6842)
\psline(5.5263,5.5263)(3.6842,7.3684)(1.8421,5.5263)
\psline(7.3684,7.3684)(5.5263,9.2105)(3.6842,7.3684)
\psline(9.2105,9.2105)(7.3684,11.053)(5.5263,9.2105)
\psline(11.053,11.053)(9.2105,12.895)(7.3684,11.053)
\psline(12.895,12.895)(11.053,14.737)(9.2105,12.895)
\psline(14.737,14.737)(12.895,16.579)(11.053,14.737)
\psline(16.579,16.579)(14.737,18.421)(12.895,16.579)
\psline(18.421,18.421)(16.579,20.263)(14.737,18.421)
\psline(20.263,20.263)(18.421,22.105)(16.579,20.263)
\psline(22.105,22.105)(20.263,23.947)(18.421,22.105)
\psline(23.947,23.947)(22.105,25.789)(20.263,23.947)
\psline(25.789,25.789)(23.947,27.632)(22.105,25.789)
\psline(27.632,27.632)(25.789,29.474)(23.947,27.632)
\psline(29.474,29.474)(27.632,31.316)(25.789,29.474)
\psline(31.316,31.316)(29.474,33.158)(27.632,31.316)
\psline(33.158,33.158)(31.316,35.000)(29.474,33.158)
\psline(35.000,35.000)(33.158,36.842)(31.316,35.000)
\psline(36.842,36.842)(35.000,38.684)(33.158,36.842)
\psline(38.684,38.684)(36.842,40.526)(35.000,38.684)
\psline(40.526,40.526)(38.684,42.368)(36.842,40.526)
\psline(42.368,42.368)(40.526,44.211)(38.684,42.368)
\psline(44.211,44.211)(42.368,46.053)(40.526,44.211)
\psline(46.053,46.053)(44.211,47.895)(42.368,46.053)
\psline(47.895,47.895)(46.053,49.737)(44.211,47.895)
\psline(49.737,49.737)(47.895,51.579)(46.053,49.737)
\psline(51.579,51.579)(49.737,53.421)(47.895,51.579)
\psline(53.421,53.421)(51.579,55.263)(49.737,53.421)
\psline(55.263,55.263)(53.421,57.105)(51.579,55.263)
\psline(57.105,57.105)(55.263,58.947)(53.421,57.105)
\psline(58.947,58.947)(57.105,60.789)(55.263,58.947)
\psline(60.789,60.789)(58.947,62.632)(57.105,60.789)
\psline(62.632,62.632)(60.789,64.474)(58.947,62.632)
\psline(64.474,64.474)(62.632,66.316)(60.789,64.474)
\psline(0.00000,3.6842)(-1.8421,5.5263)(-3.6842,3.6842)
\psline(1.8421,5.5263)(0.00000,7.3684)(-1.8421,5.5263)
\psline(3.6842,7.3684)(1.8421,9.2105)(0.00000,7.3684)
\psline(5.5263,9.2105)(3.6842,11.053)(1.8421,9.2105)
\psline(7.3684,11.053)(5.5263,12.895)(3.6842,11.053)
\psline(9.2105,12.895)(7.3684,14.737)(5.5263,12.895)
\psline(11.053,14.737)(9.2105,16.579)(7.3684,14.737)
\psline(12.895,16.579)(11.053,18.421)(9.2105,16.579)
\psline(14.737,18.421)(12.895,20.263)(11.053,18.421)
\psline(16.579,20.263)(14.737,22.105)(12.895,20.263)
\psline(18.421,22.105)(16.579,23.947)(14.737,22.105)
\psline(20.263,23.947)(18.421,25.789)(16.579,23.947)
\psline(22.105,25.789)(20.263,27.632)(18.421,25.789)
\psline(23.947,27.632)(22.105,29.474)(20.263,27.632)
\psline(25.789,29.474)(23.947,31.316)(22.105,29.474)
\psline(27.632,31.316)(25.789,33.158)(23.947,31.316)
\psline(29.474,33.158)(27.632,35.000)(25.789,33.158)
\psline(31.316,35.000)(29.474,36.842)(27.632,35.000)
\psline(33.158,36.842)(31.316,38.684)(29.474,36.842)
\psline(35.000,38.684)(33.158,40.526)(31.316,38.684)
\psline(36.842,40.526)(35.000,42.368)(33.158,40.526)
\psline(38.684,42.368)(36.842,44.211)(35.000,42.368)
\psline(40.526,44.211)(38.684,46.053)(36.842,44.211)
\psline(42.368,46.053)(40.526,47.895)(38.684,46.053)
\psline(44.211,47.895)(42.368,49.737)(40.526,47.895)
\psline(46.053,49.737)(44.211,51.579)(42.368,49.737)
\psline(47.895,51.579)(46.053,53.421)(44.211,51.579)
\psline(49.737,53.421)(47.895,55.263)(46.053,53.421)
\psline(51.579,55.263)(49.737,57.105)(47.895,55.263)
\psline(53.421,57.105)(51.579,58.947)(49.737,57.105)
\psline(55.263,58.947)(53.421,60.789)(51.579,58.947)
\psline(57.105,60.789)(55.263,62.632)(53.421,60.789)
\psline(-1.8421,5.5263)(-3.6842,7.3684)(-5.5263,5.5263)
\psline(0.00000,7.3684)(-1.8421,9.2105)(-3.6842,7.3684)
\psline(1.8421,9.2105)(0.00000,11.053)(-1.8421,9.2105)
\psline(3.6842,11.053)(1.8421,12.895)(0.00000,11.053)
\psline(5.5263,12.895)(3.6842,14.737)(1.8421,12.895)
\psline(7.3684,14.737)(5.5263,16.579)(3.6842,14.737)
\psline(9.2105,16.579)(7.3684,18.421)(5.5263,16.579)
\psline(11.053,18.421)(9.2105,20.263)(7.3684,18.421)
\psline(12.895,20.263)(11.053,22.105)(9.2105,20.263)
\psline(14.737,22.105)(12.895,23.947)(11.053,22.105)
\psline(16.579,23.947)(14.737,25.789)(12.895,23.947)
\psline(18.421,25.789)(16.579,27.632)(14.737,25.789)
\psline(20.263,27.632)(18.421,29.474)(16.579,27.632)
\psline(22.105,29.474)(20.263,31.316)(18.421,29.474)
\psline(23.947,31.316)(22.105,33.158)(20.263,31.316)
\psline(25.789,33.158)(23.947,35.000)(22.105,33.158)
\psline(27.632,35.000)(25.789,36.842)(23.947,35.000)
\psline(29.474,36.842)(27.632,38.684)(25.789,36.842)
\psline(31.316,38.684)(29.474,40.526)(27.632,38.684)
\psline(33.158,40.526)(31.316,42.368)(29.474,40.526)
\psline(35.000,42.368)(33.158,44.211)(31.316,42.368)
\psline(36.842,44.211)(35.000,46.053)(33.158,44.211)
\psline(38.684,46.053)(36.842,47.895)(35.000,46.053)
\psline(40.526,47.895)(38.684,49.737)(36.842,47.895)
\psline(42.368,49.737)(40.526,51.579)(38.684,49.737)
\psline(44.211,51.579)(42.368,53.421)(40.526,51.579)
\psline(46.053,53.421)(44.211,55.263)(42.368,53.421)
\psline(47.895,55.263)(46.053,57.105)(44.211,55.263)
\psline(49.737,57.105)(47.895,58.947)(46.053,57.105)
\psline(51.579,58.947)(49.737,60.789)(47.895,58.947)
\psline(-3.6842,7.3684)(-5.5263,9.2105)(-7.3684,7.3684)
\psline(-1.8421,9.2105)(-3.6842,11.053)(-5.5263,9.2105)
\psline(0.00000,11.053)(-1.8421,12.895)(-3.6842,11.053)
\psline(1.8421,12.895)(0.00000,14.737)(-1.8421,12.895)
\psline(3.6842,14.737)(1.8421,16.579)(0.00000,14.737)
\psline(5.5263,16.579)(3.6842,18.421)(1.8421,16.579)
\psline(7.3684,18.421)(5.5263,20.263)(3.6842,18.421)
\psline(9.2105,20.263)(7.3684,22.105)(5.5263,20.263)
\psline(11.053,22.105)(9.2105,23.947)(7.3684,22.105)
\psline(12.895,23.947)(11.053,25.789)(9.2105,23.947)
\psline(14.737,25.789)(12.895,27.632)(11.053,25.789)
\psline(16.579,27.632)(14.737,29.474)(12.895,27.632)
\psline(18.421,29.474)(16.579,31.316)(14.737,29.474)
\psline(20.263,31.316)(18.421,33.158)(16.579,31.316)
\psline(22.105,33.158)(20.263,35.000)(18.421,33.158)
\psline(23.947,35.000)(22.105,36.842)(20.263,35.000)
\psline(25.789,36.842)(23.947,38.684)(22.105,36.842)
\psline(27.632,38.684)(25.789,40.526)(23.947,38.684)
\psline(29.474,40.526)(27.632,42.368)(25.789,40.526)
\psline(31.316,42.368)(29.474,44.211)(27.632,42.368)
\psline(33.158,44.211)(31.316,46.053)(29.474,44.211)
\psline(35.000,46.053)(33.158,47.895)(31.316,46.053)
\psline(36.842,47.895)(35.000,49.737)(33.158,47.895)
\psline(38.684,49.737)(36.842,51.579)(35.000,49.737)
\psline(40.526,51.579)(38.684,53.421)(36.842,51.579)
\psline(42.368,53.421)(40.526,55.263)(38.684,53.421)
\psline(44.211,55.263)(42.368,57.105)(40.526,55.263)
\psline(46.053,57.105)(44.211,58.947)(42.368,57.105)
\psline(-5.5263,9.2105)(-7.3684,11.053)(-9.2105,9.2105)
\psline(-3.6842,11.053)(-5.5263,12.895)(-7.3684,11.053)
\psline(-1.8421,12.895)(-3.6842,14.737)(-5.5263,12.895)
\psline(0.00000,14.737)(-1.8421,16.579)(-3.6842,14.737)
\psline(1.8421,16.579)(0.00000,18.421)(-1.8421,16.579)
\psline(3.6842,18.421)(1.8421,20.263)(0.00000,18.421)
\psline(5.5263,20.263)(3.6842,22.105)(1.8421,20.263)
\psline(7.3684,22.105)(5.5263,23.947)(3.6842,22.105)
\psline(9.2105,23.947)(7.3684,25.789)(5.5263,23.947)
\psline(11.053,25.789)(9.2105,27.632)(7.3684,25.789)
\psline(12.895,27.632)(11.053,29.474)(9.2105,27.632)
\psline(14.737,29.474)(12.895,31.316)(11.053,29.474)
\psline(16.579,31.316)(14.737,33.158)(12.895,31.316)
\psline(18.421,33.158)(16.579,35.000)(14.737,33.158)
\psline(20.263,35.000)(18.421,36.842)(16.579,35.000)
\psline(22.105,36.842)(20.263,38.684)(18.421,36.842)
\psline(23.947,38.684)(22.105,40.526)(20.263,38.684)
\psline(25.789,40.526)(23.947,42.368)(22.105,40.526)
\psline(27.632,42.368)(25.789,44.211)(23.947,42.368)
\psline(29.474,44.211)(27.632,46.053)(25.789,44.211)
\psline(31.316,46.053)(29.474,47.895)(27.632,46.053)
\psline(33.158,47.895)(31.316,49.737)(29.474,47.895)
\psline(-7.3684,11.053)(-9.2105,12.895)(-11.053,11.053)
\psline(-5.5263,12.895)(-7.3684,14.737)(-9.2105,12.895)
\psline(-3.6842,14.737)(-5.5263,16.579)(-7.3684,14.737)
\psline(-1.8421,16.579)(-3.6842,18.421)(-5.5263,16.579)
\psline(0.00000,18.421)(-1.8421,20.263)(-3.6842,18.421)
\psline(1.8421,20.263)(0.00000,22.105)(-1.8421,20.263)
\psline(3.6842,22.105)(1.8421,23.947)(0.00000,22.105)
\psline(5.5263,23.947)(3.6842,25.789)(1.8421,23.947)
\psline(7.3684,25.789)(5.5263,27.632)(3.6842,25.789)
\psline(9.2105,27.632)(7.3684,29.474)(5.5263,27.632)
\psline(11.053,29.474)(9.2105,31.316)(7.3684,29.474)
\psline(12.895,31.316)(11.053,33.158)(9.2105,31.316)
\psline(14.737,33.158)(12.895,35.000)(11.053,33.158)
\psline(16.579,35.000)(14.737,36.842)(12.895,35.000)
\psline(18.421,36.842)(16.579,38.684)(14.737,36.842)
\psline(20.263,38.684)(18.421,40.526)(16.579,38.684)
\psline(22.105,40.526)(20.263,42.368)(18.421,40.526)
\psline(23.947,42.368)(22.105,44.211)(20.263,42.368)
\psline(25.789,44.211)(23.947,46.053)(22.105,44.211)
\psline(27.632,46.053)(25.789,47.895)(23.947,46.053)
\psline(29.474,47.895)(27.632,49.737)(25.789,47.895)
\psline(-9.2105,12.895)(-11.053,14.737)(-12.895,12.895)
\psline(-7.3684,14.737)(-9.2105,16.579)(-11.053,14.737)
\psline(-5.5263,16.579)(-7.3684,18.421)(-9.2105,16.579)
\psline(-3.6842,18.421)(-5.5263,20.263)(-7.3684,18.421)
\psline(-1.8421,20.263)(-3.6842,22.105)(-5.5263,20.263)
\psline(0.00000,22.105)(-1.8421,23.947)(-3.6842,22.105)
\psline(1.8421,23.947)(0.00000,25.789)(-1.8421,23.947)
\psline(3.6842,25.789)(1.8421,27.632)(0.00000,25.789)
\psline(5.5263,27.632)(3.6842,29.474)(1.8421,27.632)
\psline(7.3684,29.474)(5.5263,31.316)(3.6842,29.474)
\psline(9.2105,31.316)(7.3684,33.158)(5.5263,31.316)
\psline(11.053,33.158)(9.2105,35.000)(7.3684,33.158)
\psline(12.895,35.000)(11.053,36.842)(9.2105,35.000)
\psline(14.737,36.842)(12.895,38.684)(11.053,36.842)
\psline(16.579,38.684)(14.737,40.526)(12.895,38.684)
\psline(18.421,40.526)(16.579,42.368)(14.737,40.526)
\psline(20.263,42.368)(18.421,44.211)(16.579,42.368)
\psline(22.105,44.211)(20.263,46.053)(18.421,44.211)
\psline(23.947,46.053)(22.105,47.895)(20.263,46.053)
\psline(-11.053,14.737)(-12.895,16.579)(-14.737,14.737)
\psline(-9.2105,16.579)(-11.053,18.421)(-12.895,16.579)
\psline(-7.3684,18.421)(-9.2105,20.263)(-11.053,18.421)
\psline(-5.5263,20.263)(-7.3684,22.105)(-9.2105,20.263)
\psline(-3.6842,22.105)(-5.5263,23.947)(-7.3684,22.105)
\psline(-1.8421,23.947)(-3.6842,25.789)(-5.5263,23.947)
\psline(0.00000,25.789)(-1.8421,27.632)(-3.6842,25.789)
\psline(1.8421,27.632)(0.00000,29.474)(-1.8421,27.632)
\psline(3.6842,29.474)(1.8421,31.316)(0.00000,29.474)
\psline(5.5263,31.316)(3.6842,33.158)(1.8421,31.316)
\psline(7.3684,33.158)(5.5263,35.000)(3.6842,33.158)
\psline(9.2105,35.000)(7.3684,36.842)(5.5263,35.000)
\psline(11.053,36.842)(9.2105,38.684)(7.3684,36.842)
\psline(12.895,38.684)(11.053,40.526)(9.2105,38.684)
\psline(14.737,40.526)(12.895,42.368)(11.053,40.526)
\psline(16.579,42.368)(14.737,44.211)(12.895,42.368)
\psline(18.421,44.211)(16.579,46.053)(14.737,44.211)
\psline(20.263,46.053)(18.421,47.895)(16.579,46.053)
\psline(22.105,47.895)(20.263,49.737)(18.421,47.895)
\psline(-12.895,16.579)(-14.737,18.421)(-16.579,16.579)
\psline(-11.053,18.421)(-12.895,20.263)(-14.737,18.421)
\psline(-9.2105,20.263)(-11.053,22.105)(-12.895,20.263)
\psline(-7.3684,22.105)(-9.2105,23.947)(-11.053,22.105)
\psline(-5.5263,23.947)(-7.3684,25.789)(-9.2105,23.947)
\psline(-3.6842,25.789)(-5.5263,27.632)(-7.3684,25.789)
\psline(-1.8421,27.632)(-3.6842,29.474)(-5.5263,27.632)
\psline(0.00000,29.474)(-1.8421,31.316)(-3.6842,29.474)
\psline(1.8421,31.316)(0.00000,33.158)(-1.8421,31.316)
\psline(3.6842,33.158)(1.8421,35.000)(0.00000,33.158)
\psline(5.5263,35.000)(3.6842,36.842)(1.8421,35.000)
\psline(7.3684,36.842)(5.5263,38.684)(3.6842,36.842)
\psline(9.2105,38.684)(7.3684,40.526)(5.5263,38.684)
\psline(11.053,40.526)(9.2105,42.368)(7.3684,40.526)
\psline(12.895,42.368)(11.053,44.211)(9.2105,42.368)
\psline(14.737,44.211)(12.895,46.053)(11.053,44.211)
\psline(16.579,46.053)(14.737,47.895)(12.895,46.053)
\psline(18.421,47.895)(16.579,49.737)(14.737,47.895)
\psline(-14.737,18.421)(-16.579,20.263)(-18.421,18.421)
\psline(-12.895,20.263)(-14.737,22.105)(-16.579,20.263)
\psline(-11.053,22.105)(-12.895,23.947)(-14.737,22.105)
\psline(-9.2105,23.947)(-11.053,25.789)(-12.895,23.947)
\psline(-7.3684,25.789)(-9.2105,27.632)(-11.053,25.789)
\psline(-5.5263,27.632)(-7.3684,29.474)(-9.2105,27.632)
\psline(-3.6842,29.474)(-5.5263,31.316)(-7.3684,29.474)
\psline(-1.8421,31.316)(-3.6842,33.158)(-5.5263,31.316)
\psline(0.00000,33.158)(-1.8421,35.000)(-3.6842,33.158)
\psline(1.8421,35.000)(0.00000,36.842)(-1.8421,35.000)
\psline(3.6842,36.842)(1.8421,38.684)(0.00000,36.842)
\psline(5.5263,38.684)(3.6842,40.526)(1.8421,38.684)
\psline(7.3684,40.526)(5.5263,42.368)(3.6842,40.526)
\psline(9.2105,42.368)(7.3684,44.211)(5.5263,42.368)
\psline(11.053,44.211)(9.2105,46.053)(7.3684,44.211)
\psline(12.895,46.053)(11.053,47.895)(9.2105,46.053)
\psline(-16.579,20.263)(-18.421,22.105)(-20.263,20.263)
\psline(-14.737,22.105)(-16.579,23.947)(-18.421,22.105)
\psline(-12.895,23.947)(-14.737,25.789)(-16.579,23.947)
\psline(-11.053,25.789)(-12.895,27.632)(-14.737,25.789)
\psline(-9.2105,27.632)(-11.053,29.474)(-12.895,27.632)
\psline(-7.3684,29.474)(-9.2105,31.316)(-11.053,29.474)
\psline(-5.5263,31.316)(-7.3684,33.158)(-9.2105,31.316)
\psline(-3.6842,33.158)(-5.5263,35.000)(-7.3684,33.158)
\psline(-1.8421,35.000)(-3.6842,36.842)(-5.5263,35.000)
\psline(0.00000,36.842)(-1.8421,38.684)(-3.6842,36.842)
\psline(1.8421,38.684)(0.00000,40.526)(-1.8421,38.684)
\psline(3.6842,40.526)(1.8421,42.368)(0.00000,40.526)
\psline(5.5263,42.368)(3.6842,44.211)(1.8421,42.368)
\psline(7.3684,44.211)(5.5263,46.053)(3.6842,44.211)
\psline(9.2105,46.053)(7.3684,47.895)(5.5263,46.053)
\psline(-18.421,22.105)(-20.263,23.947)(-22.105,22.105)
\psline(-16.579,23.947)(-18.421,25.789)(-20.263,23.947)
\psline(-14.737,25.789)(-16.579,27.632)(-18.421,25.789)
\psline(-12.895,27.632)(-14.737,29.474)(-16.579,27.632)
\psline(-11.053,29.474)(-12.895,31.316)(-14.737,29.474)
\psline(-9.2105,31.316)(-11.053,33.158)(-12.895,31.316)
\psline(-7.3684,33.158)(-9.2105,35.000)(-11.053,33.158)
\psline(-5.5263,35.000)(-7.3684,36.842)(-9.2105,35.000)
\psline(-3.6842,36.842)(-5.5263,38.684)(-7.3684,36.842)
\psline(-1.8421,38.684)(-3.6842,40.526)(-5.5263,38.684)
\psline(0.00000,40.526)(-1.8421,42.368)(-3.6842,40.526)
\psline(1.8421,42.368)(0.00000,44.211)(-1.8421,42.368)
\psline(3.6842,44.211)(1.8421,46.053)(0.00000,44.211)
\psline(5.5263,46.053)(3.6842,47.895)(1.8421,46.053)
\psline(-20.263,23.947)(-22.105,25.789)(-23.947,23.947)
\psline(-18.421,25.789)(-20.263,27.632)(-22.105,25.789)
\psline(-16.579,27.632)(-18.421,29.474)(-20.263,27.632)
\psline(-14.737,29.474)(-16.579,31.316)(-18.421,29.474)
\psline(-12.895,31.316)(-14.737,33.158)(-16.579,31.316)
\psline(-11.053,33.158)(-12.895,35.000)(-14.737,33.158)
\psline(-9.2105,35.000)(-11.053,36.842)(-12.895,35.000)
\psline(-7.3684,36.842)(-9.2105,38.684)(-11.053,36.842)
\psline(-5.5263,38.684)(-7.3684,40.526)(-9.2105,38.684)
\psline(-3.6842,40.526)(-5.5263,42.368)(-7.3684,40.526)
\psline(-1.8421,42.368)(-3.6842,44.211)(-5.5263,42.368)
\psline(0.00000,44.211)(-1.8421,46.053)(-3.6842,44.211)
\psline(-22.105,25.789)(-23.947,27.632)(-25.789,25.789)
\psline(-20.263,27.632)(-22.105,29.474)(-23.947,27.632)
\psline(-18.421,29.474)(-20.263,31.316)(-22.105,29.474)
\psline(-16.579,31.316)(-18.421,33.158)(-20.263,31.316)
\psline(-14.737,33.158)(-16.579,35.000)(-18.421,33.158)
\psline(-12.895,35.000)(-14.737,36.842)(-16.579,35.000)
\psline(-11.053,36.842)(-12.895,38.684)(-14.737,36.842)
\psline(-9.2105,38.684)(-11.053,40.526)(-12.895,38.684)
\psline(-7.3684,40.526)(-9.2105,42.368)(-11.053,40.526)
\psline(-5.5263,42.368)(-7.3684,44.211)(-9.2105,42.368)
\psline(-3.6842,44.211)(-5.5263,46.053)(-7.3684,44.211)
\psline(-1.8421,46.053)(-3.6842,47.895)(-5.5263,46.053)
\psline(-23.947,27.632)(-25.789,29.474)(-27.632,27.632)
\psline(-22.105,29.474)(-23.947,31.316)(-25.789,29.474)
\psline(-20.263,31.316)(-22.105,33.158)(-23.947,31.316)
\psline(-18.421,33.158)(-20.263,35.000)(-22.105,33.158)
\psline(-16.579,35.000)(-18.421,36.842)(-20.263,35.000)
\psline(-14.737,36.842)(-16.579,38.684)(-18.421,36.842)
\psline(-12.895,38.684)(-14.737,40.526)(-16.579,38.684)
\psline(-11.053,40.526)(-12.895,42.368)(-14.737,40.526)
\psline(-9.2105,42.368)(-11.053,44.211)(-12.895,42.368)
\psline(-7.3684,44.211)(-9.2105,46.053)(-11.053,44.211)
\psline(-5.5263,46.053)(-7.3684,47.895)(-9.2105,46.053)
\psline(-25.789,29.474)(-27.632,31.316)(-29.474,29.474)
\psline(-23.947,31.316)(-25.789,33.158)(-27.632,31.316)
\psline(-22.105,33.158)(-23.947,35.000)(-25.789,33.158)
\psline(-20.263,35.000)(-22.105,36.842)(-23.947,35.000)
\psline(-18.421,36.842)(-20.263,38.684)(-22.105,36.842)
\psline(-16.579,38.684)(-18.421,40.526)(-20.263,38.684)
\psline(-14.737,40.526)(-16.579,42.368)(-18.421,40.526)
\psline(-12.895,42.368)(-14.737,44.211)(-16.579,42.368)
\psline(-11.053,44.211)(-12.895,46.053)(-14.737,44.211)
\psline(-9.2105,46.053)(-11.053,47.895)(-12.895,46.053)
\psline(-7.3684,47.895)(-9.2105,49.737)(-11.053,47.895)
\psline(-27.632,31.316)(-29.474,33.158)(-31.316,31.316)
\psline(-25.789,33.158)(-27.632,35.000)(-29.474,33.158)
\psline(-23.947,35.000)(-25.789,36.842)(-27.632,35.000)
\psline(-22.105,36.842)(-23.947,38.684)(-25.789,36.842)
\psline(-20.263,38.684)(-22.105,40.526)(-23.947,38.684)
\psline(-18.421,40.526)(-20.263,42.368)(-22.105,40.526)
\psline(-16.579,42.368)(-18.421,44.211)(-20.263,42.368)
\psline(-14.737,44.211)(-16.579,46.053)(-18.421,44.211)
\psline(-12.895,46.053)(-14.737,47.895)(-16.579,46.053)
\psline(-11.053,47.895)(-12.895,49.737)(-14.737,47.895)
\psline(-9.2105,49.737)(-11.053,51.579)(-12.895,49.737)
\psline(-29.474,33.158)(-31.316,35.000)(-33.158,33.158)
\psline(-27.632,35.000)(-29.474,36.842)(-31.316,35.000)
\psline(-25.789,36.842)(-27.632,38.684)(-29.474,36.842)
\psline(-23.947,38.684)(-25.789,40.526)(-27.632,38.684)
\psline(-22.105,40.526)(-23.947,42.368)(-25.789,40.526)
\psline(-20.263,42.368)(-22.105,44.211)(-23.947,42.368)
\psline(-18.421,44.211)(-20.263,46.053)(-22.105,44.211)
\psline(-16.579,46.053)(-18.421,47.895)(-20.263,46.053)
\psline(-14.737,47.895)(-16.579,49.737)(-18.421,47.895)
\psline(-31.316,35.000)(-33.158,36.842)(-35.000,35.000)
\psline(-29.474,36.842)(-31.316,38.684)(-33.158,36.842)
\psline(-27.632,38.684)(-29.474,40.526)(-31.316,38.684)
\psline(-25.789,40.526)(-27.632,42.368)(-29.474,40.526)
\psline(-23.947,42.368)(-25.789,44.211)(-27.632,42.368)
\psline(-22.105,44.211)(-23.947,46.053)(-25.789,44.211)
\psline(-20.263,46.053)(-22.105,47.895)(-23.947,46.053)
\psline(-18.421,47.895)(-20.263,49.737)(-22.105,47.895)
\psline(-33.158,36.842)(-35.000,38.684)(-36.842,36.842)
\psline(-31.316,38.684)(-33.158,40.526)(-35.000,38.684)
\psline(-29.474,40.526)(-31.316,42.368)(-33.158,40.526)
\psline(-27.632,42.368)(-29.474,44.211)(-31.316,42.368)
\psline(-25.789,44.211)(-27.632,46.053)(-29.474,44.211)
\psline(-23.947,46.053)(-25.789,47.895)(-27.632,46.053)
\psline(-22.105,47.895)(-23.947,49.737)(-25.789,47.895)
\psline(-20.263,49.737)(-22.105,51.579)(-23.947,49.737)
\psline(-35.000,38.684)(-36.842,40.526)(-38.684,38.684)
\psline(-33.158,40.526)(-35.000,42.368)(-36.842,40.526)
\psline(-31.316,42.368)(-33.158,44.211)(-35.000,42.368)
\psline(-29.474,44.211)(-31.316,46.053)(-33.158,44.211)
\psline(-27.632,46.053)(-29.474,47.895)(-31.316,46.053)
\psline(-25.789,47.895)(-27.632,49.737)(-29.474,47.895)
\psline(-36.842,40.526)(-38.684,42.368)(-40.526,40.526)
\psline(-35.000,42.368)(-36.842,44.211)(-38.684,42.368)
\psline(-33.158,44.211)(-35.000,46.053)(-36.842,44.211)
\psline(-31.316,46.053)(-33.158,47.895)(-35.000,46.053)
\psline(-29.474,47.895)(-31.316,49.737)(-33.158,47.895)
\psline(-27.632,49.737)(-29.474,51.579)(-31.316,49.737)
\psline(-38.684,42.368)(-40.526,44.211)(-42.368,42.368)
\psline(-36.842,44.211)(-38.684,46.053)(-40.526,44.211)
\psline(-35.000,46.053)(-36.842,47.895)(-38.684,46.053)
\psline(-33.158,47.895)(-35.000,49.737)(-36.842,47.895)
\psline(-31.316,49.737)(-33.158,51.579)(-35.000,49.737)
\psline(-29.474,51.579)(-31.316,53.421)(-33.158,51.579)
\psline(-40.526,44.211)(-42.368,46.053)(-44.211,44.211)
\psline(-38.684,46.053)(-40.526,47.895)(-42.368,46.053)
\psline(-36.842,47.895)(-38.684,49.737)(-40.526,47.895)
\psline(-35.000,49.737)(-36.842,51.579)(-38.684,49.737)
\psline(-42.368,46.053)(-44.211,47.895)(-46.053,46.053)
\psline(-40.526,47.895)(-42.368,49.737)(-44.211,47.895)
\psline(-38.684,49.737)(-40.526,51.579)(-42.368,49.737)
\psline(-36.842,51.579)(-38.684,53.421)(-40.526,51.579)
\psline(-44.211,47.895)(-46.053,49.737)(-47.895,47.895)
\psline(-42.368,49.737)(-44.211,51.579)(-46.053,49.737)
\psline(-40.526,51.579)(-42.368,53.421)(-44.211,51.579)
\psline(-38.684,53.421)(-40.526,55.263)(-42.368,53.421)
\psline(-46.053,49.737)(-47.895,51.579)(-49.737,49.737)
\psline(-44.211,51.579)(-46.053,53.421)(-47.895,51.579)
\psline(-42.368,53.421)(-44.211,55.263)(-46.053,53.421)
\psline(-47.895,51.579)(-49.737,53.421)(-51.579,51.579)
\psline(-46.053,53.421)(-47.895,55.263)(-49.737,53.421)
\psline(-49.737,53.421)(-51.579,55.263)(-53.421,53.421)
\psline(-47.895,55.263)(-49.737,57.105)(-51.579,55.263)
\psline(-51.579,55.263)(-53.421,57.105)(-55.263,55.263)
\psline(-49.737,57.105)(-51.579,58.947)(-53.421,57.105)
\psline(-53.421,57.105)(-55.263,58.947)(-57.105,57.105)
\psline(-51.579,58.947)(-53.421,60.789)(-55.263,58.947)
\psline(-55.263,58.947)(-57.105,60.789)(-58.947,58.947)
\psline(-53.421,60.789)(-55.263,62.632)(-57.105,60.789)
\psline(-57.105,60.789)(-58.947,62.632)(-60.789,60.789)
\psline(-58.947,62.632)(-60.789,64.474)(-62.632,62.632)
\psline(-60.789,64.474)(-62.632,66.316)(-64.474,64.474)
\psline(-62.632,66.316)(-64.474,68.158)(-66.316,66.316)
\psline(-64.474,68.158)(-66.316,70.000)(-68.158,68.158)
\psline(-66.316,70.000)(-68.158,71.842)(-70.000,70.000)
\endpspicture}{Les mesures de Plancherel et de Gelfand ont la même asymptotique au premier ordre, mais les mesures de Gelfand ont des fluctuations différentes, et plus grandes.}{Asymptotiques comparées des mesures de Plancherel et de Gelfand}

Clairement, les formes limites semblent identiques, et nous donnerons dans la section \ref{asymptoticgelfand} une preuve rigoureuse de ce fait, et une explication intuitive. En revanche, les fluctuations du diagramme par rapport à sa forme limite semblent plus grandes dans le cas des mesures de Gelfand que dans le cas des mesures de Plancherel --- ce point est nettement moins clair sur le dessin, et éventuellement difficile à constater à l'oeil nu\footnote{Pour rendre apparent le phénomène, on peut colorier l'aire comprise entre le diagramme et sa forme limite ; il est alors tout à fait clair que l'aire coloriée est plus grande dans le cas des mesures de Gelfand.} ; comme nous le verrons par la suite, c'est parce que les fluctuations sont seulement $\sqrt{2}$ fois plus grandes. De nouveau, nous donnerons une preuve rigoureuse de ces observations dans la section \ref{asymptoticgelfand}. Notons que les mesures de Gelfand, et plus généralement les \textbf{$\beta$-mesures de Plancherel}, ont déjà été étudiées dans \cite{BR01} ; en particulier, J. Baik et E. Rains ont établi un résultat analogue au théorème \ref{bdj} pour la longueur de la plus longue ligne ou colonne d'une partition sous la mesure de Gelfand (il suffit de remplacer le GUE par le GOE dans l'énoncé du théorème). Ainsi, la vraie nouveauté dans notre travail est la détermination des fluctuations globales. Les raisonnements de ce chapitre seront sensiblement identiques à ceux effectués dans le chapitre précédent, mais ils mettront aussi en jeu l'expression asymptotique du nombre d'\textbf{involutions} de taille $n$ (voir le paragraphe \ref{stirlinvolution}), et le \textbf{degré de Kerov} des observables de diagrammes. À un terme constant et un facteur multiplicatif $\sqrt{2}$ près, nous retrouverons le processus gaussien généralisé de Kerov ; ainsi, ce processus semble jouer un rôle universel dans le contexte des modèles de partitions aléatoires issus de la théorie des représentations.\bigskip

\section{Modèles de Gelfand des groupes symétriques}\label{gelfandmodel}
Si $G$ est un groupe fini, on appelle \textbf{modèle de Gelfand} de $G$ une représentation de $G$ dans laquelle toute représentation irréductible $\lambda \in \widehat{G}$ apparaît exactement une fois. Dans le cas des groupes symétriques, une description combinatoire d'une telle représentation a été proposée par Adin, Postnikov et Roichman dans \cite{APR07} ; la suite de cette section est consacrée à l'exposé de cette construction, qui a depuis été généralisée au cas des produits en couronne $(\Z/r\Z)\wr \sym_{n}$ et des groupes de réflexions complexes $G(r,p,n)$ (voir \cite{APR08,CF10}). \bigskip\bigskip

Si $M$ est un modèle de Gelfand du groupe symétrique $\sym_{n}$, alors sa dimension est nécessairement égale à $\sum_{\lambda \in \ym_{n}}\dim \lambda$. Or :
\begin{lemma}[Dimension d'un modèle de Gelfand du groupe symétrique]
Pour tout entier $n$, la quantité $\sum_{\lambda \in \ym_{n}}\dim \lambda$ est aussi le nombre de permutations de $\sym_{n}$ qui vérifient $\sigma^{2}=\id_{\lle 1,n\rre}$, c'est-à-dire le nombre d'involutions de taille $n$. \end{lemma}
\begin{proof} Notons $\mathfrak{I}_{n}$ l'ensemble des involutions de taille $n$. Si $\sigma$ est une involution, alors $\sigma=\sigma^{-1}$, donc la paire de tableaux standards $(P(\sigma),Q(\sigma)=P(\sigma^{-1}))$ attachée à $\sigma$ par la correspondance RSK vérifie l'identité $P=Q$. Réciproquement, comme la correspondance RSK est une bijection, si $P=Q$, alors comme $$\mathrm{RSK}(\sigma)=(P,Q) \iff \mathrm{RSK}(\sigma^{-1})=(Q,P)\,,$$ on a $\sigma=\sigma^{-1}$ ; ainsi, $\sigma$ est une involution si et seulement si les deux tableaux standards qui lui sont associés sont identiques. Par suite, la correspondance RSK établit une bijection entre tableaux standards de taille $n$ et involutions de taille $n$, et 
$$\card\,\mathfrak{I_{n}}=\sum_{\lambda \in \ym_{n}}\card \,\mathrm{Std}(\lambda) = \sum_{\lambda \in \ym_{n}}\dim \lambda\,.\vspace{-5mm}$$
\end{proof}\bigskip

Compte tenu de cette identité, il est naturel de rechercher un modèle de Gelfand de $\sym_{n}$ sous la forme d'une action de $\sym_{n}$ sur l'espace $\C\mathfrak{I}_{n}$ des combinaisons linéaires formelles d'involutions de taille $n$. Si $v \in \mathfrak{I}_{n}$, notons $B_{v}$ l'élément
de base correspondant dans $\C\mathfrak{I}_{n}$. On rappelle que les \textbf{descentes} d'une permutation $\sigma \in \sym_{n}$ sont les indices $i \in \lle 1,n-1\rre $ tels que $\ell(\sigma \,s_{i}) < \ell(\sigma)$, où $s_{i}=(i,i+1)$ est l'une des transpositions élémentaires engendrant le groupe de Coxeter $\sym_{n}$. Nous noterons $D(\sigma)$ l'ensemble des descentes de $\sigma$. Alors :
\begin{proposition}[Description combinatoire d'un modèle de Gelfand du groupe symétrique, \cite{Mel10}]\label{gmsym}
Soit $\rho : \sym_{n} \to \GL(\C\mathfrak{I}_{n})$ la représentation linéaire complexe définie par :
$$s_{i}\cdot B_{v} = \begin{cases}
- B_{v}&\text{si }s_{i}\,v\,s_{i}=v \text{ et }i\in D(v),  \\
B_{s_{i}vs_{i}}&\text{sinon}. 
\end{cases}$$
Cette représentation est bien définie, et constitue un modèle de Gelfand du groupe symétrique.
\end{proposition}
\noindent La preuve de la proposition est liée au fait suivant : pour tout groupe fini $G$, on peut calculer explicitement la somme de caractères irréductibles $\sum_{\lambda \in \widehat{G}} \,\varsigma^{\lambda}(g)$, voir \cite[chapitre 4]{Isa94}. Dans le cas particulier où toutes les représentations de $G$ sont réelles (ce qui est le cas pour tout groupe de Coxeter), on a ainsi :
$$\forall g \in G,\,\,\sum_{\lambda \in \widehat{G}} \varsigma^{\lambda}(g)=\card\{h \in G\,\,|\,\,h^{2}=g\}\,.$$
Or, la représentation décrite dans la proposition \ref{gmsym} a exactement ce caractère, voir \cite{APR07} ; ceci démontre la proposition, et on retrouve au passage l'identité $I_{n}=\card \mathfrak{I}_{n}=\sum_{\lambda \in \ym_{n}}\dim \lambda$ en évaluant la formule qui précède en $g=1$.\bigskip
\bigskip

La \textbf{mesure de Gelfand} $\Gel_{n}$ est la mesure de probabilité sur les partitions de taille $n$ associée à la représentation de $\sym_{n}$ décrite par la proposition \ref{gmsym} ; ainsi, $\Gel_{n}[\lambda]=\dim \lambda/I_{n}$, où $I_{n}$ est le nombre d'involutions de taille $n$. C'est aussi la mesure image de la mesure uniforme sur $\mathfrak{I}_{n}$ par la correspondance RSK. Comme d'habitude, on va s'intéresser à l'asymptotique des partitions sous ces mesures, en utilisant l'algèbre des observables de diagrammes pour mener une <<~technique de moments non commutatifs~>>. Si $\chi^{\lambda}$ est un caractère irréductible normalisé de $\sym_{n}$ tiré aléatoirement suivant $\Gel_{n}$, alors pour toute permutation $\sigma$, 
$$\Gel_{n}[\chi^{\lambda}(\sigma)]=\frac{1}{I_{n}}\,\sum_{\lambda \in \ym_{n}} (\dim \lambda) \,\chi^{\lambda}(\sigma)=\frac{1}{I_{n}}\,\sum_{\lambda \in \ym_{n}} \varsigma^{\lambda}(\sigma)=\frac{\tr\rho(\sigma)}{\tr\rho(\id)}=\frac{\card\{\tau \in \sym_{n}\,\,|\,\,\tau^{2}=\sigma\}}{\card\{\tau \in \sym_{n}\,\,|\,\,\tau^{2}=\id\}}\,.$$
Ceci va permettre de calculer facilement l'espérance des caractères centraux $\varSigma_{\mu}$ sous les mesures de Gelfand.

\section{Décompte des involutions et des racines carrées dans le groupe symétrique}\label{stirlinvolution}
Soit $\mu=1^{m_{1}}\,2^{m_{2}}\,\cdots\, s^{m_{s}}$ une partition de taille $n$, et $\sigma_{\mu} \in \sym_{n}$ une permutation de type cyclique $\mu$. Compte tenu de ce qui précède, il est utile de savoir dénombrer le nombre de permutations $\tau$ telles que $\sigma_{\mu}=\tau^{2}$. Commen\c cons par le cas où $\mu=1^{n}$ et $\sigma_{\mu}=\id_{\lle 1,n\rre}$, c'est-à-dire par l'énumération des involutions. Une permutation $\tau$ est une involution si et seulement si c'est un produit de transpositions à supports disjoints, donc si et seulement si elle est de type cyclique $1^{n-2k}\,2^{k}$. Par suite :
$$I_{n}= \sum_{k=0}^{\lfloor\frac{n}{2}\rfloor}\card C_{1^{n-2k}\,2^{k}}=\sum_{k=0}^{\lfloor\frac{n}{2}\rfloor}\frac{n!}{z_{1^{n-2k}\,2^{k}}}=\sum_{k=0}^{\lfloor\frac{n}{2}\rfloor}\frac{n!}{k!\, n-2k!\,2^{k}}\,.$$
Dans le cas général, la permutation $\sigma_{\mu}$ n'a pas forcément de racines carrées. En effet, si $\tau$ est de type cyclique $1^{k_{1}}\,2^{k_{2}}\,\cdots\,s^{k_{s}}$, alors :\vspace{2mm}
\begin{enumerate}
\item Le carré de tout cycle impair est un cycle impair de même longueur, donc les parts impaires $1^{k_{1}}\,3^{k_{3}}\,5^{k_{5}}\,\cdots$ sont conservées par l'élévation au carré $\tau \mapsto \tau^{2}$. \vspace{2mm}
\item Le carré de tout cycle pair de longueur $2l$ est un produit de deux $l$-cycles disjoints, donc les parts paires $2^{k_{2}}\,4^{k_{4}}\,6^{k_{6}}\,\cdots$ deviennent $1^{2k_{2}}\,2^{2k_{4}}\,3^{2k_{6}}\,\cdots$ par élévation au carré.\vspace{2mm}
\end{enumerate}
On conclut que le type cyclique de $\tau^{2}$ est :
$$1^{k_{1}+2k_{2}}\,2^{2k_{4}}\,3^{k_{3}+2k_{6}}\,4^{2k_{8}}\,\cdots\,(2l+1)^{k_{2l+1}+2k_{4l+2}}\,(2l+2)^{2k_{4l+4}}\,\cdots$$
Ainsi, $\sigma_{\mu}$ admet des racines carrées si et seulement si $m_{2i}$ est pair pour tout entier pair $2i$. Dans ce cas, pour choisir une racine carrée de $\sigma_{\mu}$, il faut :\vspace{2mm}
\begin{enumerate}
\item Pour chaque entier impair $i$, choisir les entiers $k_{i}$ et $k_{2i}$ tels que $m_{i}=k_{i}+2k_{2i}$. Puis, choisir un appariement partiel de $2k_{2i}$ des $i$-cycles de $\tau$, et pour chaque paire $(c_{1},c_{2})$ de cet appariement, choisir un $2i$-cycle dont le carré est $c_{1}c_{2}$ ; il y a à chaque fois $i$ cycles de ce type.\vspace{2mm}
\item Pour chaque entier pair $2i$, choisir simplement un appariement complet des $m_{2i}$-cycles de longueur $2i$ de $\tau$, et pour chaque paire $(c_{1},c_{2})$ de cet appairement, choisir un $4i$-cycle dont le carré est $c_{1}c_{2}$ ; il y a $2i$ cycles de ce type.\vspace{2mm}
\end{enumerate}
Le nombre d'appariements de $2k$ éléments parmi $i$ est $\frac{i!}{k!\,i-2k!\,2^{k}}$ ; on en déduit la formule suivante pour le nombre de racines carrées d'une permutation dans $\sym_{n}$.
\begin{proposition}[Nombre de racines carrées d'une permutation]
Si $\sigma$ a pour type cyclique $\mu=1^{m_{1}}\,2^{m_{2}}\,\cdots\,s^{m_{s}}$, alors $\card\{\tau \in \sym_{n}\,\,|\,\,\tau^{2}=\sigma\}=\prod_{i=1}^{s }f(i,m_{i})$, avec :
$$f(i,m)=\begin{cases}0&\text{si }i\text{ est pair et }m\text{ est impair},\\
\frac{m!}{m/2!}\,\left(\frac{i}{2}\right)^{m/2}&\text{si }i\text{ et }m\text{ sont pairs},\\
\sum_{k=0}^{\lfloor \frac{m}{2}\rfloor}\frac{m!}{m-2k!\,k!}\, \left(\frac{i}{2}\right)^{k}&\text{si }i\text{ est impair}.\end{cases}$$ 
\end{proposition}
\noindent Notons que $f(1,n)$ est bien le nombre d'involutions $I_{n}$ pour tout entier $n$.\medskip

\begin{example}
Déterminons le nombre de racines carrées de $\sigma=(1)(2)(3)(4,5)(6,7)(8,9,10)$ dans $\sym_{10}$. Si $\tau^{2}=\sigma$, alors les trois points fixes $1,2,3$ de $\sigma$ correspondent soit à une composante $(1)(2)(3)$ dans $\tau$, soit à l'une des transpositions $(1,2)(3)$, $(1,3)(2)$ et $(1)(2,3)$ ; ceci donne $f(1,3)=4$ choix possibles. Le produit de transpositions $(4,5)(6,7)$ provient de l'un des deux $4$-cycles $(4,6,5,7)$ et $(4,7,5,6)$, d'où $f(2,2)=2$ possibilités. Finalement, le $3$-cycle $(8,9,10)$ dans $\sigma$ provient nécessairement du $3$-cycle $(8,10,9)$ dans $\tau$, d'où $f(3,1)=1$ possibilité. Ainsi, le nombre de racines carrées de $\sigma$ dans $\sym_{10}$ est
$$f(1,3)\,f(2,2)\,f(3,1)=8\,.$$
\end{example}
\bigskip\bigskip

Fixons alors un caractère central $\varSigma_{\mu}$ avec $|\mu|=k$, et un entier $n \geq k$. L'espérance du caractère central sous la mesure de Gelfand $\Gel_{n}$ s'écrit :
\begin{align*}\Gel_{n}[\varSigma_{\mu}]&=n^{\downarrow |\mu|} \,\Gel_{n}[\chi^{\lambda}(\mu\sqcup 1^{n-|\mu|})]=n^{\downarrow |\mu|} \,\Gel_{n}[\chi^{\lambda}(1^{n-|\mu|+m_{1}}\,2^{m_{2}}\,\cdots \,s^{m_{s}})]\\
&=n^{\downarrow |\mu|} \,\frac{\card\{\tau \in \sym_{n}\,\,|\,\,\tau^{2}=\sigma_{\mu}\}}{I_{n}}=n^{\downarrow |\mu|}\,\frac{I_{n-|\mu|+m_{1}}}{I_{n}}\,\left(\prod_{i=2}^{s}f(i,m_{i})\right).
\end{align*}
Il est alors important de connaître un équivalent de $I_{n}$ lorsque $n$ tend vers l'infini. Pour commencer, déterminons la série génératrice exponentielle des nombres $I_{n}$. Elle s'écrit :
$$G(z)=\sum_{n=0}^{\infty}\frac{I_{n}}{n!}\,z^{n}=\sum_{n=0}^{\infty}\sum_{k=0}^{\lfloor \frac{n}{2}\rfloor}\frac{1}{n-2k!\,k!}\,z^{n-2k}\,\left(\frac{z^{2}}{2}\right)^{k}=\sum_{p=0}^{\infty}\sum_{k=0}^{\infty}\frac{1}{p!\,k!}\,z^{p}\,\left(\frac{z^{2}}{2}\right)^{k}=\exp\left(z+\frac{z^{2}}{2}\right)\,.$$
En écrivant $n!=\int_{0}^{\infty} \E^{-t}\,t^{n}\,dt=\int_{0}^{\infty} \E^{-t+n\,\log t}\,\,dt$ et en utilisant la méthode de Laplace, on peut démontrer le célèbre équivalent de Stirling
$$n! \simeq_{n \to \infty}\left(\frac{n}{\E}\right)^{n}\,\sqrt{2\pi n}\,.$$
D'autre part, étant donnée une série génératrice de rayon de convergence non nulle, il existe de nombreuses techniques analytiques donnant l'asymptotique des coefficients de la série, voir en particulier \cite[chapitre 8]{FS09} et \cite[chapitre 5]{Wil06}. Lorsque le rayon de convergence est infini, il convient d'utiliser la \textbf{méthode d'Hayman}, voir \cite[\S5.4]{Wil06}. Nous la détaillons ci-après dans le cas où $G(z)=\exp(z+z^{2}/2)$, voir \cite[chapitre 8, p. 559]{FS09}. Pour tout rayon $R$, on peut retrouver $a_{n}=I_{n}/n!$ par la formule de Cauchy : 
$$a_{n}=\frac{1}{2\pi R^{n}}\int_{-\pi}^{\pi} G(R\E^{\I \theta})\,\E^{-n \I \theta}\,d\theta=\frac{1}{2\pi R^{n}}\int_{-\pi}^{\pi} \exp\left(R\E^{\I \theta}+\frac{R^2}{2}\,\E^{2\I \theta}-n \I \theta\right)d\theta\,.$$
Notons $f(\theta)$ la fonction dans l'intégrande. Le module $|f(\theta)|$ est égal à :
$$\exp\left(R\cos\theta + \frac{R^{2}}{2}\cos 2\theta\right) \leq \exp\left(R+\frac{R^{2}}{2}\right)=f(R)\,.$$
Le module maximum vaut donc $m(R)=f(R)$, et $m(R)/R^{n}$ s'écrit :
$$\frac{m(R)}{R^{n}}=\exp\left(R+\frac{R^{2}}{2}-n\log R\right).$$
À $n$ fixé, $m(R)$ est minimum lorsque $R+R^{2}=n$, c'est-à-dire lorsque $R=\frac{\sqrt{1+4n}-1}{2}\simeq \sqrt{n}$. Dans ce qui suit, nous supposerons que $R$ a cette valeur. Le développement limité de $f(R\E^{\I\theta})$ autour de $\theta=0$ s'écrit alors :
\begin{align*}f(R\E^{\I \theta})&=\exp\left(R+R\I\theta-\frac{R\theta^{2}}{2}+R\,o(\theta^{2})+\frac{R^{2}}{2}+R^{2}\I\theta-R^{2}\theta^{2}+R^{2}\,o(\theta^{2})-n \I\theta\right)\\
&=m(R)\,\exp\left(-\left(\frac{R}{2}+R^{2}\right)\,(1+o(1))\,\theta^{2}\right).\end{align*}
Notons que la fonction $\eps(\theta)=o(1)$ est une fonction de $\theta$ qui converge vers $0$ uniformément en $R$ : ainsi,
$$\forall \eps>0,\,\, \exists \eta>0,\,\,\forall R,\,\,|\theta|\leq \eta \Rightarrow |\eps(\theta)| \leq \eps\,.$$
En découpant l'intégrale en trois parties correspondant à des intervalles $[-\pi,-\eta]$, $[-\eta,\eta]$ et $[\eta,\pi]$, on voit donc que :
\begin{align*}\bigg| a_{n}-\frac{m(R)}{2\pi\,R^{n}}\,\int_{-\eta}^{\eta} \exp\left(-\left(\frac{R}{2}+R^{2}\right)(1+\eps(\theta))\,\theta^{2}\right)&d\theta\bigg|\leq C\,\frac{m(R)}{R^{n}}\,\exp\left(R\,\cos \eta + \frac{R^{2}}{2}\,\cos 2\eta\right)\\
&\leq \frac{C}{R^{n}}\,\exp\left(-\left(\frac{R}{2}+R^{2}\right)\eta^{2}+R^{2}\eta^{2}\,\eps(\eta)\right).\end{align*}
Prenons alors $\eta=n^{-1/3}$ ; notons que $\eta$ tend vers zéro, mais $R\eta$ tend vers l'infini ; on a aussi $R\eta^{2}$ qui tend vers zéro. Le majorant est alors plus petit que 
$$\frac{C'}{R^{n}}\,\exp\left(-(R\eta)^{2}+(R\eta)^{2}\,\eps(\eta)\right)\,,$$
et l'intégrale restante est équivalente à 
$$\frac{m(R)}{2\pi\,R^{n+1}}\,\int_{-\infty}^{\infty}\E^{-\nu^{2}}\,d\nu=\frac{m(R)}{2\sqrt{\pi}\,R^{n+1}}=\frac{\exp(R+\frac{R^{2}}{2})}{2\sqrt{\pi}\,R^{n+1}}\,. $$
Ainsi, $a_{n}$ est équivalent à $\frac{1}{2\sqrt{\pi}\,R^{n+1}}\,\exp(R+\frac{R^{2}}{2})$, et comme $R=\frac{\sqrt{1+4n}-1}{2}$, on obtient :
$$a_{n} \simeq_{n \to \infty} \frac{1}{\sqrt{2}}\,\left(\frac{\E}{n}\right)^{\frac{n}{2}}\,\frac{\exp(\sqrt{n}-1/4)}{\sqrt{2\pi n}}\qquad;\qquad I_{n}\simeq_{n \to \infty} \frac{1}{\sqrt{2}}\,\left(\frac{n}{\E}\right)^{\frac{n}{2}}\,\exp(\sqrt{n}-1/4)\,.$$
En injectant cet équivalent dans l'expression précédemment donnée pour $\Gel_{n}[\varSigma_{\mu}]$, on conclut que :
\begin{proposition}[Asymptotique des espérances des caractères centraux]\label{asymptoticobsgelfand}
Lorsque $n$ tend vers l'infini,
$$\Gel_{n}[\varSigma_{\mu}]\simeq \left(\prod_{i=2}^{s}f(i,m_{i})\right)\,n^{\frac{|\mu|+m_{1}(\mu)}{2}}\,.$$
\end{proposition}
\noindent En effet, le terme $n^{\downarrow |\mu|}$ peut être remplacé par $n^{|\mu|}$, et 
$$\frac{I_{n-k}}{I_{n}}\simeq \left(\frac{n-k}{n}\right)^{\frac{n-k}{2}}\,\left(\frac{\E}{n}\right)^{\frac{k}{2}}\,\E^{\sqrt{n-k}-\sqrt{n}}\simeq \frac{1}{n^{\frac{k}{2}}}\,.$$
avec $k=|\mu|-m_{1}(\mu)$. Cette expression asymptotique va permettre d'appliquer notre méthode usuelle de moments ; notons que pour toute observable de diagrammes $f \in \obs$, on a $\Gel_{n}[f]=O(n^{\deg_{\mathrm{K}}(f)/2})$, car ceci est vrai sur la base linéaire des caractères centraux --- rappelons que le degré de Kerov des observables est défini par $\deg_{\mathrm{K}}(\varSigma_{\mu})=|\mu|+m_{1}(\mu)$ pour toute partition $\mu$. Comme le degré de Kerov est toujours inférieur au poids des observables, on a aussi $\Gel_{n}[f]=O(n^{\mathrm{wt}(f)/2})$ pour toute observable $f$.
\bigskip

\section{Asymptotique des caractères centraux et des formes des diagrammes}\label{asymptoticgelfand}
La proposition \ref{asymptoticobsgelfand} permettra d'établir la convergence des formes des diagrammes sous les mesures de Gelfand vers la courbe $\Omega$ du théorème \ref{firstasymptoticplancherel}. D'autre part, elle permet de déterminer la distribution asymptotique précise des caractères centraux des cycles $\varSigma_{k}$. Sous les mesures de Plancherel, rappelons qu'on a asymptotiquement
$$\forall k\geq 2,\,\,\frac{\varSigma_{k}(\lambda)}{n^{k/2}}\simeq \sqrt{k}\,\xi_{k}\,,$$
où les $\xi_{k}$ sont des gaussiennes indépendantes centrées de variance $1$. Le théorème suivant est l'analogue de ce résultat pour les mesures de Gelfand :
\begin{theorem}[Distribution asymptotique des caractères centraux des cycles sous les mesures de Gelfand, \cite{Mel10}]\label{nolove}
Sous les mesures de Gelfand, le vecteur des caractères centraux renormalisés $(\varSigma_{k}(\lambda)/n^{k/2})_{k \geq 2}$ converge en lois fini-dimensionnelles vers un vecteur gaussien $(\sqrt{2k}\,\xi_{k}+e_{k})_{k \geq 2}$, où les $\xi_{k}$ sont des gaussiennes standards indépendantes, et $e_{k}$ vaut $0$ si $k$ est pair, et $1$ si $k$ est impair. Ainsi,
$$\lim_{n \to \infty }\frac{\varSigma_{k}(\lambda)}{n^{k/2}} \simeq \mathcal{N}(e_{k},2k)\,.$$
\end{theorem}
\begin{lemma}[Degré de Kerov et produits de caractères centraux]\label{rageofposeidon}
Si $k$ et $l$ sont deux entiers distincts plus grands que $2$, alors $\varSigma_{k,l}$ est la composante de plus haut degré de Kerov de $\varSigma_{k}\,\varSigma_{l}$. D'autre part, si $k \geq 2$ et $m \geq 1$, alors la composante de plus haut degré de Kerov de $(\varSigma_{k})^{m}$ est 
$$\sum_{p=0}^{\lfloor \frac{m}{2}\rfloor} \frac{m!}{m-2p!\,p!}\,\left(\frac{k}{2}\right)^{p}\,\varSigma_{1^{kp}\,k^{m-2p}}\,.$$
\end{lemma}
\begin{proof}
La première partie du lemme est démontrée dans \cite[proposition 4.13]{IO02}; la seconde partie est essentiellement équivalente aux propositions 4.11, 4.12, 6.2 et 6.3 du même article, et peut être démontrée de manière purement combinatoire en utilisant l'algèbre des permutations partielles. Pour commencer, remarquons que le degré de Kerov défini sur $\obs \simeq \alg_{\infty}$ peut être relevé à l'algèbre $\blg_{\infty}$ des permutations partielles en posant :
$$\deg_{\mathrm{K}}(\sigma,S)=\card(\mathrm{Fix}(\sigma)\cap S) + \card S\,.$$
Cette définition fournit bien une filtration d'algèbre sur $\blg_{\infty}$, car 
\begin{align*}\deg_{\mathrm{K}}((\sigma,S)(\tau,T))&=\deg_{\mathrm{K}}(\sigma\tau,S\cup T)\\
&=\big\{\card(\mathrm{Fix}(\sigma\tau)\cap (S\cup T))-\card S\cap T\big\}+\card S+\card T\,,\end{align*}
et le terme entre crochets est plus petit que $\card (\mathrm{Fix}(\sigma) \cap S)+\card(\mathrm{Fix}(\tau)\cap T)$. En effet, si $x \in S \cap T$ vérifie $\sigma\tau(x)=x$, alors :
\begin{enumerate}
\item L'élément $x$ peut être dans $S\setminus T$, mais dans ce cas $\tau(x)=x$, donc $\sigma(x)=x$, et ainsi $x \in \mathrm{Fix}(\sigma)\cap S$.
\item De même, $x$ peut être dans $T \setminus S$, et dans ce cas il est aussi dans $\mathrm{Fix}(\tau)\cap T$.
\item Enfin, si les deux hypothèses précédentes ne sont pas vérifiées, alors $x$ appartient à $S\cap T$.
\end{enumerate}
On conclut que $\card(\mathrm{Fix}(\sigma\tau)\cap (S\cup T) )\leq \card (\mathrm{Fix}(\sigma)\cap S) + \card(\mathrm{Fix}(\tau)\cap T)+\card(S\cap T)$, de sorte que $\deg_{\mathrm{K}}((\sigma,S)(\tau,T))\leq \deg_{\mathrm{K}}(\sigma,S)+\deg_{\mathrm{K}}(\tau,T)$. En tant qu'élément de l'algèbre des permutations partielles, 
$$(\varSigma_{k})^{m}=\sum (c_{1},S_{1})(c_{2},S_{2})\cdots(c_{m},S_{m})$$
où la somme est prise sur les $m$-tuples de $k$-arrangements $(a_{i1},\ldots,a_{ik})$ donnant un cycle $c_{i}=(a_{i1},\ldots,a_{ik})$ et un support  $S_{i}=\{a_{i1},\ldots,a_{ik}\}$. Comme $\deg_{\mathrm{K}}(\varSigma_{k})=k$, $\deg_{\mathrm{K}}((\varSigma_{k})^{m}) \leq km$, et cette inégalité est une égalité, car $km$ est le degré de Kerov d'un produit  $(c_{1},S_{1})\cdots(c_{m},S_{m})$ de $k$-cycles disjoints. Ceci étant, déterminons tous les produits qui ont ce degré maximal. Si $\deg_{\mathrm{K}}((c_{1},S_{1})\cdots(c_{m},S_{m}))=km$, alors comme le degré de Kerov est une filtration d'algèbre, ceci impose :
\begin{align*}
\deg_{\mathrm{K}}((c_{1},S_{1})(c_{2},S_{2}))&=2k\\
\deg_{\mathrm{K}}((c_{1},S_{1})(c_{2},S_{2})(c_{3},S_{3}))&=3k\\
\vdots\qquad\qquad\qquad\qquad & \vdots\\
\deg_{\mathrm{K}}((c_{1},S_{1})(c_{2},S_{2})\cdots (c_{i},S_{i}))&=ik\\
\vdots\qquad\qquad\qquad\qquad & \vdots
\end{align*}
Considérons les deux premiers termes $(c_{1},S_{1})$ et $(c_{2},S_{2})$. On a montré que :
\begin{align*}2k&=\deg_{\mathrm{K}}((c_{1},S_{1})(c_{2},S_{2}))=\big\{\card\!((\mathrm{Fix}(c_{1}c_{2})\cap (S_{1}\cup S_{2})) -\card S_{1}\cap S_{2}\big\}+\card S_{1}+\card S_{2}\\
&=2k+\big\{\card((\mathrm{Fix}(c_{1}c_{2})\cap (S_{1}\cup S_{2})) -\card S_{1}\cap S_{2}\big\}\,.\end{align*}
On a donc $\card(\mathrm{Fix}(c_{1}c_{2})\cap(S_{1}\cup S_{2}))=\card S_{1}\cap S_{2}$, et ceci impose que tout point de l'intersection $S_{1}\cap S_{2}$ soit un point fixe de $c_{1}c_{2}$. Soit $u$ l'exposant de $c_{1}c_{2}$, de sorte que $(c_{1}c_{2})^{u}=\id$. On raisonne par l'absurde en supposant que $S_{1}$ et $S_{2}$ ont une intersection non vide, mais sont deux parties de cardinal $k$ qui ne sont pas égales. Notons alors  $c_{1}=(a_{1},\ldots,a_{k})$ et $c_{2}=(b_{1},\ldots,b_{k})$ ; par hypothèse, l'un des $b_{j}$ est dans $S_{2}$, mais pas dans $S_{1}$. On itère $c_{1}c_{2}$ sur cet élément. Si $b_{j+1}$  n'est pas dans $S_{1}$, alors $c_{1}c_{2}(b_{j})=c_{1}(b_{j+1})=b_{j+1}$, et on peut recommencer avec $b_{j+1}$, $b_{j+2}$, \emph{etc.} On construit ainsi une chaîne $b_{j},b_{j+1}, \ldots,b_{J-1}$ dont tous les éléments sont dans $S_{2}\setminus S_{1}$, et vérifient donc $c_{1}c_{2}(b_{p})=b_{p+1}$. Comme $S_{1}\cap S_{2} \neq \emptyset$, il doit exister un $b_{J}$ qui est dans $S_{1}\cap S_{2}$. Alors, $a_{i}=c_{1}(b_{J})$ est un certain $(c_{1}c_{2})^{s}(b_{j})$, et est dans $S_{1}$. Il ne peut pas être dans $S_{1}\cap S_{2}$ : sinon, ce serait un point fixe de $c_{1}c_{2}$, et tous les $(c_{1}c_{2})^{t}(a_{i})$ seraient égaux à $a_{i}$, en particulier $(c_{1}c_{2})^{u-s}(a_{i})=(c_{1}c_{2})^{u}(b_{j})=b_{j}$. Ceci est absurde car $b_{j}$ n'est pas dans $S_{1}$. Ainsi, $a_{i} \in S_{1}\setminus S_{2}$, et en continuant à itérer $c_{1}c_{2}$, on obtient $a_{i+1}$, $a_{i+2}$, \emph{etc.} jusqu'à un élment $a_{I}$ qui est dans $S_{1}$ et dans $S_{2}$ (voir la figure \ref{absurd}). Cet élément $a_{I}$ est un certain $(c_{1}c_{2})^{s'}(b_{j})$, et est un point fixe de $c_{1}c_{2}$. Alors, en prenant $a_{I}=(c_{1}c_{2})^{t'}(a_{I})$ avec $t'=u-s'$, on retrouve $b_{j}$, ce qui est absurde car  $b_{j}$ n'est pas dans $S_{1}\cap S_{2}$.
\figcapt{\psset{unit=1mm}\pspicture(0,0)(120,33)
\rput(0,30){$b_{j}$}
\rput(0,25){$S_{2}\setminus S_{1}$}
\psline{->}(4,30)(16,30)
\rput(10,28){$c_{1}c_{2}$}
\rput(20,30){$b_{j+1}$}
\rput(20,25){$S_{2}\setminus S_{1}$}
\psline{->}(24,30)(36,30)
\rput(30,28){$c_{1}c_{2}$}
\rput(40,30){$\cdots$}
\psline{->}(44,30)(56,30)
\rput(50,28){$c_{1}c_{2}$}
\rput(60,30){$b_{J-1}$}
\rput(60,25){$S_{2}\setminus S_{1}$}
\psline{->}(64,30)(76,30)
\rput(70,28){$c_{2}$}
\rput(80,30){$b_{J}$}
\rput(80,25){$S_{1} \cap S_{2}$}
\psline{->}(80,20)(80,10)
\rput(77,15){$c_{1}$}
\rput(80,5){$a_{i}$}
\rput(80,0){$S_{1} \setminus S_{2}$}
\psline{->}(84,5)(96,5)
\rput(90,3){$c_{1}c_{2}$}
\rput(100,5){$\cdots$}
\psline{->}(104,5)(116,5)
\rput(110,3){$c_{1}c_{2}$}
\rput(120,5){$a_{I}$}
\rput(120,0){$S_{1} \cap S_{2} \,???$}
\endpspicture}{Si $S_{1} \neq S_{2}$ et $S_{1}\cap S_{2}\neq \emptyset$, alors l'identité $\deg_{\mathrm{K}}((c_{1},S_{1})(c_{2},S_{2}))=2k$ fournit une contradiction. Les produits de $k$-cycles de degré de Kerov maximal $2k$ sont donc ceux pour lesquels les supports sont disjoints ou identiques.\label{absurd}}{Produits de $k$-cycles de degré de Kerov $2k$ à supports disjoints ou identiques}

Ainsi, si $(c_{1},S_{1})(c_{2},S_{2})$ a pour degré de Kerov $2k$, alors soit les cycles sont disjoints, soit $S_{1}=S_{2}$. Dans ce cas, le nombre de points fixes de $c_{1}c_{2}$ doit être égal à $k$, et ceci impose $c_{1}=c_{2}^{-1}$. Ainsi, $c_{2}$ doit être soit disjoint de $c_{1}$, soit égal à $c_{1}^{-1}$. Pour les mêmes raisons, $c_{3}$ doit être soit disjoint de $c_{1}$ et de $c_{2}$, soit égal à l'inverse d'un de ces cycles, en supposant dans ce cas qu'on n'a pas déjà $c_{2}=c_{1}^{-1}$. Par récurrence sur $m$, on conclut qu'un cycle $c_{i}$ mis en jeu dans un produit de degré de Kerov $km$ est soit disjoint de tous les autres cycles, soit apparié avec un (et exactement un) autre cycle dont il est l'inverse. Le type cyclique d'un produit de degré de Kerov $km$ est donc toujours de la forme $1^{kp}\,k^{m-2p}$, où $p$ est le nombre de paires de cycles inverses l'un de l'autre. L'entier $p$ étant fixé, il y aura un facteur $$\frac{m!}{m-2p!\,p!\,2^{p}}$$
devant le terme $\varSigma_{1^{kp}\,k^{m-2p}}$  correspondant au nombre d'appariements partiels avec $p$ paires ; et aussi un facteur $k^{p}$, car étant donné un arrangement $(a_{i1},\ldots,a_{ik})$ donnant un cycle $c_{i}$, il y a exactement $k$ arrangements distincts  $(a_{j1},\ldots,a_{jk})$ correspondant au cycle inverse (les $k$ permutations cycliques d'une écriture du cycle inverse). On conclut que le terme de plus haut degré de Kerov de $(\varSigma_{k})^{m}$ est bien 
$$\sum_{p=0}^{\lfloor \frac{m}{2}\rfloor} \frac{m!}{m-2p!\,p!}\,\left(\frac{k}{2}\right)^{p}\,\varSigma_{1^{kp}\,k^{m-2p}}\,.\vspace{-5mm}$$
\end{proof}
\begin{example}
$(\varSigma_{4})^{4}=\varSigma_{4^{4}}+24\,\varSigma_{1^{4}\, 4^{2}}+48\,\varSigma_{1^{8}}+(\text{termes de degré strictement inférieur à }16)$.
\end{example}\bigskip

\begin{proof}[Preuve du théorème \ref{nolove}, première partie]
Soit $k$ un entier pair, et $m$ un entier impair. Alors, le terme $\varSigma_{1^{kp}\,k^{m-2p}}$ mis en jeu dans le développement de $(\varSigma_{k})^{m}$ a une espérance nulle sous la mesure de Gelfand, car $m-2p$ est impair, et les permutations qui ont ce type cyclique n'ont donc pas de racines carrées. Par conséquent, $\Gel_{n}[(\varSigma_{k})^{m}]$ est dans ce cas égale à l'espérance d'une observable de degré de Kerov inférieur à $km-1$, donc :
$$\Gel_{n}\left[\left(\frac{\varSigma_{k}}{n^{k/2}}\right)^{m}\right]=\frac{O(n^{km-1/2})}{n^{km/2}}=O(n^{-1/2})\to 0\,.$$
Supposons maintenant $k$ et $m$ pairs. Alors, les espérances des termes de degré de Kerov inférieur à $km-1$ sont négligeables dans l'asymptotique de $\Gel_{n}\left[\left(\frac{\varSigma_{k}}{n^{k/2}}\right)^{m}\right]$, donc 
\begin{align*}
\Gel_{n}\left[\left(\frac{\varSigma_{k}}{n^{k/2}}\right)^{m}\right]&\simeq \sum_{p=0}^{\lfloor \frac{m}{2}\rfloor} \frac{m!}{m-2p!\,p!}\,\left(\frac{k}{2}\right)^{p}\,\Gel_{n}\left[\frac{\varSigma_{1^{kp}\,k^{m-2p}} }{n^{km/2}}\right]=\sum_{p=0}^{\lfloor \frac{m}{2}\rfloor} \frac{m!}{m-2p!\,p!}\,\left(\frac{k}{2}\right)^{p}\,f(k,m-2p)\\
&\simeq\left(\frac{k}{2}\right)^{m/2}\,\sum_{p=0}^{\lfloor \frac{m}{2}\rfloor} \frac{m!}{m/2-p!\,p!}=k^{m/2}\,\frac{m!}{m/2!}=(2k)^{m/2}\,(m-1!!)\,.
\end{align*}
Mais on a vu dans le chapitre \ref{matrix} que les doubles factorielles étaient justement les moments pairs d'une variable gaussienne standard centrée, les moments impairs étant pour leur part nuls. Comme les variables gaussiennes sont caractérisées par leurs moments, on conclut à la convergence en loi 
$$\frac{\varSigma_{k}}{n^{k/2}} \to \mathcal{N}(0,2k) $$
sous les mesures de Gelfand. Le cas où $k$ est impair peut être traité de fa\c con tout à fait similaire. Pour les mêmes raisons que précédemment, on peut négliger les termes qui ne sont pas de degré de Kerov maximal dans l'asymptotique des moments de $\varSigma_{k}$, et ainsi,
\begin{align*}
\Gel_{n}\left[\left(\frac{\varSigma_{k}}{n^{k/2}}\right)^{m}\right]&\simeq \sum_{p=0}^{\lfloor \frac{m}{2}\rfloor} \frac{m!}{m-2p!\,p!}\,\left(\frac{k}{2}\right)^{p}\,\Gel_{n}\left[\frac{\varSigma_{1^{kp}\,k^{m-2p}} }{n^{km/2}}\right]=\sum_{p=0}^{\lfloor \frac{m}{2}\rfloor} \frac{m!}{m-2p!\,p!}\,\left(\frac{k}{2}\right)^{p}\,f(k,m-2p)\\
&\simeq \sum_{p=0}^{\lfloor \frac{m}{2}\rfloor}\sum_{q=0}^{\lfloor \frac{m-2p}{2}\rfloor} \frac{m!}{p!\,q!\,m-2(p+q)!}\,\left(\frac{k}{2}\right)^{p+q}=\sum_{r=0}^{\lfloor \frac{m}{2}\rfloor} \frac{m!}{m-2r!} \,\left(\frac{k}{2}\right)^{r}\,\left(\sum_{p+q=r}\frac{1}{p!\,q!}\right)\\
&\simeq \sum_{r=0}^{\lfloor\frac{m}{2}\rfloor} \frac{m!}{m-2r!\,r!}\,k^{r}\,.
\end{align*}
La fonction génératrice correspondant à ces moments limites est donnée par :
$$\sum_{m=0}^{\infty} \sum_{r=0}^{\lfloor\frac{m}{2}\rfloor} \frac{1}{m-2r!\,r!}\,k^{r}\,z^{m}=\sum_{m=0}^{\infty} \sum_{r=0}^{\lfloor\frac{m}{2}\rfloor} \frac{1}{m-2r!\,r!}\,(kz^{2})^{r}\,z^{m-2r}=\exp\left(z+kz^{2}\right),$$
et c'est donc la fonction génératrice $\esper[\E^{zX}]$ d'une variable gaussienne de loi $\mathcal{N}(1,2k)$. Ainsi, on a la convergence en loi
$$\frac{\varSigma_{k}}{n^{k/2}} \to \mathcal{N}(1,2k) $$
sous les mesures de Gelfand, et on a établi la convergence (pour le moment non jointe) de chaque caractère central renormalisé vers une variable gaussienne. \end{proof}\bigskip

Pour démontrer l'indépendance asymptotique des caractères centraux renormalisés, on doit de nouveau utiliser la théorie des cumulants d'observables. Ainsi, nous allons montrer l'annulation asymptotique des cumulants d'ordre plus grand que $3$ :
\begin{lemma}[Convergence jointe des caractères centraux renormalisés]\label{spacebound}
$$\forall r\geq 3,\,\,\,\forall l_{1},\ldots,l_{r} \geq 2,\,\,\,k\left(\frac{\varSigma_{l_{1}}}{n^{l_{1}/2}},\frac{\varSigma_{l_{2}}}{n^{l_{2}/2}},\ldots,\frac{\varSigma_{l_{r}}}{n^{l_{r}/2}}\right)=o(1)\,.$$
\end{lemma}
\noindent Comme les $\varSigma_{l_{i}}$ sont de degré de Kerov $l_{i}$, on sait déjà que l'expression précédente est un $O(1)$. En effet : 
$$k(\varSigma_{l_{1}},\ldots,\varSigma_{l_{r}})=\sum_{\pi \in \mathfrak{Q}(\lle 1,r\rre)} \mu(\pi)\,\prod_{\pi_{j} \in \pi}\esper\left[\prod_{i \in \pi_{j}} \varSigma_{l_{i}}\right]  = \sum_{\pi \in \mathfrak{Q}(\lle 1,r\rre)} \prod_{\pi_{j} \in \pi} O\left(n^{\frac{\sum_{i \in \pi_{j}} l_{i}}{2}}\right) = O\left(n^{\frac{l_{1}+\cdots+l_{r}}{2}}\right).$$
On doit donc seulement gagner un ordre de grandeur. Rappelons que le cumulant standard peut être exprimé en fonction des cumulants disjoints et des cumulants identité :
$$k(\varSigma_{l_{1}},\ldots,\varSigma_{l_{r}})=\sum_{\pi \in \mathfrak{Q}(\lle 1,r\rre)} k^{\bullet}\left(k^{\id}(\varSigma_{l_{i} \in \pi_{1}}), \ldots, k^{\id}(\varSigma_{l_{i} \in \pi_{s}}) \right).$$
\begin{lemma}[Degré de Kerov d'un cumulant identité de caractère centraux]\label{blopblop}
Soient $l_{1},\ldots,l_{t}$ des entiers plus grands que $2$. Si ces entiers ne sont pas tous égaux, alors :
$$\deg_{\mathrm{K}}\left(k^{\id}(\varSigma_{l_{1}},\ldots,\varSigma_{l_{t}})\right)\leq l_{1}+\cdots+l_{t}-1\,.$$
Si ces entiers sont tous égaux et si $t\geq 3$, alors 
$$\deg_{\mathrm{K}}\left(k^{\id}(\varSigma_{l},\ldots,\varSigma_{l})\right)\leq tl-1\,.$$
\end{lemma}
\begin{proof}
Le lemme \ref{sanitarium} donne une description explicite des produits de cycles intervenant dans le cumulant identité :
$$k^{\id}(\varSigma_{l_{1}},\ldots,\varSigma_{l_{t}})= \sum_{\substack{\forall j,\,\,A_{j} \in \arra(l_{j})\\ \pi(A_{1},\ldots,A_{t})=\{\lle 1,t \rre\}}} C(A_{1})\,C(A_{2})\,\cdots\,C(A_{t})\,.$$
où $\arra(l_{j})$ est l'ensemble des $l_{j}$-arrangements d'entiers, et où la restriction $\pi(A_{1},\ldots,A_{t})=\{\lle 1,t \rre\}$ impose que ces arrangements <<~se recouvrent les uns les autres~>>. Dans le premier cas, comme les $l_{i}$ ne sont pas tous égaux, quitte à les permuter, on peut trouver un entier $u<t$ tels que les $l_{1},\ldots,l_{u}$ soient identiques, et distincts de tous les $l_{u+1},\ldots,l_{t}$. Considérons alors un produit de cycles $C(A_{1})\,\cdots\,C(A_{t})$ apparaissant dans la somme. Si le produit partiel $C(A_{1})\,\cdots\,C(A_{u})$ n'est pas de degré de Kerov maximal $l_{1}+\cdots+l_{u}=u\,l_{1}$, alors comme le degré de Kerov est une filtration sur l'algèbre des permutations partielles, on a bien 
$$\deg_{\mathrm{K}}(C(A_{1})\,\cdots\,C(A_{t}))\leq l_{1}+\cdots+l_{t}-1\,.$$
Dans le cas contraire, d'après le lemme \ref{rageofposeidon}, le type cyclique du produit partiel des $u$ premiers cycles est forcément un $l_{1}^{u-2k}\,1^{l_{1}k}$. D'autre part, il existe un indice $v \in \lle u+1,t\rre$ tel que le cycle $C(A_{v})$ ait son support intersectant le support $S$ du produit des $u$ premiers cycles ; prenons le premier d'entre eux. Notons que la longueur $l_{v}$ n'intervient pas dans le type type cyclique $l_{1}^{u-2k}\,1^{l_{1}k}$. Or, d'après la première partie du lemme \ref{rageofposeidon}, si deux permutations partielles ont des types cycliques sans parts communes, alors soit elles ont des supports disjoints, soit leur produit est de degré de Kerov strictement inférieur à la somme des degré de Kerov. Comme les supports s'intersectent, on conclut que le produit partiel $C(A_{1})\,\cdots\,C(A_{v})$ est de degré de Kerov inférieur à $l_{1}+\cdots+l_{v}-1$, et en rajoutant les derniers cycles, on obtient :
$$\deg_{\mathrm{K}}(C(A_{1})\,\cdots\,C(A_{t}))\leq l_{1}+\cdots+l_{t}-1\,.$$
Tous les produits de cycles intervenant dans la somme vérifient donc bien l'inégalité annoncée.\bigskip

Dans le second cas, considérant de même un produit de $l$-cycles $C(A_{1})\,\cdots\,C(A_{t})$ apparaissant dans la somme, comme $t\geq 3$, il y a forcément un arrangement $A_{i}$ qui intersecte deux autres arrangements $A_{j_{1}}$ et $A_{j_{2}}$. En effet, la condition de recouvrement des arrangements impose que le graphe de la relation $i \sim j \iff A_{i}\cap A_{j}\neq \emptyset$ soit connexe, donc ait plus de $t-1$ arêtes. Alors, la somme des valences des sommets $1,2,\ldots,t$ est plus grande que $2(t-1)>t$ pour $t\geq 3$, ce qui impose qu'un sommet ait une valence supérieure à $2$. Ainsi, fixons trois entiers $i,j_{1},j_{2}$ tels que $A_{i}\cap A_{j_{1}}\neq \emptyset$ et $A_{i}\cap A_{j_{2}}\neq \emptyset$ ; on ne perd pas de généralité en supposant que $i<j_{1}<j_{2}$, et que si $u \in \lle i+1,j_{1}-1\rre \cup \lle j_{1}+1,j_{2}-1\rre$, alors $A_{i}\cap A_{u}=\emptyset$ et $A_{j_{1}}\cap A_{u}=\emptyset$. Dans ces conditions, on peut permuter les cycles et regrouper le facteur $$c(A_{i})\,c(A_{j_{1}})\,c(A_{j_{2}})\,,$$ et il suffit de montrer que ce facteur est de degré inférieur à $3l-1$. Si le produit $c(A_{i})\,c(A_{j})$ est de degré de Kerov inférieur à $2l-1$, c'est évident. Sinon, comme $A_{i} \cap A_{j_{1}}\neq \emptyset$, la condition $\deg_{\mathrm{K}} (c(A_{i})\,c(A_{j_{1}}))=2l$ impose que les deux premiers cycles soient inverses l'un de l'autre, et le produit $c(A_{i})\,c(A_{j_{1}})$ est l'identité sur le support de $A_{i}$. Mais dans ce cas, comme $A_{i}\cap A_{j_{2}}\neq \emptyset$, notant $t\geq 1$ le cardinal de cette intersection, le produit $c(A_{i})\,c(A_{j_{1}})\,c(A_{j_{2}})$ est un $l$-cycle sur un support de taille $2l-t$, donc a un degré de Kerov égal à $(2l-t)+(l-t)=3l-2t<3l$. 
\end{proof}
\bigskip

\begin{lemma}[Formule d'inversion de M\"obius pour les partitions d'ensembles]\label{joeystarr}
Soit $F$ une fonction sur les couples d'entiers strictement positifs, et $(1,\ldots,1,2,\ldots,2,\ldots,s,\ldots,s)$ une suite de $r$ entiers avec $r_{1}\geq 1$ entiers $1$, $r_{2}\geq 1$ entiers $2$, \emph{etc}. Si $\pi \in \mathfrak{Q}(\lle 1,r\rre)$ est une partition d'ensemble de parts $\pi_{1}\sqcup \pi_{2}\sqcup \cdots \sqcup \pi_{\ell(\pi)}$, on note $r_{ij}$ le nombre d'entiers $i$ tombant dans $\pi_{j}$ ; en particulier, $|\pi_{j}|=\sum_{i=1}^{s} r_{ij}$ et $r_{i}=\sum_{j=1}^{\ell(\pi)} r_{ij}$.  Supposons $s \geq 2$. Alors,
$$S(r_{1},\ldots,r_{s})=\sum_{\pi \in \mathfrak{Q}(\lle 1,r\rre)} (-1)^{\ell(\pi)-1}\,(\ell(\pi)-1)!\,\prod_{j=1}^{\ell(\pi)}\prod_{r_{ij}\geq 1} F(i,r_{ij})=0\,.$$
\end{lemma}
\begin{example}
Supposons $r=4$. On a alors 15 partitions d'ensemble : $\{1,2,3,4\}$, $\{1,2,3\}\sqcup\{4\}$, $\{1,2,4\}\sqcup\{3\}$, $\{1,3,4\}\sqcup\{2\}$, $\{2,3,4\}\sqcup\{1\}$, $\{1,2\}\sqcup\{3,4\}$, $\{1,3\}\sqcup\{2,4\}$, $\{1,4\}\sqcup\{2,3\}$, $\{1,2\}\sqcup\{3\}\sqcup\{4\}$, $\{1,3\}\sqcup\{2\}\sqcup\{4\}$, $\{1,4\}\sqcup\{2\}\sqcup\{3\}$, $\{2,3\}\sqcup\{1\}\sqcup\{4\}$, $\{2,4\}\sqcup\{1\}\sqcup\{3\}$, $\{3,4\}\sqcup\{1\}\sqcup\{2\}$ et $\{1\}\sqcup\{2\}\sqcup\{3\}\sqcup\{4\}$. Si $s=2$ et si la suite est $(1,1,2,2)$, alors la somme correspondante $S(2,2)$ s'écrit 
{\small \begin{align*}
&F(1,2)\,F(2,2)-F(1,2)\,F(2,1)^{2}-F(1,2)\,F(2,1)^{2}-F(1,1)^{2}\,F(2,2)-F(1,1)^{2}\,F(2,2)\\
&-F(1,2)\,F(2,2)-F(1,1)^{2}\,F(2,1)^{2}-F(1,1)^{2}\,F(2,1)^{2}+2\,F(1,2)\,F(2,1)^{2}+2\,F(1,1)^{2}\,F(2,1)^{2}\\
&+2\,F(1,1)^{2}\,F(2,1)^{2}+2\,F(1,1)^{2}\,F(2,1)^{2}+2\,F(1,1)^{2}\,F(2,1)^{2}+2\,F(1,1)^{2}\,F(2,2)-6\,F(1,1)^{2}\,F(2,1)^{2}
\end{align*}}
ce qui donne bien $0$.
\end{example}
\begin{proof}
La preuve suivante est due à M. Sage\footnote{<<~Entre minuit et une heure du matin...~>>}. On raisonne par récurrence sur les $r_{i}$. Si $r_{i}=1$ pour tout $i$, alors le multi-ensemble considéré est simplement $\{1,2,\ldots,r\}$, et la somme $S(r_{1},\ldots,r_{s})= S(1,\ldots,1)$ s'écrit :
$$S(1^{r})=\sum_{\pi \in \mathfrak{Q}_{r}} \mu(\pi,\lle 1,r\rre)\,\prod_{j=1}^{\ell(\pi)} \prod_{r_{ij} \geq 1}F(i,r_{ij})=\left(\sum_{\pi \in \mathfrak{Q}_{r}} \mu(\pi,\lle 1,r\rre)\right)\,\left(\prod_{i=1}^{r}F(i,1)\right)=0$$
car $r\geq 2$, et la somme ci-dessus est le produit de convolution $$\zeta *\mu \left(\sqcup_{i=1}^{r}\{i\},\lle 1,r\rre\right)=\delta\left(\sqcup_{i=1}^{r}\{i\},\lle 1,r\rre\right)$$ dans l'algèbre d'incidence du treillis $\mathfrak{Q}(\lle 1,r\rre)$, donc vaut $0$ car $r\geq 2$ et $\sqcup_{i=1}^{r} \{i\} \neq \lle 1,r\rre$ dans $\mathfrak{Q}(\lle 1,r\rre)$ --- on renvoie à la section \ref{fiber} pour des précisions sur les fonctions de M\"obius, les algèbres d'incidence d'ensembles ordonnés et les notations employées ici.\bigskip

Supposons maintenant que l'un des entiers $i \in \lle 1,s\rre$ apparaît plus d'une fois ; à permutation près, on peut supposer que c'est $i=1$, et on notera donc $r_{1}+1$ le nombre d'occurrences de cet entier, et $r+1=(r_{1}+1)+r_{2}+\cdots+r_{s}$ le cardinal total de l'ensemble. Par hypothèse de récurrence, étant données des variables indépendantes $(X_{i,q})_{i,q\geq 1}$, le polynôme
$$\sum_{\pi \in \mathfrak{Q}_{r}} \mu(\pi,\lle 1,r\rre)\,\prod_{j=1}^{\ell(\pi)}\prod_{r_{ij}\geq 1}X_{i,r_{ij}}$$
est formellement nul. Si $\pi$ est une partition d'ensemble, notons $v_{i,q}(\pi)$ le nombre de parts de $\pi$ qui contiennent exactement $q$ symboles $i$ ; alors, par hypothèse de récurrence,
$$P((X_{i,q})_{i,q \geq 1},r_{1},\ldots,r_{s})=\sum_{\pi \in \mathfrak{Q}_{r}}\mu(\pi)\,\prod_{i=1}^{s}\prod_{q \geq 1}X_{i,q}^{v_{i,q}(\pi)}=0$$
en tant que polynôme, avec $\mu(\pi)=\mu(\pi,\lle1,r\rre)=(-1)^{\ell(\pi)-1}\,(\ell(\pi)-1)!\,$. Ceci étant, une partition $\pi$ dans $\mathfrak{Q}_{r+1}$ du multi-ensemble $(1^{r_{1}+1},2^{r_{2}},\ldots,s^{r_{s}})$ peut être obtenue à partir d'une partition $\pi'$ du multi-ensemble $(1^{r_{1}},2^{r_{2}},\ldots,s^{r_{s}})$ soit en rajoutant $1$ à l'une des parts déjà existante de $\pi'$, soit en mettant $1$ tout seul dans une nouvelle part.\vspace{2mm}
\begin{enumerate}
\item Dans le premier cas, la fonction de M\"obius $\mu(\pi)$ reste égale à $\mu(\pi')$, et le poids $\omega(\pi')=\prod_{j=1}^{\ell(\pi)}\prod_{r_{ij}\geq 1}F(i,r_{ij})$ est multiplié par un facteur $F(1,a_{j}+1)/F(1,a_{j})$, où $a_{j}$ est le nombre de $1$ dans la part $\pi_{j}$ de $\pi'$ à laquelle on adjoint le $1$.\vspace{2mm}
\item Dans le second cas, la fonction de M\"obius $\mu(\pi)$ vaut $-\ell(\pi')\,\mu(\pi')$, et le poids $\omega(\pi)$ est égal à $F(1,1)\times \omega(\pi')$.\vspace{2mm}
\end{enumerate}
Par conséquent, la somme $S$ vérifie la relation de récurrence 
$$S(r_{1}+1,\ldots,r_{s})=\sum_{\pi' \in \mathfrak{Q}_{r}}\sum_{j=1}^{\ell(\pi')}\frac{F(1,a_{j}+1)}{F(1,a_{j})}\,\mu(\pi')\,\omega(\pi')-\sum_{\pi' \in \mathfrak{Q}_{r}}F(1,1)\,\ell(\pi')\,\mu(\pi')\,\omega(\pi')\,.$$
La première partie s'écrit aussi 
$$\sum_{q\geq 1}\sum_{\pi' \in \mathfrak{Q}_{r}}\frac{F(1,q+1)}{F(1,q)}\,v_{1,q}(\pi')\,\mu(\pi')\omega(\pi')=\sum_{q\geq 1} F(1,q+1)\left.\frac{\partial P(r_{1},\ldots,r_{s})}{\partial X_{1,q}}\right|_{(X_{i,q})_{i,q \geq 1}=(F(i,q))_{i,q\geq 1}}=0$$
car le polynôme est par hypothèse de récurrence formellement nul. De même, la seconde partie est proportionnelle à :
$$\sum_{q\geq 1}\sum_{\pi'\in \mathfrak{Q}_{r}} v_{1,q}(\pi')\,\mu(\pi')\,\omega(\pi')= \sum_{q\geq 1} F(1,q)\,\left.\frac{\partial P(r_{1},\ldots,r_{s})}{\partial X_{1,q}}\right|_{(X_{i,q})_{i,q \geq 1}=(F(i,q))_{i,q\geq 1}}=0\,.$$
On conclut que $S(r_{1}+1,\ldots,r_{s})=0$ est encore nulle. 
\end{proof}\bigskip

\begin{proof}[Démonstration du lemme \ref{spacebound}]
Compte tenu de la convergence simple (non jointe) vers des gaussiennes, le cas où tous les $l_{i}$ sont identiques est déjà traité ; on peut donc supposer que certains des $l_{i}$ sont distincts. Nous noterons d'autre part $L=l_{1}+\cdots+l_{r}$. Comme le produit disjoint $\bullet$ des observables est compatible avec le degré de Kerov, dans la décomposition du cumulant standard $k(\varSigma_{l_{1}},\ldots,\varSigma_{l_{r}})$ en cumulants disjoints, si pour une partition d'ensembles $\pi$ l'un des cumulants identité est de degré de Kerov strictement inféireur à la somme des $l_{i}$, $i \in \pi_{j}$, alors le cumulant disjoint correspondant
$$k_{\pi}= k^{\bullet}\left(k^{\id}(\varSigma_{l_{i} \in \pi_{1}}), \ldots, k^{\id}(\varSigma_{l_{i} \in \pi_{s}}) \right)$$
est un $O(n^{(L-1)/2})$. D'après la première partie du lemme \ref{blopblop}, il reste donc à traiter le cas des cumulants disjoints $k_{\pi}$, où $\pi$ est une partition d'ensembles telle que pour tous indices $i_{1}$ et $i_{2}$ dans la même part de $\pi$, $l_{i_{1}}=l_{i_{2}}$. Autrement dit, on doit évaluer des cumulants disjoints du type :
$$k^{\bullet}\left(k^{\id}(\varSigma_{m_{1}},\ldots,\varSigma_{m_{1}}),\ldots,k^{\id}(\varSigma_{m_{s}},\ldots,\varSigma_{m_{s}})\right)\,.$$
Ensuite, si l'un des cumulants identité $k^{\id}(\varSigma_{m_{j}},\ldots,\varSigma_{m_{j}})$ contient $n_{j}$ termes avec $n_{j}\geq 3$, alors la seconde partie du lemme \ref{blopblop} s'applique, et $k_{\pi}$ est donc de nouveau un $O(n^{(L-1)/2})$. On peut donc encore restreindre l'ensemble des partitions d'ensembles à étudier, et supposer que tous les $n_{j}$ valent $1$ ou $2$. Bien sûr, si $n_{j}=1$, alors $k^{\id}(\varSigma_{m_{j}})=\varSigma_{m_{j}}$. D'autre part, d'après la seconde partie du lemme \ref{rageofposeidon}, le terme de degré de Kerov maximal $2m_{j}$ dans $k^{\id}(\varSigma_{m_{j}},\varSigma_{m_{j}})$ est $m_{j}\,\varSigma_{1^{m_{j}}}$ ; et on peut bien sûr négliger les autres termes. Au final, le cumulant de l'énoncé que l'on souhaite évaluer est à un $O(n^{(L-1)/2})$ près égal à une somme de cumulants disjoints du type 
$$k^{\bullet}(\varSigma_{m_{1}},\ldots,\varSigma_{m_{u}},\varSigma_{1^{m_{u+1}}},\ldots,\varSigma_{1^{m_{u+v}}})\,,$$
avec $m_{1}+\cdots+m_{u}+2(m_{u+1}+\cdots+m_{u+v})=L$. Alors, le lemme \ref{joeystarr} montre que ces cumulants disjoints sont encore des $O(n^{(L-1)/2})$. En effet, réordonnons les caractères centraux dans un tel cumulant disjoint pour l'écrire sous la forme :
$$k^{\bullet}\left(\underbrace{\varSigma_{2},\ldots,\varSigma_{2}}_{r_{1}\text{ termes}},\ldots,\underbrace{\varSigma_{s+1},\ldots,\varSigma_{s+1}}_{r_{s}\text{ termes}},\underbrace{\varSigma_{1^{2}},\ldots,\varSigma_{1^{2}}}_{r_{s+1}\text{ termes}},\ldots,\underbrace{\varSigma_{1^{t+1}},\ldots,\varSigma_{1^{t+1}}}_{r_{s+t}\text{ termes}}\right).$$
On a $2r_{1}+\cdots+(s+1)r_{s}+2(2r_{s+1}+\cdots+(t+1)r_{s+t})=L$ ; d'autre part, l'hypothèse faite au début de la preuve du lemme garantit que $s+t\geq 2$. Notons alors $F$ la fonction définie par :
$$F(i,r)=\begin{cases} f(i+1,r)\,n^{\frac{(i+1)r}{2}} &\text{si }i \in \lle 1,s\rre,\\
n^{(i+1-s) r}& \text{si }i \in \lle s+1,s+t \rre.
\end{cases}$$
Alors, compte tenu de l'expression des espérances $\Gel_{n}[\varSigma_{\mu}]$ donnée par la proposition \ref{asymptoticobsgelfand}, et en utilisant la formule de M\"obius
$$k(X_{1},\ldots,X_{r})=\sum_{\pi \in \mathfrak{Q}(\lle 1,r\rre)} (-1)^{\ell(\pi)-1}\,(\ell(\pi)-1)! \, \prod_{\pi_{j }\in \pi} \esper\left[\prod_{i \in \pi_{j}}X_{j}\right]\,,$$ 
on voit que le cumulant disjoint a pour terme de degré $L/2$ en $n$ 
$$\sum_{\pi \in \mathfrak{Q}(\lle 1, s+t\rre)}(-1)^{\ell(\pi)-1}\,(\ell(\pi)-1)! \, \prod_{j=1}^{\ell(\pi)}\, \prod_{r_{ij}\geq 1} F(i,r_{ij})\,,$$
étant entendu que par rapport à la formule exacte la différence est de degré inférieur en $n$. Le lemme \ref{joeystarr} montre que ce coefficient est nul ; ainsi, le terme de plus haut degré en $n$ est d'ordre inférieur, ce qui achève la preuve.
\end{proof}\bigskip

\begin{proof}[Preuve du théorème \ref{nolove}, seconde partie]
D'après le lemme \ref{spacebound}, les caractères centraux renormalisés convergent conjointement vers un vecteur gaussien ; il suffit dès lors de montrer que les covariances tendent vers $0$, c'est-à-dire que 
$$\frac{\Gel_{n}[\varSigma_{k}\,\varSigma_{l}]-\Gel_{n}[\varSigma_{k}]\,\Gel_{n}[\varSigma_{l}]}{n^{\frac{k+l}{2}}} \to 0\,$$
si $k\neq l$. Mais dans la preuve du lemme \ref{spacebound}, l'hypothèse $r\geq 3$ a seulement servi dans le cas exclus tout au début, c'est-à-dire lorsque tous les $l_{i}$ étaient identiques. Ici, $k \neq l$, donc la preuve s'adapte et on a bien convergence des covariances vers zéro, et indépendance asymptotique.
\end{proof}\bigskip

Le théorème \ref{nolove} décrit le comportement asymptotique des observables décomposées dans la base algébrique des caractères centraux ; en changeant de base comme dans la section \ref{schurplus}, nous allons maintenant en déduire l'asymptotique des formes $\lambda^{*}$ des diagrammes de Young renormalisés. Pour commencer, établissons l'équivalent du théorème \ref{firstasymptoticplancherel} pour les mesures de Gelfand :
\begin{theorem}[Formes limites des diagrammes sous les mesures de Gelfand, \cite{LS77,Mel10}]\label{firstasymptoticgelfand}
Le diagramme continu $\Omega$ est la forme limite des diagrammes de Young tirés suivant les mesures de Gelfand: pour tout $\eps >0$, $\Gel_{n}[\|\lambda^{*}-\Omega\|_{\infty} \geq \eps]$ tend vers $0$.
\end{theorem}
\begin{proof}
Comme dans le chapitre \ref{plancherel}, on montre que les cumulants libres convergent en probabilité :
$$\forall j\geq 2,\,\,\,R_{j}(\lambda^{*}) \to \begin{cases} 1 &\text{si }j=2,\\
0 &\text{sinon}.
\end{cases}
$$
Pour tout $k \geq 2$, on écrit $R_{k+1}=\varSigma_{k}+h$ avec $\mathrm{wt}(h) \leq k$. Alors,
$$R_{k+1}(\lambda^{*})=n^{-\frac{k+1}{2}}\,R_{k+1}(\lambda) = n^{-\frac{k+1}{2}}\,\varSigma_{k}(\lambda) + n^{-\frac{k+1}{2}}\,h(\lambda)\,.$$
Pour tout observable de diagrammes $f$, $\Gel_{n}[f]=O(n^{\mathrm{wt}(f)/2})$, donc le second terme vérifie
$$\Gel_{n}\left[\left(n^{-\frac{k+1}{2}}\,h(\lambda)\right)^{2}\right]=n^{-(k+1)}\,\,\Gel_{n}\left[h(\lambda)^{2}\right]=O(n^{-1})\to 0$$
et converge en probabilité vers $0$. D'autre part, pour les caractères centraux $\varSigma_{k\geq 2}$, $(\varSigma_{k})^{2}$ est de degré de Kerov $2k$, donc en utilisant l'estimation plus précise des espérances d'observables données par le degré de Kerov, on voit que 
$$\Gel_{n}\left[\left(n^{-\frac{k+1}{2}}\,\varSigma_{k}(\lambda)\right)^{2}\right]=O(n^{-1})\to 0\,. $$
Ainsi, le premier terme converge aussi en probabilité vers $0$, donc $R_{k+1}$ converge bien vers $0$ si $k \geq 2$. D'autre part, $R_{2}(\lambda^{*})=\varSigma_{1}(\lambda^{*})=1$ pour tout diagramme de Young, donc $R_{2}(\lambda^{*})$ converge bien en probabilité vers la constante $1$. Le même argument que pour les mesures de Plancherel --- c'est-à-dire, le contenu de la proposition \ref{strongweak} --- permet de conclure quant à la convergence uniforme en probabilité des diagrammes renormalisés vers $\Omega$.
\end{proof}\bigskip

De nouveau, on a en réalité convergence en probabilité au sens ultra-fort, c'est-à-dire qu'en plus de la convergence uniforme en probabilité, les bornes $a(\lambda^{*})$ et $b(\lambda^{*})$ des supports des diagrammes renormalisés convergent en probabilité vers les bornes $a(\Omega)=-2$ et $b(\Omega)=2$. De fa\c con équivalente, la longueur $\ell_{-,n}$ d'un plus long sous-mot décroissant, et la longueur $\ell_{+,n}$ d'un plus long sous-mot croissant dans une involution de taille $n$ tirée au hasard uniformément vérifient :
$$\lim_{n \to \infty}\frac{\ell_{-,n}}{\sqrt{n}} =\lim_{n \to \infty}\frac{\ell_{+,n}}{\sqrt{n}} = 2$$
les limites s'entendant en probabilité. Ce point est démontré dans l'article \cite{BR01}, et la distribution asymptotique de $\ell_{-,n}-2\sqrt{n}$ et  $\ell_{+,n}-2\sqrt{n}$ y est décrite par un résultat analogue au théorème \ref{bdj}, et mettant de nouveau en jeu les équations de Painlevé II. Notons que si l'on interprète une involution par un appariement partiel des entiers de $\lle 1,n\rre$, alors ce résultat se traduit par des estimations sur les nombres de croisements et d'imbrications d'appariements aléatoires, voir \cite{CDDSY05}.\bigskip

La correspondance des formes limites pour les mesures de Plancherel et les mesures de Gelfand n'est en réalité pas très étonnante, à condition d'utiliser la méthode originale de preuve du théorème de Logan-Shepp-Kerov-Vershik, c'est-à-dire celle de l'article \cite{LS77}. Ainsi, en utilisant la formule des crochets, Logan et Shepp ont construit une fonctionnelle énergie $I$ sur les diagrammes de Young continus telle que, pour $\lambda$ partition de taille $n$, $(\dim \lambda)^{2}$ ait un comportement asymptotique grossièrement équivalent à $\exp(-nI(\lambda^{*}))$. Or, la courbe $\Omega$ est l'unique diagramme continu qui maximise $I$. Par conséquent, il est intuitivement évident que la mesure de Plancherel charge surtout les diagrammes qui sont proches de cette courbe (à renormalisation près). Cet argument variationnel peut être rendu rigoureux, voir \cite{LS77}. Maintenant, si l'on remplace $(\dim \lambda)^{2}$ par $\dim \lambda$, on peut clairement utiliser le même type d'argument, en divisant simplement l'énergie par $2$. En réalité, on pourrait adapter la preuve au cas plus général des \textbf{$\beta$-mesures de Plancherel}, qui sont les mesures définies par :
$$M_{n}^{\beta}[\lambda]=\frac{(\dim \lambda)^{\beta}}{\sum_{\mu \in \ym_{n}} (\dim \mu)^{\beta}},\quad\beta>0\,.$$
Notons qu'on retrouve les mesures de Plancherel pour $\beta=2$, et les mesures de Gelfand pour $\beta=1$. Ces $\beta$-mesures sont les analogues pour les partitions des $\beta$-ensembles, qui sont les distributions ponctuelles de densité
$$d\proba[x_{1},\ldots,x_{n}]=\frac{1}{Z_{n,\beta}}\,|\Delta(x_{1},\ldots,x_{n})|^{\beta}\,\E^{-\frac{1}{2}\sum_{i=1}^{n}(x_{i})^{2}}\,dx_{1}\cdots dx_{n}\,,$$
voir par exemple \cite{DE02}. En particulier, on conjecture que la correspondance de Baik-Deift-Johansson s'étende à ces $\beta$-ensembles : ainsi, la déviation de la plus grande valeur d'un point d'un $\beta$-ensemble par rapport à la valeur $2\sqrt{n}$ est probablement équivalente en loi (après renormalisation) à la déviation de la plus grande part d'une partition sous la $\beta$-mesure de Plancherel par rapport à la valeur $2\sqrt{n}$ (on sait que c'est le cas lorsque $\beta=1,2$).
\bigskip\bigskip

Ceci étant, revenons à notre étude asymptotique des mesures de Gelfand, et introduisons comme dans la section \ref{schurplus} la déviation $\Delta_{\lambda}(s)=\lambda^{*}(s)-\Omega(s)$, et les polynômes de Chebyshev de seconde espèce $u_{k}(X)$, avec la condition de normalisation 
$$u_{k}(2\cos \theta)=\frac{\sin (k+1)\theta}{\sin \theta}\,.$$
On note $W_{k}$ la variable aléatoire $\frac{\varSigma_{k}}{n^{k/2}}$ ; pour $k\geq 2$, elle converge vers une gaussienne $\mathcal{N}(e_{k},2k)$. Le lemme suivant est l'exact analogue du lemme \ref{recovery} pour les mesures de Gelfand :
\begin{lemma}[Développement gaussien des moments de la déviation d'un diagramme]
Le $k$-ième moment de la déviation 
$$\frac{\sqrt{n}}{2}\int_{\R}s^{k}\,(\lambda^{*}(s)-\Omega(s))\,ds\,.$$
est égal à $\sqrt{n}\,\frac{\tilp_{k+2}(\lambda^{*})-\tilp_{k+2}(\Omega)}{(k+1)(k+2)}$, et c'est aussi
$$\sum_{j=0}^{\lfloor \frac{k-1}{2}\rfloor} \frac{k!}{k+1-j!\,j!}\,W_{k+1-2j}\,,$$
plus une variable aléatoire qui tend en probabilité vers $0$ sous les mesures de Gelfand.
\end{lemma}
\begin{proof}
La première partie du lemme a déjà été démontrée à la fin du chapitre \ref{plancherel} ; on renvoie sinon à \cite[proposition 7.2]{IO02}. Rappelons que les moments entrela\c cés du diagramme continu $\Omega$ sont :
$$\tilp_{2k}(\Omega)=\binom{2k}{k}\qquad;\qquad\tilp_{2k+1}(\Omega)=0$$
voir la proposition \ref{momentwignerdeform}. D'autre part, il existe un résultat analogue au lemme \ref{otherchangeofbasis} pour le degré de Kerov :
$$\tilp_{k}=\sum_{j=0}^{\lfloor \frac{k-3}{2}\rfloor}\frac{k^{\downarrow j+1}}{j!}\,\varSigma_{k-1-2j}\,(\varSigma_{1})^{j} +\begin{cases}\binom{k}{k/2}\,(\varSigma_{1})^{k/2}&\text{si }k \text{ est pair},\\
0&\text{sinon},\end{cases}$$
plus une observable de degré de Kerov inférieur à $k-2$, voir \cite[proposition 7.3]{IO02}. Par conséquent,
$$\frac{\tilp_{k+2}(\lambda^{*})}{(k+1)(k+2)}=\sum_{j=0}^{\lfloor \frac{k-1}{2}\rfloor} \frac{k!}{k+1-j!\,j!}\,\frac{\varSigma_{k+1-2j}(\lambda)}{n^{\frac{k+2-2j}{2}}}+\begin{cases} \binom{k+2}{(k+2)/2}&\text{si }k\text{ est pair},\\0 &\text{if }k\text{ est impair},\end{cases}$$ 
plus $n^{-\frac{k+2}{2}}$ fois une observable de degré de Kerov au plus $k$. En soustrayant $\tilp_{k+2}(\Omega)$, on se débarrasse du second terme à droite, et en multipliant par $\sqrt{n}$, on obtient finalement $$\sum_{j=0}^{\lfloor \frac{k-1}{2}\rfloor} \frac{k!}{k+1-j!\,j!}\,W_{k+1-2j}\,$$
plus $n^{-\frac{k+1}{2}}$ fois une observable de degré de Kerov plus petit que $k$, donc d'ordre de grandeur au plus $n^{\frac{k}{2}}$. Ce terme résiduel est donc négligeable à l'infini, et le lemme est établi.
\end{proof}\bigskip

\begin{lemma}[Polynômes de Chebyshev de second espèce et déviation d'un diagramme]
Pour tout entier $k$, l'observable 
$$\Upsilon_{k}=\frac{\sqrt{n}}{2}\int u_{k}(s)\,\Delta_{\lambda}(s)\,ds$$ est égale à $\frac{W_{k+1}}{k+1}$, plus une observable qui sous les mesures de Gelfand tend en probabilité vers $0$ lorsque $n$ tend vers l'infini.
\end{lemma}
\begin{proof}
On rappelle que le développement explicite du polynôme $u_{k}(X)$ est 
$$u_{k}(X)=\sum_{m=0}^{\lfloor \frac{k}{2}\rfloor} (-1)^{m}\binom{k-m}{m}\,X^{k-2m}\,.$$ 
Par conséquent, si $\Theta_{k}$ désigne l'observable du lemme précédent, alors 
$$\Upsilon_{k}= \sum_{m=0}^{\lfloor\frac{k}{2}\rfloor} (-1)^{m}\,\binom{k-m}{m}\,\Theta_{k-2m}\qquad;\qquad \Theta_{k} \simeq \frac{1}{k+1}\sum_{j=0}^{\lfloor \frac{k}{2}\rfloor} \binom{k+1}{j}\,X_{k+1-2j}\,.$$
Les mêmes arguments que dans la preuve du théorème \ref{almosttheorem} montrent que ces relations sont essentiellement inverses l'une de l'autre, et ainsi, $\Upsilon_{k}\simeq \frac{W_{k+1}}{k+1}$, qui est asymptotiquement une gaussienne (pas forcément centrée).
\end{proof}\bigskip

\begin{theorem}[Déviation des diagrammes par rapport à leur forme limite sous les mesures de Gelfand, \cite{Mel10}]\label{scaledkerovgaussian}
Soit $(\xi_{k})_{k \geq 2}$ une famille de variables gaussiennes indépendantes, centrées et toutes de variance $1$. Sous les mesures de Gelfand et au sens des distributions tempérées sur l'intervalle $[-2,2]$, le processus $\sqrt{n}\,\Delta_{\lambda}(s)$ converge en loi vers le processus gaussien généralisé 
$$\frac{1}{2}-\frac{\sqrt{4-s^{2}}}{\pi}+\sqrt{2}\,\Delta(s)=\frac{1}{2}-\frac{\sqrt{4-s^{2}}}{\pi}+\sqrt{2}\,\left(\frac{2}{\pi}\sum_{k=2}^{\infty}\frac{\xi_{k}}{\sqrt{k}}\,\sin(k\theta)\right),$$
avec $s=2\cos \theta$. À une fonction déterministe près, le processus gaussien de Kerov (multiplié par $\sqrt{2}$) décrit donc encore les fluctuations des diagrammes par rapport à leur forme limite.
\end{theorem}
\begin{proof}
Si $f$ est une fonction de $\mathscr{C}^{\infty}([-2,2])$, on peut la développer dans la base orthonormale de $\mathscr{L}^{2}([-2,2],m_{\Omega})$ constituée par les polynômes de Chebyshev. Ainsi, 
$$f(s)=\sum_{k=0}^{\infty} \left(\int_{-2}^{2}f(s)\,u_{k}(s)\, \frac{\sqrt{4-s^{2}}}{2\pi}\,ds\right)\,u_{k}(s)\,.$$
La fonctionnelle linéaire de la fonction $f$ associée à la déviation renormalisée $\frac{\sqrt{n}}{2}\,\Delta_{\lambda}(s)$ s'écrit donc asymptotiquement
$$\int_{-2}^{2}f(s)\,\left(\sum_{k=1}^{\infty} \frac{W_{k+1}}{k+1}\,u_{k}(s)\, \frac{\sqrt{4-s^{2}}}{2\pi}\right)\,ds\,,$$
et ainsi, en tant que distribution, $\sqrt{n}\,\Delta_{\lambda}(s)$ converge faiblement vers le processus gaussien généralisé
$$\frac{1}{\pi}\sum_{k=1}^{\infty}\frac{\mathcal{N}(e_{k+1},2(k+1))}{k+1}\,u_{k}(s)\,\sqrt{4-s^{2}}=\frac{2}{\pi}\sum_{k=2}^{\infty}\frac{\mathcal{N}(e_{k},2k)}{k}\,\sin (k\theta)\,,$$ avec $s=2\cos\theta$ et $\theta \in [0,\pi].$ Décomposons les variables $\mathcal{N}(e_{k},2k)$ sous la forme $e_{k}+\sqrt{2k}\,\xi_{k}$ avec les $\xi_{k}$ comme dans l'énoncé du théorème. Alors, le processus limite est somme des deux contributions :
$$\frac{2}{\pi}\sum_{k=1}^{\infty} \frac{\sin(2k+1)\theta}{2k+1}= \frac{1}{2}-\frac{2\sin\theta}{\pi} =\frac{1}{2}-\frac{\sqrt{4-s^{2}}}{\pi}\qquad\!\!\text{et}\qquad \!\!\frac{2\sqrt{2}}{\pi}\left(\sum_{k=2}^{\infty}\frac{\zeta_{k}}{\sqrt{k}}\,\sin (k\theta)\right) = \sqrt{2}\,\Delta(s)\,,$$
où l'on reconnaît dans le second terme le processus $\Delta(s)$ des théorèmes \ref{secondasymptoticplancherel} et \ref{schurweylplus}.
\end{proof}\bigskip

Ainsi, le processus gaussien de Kerov décrit les fluctuations des formes des diagrammes dans au moins deux autres cadres que celui des mesures de Plancherel : le cadre des mesures de Schur-Weyl et le cadre des mesures de Gelfand. On peut raisonnablement conjecturer que pour tout paramètre $\beta$, les fluctuations des formes des diagrammes de Young tirés aléatoirement selon les $\beta$-mesures de Plancherel sont encore décrites par un multiple\footnote{On peut conjecturer que le facteur multiplicatif est $\sqrt{\frac{2}{\beta}}$.} de ce processus, plus éventuellement une fonction déterministe. L'idéal serait de disposer de théorèmes universels pour les partitions semblables à ce que l'on connaît pour les modèles de matrices aléatoires (voir en particulier \cite{LP09}) ; d'après les deux chapitres de cette partie, il semble qu'on puisse établir un théorème central limite semblable à celui de Kerov dès que les caractères centraux cycliques renormalisés convergent conjointement vers un vecteur gaussien.

\pagestyle{empty}

\clearpage
~
\clearpage
~

\pagestyle{fancy}
\fancyhead{}
\fancyfoot{}
\fancyfoot[C]{\thepage}
\renewcommand{\chaptermark}[1]{\markboth{\chaptername\ \thechapter.\ #1.}{}} 
\renewcommand{\sectionmark}[1]{\markright{\thesection.\ #1.}} 
\fancyhead[RO]{\rightmark}
\fancyhead[LE]{\leftmark}
\setlength{\headheight}{15.5pt}

\part{Identités génériques dans les algèbres de groupes}
Les trois premières parties du mémoire ont mis en évidence les propriétés de concentration des mesures issues de la théorie des représentations des groupes symétriques $\sym_{n}$, de leurs algèbres d'Hecke $\IH_{q}(\sym_{n})$ et des groupes de Chevalley finis $\GL(n,\For_{q})$ ou $\Sp(2n,\For_{q})$ ; de plus, on a presque toujours un comportement asymptotique gaussien. L'un des arguments centraux dans les preuves de ces résultats était la propriété de factorisation asymptotique des caractères centraux dans l'algèbre des observables de diagrammes $\obs$ :
\begin{align*}\varSigma_{\mu}\,*\,\varSigma_{\nu}&=\varSigma_{\mu\sqcup \nu} + (\text{termes de degré strictement inférieur})\\
&=\varSigma_{\mu\sqcup \nu} + (\text{termes de poids strictement inférieur})\\
&=\varSigma_{\mu\sqcup \nu} + (\text{termes de $\alpha$-degré strictement inférieur})\,.
\end{align*}
De plus, les covariances des lois gaussiennes limites sont liées aux termes suivants des développements de $\varSigma_{\mu}\,\varSigma_{\nu}$. \bigskip

Cette propriété de factorisation asymptotique est tout à fait évidente si l'on interprète les $\varSigma_{\mu}$ comme éléments de l'\textbf{algèbre des permutations partielles} (voir \cite{IK99}) : en effet, dans un <<~grand~>> groupe symétrique $\sym_{n}$, une permutation de type $\mu \sqcup 1^{n- |\mu|}$ et une permutation de type $\nu \sqcup 1^{n-|\nu|}$ ont une très grande chance d'avoir leurs supports essentiels disjoints, donc la contribution principale d'un produit $\varSigma_{\mu}\,\varSigma_{\nu}$ est bien $\varSigma_{\mu \sqcup \nu}$. Notons d'autre part que l'existence de l'algèbre des permutations partielles permet une preuve très simple du théorème de Farahat-Higman \ref{farahathigman}.\bigskip
\bigskip

Dans cette dernière partie, nous cherchons à généraliser ces raisonnements au cas d'autres groupes ou algèbres, et à construire des \textbf{algèbres d'Ivanov-Kerov} pour la famille des algèbres d'Hecke de type A $(\IH_{q}(\sym_{n}))_{n \geq 1}$ et la famille des groupes linéaires finis $(\GL(n,\For_{q}))_{n\geq 1}$. Le premier cas que nous traitons est celui des algèbres d'Hecke (chapitre \ref{badbeat}), et après avoir rappelé les détails de la construction d'Ivanov-Kerov, nous démontrons l'exact analogue du théorème \ref{farahathigman} pour les algèbres $\IH_{q}(\sym_{n})$, les classes de conjugaison de $Z(\C\sym_{n})$ étant remplacées par les \textbf{éléments de Geck-Rouquier} dans $Z(\IH_{q}(\sym_{n}))$, \emph{cf.} le théorème \ref{heckefarahathigman} --- c'est le résultat principal de l'article \cite{Mel10c}. Notons que la preuve est nettement plus alambiquée que dans le cas du groupe symétrique ; c'est parce que les centres des algèbres d'Hecke ont une structure combinatoire beaucoup plus riche que les centres des algèbres des groupes symétriques (voir \cite{Las06,Jon90,Fra99,GR97}).\bigskip\bigskip

Les algèbres de permutations partielles et leurs analogues Hecke rentrent dans le cadre plus général des \textbf{fibrés de semi-groupes par des semi-treillis} ; dans le chapitre \ref{bundle}, nous donnons un autre exemple important relié au problème des nombres de Hurwitz, et nous expliquons la construction générale. L'idée est de combiner les constructions suivantes :\vspace{2mm}
\begin{enumerate}
\item \'Etant donné un semi-groupe $M$ et un semi-treillis $L$ qui <<~fibre~>> les éléments de $M$, on peut construire un semi-groupe fibré $M \times_{I} L$. La théorie des représentations de ce nouveau semi-groupe est aisée à déterminer à partir de celle de $M$ et de ses sous-semi-groupes.\vspace{2mm}
\item \'Etant donnée une tour de semi-groupes $(M_{n})_{n \in \N}$ et une famille croissante $(L_{n})_{n\in \N}$ de semi-treillis fibrant les $L_{n}$, on peut sous certaines conditions construire une limite projective des $M_{n}\times_{I_{n}} L_{n}$, qui vient s'ajouter à la construction naturelle de limite injective des $M_{n}$, et joue un rôle dual.\vspace{2mm}
\end{enumerate}
À la fin du chapitre \ref{bundle}, nous présentons des conjectures concernant les produits de classes de conjugaison de $\GL(n,\For_{q})$. D'après le chapitre \ref{general}, ces classes sont indexées par des polypartitions ; on conjecture que les produits de classes complétées $C_{\bbmu \rightarrow n}$ vérifient des relations du type
$$C_{\bblambda \rightarrow n}\,*\,C_{\bbmu \rightarrow n} = \sum_{\|\bbnu\|\leq \|\bblambda\|+\|\bbmu\|} a_{\bblambda\bbmu}^{\bbnu}(n,q)\,C_{\bbnu \rightarrow n}\,$$
avec les constantes de structure $a_{\bblambda\bbmu}^{\bbnu}$ dans $\Q(q^{n},q)$, et qu'il existe une algèbre de type Ivanov-Kerov impliquée dans ces identités génériques. Les éléments d'une telle algèbre limite projective $\obs(\GL(\For_{q}))$ seraient également des candidats naturels pour des \textbf{observables de polydiagrammes}, ce qui permettrait une étude asymptotique fine des caractères des groupes $\GL(n,\For_{q})$.

\chapter{Algèbres d'Ivanov-Kerov et d'Hecke-Ivanov-Kerov}\label{badbeat}

Dans ce chapitre, nous donnons une preuve très simple du théorème de Farahat-Higman \ref{farahathigman}, et nous en démontrons un analogue dans le contexte des algèbres d'Hecke de type A (\emph{cf.} \cite{Mel10c}). Ces résultats constituent ce que nous appellerons des \textbf{identités génériques}, c'est-à-dire des identités entre éléments d'algèbres (de groupes, ou générales) $A_{n}$ qui mettent en jeu les <<~mêmes~>> éléments pour tout $n$, et avec des coefficients dépendant d'une manière explicite de $n$.

\begin{example} Dans les algèbres des groupes symétriques $\C\sym_{n}$, si $C_{1 \rightarrow n}$ désigne la somme (formelle) de toutes les transpositions et si $C_{2 \rightarrow n}$ (resp., $C_{(1,1)\rightarrow n}$) désigne la somme de tous les $3$-cycles (resp., de tous les produits de deux transpositions disjointes), alors 
$$C_{1 \rightarrow n}\,*\,C_{1\rightarrow n}=2\, C_{(1,1) \rightarrow n} + 3 \,C_{2 \rightarrow n} + \frac{n(n-1)}{2}\,\id_{\lle 1,n\rre}\,.$$
En effet, parmi les $\left(\frac{n(n-1)}{2}\right)^{2}$ éléments du terme de gauche, $\frac{n(n-1)(n-2)(n-3)}{4}$ sont des produits disjoints, et donnent l'un des produits de deux transpositions de $C_{(1,1) \rightarrow n}$, chacun de ses produits $\tau_{1}\tau_{2}$ étant atteint deux fois (par $\tau_{1}\times \tau_{2}$ et $\tau_{2}\times \tau_{1}$) ; $n(n-1)(n-2)$ sont des produits de transpositions avec un point commun, donc donnent un $3$-cycle, qui peut s'écrire de trois fa\c cons comme produit de transpositions ; et $\frac{n(n-1)}{2}$ sont des carrés de transpositions, donc donnent un terme $\id_{\lle 1,n\rre}$.
\end{example}\bigskip

L'idée des preuves de ces résultats est de construire une algèbre $A^{\infty}$ qui se projette sur toutes les algèbres $A_{n}$, et dans laquelle on peut effectuer les calculs ; alors, toute identité dans $A^{\infty}$ donne une identité valable dans tous les $A_{n}$. Malheureusement, lorsque les $A_{n}$ forment une tour d'algèbres --- \emph{i.e.}, $$\cdots \subset A_{n}\subset A_{n+1}\subset \cdots$$ pour tout $n$ --- l'objet universel naturel est la limite injective $A_{\infty}=\varinjlim_{n \to \infty}A_{n}$, et il n'existe pas de limite projective (qui serait un candidat naturel pour $A^{\infty}$). L'\textbf{algèbre d'Ivanov-Kerov} (\emph{cf.} \cite{IK99}) permet de contourner cette difficulté, et l'essentiel de cette dernière partie est consacrée à des généralisations de cette belle construction, qui est rappelée dans la section \ref{ivanovkerov}.  \bigskip\bigskip

Ainsi, le résultat nouveau principal de ce chapitre est le théorème \ref{heckefarahathigman}, qui assure que le théorème de Farahat-Higman reste vrai dans les centres des algèbres d'Hecke, à condition de remplacer les classes de conjugaison par les \textbf{éléments de Geck-Rouquier} (voir \cite{GR97}), et les coefficients polynomiaux en $n$ par des coefficients polynomiaux en $n$ et en $q$ (ceci avait été conjecturé par A. Francis et W. Wang dans \cite{FW09}). Ces éléments de Geck-Rouquier ont une définition implicite, et on ne dispose pas d'expression explicite \emph{a priori} ;  il est donc très difficile de construire une algèbre <<~limite projective~>> les contenant tous. Fort heureusement, les centres des algèbres d'Hecke de type A admettent d'autres bases plus faciles à manipuler, en particulier des bases constituées de \textbf{normes}. Nous expliquons ceci dans la section \ref{lascoux}, en nous appuyant pour l'essentiel sur \cite{Las06} ; ainsi, la base canonique de notre \textbf{algèbre d'Hecke-Ivanov-Kerov} sera constituée de \textbf{normes génériques}, voir la section \ref{genericnorm}. Des exemples d'identités génériques dans les centres des algèbres $\C\sym_{n}$ et $\IH_{q}(\sym_{n})$ seront donnés tout au long du chapitre. \bigskip

\section{Permutations partielles et permutations composées}\label{ivanovkerov}
Dans tout ce qui suit, si $\mu$ est une partition et $n$ est un entier plus grand que $|\mu|+\ell(\mu)$, nous noterons $\mu\rightarrow n$ la partition complétée de taille $n$ obtenue en rajoutant $1$ aux $n-|\mu|$ premières parts (éventuellement nulles) de $\mu$. D'autre part, si $\mu=(\mu_{1},\ldots,\mu_{r})$, alors $\mu-1$ désigne la partition $(\mu_{1}-1,\ldots,\mu_{r}-1)$, étant entendu qu'on supprime ensuite les éventuelles parts nulles. Ainsi, pour toute partition, $\mu\sqcup 1^{n-\mu}= (\mu-1)\rightarrow n$. On adopte les mêmes notations avec $\mu+1$.\bigskip\bigskip

Si $n\geq 1$, on rappelle qu'une \textbf{permutation partielle} de taille $n$ est la donnée d'une partie $S \subset \lle 1,n\rre$ et d'une permutation $\sigma \in \sym(S)$, qu'on peut également voir comme une permutation dans $\sym_{n}$ laissant fixes tous les points hors de $S$ (mais éventuellement aussi des points dans $S$). Par exemple,
$$\sigma=(1,2,5)\,,\,\,\,S=\{1,2,4,5,7\}\subset \lle 1,8\rre $$
est une permutation partielle de taille $8$. Ainsi que nous l'avons expliqué dans la section \ref{centralcharacter}, les permutations partielles peuvent être multipliées entre elles par la règle
$$ (\sigma,S)(\tau,T)=(\sigma\tau,S\cup T)\,,$$
et l'algèbre complexe $\blg_{n}$ du monoïde ainsi constituée se projette naturellement sur $\C\sym_{n}$ en oubliant les supports des permutations partielles (nous noterons $\pi_{n}$ cette projection). Enfin, en autorisant des supports dans l'ensemble des parties finies de $\N$, on peut construire une limite projective $\blg_{\infty}=\varprojlim_{n\to \infty} \blg_{n}$ dans la catégorie des algèbres filtrées, et relativement aux morphismes d'algèbres 
$$\phi_{N,n} (\sigma,S)= \begin{cases} (\sigma,S)& \text{si }S\subset \lle 1,n\rre, \\0 &\text{sinon.} \end{cases}$$
de $\blg_{N}$ dans $\blg_{n \leq N}$. Les éléments de $\blg_{\infty}$ sont les combinaisons linéaires formelles (éventuellement infinies) de permutations partielles de degrés restant bornés.
\bigskip\bigskip

Un premier point nouveau qui mérite d'être précisé est la théorie des représentations des algèbres de permutations partielles $\blg_{n}$. Ainsi :
\begin{proposition}[Théorie des représentations des algèbres de permutations partielles, \cite{IK99}]\label{reprpartialpermutation}
L'algèbre $\blg_{n}$ des permutations partielles de taille $n$ est isomorphe à la somme directe 
$$\bigoplus_{S \subset \lle 1,n\rre}\C\sym(S)\,.$$
Elle est donc semi-simple, et ses classes de modules simples sont indexées par les partitions de taille inférieure à $n$.
\end{proposition}
\begin{proof}
Nous suivons la preuve de \cite[proposition 3.2]{IK99}. Si $\sigma\in \sym(S)$, nous noterons $\sigma_{S}=[0,\ldots,0,\sigma,0,\ldots,0]$ le vecteur de $\bigoplus_{S \subset \lle 1,n\rre}\C\sym(S)$ qui a toutes ses coordonnées nulles dans les sous-espaces $\C\sym(T)$ avec $T\neq S$, et pour coordonnée $\sigma$ dans $\C\sym(S)$. Lorsque $S$ parcourt $\Part(\lle 1,n\rre)$ et $\sigma$ parcourt $\sym(S)$, $(\sigma_{S})$ forme une base linéaire de $\bigoplus_{S \subset \lle 1,n\rre}\C\sym(S)$. Le point difficile de la proposition est la détermination de l'isomorphisme : en effet, l'application linéaire
$$(\sigma,S) \in \blg_{n} \mapsto \sigma_{S} \in \bigoplus_{S \subset \lle 1,n\rre}\C\sym(S)$$
est seulement un isomorphisme d'espaces vectoriels. Considérons plutôt l'application linéaire définie par 
$$\psi : (\sigma,S)\mapsto \sum_{S\subset T} \sigma_{T}\,.$$
Ceci a bien un sens, car si une permutation $\sigma$ laisse fixe tous les points hors de $S$, alors \emph{a fortiori} elle laisse fixe tous les points hors d'une partie $T$ contenant $S$ ; ainsi, $\sigma$ est bien un élément de $\sym(T)$ pour tout $T \supset S$. On montre sans mal que $\psi$ est bien un morphisme d'algèbres :
\begin{align*}\psi((\sigma,S))\,\psi((\tau,T))&=\sum_{\substack{S \subset S'\\T \subset T'}} \sigma_{S'}\,\tau_{T'}=\sum_{\substack{S\subset U\\ T \subset U}}(\sigma\tau)_{U}=\sum_{\substack{S\cup T \subset  U}}(\sigma\tau)_{U}\\
&=\psi((\sigma\tau,S\cup T))=\psi((\sigma,S)(\tau,T))\,.\end{align*}
De plus, si $\mu(S,T)=(-1)^{\card T - \card S}$, alors l'application linéaire
$$\theta : \sigma_{S} \mapsto \sum_{S\subset T}\mu(S,T)\,(\sigma,T)$$
est la réciproque de $\psi$. En effet,
$$\theta(\psi(\sigma,S)))=\theta\left(\sum_{S\subset T}\sigma_{T} \right)=\sum_{S \subset T \subset U} \mu(T,U)\, (\sigma,U)=\sum_{S \subset U} \left( \sum_{T \,\,|\,\,S\subset T\subset U} \mu(T,U)\right)(\sigma,U)\,,$$
et si $k=\card U-\card S$, alors le coefficient entre parenthèses s'écrit $\sum_{i=0}^{k} \binom{k}{i} (-1)^{k-i}=(1-1)^{k}$ en réunissant les parts $T$ suivants leur cardinalité. Tous les coefficients sont donc nuls, sauf celui de $(\sigma,S)$, qui vaut $1$ ; on a donc bien $\theta(\psi((\sigma,S)))=(\sigma,S)$, ce qu'il fallait démontrer. Ainsi,
$$\blg_{n}\simeq \bigoplus_{S \subset \lle 1,n\rre}\C\sym(S) \simeq \bigoplus_{k=0}^{n}\bigoplus_{S \in \Part_{k}(\lle 1,n\rre)}\bigoplus_{\lambda \in \ym_{k}} \hendo(V_{\lambda})\,,$$
donc $\blg_{n}$ est semi-simple et a bien ses classes de modules simples indexées par les partitions de taille inférieure à $n$.
\end{proof}
\noindent On donnera plus loin des résultats (et des preuves) analogues pour les algèbres (d'Hecke) de permutations composées, pour les algèbres de permutations scindées et plus généralement pour les algèbres de fibrés de semi-groupes par des semi-treillis, voir la proposition \ref{fiberrepresentationtheory}.
\bigskip\bigskip

Comme $\blg_{\infty}$ se projette sur les $\blg_{n}$ puis sur les $\C\sym_{n}$, cette algèbre peut être vue comme une limite projective\footnote{Dans ce chapitre et le suivant, nous utiliserons le terme <<~limite projective~>> de fa\c con quelque peu abusive ; par exemple, au sens strict du terme, $\blg_{\infty}$ est la limite projective des $\blg_{n}$, mais nous dirons aussi que c'est une limite projective des algèbres $\C\sym_{n}$. } des algèbres des groupes symétriques. Voyons maintenant comment construire une limite projective des centres de ces algèbres. Malheureusement, le centre de l'algèbre $\blg_{n}$ n'a que peu de rapports avec $Z(\C\sym_{n})$ ; d'après ce qui précède, on peut juste dire qu'il est isomorphe à $\bigoplus_{S\subset \lle 1,n\rre}Z(\C\sym(S))$. Par contre, on peut considérer la sous-algèbre $\alg_{n}\subset \blg_{n}$ constituée des invariants pour l'action du groupe symétrique par conjugaison sur les permutations partielles :
$$\sigma\cdot(\tau,T)=(\sigma\tau\sigma^{-1},\sigma(T))\,.$$
\'Etendons cette action en une action linéaire de $\sym_{n}$ sur $\blg_{n}$, et notons $\alg_{n}=(\blg_{n})^{\sym_{n}}$ la sous-algèbre\footnote{Comme l'action est compatible avec le produit des permutations partielles, on obtient bien une sous-algèbre.} des invariants pour cette action. Si deux permutations partielles $(\sigma,S)$ et $(\tau,T)$ sont conjuguées par $\sym_{n}$, alors on a évidemment $k=\card S=\card T$, puis, le type cyclique de $\sigma \in \sym(S)$ doit être le même que le type cyclique de $\tau \in \sym(T)$. Par conséquent, les classes de conjugaison de permutations partielles sont elles aussi en bijection avec l'ensemble des partitions de taille $k$ inférieure à $n$, et une base de $\alg_{n}$ est donnée par les :
$$A_{\mu,n}=\sum_{\substack{S \in \Part(\lle 1,n\rre)\\ \card S=|\mu|}}\sum_{\substack{\sigma \in \sym(S)\\ t(\sigma)=\mu}} (\sigma,S)\,.$$
Si $A_{\mu}$ désigne la même somme formelle de permutations partielles avec $S$ n'importe quelle partie finie de $\N$, alors on a d'une part 
$$\phi_{\infty,n}(A_{\mu})=\begin{cases} A_{\mu,n}&\text{si }n\geq |\mu|,\\
0 &\text{sinon,}\end{cases}$$ 
et d'autre part, $A_{\mu}$ est un multiple de l'élément $\varSigma_{\mu}$ défini dans la section \ref{centralcharacter}. 
\begin{definition}[Algèbre d'Ivanov-Kerov]
L'algèbre d'Ivanov-Kerov des groupes symétriques est la sous-algèbre de $\blg_{\infty}$ engendrée (linéairement) par les $A_{\mu}$ (ou les $\varSigma_{\mu}$) ; cette algèbre est aussi
$$\alg_{\infty}=(\blg_{\infty})^{\sym_{\infty}}=\varprojlim_{n \to \infty} \alg_{n}\,.$$
Elle est commutative et graduée par le degré $\deg(A_{\mu})=|\mu|$, et aussi par $\mathrm{deg}_{-1/2}(A_{\mu})=|\mu|-\ell(\mu)$.
\end{definition}
\begin{proof}
Ces faits affirmés dès le chapitre \ref{tool} sont en réalité relativement aisés à démontrer. Notons $\alg_{\infty}$ l'espace vectoriel engendré par les $A_{\mu}$ dans $\blg_{\infty}$. Pour les mêmes raisons que dans le cas fini, $\alg_{\infty}$ est égal à l'espace des invariants $(\blg_{\infty})^{\sym_{\infty}}$. Mais comme l'action de $\sym_{\infty}$ est compatible avec le produit, cet espace d'invariants est une sous-algèbre ; ainsi, $\alg_{\infty}$ est bien une sous-algèbre, et elle est graduée par $\deg(A_{\mu})=|\mu|$, car ceci est la restriction de la filtration canonique de $\blg_{\infty}$ aux éléments de $\alg_{\infty}$. D'autre part, comme $\phi_{\infty,n}(A_{\mu})=A_{\mu,n}$, cette algèbre est également la limite projective des algèbres graduées $\alg_{n}=(\blg_{n})^{\sym_{n}}$. Pour le caractère commutatif, on peut utiliser la définition des $A_{\mu}$ pour montrer qu'ils commutent avec toutes les permutations partielles ; ainsi, $\alg_{\infty}\subset Z(\blg_{\infty})$, mais l'inclusion est stricte. Finalement, $\deg_{-1/2}$ est bien une autre graduation d'algèbre sur $\alg_{\infty}$, car si une permutation partielle $(\eps,E)=(\sigma,S)(\tau,T)$ de type $\rho$ apparaît dans un produit $A_{\mu}\,A_{\nu}$ avec $t(\sigma,S)=\mu$ et $t(\tau,T)=\nu$, alors le nombre minimal de transpositions nécessaires pour écrire $\sigma$ est $|\mu|-\ell(\mu)$, et ce quelque soit la taille $|\mu|$ de $S$ ; de même pour $\tau$ avec $|\nu|-\ell(\nu)$ et pour $\eps$ avec $|\rho|-\ell(\rho)$. Comme $\eps=\sigma \tau$, on en déduit 
$$|\rho|-\ell(\rho) \leq |\mu|-\ell(\mu)+|\nu|-\ell(\nu)\,,$$
ce que l'on voulait démontrer.
 \end{proof}
 \noindent Nous noterons $g_{\mu\nu}^{\rho}$ les coefficients de structure de $\alg_{\infty}$ dans la base $A_{\mu}$ ; ainsi, $A_{\mu}\,A_{\nu}=\sum_{\rho}g_{\mu\nu}^{\rho}\,A_{\rho}$.
\bigskip\bigskip

Maintenant, remarquons que les $A_{\mu,n}$ ont pour projetés sur les algèbres $\C\sym_{n}$ des multiples des classes de conjugaisons\footnote{On convient que $C_{\lambda\rightarrow n}=0$ si $\lambda+\ell(\lambda)>n$.} $C_{(\mu-1 )\rightarrow n}$. Plus précisément, 
$$\pi_{n}(A_{\mu,n})=\binom{n-|\mu|+m_{1}(\mu)}{m_{1}(\mu)}C_{(\mu-1)\rightarrow n}\,,$$
car une permutation partielle de type $\mu$ et de taille $n$ est entièrement déterminée par la donnée d'une permutation $\sigma$ de type cyclique $\mu \sqcup 1^{n-\mu}$ et d'un marquage de $m_{1}(\mu)$ points fixes dans le support de la permutation partielle (les autres points du support sont déjà déterminés par les cycles de $\sigma$). En particulier, si $\mu$ n'a pas de parts égale à $1$, alors $\pi_{n}(A_{\mu,n})=C_{(\mu-1)\rightarrow n}$. Notons $\proj_{n}=\pi_{n} \circ \phi_{\infty,n}$ la projection de $\blg_{\infty}$ vers $\C\sym_{n}$, qu'on restreint en une projection de $\alg_{\infty}$ vers $Z(\C\sym_{n})$. Alors, étant données deux partitions $\lambda$ et $\mu$, on peut écrire : 
\begin{align*}
C_{\lambda \rightarrow n}\,*\,C_{\mu\rightarrow n} &=\proj_{n}(A_{\lambda+1})\,\proj_{n}(A_{\mu+1})=\proj_{n}(A_{\lambda+1}\,A_{\mu+1})=\proj_{n}\left( \sum_{|\rho|-\ell(\rho)\leq |\lambda|+|\mu|} g_{\lambda+1,\mu+1}^{\rho}\, A_{\rho}\right)\\
&=\sum_{|\rho|-\ell(\rho)\leq |\lambda|+|\mu|} g_{\lambda+1,\mu+1}^{\rho}\binom{n-|\rho|+m_{1}(\rho)}{m_{1}(\rho)}\,C_{(\rho-1)\rightarrow n}\,.
\end{align*}
Les partitions $\rho'=\rho-1$ que l'on obtient ont une taille $|\rho|'=|\rho|-\ell(\rho)\leq |\mu|+|\nu|$, et les coefficients devant les classes $C_{\rho'\rightarrow n}$ sont bien des polynômes en $n$. Ainsi, l'algèbre d'Ivanov-Kerov fournit une preuve très simple du théorème \ref{farahathigman}.
\begin{example}
L'identité générique présentée dans l'introduction du chapitre découle de l'identité suivante dans $\alg_{\infty}$ : $A_{(2)}*A_{(2)}=2\,A_{(2,2)}+3\,A_{(3)}+A_{(1,1)}$.
\end{example}
\bigskip\bigskip

L'objectif du chapitre est d'étendre ce résultat et ces méthodes dans la direction des algèbres d'Hecke de type A. En particulier, il convient de se demander s'il est possible de définir des <<~algèbres d'Hecke de permutations partielles~>> qui seraient des déformations des algèbres $\blg_{n}$. La réponse est non, et ceci est la source de nombreuses difficultés supplémentaires : en particulier, la base canonique de l'algèbre d'Hecke-Ivanov-Kerov que nous construirons dans la section \ref{genericnorm} ne sera pas directement reliée aux classes de conjugaison dans les centres $Z(\IH_{q}(\sym_{n}))$ --- ou plutôt les éléments de Geck-Rouquier, voir la section \ref{lascoux}. Essentiellement, c'est parce que la notion de support se comporte mal vis-à-vis de la structure de groupe de Coxeter de $\sym_{n}$ ; l'exemple qui suit illustre ce point.

\begin{example}
Dans $\sym_{4}$, considérons la permutation $(1,3)$ ; ses supports admissibles sont $\{1,3\}$, $\{1,2,3\}$, $\{1,3,4\}$ et $\{1,2,3,4\}$. D'autre part, les décompositions réduites de $(1,3)$ sont de longueur $3$, et il s'agit de :
$$(1,3)=s_{2}s_{1}s_{2}=s_{1}s_{2}s_{1}\,.$$
Si l'on attache des supports $S_{1}$ et $S_{2}$ aux transpositions $s_{1}$ et $s_{2}$, les seuls supports admissibles pour $(1,3)$ que l'on obtient dans les produits $(s_{2},S_{2})(s_{1},S_{1})(s_{2},S_{2})=(s_{1},S_{1})(s_{2},S_{2})(s_{1},S_{1})$ sont les supports <<~connexes~>> $\{1,2,3\}$ et $\{1,2,3,4\}$. Or, les éléments de base $T_{\sigma}$ des algèbres d'Hecke $\IH_{q}(\sym_{n})$ sont définis à partir de ces décompositions réduites ; par conséquent, on ne peut pas attacher n'importe quel support à un élément $T_{\sigma}$, et il n'y a pas d'équivalent Hecke de l'algèbre des permutations partielles $\blg_{n}$.
\end{example}
\bigskip\bigskip

L'obstruction constatée précédemment amène à considérer un autre modèle de permutation fibrées, à savoir, les permutations composées. On rappelle qu'une \textbf{composition} de taille $n$ est une suite finie $c=(c_{1},\ldots,c_{r})$ de nombres entiers strictement positifs tels que $\sum_{i=1}^{r} c_{i}=n$. On note $\comp_{n}$ l'ensemble des compositions de taille $n$, et $\comp=\bigsqcup_{n \in \N}\comp_{n}$ l'ensemble de toutes les compositions. Les \textbf{descentes} d'une composition $c=(c_{1},\ldots,c_{r})$ sont les entiers
$$c_{1},c_{1}+c_{2},c_{1}+c_{2}+c_{3},\ldots,c_{1}+\cdots+c_{r}\,.$$ 
Ces entiers caractérisent clairement la composition, et ils forment une partie de $\lle 1,n\rre$ qui contient $n$ ; par conséquent, le cardinal de $\comp_{n}$ est $2^{n-1}$. Notons que l'on inclut $|c|=n$ dans l'ensemble des descentes de $c$ ; cette convention n'est pas la plus courante. \bigskip
\bigskip

D'autre part, une \textbf{partition d'ensembles} de taille $n$ est la donnée d'une partition $\lle 1,n\rre=\bigsqcup_{i=1}^{r} \Pi_{i}$ en parts non vides. Les partitions d'ensembles sont ordonnées par la relation de raffinement :
$$\Pi \preceq \Phi \iff \forall \Pi_{i} \in \Pi,\,\,\exists \Phi_{j} \in \Phi,\,\,\Pi_{i}\subset \Phi_{j}$$
et il est bien connu que l'ensemble $\mathfrak{Q}_{n}$ des partitions d'ensembles de taille $n$ est un \textbf{treillis} pour cet ordre non total, c'est-à-dire que deux partitions $\Pi$ et $\Phi$ ont toujours une borne supérieure $\Pi \vee \Phi$ et une borne inférieure $\Pi \wedge \Phi$. En effet, une partition d'ensembles $\Pi$ de taille $n$ peut être vue comme une relation d'équivalence $R_{\Pi}$ sur $\lle 1,n\rre$ (de classes les parts de $\Pi$), et dans ce cas, $\Pi \wedge \Phi$ est la partition d'ensembles associée à la relation $R_{\Pi}\textsc{ et }R_{\Phi}$, et $\Pi \vee \Phi$ est la partition d'ensembles associée à la clôture transitive de $R_{\Pi} \textsc{ ou } R_{\Phi}$ en une relation d'équivalence.\bigskip
\bigskip

On peut associer à toute composition de taille $n$ une partition d'ensembles de $\lle 1,n\rre$ en intervalles : ainsi, on notera 
$$\Pi(c)=\lle 1,c_{1}\rre\sqcup \lle c_{1}+1,c_{1}+c_{2}\rre \sqcup \cdots \sqcup \lle c_{1}+\cdots+c_{r-1}+1,n\rre\,.$$
D'autre part, on peut associer à toute permutation de taille $n$ une partition d'ensembles $\orb(\sigma)$ dont les parts sont les orbites de $\sigma$. Alors, on appelle \textbf{permutation composée} de taille $n$ la donnée d'une permutation $\sigma \in \sym_{n}$ et d'une composition $c\in \comp_{n}$ telle que $\orb(\sigma)\preceq \Pi(c)$. Ceci revient à dire que $\sigma$ appartient au sous-groupe de Young $\sym_{c}$, et si $(\sigma,c)$ est une permutation composée et $j=\sigma(i)$, alors il n'y a pas de descente de $c$ entre $i$ et $j$.
\begin{example}
Le couple $(32154867,(5,3))$ est une permutation composée de taille $8$. On notera cette permutation composée $32154|867$.
\end{example}
Le produit de deux permutations composées est défini par $(\sigma,c)\,(\tau,d)=(\sigma\tau,c\vee d)$, où $c\vee d$ est l'unique composition telle que $\Pi(c\vee d)=\Pi(c)\vee \Pi(d)$ dans le treillis des partitions d'ensembles. Par exemple,
$$ 12|435|687 \times 321|54|867=42153|768\,.$$
On obtient ainsi un monoïde d'élément neutre $1|2|3|\cdots|n$, et nous noterons $\dlg_{n}$ l'algèbre complexe de ce monoïde. Elle se projette évidemment sur $\C\sym_{n}$ par le morphisme d'oubli de la composition.
\bigskip\bigskip

Par rapport aux algèbres de permutations partielles $\blg_{n}$, les algèbres de permutations composées $\dlg_{n}$ ont l'avantage d'admettre une présentation de type Coxeter. Plus précisément, notons $S_{i}$ la permutation composée $1|2|\cdots|i-1|i+1,i|i+2|\cdots|n$, et $I_{i}$ la permutation composée $1|2|\cdots|i-1|i,i+1|i+2|\cdots|n$. Alors, les $S_{i}$ et les $I_{i}$ entretiennent les relations suivantes :
\begin{align*}
&\forall i,\,\,\,(S_{i})^{2}=I_{i}\\
&\forall i,\,\,\,S_{i}S_{i+1}S_{i}=S_{i+1}S_{i}S_{i+1}\\
&\forall |j-i|\geq 2,\,\,\,S_{i}S_{j}=S_{j}S_{i}\\
&\forall i,j,\,\,\,S_{i}I_{j}=I_{j}S_{i} \\
&\forall i,j,\,\,\,I_{i}I_{j}=I_{j}I_{i}\\
&\forall i,\,\,\,S_{i}I_{i}=S_{i}\\
&\forall i,\,\,\,(I_{i})^{2}=I_{i}\,.
\end{align*}
Nous verrons dans quelques instants que ces relations constituent une présentation de l'algèbre $\dlg_{n}$. Plus généralement, notons $\IH(\dlg_{n})$ la $\C(q)$-algèbre engendrée par des éléments $(S_{i},I_{i})_{i \in \lle 1,n-1\rre}$ qui vérifient les six dernières relations de la liste précédente, et la relation quadratique 
$$(S_{i})^{2}=(q-1)\,S_{i}+q\,I_{i}\,.$$
Nous dirons que $\IH(\dlg_{n})$ est l'\textbf{algèbre d'Hecke (générique) des permutations composées} de taille $n$.
\begin{theorem}[Algèbre d'Hecke des permutations composées, \cite{Mel10c}]\label{heckecomposedpermutation}
L'algèbre $\IH(\dlg_{n})$ spécialise en $\dlg_{n}$ lorsque $q=1$ ; en $\IH(\sym_{n})$ lorsque $I_{1}=I_{2}=\cdots=I_{n-1}=1$ ; et en l'algèbre d'ordre inférieur $\IH(\dlg_{m})$ lorsque $m<n$ et $I_{m}=\cdots=I_{n-1}=0$ et $S_{m}=\cdots=S_{n-1}=0$. L'algèbre $\IH(\dlg_{n})$ admet une $\C(q)$-base $(T_{(\sigma,c)})$ indexée par les permutations composées de taille $n$, et il existe un isomorphisme de $\C(q)$-algèbres entre 
$$\IH(\dlg_{n})\quad\text{et}\quad\bigoplus_{c \in \comp_{n}} \IH(\sym_{c})\,,$$
où $\IH(\sym_{c})=\IH(\sym_{c_{1}}) \otimes \IH(\sym_{c_{2}}) \otimes \cdots \otimes \IH(\sym_{c_{r}})$ est la sous-algèbre de Young de $\IH(\sym_{n})$ associée à la composition $c$.
\end{theorem}
\begin{proof}
Si $c$ est une composition de taille $n$, nous appellerons \textbf{code} de $c$ l'ensemble des entiers de $\lle 1,n\rre$ qui ne sont pas des descentes de $c$. Par exemple, le code de $(3,2,3)$ est $\{1,2,4,6,7\}$. Si $c$ est une composition de taille $n$, notons $I_{c}$ le produit des $I_{i}$ avec $i$ dans le code de $c$ ; ainsi,
$$I_{(3,2,3)}=I_{1}I_{2}I_{4}I_{6}I_{7}\,.$$
Ces éléments sont des idempotents centraux, et $I_{c}$ correspond à la permutation composée $(\id_{\lle 1,n\rre},c)$. D'autre part, si $\sigma\in \sym_{n}$ admet pour décomposition réduite $\sigma=s_{i_{1}}s_{i_{2}}\cdots s_{i_{r}}$, notons comme dans le cas des algèbres d'Hecke
$$T_{\sigma}=S_{i_{1}}S_{i_{2}}\cdots S_{i_{r}}\,.$$
Finalement, si $(\sigma,c)$ est une permutation composée, notons $T_{(\sigma,c)}$ le produit $T_{\sigma}I_{c}$ dans l'algèbre $\IH(\dlg_{n})$. D'après le théorème de Matsumoto, si $\sigma$ admet deux expressions réduites $s_{i_{1}}s_{i_{2}}\cdots s_{i_{r}}$ et $s_{j_{1}}s_{j_{2}}\cdots s_{j_{r}}$, alors il est toujours possible de passer d'une expression à l'autre par des mouvements de tresse $s_{i}\,s_{i+1}\,s_{i}\leftrightarrow s_{i+1}\,s_{i}\,s_{i+1}$ et des commutations $s_{i}\,s_{j} \leftrightarrow s_{j}\,s_{i}$ lorsque $|j-i|\geq 2$. Comme les $S_{i}$ vérifient les mêmes relations de tresse dans l'algèbre $\IH(\dlg_{n})$, on conclut que les éléments $T_{\sigma}$ ne dépendent pas du choix d'expressions réduites. Maintenant, considérons un produit arbitraire $\Pi$ d'éléments $S_{i}$ et $I_{i}$ (dans n'importe quel ordre), et montrons qu'on peut toujours l'écrire comme combinaison linéaire d'éléments $T_{\sigma,c}$. Comme les $I_{i}$ sont des idempotents centraux, on peut toujours se ramener au cas où 
$$\Pi=S_{i_{1}}S_{i_{2}}\cdots S_{i_{p}}\,I_{c}\,,$$
où $c$ est une composition de taille $n$, et $s_{i_{1}}s_{i_{2}}\cdots s_{i_{p}}$ n'est pas \emph{a priori} une expression réduite. De plus, comme $S_{i}\,I_{i}=S_{i}$, on peut également supposer que le code de $c$ contient $\{i_{1},i_{2},\ldots,i_{p}\}$. Supposons l'expression $s_{i_{1}}s_{i_{2}}\cdots s_{i_{p}}$ non réduite ; alors, en utilisant des mouvements de tresse et des commutations, on peut transformer l'expression en un expression avec deux lettres consécutives identiques, \emph{i.e.},
$$\sigma=s_{j_{1}} \cdots s_{j_{k}}s_{j_{k+1}}\cdots s_{j_{p}}=s_{j_{1}} \cdots s_{j_{k-1}}s_{j_{k+2}}\cdots s_{j_{p}} \quad \text{car }j_{k}=j_{k+1}.$$
On peut appliquer les mêmes mouvements dans $\IH(\dlg_{n})$, et ainsi, $\Pi=S_{j_{1}}\cdots S_{j_{k}}S_{j_{k+1}}\cdots S_{j_{p}}\,I_{c}$ avec $j_{k}=j_{k+1}$ ; notons que le code de $c$ contient toujours $\{j_{1},j_{2},\ldots,j_{p}\}$. La relation quadratique dans $\IH(\dlg_{n})$ donne alors :
\begin{align*}
\Pi&=(q-1)\big\{S_{j_{1}}\cdots S_{j_{k-1}}S_{j_{k}}S_{j_{k+2}}\cdots S_{j_{p}}\,I_{c}\big\}+q\big\{S_{j_{1}}\cdots S_{j_{k-1}}I_{j_{k}}S_{j_{k+2}}\cdots S_{j_{p}}\,I_{c}\big\}\\
&=(q-1)\big\{S_{j_{1}}\cdots S_{j_{k-1}}S_{j_{k}}S_{j_{k+2}}\cdots S_{j_{p}}\,I_{c}\big\}+q\big\{S_{j_{1}}\cdots S_{j_{k-1}}S_{j_{k+2}}\cdots S_{j_{p}}\,I_{c}\big\}
\end{align*}
car $I_{j_{k}}I_{c}=I_{c}$. Par récurrence sur $p$, on conclut que $\Pi$ est une $\Z[q]$-combinaison linéaire d'éléments $T_{(\sigma,c)}$, avec la même composition $c$ pour toutes les permutations $\sigma$ de la combinaison linéaire. Ainsi, les produits réduits $T_{(\sigma,c)}$ engendrent linéairement l'algèbre $\IH(\dlg_{n})$.\bigskip

Si $c$ est une composition de taille $n$, notons $\psi_{c}$ le morphisme de $\C(q)$-algèbres  entre $\IH(\dlg_{n})$ et $\IH(\sym_{c})$ défini par :
$$\psi_{c}(S_{i})=\begin{cases} S_{i} &\text{si }i\text{ est dans le code de }c,\\
0 &\text{sinon,}
\end{cases}
\qquad;\qquad
\psi_{c}(I_{i})=\begin{cases} 1 &\text{si }i\text{ est dans le code de }c,\\
0 &\text{sinon.}
\end{cases}$$
Les éléments $\psi_{c}(S_{i})$ et $\psi_{c}(I_{i})$ vérifient dans $\IH(\sym_{c})$ les relations des éléments $S_{i}$ et $I_{i}$ dans l'algèbre $\IH(\dlg_{n})$ ; par conséquent, on peut bien définir un tel morphisme de $\C(q)$-algèbres $\psi_{c} : \IH(\dlg_{n})\to\IH(\sym_{c})$, et on voit facilement que :
$$\psi_{c}(T_{(\sigma,b)})=\begin{cases}T_{\sigma} &\text{si }\Pi(b) \preceq \Pi(c),\\
0&\text{sinon.}
\end{cases}$$
Considérons alors la somme directe de morphismes $\psi = \bigoplus_{c \in \comp_{n}}\psi_{c} : \IH(\dlg_{n}) \to \bigoplus_{c \in \comp_{n}}\IH(\sym_{c})$. Les vecteurs de base $[0,0,\ldots,(T_{\sigma}\in \IH(\sym_{c})),\ldots,0]$ de $\IH(\sym_{\comp_{n}})=\bigoplus_{c \in \comp_{n}}\IH(\sym_{c})$ seront notés $T_{\sigma \in \sym_{c}}$. D'après ce qui précède,
$$\psi(T_{(\sigma,c)})=\sum_{c \preceq d} T_{\sigma \in \sym_{d}}$$
pour toute permutation composée $(\sigma,c)$. Or, les compositions forment un sous-treillis (hypercube) du treillis des partitions d'ensembles, et la fonction de M\"obius\footnote{On renvoie au paragraphe \ref{fiber} pour des précisions sur cette terminologie, et aussi à \cite{Rot64}.} de ce sous-treillis est :
$$\mu(c,d)=(-1)^{\ell(c)-\ell(d)}\,.$$
Par suite, $\psi$ est un morphisme d'algèbres surjectif, puisque
$$\psi\left( \sum_{c \preceq d} \mu(c,d)\,T_{(\sigma,d)}\right)=T_{\sigma\in\sym_{c}}$$
pour toute permutation composée $(\sigma,c)$. Or, la dimension de $\IH(\sym_{\comp_{n}})$ est  exactement le nombre de permutations composées de taille $n$ ; par conséquent, le rang de la famille $(T_{(\sigma,c)})$ dans $\IH(\dlg_{n})$ est plus grand que ce nombre, et on conclut que $(T_{(\sigma,c)})$ forme une base de $\IH(\dlg_{n})$ lorsque $(\sigma,c)$ parcourt l'ensemble des permutations composées de taille $n$. De plus, $\psi$ est bien un isomorphisme d'algèbres.\bigskip

Finalement, étudions les spécialisations évoquées dans l'énoncé du théorème. Si $q$ est spécialisé en $1$, alors $T_{(\sigma,c)} \mapsto (\sigma,c)$ est un isomorphisme d'espaces vectoriels entre $\IH_{1}(\dlg_{n})$ et $\dlg_{n}$, et c'est un isomorphisme d'algèbres car il envoie une famille génératrice sur une famille génératrice vérifiant les mêmes relations. Ainsi, la liste des relations entre les $S_{i}$ et les $I_{i}$ donnée plus haut était bien une présentation de $\dlg_{n}$. Les autres spécialisations $\pi_{n} : \IH(\dlg_{n})\to \IH(\sym_{n})$ et $\phi_{n,m} : \IH(\dlg_{n})\to \IH(\dlg_{m})$ sont tout à fait triviales.
\end{proof}
\bigskip
\bigskip

Comme les spécialisations $\phi_{N\geq n} : \IH(\dlg_{N}) \to \IH(\dlg_{n})$ forment un système dirigé de morphismes d'algèbres, il existe une limite projective $\IH(\dlg_{\infty})=\varprojlim_{n \to \infty}\IH(\dlg_{n})$ dans la catégorie des algèbres\footnote{Par rapport à la construction d'Ivanov et Kerov, notons qu'on ne prend pas la limite projective dans la catégorie des algèbres filtrées ; en particulier, on ne pourra pas utiliser directement de filtration dans $\IH(\dlg_{\infty})$. Fort heureusement, il existera une filtration d'algèbre sur l'équivalent Hecke $\IH(\alg_{\infty})\subset \IH(\dlg_{\infty})$ de l'algèbre d'Ivanov-Kerov, mais ceci résultera d'un résultat de A. Francis et W. Wang, \emph{cf.} la section \ref{genericnorm}.} ; ses éléments sont les combinaisons linéaires formelles (éventuellement infinies) de $T_{(\sigma,c)}$ avec $\sigma$ permutation finie dans $\sym_{\infty}$, et $c$ composition infinie compatible avec $\sigma$ et telle que $c$ a une infinité de parts égales à $1$. Ainsi, on a construit une algèbre qui se projette sur toutes les algèbres d'Hecke $\IH(\sym_{n})$ (on notera comme précédemment $\proj_{n}$ les projections), et dans laquelle on peut espérer être en mesure d'établir des identités génériques. 
\bigskip\bigskip

Pour trouver dans $\IH(\dlg_{n})$ un équivalent de la sous-algèbre $\alg_{n}\subset \blg_{n}$, il conviendra de connaître la théorie des centres des algèbres d'Hecke de type A ; nous la rappellerons dans la section suivante. Avant cela, concluons ce paragraphe en étudiant les projections $\proj_{n}=\pi_{n} \circ \phi_{\infty,n} : \IH(\dlg_{\infty})\to \IH(\sym_{n})$. Deux éléments distincts $x$ et $y$ de l'algèbre d'Hecke des permutations composées peuvent avoir les mêmes projections $\proj_{n}(x)=\proj_{n}(y)$ dans toutes les algèbres d'Hecke : par exemple, c'est le cas de
$$T[21|34|5|6|\cdots]=S_{1}I_{1}I_{3} \quad\text{et}\quad T[2134|5|6|\cdots]=S_{1} I_{1}I_{2}I_{3}$$
qui se projettent sur $T_{(1,2)}$ si $n \geq 4$ et sur $0$ sinon. Néanmoins, le résultat devient vrai si l'on se restreint à la sous-algèbre $\IH(\dlg_{\infty}')$ engendrée par les $T_{(\sigma,c)}$ avec $c=(k,1^{\infty})$ :
\begin{proposition}[Séparation des vecteurs dans $\IH(\dlg_{\infty}')$]
Soit $\IH(\dlg_{n}')$ le sous-espace vectoriel de $\IH(\dlg_{n})$ engendré par les $T_{(\sigma,c)}$ avec $c$ de la forme $(k,1^{n-k})$. Ce sous-espace est une sous-algèbre, et dans la limite projective $\IH(\dlg_{\infty}') \subset \IH(\dlg_{\infty})$, les projections $\proj_{n}$ séparent les vecteurs : 
$$\forall x,y \in \IH(\dlg_{\infty}'),\quad (\forall n\in \N,\,\,\proj_{n}(x)=\proj_{n}(y)) \iff (x=y)\,.$$
\end{proposition}
\begin{proof}
Le supremum de deux compositions $(k,1^{n-k})$ et $(l,1^{n-l})$ est $(m,1^{n-m})$ avec $m=\max(k,l)$. Par conséquent, $\IH(\dlg_{n}')$ est bien une sous-algèbre de $\IH(\dlg_{n})$. Un élément $x$ de la limite projective s'écrit de manière unique :
$$x=\sum_{k=0}^{\infty}\sum_{\sigma \in \sym_{k}} a_{\sigma,k}(x)\,T_{\sigma,(k,1^{\infty})}\,.$$
Supposons que $x$ et $y$ admettent les mêmes projections dans toutes les algèbres d'Hecke $\IH(\sym_{n})$, et fixons une permutation $\sigma$. Il existe un entier $k$ minimal tel que $\sigma \in \sym_{k}$, et $a_{\sigma,k}(x)$ est le coefficient de $T_{\sigma}$ dans $\proj_{k}(x)$. Par conséquent, $a_{\sigma,k}(x)=a_{\sigma,k}(y)$. De même, le coefficient de $T_{\sigma}$ dans $\proj_{k+1}(x)$ est $a_{\sigma,k}(x)+a_{\sigma,k+1}(x)$, donc on a également $a_{\sigma,k}(x)+a_{\sigma,k+1}(x)=a_{\sigma,k}(y)+a_{\sigma,k+1}(y)$, et compte tenu de ce qui précède, $a_{\sigma,k+1}(x)=a_{\sigma,k+1}(y)$. Par récurrence sur $l$ et en utilisant le même argument, on voit donc que $a_{\sigma,k+l}(x)=a_{\sigma,k+l}(y)$ pour tout $l$. Ainsi, $x=y$, et on a montré que les projections séparaient les vecteurs de $\IH(\dlg_{\infty}')$.
\end{proof}
\noindent La sous-algèbre $\IH(\dlg_{\infty}')$ est bien plus simple que $\IH(\dlg_{\infty})$, mais sa présentation est moins naturelle ; c'est la raison pour laquelle nous avons commencé par introduire l'algèbre d'Hecke générale des permutations composées.
\bigskip

\section{Bases du centre d'une algèbre d'Hecke de type A}\label{lascoux}
Dans l'algèbre d'Hecke $\IH(\sym_{n})$, les classes de conjugaison $C_{\lambda}=\sum_{t(\sigma)=\lambda }T_{\sigma}$ ne sont plus forcément des éléments centraux ; par exemple, dans $\IH(\sym_{4})$, 
\begin{align*}T_{1}\,C_{(2,2)}&=T_{1}\,(T_{2143}+T_{3412}+T_{4321})=(q-1)\,T_{2143}+q\,T_{1243}+T_{3421}+(q-1)\,T_{4321}+q\,T_{4312}\,;\\
C_{(2,2)}\,T_{1}&=(T_{2143}+T_{3412}+T_{4321})\,T_{1}=(q-1)\,T_{2143}+q\,T_{1243}+T_{4312}+(q-1)\,T_{4321}+q\,T_{3421}\,.
\end{align*}
Le $q$-analogue du théorème de Farahat-Higman doit donc être énoncé pour d'autres éléments formant une base du centre $Z(\IH(\sym_{n}))$. En réalité, ce centre admet de nombreuses bases intéressantes, en particulier la base des \textbf{éléments de Geck-Rouquier} (\cite{GR97,Fra99,FW09}) et la base des \textbf{normes} (\cite{Jon90}). Nous rappelons ici ces points, en suivant pour l'essentiel \cite{Las06}, qui a l'intérêt de présenter très clairement les relations de changement de base.\bigskip
\bigskip

Les éléments de Geck-Rouquier sont des éléments $\Gamma_{\lambda}$ indexés par les partitions $\lambda \in \ym_{n}$, et caractérisés par les propriétés suivantes :
\begin{enumerate}[(i)]
\item Pour tout entier $n$, la famille $(\Gamma_{\lambda})_{\lambda \in \ym_{n}}$ est une base linéaire de $Z(\IH(\sym_{n}))$, qui spécialise en $(C_{\lambda})_{\lambda \in \ym_{n}}$ lorsque $q$ est spécialisé en la valeur $1$.
\item La différence $\Gamma_{\lambda}-C_{\lambda}$ ne met en jeu aucun élément $T_{\sigma}$ tel que $\sigma$ soit de longueur minimal dans sa classe de conjugaison.
\end{enumerate}
\begin{example}
Les éléments de Geck-Rouquier de taille $3$ sont :
$$\Gamma_{3}=T_{231}+T_{312}+(q-1)q^{-1}\,T_{321}\qquad;\qquad \Gamma_{2,1}=T_{213}+T_{132}+q^{-1}\,T_{321}\qquad;\qquad \Gamma_{1,1,1}=T_{123}$$
et ils forment une base de $Z(\IH(\sym_{3}))$. En taille supérieure, un exemple intéressant est constitué par :
\begin{align*}
\Gamma_{3,1}&=T_{1342}+T_{1423}+T_{2314}+T_{3124}+q^{-1}\,(T_{2431}+T_{4132}+T_{3214}+T_{4213})\\
&\quad+(q-1)q^{-1}\,(T_{1432}+T_{3214})+(q-1)q^{-2}\,(T_{3421}+T_{4312}+2\,T_{4231})+(q-1)^{2}q^{-3}\,T_{4321}
\end{align*}
Les termes de coefficient $1$ sont les quatre $3$-cycles de longueur minimale dans $\sym_{4}$ ; les termes dont le coefficient spécialise en $1$ lorsque $q=1$ sont les huit $3$-cycles de $\sym_{4}$ ; et tous les autres termes ne sont pas de longueur minimale dans leurs classes de conjugaison respectives, et ont leurs coefficients qui spécialisent en $0$ lorsque $q=1$.
\end{example}\bigskip

Notons $\Q_{\Z}[X]$ l'ensemble des polynômes à coefficients rationnels, et qui prennent des valeurs entières en les entiers : $\forall n \in \Z,\,\,P(n) \in \Z$. Il est bien connu que $\Q_{\Z}[X]$ est le $\Z$-module engendré par les coefficients binomiaux $\binom{X}{k}$ avec $k \geq0$. L'objectif du chapitre est la démonstration du théorème suivant :
\begin{theorem}[Hecke-Farahat-Higman, \cite{FW09,Mel10c}]\label{heckefarahathigman}
Le centre $Z(\IH(\sym_{n}))$ de l'algèbre d'Hecke est gradué par le degré $\deg \Gamma_{\lambda \rightarrow n}=|\lambda|$, c'est-à-dire que pour toutes partitions $\lambda$ et $\mu$,
$$\Gamma_{\lambda \rightarrow n}\,\Gamma_{\mu \rightarrow n}=\sum a_{\lambda\mu}^{\tau}(n,q)\,\Gamma_{\tau \rightarrow n}$$
la somme étant restreinte aux partitions $\tau$ telles que $|\tau|\leq |\lambda|+|\mu|$. De plus, les coefficients $a_{\lambda\mu}^{\tau}(n,q)$ sont des polynômes en $n$, $q$ et $q^{-1}$, et plus précisément, des éléments de $\Q_{\Z}[n]\otimes \Z[q,q^{-1}]$.
\end{theorem}\bigskip\bigskip

La première partie du théorème (le caractère gradué) a été prouvée par Francis et Wang ; nous rappellerons plus loin leur preuve. Les coefficients $a_{\lambda\mu}^{\tau}(n,q)$ peuvent pour leur part être étudiés à l'aide de l'algèbre d'Hecke des permutations composées, et plus précisément une sous-algèbre de celle-ci qui constitue un analogue Hecke de l'algèbre d'Ivanov-Kerov ; le reste du chapitre est consacré à la construction de cette algèbre. Avant cela, évoquons d'autres bases importantes des centres des algèbres d'Hecke de type A. Si $c$ est une composition de taille $n$ et $\sym_{c}$ est le sous-groupe de Young correspondant, les ensembles de classes $\sym_{c}\backslash \sym_{n}$ et $\sym_{n}/\sym_{c}$ admettent des représentants uniques de longueur minimale, appelés \textbf{représentants distingués} (voir \cite[\S2.1]{GP00} --- ce fait est même vrai pour les ensembles de doubles classes paraboliques). Dans ce qui suit, nous travaillerons plutôt avec les classes à droite, et dans ce contexte, les représentants distingués de $\sym_{c}\backslash \sym_{n}$  sont exactement les permutations dont les reculs\footnote{Un recul d'une permutation $\sigma(1)\sigma(2)\cdots\sigma(n)$ est une lettre $\sigma(i)$ telle que la lettre $\sigma(i)+1$ soit placée avant dans le mot. Par exemple, dans la permutation $6213475$, les reculs sont $1$ et $5$.} sont contenus dans l'ensemble des descentes de $c$. Par exemple, si $c=(2,3)$, alors
$$\sym_{(2,3)}\backslash \sym_{5}=\{12345,13245,13425,13452,31245,31425,31452,34125,34152,34512 \}=12 \sqcup\hspace{-1.4mm}\sqcup\, 345\,,$$
où $\sqcup\hspace{-0.85mm}\sqcup$ est le produit de mélange des mots ; en effet, le seul recul autorisé est entre $2$ et $3$.\bigskip
\bigskip

Soit $c$ une composition de taille $n$. Le morphisme de normalisation de $\IH(\sym_{c})$ vers $\IH(\sym_{n})$ est défini par :
$$N_{c}(h)=\sum_{\omega \in \sym_{c}\backslash \sym_{n}} q^{-\ell(\omega)}\,T_{\omega^{-1}}\,h\,T_{\omega}\,.$$
On peut montrer que si $h$ est central dans $\IH(\sym_{c})$, alors $N_{c}(h)$ est central dans $\IH(\sym_{n})$, voir \cite{Jon90}. Ceci permet de construire inductivement suffisamment d'éléments centraux dans $Z(\IH(\sym_{n}))$ pour obtenir des bases linéaires. En particulier, si 
$$N_{c}=N_{c}(1)=\sum_{\omega \in \sym_{c}\backslash \sym_{n}} q^{-\ell(\omega)}\,T_{\omega^{-1}}\,T_{\omega}\,,$$
alors $N_{c}$ est un élément central de l'algèbre d'Hecke, et on peut montrer que $N_{c}$ ne dépend de l'ordre des parts de $c$. De plus, les normes $N_{\lambda}$ avec $\lambda \in \ym_{n}$ forment une base du centre $Z(\IH(\sym_{n}))$. On doit à A. Lascoux\footnote{Le résultat a aussi été démontré à peu près au même moment par A. Francis et L. Jones, mais leur preuve est considérablement plus difficile ; voir \cite{FJ05}.} le résultat suivant :
\begin{proposition}[Changement de base entre éléments de Geck-Rouquier et normes, \cite{Las06}]\label{devilascoux}
Dans $Z(\IH(\sym_{n}))$, la base des éléments de Geck-Rouquier est reliée à la base des normes par la relation :
$$(\Gamma_{\lambda})_{\lambda \in \ym_{n}}=D\cdot P_{M,E}\cdot (N_{\mu})_{\mu \in \ym_{n}}$$
où $D$ est la matrice diagonale de coefficients $(q/q-1)^{n-\ell(\lambda)}$, et où $P_{M,E}$ est la matrice de changement de base entre fonctions monomiales et fonctions élémentaires de degré $n$.
\end{proposition}
\begin{example}
Les normes en taille $n=3$ sont :
\begin{align*}&N_{3}=T_{123}\qquad;\qquad N_{2,1}=3\,T_{123}+(q-1)q^{-1}\,(T_{213}+T_{132})+(q-1)q^{-2}\,T_{321}\qquad;\\
&N_{1,1,1}=6\,T_{123}+3(q-1)q^{-1}\,(T_{213}+T_{132})+(q-1)^{2}q^{-2}\,(T_{231}+T_{312})+(q^{3}-1)q^{-3}\,T_{321}
\end{align*}
et on peut vérifier que $\Gamma_{3}=q^{2}\,(q-1)^{-2}\,(3\,N_{3}-3\,N_{2,1}+N_{1,1,1})$, ce qui correspond à la relation $m_{3}=3\,e_{3}-3\,e_{2,1}+e_{1,1,1}$ dans l'algèbre des fonctions symétriques.
\end{example} 
\bigskip

\section[Normes génériques et construction d'une algèbre d'Hecke-Ivanov-Kerov]{Normes génériques et construction d'une algèbre d'Hecke-\\
-Ivanov-Kerov}\label{genericnorm}
La définition implicite des éléments de Geck-Rouquier donnée dans la section précédente ne permet pas de les décomposer simplement sur la base $(T_{\sigma})_{\sigma \in \sym_{n}}$, et en réalité, la fa\c con la plus simple de les calculer est de partir des normes et d'utiliser la proposition \ref{devilascoux}. Pour cette raison, il est beaucoup plus simple de contruire des \textbf{normes génériques} dans $\IH(\dlg_{\infty})$ que des éléments de Geck-Rouquier génériques (qui seraient des $q$-analogues des éléments $A_{\lambda}$ dans l'algèbre d'Ivanov-Kerov $\alg_{\infty}$). Si $c=(c_{1},\ldots,c_{r})$ est une composition de taille inférieure à $n$, notons $c \uparrow n$ la composition $(c_{1},\ldots,c_{r},n-|c|)$. On note d'autre part $H_{c}=I_{1}I_{2}\cdots I_{|c|-1}$ dans $\IH(\dlg_{\infty})$ ou dans $\IH(\dlg_{n})$ pour $n\geq |c|$. 
\begin{lemma}[Normes dans les algèbres d'Hecke de permutations composées]\label{whitetrashparty}
Si $|c|\leq n$, notons $M_{c,n}$ l'élément de $\IH(\dlg_{n})$ défini par :
$$
M_{c,n}=\sum_{\omega \in \sym_{c \uparrow n}\backslash\sym_{n}} q^{-\ell(\omega)}\,T_{\omega^{-1}}\,T_{\omega}\,H_{c}
$$
où les $T_{\omega}$ sont considérés comme des éléments de $\IH(\dlg_{n})$. On étend cette définition par $M_{c,n}=0$ si $|c|>n$. Alors, les $M_{c,n}$ vérifient la relation de compatibilité $\phi_{N,n}(M_{c,N})=M_{c,n}$, et leurs projections sur les algèbres d'Hecke sont les normes $N_{c\uparrow n}$ :
$$\pi_{n}(M_{c,n})=\begin{cases}N_{c\uparrow n}&\text{si }|c|\leq n,\\
0&\text{sinon.}\end{cases}$$
De plus, pour tout $n$ et tout $c$, $M_{c,n}$ appartient à $\IH(\dlg_{n}')$.
\end{lemma}
\begin{proof}
Supposons $|c|\leq n-1$ ; d'après ce qui précède, l'ensemble $R_{c,n}$ des mots-permutations $\omega$ de taille $n$ dont les reculs sont dans l'ensemble des descentes de $c$ constitue un système de représentants distingués de $\sym_{c \uparrow n}\backslash \sym_{n}$. Soit $\omega \in R_{c,n}$ un tel mot. Si $\omega(n)\neq n$, alors $T_\omega$ met en jeu l'élément $S_{n-1}$, donc $\phi_{n,n-1}(q^{-\ell(\omega)}\,T_{\omega^{-1}}\,T_{\omega})=0$. Au contraire, si $\omega(n)=n$, alors aucune décomposition réduite de $T_{\omega}$ ne met en jeu $S_{n-1}$, donc le terme correspondant dans la somme $M_{c,n}$ est préservé par $\phi_{n,n-1}$. Par conséquent, $\phi_{n,n-1}(M_{c,n})$ est la même somme que $M_{c,n}$, mais portant sur l'ensemble de mots $\omega \in R_{c,n}$ tels $\omega(n)=n$ ; autrement dit, sur $R_{c,n-1}$. Ainsi,
$$\phi_{n,n-1}(M_{c,n})=M_{c,n-1}$$
si $|c|\leq n-1$. Si $|c|=n$, alors $M_{c,n-1}=0$, et $\phi_{n,n-1}(M_{c,n})=0$, car $\phi_{n,n-1}(J_{c})=0$ ; et si $|c|>n$, alors $M_{c,n}$ et $M_{c,n-1}$ sont tous deux nuls, donc on a encore dans ce cas $\phi_{n,n-1}(M_{c,n})=M_{c,n-1}$. Finalement, comme 
$$\phi_{N,n}=\phi_{N,N-1}\circ \phi_{N-2,N-1}\circ \cdots \circ \phi_{n+1,n}\,,$$
on conclut que pour tous $N \geq n$ et toute composition $c$, $\phi_{N,n}(M_{c,N})=M_{c,n}$.\bigskip

Montrons maintenant que $M_{c,n} \in \IH(\dlg_{n}')$. Le résultat est trivial si $|c|>n$, car dans ce cas $M_{c,n}=0$ ; il l'est aussi si $|c|=n$, car dans ce cas $H_{c}=I_{(n)}$, et $d=(n)$ pour toute permutation composée $(\sigma,d)$ intervenant dans $M_{c,n}$. Supposons donc $|c|\leq n-1$. Comme les éléments de $\sym_{d}\backslash \sym_{|d|}$ peuvent être décrits par un produit de mélange, tout représentant distingué $\omega$ de $\sym_{c \uparrow n}\backslash \sym_{c}$ est le mélange d'un représentant distingué $\omega_{c}$ de $\sym_{c}\backslash \sym_{|c|}$ et du mot $|c|+1,|c|+2,\ldots,n$. Par exemple, $5613724$ est un représentant distingué d'une classe à droite dans $\sym_{(2,2,3)}\backslash \sym_{7}$, et c'est un mélange de $567$ avec le représentant distingué $1324$ d'une classe à droite de $\sym_{(2,2)}\backslash \sym_{4}$. Notons $s_{i_{1}}s_{i_{2}}\cdots s_{i_{r}}$ une expression réduite de $\omega_{c}$, et $j_{|c|+1},\ldots,j_{n}$ les positions de $|c|+1,\ldots,n$ dans le mot $\omega$. Alors, on peut montrer qu'une expression réduite de $\omega$ est :
$$ s_{i_{1}}s_{i_{2}}\cdots s_{i_{r}}\,\,(s_{|c|}s_{|c|-1}\cdots s_{j_{|c|+1}})\,(s_{|c|+1}s_{|c|}\cdots s_{j_{|c|+2}})\cdots (s_{n-1}s_{n-2}\cdots s_{j_{n}})\,.$$
Par exemple, $s_{2}$ est une expression réduite de $1324$, et 
$$s_{2}\,\,(s_{4}s_{3}s_{2}s_{1})\,(s_{5}s_{4}s_{3}s_{2})\,(s_{6}s_{5})$$
est une expression réduite de $5613724$. De ceci, on déduit que $T_{\omega} H_{c}=T_{(\omega,(k,1^{n-k}))}$, où $k$ est le plus grand entier de $\lle |c|+1,n\rre$ tel que $j_{k}<k$ --- on prend $k=|c|$ si $\omega$ et $\omega_{c}$ représentent la même permutation. Puis, la multiplication par $T_{\omega^{-1}}$ ne peut plus grossir plus la composition, donc $q^{-\ell(\omega)}\,T_{\omega^{-1}}\,T_{\omega}\,H_{c}$ est une combinaison linéaire d'éléments $T_{\tau,(k,1^{n-k})}$. AInsi, on a bien montré que $M_{c,n}$ appartenait à $\IH(\dlg_{n}')$. Finalement, comme la projection $\pi_{n}$ conserve les $T_{\omega}$ et envoie $H_{c}$ sur $1$, il est évident que $\pi_{n}(M_{c,n})=N_{c \uparrow n}$ si $|c|\leq n$. Dans ce qui suit, on convient que $N_{c \uparrow n}=0$ si $|c|<n$. Alors, on a également $\pi_{n}(M_{c,n})=N_{c \uparrow n}$ si $|c|<n$, car les deux termes de l'équation sont nuls.
\end{proof}
\bigskip

On appelle \textbf{norme générique} la somme infinie $M_{c} =\sum q^{-\ell(\omega)}\,T_{\omega^{-1}}\,T_{\omega}\,H_{c}$ effectuée sur toutes les permutations $\omega \in \sym_{\infty}$ qui ont leurs reculs dans l'ensemble des descentes de $c$. Ainsi, $M_{c}$ est l'unique élément de $\IH(\dlg_{\infty})$ dont les projections dans les algèbres $\IH(\dlg_{n})$ sont les $M_{c,n}$, et d'après le lemme \ref{whitetrashparty}, $M_{c}$ appartient à $\IH(\dlg_{\infty}')$, et est l'unique élément de cette sous-algèbre tel que $\proj_{n}(M_{c})=N_{c \uparrow n}$ pour tout $n$. Comme les $N_{c \uparrow n}$ ne dépendent pas de l'ordre des parts de $c$, il en va de même pour les $M_{c}$ ; ainsi, on considérera seulement dans la suite des éléments $M_{\lambda}$ avec $\lambda \in \ym$.
\begin{example}
L'élément $M_{(2)}$ admet les projections suivantes dans les premières algèbres d'Hecke de permutations composées :
\begin{align*}
&M_{(2),1}=0\qquad;\qquad M_{(2),2}=T_{12}\qquad;\\
&M_{(2),3}=T_{12|3}+2T_{123}+(1-q^{-1})\,(T_{132}+T_{213})+(q^{-1}-q^{-2})\,T_{321}\,.
\end{align*}
\end{example}
\bigskip
\bigskip

L'idée est maintenant d'utiliser les $M_{\lambda}$ pour établir des identités génériques entre les normes $N_{\lambda}$, puis d'utiliser une propriété polynomiale des matrices de passage $P_{M,E}$. Dans ce qui suit, nous noterons $(S_{i})^{-1}$ l'élément de $\IH(\dlg_{\infty})$ défini par $$(S_{i})^{-1}=q^{-1}\,S_{i}+(q^{-1}-1)\,I_{i}\,.$$ Le produit $S_{i}\,(S_{i})^{-1}=(S_{i})^{-1}\,S_{i}$ est égal à $I_{i}$, et par projection dans les algèbres d'Hecke, on retrouve $T_{i}\,(T_{i})^{-1}=1$ dans $\IH(\sym_{n})$ si $i<n$.
\begin{theorem}[Algèbre des normes génériques, \cite{Mel10c}]
Soit $\IH(\clg_{\infty})$ le sous-espace vectoriel de l'algèbre $\IH(\dlg_{\infty}')$ constitué des vecteurs $x$ tels que $S_{i}\,x\,(S_{i})^{-1}=I_{i}\,x$ pour tout $i$. Ce sous-espace vectoriel est aussi l'ensemble des combinaisons linéaires formelles de normes génériques $M_{\lambda}$, et c'est une sous-algèbre, qu'on appellera algèbre des normes génériques.
\end{theorem}\medskip

\begin{proof}
Si $I_{i}\,x=S_{i}\,x\,(S_{i})^{-1}$ et $I_{i}\,y=S_{i}\,y\,(S_{i})^{-1}$, alors
$$I_{i}\,xy=I_{i}\,x\,I_{i}\,y=S_{i}\,x\,(S_{i})^{-1}S_{i}\,y\,(S_{i})^{-1}=S_{i}\,x\,I_{i}\,y\,(S_{i})^{-1}=S_{i}\,xy\,(S_{i})^{-1}\,,$$
donc les éléments qui <<~commutent~>> avec $S_{i}$ dans $\IH(\dlg_{\infty})$ forment une sous-algèbre. En tant qu'intersection de sous-algèbres, le sous-espace $\IH(\clg_{\infty})\subset \IH(\dlg_{\infty}')\subset \IH(\dlg_{\infty})$ est donc une sous-algèbre. Montrons que les $M_{\lambda}$ appartiennent à $\IH(\clg_{\infty})$. Pour tout indice $i$, notons $\IH(\dlg_{\infty}^{(i)})$ la sous-algèbre de $\IH(\dlg_{\infty})$ engendrée linéairement par les $T_{(\sigma,c)}$ avec $c$ composition de la forme $(k,1^{\infty})\vee (1^{i-1},2,1^{\infty})$. Pour les mêmes raisons que dans le cas de $\IH(\dlg_{\infty}')$, les projections $\proj_{n}$ séparent les vecteurs de $\IH(\dlg_{\infty}^{(i)})$. Fixons alors un indice $i$ et une partition $\lambda \in \ym$. Les éléments $I_{i}\,M_{\lambda}$ et $S_{i}\,M_{\lambda}\,(S_{i})^{-1}$ appartiennent à $\IH(\dlg_{\infty}^{(i)})$, et leurs projections sont les mêmes, car pour tout $n$, $\proj_{n}(M_{\lambda})$ est une norme et en particulier un élément central de l'algèbre d'Hecke. Ainsi,  $I_{i}\,M_{\lambda}=S_{i}\,M_{\lambda}\,(S_{i})^{-1}$ pour tout $i$, donc $M_{\lambda}$ appartient à $\IH(\clg_{\infty})$. 
Fixons maintenant un élément $x$ de $\IH(\clg_{\infty})$ ; pour tout $x$, $\proj_{n}(x)= T_{i}\, \proj_{n}(x)\,(T_{i})^{-1}$, donc $\proj_{n}(x) \in Z(\IH(\sym_{n}))$ et est une combinaison linéaire de normes $N_{\lambda}$ :
$$\forall n \in \N,\,\,\,\proj_{n}(x)=\sum_{\lambda \in \ym_{n}} a_{\lambda}(x)\,N_{\lambda}\,.$$
D'après ce qui précède, on a le même résultat pour toute différence $x-\sum b_{\lambda}M_{\lambda}$. On peut donc construire par récurrence sur $n$ une combinaison linéaire infinie $S_{\infty}$ de $M_{\lambda}$ qui a les mêmes projections que $x$ :
\begin{align*} &\proj_{1}(x)=\sum_{|\lambda|=1} b_{\lambda}\,N_{\lambda} \quad\Rightarrow\quad \proj_{1}\bigg(x-\sum_{|\lambda|=1} b_{\lambda}\,M_{\lambda}\bigg)=0,\,\,\,S_{1}=\sum_{|\lambda|=1} b_{\lambda}\,M_{\lambda}\\
&\proj_{2}\left(x-S_{1}\right)=\sum_{|\lambda|=2} b_{\lambda}\,N_{\lambda} \quad\Rightarrow\quad \proj_{1,2}\bigg(x-\sum_{|\lambda|\leq 2} b_{\lambda}\,M_{\lambda}\bigg)=0,\,\,\,S_{2}=\sum_{|\lambda|\leq 2} b_{\lambda}\,M_{\lambda}\\
&\quad\vdots\\
&\proj_{n+1}\left(x-S_{n}\right)=\sum_{|\lambda|=n+1}b_{\lambda}\,N_{\lambda} \quad\Rightarrow\quad S_{n+1}=S_{n}+\sum_{|\lambda|=n+1}b_{\lambda}\,M_{\lambda}=\sum_{|\lambda|\leq n+1}b_{\lambda}\,M_{\lambda}\,.
\end{align*}
Alors, $S_{\infty}=\sum_{\lambda \in \ym}b_{\lambda}M_{\lambda}$ est dans $\IH(\clg_{\infty})$ et a les mêmes projections que $x$, donc est égal à $x$. On conclut que l'algèbre $\IH(\clg_{\infty})$ est bien l'algèbre des combinaisons linéaires formelles (éventuellement infinies) de normes génériques $M_{\lambda}$. 
\end{proof}
\noindent En particulier, un produit $M_{\lambda}\,*\,M_{\mu}$ dans $\IH(\clg_{\infty})$ est forcément une combinaison linéaire éventuellement infinie de $M_{\nu}$ :
$$\forall \lambda,\mu,\,\,\,M_{\lambda}\,*\,M_{\mu}=\sum \,g_{\lambda\mu}^{\nu} \,M_{\nu}$$
et comme les normes $N_{\lambda}$ sont définies sur $\Z[q,q^{-1}]$, les coefficients $g_{\lambda\mu}^{\nu}$ sont des éléments de $\Z[q,q^{-1}]$ (et même des polynômes symétriques en $q$ et $q^{-1}$). Il reste à montrer que la somme précédente est finie et porte en fait sur les partitions $\nu$ telles que $|\nu|\leq |\lambda|+|\mu|$ ; nous utiliserons pour ceci un argument de Francis et Wang (\cite{FW09}).\bigskip

\begin{example}
Le calcul des produits $(N_{(1)\uparrow n})^{2}$ se déduit de l'identité $$M_{1}\,*\,M_{1}=M_{1}+(q+1+q^{-1})\,M_{1,1}-(q+2+q^{-1})\,M_{2}$$ dans l'algèbre des normes génériques. On en déduit par exemple :
\begin{align*}(N_{1,1})^{2}&=(q+2+q^{-1})\,(N_{1,1}-N_{2})\,\,; \\
(N_{2,1})^{2}&=(q+1+q^{-1})\,(N_{1,1,1}-N_{2,1})\,\,;\\
(N_{3,1})^{2}&=N_{3,1}+(q+1+q^{-1})\,N_{2,1,1}-(q+2+q^{-1})\,N_{2,2}\,.
 \end{align*}
\end{example}
\bigskip
\bigskip

Dans ce qui suit, on convient que $e_{\lambda \uparrow n}=0$ si $|\lambda|>n$, et que $m_{\lambda \rightarrow n}=0$ si $|\lambda|+\ell(\lambda)>n$. Le point clef de notre raisonnement est la propriété très simple suivante :
\begin{lemma}[Caractère polynomial des matrices de passage $P_{M,E}$ et $P_{E,M}$]\label{evilwhistler}
Il existe des polynômes $P_{\lambda\mu}(n) \in \Q_{\Z}[n]$ et $Q_{\lambda\mu}(n) \in \Q_{\Z}[n]$ tels que 
$$\forall \lambda,n,\qquad m_{\lambda \rightarrow n} = \sum_{\mu' \leq_{d}\, \lambda} P_{\lambda\mu}(n) \,\,e_{\mu \uparrow n}\quad\mathrm{et}\quad 
                                          e_{\lambda \uparrow n} = \sum_{\mu \leq_{d}\, \lambda'} Q_{\lambda\mu}(n) \,\,m_{\mu \rightarrow n}\,,$$
où $\leq_{d}$ est la relation de domination des partitions (\emph{cf.} \cite[\S1.1]{Mac95}).
\end{lemma}
\begin{example}
Développons $m_{(2,1)\rightarrow n}$ dans la base des fonctions élémentaires, et $e_{(2,1)\uparrow n}$ dans la base des fonctions monomiales :
\begin{align*}m_{2,1\rightarrow n}&=e_{2,1 \uparrow n}-3\,e_{3\uparrow n}-(n-3)\,e_{1,1 \uparrow n} +(2n-8)\,e_{2 \uparrow n} +(2n-5)\,e_{1 \uparrow n}- n(n-4)\,e_{\uparrow n}\,\,;\\
e_{2,1\uparrow n}&=\frac{n(n-1)(n-2)}{2}\,m_{\rightarrow n}+\frac{(n-2)(3n-7)}{2}\,m_{1 \rightarrow n}+(3n-10)\,m_{1,1\rightarrow n}+3\,m_{1,1,1\rightarrow n}\\
&+(n-3)\,m_{2 \rightarrow n}+m_{2,1 \rightarrow n}\,.
\end{align*}
\end{example}\bigskip

La preuve du lemme est de nature algorithmique : on développe les $e_{\lambda \rightarrow n}$ sur un nombre suffisant de variables, et on réunit les monômes. Détaillons par exemple le calcul de $e_{2\uparrow n}=e_{n-2,2}$. Par définition,
$$e_{n-2,2}(X)=\sum_{\substack{i_{1}\neq i_{2}\neq \cdots \neq i_{n-2}\\ i_{n-1}\neq i_{n}}} x_{i_{1}}x_{i_{2}}\cdots x_{i_{n}}\,.$$
Si les indices $i_{n-1}$ et $i_{n}$ n'appartiennent pas à $\{i_{1},\ldots,i_{n-2}\}$, alors on reconnaît un terme $m_{1^{n}}$, et il y a $\binom{n}{2}$ termes de ce type, qui correspondent aux choix des deux indices $i_{n-1}$ et $i_{n}$ parmi $n$ indices. Si l'un des deux indices $i_{n-1}$ ou $i_{n}$ est dans $\{i_{1},\ldots,i_{n-2}\}$ (disons $i_{n-1}$), alors on reconnaît un terme $m_{21^{n-2}}$, et il y a $(n-2)$ termes de ce type, qui correspondent au choix de l'indice $i_{n}$ parmi les $n-2$ indices donnant un terme de degré $1$. Finalement, si les deux indices $i_{n-1}$ et $i_{n}$ sont dans $\{i_{1},\ldots,i_{n-2}\}$, alors on reconnaît un terme $m_{2^{2}1^{n-4}}$. On conclut que 
$$e_{2\uparrow n}=\frac{n(n-1)}{2}\,m_{\rightarrow n}+(n-2)\,m_{1\rightarrow n}+m_{1,1\rightarrow n}\,,$$
et il est clair que les mêmes arguments s'appliquent à toute fonction élémentaire $e_{\lambda \uparrow n}$. Les coefficients qui apparaissent sont des combinaisons linéaires de coefficients binomiaux, donc des éléments de $\Q_{\Z}[n]$ ; ceci démontre l'existence des polynômes $Q_{\lambda\mu}(n) \in \Q_{\Z}[n]$. Un raisonnement similaire permet de traiter le cas des $P_{\lambda\mu}$ ; on renvoie d'autre part à \cite[\S1.6, exemple 4.c]{Mac95} pour l'étude des coefficients des matrices de Kotska $K_{\lambda,\mu\rightarrow n}$.\bigskip
\bigskip

Comme dans le cas du groupe symétrique, on définit les $q$-éléments de Jucys-Murphy par la formule :
$$J_{k}=\sum_{j<k}q^{j-k+1} T_{(j,k)} = q^{-(k-1)}T_{(1,k)}+q^{-(k-2)}T_{(2,k)}+\cdots+T_{(k-1,k)}\,.$$
Ils commutent et engendrent une sous-algèbre abélienne maximale $\IH(GZ(n)) \subset \IH(\sym_{n})$. De plus, les polynômes symétriques en les $q$-éléments de Jucys-Murphy forment le centre de $\IH(\sym_{n})$, et comme dans la section \ref{jucysmurphy}, on a :
$$e_{r}(J_{1},\ldots,J_{n})=\sum_{|\lambda|=r} \Gamma_{\lambda\rightarrow n}\,.$$
Ces résultats ont été conjecturés par R. Dipper et G. James, et ils sont démontrés en particulier dans \cite{FG06}. Il sont liés au développement suivant des éléments de Geck-Rouquier :
$$\Gamma_{\lambda \rightarrow n} = m_{\lambda}(J_{1},\ldots,J_{n})+(\text{combinaison linéaire de $m_{\mu}(J_{1},\ldots,J_{n})$ avec } |\mu|<|\lambda|)\,.$$
D'autre part, chaque $m_{\lambda}(J_{1},\ldots,J_{n})$ est une combinaison linéaire (à coefficients linéaires dépendant de $n$) de termes $\Gamma_{\mu \rightarrow n}$ avec $|\mu|\leq |\lambda|$, \emph{cf.} \cite[lemme 2.3]{FW09}.
\begin{lemma}[Francis-Wang]
Pour toutes partitions $\lambda$ et $\mu$ et tout entier $n$,  $\Gamma_{\lambda \rightarrow n}\,\Gamma_{\mu \rightarrow n}$ est une combinaison linéaire d'éléments de Geck-Rouquier $\Gamma_{\tau \rightarrow n}$ avec $|\tau|\leq |\lambda|+|\mu|$. De même, $N_{\lambda \uparrow n}\,N_{\mu \uparrow n}$ est une combinaison linéaire de normes $N_{\tau \uparrow n}$ avec $|\tau| \leq |\lambda|+|\mu|$.
\end{lemma}
\begin{proof}
Si l'on exprime $\Gamma_{\lambda \rightarrow n}\,\Gamma_{\mu\rightarrow n}$ en termes des $q$-éléments de Jucys-Murphy, on obtient un terme $(m_{\lambda}m_{\mu})(J_{1},\ldots,J_{n})$, plus des polynômes symétriques en les $J_{k}$ de degré strictement inférieur à $|\lambda|+|\mu|$. Ainsi, $\Gamma_{\lambda \rightarrow n}\,\Gamma_{\mu\rightarrow n}$ est un polynôme symétrique en les $J_{k}$ de degré exactement $|\lambda|+|\mu|$, et en décomposant ce polynôme sur la base des fonctions monomiales, compte tenu du lemme de \cite{FW09} précité, on conclut que :
$$\Gamma_{\lambda\rightarrow n}\,*\,\Gamma_{\mu\rightarrow n}= \sum_{|\tau|\leq |\lambda|+|\mu|} a_{\lambda\mu}^{\tau}(n,q) \,\Gamma_{\tau \rightarrow n}\,,$$
ce qui constitue la première partie du théorème \ref{heckefarahathigman}. En combinant ce résultat et le lemme \ref{evilwhistler}, on obtient le même résultat pour les produits de normes $N_{\lambda\uparrow n}\,N_{\mu \uparrow n}$, car d'après la proposition \ref{devilascoux}, la matrice de passage entre les éléments de Geck-Rouquier et les normes est à une matrice diagonale près $P_{M,E}$.
\end{proof}
\noindent Comme les projections des $M_{\lambda}$ sont les $N_{\lambda \uparrow n}$, on en déduit facilement :
\begin{corollary}[Produits de normes génériques]
Le produit de deux normes génériques $M_{\lambda}$ et $M_{\mu}$ est une combinaison linéaire finie $\sum \,g_{\lambda\mu}^{\tau} \,M_{\tau}$ portant sur les partitions $\tau$ telles que $|\tau|\leq |\lambda|+|\mu|$. Ici, $g_{\lambda\mu}^{\tau}$ est un élément de $\Z[q,q^{-1}]$.
\end{corollary}
\begin{remark}
Une autre justification peut être donnée en introduisant les $q$-éléments de Jucys-Murphy génériques. Notons $L_{k}$ l'élément de $\IH(\dlg_{\infty}')$ défini par :
$$L_{k}=I_{(k)}\left(1+\frac{q-1}{q}\,J_{k}\right)=I_{(k)}\left(1+\frac{q-1}{q}\,\sum_{j<k}q^{j-k+1}\,T_{(j,k)}\right).$$
Ces éléments commutent, et on peut montrer qu'une norme générique $M_{\lambda}$ est un polynôme symétrique de degré $|\lambda|$ en les $q$-éléments de Jucys-Murphy génériques $L_{1},L_{2},\ldots,L_{n},\ldots$
\end{remark}
\bigskip
\bigskip

Nous pouvons finalement construire l'algèbre d'Hecke-Ivanov-Kerov. Notons $\IH(\alg_{\infty})$ le sous-espace vectoriel de $\IH(\clg_{\infty})$ constitué des combinaisons linéaires finies de normes génériques $M_{\lambda}$. D'après ce qui précède, $\IH(\alg_{\infty})$ est en fait une sous-algèbre, que nous appellerons \textbf{algèbre d'Hecke-Ivanov-Kerov}. Cette algèbre explique le caractère polynomial des coefficients $a_{\lambda\mu}^{\tau}(n,q)$. En effet,
\begin{align*}
\Gamma_{\lambda\rightarrow n}\,*\,\Gamma_{\mu \rightarrow n}&=\left(\frac{q}{q-1}\right)^{|\lambda|+|\mu|}\,\sum_{\rho,\sigma}P_{\lambda\rho}(n)\,P_{\mu\sigma}(n)\,\,N_{\rho\uparrow n}\,*\,N_{\sigma \uparrow n} \\
&=\left(\frac{q}{q-1}\right)^{|\lambda|+|\mu|}\,\sum_{\rho,\sigma,\tau}\,P_{\lambda\rho}(n)\,P_{\mu\sigma}(n)\,g_{\rho\sigma}^{\tau}\,N_{\tau \uparrow n} \\
&=\sum_{\rho,\sigma,\tau,\nu}\left(\frac{q}{q-1}\right)^{|\lambda|+|\mu|-|\nu|}\,P_{\lambda\rho}(n)\,P_{\mu\sigma}(n)\,g_{\rho\sigma}^{\tau}(q)\,Q_{\tau\nu}(n)\,\Gamma_{\nu \rightarrow n}=\sum_{\nu} a_{\lambda\mu}^{\nu}(n,q)\,\Gamma_{\nu \rightarrow n}
\end{align*}
avec $a_{\lambda\mu}^{\nu}(n,q)=(q/(q-1))^{|\lambda|+|\mu|-|\nu|}\,(P^{\otimes 2}(n)\,g(q)\,Q(n))_{\lambda\mu}^{\nu}$ en notation tensorielle. La seconde partie du théorème \ref{heckefarahathigman} est donc établie.
\begin{example}
En remarquant que $m_{1\rightarrow n}=e_{1\uparrow n} - n\,e_{\uparrow n}$ et que $e_{1\uparrow n} = n\,m_{\rightarrow n}+m_{1\rightarrow n}$, on obtient :
$$(\Gamma_{(1) \rightarrow n})^{2}=\binom{n}{2}q\,\Gamma_{(0) \rightarrow n}+(n-1)(q-1)\,\Gamma_{(1) \rightarrow n}+(q+q^{-1})\,\Gamma_{(1,1)\rightarrow  n}+(q+1+q^{-1})\,\Gamma_{(2)\rightarrow n}\,.$$
\end{example}\bigskip

Pour conclure, notons que notre algèbre $\IH(\alg_{\infty})$ ne contient pas d'<<~éléments de Geck-Rouquier génériques~>> ; pour cette raison, le calcul explicite des coefficients $a_{\lambda\mu}^{\tau}(n,q)$ est relativement difficile (en particulier, on ne connaît pas d'interprétation combinatoire de ces coefficients semblable à celle donnée par \cite[proposition 6.2]{IK99} dans le cas $q=1$). Une voie qui n'a pas été explorée dans cette direction est l'utilisation de l'algorithme de Francis (voir \cite{Fra99}) qui permet de calculer les éléments de Geck-Rouquier $\Gamma_{\lambda}$ à partir des classes $C_{\lambda}$ en effectuant des transformations récursives sur les combinaisons linéaires de $T_{\sigma}$ qui apparaissent. Ainsi, il est possible qu'il existe une algèbre d'éléments de Geck-Rouquier génériques qui fournisse une preuve plus directe du théorème \ref{heckefarahathigman} ; notre construction a pour sa part l'avantage de rentrer dans le cadre plus général présenté dans le chapitre suivant.

\chapter{Fibrés de semi-groupes et limites projectives}\label{bundle}

Dans le chapitre précédent, nous avons expliqué comment utiliser des permutations fibrées (par un support ou par une composition) pour construire des algèbres se projetant sur toutes les algèbres de groupes $\C\sym_{n}$, ou sur toutes les algèbres d'Iwahori-Hecke $\IH(\sym_{n})$. En réalité, on peut effectuer ce type de construction dans un contexte très général. Ainsi, si $(M,\cdot)$ est un semi-groupe fini et si $(L,\vee)$ est un semi-treillis fini, alors on peut sous certaines conditions construire un \textbf{semi-groupe fibré} $M \times_{I} L$ qui se projette sur $M$, et dont l'algèbre complexe vérifie des propriétés analogues aux algèbres de permutations partielles ou aux algèbres (d'Hecke) de permutations composées. De plus, étant données une tour de semi-groupes $(M_{n})_{n \in \N}$ et une famille croissante $(L_{n})_{n \in \N}$ de semi-treillis fibrant les $M_{n}$, on peut (de nouveau sous certaines conditions) utiliser la structure d'ensemble ordonné de $L=\bigcup_{n \in \N} L_{n}$ pour construire une <<~limite projective~>> des algèbres des semi-groupes $M_{n}$. Ainsi, pour construire des analogues de l'algèbre d'Ivanov-Kerov pour d'autres familles de groupes $(G_{n})_{n\in \N}$, il suffit de trouver des semi-treillis fibrant les groupes et munis d'actions naturelles de ceux-ci ; nous expliquons tout ceci dans les sections \ref{fiber} et \ref{chain}.
\bigskip\bigskip

La notion de fibrés de semi-groupes par des semi-treillis (ou plus généralement, de fibrés de semi-groupes par d'autres semi-groupes) n'est sans doute pas nouvelle, mais nous n'avons pas trouvé de références écrites sur cette construction, et les seuls exemples que nous connaissons (en dehors de ceux construits dans cette thèse) sont :\vspace{2mm}
\begin{enumerate}
\item les semigroupes de permutations partielles d'Ivanov et Kerov ;\vspace{2mm}
\item et les semigroupes de permutations <<~colorées~>> utilisés dans des problèmes d'énumération de factorisations de permutations, et en particulier le problème des \textbf{nombres de Hurwitz}.\vspace{2mm}
\end{enumerate}
Ainsi, dans la section \ref{hurwitz}, nous utilisons la théorie des fibrés de semi-groupes pour apporter un éclairage algébrique au problème des nombres de Hurwitz, que nous interprétons comme coefficients de structure d'algèbres de \textbf{permutations scindées}. En particulier, nous proposons un algorithme efficace en genre pour le calcul de ces nombres --- ce travail a été effectué en collaboration avec M. Sage. Ces méthodes étaient certainement déjà connues, mais nous avons jugé intéressant de les présenter dans le cadre algébrique nouveau (et particulièrement clair) des algèbres de permutations fibrées. Finalement, dans la section \ref{polyobs}, nous présentons des conjectures concernant les produits de classes de conjugaison dans les centres $Z(\C\GL(n,\For_{q}))$ des algèbres des groupes linéaires finis, et une tentative d'attaque du problème utilisant la notion d'\textbf{isomorphisme composé}.\bigskip

\section{Fibration d'un semi-groupe par un semi-treillis}\label{fiber}
On rappelle qu'un \textbf{semi-groupe} $M$ est un ensemble muni d'une loi de composition $\cdot$ qui est associative, c'est-à-dire que $m\cdot(n\cdot o)=(m \cdot n)\cdot o$ pour tous les éléments du semi-groupe. On parle de \textbf{monoïde} si $M$ admet un élément neutre (à droite et à gauche), et de groupe si de plus tout élément est inversible. D'autre part, un \textbf{semi-treillis} $L$ est un ensemble muni d'un ordre partiel $\leq$, et tel que la borne supérieure $x \vee y$ de deux éléments existe toujours. Alternativement, un semi-treillis peut être vu comme un semi-groupe de loi de composition $\vee$ qui est commutative ($\forall x,y,\,\,x \vee y=y\vee x$) et idempotente ($\forall x,\,\,x\vee x=x$). \bigskip
\bigskip

Un \textbf{sous-semi-groupe} $N$ d'un semi-groupe $(M,\cdot)$ est une partie $N$ qui est stable pour la loi de composition $\cdot$, c'est-à-dire que si $x$ et $y$ appartiennent à $N$, alors $x \cdot y$ appartient aussi à $N$ (on autorise en particulier le sous-semi-groupe vide). Nous noterons $\mathcal{S}(M)$ l'ensemble des sous-semi-groupes de $M$. D'autre part, un \textbf{idéal} $I$ d'un semi-treillis $(L,\vee)$ est une partie $I$ qui est close pour l'ordre $\leq$, c'est-à-dire que si $x \in I$ et $x\leq y$, alors $y \in I$. En particulier, $I$ est une partie close pour l'opération $\vee$. Nous noterons $\mathcal{I}(L)$ l'ensemble des idéaux de $L$. L'opération $\vee$ s'étend à $\mathcal{I}(L)$ : en effet, si $I_{1}$ et $I_{2}$ sont deux idéaux et si 
$$I_{1}\vee I_{2}= \{x\vee y\,\,|\,\,x \in I_{1},\,\,y \in I_{2}\}\,,$$
alors $I_{1}\vee I_{2}$ est un idéal dans $\mathcal{I}(L)$.
\bigskip

\begin{definition}[Fibré d'un semi-groupe par un semi-treillis]
Un fibré d'un semi-groupe $M$ par un semi-treillis $L$ est la donnée d'une application $m \in M \mapsto I_{m} \in \mathcal{I}(L)$ telle que 
$$\forall m_{1},m_{2} \in M,\,\,\,I_{m_{1}}\vee I_{m_{2}}\subset I_{m_{1}m_{2}}\,.$$
\end{definition}
\noindent \'Etant donné un fibré $(M,L,I)$ et un élément $l \in L$, notons $N_{l}$ la partie de $M$ définie par $N_{l}=\{m \in M\,\,|\,\,l \in I_{m}\}$. Cette partie est un sous-semi-groupe de $M$. En effet, si $m_{1}$ et $m_{2}$ appartiennent à $N_{l}$, alors $l \in I_{m_{1}}$ et $l\in I_{m_{2}}$, donc $$l=l \vee l \in I_{m_{1}}\vee I_{m_{2}}\subset I_{m_{1}m_{2}}\,,$$ et $m_{1}m_{2}\in N_{l}$. De plus, l'application $l \in L \mapsto N_{l} \in \mathcal{S}(M)$ est croissante pour l'inclusion des sous-semi-groupes. En effet, si $l_{1}\leq l_{2}$, alors si $m$ est un élément de $M$ tel que $l_{1}\in I_{m}$, on a aussi $l_{2} \in I_{m}$ ; par conséquent, $N_{l_{1}}\subset N_{l_{2}}$. Réciproquement, considérons une application $l \in L \mapsto N_{l} \in \mathcal{S}(M)$ qui est croissante, et montrons que les parties
$$I_{m}= \{l \in L\,\,|\,\,m \in N_{l} \}$$
sont des idéaux de $L$ vérifiant la condition de fibré. Si $l_{1}\in I_{m}$ et $l_{2} \geq l_{1}$, alors $m \in N_{l_{1}}\subset N_{l_{2}}$, donc $l_{2} \in I_{m}$ ; les parties $I_{m}$ sont donc bien des idéaux de $L$. Ensuite, si $l_{1} \in I_{m_{1}}$ et $l_{2} \in I_{m_{2}}$, alors 
$$m_{1} \in N_{l_{1}} \subset N_{l_{1}\vee l_{2}}\quad;\quad m_{2} \in N_{l_{2}} \subset N_{l_{1}\vee l_{2}}$$
et comme $N_{l_{1}\vee l_{2}}$ est un sous-semi-groupe, $m_{1}m_{2}$ appartient à $N_{l_{1}\vee l_{2}}$,  donc $l_{1}\vee l_{2}\in I_{m_{1}m_{2}}$. Ainsi, $I_{m_{1}}\vee I_{m_{2}}\subset I_{m_{1}m_{2}}$, c'est-à-dire que la condition de fibré est respectée. On conclut qu'un fibré de $M$ par $L$ peut alternativement être donné par une application
$l\in L \mapsto N_{l}\in \mathcal{S}(M)$ croissante pour l'ordre de $L$ et l'inclusion des sous-semi-groupes.
\bigskip\bigskip

Expliquons maintenant la terminologie de fibrés. On considère un fibré $(M,L,I)$, et on note $M \times_{I} L$ l'union disjointe $\bigsqcup_{m \in M}\{m\}\times I_{m}$. Ainsi,
$$M \times_{I} L = \bigsqcup_{m \in M}\{m\}\times I_{m}=\{(m,l) \,\,|\,\, l \in I_{m}\}= \{(m,l) \,\,|\,\, m \in N_{l}\} = \bigsqcup_{l \in L} N_{l}\times \{l\}\,.$$
L'ensemble $M\times_{I} L$ est muni d'une structure naturelle de semi-groupe pour le produit :
$$(m_{1},l_{1}) \times (m_{2},l_{2})=(m_{1}m_{2},l_{1}\vee l_{2})\,.$$
En effet, par hypothèse, si $l_{1} \in I_{m_{1}}$ et $l_{2}\in I_{m_{2}}$, alors $l_{1}\vee l_{2} \in I_{m_{1}m_{2}}$, donc le produit est bien dans $M\times_{I} L$ ; l'associativité du produit est évidente. De plus, la projection $$(m,l )\in M\times_{I}L \mapsto m \in M$$ est un morphisme de semi-groupes, et l'image réciproque de $m$ par cette projection s'identifie à $I_{m}$. On dira donc que $I_{m}$ est la \textbf{fibre} de $M \times_{I} L$ au-dessus de $m$.\bigskip

\begin{examples}
Notons $\sym_{n}^{p}$ le semi-groupe des permutations partielles de taille $n$, et $\sym_{n}^{c}$ le semi-groupe des permutations composées de taille $n$. Dans les deux cas, il s'agit de fibrés du groupe $\sym_{n}$. Ainsi, $\sym_{n}^{p}$ est le fibré de $\sym_{n}$ par le treillis $\mathfrak{P}_{n}=\mathfrak{P}(\lle 1,n\rre)$ des parties de $\lle 1,n\rre$ (muni de la relation d'ordre inclusion), avec :
\begin{align*}I_{\sigma \in \sym_{n}} &= \big\{\text{parties $A$ contenant le support essentiel de }\sigma, \text{ \emph{i.e.}, les cycles non triviaux}\big\}\,\,;\\
N_{A \in \mathfrak{P}_{n}}&= \sym(A)\,.\end{align*}
De même, $\sym_{n}^{c}$ est le fibré de $\sym_{n}$ par le treillis $\comp_{n}$ des compositions de taille $n$ (muni de la relation d'ordre de raffinement), avec :
\begin{align*}
I_{\sigma \in \sym_{n}}&=\big\{\text{compositions $c$ telles que }\Pi(c) \succeq \orb(\sigma)\big\}\,\,;\\
N_{c \in \comp_{n}}&= \sym_{c}=\sym_{c_{1}}\times \sym_{c_{2}}\times \cdots \times \sym_{c_{r}}\,.
\end{align*}
Nous verrons dans la section \ref{hurwitz} que $\sym_{n}^{c}$ est un sous-fibré du semi-groupe des permutations scindées, obtenu en rempla\c cant le treillis $\comp_{n}$ par le treillis $\mathfrak{Q}_{n}$ des partitions d'ensembles de $\lle 1,n\rre$. 
\end{examples}\medskip

\begin{remark}
Dans les deux cas présentés ci-dessus, le fibré obtenu est en fait un monoïde. Une condition suffisante pour qu'un fibré $M \times_{I} L$ soit un monoïde est la suivante : si $M$ a un neutre $e_{M}$ et $L$ a un élément minimal $0$, et si $I_{e_{M}}=L$, alors $(e_{M},0)$ est un élément neutre dans le semi-groupe fibré.
\end{remark}\medskip

\begin{example}
Soit $G$ un groupe et $E$ un $G$-ensemble, c'est-à-dire que $G$ agit sur $E$. On fixe un ensemble $L$ de parties de $E$ qui est stable par intersections, et tel que si $A \in L$ et $g \in G$, alors $g(A) \in L$. Par exemple, on peut imaginer que $E$ est un espace topologique sur lequel $G$ agit continuement ; alors, on peut prendre pour $L$ l'ensemble des fermés de $E$. En particulier, si $E$ et $G$ sont des ensembles finis, on prendra systématiquement $L=\mathfrak{P}(E)$. Alors, si l'on munit $L$ de l'ordre inverse de l'inclusion, l'application 
$$N : A \in L \mapsto \mathrm{Fix}(A) \subset G$$
vérifie les conditions de fibrés ; les éléments du fibrés sont les couples $(g,A)$ avec $g(a)=a$ pour tout $a$ dans $A$. Le semi-groupe des permutations partielles est un fibré de ce type, car le sous-groupe $\sym(A\subset \lle 1,n\rre)$ peut être vu comme le fixateur du complémentaire de $A$.
\end{example}
\bigskip\bigskip

Dans ce qui suit, on suppose $M$ et $L$ finis. En particulier, $L$ a un plus grand élément, qu'on notera $\infty$. Si $(E,\leq)$ est un ensemble partiellement ordonné fini, son \textbf{algèbre d'incidence} $\mathscr{I}(E)$ est l'ensemble des fonctions à valeurs entières définies sur les couples $a\leq b$ d'éléments de $E$, avec pour loi de multiplication 
$$(f * g)(a,b)=\sum_{a \leq x \leq b}f(a,x)\,g(x,b)\,.$$
En considérant le produit tensoriel $\mathscr{I}(E)\otimes_{\Z} A$, on peut définir l'algèbre d'incidence des fonctions à valeurs dans n'importe quel groupe abélien $A$. Le neutre de l'algèbre $\mathscr{I}(E)$ est la fonction de Dirac :
$$\delta(a,b)=\begin{cases} 1 & \text{si }a=b,\\
0 &\text{sinon.}\end{cases}
$$
Soit $\zeta(a,b)$ la fonction constante égale à $1$. Cette fonction est inversible dans $\mathscr{I}(E)$, d'inverse la \textbf{fonction de M\"obius} de $E$, \emph{cf.} \cite{Rot64}. Cette fonction $\mu(a,b)$ peut être définie récursivement par la relation :
$$\mu(a,b)=\begin{cases} 1&\text{si }a=b,\\
-\sum_{a\leq c< b} \mu(a,c)&\text{sinon}.
\end{cases}$$
Alors, si $f(a,b)=\sum_{a \leq c \leq b}g(a,c)$, alors $g(a,b)=\sum_{a \leq c \leq b} f(a,c)\,\mu(c,b)$. L'existence de cette fonction permet de généraliser les propositions \ref{reprpartialpermutation} et \ref{heckecomposedpermutation} du chapitre précédent :

\begin{proposition}[Théorie des représentations d'un fibré de semi-groupes]\label{fiberrepresentationtheory}
L'algèbre du semi-groupe fibré $\Z[M \times_{I} L]$ est abstraitement isomorphe à la somme directe
$$\bigoplus_{l \in L}\Z[N_{l}]\,.$$
En particulier, si les $N_{l}$ sont des groupes finis, alors $\C[M \times_{I} L]$ est une algèbre semi-simple dont les blocs correspondent à des blocs des $N_{l}$.
\end{proposition}
\begin{proof}
La preuve est essentiellement la même que précédemment. Si $m \in N_{l} \subset M$, on note $m_{N_{l}}$ le vecteur de $\bigoplus_{l \in L} \Z[N_{l}]$ dont la coordonnée dans $\Z[N_{l}]$ est $m$, et dont les autres coordonnées sont nulles. Lorsque $l$ parcourt $L$ et $m$ parcourt $N_{l}$, les $m_{N_{l}}$ forment une base de la somme directe. L'isomorphisme entre l'algèbre du fibré et la somme directe est défini par 
$$\psi : (m,l) \mapsto \sum_{l'\geq l}m_{N_{l'}}\,. $$
De nouveau, ceci fait sens, car l'application $N$ est croissante : si $(m,l) \in M \times_{I} L$ et $l' \geq l$, alors $m \in N_{l}$ et $N_{l} \leq N_{l'}$, donc $m \in N_{l'}$ et $m_{N_{l'}}$ existe. On vérifie trivialement que $\psi$ est compatible avec le produit du fibré. L'application réciproque de $\psi$ est l'extension $\Z$-linéaire de la correspondance 
$$\psi^{-1} : m_{N_{l}} \mapsto \sum_{l'\geq l} \mu(l,l')\,(m,l')\,.$$
En effet, pour tout élément $(m,l)$, on a bien :
\begin{align*}\psi^{-1}(\psi(m,l))&=\psi^{-1}\left(\sum_{l' \geq l} m_{N_{l'}}\right)=\sum_{l' \geq l}\psi^{-1}(m_{N_{l'}})=\sum_{l''\geq l'\geq l} \mu(l',l'')\,(m,l'')\\
&=\sum_{l'' \geq l} (\zeta*\mu)(l,l'')\,(m,l'')=\sum_{l'' \geq l} \delta(l,l'')\,(m,l'')=(m,l)\,.
\end{align*}
\noindent Par exemple, si $\dlg_{n}$ est l'algèbre des permutations composées de taille $n$, alors on retrouve immédiatement l'isomorphisme $\dlg_{n}\simeq \bigoplus_{c \in \comp_{n}} \C\sym_{c}$ ; nous avons vu dans le chapitre précédent que l'on avait en fait un isomorphisme entre les versions quantifiées $\IH(\dlg_{n})$ et $\bigoplus_{c\in \comp_{n}}\IH(\sym_{c})$.
\end{proof}\bigskip

Si $G$ est un groupe qui agit sur le semi-groupe $M$ et sur le semi-treillis $L$ de fa\c con compatible avec les lois de $M$ et $L$ et avec la fibration $I$, alors on peut comme dans le cas des permutations partielles prolonger l'action de $G$ au fibré $M \times_{I} L$, et restreindre la projection $\Z[M \times_{I} L]\to \Z[M]$ en une projection $\Z[M \times_{I } L]^{G}  \to \Z[M]^{G}$. Dans ce qui suit, on suppose que les actions de $G$ sur $M$ et $L$ vérifient :
\begin{align*}\forall g \in G&,\,\,\,\forall m_{1},m_{2} \in M,\,\,\,g\cdot(m_{1}m_{2})=(g\cdot m_{1})(g\cdot m_{2})\\
\forall g \in G&,\,\,\,\forall l_{1},l_{2} \in L,\,\,\,g\cdot(l_{1}\vee l_{2})=(g\cdot l_{1})\vee (g\cdot l_{2})\\
\forall g \in G&,\,\,\,\forall (m,l) \in M \times_{I} L,\,\,\,g\cdot (I_{m})=I_{g\cdot m}\,.
\end{align*}
La dernière hypothèse revient à dire que $I$ est $G$-équivariante pour l'action de $G$ sur $M$ et pour l'action de $G$ sur $\mathcal{I}(L)$. Alors, $G$ agit sur le fibré $M \times_{I} L$ par 
$$g\cdot(m,l)=(g\cdot m,g\cdot l)\,.$$
En effet, comme $l \in I_{m}$, $g \cdot l$ appartient bien à $g\cdot I_{m}=I_{g\cdot m}$ compte tenu de l'hypothèse d'équivariance. D'autre part, les deux premières hypothèses montrent que les sous-espaces d'invariants $\Z[M \times_{I}L ]^{G}$ et $\Z[M]^{G}$ sont des sous-algèbres de $\Z[M \times_{I} L]$ et de $\Z[M]$, et par définition de l'action de $G$ sur le fibré $M \times_{I}L$, la projection $\pi : \Z[M \times_{I} L] \to \Z[M]$ est un morphisme d'algèbres $G$-équivariant. En particulier, $\pi(\Z[M \times_{I} L])^{G} \subset \Z[M]^{G}$, et cette inclusion est en réalité une égalité. En effet, si $c=\sum_{m \in M} a_{m}\,m$ appartient à $\Z[M]^{G}$, alors $(c,\infty)=\sum_{m \in M}a_{m}\,(m,\infty)$ appartient à $\Z[M\times_{I}L]^{G}$, car $g\cdot \infty=\infty$ pour tout $g\in G$ ; et $\pi(c,\infty)=c$. On conclut que :
\begin{proposition}[Fibré de semi-groupes et actions de groupes]
Si $G$ est un groupe qui agit sur $M$ et $L$ de fa\c con compatible avec les lois de $M$ et $L$ et avec la fibration $I$, alors $G$ agit sur le fibré $M\times_{I}L$, et l'espace d'invariants $\Z[M \times_{I} L]^{G}\subset \Z[M \times_{I}L]$ est une sous-algèbre qui se projette surjectivement sur la sous-algèbre d'invariants $\Z[M]^{G}$.
\end{proposition}
\noindent En particulier, si $M=G$ et si l'action de $G$ sur lui-même est la conjugaison, alors une sous-algèbre commutative de l'algèbre du fibré $G\times_{I}L$ se projette sur le centre de l'algèbre de groupe $\Z[G]^{G}=Z(\Z G)$.
\bigskip

\section{Chaînes et construction de limites projectives}\label{chain}
Dans cette section, on considère un semi-groupe fibré (éventuellement infini) $M \times_{I} L$, et on suppose que $L$ admet un élément minimal $l_{0}$, et une \textbf{chaîne exhaustive} (infinie) $$C=(l_{0} \leq l_{1}\leq l_{2}\leq \cdots)\,,$$ c'est-à-dire une suite croissante telle que pour tout $l \in L$, il existe un indice $n$ tel que $l \leq l_{n}$. Les $N_{n}=N_{l_{n}}$ forment une suite croissante de sous-semi-groupes de $M$, et tout élément $m \in M$ est dans un $N_{n}$. En effet, si $m \in M$, choisissant $l \in I_{m}$, un certain $l_{n}$ est plus grand que $l$, donc dans l'idéal $I_{m}$ ; par suite, $m \in N_{n}$. Ainsi :
$$L=\bigcup_{n \in \N}\!\!\!\uparrow \,[l_{0},l_{n}]\qquad;\qquad M=\bigcup_{n \in \N}\!\!\!\uparrow \,N_{n}\,.$$
Pour $n \in \N$, notons $L_{n}$ l'intervalle $[l_{0},l_{n}]$, et pour $m \in N_{n}$, $I_{m,n}=I_{m}\cap L_{n}$. L'intervalle $L_{n}$ est un sous-treillis de $L$, et $I_{m,n}$ est un idéal de $L_{n}$. La chaîne $C$ fournit donc une suite de fibrés de semi-groupes $(N_{n} \times_{I_{n}} L_{n})_{n \in \N}$.
\begin{example}
Les permutations partielles rentrent parfaitement dans le cadre exposé ci-dessus, avec :
\begin{align*}&L=\Part_{\mathrm{finies}}(\N^{*})\qquad ;\qquad L_{n}=\Part(\lle 1,n\rre)\qquad;\\
&M=\sym_{\infty}\qquad\qquad\quad;\qquad N_{n}=\sym_{n}\qquad;\\
&I_{\sigma}=\{\text{parties de }\N^{*}\text{ contenant le support essentiel de }\sigma\}\qquad;\\
&I_{\sigma,n}=\{\text{parties de }\lle 1,n\rre\text{ contenant le support essentiel de }\sigma\}\,.
\end{align*}
\end{example}\bigskip

\noindent Si $n_{1} \geq n_{2}$, notons $\phi_{n_{1},n_{2}}$ la projection :
\begin{align*}
\phi_{n_{1},n_{2}} : \Z\left[N_{n_{1}}\times_{I_{n_{1}}}L_{n_{1}} \right]&\to \Z\left[N_{n_{2}}\times_{I_{n_{2}}}L_{n_{2}} \right]\\
(m,l) &\mapsto\begin{cases} (m,l) &\text{si } l \in L_{n_{2}},\\
0&\text{sinon.}
\end{cases}
\end{align*}
Cette application est bien définie, car si $l \in L_{n_{2}}$, alors $l \leq l_{n_{2}}$ et $l \in I_{m}$, donc $l_{n_{2}} \in I_{m}$ et $(m,l) \in N_{n_{2}} \times_{I_{n_{2}}} L_{n_{2}}$. De plus, on vérifie sans mal que $\phi_{n_{1},n_{2}}$ est un morphisme d'algèbres, et si $n_{1}\geq n_{2}\geq n_{3}$, alors on a évidemment la condition de compatibilité 
$$\phi_{n_{2},n_{3}}\circ \phi_{n_{1},n_{2}}=\phi_{n_{1},n_{3}}\,.$$
Ainsi, on a une famille dirigée de morphismes d'algèbres, ce qui permet d'envisager la limite projective $\varprojlim_{n \to \N} \Z[N_{n} \times_{I_{n}} L_{n}]$. On doit en réalité distinguer deux limites projectives : l'une dans la catégorie des algèbres, et l'autre dans la catégorie des algèbres filtrées. Si $(m,l) \in M\times_{I} L$, notons $\deg(m,l)$ le plus petit entier $n$ tel que $l \in L_{n}$. Cette graduation s'étend en une filtration d'algèbres sur $\Z[M\times_{I} L]$ et sur les $\Z[N_{n} \times_{I_{n}} L_{n}]$.
\begin{proposition}[Chaînes et limites projectives d'algèbres de semi-groupes fibrés]
On suppose que tous les sous-semi-groupes $N_{n}$ et tous les intervalles $L_{n}$ sont finis. Dans la catégorie des algèbres, $\varprojlim_{n \to \infty} \Z[N_{n} \times_{I_{n}} L_{n}]=\Z[[M \times_{I} L]]$. Dans la catégorie des algèbres filtrées, $\varprojlim_{n \to \infty} \Z[N_{n} \times_{I_{n}} L_{n}]=\Z[M \times_{I} L]$.
\end{proposition}
\begin{proof}
Les hypothèses réalisées sur la chaîne assure que tout élément $(m,l)$ de $M \times_{I} L$ est dans un certain $N_{n} \times_{I_{n}} L_{n}$. D'autre part, l'hypothèse de finitude des $N_{n}$ et des $L_{n}$ assure l'existence de l'algèbre de séries formelles $\Z[[M \times_{I} L]]$. Si l'on n'impose pas le degré fini pour les éléments de la limite projective $\varprojlim_{n \to \infty} \Z[N_{n} \times_{I_{n}} L_{n}]$, alors on obtient n'importe quelle combinaison formelle d'éléments fibrés $(m,l)$, c'est-à-dire n'importe quel élément de l'algèbre $\Z[[M\times_{I} L]]$. À l'inverse, si l'on impose un degré fini, on obtient n'importe quelle combinaison formelle finie d'éléments fibrés $(m,l)$, c'est-à-dire n'importe quel élément de $\Z[M \times_{I}L]$.
\end{proof}

\noindent Si l'on considère d'autres filtrations d'algèbres sur les $\Z[N_{n} \times_{I_{n}} L_{n}]$ et sur $\Z[[M \times_{I} L]]$, on peut construire d'autres limites projectives inclues dans $\Z[[M \times_{I} L]]$. Par exemple, dans la catégorie des algèbres filtrées, la limite projective des $\C\sym_{n}^{p}$ équipées du degré de Kerov, ou du degré canonique, est l'algèbre $\blg_{n}$ des permutations partielles.\bigskip\bigskip

Inversement, considérons une famille croissante $N_{0} \subset N_{1}\subset N_{2}\subset \cdots$ de semi-groupes. Les $N_{n}$ ont une limite directe naturelle $M= \bigcup_{n \in \N}\uparrow N_{n}$, mais la construction décrite ci-dessus permet également de construire une <<~limite projective~>> des algèbres $\Z[N_{n}]$, c'est-à-dire une algèbre générique qui se projette sur toutes ces algèbres. En effet, il suffit de trouver un semi-treillis $L$, une fibration $I: N \mapsto \mathcal{I}(L)$ et une chaîne exhaustive $C=(l_{0}\leq l_{1}\leq \cdots)$ telle que $N_{n}=N_{l_{n}}$ pour tout $n$ ; on a alors des projections $\phi_{\infty,n} : \Z[M \times_{I} L ]\to \Z[N_{n} \times_{I_{n}} L_{n}]$, et aussi des projections $\pi_{n} : \Z[N_{n} \times_{I_{n}} L_{n}]\to \Z[N_{n}]$, donc des projections $\proj_{n} : \Z[M \times_{I} L] \to \Z[N_{n}]$ pour tout $n$.  Nous utiliserons ce principe dans la section \ref{polyobs} pour construire une algèbre se projetant sur toutes les algèbres de groupes $\C\GL(n,\For_{q})$.

\bigskip
\bigskip

\section{Permutations scindées et nombres de Hurwitz}\label{hurwitz}

Le concept de permutations fibrées est également intéressant en dehors du problème de construction de limites projectives : en particulier, il est particulièrement utile pour le décompte de factorisations de permutations vérifiant certaines propriétés, et par exemple pour le calcul des nombres de Hurwitz. Si $n$ et $r$ sont deux entiers positifs, on appelle \textbf{constellation} (\emph{cf.} \cite[chapitre 1]{LZ04}) de degré $n$ et de longueur $r$ une famille $(\sigma^{(1)},\ldots,\sigma^{(r)})$ de permutations de $\sym_{n}$ telle que $\sigma^{(1)}\sigma^{(2)}\cdots\sigma^{(r)}=\id_{\lle 1,n\rre}$, et telle que le groupe engendré par les $\sigma^{(i)}$ agisse transitivement sur $\lle 1,n\rre$. Par exemple,
$$(2,4,3)(1,3)(2,4)(1,3,2)=\id_{\lle 1,4\rre}$$
est une constellation de longueur $4$ et de degré $4$. Le \textbf{type cyclique} d'une constellation est la famille $(\mu^{(1)},\ldots,\mu^{(r)})$ des types cycliques des pemutations $\sigma^{(i)}$ la composant ; ainsi, chaque $\mu^{(i)}$ est un élément de $\ym_{n}$. Réciproquement, étant donnée une famille $(\mu^{(1)},\ldots,\mu^{(r)})$ de partitions de même taille $n$, le \textbf{nombre de Hurwitz} de ces partitions est défini par :
$$H_{n}(\mu^{(1)},\ldots,\mu^{(r)})=\frac{1}{n!}\,\card\big\{\text{constellations de type }(\mu^{(1)},\ldots,\mu^{(r)})\text{ dans }\sym_{n}\big\}\,.$$
Plus généralement, considérons un entier $n$, des partitions $\mu^{(1)},\ldots,\mu^{(r)}$ de tailles plus petites que $n$, et un entier $g$ tel que
$$2(n+g-1)\geq \sum_{i=1}^{r}|\mu^{(i)}|-\ell(\mu^{(i)})=\sum_{i=1}^{r} R(\mu^{(i)})\,.$$
\begin{definition}[Nombres de Hurwitz]
Le nombre de Hurwitz $H_{n,g}(\mu^{(1)},\ldots,\mu^{(r)})$ est défini par :
$$H_{n,g}(\mu^{(1)},\ldots,\mu^{(r)})=H_{n}\bigg(\mu^{(1)}\sqcup 1^{n-|\mu^{(1)}|},\ldots,\mu^{(r)}\sqcup 1^{n-|\mu^{(r)}|},\underbrace{21^{n-2},21^{n-2},\ldots,21^{n-2}}_{k\text{ transpositions}}\bigg)$$
avec $k=2(n+g-1)-\sum_{i=1}^{r}R(\mu^{(i)})$.
\end{definition}
\begin{example}
Pour une seule partition $\mu$, $H_{n,g}(\mu)$ est aussi égal à $c_{\mu\sqcup1^{n-|\mu|}}/n!=1/z_{\mu\sqcup 1^{n-|\mu|}}$ fois le nombre de factorisations de la permutation standard de type $\mu$
$$\sigma_{\mu}=(1,2,\ldots,\mu_{1})(\mu_{1}+1,\ldots,\mu_{1}+\mu_{2})\,\cdots\,(\mu_{1}+\cdots+\mu_{r-1}+1,\ldots,\mu_{1}+\cdots+\mu_{r})$$
en produit de $k(\mu)=2(n+g-1)-R(\mu)$ transpositions de $\sym_{n}$ agissant transitivement. On en déduit les valeurs suivantes pour les premiers nombres de Hurwitz :
\begin{align*}
&H_{2,g}(2)=1/2\quad;\\
&H_{3,0}(3)=1\quad;\quad H_{3,1}=9\quad;\quad H_{3,2}(3)=81 \quad;\quad H_{3,g}=9^{g}\quad;\\
&H_{3,0}(2)=4\quad;\quad H_{3,1}(2)=40 \quad;\quad H_{3,2}(2)=364 \quad;\quad H_{3,g}=(9^{g+1}-1)/2\,.
\end{align*}
Nous expliquerons plus loin pourquoi $H_{n,g}(\mu)$ peut toujours être développé en combinaison linéaire rationnelle de puissances d'entiers $e^{g}$.
\end{example}\bigskip

Ces nombres comptent également des classes d'équivalence de revêtements ramifiés marqués de la sphère $\CP^{1}$. On rappelle qu'un \textbf{revêtement ramifié} de la sphère de degré $n$ est la donnée d'une surface de Riemann $X$ et d'une application holomorphe $f : X \to \CP^{1}$ telle que :\vspace{2mm}
\begin{enumerate}
\item Il existe une partie finie $S \subset X$ telle que $f: X \setminus S \mapsto \CP^{1} \setminus f(S)$ soit un revêtement holomorphe de degré $n$, c'est-à-dire une application localement équivalente à l'identité $z \mapsto z$ (dans des cartes holomorphes), et qui atteint chaque point $n$ fois.\vspace{2mm}
\item Autour de chaque $s\in S$, l'application $f$ est localement équivalente (dans des cartes holomorphes) à l'application $z \mapsto z^{d}$ avec $d\geq 1$.
\end{enumerate}
\vspace{-1mm}
\figcapt{\psset{unit=1mm}\pspicture(0,0)(120,60)
\psline(0,15)(15,0)(115,0)(100,15)
\psdots(35,10)(75,10)
\psellipse(35,10)(10,5)
\psellipse(75,10)(10,5)
\psellipse(35,25)(10,5)
\psellipse*[linecolor=white](35,29)(10,5)
\psellipse(35,29)(10,5)
\psellipse*[linecolor=white](35,33)(10,5)
\psellipse(35,33)(10,5)
\psellipse(35,50)(10,5)
\psline(35,33)(35,44.5)
\psline{->}(35,19.5)(35,11)
\psline{->}(75,21.8)(75,11)
\psline(75,44.5)(75,33)(68,31)
\psline(75,24)(68,22)
\psline[linestyle=dashed](75,33)(75,24)
\pscurve(68,31)(65,33)(85,32)(68,25)(65,27)(85,27)(68,22)
\psellipse(75,50)(10,5)
\psline(20,55)(18,53)(18,38.5)(16.5,37)(18,35.5)(18,22)(20,20)
\rput(6,37){$n$ feuillets}
\psline(90,38)(92,36)(92,31.5)(93.5,30)(92,28.5)(92,24)(90,22)
\rput(105,32){point de}
\rput(105,28.5){ramification}
\rput(33,10){$u$}
\rput(77,10){$t$}
\endpspicture}{Aspect local d'un revêtement ramifié de la sphère au-dessus d'un point régulier $u \notin f(S)$ et d'un point critique $t \in f(S)$.}{Aspect local d'un revêtement ramifié de la sphère}

L'entier $d(s)$ est appelé \textbf{degré de ramification} du revêtement en $s$, et l'\textbf{indice de ramification} est $i(s)=d(s)-1$. Le point $s$ est appelé \textbf{point de ramification} si $i(s)\geq 1$, et point de ramification simple si $i(s)=1$. Si $t \in f(S)$, le \textbf{type de ramification} de $t$ est la partition $\mu(t)=(d(s_{1}),\ldots,d(s_{l}))$ dont les parts sont les degrés de ramification des points $s_{i}$ dans la \textbf{fibre} de $t$, c'est-à-dire les points tels que $f(s_{i})=t$. En utilisant des triangulations et la formule d'Euler, on peut montrer que les indices de ramification et le genre de la surface $X$ sont reliés par la \textbf{formule de Riemann-Hurwitz} :
$$2-2g=\chi(X)=n\chi(\CP^{1})-\sum_{s\in S}i(s)=2n-\sum_{s \in S}i(s)=2n-\sum_{t \in f(S)} R(\mu(t))\,.$$
Étant données des partitions $\mu^{(1)},\ldots,\mu^{(r)}$ de tailles plus petites que $n$, on appelle \textbf{revêtement ramifié marqué} de $\CP^{1}$ de degré $n$, de genre $g$ et de type $(\mu^{(1)},\ldots,\mu^{(r)})$ la donnée :\vspace{2mm}
\begin{itemize} 
\item d'un revêtement ramifié $f : X \to \CP^{1}$ de degré $n$, \vspace{2mm}
\item de $r$ points critiques $t^{(1)},\ldots,t^{(r)}$, ces points ayant pour types de ramification respectifs $\mu^{(1)}\sqcup 1^{n-|\mu^{(1)}|}, \ldots,\mu^{(r)}\sqcup 1^{n-|\mu^{(r)}|}$,\vspace{2mm}
\item dans chaque fibre d'un point critique $t^{(i)}$, de $m_{1}(\mu^{(i)})$ points spéciaux marqués parmi les $n-|\mu^{(i)}|+m_{1}(\mu^{(i)})$ points non ramifiés,\vspace{2mm}
\item de $k=2(n+g-1)-\sum_{i=1}^{r}R(\mu^{(i)})$ points critiques supplémentaires $t^{(r+1)},\ldots,t^{(r+k)}$, chaque point étant de type de ramification $21^{n-2}$, c'est-à-dire qu'il y a un unique point de ramification simple dans la fibre.\vspace{2mm}
\end{itemize}
Par la formule de Riemann-Hurwitz, un tel revêtement est forcément de genre $g$. On compte ces revêtements avec le poids $1/\card \mathrm{Aut}(X,\CP^{1})$, où $\mathrm{Aut}(X,\CP^{1})$ désigne le groupe (fini) des automorphismes du revêtement ramifié $f : X \to \CP^{1}$. Alors :
\begin{proposition}[Nombres de Hurwitz et revêtements ramifiés marqués]
Le nombre de revêtements ramifiés marqués de degré $n$, genre $g$ et type $(\mu^{(1)},\ldots,\mu^{(r)})$ est exactement $H_{n,g}(\mu^{(1)},\ldots,\mu^{(r)})$.
\end{proposition}
\begin{proof}
Les arguments sont du même type que ceux donnés dans le paragraphe \ref{bdjgeom} pour les revêtements de Jucys-Murphy, et reposent sur la notion de \textbf{monodromie} dans les fibres d'un revêtement ramifié. Fixons un point régulier $z \in \CP^{1}$, et notons $t^{(1)},\ldots,t^{(r+k)}$ les points critiques d'un revêtement ramifié marqué. On dessine des boucles $\gamma^{(1)},\ldots,\gamma^{(r+k)}$ basées en $z$ et contournant les points $t^{(1)},\ldots,t^{(r+k)}$, \emph{cf.} la figure \ref{monodromy}.
\figcapt{\psset{unit=1mm}\pspicture(0,0)(100,50)
\pscircle(25,20){20}
\psdots(19.7,20)(35,15)(35,25)(15,31)
\pscurve{->}(20.5,19)(37,13)(37,16)(20.5,19.5)
\pscurve{->}(21,20.5)(37,24)(36,27)(20.5,21)
\pscurve{->}(20,21)(17,32)(13,32)(19,20.5)
\rput(25,29){$\cdots$}
\rput(18,19){$z$}
\rput(38,11){$\gamma^{(1)}$}
\rput(40,22.5){$\gamma^{(2)}$}
\rput(22.5,33.5){$\gamma^{(r+k)}$}
\pscurve(56,20)(54.5,30)(52,40)
\pscurve(52,40)(58,41)(64,43)
\pscurve(56,20)(62,21)(68,23)
\pscurve(71,20)(69.5,30)(67,40)
\pscurve(67,40)(73,41)(79,43)
\pscurve(71,20)(77,21)(83,23)
\psline[border=1mm,bordercolor=white]{->}(54,28)(36,25)
\psline[linestyle=dashed](55.5,28.25)(75,31.5)
\psdots(75,31.5)(60,29)
\psarc[linecolor=blue]{->}(60,29){3}{210}{5}
\psarc[linecolor=blue]{->}(75,31.5){3}{0}{170}
\psdots[linecolor=blue](57,28.5)(72,31)
\rput(75,28.5){\textcolor{blue}{$y_{\sigma^{(2)}(i)}$}}
\rput(57,30.5){\textcolor{blue}{$y_{i}$}}
\endpspicture}{Lien entre constellations et revêtements ramifiés marqués donné par la monodromie dans une fibre au-dessus d'un point régulier.\label{monodromy}}{Lien entre constellations et revêtements ramifiés donné par la monodromie}
Dans $\CP^{1}\setminus \{t^{(1)},\ldots,t^{(r+k)}\}$, le chemin concaténé $\Gamma=\gamma^{(1)}\gamma^{(2)}\cdots\gamma^{(r+k)}$ est homotope à $z$ --- il suffit de faire passer les boucles $\gamma^{(i)}$ de l'autre côté de la sphère. On relève les chemins $\gamma^{(i)}$ à $X \setminus S$ : ainsi, il existe des chemins $\delta^{(i)}$ sur $X \setminus S$ tels que $f(\delta^{(i)})=\gamma^{(i)}$. Si l'on suit une fois l'une de ces boucles $\delta^{(i)}$, on obtient une permutation $\sigma^{(i)}$ des préimages $y_{1},\ldots,y_{n}$ de $z$, et par hypothèse sur le revêtement ramifié, le type de cette permutation est $\mu^{(i)}\sqcup 1^{n-|\mu^{(i)}|}$ si $i \leq r$, et $21^{n-2}$ si $i \geq r+1$. Comme $\Gamma$ est homotope à $z$, son relevé $\delta^{(1)}\cdots \delta^{(r+k)}$ induit l'identité sur la fibre de $z$ ; ainsi, 
$$\sigma^{(1)}\cdots \sigma^{(r+k)}=\id_{\lle 1,n\rre}\,.$$
D'autre part, comme la surface $X$ est connexe, le groupe engendré par ces permutations agit transitivement sur $\lle1,n\rre$. Ainsi, tout revêtement ramifié marqué de degré $n$, genre $g$ et type $(\mu^{(1)},\ldots,\mu^{(r)})$ fournit une constellation de même type, et une analyse approfondie des symétries de cette correspondance conduit à l'identité annoncée.
\end{proof}
 \bigskip
 
 Dans la suite de cette section, nous expliquerons comment les permutations scindées, qui sont un type de permutations fibrées, permettent de calculer efficacement les nombres de Hurwitz généraux $H_{n,g}(\mu^{(1)},\ldots,\mu^{(r)})$, et plus particulièrement les nombres de Hurwitz d'une partition $H_{n,g}(\mu)$. Si $g=0$ et $\mu=(k)$, une formule exacte est connue depuis Hurwitz lui-même, voir \cite{Hur02}. Plus récemment, un lien a été établi entre le problème des nombres d'Hurwitz et la théorie de l'intersection des espaces de modules des courbes\footnote{Cette \textbf{théorie de Gromov-Witten} est aussi reliée à des modèles $2$-dimensionnels de la gravité quantique, voir \cite{Witten91,Zvon05}.}, voir par exemple \cite{ELSV,FP00,OP01,OP02}. En particulier, il existe une formule exacte, la \textbf{formule ELSV}, qui exprime les nombres de Hurwitz d'une partition en termes d'intégrales de formes différentielles sur les espaces de modules. Notons $\mathscr{M}_{n,g}$ l'espace de modules des courbes lisses de genre $g$ avec $n$ points marqués\footnote{On supposera $g\geq 1$ ou $g=0$ et $n\geq 3$, car sinon $\mathscr{M}_{n,g}$ n'est pas correctement défini.} $x_{1},\ldots,x_{n}$. Cet espace est un orbifold complexe (ou champ algébrique) de dimension complexe $3g-3+n$, et il admet une compactification $\mathscr{M}_{n,g}^{c}$ due à Deligne et Mumford et dont les points correspondent aux courbes stables de genre $g$ et avec $n$ points marqués. Sur une variété complexe compacte $X$ de faisceau structural $\obs$, rappelons que la suite exacte de faisceaux
 $$\begin{CD}
 1 @>>> \underline{\Z} @>>>\obs @>{\exp(2\I\pi\cdot)}>> \obs^{\times}@>>> 1
 \end{CD}$$
fournit une suite exacte longue entre les espaces de cohomologie $H^{k}(X,\Z)$, $H^{k}(X,\obs)$ et $H^{k}(X,\obs^{\times})$. En particulier, il existe une flèche canonique $H^{1}(X,\obs^{\times}) \to H^{2}(X,\Z)$ appelée \textbf{première classe de Chern}, et notée $c_{1}$. De plus, les classes de cohomologie de $H^{1}(X,\obs^{\times})$ correspondent aux classes d'équivalence conforme de fibrés en droites sur $X$, voir \cite[chapitre 1]{GH78} ; on peut donc considérer la classe de Chern de tout fibré en droites holomorphe sur $X$. Pour un orbifold compact tel que $\mathscr{M}_{n,g}^{c}$, la situation est semblable, à ceci près que la classe de Chern d'un fibré $c_{1}(L)$ prend maintenant ses valeurs dans l'espace de cohomologie rationnelle $H^{2}(\mathscr{M}_{n,g}^{c},\Q)$. D'autre part, on peut définir de fa\c con analogue les classes de Chern d'ordre supérieur $c_{k\geq 2}(L) \in H^{2k}(\mathscr{M}_{n,g}^{c},\Q)$. Considérons alors les \textbf{fibrés en droites tautologiques} $L_{i}$ de fibres
$$L_{i,(C;x_{1},\ldots,x_{n})}=T^{*}_{x_{i}}C\,,$$
 et le \textbf{fibré de Hodge} $\Lambda_{n,g}$, qui est le fibré vectoriel de dimension $2g$ dont les sections sont les $1$-formes holomorphes sur $\mathscr{M}_{n,g}^{c}$. On note finalement $\psi_{i}=c_{1}(L_{1})$, et $c(\Lambda^{\vee}_{n,g})$ la classe de Chern totale du fibré de Hodge dual $\Lambda^{\vee}_{n,g}$, \emph{i.e.}, la somme alternée de toutes les classes de Chern $c_{k}(\Lambda^{\vee}_{n,g})$. Alors :
 \begin{theorem}[Formule ELSV, \cite{ELSV}]
Si $\mu \in \ym_{n}$ et $\mathrm{Aut}(\mu)=\prod_{i\geq 1}m_{i}(\mu)!$, alors
$$ H_{n,g}(\mu)=\frac{k(\mu)!}{\mathrm{Aut}(\mu)}\,\int_{\mathscr{M}_{n,g}^{c}} \frac{c(\Lambda^{\vee}_{n,g})}{(1-\mu_{1}\psi_{1})(1-\mu_{2}\psi_{2})\cdots(1-\mu_{r}\psi_{r})}\,,$$
 où la fraction rationnelle de classes de cohomologie doit être interprétée comme une série formelle, et l'intégrale d'une classe de cohomologie $\xi$ vaut $0$ si $\xi$ n'est pas dans $H^{6g-6+2n}(\mathscr{M}_{n,g}^{c},\Q)$.
 \end{theorem} 
\noindent La preuve de ce résultat met en jeu l'application de Lyashko-Looijenga et ses extensions à des fibrés en cônes sur $\mathscr{M}_{n,g}^{c}$, en particulier l'espace de Hurwitz ; on renvoie à l'article original \cite{ELSV} pour plus de précisions. Notre approche est beaucoup moins ambitieuse : on cherche seulement à calculer efficacement (et exactement\footnote{Des expressions asymptotiques en $n$ et mettant en jeu les invariants de Gromov-Witten des espaces de modules ont été proposées par Kazarian et Zvonkine, voir par exemple \cite{Zvon04}.}) les nombres $H_{n,g}(\mu)$.
 \bigskip\bigskip

La difficulté principale pour le calcul des nombres de Hurwitz est la manipulation de l'hypothèse de transitivité ; sans celle-ci, il existe une formule très simple pour le nombre de factorisations remontant à Frobenius et mettant en jeu les caractères des groupes symétriques. Appelons \textbf{permutation scindée} de taille $n$ la donnée d'une paire $(\sigma,\pi)$, où $\sigma \in \sym_{n}$ est une permutation et $\pi \in \mathfrak{Q}_{n}$ est une partition d'ensemble de $\lle 1,n\rre$ qui est moins fine que la partition $\orb(\sigma)$ donnée par les orbites de $\sigma$. Ainsi, $j=\sigma(i)$ implique que $i$ et $j$ sont dans la même part de $\pi$.
\begin{example}
Le couple $s=((1,2)(3,4),\{1,2,6\}\sqcup\{3,4\}\sqcup\{5\})$ est une permutation scindée de taille $6$.
\end{example}\medskip
\noindent Notons $\sym_{n}^{s}$ l'ensemble des permutations scindées de taille $n$. Le nombre de permutations scindées $(\sigma,\pi)$ de taille $n$ s'écrit :
$$|\sym_{n}^{s}|=\sum_{\pi \in \mathfrak{Q}_{n}} \prod_{\pi_{i}\in \pi} |\pi_{i}|!=\sum_{\sigma \in \sym_{n}} B_{\text{nombre de cycles de }\sigma}=\sum_{\lambda \in \ym_{n}} |C_{\lambda}|\,B_{\ell(\lambda)}=\sum_{\lambda \in \ym_{n}} \frac{n!\,B_{\ell(\lambda)}}{z_{\lambda}}$$
où $B_{k}$ est le $k$-ième nombre de Bell --- le cardinal de $\mathfrak{Q}_{k}$ --- et correspond au nombre de fa\c cons de réunir des parts de $\orb(\sigma)$ pour obtenir $\pi$. Ainsi, les premiers cardinaux s'écrivent :
$$1,3,13,73,501,4051,37633,4596553,58941091,\ldots$$\bigskip

Dans le treillis des partitions d'ensembles, si $\sigma,\tau \in \sym_{n}$, alors $\orb(\sigma\tau) \preceq \orb(\sigma)\vee \orb(\tau)$. En effet, si $i$ et $j$ sont dans la même part de $\orb(\sigma\tau)$, alors 
$$j=\underbrace{\sigma\tau\sigma\tau\cdots \sigma\tau}_{k\text{ termes}}(i)\,,$$
donc il existe une chaîne $i_{0}=i,i_{1},\ldots,i_{2m}=j$ telle que $i_{l+1}=\sigma(i_{l})$ où $i_{l+1}=\tau(i_{l})$ pour tout $l$. Alors, pour la relation d'équivalence $\sim$ associée à $\orb(\sigma)\vee \orb(\tau)$, 
$$i=i_{0}\sim i_{1}\sim i_{2} \sim\cdots\sim i_{2m}=j\,,$$
donc on a bien $\orb(\sigma\tau) \preceq \orb(\sigma)\vee \orb(\tau)$. Par conséquent, l'application
\begin{align*}I : \sym_{n} &\to \mathcal{I}(\mathfrak{Q}_{n})\\
\sigma &\mapsto \big[\orb(\sigma),\{1,2,\ldots,n\}\big]
\end{align*}
vérifie la condition de fibré de semi-groupe, donc l'ensemble $\sym_{n}^{s}=\sym_{n}\times_{I} \mathfrak{Q}_{n}$ est un fibré\footnote{Remarquons que les permutations composées $(\sigma,c) \in \sym_{n}^{c}$ du chapitre précédent forment un sous-semi-groupe fibré de $\sym_{n}^{s}$.} sur $\sym_{n}$ avec pour loi multiplicative :
$$ (\sigma_{1},\pi_{1})\cdot(\sigma_{2},\pi_{2})= (\sigma_{1}\sigma_{2},\pi_{1}\vee \pi_{2})\,.$$
\begin{example}
Si $s=((1,2)(3,4),\{1,2,6\}\sqcup\{3,4\}\sqcup\{5\})$ et $t=((2,1,6),\{1,2,6\}\sqcup \{4,5\}\sqcup\{3\})$, alors 
$$s \cdot t = \big((1,6)(3,4),\{1,2,6\}\sqcup \{3,4,5\}\big)\,.$$
\end{example}\medskip

\noindent Compte tenu de la théorie développée dans la section \ref{fiber}, l'algèbre de semi-groupe $\C\sym_{n}^{s}$ se projette canoniquement sur $\C\sym_{n}$, et elle est isomorphe à la somme directe $\bigoplus_{\pi \in \mathfrak{Q}_{n}}\C\sym_{\pi}$, où $\sym_{\pi} \subset \sym_{n}$ désigne le sous-groupe produit des $\sym(\pi_{i})$ si $\pi=\pi_{1}\sqcup \cdots \sqcup \pi_{\ell(\pi)}$. L'isomorphisme s'écrit 
$$\psi(\sigma,\pi)=\sum_{\nu \succeq \pi}(\sigma_{\nu} \in \C\sym_{\nu})\,,$$
et sa réciproque met en jeu la fonction de M\"obius du réseau $\mathfrak{Q}_{n}$ :
$$\psi^{-1}(\sigma_{\pi})=\sum_{\nu \succeq \pi} \mu(\pi,\nu) \,(\sigma,\nu)$$
avec $\mu(\pi,\nu)=(-1)^{r-s}\,(2!)^{r_{3}}\,(3!)^{r_{4}}\,\cdots\,(n-1!)^{r_{n}}$, où $r$ est le nombre de parts de $\nu$, $s$ est le nombre de parts de $\pi$, et $r_{i}$ est le nombre de parts de $\nu$ qui se scindent en $i$ parts de $\pi$. Ceci permettra plus tard de donner une formule de type Frobenius pour les nombres de Hurwitz. \bigskip
\bigskip

Le groupe $\sym_{n}$ agit sur lui-même par conjugaison, et sur $\mathfrak{Q}_{n}$ en prenant les images des parties des partitions d'ensembles. Ces actions vérifient les axiomes donnés à la fin de la section \ref{fiber}, d'où une action sur le semi-groupe fibré $\sym_{n}^{s}$ :
$$\rho\cdot(\sigma,\pi)=(\rho\sigma\rho^{-1},\rho(\pi))$$
telle que l'espace des invariants dans $\C\sym_{n}^{s}$ forme une sous-algèbre se projetant sur $Z(\C\sym_{n})$. Deux permutations scindées $(\sigma,\pi)$ et $(\sigma',\pi')$ sont conjuguées sous cette action si et seulement si :\vspace{2mm}
\begin{enumerate}
\item Les partitions d'ensembles $\pi$ et $\pi'$ ont même profil, c'est-à-dire que les partitions d'entiers $\lambda$ et $\lambda'$ obtenues en ordonnant les tailles des parts de $\pi$ et de $\pi'$ sont identiques.\vspace{2mm}
\item Il existe une correspondance $\pi_{i} \leftrightarrow \pi_{i}'$ préservant les tailles des parts, et telle que pour tout indice $i$, les restrictions $\sigma_{i}$ et $\sigma_{i}'$ de $\sigma$ à la part $\pi_{i}$ et de $\sigma'$ à la part $\pi_{i}'$ aient même type cyclique.\vspace{2mm}
\end{enumerate} 
Par conséquent, les classes de $\sym_{n}$-conjugaison de permutations scindées sont indexées par les \textbf{multipartitions} de taille $n$, c'est-à-dire les multi-ensembles $\bbmu=\{\mu^{(1)},\ldots,\mu^{(s)}\}$ de partitions d'entiers tels que $\|\bbmu\|=\sum_{i=1}^{s}|\mu^{(i)}|=n$. L'ensemble des multipartitions de taille $n$ sera noté $\mathscr{Z}_{n}$ ; les cardinaux des premiers ensembles de multipartitions sont
$$1,3,6,14,27,58,111,223,424,817,\ldots$$
Le profil d'une multipartition $\bbmu=\{\mu^{(1)},\ldots,\mu^{(s)}\}$ est la partition $|\bbmu|$ dont les parts sont les tailles $|\mu^{(1)}|,\ldots,|\mu^{(s)}|$. En réunissant les multipartitions en fonction de leur profil, on peut montrer sans difficulté le lemme suivant :
\begin{lemma}[Majoration du nombre de multipartitions]\label{numberofmultipartitions}
Il existe une constante $B>1$ telle que $\card \mathscr{Z}_{n} \leq B^{n}$ pour tout $n$.
\end{lemma}
\begin{proof}
Asymptotiquement, le nombre de partitions de taille $n$ vérifie 
$$\card \ym_{n} \simeq \frac{\E^{\pi\sqrt{2n/3}}}{4\sqrt{3}\,n}\,,$$
voir \cite[chapitre VIII. 6]{FS09}. En particulier, $\card \ym_{n}$ est toujours majoré par $A^{\sqrt{n}}$ pour une certaine constante $A$. Maintenant,
\begin{align*}\card \mathscr{Z}_{n}&=\sum_{\lambda \in \ym_{n}} \card \{\bblambda \in \mathscr{Z}_{n}\,\,|\,\,|\bblambda|=\lambda\}\leq \sum_{\lambda \in \ym_{n}} \prod_{i=1}^{\ell(\lambda)}\card \ym_{\lambda_{i}} \leq \sum_{\lambda \in \ym_{n}} \prod_{i=1}^{\ell(\lambda)} A^{\sqrt{\lambda_{i}}}\\
&\leq \sum_{\lambda \in \ym_{n}} \prod_{i=1}^{\ell(\lambda)} A^{\lambda_{i}} \leq A^{n} \sum_{\lambda \in \ym_{n}} 1 \leq A^{n+\sqrt{n}} \leq (A^{2})^{n}\,,
\end{align*}
d'où le résultat avec $B=A^{2}$.
\end{proof}\bigskip

Par analogie avec ce qui a été fait dans la section \ref{ivanovkerov}, notons $\blg_{n}=\C\sym_{n}^{s}$, et $\alg_{n}$ la sous-algèbre des invariants, c'est-à-dire $(\C\sym_{n}^{s})^{\sym_{n}}$. Une base de $\alg_{n}$ est formée par les classes de pemutations scindées
$$C_{\bbmu}=\sum_{t(\sigma,\pi)=\bbmu} (\sigma,\bbmu)$$
indexées par les multipartitions $\bbmu \in \mathscr{Z}_{n}$. Si $\bbmu=\{\mu^{(1)},\ldots,\mu^{(s)}\}$ est une multipartition de taille $n$, notons $\bbmu^{\sqcup}=\mu^{(1)}\sqcup \mu^{(2)}\sqcup \cdots\sqcup \mu^{(s)}$ ; c'est une partition de $n$. L'image de $C_{\bbmu}$ par la projection $\pi_{n} : \blg_{n} \to \C\sym_{n}$ est un multiple de la classe de conjugaison $C_{\bbmu^{\sqcup}}$ ; ceci implique en particulier l'identité $\pi_{n}(\alg_{n})=\Z(\C\sym_{n})$, et le fait que $\alg_{n}$ est une sous-algèbre commutative. Réciproquement, étant donnée une partition $\mu \in \ym_{k \leq n}$, notons $\mu^{s}$ la multipartition de taille $n$ définie par 
$$\mu^{s}= \{(\mu_{1}),\ldots,(\mu_{\ell(\mu)}), (1),\ldots,(1)\}\,.$$
C'est la multipartition <<~la plus scindée possible~>> parmi celles telles que $\bbmu^{\sqcup}=\mu$.
\begin{theorem}[Nombres de Hurwitz et constantes de structure de l'algèbre des invariants]
Le nombre de Hurwitz $H_{n,g}(\mu^{(1)},\ldots,\mu^{(r)})$ est égal à $1/n!$ fois le coefficient de $C_{\{(1^{n)}\}}$ dans
$$C_{(\mu^{(1)})^{s}}\,C_{(\mu^{(2)})^{s}}\,\cdots\,C_{(\mu^{(r)})^{s}}\,(C_{(2)^{s}})^{k}\,.$$
\end{theorem}
\begin{proof}
Si $\mu$ est une partition de taille inférieure à $n$, la classe $C_{\mu^{s}}$ est simplement égale à 
$$\sum_{t(\sigma)=\mu\sqcup 1^{n-|\mu|}} (\sigma,\orb(\sigma))\,,$$
et sa projection sur $Z(\C\sym_{n})$ est $\pi_{n}(C_{\mu^{s}})=C_{\mu\sqcup 1^{n-|\mu|}}$. D'autre part, la classe $C_{\{(1^{n})\}}$ contient pour seule permutation scindée $(\id_{\lle 1,n\rre}, \lle 1,n\rre)$. Par conséquent, le coefficient de l'énoncé est égal à $1/n!$ fois le nombre de familles de permutations $(\sigma^{(1)},\ldots,\sigma^{(r)},\tau^{(1)},\ldots,\tau^{(k)})$ telles que :
\begin{align*}&\sigma^{(1)}\cdots\sigma^{(r)}\tau^{(1)}\cdots\tau^{(k)}=\id_{\lle 1,n\rre}\qquad;\qquad t(\sigma^{(i)})=\mu^{(i)} \sqcup 1^{n-|\mu^{(i)}|}\qquad;\qquad t(\tau^{(j)})=21^{n-2}\quad;\\
&\orb(\sigma^{(1)})\vee \cdots\vee \orb(\sigma^{(r)})\vee \orb(\tau^{(1)})\vee \cdots \vee \orb(\tau^{(k)})=\lle 1,n\rre\,.
\end{align*}
Il suffit donc de se convaincre que la dernière condition est équivalente à la transitivité du groupe engendré par les $\sigma^{(i)}$ et les $\tau^{(j)}$. Pla\c cons les entiers de $\lle 1,n\rre$ sur un cercle, et pour chaque permutation $\sigma$, dessinons le graphe $G(\sigma)$ dont les arêtes sont les $(i,j)$ avec $j=\sigma(i)$, voir la figure \ref{transitivity}.
\figcapt{\psset{unit=1mm}\pspicture(0,0)(80,40)
\pscircle[linewidth=0.25pt,linestyle=dashed](20,20){20}
\rput(53.7,31){$n=6$}
\rput(60,26){$\sigma^{(1)}=(1,5,6)$}
\rput(58.5,22){$\sigma^{(2)}=(2,4)$}
\rput(60,18){$\sigma^{(3)}=(2,6,3)$}
\psline[linecolor=blue]{->}(40,20)(10,2.68)
\psline[linecolor=blue]{->}(10,2.68)(30,2.68)
\psline[linecolor=blue]{->}(30,2.68)(40,20)
\psline[linecolor=violet]{<->}(30,37.32)(0,20)
\psline[linecolor=red]{->}(10,37.32)(30,37.32)
\psline[linecolor=red]{->}(30,37.32)(30,2.68)
\psline[linecolor=red]{->}(30,2.68)(10,37.32)
\psdots(40,20)(0,20)(10,37.32)(30,37.32)(10,2.68)(30,2.68)
\rput(42,20){$1$}
\rput(33,37.32){$2$}
\rput(7,37.32){$3$}
\rput(-2,20){$4$}
\rput(33,2.68){$6$}
\rput(7,2.68){$5$}
\endpspicture}{La transitivité du sous-groupe engendré par des permutations $\sigma \in S$ est équivalente à la connexité du graphe $\bigcup_{\sigma \in S} G(\sigma)$, et donc à la confition $\bigvee_{\sigma \in S} \,\orb(\sigma)=\lle 1,n\rre$.\label{transitivity}}{Transitivité du sous-groupe engendré par des permutations $\sigma \in S$}\bigskip

\noindent La transitivité du groupe engendré par les permutations est clairement équivalente à la con\-ne\-xité du graphe $(\bigcup_{i=1}^{r}G(\sigma^{(i)}))\cup(\bigcup_{j=1}^{k}G(\tau^{(j)}))$. Or, si ce graphe est connexe, alors tous les entiers de l'intervalle $\lle 1,n\rre$ sont en relation pour l'équivalence associée à la partition $\orb(\sigma^{(1)})\vee \cdots\vee \orb(\sigma^{(r)})\vee \orb(\tau^{(1)})\vee \cdots \vee \orb(\tau^{(k)})$, et réciproquement. 
\end{proof}\bigskip

\begin{corollary}[Nombres de Hurwitz d'une partition]\label{hurwitzonepartition}
Le nombre de Hurwitz $H_{n,g}(\mu)$ est égal à $1/z_{\mu \sqcup 1^{n-|\mu|}}$ fois le coefficient de $C_{\{\mu\sqcup1^{n-|\mu|}\}}$ dans $(C_{(2)^{s}})^{k}$.
\end{corollary}
\begin{proof}
Pour une seule partition, l'identité $\sigma\tau^{(1)}\cdots \tau^{(k)}=\id_{\lle 1,n\rre}$ se réécrit sous la forme $\sigma=\tau^{(k)}\cdots \tau^{(1)}$, et par conséquent, le groupe engendré $\langle \sigma,\tau^{(1)},\ldots,\tau^{(k)}\rangle$ est simplement $\langle \tau^{(1)},\ldots,\tau^{(k)}\rangle$. Ce sous-groupe est transitif si et seulement si $\bigvee \orb(\tau^{(j)})=\lle 1,n\rre$. Par conséquent, $n!\,H_{n,g}(\mu)$ est égal au nombre de termes dans $(C_{(2)^{s}})^{k}$ qui sont égal à un $(\sigma,\lle 1,n\rre)$ avec $t(\sigma)=\mu\sqcup1^{n-|\mu|}$. Comme $\alg_{n}$ est une sous-algèbre de $\blg_{n}$, ce nombre est aussi $\card C_{\mu \sqcup 1^{n-|\mu|}}$ fois le coefficient de $C_{\{\mu\sqcup1^{n-|\mu|}\}}$.
\end{proof}
\begin{example}
Dans $\alg_{3}$, les puissances successives de $C_{\{(2),(1)\}}$ sont :
\begin{align*}(C_{(2)^{s}})^{2m}&=C_{\{(1,1),(1)\}}+(3^{2m-1}-3)\,C_{\{(1,1,1)\}}+3^{2m-1}\,C_{\{(3)\}}\quad;\\
(C_{(2)^{s}})^{2m+1}&=C_{\{(2),(1)\}}+(3^{2m}-1)\,C_{\{(2,1)\}}\,.
\end{align*}
On retrouve ainsi les formules données précédemment pour $H_{3,g}(\mu)$.
\end{example}
 \bigskip
\bigskip

Compte tenu des deux propositions précédentes, le problème des nombres de Hurwitz se ramène à la compréhension des puissances de la classe $C_{(2)^{s}}$ dans $\alg_{n}$. En réalité, il va seulement s'agir de diagonaliser une matrice de taille $N \times N$, où $N=\card \mathscr{Z}_{n} \leq B^{n}$. Pour commencer, précisons quelque peu la structure des algèbres $\blg_{n}$. D'après ce qui précède, $\blg_{n}$ est isomorphe à la somme directe des $\C\sym_{\pi}$, $\pi$ parcourant $\mathfrak{Q}_{n}$ ; c'est donc une algèbre complexe semisimple, c'est-à-dire une somme directe de blocs qui sont des algèbres de matrices. Pour tout tableau standard $T$ de taille $n$, notons $e(T)$ l'\textbf{idempotent de Young} de $\C\sym_{n}$ (voir \cite{JK81}), qui est le produit :\vspace{2mm}
\begin{itemize}
\item de la somme de tous les éléments du sous-groupe de Young associé aux lignes de $T$\vspace{2mm}
\item et de la somme alternée de tous les éléments du sous-groupe de Young associé aux colonnes de $T$.
\end{itemize}\medskip

\begin{example}
Si $T=\young(25,134)$\,, alors la première somme est
$$R(T)=\big(\id+(1,3)+(1,4)+(3,4)+(1,3,4)+(1,4,3)\big)\,*\,\big(\id+(2,5)\big)\,,$$
et la seconde somme est
$$C(T)=(\id -(1,2))\,*\,(\id-(3,5))\,.$$
L'idempotent de Young est le produit $e(T)=R(T)\,C(T)$.
\end{example}
\noindent Les éléments $e(T)$ sont à des constantes multiplicatives près des idempotents dans l'algèbre $\C\sym_{n}$ : $e(T)^{2}=c(T)\,e(T)$ avec $c(T)=c(\lambda)=n!/\dim\lambda$. Plus précisément, les $f(T)=e(T)/c(T)$ forment un \textbf{système complet d'idempotents orthogonaux et indécomposables}, c'est-à-dire que :\vspace{2mm}
\begin{itemize} 
\item $1=\sum_{\lambda \in \ym_{n}}\sum_{T \in \mathrm{Std}(\lambda)} f(T)$ (complétude).\vspace{2mm}
\item $f(T)\,f(U)=0$ si $T \neq U$ (orthogonalité).\vspace{2mm}
\item chaque $f(T)$ ne peut être écrit comme somme d'idempotents orthogonaux non nuls (indécomposabilité).\vspace{2mm}
\end{itemize}
Si l'on réunit les idempotents de Young en fonction des formes $\lambda$ des tableaux $T$, alors on obtient un système complet d'\textbf{idempotents centraux}, \emph{i.e.}, $f(\lambda)=\sum_{T \in \mathrm{Std}(\lambda)}f(T)$ est dans $Z(\C\sym_{n})$ pour toute partition $\lambda$. Ces idempotents de Young fournissent l'isomorphisme d'algèbres $\C\sym_{n} \simeq \bigoplus_{\lambda \in \ym_{n}}\hendo(V_{\lambda})$ :
$$\hendo(V^{\lambda})=f(\lambda) \,\C\sym_{n}\,f(\lambda)=f(\lambda) \, \C\sym_{n} = \C\sym_{n}\,f(\lambda) = \bigoplus_{T,U \in \mathrm{Std}(\lambda)} f(T)\,\C\sym_{n}\,f(U)\,.$$
Pour tout tableau standard $U$ de forme $\lambda$, l'idéal à gauche $\C\sym_{n}\,f(U)$ est l'un des $\dim\lambda$ modules irréductibles de type $\lambda$ intervenant dans la décomposition $\C\sym_{n}\simeq \bigoplus_{\lambda \in \ym_{n}} (\dim \lambda)\,V^{\lambda}$, et les $f(T)\,\C\sym_{n}\,f(U)$ forment la \textbf{base de Young} du module. Dans ce qui suit, nous noterons $F(\lambda)=f(\lambda)\,\C\sym_{n}$ le bloc de type $\lambda$ de l'algèbre du groupe symétrique.\bigskip\bigskip

Pour un sous-groupe de Young $\sym_{\pi}=\sym(\pi_{1})\times \cdots \times \sym(\pi_{r})$ associé à une partition d'ensembles $\pi \in \mathfrak{Q}_{n}$, la décomposition de l'algèbre de groupes s'écrit 
$$\C\sym_{\pi}=\bigotimes_{i=1}^{r} \C\sym(\pi_{i})=\bigotimes_{i=1}^{r}\left(\bigoplus_{\lambda^{(i)} \in \ym_{|\pi_{i}|}} f_{\pi_{i}}(\lambda^{(i)}) \,\C\sym(\pi_{i}) \right)=\!\!\!\!\!\bigoplus_{\lambda^{(1)}\in \ym_{|\pi_{1}|},\ldots,\lambda^{(r)} \in \ym_{|\pi_{r}|}}\!\!\! \left(\bigotimes_{i=1}^{r} f_{\pi_{i}}(\lambda^{(i)})\right)\C\sym_{\pi}\,,$$
où par $f_{\pi_{i}}(\lambda^{(i)})$ on entend l'idempotent central de type $\lambda^{(i)}$ dans l'algèbre de groupes $\C\sym(\pi_{i})$. Ainsi, si $\Lambda=(\lambda^{(1)},\ldots,\lambda^{(r)}) \in \prod_{i=1}^{r}\ym_{|\pi_{i}|}$ et $f(\pi,\Lambda)=\bigotimes_{i=1}^{r} f_{\pi_{i}}(\lambda^{(i)})$, alors $\C\sym_{\pi}$ est isomorphe à la somme directe des blocs $F(\pi,\Lambda)=f(\pi,\Lambda)\,\C\sym_{\pi}$, et un tel $F(\pi,\lambda)$ est isomorphe au produit tensoriel $\bigotimes_{i=1}^{r}F(\lambda^{(i)})$. Comme $\blg_{n}\simeq \bigoplus_{\pi \in \mathfrak{Q}_{n}} \C\sym_{\pi}$, on conclut que :
\begin{lemma}[Décomposition de l'algèbre des permutations scindées]
On a la décomposition en blocs d'algèbres matricielles
$$\blg_{n}\simeq \bigoplus_{\pi,\Lambda}F(\pi,\Lambda)\,,$$
où $\pi$ parcourt $\mathfrak{Q}_{n}$, et $\Lambda$ parcourt $\ym_{\pi}=\prod_{i=1}^{\ell(\pi)}\ym_{|\pi_{i}|}$.
\end{lemma}
\noindent On peut ensuite réunir les blocs $F(\pi,\Lambda)$ en fonction de la multipartition $\bblambda=\{\lambda^{(1)},\ldots,\lambda^{(r)}\} \in \mathscr{Z}_{n}$ dont les partitions sont celles apparaissant dans la suite de partitions $\Lambda$. Ainsi, le nombre de blocs $F(\pi,\Lambda)$ de type $\bblambda \in \mathscr{Z}_{n}$ est
$$b(\bblambda)=\frac{n!}{|\bblambda|!\,\prod_{\lambda}m_{\lambda}(\bblambda)!}\,,$$
où $|\bblambda|!=\prod_{i=1}^{r}|\lambda^{(i)}|!$ si $\bblambda=\{\lambda^{(1)},\ldots,\lambda^{(r)}\}$, et $m_{\lambda}(\bblambda)$ est le nombre de partitions $\lambda$ dans la multipartition $\bblambda$. En effet, étant fixée une multipartition $\bblambda$, pour choisir un couple $(\pi,\Lambda)$ correspondant, il faut :\vspace{2mm}
\begin{enumerate}
\item En supposant $|\lambda^{(1)}|\geq \cdots \geq |\lambda^{(r)}|$, choisir une permutation $\Lambda=(\lambda^{\sigma(1)},\ldots,\lambda^{\sigma(r)})$ qui conserve la décroissance des tailles de partitions ; il y a
$$\mathrm{Aut}(\bblambda)=\prod_{k\geq1} \binom{m_{k}(|\bblambda|)}{\{m_{\lambda}(\bblambda)\}_{|\lambda|=k}}=\frac{\prod_{k\geq 1}m_{k}(|\bblambda|)!}{\prod_{\lambda} m_{\lambda}(\bblambda)!}\quad\text{possibilités}.$$
\item Puis, choisir une partition d'ensembles $\pi$ qui a pour profil $(|\pi_{1}|\geq \cdots \geq |\pi_{r}|)=(|\lambda^{\sigma(1)}|\geq \cdots \geq |\lambda^{\sigma(r)}|)$ ; il y a 
$$B(n,|\bblambda|)=\binom{n}{|\lambda^{(1)}|\cdots |\lambda^{(r)}|}\,\frac{1}{\prod_{k\geq 1}m_{k}(|\bblambda|)!}\quad\text{possibilités}.$$
\end{enumerate}
Ensuite, $b(\bblambda)=B(n,\bblambda)\,\mathrm{Aut}(\bblambda)$, et ainsi, $\blg_{n}\simeq \bigoplus_{\bblambda \in \mathscr{Z}_{n}} F(\bblambda)^{\oplus b(\bblambda)}$, où $F(\bblambda)$ désigne un bloc de type $\bblambda$.
\begin{example}
Si $n=4$, la décomposition de $\blg_{4}$ en blocs donne :\vspace{2mm}
\begin{itemize}
\item $B(4,(4))=1$ bloc pour les multipartitions $\{(1,1,1,1)\}$, $\{(2,1,1)\}$, $\{(2,2)\}$, $\{(3,1)\}$  et $\{(4)\}$. Les dimensions correspondantes sont $1$, $9$, $4$, $9$ et $1$.\vspace{2mm}
\item $B(4,(3,1))=4$ blocs pour les multipartitions $\{(3),(1)\}$, $\{(2,1),(1)\}$ et $\{(1,1,1),(1)\}$. Les dimensions de ces blocs sont $1$, $4$ et $1$.\vspace{2mm}
\item $B(4,(2,2))=3$ blocs pour les multipartitions $\{(2),(2)\}$ et $\{(1,1),(1,1)\}$, et $2 \times 3=6$ blocs pour la multipartition $\{(2),(1,1)\}$, car le groupe d'automorphismes de celle-ci est d'ordre $2$. Tous ces blocs sont de dimension $1$.\vspace{2mm}
\item $B(4,(2,1,1))=6$ blocs pour les multipartitions $\{(2),(1),(1)\}$ et $\{(1,1),(1),(1)\}$ ; ils sont de dimension $1$.\vspace{2mm}
\item et $B(4,(1,1,1,1))=1$ bloc de dimension $1$ pour la multipartition $\{(1),(1),(1),(1)\}$.  \vspace{2mm}
\end{itemize}
On retrouve bien $\dim \blg_{4}=73=1+9+4+9+1+4+16+4+3+3+6+6+6+1$.
\end{example}\bigskip
\bigskip

Revenons maintenant au problème du calcul des puissances de la classe des transpositions scindées $C_{(2)^{s}}$. Si $(a,b)$ est une transposition, on note
$$(a,b)^{s}=\left((a,b),\{a,b\}\sqcup \bigsqcup_{c\neq a,b}\{c\}\right)$$
la transposition scindée maximalement correspondante ; par définition, $C_{(2)^{s}}$ est la somme des $(a,b)^{s}$ pour $1\leq a <b \leq n$. Pour toute multipartition $\bbmu$, notons $\widetilde{C}_{\bbmu}=C_{\bbmu}/\card C_{\bbmu}$. Si $\proj_{\alg_{n}}$ désigne la projection 
$$x \in \blg_{n} \mapsto \frac{1}{n!} \sum_{\sigma \in \sym_{n}} \sigma\cdot x \in \alg_{n}\,,$$
alors la projection de tout élément $(\sigma,\pi)$ de type $\bbmu$ est $\widetilde{C}_{\bbmu}$. Par conséquent, pour toute classe $\bbmu$ et tout entier $m$,
$$C_{(2)^{s}}\,\widetilde{C}_{\bbmu}=\proj_{\alg_{n}}(C_{(2)^{s}})\,\proj_{n}(\sigma,\pi)=\proj_{\alg_{n}}\left(\sum_{a<b}(a,b)^{s}\right)\,\proj_{n}(\sigma,\pi)=\proj_{\alg_{n}}\left(\sum_{a<b}(a,b)^{s}(\sigma,\pi)\right)$$
si $(\sigma,\pi)$ est dans la classe $C_{\bbmu}$. Pour décomposer un produit $C_{(2)}^{s}\,\widetilde{C}_{\bbmu}$ dans la base des $\widetilde{C}_{\bblambda}$ de $\alg_{n}$, il suffit donc de déterminer les types des produits $(a,b)^{s}(\sigma,\pi)$ avec $(\sigma,\pi)$ permutation scindée fixée dans la classe de type $\bbmu$ :\vspace{2mm}
\begin{enumerate}
\item Si $a$ et $b$ sont dans des parts différentes $\pi_{i}$ et $\pi_{j}$ et dans des cycles de longueurs respectives $\mu^{(i)}_{k}$ et $\mu^{(j)}_{l}$, alors le type du produit est la multipartition $\bbmu[\mu^{(i)}\sqcup \mu^{(j)},\mu^{(i)}_{k}+\mu^{(j)}_{l}]$ obtenue à partir de $\bbmu$ en rempla\c cant les deux partitions $\mu^{(i)}$ et $\mu^{(j)}$ par leur union disjointe, puis les parts $\mu^{(i)}_{k}$ et $\mu^{(j)}_{l}$ de cette nouvelle partition par leur somme.\vspace{2mm}
\item Si $a$ et $b$ sont dans une même part $\pi_{i}$ et dans des cycles de longueur $\mu^{(i)}_{k}$ et $\mu^{(i)}_{l}$, alors le type du produit est la multipartition $\bbmu[\mu^{(i)}_{k}+\mu^{(i)}_{l}]$ obtenue à partir de $\bbmu$ en rempla\c cant dans la partition $\mu^{(i)}$ les parts   $\mu^{(i)}_{k}$ et $\mu^{(i)}_{l}$ par leur somme.\vspace{2mm}
\item Enfin, si $a$ et $b$ sont dans une même part $\pi_{i}$ et un même cycle de longueur $\mu^{(i)}_{k}$, alors le type du produit est la multipartition $\bbmu[\mu^{(i)}_{k}=d+d']$ obtenue en rempla\c cant dans la partition $\mu^{(i)}$ de $\bbmu$ la part $\mu^{(i)}_{k}$ par les deux distances $d$ et $d'$ entre $a$ et $b$ dans le cycle.\vspace{2mm}
\end{enumerate}
En énumérant les occurrences de chaque cas, on conclut que :
\begin{lemma}[Matrice de la multiplication par la classe des transpositions scindées]\label{matrixsplittransposition}
Pour toute multipartition $\bbmu=\{\mu^{(1)},\ldots,\mu^{(r)}\} \in \mathscr{Z}_{n}$, le produit $C_{(2)^{s}}\,*\,\widetilde{C}_{\bbmu}$ est égal à 
\begin{align*}&\sum_{1 \leq i<j \leq r}\sum_{k=1}^{\ell(\mu^{(i)})}\sum_{l=1}^{\ell(\mu^{(j)})} \mu^{(i)}_{k}\mu^{(j)}_{l}\,\widetilde{C}_{\bbmu[\mu^{(i)}\sqcup \mu^{(j)},\mu^{(i)}_{k}+\mu^{(j)}_{l}]}+\sum_{i=1}^{r} \sum_{1 \leq k <l \leq \ell(\mu^{(i)}} \mu^{(i)}_{k}\mu^{(i)}_{l}\,\widetilde{C}_{\bbmu[\mu^{(i)}_{k}+\mu^{(i)}_{l}]}\\
&+\sum_{i=1}^{r}\sum_{k=1}^{\ell(\mu^{(i)})}\sum_{d=1}^{\mu^{(i)}_{k}}\frac{\mu^{(i)}_{k}}{2}\,\widetilde{C}_{\bbmu[\mu^{(i)}_{k}=d+d']}\, .\end{align*}
\end{lemma}
\begin{example}
Dans la base $(\widetilde{C}_{\bbmu})_{\bbmu}$ de $\alg_{3}$ indexée par les multipartitions $\{(3)\}$, $\{(2,1)\}$, $\{(1,1,1)\}$, $\{(2),(1)\}$, $\{(1,1),(1)\}$ et $\{(1),(1),(1)\}$ (dans cet ordre), la matrice de la multiplication par $C_{(2)^{s}}$ est
$$\begin{pmatrix}
\text{\tiny 0}& 2 & \text{\tiny 0} & 2 & \text{\tiny 0} & \text{\tiny 0} \\ 
3 & \text{\tiny 0} & 3 & \text{\tiny 0} & 2 & \text{\tiny 0} \\ 
\text{\tiny 0} & 1 & \text{\tiny 0} & \text{\tiny 0} & \text{\tiny 0} & \text{\tiny 0} \\ 
\text{\tiny 0} & \text{\tiny 0} & \text{\tiny 0} & \text{\tiny 0} & 1 & 3 \\ 
\text{\tiny 0} &\text{\tiny 0}  &\text{\tiny 0}  & 1 & \text{\tiny 0} & \text{\tiny 0} \\ 
\text{\tiny 0}&\text{\tiny 0}  &\text{\tiny 0}  &\text{\tiny 0}  &\text{\tiny 0}  &\text{\tiny 0} 
\end{pmatrix}.
$$
\end{example}
\noindent En regardant la colonne correspondant à la multipartition $\{(1)^{n}\}$ dans les puissances de la matrice $M$ précédemment décrite, on obtient les nombres de Hurwitz d'une partition compte tenu du corollaire \ref{hurwitzonepartition}. Il reste donc à calculer efficacement les puissances de la matrice du lemme \ref{matrixsplittransposition}. Si $\bblambda=\{\lambda^{(1)},\ldots,\lambda^{(r)}\}$ est une multipartition, on note $c(\bblambda)=\sum_{i=1}^{r}c(\lambda^{(i)})$ son contenu total.
\begin{proposition}[Diagonalisation de la matrice de la classe des transpositions scindées]\label{diagosplittransposition}
Soit $M$ la matrice décrite par le lemme \ref{matrixsplittransposition}, c'est-à-dire la matrice de la multiplication par $C_{(2)^{s}}$ dans $\alg_{n}$. Cette matrice est diagonalisable, et son spectre est l'ensemble des contenus $c(\bblambda)$ avec $\bblambda$ parcourant $\mathscr{Z}_{n}$.
\end{proposition}
\begin{proof}
Fixons dans ce qui suit l'entier $n$. Dans la section \ref{jucysmurphy}, nous avons vu que les éléments de Jucys-Murphy $J_{1},\ldots,J_{n}$ agissaient diagonalement sur la base de Young du module $V^{\lambda}$, les valeurs propres étant les contenus des cases de $\lambda$. Or, la classe des transpositions $C_{21^{n-2}}$ est égale à la somme $p_{1}(J_{1},\ldots,J_{n})=J_{1}+J_{2}+\cdots+J_{n}$ ; son action sur $V^{\lambda}$ est donc la multiplication par la somme $c(\lambda)$ de tous les contenus des cases de $\lambda$. Notons $\proj_{\pi}$ la projection $\blg_{n}\to \C\sym_{\pi}$ tel que l'isomorphisme $\psi : \blg_{n}\to \sum_{\pi \in \mathfrak{Q}_{n}}$ soit égal à $\sum_{\pi \in \mathfrak{Q}_{n}} \proj_{\pi}$ ; $\proj_{\Lambda}$ la projection de $\C\sym_{\pi}$ sur le bloc $F(\pi,\Lambda)$ ; et $\proj_{\pi,\Lambda}$ la composée $\proj_{\Lambda}\circ \proj_{\pi}$. Comme 
$$\proj_{\pi}((a,b)^{s})=\begin{cases}(a,b)&\text{si $a$ et $b$ sont dans la même part de $\pi$},\\
0&\text{sinon},\end{cases}$$
l'image $\proj_{\pi}(C_{(2)^{s}})$ de la classe des transpositions scindées sur $\C\sym_{\pi}$ est simplement la somme $C_{21^{n-2},\pi}$ de toutes les transpositions de $\C\sym_{\pi}$. La multiplication par $C_{(2)^{s}}$ s'écrit donc dans $\blg_{n}$ :
$$C_{(2)^{s}}\,*\,x=\sum_{\pi \in \mathfrak{Q}_{n}} \proj_{\pi}(C_{(2)^{s}})\,*\,\proj_{\pi}(x)=\!\!\sum_{\pi \in \mathfrak{Q}_{n},\Lambda \in \ym_{\pi}}\! C_{21^{n-2},\pi}\,*\,\proj_{\pi,\Lambda}(x)=\!\!\sum_{\pi \in \mathfrak{Q}_{n},\Lambda \in \ym_{\pi}} \!c(\Lambda)\,\proj_{\pi,\Lambda}(x)\,.$$
Ainsi, l'endomorphisme $m$ associé à la multiplication par $C_{(2)^{s}}$ dans $\blg_{n}$ est diagonalisable, avec pour valeurs propres les contenus de multipartitions. Or, $m$ laisse stable $\alg_{n}$, car $\alg_{n}$ est une sous-algèbre de $\blg_{n}$. La restriction de $m$ à $\alg_{n}$ est donc encore diagonalisable, et la matrice $M$ est bien diagonalisable de spectre les contenus des multipartitions.
\end{proof}
\begin{remark}
En réalité, on peut même préciser les multiplicités des valeurs propres ; ainsi, le spectre avec multiplicités de $M$ est exactement $\{c(\bblambda)\,\,|\,\,\bblambda \in \mathscr{Z}_{n}\}$. Par exemple, les six multipartitions de taille $3$ ont pour contenus $3,0,-3,1,-1,0$, donc le spectre avec multiplicités de la matrice $M$ en taille $3$ est $\{-3,-1,0^{2},1,3\}$.
\end{remark}\bigskip

En combinant les lemmes précédents, on peut finalement donner un algorithme efficace pour le calcul des nombres de Hurwitz :
\begin{algorithm}[Calcul des nombres de Hurwitz d'une partition]
Pour calculer un nombre de Hurwitz d'une partition $H_{n,g}(\mu)$, il suffit de :
\begin{enumerate}
\item Lister toutes les multipartitions de taille $n$, et écrire la matrice $M$ de la multiplication par $C_{(2)^{s}}$, voir le lemme \ref{matrixsplittransposition}.
\item Diagonaliser la matrice $M$ ; comme on connaît à l'avance ses valeurs propres (\emph{cf.} la proposition \ref{diagosplittransposition}), ceci revient à résoudre des systèmes d'équations linéaires.
\item Calculer la puissance $(2(n+g-1)-R(\mu))$-ième de $M$, ce qui est facile une fois qu'elle est diagonalisée.
\end{enumerate}
Le nombre $H_{n,g}(\mu)$ se situe alors à l'intersection de la ligne indexée par $\{\mu\sqcup1^{n-|\mu|}\}$ et de la colonne indexée par $\{(1)^{n}\}$. Comme $\card \mathscr{Z}_{n}=O(B^{n})$ par le lemme \ref{numberofmultipartitions}, notre algorithme a pour complexité $O(C^{n})$ pour une certaine constante $C>1$.
\end{algorithm}
\begin{remark}
L'entier $n$ étant fixé, une fois connue une base de diagonalisation explicite de $M$, le calcul des puissances de $M$ est quasi immédiat, et on a dès lors accès simultanément à tous les nombres de Hurwitz $H_{n,g}(\mu)$. En particulier, l'algorithme est de complexité quasiment indépendante de $g$ (le genre est seulement mis en jeu dans le calcul de la puissance d'une matrice diagonale...).
\end{remark}\medskip
\begin{example}
L'algorithme précédent a été implémenté en \texttt{sage} ; il permet de calculer en quelques secondes $H_{9,3}([3,2])=11335243138639147728000$. Il n'est pas beaucoup plus dur de calculer $H_{10,100}([3,3])$ :
{\small\begin{align*}&78209797946099221469380408333253658389335110778578102493417366937278419420971892637983710 \\
&75560582522421501772573340373051838027863257564920539419318289349146733779503133393782164 \\
&00502995632992349968406352652755255329660159383909006457131068007080172851654851060277221 \\
&485502282528772332192003548685671573635386956399466111869724001404563147200000 \,.\end{align*}}
\end{example}\bigskip

Pour conclure, donnons une version plus théorique de l'algorithme précédemment décrit. Comme les isomorphismes $\psi$ et $\psi^{-1}$ ont une expression explicite, on peut décrire précisément les projections $\proj_{\pi} : \blg_{n} \to \C\sym_{\pi}$, assurant la décomposition en blocs $C\sym_{\pi}$. Ensuite, l'idempotent central $f(\lambda)$ peut être exprimé en fonction du caractère irréductible $\varsigma^{\lambda}$ du groupe symétrique :
$$f(\lambda)=\frac{\dim\lambda}{n!}\,\sum_{\sigma \in \sym_{n}}\varsigma^{\lambda}(\sigma)\,\sigma\,.$$
Si $\bblambda=\{\lambda^{(1)},\ldots,\lambda^{(r)}\}$ est une multipartition, fixons un ordre arbitraire sur les partitions la composant, et notons :
\begin{align*}
&\sym_{\bblambda}=\sym_{\lambda^{(1)}}\times \cdots \times \sym_{\lambda^{(r)}}\qquad;\qquad\varsigma^{\bblambda}=\varsigma^{\lambda^{(1)}}\otimes \cdots \otimes \varsigma^{\lambda^{(r)}}\qquad;\\
&\dim {\bblambda}=\dim \lambda^{(1)}\times \cdots \times \dim \lambda^{(r)} \qquad;\qquad m(\bblambda)=(-1)^{r-1}(r-1)!\,.
\end{align*}
On adopte les mêmes notations pour une suite ordonnée de partitions $\Lambda$. La description des idempotents centraux à l'aide des caractères s'étend alors aux $f(\Lambda,\pi)$ :
$$f(\Lambda,\pi)=\frac{\dim \Lambda}{|\Lambda|!}\,\sum_{\sigma \in \sym_{\pi}} \varsigma^{\Lambda}(\sigma)\,\sigma$$
et ceux-ci réalisent les projections $\proj_{\Lambda}$. En combinant ces descriptions des projecteurs, on obtient la formule abstraite suivante pour les nombres de Hurwitz :
\begin{theorem}[Formule de Frobenius pour les nombres de Hurwitz]\label{lasthurwitz}
$$H_{n,g}(\mu)=\sum_{\bblambda \in \mathscr{Z}_{n}} c(\bblambda)^{2(n+g-1)-R(\mu)}\,\frac{m(\bblambda)\,b(\bblambda)\,\dim \bblambda}{n!\,|\bblambda|!}\left(\sum_{\sigma \in \sym_{\bblambda} \cap C_{\mu\sqcup 1^{n-|\mu|}}}\varsigma^{\bblambda}(\sigma)\right).$$
\end{theorem}
\begin{example}
Cette formule abstraite explique le développement des nombres de Hurwitz en puissances $2g$-ième des contenus des multipartitions de taille $n$. Ainsi, on obtient par exemple :
$$H_{3,g}(2)=(9^{g+1}-1)/2\qquad;\qquad H_{4,g}(3)=6^{2g+2}-3^{2g+2}\,.$$
\end{example}\bigskip

\noindent La formule du théorème \ref{lasthurwitz} n'est en réalité pas nouvelle : c'est l'application d'un principe d'inclusion-exclusion aux \textbf{nombres de Hurwitz déconnectés} $H_{n,g}^{\bullet}(\mu)$, qui comptent les factorisations non transitives (ou de fa\c con équivalente, les revêtements ramifiés non connexes), et sont donnés par la formule de Frobenius usuelle
\begin{align*}
H_{n,g}^{\bullet}(\mu^{(1)},\ldots,\mu^{(r)})&=\frac{1}{n!}\,\,\card\!\left\{(\sigma^{(1)},\ldots,\sigma^{(r)},\tau^{(1)},\ldots,\tau^{(k)}) \,\,\,\bigg|\,\,\,\substack{t(\sigma^{(i)})=\mu^{(i)}\,\,;\,\,t(\tau^{(j)})=21^{n-2}\\ \sigma^{(1)}\cdots\sigma^{(r)}\tau^{(1)}\cdots\tau^{(k)}=\id_{\lle1,n\rre}}\right\}\\
&\!\!\!\!\!\!=\frac{1}{n!}\left(\prod_{i=1}^{r}\card C_{\mu^{(i)} \sqcup1^{n-|\mu^{(i)}|}}\right)\binom{n}{2}^{k}\,\sum_{\lambda \in \ym_{n}} \frac{(\dim \lambda)^{2}}{n!}\,\left(\frac{c(\lambda)}{\dim \lambda}\right)^{k}\,\left(\prod_{i=1}^{r}\chi^{\lambda}(\mu^{(i)})\right).
\end{align*}
voir par exemple \cite[appendice]{LZ04}. Néanmoins, l'algèbre des permutations scindées a permis de manipuler rigoureusement les symétries des énumérations, et elle a fourni un algorithme efficace pour le calcul explicite.
\bigskip
\bigskip

\section{Conjectures sur les produits de classes de conjugaison dans $Z(\C\GL(n,\For_{q}))$}\label{polyobs}

Les constructions effectuées au début du chapitre dans le contexte très général des fibrés de semi-groupes ont été réalisées en vu de conjectures sur les produits de classes de conjugaison dans les algèbres $\C\GL(n,\For_{q})$. Rappelons que ces classes de conjugaison sont indexées par les polypartitions de taille $n$, c'est-à-dire les fonctions $\bbmu : M /\Gal \to \ym$ qui associent des partitions $\bbmu(\phi)$ à des $\For_{q}$-polynômes irréductibles $\phi \in M/\Gal$, de telle fa\c con que
$$\|\bbmu\|=\sum_{\phi \in M/\Gal} (\deg \phi)\,|\bbmu(\phi)|=n\,,$$
voir la section \ref{jordanfrobenius}. Si $\bbmu$ est une polypartition telle que $\|\bbmu\|+\ell(\bbmu(X-1))\leq n$, on note $\bbmu \rightarrow n$ la \textbf{polypartition complétée} définie par 
\begin{align*}(\bbmu\rightarrow n)(\phi \neq X-1)&=\bbmu(\phi)\\
(\bbmu\rightarrow n)(X-1)&=\bbmu(X-1)\rightarrow \left(n-\!\!\sum_{\phi \neq X-1} (\deg\phi)\,|\bbmu(\phi)|\right).
\end{align*}
en utilisant les notations du début de la section \ref{ivanovkerov}. Ainsi, si l'on regarde la $\For_{3}$-polypartition $\bbmu=\{X^{2}+1:(2,1)\,;\,(X-1):(1,1)\}$, alors $\bbmu \rightarrow n$ fait sens si $n \geq 8+2=10$, et par exemple,
$$\bbmu\rightarrow 15=\{X^{2}+1:(2,1)\,;\,(X-1):(2,2,1,1,1,1,1)\}\,.$$
La complétion définit une bijection entre les polypartitions de taille $n$ et les polypartitions $\bbmu$ telles que $\|\bbmu\|+\ell(\bbmu(X-1)) \leq n$. D'autre part, remarquons que si $g$ est un isomorphisme de type $\bbmu \rightarrow n$ (au sens de la section \ref{jordanfrobenius}), alors la codimension 
$$\codim \mathrm{Fix}(g)=n-\dim\{x \in (\For_{q})^{n}\,\,|\,\, g(x)=x\}$$
est exactement $\|\bbmu\|$. En effet, on peut supposer à conjugaison près que $g=J_{\bbmu \rightarrow n}$ est une matrice compagnon de Jordan, et si l'on écrit les blocs de cette matrice sous la forme
$$\begin{pmatrix}
J_{\bbmu(\phi_{1})}&          &                               &                                     &          &                                  &                    \\
                              &\ddots &                               &                                     &          &                                  &                    \\
                              &           &J_{\bbmu(\phi_{r})} &                                    &          &                                  &                    \\
                              &           &                               & J_{(X-1)^{\mu_{1}+1}} &          &                                  &                    \\
                              &           &                               &                                     &\ddots&                                  &                     \\
                              &           &                               &                                     &          &J_{(X-1)^{\mu_{l}+1}}&                     \\
                              &           &                               &                                     &          &                                  &I_{n-\|\bbmu\|-\ell(\bbmu(X-1))}
\end{pmatrix}$$
avec $\mu=\bbmu(X-1)$, alors le sous-espace propre $\mathrm{Fix}(g)$ est exactement celui engendré par les $n-\|\bbmu\|-\ell(\mu)$ derniers vecteurs de base, plus chaque premier vecteur d'un bloc $J_{(X-1)^{\mu_{i}+1}}$, d'où :
$$\dim \mathrm{Fix}(g)=n-\|\bbmu\|-\ell(\mu)+\ell(\mu)=n-\|\bbmu\|\,.$$
Or, étant donnés deux isomorphismes $g$ et $h$, on a bien sûr $\mathrm{Fix}(gh)\subset \mathrm{Fix}(g)\cap \mathrm{Fix}(h)$, donc, au niveau des codimensions,
$$\codim \mathrm{Fix}(gh) \leq \codim \big(\mathrm{Fix}(g) \cap \mathrm{Fix}(h)\big) \leq \codim \mathrm{Fix}(g)+\codim \mathrm{Fix}(h)\,.$$
On en déduit le caractère filtré des centres des algèbres $\C\GL(n,\For_{q})$ :
\begin{proposition}[Filtration sur les centres $Z(\C\GL(n,\For_{q}))$]
Pour tout entier $n$ et toutes polypartitions $\bblambda$ et $\bbmu$, on a dans $Z(\C\GL(n,\For_{q}))$ une identité du type 
$$C_{\bblambda \rightarrow n}\,*\,C_{\bbmu \rightarrow n} = \sum_{\|\bbnu\|\leq \|\bblambda\|+\|\bbmu\|} a_{\bblambda\bbmu}^{\bbnu}(n,q)\,C_{\bbnu \rightarrow n}\,$$
avec les $a_{\bblambda\bbmu}^{\bbnu}(n,q)$ dans $\N$, et étant entendu que la classe $C_{\bblambda \rightarrow n}$ vaut $0$ si la polypartition $\bblambda$ satisfait $\|\bblambda\|+\ell(\bblambda(X-1))>n$.
\end{proposition}
\begin{proof}
Le terme de gauche est une combinaison linéaire (à coefficients dans $\N$) de produits d'isomorphismes $gh$ avec 
$$\codim \mathrm{Fix}(gh)\leq \codim \mathrm{Fix}(g)+\codim \mathrm{Fix}(h)= \|\bblambda\|+\|\bbmu\|\,.$$ 
D'autre part, comme $Z(\C\GL(n,\For_{q}))$ est stable par produit et engendré linéairement par les classes de conjugaison, ce produit est une combinaison linéaire (à coefficients dans $\N$) de classes de conjugaison $C_{\bbnu \rightarrow n}$. D'après ce qui précède, tous les $\bbnu$ intervenant satisfont $\|\bbnu\| \leq \|\bblambda\|+\|\bbmu\|$.
\end{proof}\bigskip

Des expériences numériques sur des cas simples --- par exemple, un produit de classes de degré $1$ sur $\For_{3}$ ou $\For_{5}$ --- laissent penser que la conjecture suivante est vraie :
\begin{conjecture}[Dépendance en $n$ et $q$ des constantes de structure dans les centres des algèbres $\C\GL(n,\For_{q})$]\label{lastconjecture}
Les coefficients de structure $a_{\bblambda\bbmu}^{\bbnu}(n,q)$ sont toujours des éléments de $\Q(q^{n},q)$. 
\end{conjecture} 
\noindent Plus précisément, il n'est pas difficile de se convaincre (au moins sur des cas simples) que ces coefficients correspondent à des énumérations d'arrangements des sous-espaces cycliques des éléments des classes $C_{\bblambda \rightarrow n}$ et $C_{\bbmu \rightarrow n}$ ; comme les classes de conjugaison (voir la proposition \ref{cardclass}) et les arrangements d'hyperplans vérifiant certaines conditions d'incidence sont énumérés par des fractions rationnelles\footnote{Comme la caractéristique $q$ du corps de base est fixée, la seule vraie variable est en fait $q^{n}$, donc on pourrait énoncer la conjecture avec $\Q[q^{n}]$ à la place de $\Q(q^{n},q)$.} en $q$ et $q^{n}$, ceci explique intuitivement la conjecture. Une preuve rigoureuse nécessiterait néanmoins des arguments à la Ivanov-Kerov. Ainsi, s'il existait une algèbre graduée $\alg_{\infty}(\For_{q})$ engendrée linéairement par des classes $A_{\bbmu}$ avec $\bbmu \in \ym(\For_{q})$, et se projetant sur tous les centres $Z(\C\GL(n,\For_{q}))$ avec 
$$\proj_{n}(A_{\bbmu})=f_{\bbmu}(q^{n},q)\,C_{\bbmu \rightarrow n}\quad;\quad \deg A_{\bbmu}=\|\bbmu\|\,,$$
alors on aurait une preuve agréable du caractère rationnel en $q$ et $q^{n}$ des $a_{\bblambda\bbmu}^{\bbnu}(n,q)$. La fin de ce chapitre est consacré à la construction d'une algèbre $\blg_{\infty}(\For_{q})$ qui se projette sur tous les $\C\GL(n,\For_{q})$. On espère qu'une sous-algèbre de $\blg_{\infty}(\For_{q})$ permette d'établir le théorème \ref{lastconjecture}, mais ceci reste encore à l'état de conjecture.\bigskip
\bigskip

Par rapport au cas des permutations partielles, une difficulté majeure est qu'il n'est pas possible de définir un semi-groupe des <<~isomorphismes partiels~>> de taille $n$ dont les éléments seraient les couples $(g,V)$ avec $V$ sous-espace de $(\For_{q})^{n}$ et $g \in \GL(V)$. En effet, étant donnés deux tels couples $(g,V)$ et $(h,W)$, on aimerait alors définir leur produit comme un isomorphisme partiel de l'espace $V+W$, mais contrairement au cas $q=1$, il n'y a pas de fa\c con canonique de compléter $g$ et $h$ en isomorphismes de cet espace --- il faudrait spécifier des supplémentaires de $V$ et de $W$ dans $V+W$ sur lesquels $g$ ou $h$ induisent l'identité. Compte tenu de cette difficulté, il convient de se raccrocher aux principes énoncés dans les deux premières sections du chapitre. Ainsi, pour construire une limite projective des algèbres $\C\GL(n,\For_{q})$, il suffit de trouver des semi-treillis $L_{0}\subset L_{1}\subset \cdots \subset L_{n}\subset \cdots$ fibrant les $\GL(n,\For_{q})$. Les \textbf{décompositions en sommes directes de sous-espaces} vont fournir de tels semi-treillis. \'Etant fixés un entier $n$ et une puissance $q$ d'un nombre premier, on appelle décomposition de $(\For_{q})^{n}$ la donnée d'une famille (non ordonnée) $\{V^{(1)},\ldots,V^{(r)})$ de $\For_{q}$-sous-espaces vectoriels tels que 
$$(\For_{q})^{n}=V^{(1)}\oplus V^{(2)}\oplus \cdots \oplus V^{(r)}\,.$$
Le \textbf{type}, ou profil d'une décomposition $\Delta$ est la partition $t(\Delta)$ obtenue en ordonnant les dimensions des espaces $V^{(i)}$ composant $\Delta$. On notera $\dec(n,\For_{q})$ l'ensemble des décompositions de $(\For_{q})^{n}$, et si $\lambda \in \ym_{n}$, $\dec(\lambda,\For_{q})$ l'ensemble des décompositions de type $\lambda$. Le groupe $\GL(n,\For_{q})$ agit sur $\dec(n,\For_{q})$ par
$$g\cdot\{V^{(1)},\ldots,V^{(r)}\}=\big\{g(V^{(1)}),\ldots,g(V^{(r)})\big\}\,,$$
et les orbites de cette action sont exactement les sous-ensembles $\dec(\lambda,\For_{q})$. Le stabilisateur d'une décomposition de type $\lambda=(\lambda_{1},\ldots,\lambda_{r})$ a pour cardinal 
$$\left(\prod_{j \geq 1}m_{j}(\lambda)!\right)\,\,\left(\prod_{i=1}^{r} \card \GL(\lambda_{i},\For_{q})\right),$$
et par conséquent, le cardinal de $\dec(n,\For_{q})$ s'écrit :
$$\card \dec(n,\For_{q})=\sum_{\lambda \in \ym_{n}}\card\dec(\lambda,\For_{q})=\sum_{\lambda \in \ym_{n}}\frac{1}{\prod_{j\geq 1}m_{j}(\lambda)!}\,\,\frac{\card \GL(n,\For_{q})}{\prod_{i=1}^{\ell(\lambda)} \card \GL(\lambda_{i},\For_{q})}\,.$$
On dira qu'une décomposition $(V^{(i)})_{i\in I}$ est plus fine qu'une décomposition $(W^{(j)})_{j \in J}$ s'il existe une partition d'ensemble $I=\bigsqcup_{j \in J}I_{j}$ telle que 
$$\forall j\in J,\,\,\,W^{(j)}=\bigoplus_{i \in I_{j}} V^{(i)}\,.$$
\begin{proposition}[Semi-treillis des décompositions en sommes directes]\label{semilatticedecomposition}
Muni de cette relation, $\dec(n,\For_{q})$ est un semi-treillis d'élément maximal la décomposition triviale $\{(\For_{q})^{n}\}$.
\end{proposition}
\begin{proof}
Dans ce qui suit, $\Delta_{U}=\{U^{(1)},\ldots,U^{(r)}\}$ et $\Delta_{V}=\{V^{(1)},\ldots,V^{(s)}\}$ sont deux décompositions dans $\dec(n,\For_{q})$, et si $\pi$ est une partition d'ensemble dans $\mathfrak{Q}_{r}$, on note $\pi(\Delta_{U})$ la décomposition plus grossière que $\Delta_{U}$ obtenue en regroupant les espaces $U^{(i)}$ dans chaque part de $\pi$. Par exemple, supposons que
$$\Delta_{U}=\For_{q}[e_{1}]\oplus \For_{q}[e_{2},e_{3}]\oplus\For_{q}[e_{4}]=U^{(1)}\oplus U^{(2)}\oplus U^{(3)}$$
dans $(\For_{q})^{4}$, où $e_{i}$ désigne le $i$-ième vecteur de la base canonique, et où $\For_{q}[(x_{i})_{i\in I}]$ désigne l'espace vectoriel engendré par la famille $(x_{i})_{i \in I}$. Alors, si $\pi=\{1,3\}\sqcup \{2\}$, 
$$\pi(\Delta_{U})=(U^{(1)}\oplus U^{(3)})\oplus U^{(2)}=\For_{q}[e_{1},e_{4}]\oplus \For_{q}[e_{2},e_{3}]\,.$$
Soit $S(\Delta_{U})\simeq \mathfrak{Q}_{r}$ l'ensemble des décompositions de $(\For_{q})^{n}$ plus grossières que $\Delta_{U}$, et $S(\Delta_{U},\Delta_{V})$ l'ensemble des décompositions de $(\For_{q})^{n}$ plus grossières que $\Delta_{U}$ et que $\Delta_{V}$. Notons que $S(\Delta_{U},\Delta_{V})$ est une partie non vide de $S(\Delta_{U})$, car $\{(\For_{q})^{n}\}$ est une décomposition plus grossière que toute autre, donc est dans $S(\Delta_{U},\Delta_{V})$. Fixons deux décompositions $\pi(\Delta_{U})$ et $\nu(\Delta_{U})$ dans $S(\Delta_{U},\Delta_{V})$. Tout sous-espace $V^{(j)}$ est contenu dans un certain $\bigoplus_{i \in \pi_{a_{j}}}U^{(i)}$, et dans un $\bigoplus_{i \in \nu_{b_{j}}}U^{(i)}$. Par conséquent\footnote{L'identité pour l'intersection des sous-espaces est vraie compte tenu de la décomposition en somme directe $(\For_{q})^{n}=\bigoplus_{i=1}^{r}U^{(i)}$.},
$$V_{j} \subset \left(\bigoplus_{i \in \pi_{a_{j}}}U^{(i)}\right) \cap\left(\bigoplus_{i \in \nu_{b_{j}}}U^{(i)}\right)=\bigoplus_{i \in \pi_{a_{j}} \cap \nu_{b_{j}}} U^{(i)}\,.$$
Ceci implique le point suivant: si $\mu=\pi \wedge \nu$ est la borne inférieure de $\pi$ et $\nu$ dans $\mathfrak{Q}_{r}$, alors tout sous-espace $V^{(j)}$ est contenu dans un sous-espace de la décomposition $\mu(\Delta_{U})$, donc $\mu(\Delta_{U})$ appartient à $S(\Delta_{U},\Delta_{V})$. Ainsi, si l'on identifie $S(\Delta_{U},\Delta_{V})$ à une partie de $\mathfrak{Q}_{r}$, alors $S(\Delta_{U},\Delta_{V})$ est stable pour l'opération $\wedge$ sur les partitions d'ensembles, et on peut considérer :
$$\pi_{\min}=\bigwedge \{\pi'\,\,|\,\,\pi'(\Delta_{U}) \in S(\Delta_{U},\Delta_{V})\}\,.$$
La décomposition $\pi_{\min}(\Delta_{U})$ est plus grossière que $\Delta_{U}$ et $\Delta_{V}$, et par construction, c'est la décomposition en sous-espaces la plus fine ayant cette propriété, car toute autre décomposition $\pi(\Delta_{U}) \in S(\Delta_{U},\Delta_{V})$ vérifie $\pi_{\min} \preceq \pi$ et donc $\pi_{\min}(\Delta_{U}) \preceq \pi(\Delta_{U})$. Ainsi, la borne supérieure $\Delta_{U} \vee \Delta_{V}=\pi_{\min}(\Delta_{U})$ existe toujours.
\end{proof}
\noindent La preuve qui précède implique également le point suivant : si $\Delta_{U}\preceq \Delta_{W}=\pi(\Delta_{U})$, alors la fonction de M\"obius $\mu(\Delta_{U},\Delta_{W})$ est simplement la fonction de M\"obius du treillis $\mathfrak{Q}_{r}$ évaluée en les partitions d'ensembles $\{1\}\sqcup \{2\}\sqcup \cdots \sqcup \{r\}$ et $\pi$. D'autre part, notons que le semi-treillis $\dec(n,\For_{q})$ n'est pas un treillis : par exemple, si $\mathcal{B}=(e_{i})_{i \in \lle 1,n\rre}$ et $\mathcal{B}'=(e_{i}')_{i \in \lle 1,n\rre}$ sont deux bases de $(\For_{q})^{n}$, alors les deux décompositions
$$\For_{q}[e_{1}]\oplus \cdots \oplus \For_{q}[e_{n}]\quad\text{et}\quad\For_{q}[e_{1}']\oplus\cdots \oplus \For_{q}[e_{n}']$$
sont minimales pour l'ordre de raffinement, mais elles n'ont de borne inférieure que si elles sont égales (à permutation des vecteurs près).\bigskip
\bigskip

On appelle \textbf{isomorphisme scindé} de taille $n$ et sur $\For_{q}$ la donnée d'un isomorphisme $g \in \GL(n,\For_{q})$ et d'une décomposition $\{U^{(1)},\ldots,U^{(r)}\}=\Delta \in \dec(n,\For_{q})$ laissée stable par $g$, c'est-à-dire que $g(U^{(i)})=U^{(i)}$ pour tout $i$. Autrement dit, un isomorphisme scindé est la donnée d'une décomposition et d'isomorphismes $g^{(i)}\in \GL(U^{(i)})$ pour tout $i$ ; on a alors $g=g^{(1)}\oplus \cdots \oplus g^{(r)}$. Les isomorphismes scindés s'agrègent en un fibré de semi-groupe $\IS(n,\For_{q})=\GL(n,\For_{q})\times_{I} \dec(n,\For_{q})$, l'application de fibration étant :
$$I_{g}=\left\{\text{décompositions de $(\For_{q})^{n}$ laissées stables par $g$}\right\}\,.$$
Il est trivial de vérifier que $I$ vérifie bien les axiomes d'une fibration : en effet, si $g$ et $h$ sont des isomorphismes et si $\Delta_{U} \in I_{g}$ et $\Delta_{V}\in I_{h}$, alors $g$ et $h$ laissent tous les deux stable la décomposition $\Delta_{U}\vee \Delta_{V}$, donc $g\circ h$ également, ce qui qui implique $I_{g} \vee I_{h} \subset I_{gh}$.\bigskip\bigskip

Dans ce qui suit, un isomorphisme scindé $(g,\Delta_{U})$ sera aussi noté sous la forme additive $(g^{(1)},U^{(1)})\oplus\cdots\oplus (g^{(r)},U^{(r)})$. Une notation particulière sera employée lorsque $\Delta_{U}$ est la décomposition standard associée à une composition $c=(c_{1},\ldots,c_{r})$ de $\lle 1,n\rre$, c'est-à-dire que $\Delta_{U}=\{U^{(1)},\ldots,U^{(r)}\}$ avec 
$$U^{(i)}=\For_{q}[e_{c_{1}+\cdots+c_{i-1}+1},\ldots,e_{c_{1}+\cdots+c_{i}}]\,.$$
Dans ce cas, nous noterons simplement $(g^{(1)},U^{(1)})\oplus \cdots \oplus (g^{(r)},U^{(r)})=g^{(1)}\oplus \cdots \oplus g^{(r)}$.
\begin{example}
Dans $\IS(5,\For_{3})$,  $(\begin{smallmatrix}1 & 2 \\ 0 & 2 \end{smallmatrix}) \oplus (1) \oplus (\begin{smallmatrix} 1 & 1 \\ 2 & 1\end{smallmatrix})$ et $\left(\begin{smallmatrix}  0 & 2 & 0 \\ 1 & 0 & 1 \\ 2 & 2 & 1\end{smallmatrix}\right)\oplus (2) \oplus (1)$ ont pour produit
$$ \begin{pmatrix}2 & 2 & 2 \\
2 & 0 & 2 \\
2 & 2 & 1 \end{pmatrix} \oplus \begin{pmatrix} 2 & 1 \\ 1 & 1\end{pmatrix}.$$
\end{example}\bigskip

Il reste maintenant à relier entre eux les semi-treillis $\dec(n,\For_{q})$ et considérer une limite projective $\GL(\infty,\For_{q}) \times_{I}\dec(\infty,\For_{q})$. Si $N \geq n$ et $\Delta_{U}=\big\{(\For_{q})^{n}=U_{1}\oplus \cdots \oplus U_{r}\big\}$ est une décomposition dans $\dec(n,\For_{q})$, on peut toujours la compléter en une décomposition $\Delta_{U}\rightarrow N$ de $\dec(N,\For_{q})$ en rajoutant les sous-espaces vectoriels $\For_{q}[e_{n+1}],\ldots,\For_{q}[e_{N}]$ :
$$\Delta_{U} \rightarrow N = \big\{(\For_{q})^{N}=U_{1}\oplus \cdots \oplus U_{r}\oplus \For_{q}[e_{n+1}]\oplus \cdots \oplus \For_{q}[e_{N}]\big\}\,.$$
On appellera décomposition de $(\For_{q})^{\infty}$ la donnée d'une famille $(\Delta_{N})_{N \geq n}$ de décompositions des $((\For_{q})^{N})_{N \geq n}$ telles que
$$\Delta_{N} = \Delta_{n} \rightarrow N$$
pour tout $N \geq n$, étant entendu qu'on identifie deux telles familles si elles sont égales à partir d'un certain rang $N$. Autrement dit, une décomposition de $(\For_{q})^{\infty}$ est la donnée d'une décomposition d'un certain $(\For_{q})^{n}$, étant entendu qu'on identifie deux décompositions de tailles distinctes si l'une peut être obtenue à partir de l'autre par le procédé de complétion décrit ci-dessus. Nous noterons $\dec(\infty,\For_{q})$ l'ensemble des décompositions de $(\For_{q})^{\infty}$.
\begin{example}
Les décompositions $\For_{q}[e_{1},e_{3}]\oplus \For_{q}[e_{2},e_{4}]$ et $\For_{q}[e_{1},e_{3}]\oplus \For_{q}[e_{2},e_{4}]\oplus \For_{q}[e_{5}]\oplus \For_{q}[e_{6}]$ correspondent à la même décomposition de $(\For_{q})^{\infty}$, car celle de taille $6$ est obtenue par complétion à partir de celle de taille $4$.
\end{example}
\bigskip

On munit $\dec(\infty,\For_{q})$ de la relation d'ordre suivante : une décomposition $\Delta_{1} \in \dec(n_{1},\For_{q})$ est plus petite qu'une décomposition $\Delta_{2} \in \dec(n_{2},\For_{q})$ si $(\Delta_{1}\rightarrow N) \preceq (\Delta_{2}\rightarrow N)$ dans $\dec(N,\For_{q})$ pour $N$ assez grand (plus grand que $\max(n_{1},n_{2})$). On vérifie sans mal que la preuve de la proposition \ref{semilatticedecomposition} s'adapte à cette définition ; ainsi, $(\dec(\infty,\For_{q}),\preceq)$ est un semi-treillis (mais cette fois-ci sans élément maximal). Une chaîne exhaustive dans $\dec(\infty,\For_{q})$ est constituée par les décompositions $\Delta_{n\geq 1}=\For_{q}[e_{1},e_{2},\ldots,e_{n}]$. En effet, on a bien $\Delta_{N}\succeq \Delta_{n}$ si $N \geq n$, car 
$$\For_{q}[e_{1},\ldots,e_{N}]\succeq \For_{q}[e_{1},\ldots,e_{n}]\oplus \For_{q}[e_{n+1}]\oplus \cdots \oplus \For_{q}[e_{N}]$$
dans $\dec(N,\For_{q})$. D'autre part, si $\Delta$ est une décomposition de $(\For_{q})^{\infty}$ donnée par une décomposition dans $\dec(n,\For_{q})$, alors $\Delta \preceq \Delta_{n}$, car $\Delta_{n}$ est maximale dans $\dec(n,\For_{q})$. On est donc dans la situation décrite\footnote{Il n'y a pas d'élément minimal $l_{0}$, mais c'est sans importance.} dans la section \ref{chain}, avec
$$L=\dec(\infty,\For_{q}) \qquad;\qquad L_{n}=\dec(n,\For_{q})=\{\Delta \preceq \Delta_{n}\}\,.$$
Considérons maintenant la limite directe $\GL(\infty,\For_{q})$ des groupes $\GL(n,\For_{q})$ pour les morphismes
\begin{align*}
i_{n \leq N} : \GL(n,\For_{q}) &\to \GL(N,\For_{q})\\
g &\mapsto (g \rightarrow N )= \begin{pmatrix}g & 0 \\
0 & I_{N-n}\end{pmatrix} ;
\end{align*}
un élément de $\GL(\infty,\For_{q})$ est donc donné par un $g \in \GL(n,\For_{q})$, étant entendu qu'on identifie deux isomorphismes si l'un est obtenu à partir de l'autre en rajoutant des $1$ sur la diagonale de la matrice. On définit comme précédemment une fibration $I : \GL(\infty,\For_{q}) \to \mathcal{I}(\dec(\infty,\For_{q}))$ par :
$$I_{g}=\left\{\text{décompositions de $(\For_{q})^{\infty}$ laissées stables par $g$}\right\}\,.$$
On obtient ainsi un fibré de semi-groupes $\IS(\infty,\For_{q})=\GL(\infty,\For_{q})\times_{I}\dec(\infty,\For_{q})$ dont les éléments sont les isomorphismes scindés de taille arbitraire, étant entendu qu'on identifie deux isomorphismes scindés si l'un peut être obtenu à partir de l'autre par l'opération
$$(h,\Theta)=(g,\Delta)\oplus(1,\For_{q}[e_{n+1}])\oplus(1,\For_{q}[e_{n+2}])\oplus\cdots\oplus(1,\For_{q}[e_{N}])=(g\rightarrow N,\Delta \rightarrow N)\,.$$
On est donc dans la situation du paragraphe \ref{chain}, avec $M=\GL(\infty,\For_{q})$ et $N_{n}=\GL(n,\For_{q})$. \bigskip\bigskip

Si $(g,\Delta)=(g^{(1)},U^{(1)})\oplus \cdots \oplus (g^{(r)},U^{(r)})$ est un isomorphisme scindé de taille arbitraire $n$, notons $\deg (g,\Delta)$ la quantité $n-\dim \mathrm{Fix}(g,\Delta)$, où $\dim \mathrm{Fix}(g,\Delta)$ est la dimension du sous-espace vectoriel 
$$\mathrm{Fix}(g,\Delta)=\mathrm{Fix}(g)=\{x \in (\For_{q})^{n}\,\,|\,\, g(x)=x\}\,.$$ Cette définition est compatible avec la complétion des isomorphismes scindés, c'est-à-dire que $\deg (g,\Delta) \rightarrow N=\deg (g,\Delta)$ pour tout $N \geq n$. Par conséquent, le degré est bien défini sur $\IS(\infty,\For_{q})$ et sur tous les $\IS(n,\For_{q})$. De plus, si $(g_{1},\Delta_{1})$ et $(g_{2},\Delta_{2})$ sont deux isomorphismes scindés, alors
$$\deg((g_{1},\Delta_{1})(g_{2},\Delta_{2}))\leq \deg(g_{1},\Delta_{1})+\deg(g_{2},\Delta_{2})$$
pour les mêmes raisons que dans le cas des isomorphismes simples. Ainsi, $\deg$ s'étend en des filtrations d'algèbres sur les $\blg_{n}(\For_{q})=\C\IS(n,\For_{q})$. On note $\blg_{\infty}(\For_{q})$ la limite projective des $\blg_{n}(\For_{q})$ dans la catégorie des algèbres filtrées et relativement aux morphismes
\begin{align*}\phi_{N,n} : \blg_{N}(\For_{q}) & \to\blg_{n}(\For_{q})\\
(g,\Delta)&\mapsto \begin{cases} (h,\Theta) & \text{si }(g,\Delta)=(h,\Theta)\rightarrow N,\\
0&\text{sinon}.
\end{cases}
\end{align*}
Notons que les projections $\phi_{N,n}$ sont légèrement plus restrictives que celles décrites dans la section \ref{chain} : ainsi, on peut avoir $\phi_{N,n}(g,\Delta)=0$ même si $\Delta \preceq \Delta_{n}$, par exemple si 
$$(g,\Delta)=(h,\Delta_{n})\oplus (\lambda,\For_{q}[e_{n+1}])\oplus \cdots \oplus (\lambda,\For_{q}[e_{N}])$$
avec $\lambda \in \For_{q}^{\times}$ différent de $1$. Néanmoins, la construction fournit de nouveau une algèbre naturelle $\blg_{\infty}(\For_{q})\subset \C[[\IS(\infty,\For_{q})]]$ se projetant sur toutes les algèbres de groupes $\C\GL(n,\For_{q})$.
\begin{conjecture}[Isomorphismes scindés et théorème de Farahat-Higman pour les groupes $\GL(n,\For_{q})$]\label{verylastconjecture}
Une sous-algèbre de $\blg_{\infty}(\For_{q})$, ou d'une algèbre construite sur le même principe, se projette sur tous les centres $Z(\C\GL(n,\For_{q}))$, et permet une démonstration agréable de la conjecture \ref{lastconjecture}.
\end{conjecture}
\bigskip

Pour tout entier $n$, il existe une action naturelle de $\GL(n,\For_{q})$ sur $\dec(n,\For_{q})$ qui est compatible avec la fibration décrite précédemment, d'où une action de $\GL(n,\For_{q})$ sur le fibré $\IS(n,\For_{q})$ :
$$g\cdot(h,\Delta)=(ghg^{-1},g(\Delta))$$
avec $g(\Delta=\{U^{(1)},\ldots,U^{(r)}\})=\{g(U^{(1)}),\ldots,g(U^{(r)})\}$. Deux isomorphismes scindés $(g,\Delta_{U})$ et $(h,\Delta_{V})$ sont conjugués sous cette action si et seulement si $\Delta_{U}$ et $\Delta_{V}$ ont même profil, et s'il existe une bijection $U^{(i)}\leftrightarrow V^{(i)}$ préservant les dimensions et telle que $g^{(i)}$ et $h^{(i)}$ ait même type $\bbmu^{(i)}$ au sens de la section \ref{jordanfrobenius}. Les classes de $\GL(n,\For_{q})$-conjugaison d'isomorphismes scindés de taille $n$ sont donc indexées par les \textbf{polypartitions multiples} de taille $n$, c'est-à-dire, les familles $\{\bbmu^{(1)},\ldots,\bbmu^{(r)}\}$ de $\For_{q}$-polypartitions telles que $\sum_{i=1}^{r} \|\bbmu^{(i)}\|=n$. De même, il existe une action naturelle de $\GL(\infty,\For_{q})$ sur le fibré $\IS(\infty,\For_{q})$, et qui est compatible avec le degré des isomorphismes scindés. Ainsi, il est naturel de considérer les algèbres d'invariants :
$$\alg_{n}(\For_{q})=(\blg_{n}(n,\For_{q}))^{\GL(n,\For_{q})}\qquad;\qquad \alg_{\infty}(\For_{q})=(\blg_{\infty}(\For_{q}))^{\GL(\infty,\For_{q})}\,.$$
La première se projette sur $Z(\C\GL(n,\For_{q}))$ et a une base $(B_{M,n})_{M}$ indexée par les polypartitions multiples $M=\{\bbmu^{(1)},\ldots,\bbmu^{(r)}\}$ de taille $n$ ; la seconde est la limite projective des $\alg_{n}(\For_{q})$, et elle a une base $(B_{M})_{M}$ indexée par les polypartitions multiples, à identification près par un procédé de complétion des polypartitions multiples semblable à ce que l'on a exposé au début de la section pour les polypartitions simples. Malheureusement, même si les projetés $\pi_{n}(B_{M,n})$ sont bien des multiples des classes $C_{\bbmu^{(1)}\sqcup \cdots \sqcup \bbmu^{(r)}}$ (avec des coefficients rationnels en $q$ et $q^{n}$), il n'est pas du tout clair que dans l'algèbre $\alg_{\infty}(\For_{q})$, le produit de deux classe $B_{M_{1}}$ et $B_{M_{2}}$ ne mette en jeu que des classes $B_{N}$ avec $\|N\|\leq \|M_{1}\|+\|M_{2}\|$. D'autre part, il serait largement préférable de manipuler des algèbres dont les bases naturelles sont indexées par les polypartitions simples (au lieu des polypartitions multiples). Il manque donc encore quelques arguments à la preuve de la conjecture \ref{verylastconjecture}.

\addtocontents{lof}{\protect\bigskip}

\pagestyle{empty}

\clearpage
~

\backmatter
\pagestyle{empty}

\renewcommand{\thechapter}{\Alph{chapter}}

\chapter{Conclusion et perspectives}
\setcounter{chapter}{3}
\setcounter{figure}{0}
\setcounter{footnote}{0}
En utilisant exclusivement des techniques d'observables de diagrammes, nous avons déterminé  le comportement asymptotique des partitions aléatoires tirées suivant les mesures de Plancherel, les $q$-mesures de Plancherel de type A ou B, les mesures d'induction parabolique, les mesures de Gelfand et les mesures de Schur-Weyl de paramètre $\alpha>1/2$, $\alpha=1/2$ ou $\alpha<1/2$. Nos résultats reposent en grande partie sur la combinatoire de l'algèbre $\obs$, qui peut être interprétée comme sous-algèbre des invariants dans l'algèbre des permutations partielles. Cette dernière construction rentre dans le cadre plus général des fibrés de semi-groupes par des semi-treillis, et on conjecture qu'une telle construction dans le contexte des groupes linéaires finis $\GL(n,\For_{q})$ fournisse une algèbre naturelle d'observables de polydiagrammes. Ceci nous mène aux problèmes (non résolus) suivants :\vspace{1.5mm}
\begin{enumerate}
\item Quels sont les semi-treillis donnant lieu à une algèbre d'Ivanov-Kerov pour les groupes $\GL(n,\For_{q})$ ? À défaut, peut-on construire une algèbre de Farahat-Higman $FH(\GL(\For_{q}))$, et quelles sont ces propriétés ?\vspace{2mm}
\item Soit $X^{\lambda}$ le caractère unipotent normalisé associé au $\GL(n,\For_{q})$-module $U^{\lambda}(q)$. On fixe une polypartition $\bbmu$, et on tire au hasard $X^{\lambda}$ suivant la $q$-mesure de Plancherel $M_{n,q}$. Quelle est la distribution asymptotique de $X^{\lambda}(g)$, où $g$ est un isomorphisme de type $\bbmu \sqcup (\mathbb{1}:1^{n-|\bbmu|})$ ?\vspace{2mm}
\item Soit $\chi^{\bblambda}$ le caractère normalisé associé au $\GL(n,\For_{q})$-module irréductible de type $\bblambda$. On tire au hasard $\chi^{\bblambda}$ suivant la $\GL(n,\For_{q})$-mesure de Plancherel évoquée à la fin du chapitre \ref{general}. Quelle est la distribution asymptotique de la variable aléatoire centrée $\chi^{\bblambda}(g)$, où $g$ est un isomorphisme de type $\bbmu \sqcup (\mathbb{1}:1^{n-|\bbmu|})$ avec $\bbmu$ polypartition fixée ? \vspace{1.5mm}
\end{enumerate}
Il est vraisemblable qu'après renormalisation, ces variables aléatoires aient un comportement asymptotique gaussien, mais ce type de résultat dépend d'une meilleure compréhension des caractères et des classes de conjugaison des groupes $\GL(n,\For_{q})$. On peut évidemment envisager les mêmes problèmes en type B, et dans ce cadre, le premier problème est :\vspace{1.5mm}
\begin{enumerate}
\setcounter{enumi}{3}
\item Comment démontrer la conjecture \ref{bqmix} ? Autrement dit, quelle est la $q$-formule de Frobenius pour les caractères des algèbres d'Hecke de type B ?\vspace{1.5mm}
\end{enumerate}
Notons qu'il existe de nombreuses preuves de la formule de Ram \ref{qfrobeniusschur} en type A, mais aucune ne semble évidente à généraliser en type B.\bigskip\medskip

Les méthodes d'observables de diagrammes sont essentiellement des méthodes de moments dans un contexte non commutatif ; en particulier, les espérances $\esper[p_{\mu}(\lambda)]$ ou $\esper[\varSigma_{\mu}(\lambda)]$ pour des modèles de partitions aléatoires sont à rapprocher des espérances $\esper[\prod_{i=1}^{r}\tr M^{\mu_{i}}]$ pour les modèles de matrices aléatoires. D'autre part, nous avons expliqué dans les chapitres \ref{determinantal} et \ref{matrix} que les méthodes de processus ponctuels déterminantaux utilisées en théorie des matrices aléatoires avaient également des applications en théorie des partitions aléatoires. Plus généralement, à peu près tous les outils naturels de la théorie des matrices aléatoires ont des analogues pour les modèles de partitions. Dès lors, il convient de se demander ce que donnent ces outils pour les nouveaux modèles présentés dans cette thèse, ou pour d'autres modèles. Mentionnons en particulier les trois problèmes suivants :\vspace{2mm}
\begin{enumerate}\setcounter{enumi}{4}
\item De nombreux résultats asymptotiques sur les grandes matrices aléatoires ont un caractère universel, c'est-à-dire qu'ils sont valables dans un contexte beaucoup plus général que le cadre hermitien gaussien usuel (voir les travaux récents \cite{Joh09,EPRSY09,ESY09}). Pour les partitions aléatoires, nous avons vu que le théorème central limite de Kerov avait cette même propriété. Dans quel cadre peut-on encore étendre ce résultat ? A-t-on d'autres propriétés universelles pour les grandes partitions aléatoires ?\vspace{2mm}
\item L'équivalence de Baik-Deift-Johansson entre le GUE et la mesure de Plancherel est également valable entre le GOE et la mesure de Gelfand, voir \cite{BR01}. Plus généralement, a-t-on une équivalence asymptotique entre les $\beta$-ensembles et les $\beta$-mesures de Plancherel ? Quelles interprétations combinatoires ou géométriques peut-on donner à une telle correspondance ?\vspace{2mm}
\item Les mesures spectrales des matrices aléatoires hermitiennes gaussiennes satisfont à un principe de grandes déviations, voir en particulier \cite{BAG97,Gui04}, et \cite{Var08} pour des généralités sur les grandes déviations. De plus, la fonction de taux sous-jacente met en jeu l'entropie non-commutative de Voiculescu des mesures de probabilité. Dans la preuve originale de la loi des grands nombres de Logan-Shepp-Kerov-Vershik (\cite{LS77}), Logan et Shepp construisent une fonction <<~énergie~>> sur les diagrammes de Young continus, qui met en jeu un analogue de cette entropie non-commutative, et telle que la forme limite minimise l'énergie. Cette énergie est-elle une fonction de taux pour un principe de grandes déviations vérifiées par les mesures de Plancherel ? A-t-on des analogues de ce principe variationnel pour d'autres modèles de partitions aléatoires ?\vspace{2mm}
\end{enumerate}
Le dernier problème peut être formulé pour une famille de déformations des mesures de Plancherel des groupes symétriques, qui est liée aux marches aléatoires sur ces groupes ; ceci mène à une généralisation de la définition \ref{defplancherel} que nous souhaitions donner dans cette conclusion. Si $G$ est un groupe fini et si $V=\bigoplus_{\lambda \in \widehat{G}} n_{\lambda}\,V^{\lambda}$ est une représentation réductible de $G$, sa mesure de Plancherel $\proba_{V}$ peut s'écrire :
$$\proba_{V}[\lambda]=\frac{n_{\lambda}\,\dim V^{\lambda}}{\dim V}=\frac{n_{\lambda}\,\varsigma^{\lambda}(e_{G})}{\varsigma^{V}(e_{G})}$$
où $\varsigma^{V}$ est le caractère non normalisé de $V$. Déformons alors $e_{G}$ et considérons un élément central $c$ suffisamment proche de $e_{G}$ dans $Z(\C G)$. On peut lui associer une nouvelle fonction 
$$\proba_{V,c}[\lambda]=\frac{n_{\lambda}\,\varsigma^{\lambda}(c)}{\varsigma^{V}(c)}$$
qui est bien définie si $c$ est suffisamment proche de $e_{G}$, car le dénominateur ne s'annule pas dans ce cas. Cette fonction $\proba_{V,c}$ a pour somme $1$ sur $\widehat{G}$ ; de plus, si l'algèbre $\Z G$ est scindée sur $\R$ (ce qui est le cas pour les groupes symétriques) et si $c$ est une combinaison linéaire à coefficients réels de classes de conjugaison, alors $\proba_{V,c}$ est bien une mesure de probabilité (réelle positive) sur $\widehat{G}$. \bigskip

Un cas particulièrement intéressant est celui où $c=\E^{tC}$ avec $C \in Z(\R G)$, étant entendu que $G$ est un groupe tel que $\Z G$ se scinde sur $\R$. Alors, on obtient une famille à un paramètre de déformations de la mesure de Plancherel de $V$ :
$$\proba_{V,C,t}[\lambda]=\frac{n_{\lambda}\,\varsigma^{\lambda}(\E^{tC})}{\varsigma^{V}(\E^{tC})}=\frac{n_{\lambda}\,\dim V^{\lambda}\,\exp(t\chi^{\lambda}(C))}{\sum_{\lambda \in \widehat{G}} n_{\lambda}\,\dim V^{\lambda}\,\exp(t\chi^{\lambda}(C))}$$
et ces mesures sont directement reliées au processus markovien sur $G$ de loi invariante par conjugaison, et de matrice de transition $C$. Dans le cas des groupes symétriques, on peut étudier ces mesures lorsque 
$$V=\C \sym_{n} \qquad ;\qquad  C=C_{21^{n-2}}-\frac{n(n-1)}{2} \,C_{1^{n}}\qquad;\qquad t \simeq \frac{\beta}{\sqrt{n}}\,, $$
ce problème étant lié à la marche aléatoire isotrope sur le graphe de Cayley de $\sym_{n}$ par rapport au système de générateurs donné par les transpositions. Les formes limites des diagrammes aléatoires renormalisés sont dans ce contexte des déformations de la courbe de Logan-Shepp-Kerov-Vershik, et il semble possible d'énoncer dans ce cadre des principes de grandes déviations (le cas $\beta=0$ correpondant aux mesures de Plancherel usuelles).\bigskip

\figcapt{\vspace{-2.5mm}
\include{french1}\vspace{-5mm}}{Diagramme de Young aléatoire de taille $n=1000$ sous la mesure de Plancherel déformée de paramètre $\beta=1/2$.}{Diagramme sous la mesure de Plancherel déformée de paramètre $\beta=1/2$}
Finalement, il est tentant de vouloir généraliser les techniques d'observables de diagram\-mes à d'autres types d'objets aléatoires planaires, par exemple les arbres planaires ou les partitions planes. L'idée générale est d'interpréter les caractères irréductibles de certaines algèbres\footnote{Par exemple, si l'on remplace dans ce mémoire $\sym_{n}$ ou $\GL(n,\For_{q})$ par $(\Z/2\Z)^{n}$, ceci mène aux résultats asymptotiques classiques liant marches aléatoires et mouvements browniens unidimensionnels.} comme observables de ces objets planaires, étant entendu que les espérances de ces variables aléatoires sous les mesures de Plancherel de représentations peuvent parfois être calculées très simplement. Ainsi, il serait très intéressant de retrouver les résultats asymptotiques connus relatifs aux partitions planes ou aux arbres planaires aléatoires par ces <<~méthodes de moments non commutatifs~>>.

\chapter*{Conclusion and open problems}
\setcounter{footnote}{0}
\foreignlanguage{english}{By using exclusively techniques of observables of diagrams, we have determined the asymptotic behaviour of partitions picked randomly according to Plancherel measures, $q$-Plancherel measures of type A or B, measures of parabolic induction, Gelfand measures and Schur-Weyl measures of parameter $\alpha>1/2$, $\alpha=1/2$ or $\alpha<1/2$. Our results mostly relie on the combinatorics of the algebra $\obs$, which can be interpreted as the subalgebra of invariants in the algebra of partial permutations. This latter construction sits in the more general setting of semilattice bundles over semigroups, and we conjecture that by using such a construction in the context of finite general linear groups $\GL(n,\For_{q})$, we could exhibit a natural algebra of observables of polydiagrams. This leads to the following questions:\vspace{2mm}
\begin{enumerate}
\item Which semilattices give rise to an Ivanov-Kerov algebra for the groups $\GL(n,\For_{q})$? Other\-wise, is it possible to construct a Farahat-Higman algebra $FH(\GL(\For_{q}))$, and what are the properties of this algebra?\vspace{2mm}
\item Let $X^{\lambda}$ be the normalized character of the unipotent $\GL(n,\For_{q})$-module $U^{\lambda}(q)$. We fix a polypartition $\bbmu$, and we choose $X^{\lambda}$ randomly according to the $q$-Plancherel measure $M_{n,q}$. What is the asymptotic distribution of $X^{\lambda}(g)$, where $g$ is an isomorphism of type $\bbmu \sqcup (\mathbb{1} : 1^{n-|\bbmu|})$?  \vspace{2mm}
\item Let $\chi^{\bblambda}$ be the normalized character of the irreducible $\GL(n,\For_{q})$-module of type $\bblambda$. We choose $\chi^{\bblambda}$ randomly according to the $\GL(n,\For_{q})$-Plancherel measure evoked at the end of Chapter \ref{general}. What is the asymptotic distribution of the centered random variable $\chi^{\bblambda}(g)$, where $g$ is an isomorphism of type  $\bbmu \sqcup (\mathbb{1} : 1^{n-|\bbmu|})$ with $\bbmu$ fixed polypartition?\vspace{2mm}
\end{enumerate}
It is likely that after renormalization, these random variables exhibit a gaussian asymptotic behaviour, but this kind of result depends on a better understanding of the characters and conjugacy classes of the groups $\GL(n,\For_{q})$. One can consider the same problems in type B, and in this setting, the first question to tackle is:\vspace{2mm}
\begin{enumerate}
\setcounter{enumi}{3}
\item How to prove Conjecture \ref{bqmix}? In other words, what is the $q$-Frobenius formula for the characters of the Hecke algebras of type B? \vspace{2mm}
\end{enumerate}
Notice that there are many proofs of Ram's formula \ref{qfrobeniusschur} in type A, but none of them seems easy to adapt in type B.}\bigskip
\bigskip

\foreignlanguage{english}{The techniques of observables of diagrams are in essence methods of moments in a noncommutative setting; in particular, the expectations $\esper[p_{\mu}(\lambda)]$ or $\esper[\varSigma_{\mu}(\lambda)]$ for models of random partitions bear some resemblance to the expectations $\esper[\prod_{i=1}^{r}\tr M^{\mu_{i}}]$ for models of random matrices. On the other hand, we have explained in Chapters \ref{determinantal} and \ref{matrix} that the methods of determinantal point processes used in random matrix theory can also be used in random partition theory. More generally, almost all the natural tools of random matrix theory have analogues for models of partitions. Therefore, it is worth asking what these tools can give for the new models presented in this thesis, or for other models. Let us mention in particular the three following problems:\vspace{2mm}
\begin{enumerate}
\setcounter{enumi}{4}
\item Many asymptotic results for large random matrices are universal, meaning that they hold in a much broader context than the usual gaussian hermitian setting (see the recent works \cite{Joh09,EPRSY09,ESY09}). For random partitions, we have seen that Kerov's central limit theorem has this same property. In which setting can we further expand this result? What are the other universal properties of large random partitions? \vspace{2mm}
\item The Baik-Deift-Johansson equivalence between the GUE and the Plancherel measure holds also between the GOE and the Gelfand measure, see \cite{BR01}. More generally, is there an asymptotic equivalence between the $\beta$-ensembles and the $\beta$-Plancherel measures? How can it be interpreted in combinatoric or geometric terms?\vspace{2mm}
\item The spectral measures of large random matrices satisfy a large deviation principle, see in particular \cite{BAG97,Gui04}, and \cite{Var08} for generalities on large deviations. Moreover, the underlying rate function involves Voiculescu's noncommutative entropy of probability measures. In the original proof of Logan-Shepp-Kerov-Vershik law of large numbers (\cite{LS77}), Logan and Shepp construct an ``energy'' of continuous Young diagrams, which involves an analogue of that noncommutative entropy, and such that the limit shape minimizes the energy. Is this energy a rate function for a large deviation principle satisfied by the Plancherel measures? Does this variational principle have analogues for other models of random partitions?   \vspace{2mm}
\end{enumerate}
The latter problem can be formulated for a family of deformations of the Plancherel measures of the symmetric groups, that is related to the random walks on these groups; it leads to a generalization of Definition \ref{defplancherel} that we wanted to present in this conclusion. If $G$ is a finite group and if $V=\bigoplus_{\lambda \in \widehat{G}}n_{\lambda}\,V^{\lambda}$ is a reducible representation of $G$, its Plancherel measure can be written as:
$$\proba_{V}[\lambda]=\frac{n_{\lambda}\,\dim V^{\lambda}}{\dim V}=\frac{n_{\lambda}\,\varsigma^{\lambda}(e_{G})}{\varsigma^{V}(e_{G})}$$
where $\varsigma^{V}$ is the non normalized character of $V$. Then, let us deform $e_{G}$ and consider a central element $c$ that is sufficiently close to $e_{G}$ in $Z(\C G)$. One can associate to $c$ a new function
$$\proba_{V,c}[\lambda]=\frac{n_{\lambda}\,\varsigma^{\lambda}(c)}{\varsigma^{V}(c)}$$
that is well-defined if $c$ is sufficiently close to $e_{G}$, because in that case the denominator does not vanish. This function $\proba_{V,c}$ sums up to $1$ on $\widehat{G}$; moreover, if  the algebra $\Z G$ is split over $\R$ (this is the case for symmetric groups) and if $c$ is a real linear combination of conjugacy classes, then $\proba_{V,c}$ is indeed a (real non-negative) probability measure on $\widehat{G}$.}\bigskip\bigskip

\foreignlanguage{english}{A very interesting case is when $c=\E^{tC}$ with $C \in \Z(\R G)$, provided that $G$ is a group such that $\Z G$ splits over $\R$. Then, one obtains a one-parameter family of deformations of the Plancherel measure of $V$:
$$\proba_{V,C,t}[\lambda]=\frac{n_{\lambda}\,\varsigma^{\lambda}(\E^{tC})}{\varsigma^{V}(\E^{tC})}=\frac{n_{\lambda}\,\dim V^{\lambda}\,\exp(t\chi^{\lambda}(C))}{\sum_{\lambda \in \widehat{G}} n_{\lambda}\,\dim V^{\lambda}\,\exp(t\chi^{\lambda}(C))}$$
and these measures are directly related to the markovian process on $G$ whose law is invariant by conjugation, and with transition matrix $C$. In the case of symmetric groups, one can study these measures when
$$V=\C \sym_{n} \qquad ;\qquad  C=C_{21^{n-2}}-\frac{n(n-1)}{2} \,C_{1^{n}}\qquad;\qquad t \simeq \frac{\beta}{\sqrt{n}}\,, $$
this problem being linked to the isotropic random walk on the Cayley graph of $\sym_{n}$ with respect to the system of generators given by the transpositions. In this context, the limit shapes of renormalized random diagrams are some deformations of Logan-Shepp-Kerov-Vershik curve, and it seems possible to enounce in this setting some large deviation principles (the case when $\beta=0$ corresponds to the usual Plancherel measures). }\bigskip

\begin{minipage}{160mm}{
\begin{center}
\footnotesize{\include{english1}}
\end{center}}
 \end{minipage} 

\normalsize
\noindent {\textsc{Figure} -- Random Young diagram of size $n=1000$ under the deformed Plancherel measure of parameter $\beta=1/2$.}\bigskip\bigskip

Finally, it is tempting to generalize the techniques of observables of diagrams to the case of other random planar objects, for instance planar trees of plane partitions. The general idea is to interpret the irreducible characters of some algebras\footnote{For instance, if one replaces in this thesis $\sym_{n}$ or $\GL(n,\For_{q})$ by $(\Z/2\Z)^{n}$, it leads to the classical asymptotic results that relate the one-dimensional random walks and brownian motions.} as observables of these planar objects, provided that the expectations of these variables under Plancherel measures of representations are sometimes very easy to compute. Thus, it would be very interesting to retrieve the known asymptotic results relative to random plane partitions or random planar trees by using these ``methods of noncommutative moments''.

\addtocontents{toc}{\protect\contentsline {chapter}{Bibliographie}{275}{part.B}
}

\bibliographystyle{alpha}
\bibliography{source}

\end{document}